\documentclass[12pt]{amsart}
\pdfoutput=1
\usepackage{amsmath,amsfonts,amssymb,amsthm,ifpdf,mathtools}
\usepackage{graphicx,color}
\usepackage{imakeidx}
\makeindex[title=Index Notation,columns=3]

\ifpdf
\usepackage[pdftex,pdfstartview=FitH,pdfborderstyle={/S/B/W 1}, %
colorlinks=true, linkcolor=blue, urlcolor=blue,citecolor=blue,%
pagebackref=true]{hyperref}
\else
\usepackage[hypertex,colorlinks]{hyperref}
\fi

\usepackage{bookmark}

\DeclareMathRadical{\sqrtsign}{symbols}{"70}{largesymbols}{"70}
\newcommand{\bb}{\mathbb}
\newcommand{\gothic}{\mathfrak}

\newcommand{\cx}{{\bb C}}
\newcommand{\half}{{\bb H}}
\newcommand{\integers}{{\bb Z}}
\newcommand{\natls}{{\bb N}}

\newcommand{\reals}{{\bb R}}
\newcommand{\proj}{{\bb P}}
        
\newlength{\figboxwidth}             
\newcommand{\makefig}[3]{
        \begin{figure}[htb]
        \refstepcounter{figure}
        \label{#2}
        \begin{center}
                #3~\\
                \smallskip
                Figure \thefigure.  #1
        \end{center}
        \medskip
        \end{figure}
}
\newcommand{\makefignocenter}[3]{
        \begin{figure}[htb]
        \refstepcounter{figure}
        \label{#2}
                #3~\\
                \smallskip

                Figure \thefigure.  #1
        \medskip
        \end{figure}
}
\newcommand{\starredsubsection}[1]{
\renewcommand{\thesubsection}{\arabic{section}.\arabic{subsection}$^\ast$}
\subsection{#1}
\renewcommand{\thesubsection}{\arabic{section}.\arabic{subsection}}
}

\renewcommand{\bold}[1]{\medskip \noindent {\bf #1 }\nopagebreak}

\newcommand{\dirsum}{\oplus}

\newcommand{\isom}{\cong}

\newcommand{\tensor}{\otimes}
\newcommand{\bs}{\backslash}
\newcommand{\cross}{\times}
\renewcommand{\Re}{\mbox{Re}\:}
\renewcommand{\Im}{\mbox{Im}\:}

\newcommand{\st}{\;\: : \;\:}         

\newcommand{\zed}{\integers}

\newcommand{\diag}{\operatorname{diag}}
\newcommand{\area}{\operatorname{area}}

\newcommand{\diam}{\operatorname{diam}}

\newcommand{\Hom}{\operatorname{Hom}}

\newcommand{\Vol}{\operatorname{Vol}}

\newcommand{\Span}{\operatorname{span}}
\newcommand{\GSpc}{\operatorname{GSpc}}
\newcommand{\Lie}{\operatorname{Lie}}
\newcommand{\Conf}{\operatorname{Conf}}

\def\@ifundefined#1#2#3%
  {\expandafter\ifx\csname#1\endcsname\relax#2\else#3\fi}

\@ifundefined{theoremstyle}{
}{
\theoremstyle{plain} 
}
\newtheorem{theorem}{Theorem}[section]
\newtheorem{prop}[theorem]{Proposition}
\newtheorem{proposition}[theorem]{Proposition}
\newtheorem{lemma}[theorem]{Lemma}

\newtheorem{corollary}[theorem]{Corollary}
\newtheorem{claim}[theorem]{Claim}

\@ifundefined{theoremstyle}{
}{
\theoremstyle{definition} 
}
\newtheorem{definition}[theorem]{Definition}

\newtheorem{remark}[theorem]{Remark}

\newcommand{\cA}{{\mathcal A}}
\newcommand{\cB}{{\mathcal B}}
\newcommand{\cC}{{\mathcal C}}
\newcommand{\cD}{{\mathcal D}}
\newcommand{\cE}{{\mathcal E}}
\newcommand{\cF}{{\mathcal F}}
\newcommand{\cG}{{\mathcal G}}
\newcommand{\cH}{{\mathcal H}}

\newcommand{\cK}{{\mathcal K}}
\newcommand{\cL}{{\mathcal L}}
\newcommand{\cM}{{\mathcal M}}
\newcommand{\cN}{{\mathcal N}}
\newcommand{\cO}{{\mathcal O}}
\newcommand{\cP}{{\mathcal P}}
\newcommand{\cQ}{{\mathcal Q}}

\newcommand{\cS}{{\mathcal S}}

\newcommand{\cU}{{\mathcal U}}
\newcommand{\cV}{{\mathcal V}}
\newcommand{\cW}{{\mathcal W}}

\newcommand{\cY}{{\mathcal Y}}

\mathchardef\GG="321D
\newcommand{\gB}{{\gothic B}}

\newcommand{\gF}{{\gothic F}}

\newcommand{\gP}{{\gothic P}}

\newcommand{\bA}{{\mathbf A}}
\newcommand{\bB}{{\mathbf B}}

\newcommand{\bE}{{\mathbf E}}
\newcommand{\bF}{{\mathbf F}}

\newcommand{\bH}{{\mathbf H}}

\newcommand{\bL}{{\mathbf L}}

\newcommand{\bP}{{\mathbf P}}

\newcommand{\bS}{{\mathbf S}}

\newcommand{\bV}{{\mathbf V}}

\newcommand{\bZ}{{\mathbf Z}}

\newcommand{\bfe}{{\mathbf e}}
\newcommand{\bff}{{\mathbf f}}

\newcommand{\bfi}{{\mathbf i}}
\newcommand{\bfj}{{\mathbf j}}

\newcommand{\bfu}{{\mathbf u}}
\newcommand{\bfv}{{\mathbf v}}
\newcommand{\bfw}{{\mathbf w}}

\DeclarePairedDelimiter{\abs}{|}{|}
\DeclarePairedDelimiter{\norm}{\|}{\|}

\newcommand{\Nbhd}{\operatorname{Nbhd}}

\newcommand{\mcc}[1]{{}} 
\newcommand{\mccc}[1]{{}} 

\definecolor{mygreen}{rgb}{0,.5,0}

\numberwithin{equation}{section}
\newcommand{\noz}{n}

\title[Invariant and stationary measures]
{Invariant and stationary measures for the $SL(2,\reals)$ action on
Moduli space}

\author{Alex Eskin}
\thanks{Research  of  the first author is partially supported  by
NSF grants DMS 0604251, DMS 0905912 and DMS 1201422}
\address{
Department of Mathematics,
University of Chicago,
Chicago, Illinois 60637, USA\\
}
\email{eskin@math.uchicago.edu}

\author{Maryam Mirzakhani}
\thanks{Research of the second author is 
partially supported by the Clay foundation and by NSF grant DMS 0804136}
\address{Department of Mathematics, 
Stanford University, Stanford CA 94305 USA \\}
\email{mmirzakh@math.stanford.edu}

\makeindex

\begin{document}

\begin{abstract}
We prove some ergodic-theoretic rigidity properties of the action of
$SL(2,\reals)$ on moduli space.  
In particular,  we show that any ergodic measure
invariant under the action of the upper triangular subgroup of
$SL(2,\reals)$ is supported on an invariant affine submanifold.

The main theorems are inspired by the
results of several authors on unipotent flows on homogeneous spaces,
and in particular by Ratner's seminal work. 
\end{abstract}

\maketitle

\vspace{-3.1pt}
\setcounter{tocdepth}{2}
\makeatletter
\def\l@subsection{\@tocline{2}{0pt}{2.5pc}{5pc}{}}
\makeatother
\tableofcontents


\section{Introduction}
\label{sec:intro}

Suppose $g \ge 1$, \index{$g$} and 
let $\alpha = (\alpha_1,\dots, \alpha_\noz)$ \index{$\alpha$} 
be a partition of $2g-2$, and 
let $\cH(\alpha)$ \index{$H($@$\cH(\alpha)$} be a stratum of Abelian differentials,
i.e.\ the space of pairs $(M,\omega)$ where $M$ \index{$M$} is a Riemann surface
and $\omega$ \index{$\omega$} is a holomorphic $1$-form on $M$ whose
zeroes have
multiplicities $\alpha_1 \dots \alpha_\noz$. The form $\omega$ defines a
canonical flat metric on $M$ with conical singularities at the zeros
of $\omega$. Thus we refer to points of $\cH(\alpha)$ as
{\em flat surfaces} or {\em translation surfaces}. For an introduction
to this subject, see the survey \cite{Zorich:survey}. 

The space $\cH(\alpha)$
admits an action of the group $SL(2,\reals)$ which generalizes the
action of $SL(2,\reals)$ on the space $GL(2,\reals)/SL(2,\zed)$ of flat
tori. In this paper we prove ergodic-theoretic rigidity
properties of this action. 

In what follows, we always replace $\cH(\alpha)$ by a finite cover
$X_0$ which is a manifold. Such a cover can be found by
e.g.\ considering a level $3$ structure (see \S\ref{sec:semi:markov}). 
However, in the introduction,
we suppress this from the notation.

Let $\Sigma \subset M$ \index{$\Sigma$} 
denote the set of zeroes of $\omega$. Let
$\{ \gamma_1, \dots, \gamma_k\}$ denote a symplectic $\zed$-basis for
the relative 
homology group $H_1(M,\Sigma, \zed)$. We can define a map $\Phi:
\cH(\alpha) \to \cx^k$ by 
\begin{displaymath}
\Phi(M,\omega) = \left( \int_{\gamma_1} \omega, \dots, \int_{\gamma_k}
  \omega \right).
\end{displaymath}
The map $\Phi$ (which depends on a choice of the basis $\{ \gamma_1,
\dots, \gamma_k\}$) is a local coordinate system on $(M,\omega)$.  
Alternatively,
we may think of the cohomology class $[\omega] \in H^1(M,\Sigma, \cx)$ as a
local coordinate on the stratum $\cH(\alpha)$. We
will call these coordinates {\em period coordinates}.

We can consider the measure $\lambda$ on $\cH(\alpha)$ which
is given by the pullback of the Lebesgue measure on $H^1(M,\Sigma,
\cx) \approx \cx^k$. The measure $\lambda$ is independent of the
choice of basis  
$\{ \gamma_1, \dots, \gamma_k\}$, and is easily seen to be
$SL(2,\reals)$-invariant. We call $\lambda$ the {\em Lebesgue} or the
{\em Masur-Veech} measure
on $\cH(\alpha)$. 

The area of a translation surface is given by 
\begin{displaymath}
a(M,\omega) = \frac{i}{2} \int_M \omega \wedge \bar{\omega}.
\end{displaymath}
A ``unit hyperboloid'' $\cH_1(\alpha)$ \index{$H$@$\cH_1(\alpha)$}
is defined as a subset of translation surfaces in $\cH(\alpha)$ of
area one. The $SL(2,\reals)$-invariant Lebesgue measure $\lambda_{(1)}$ on
$\cH_1(\alpha)$ is defined by disintegration of the
Lebesgue measure $\lambda$ on $\cH_1(\alpha)$, namely
\mccc{FIX ME}
$$
d\lambda=c \, d\lambda_{(1)} \, da.
$$
where $c$ is a constant. 
A fundamental result of Masur \cite{masur:interval} and Veech
\cite{veech:gauss}  is that 
$\lambda_{(1)}(\cH_1(\alpha)) < \infty$. In this paper, we
normalize $\lambda_{(1)}$ so that $\lambda_{(1)}(\cH_1(\alpha)) =
1$ (and so $\lambda_{(1)}$ is a probability measure). 

For a subset $\cM_1 \subset \cH_1(\alpha)$ we write
\begin{displaymath}
\reals \cM_1 = \{ (M, t \omega) \;|\; (M,\omega) \in \cM_1, \quad t \in
\reals \setminus \{0 \} \} \subset \cH(\alpha).
\end{displaymath}
\begin{definition}
\label{def:affine:measure}
An ergodic $SL(2,\reals)$-invariant probability measure $\nu_1$ on
$\cH_1(\alpha)$ is called {\em affine} if the following conditions hold:
\begin{itemize}
\item[{\rm (i)}] The support $\cM_1$ of $\nu_1$ is an 
{\em immersed submanifold} of
  $\cH_1(\alpha)$, i.e.\
there exists a manifold $\cN$ and a proper continuous
  map $f: \cN \to \cH_1(\alpha)$ so that $\cM_1 =
  f(\cN)$. The self-intersection set of $\cM_1$, i.e.\ the set of
  points of $\cM_1$ which do not have a 
  unique preimage under $f$, is a closed subset of $\cM_1$ of
  $\nu_1$-measure $0$.   Furthermore, each point in $\cN$ has a
  neighborhood $U$ such 
  that locally $\reals f(U)$ is given by a complex linear subspace defined
  over $\reals$ in the period coordinates.
\item[{\rm (ii)}] Let $\nu$ be the measure supported on $\cM = \reals
  \cM_1$ so that $d\nu = 
    d\nu_1 da$. Then each point in $\cN$ has a
      neighborhood $U$ such that the restriction of $\nu$ to $\reals
      f(U)$ is an affine  linear measure in the period
    coordinates on $\reals f(U)$, i.e.\ it is (up to normalization) the
    induced measure of the Lebesgue measure $\lambda$ to the subspace
    $\reals f(U)$.
\end{itemize}
\end{definition}

\begin{definition}
\label{def:affine:invariant:submanifold}
We say that any suborbifold $\cM_1$ for which there exists a measure
$\nu_1$ such that the pair $(\cM_1, \nu_1)$ 
satisfies (i) and (ii) is an {\em affine invariant submanifold}. 
\end{definition}
We also consider the entire stratum $\cH(\alpha)$ to be an
(improper) affine invariant submanifold. 
It follows from 
\cite[Theorem~2.2]{Eskin:Mirzakhani:Mohammadi} that
the self-intersection set of an affine invariant manifold is itself a
finite union of 
affine invariant manifolds of lower dimension. 

For many applications we need the following:
\begin{proposition}
\label{prop:countable}
Any stratum $\cH_1(\alpha)$ contains at most
countably many affine invariant submanifolds.
\end{proposition}

Proposition~\ref{prop:countable} is deduced as a consequence of some
isolation theorems in \cite{Eskin:Mirzakhani:Mohammadi}. 
This argument relies on
adapting some ideas of G.A. Margulis to the Teichm\"uller space
setting. Another proof is given by A.~Wright in
  \cite{Wright:field}, where it is proven that affine invariant
  submanifolds are always defined over a number field.

The classification of the affine invariant submanifolds is complete in
genus $2$ by the work of McMullen \cite{McMullen:Billiards}
\cite{McMullen:geodesics} \cite{McMullen:discriminant}
\cite{McMullen:decagon} \cite{McMullen:torsion}
and Calta \cite{Calta:thesis}. 
In genus $3$ or greater it is an important open problem. See 
\cite{Moeller:Variations},
\cite{Moeller:Periodic},
\cite{Moeller:Finiteness},
\cite{Moeller:Linear}, 
\cite{Bouw:Moeller},
\cite{Bainbridge:Moeller},
\cite{Hubert:Lanneau:Moeller:Ratner},
\cite{Lanneau:Nguyen:Teichmuller:curves},
\cite{Lanneau:Nguyen:periodicity},
\cite{Lanneau:Nguyen:GL2R},
\cite{Wright:field},
\cite{Wright:cylinder},
\cite{Matheus:Wright},
\cite{Nguyen:Wright},
\cite{Aulicino:Nguyen:Wright},
\cite{Filip:semisimplicity}
and
\cite{Filip:splitting}
for some results in this direction. 

\subsection{The main theorems}
\label{sec:theorems}
Let \index{$N$} \index{$A$} \index{$N$@$\bar{N}$}
\begin{displaymath}
N = \left\{ \begin{pmatrix} 1 & t \\ 0 & 1 \end{pmatrix}, t \in \reals
  \right\}, \quad A =
\left\{ \begin{pmatrix} e^t & 0 \\ 0 & e^{-t} \end{pmatrix}, t \in
  \reals \right\}, \quad \bar{N} = \left\{ \begin{pmatrix} 1 & 0 \\ t & 1
  \end{pmatrix}, t \in \reals   \right\}
\end{displaymath}
Let $r_\theta = \begin{pmatrix} \cos \theta & \sin \theta
  \\ - \sin \theta & \cos \theta \end{pmatrix}$, \index{$r_\theta$}
and let $SO(2) = \{
r_\theta \mid \theta \in [0,2\pi) \}$. \index{$SO(2)$}
Then $N$, $\bar{N}$, $A$ and $SO(2)$ 
are subgroups of $SL(2,\reals)$. Let $P = AN$
\index{$P$} denote the set of upper triangular matrices of determinant
$1$, which is a subgroup of $SL(2,\reals)$.

\begin{theorem}
\label{theorem:P:measures}
Let $\nu$ be any ergodic $P$-invariant probability measure on
$\cH_1(\alpha)$. Then $\nu$ is $SL(2,\reals)$-invariant and affine.  
\end{theorem}

The following (which uses Theorem~\ref{theorem:P:measures}) is joint
work with A.~Mohammadi and is proved
in \cite{Eskin:Mirzakhani:Mohammadi}:
\begin{theorem}
\label{theorem:closure:submanifold}
Suppose $S \in \cH_1(\alpha)$. Then, the orbit closure
$\overline{P S} = \overline{SL(2,\reals) S}$ is an affine invariant
submanifold of $\cH_1(\alpha)$. 
\end{theorem}

For the case of strata in genus $2$, the 
$SL(2,\reals)$ part of Theorem~\ref{theorem:P:measures} and 
Theorem~\ref{theorem:closure:submanifold} were
proved using a different method by Curt McMullen \cite{McMullen:SL2R}. 
\noindent
\medskip

The proof of Theorem~\ref{theorem:P:measures}
uses extensively entropy and conditional measure 
techniques developed in the context of  homogeneous
spaces (Margulis-Tomanov \cite{Margulis:Tomanov:Ratner},
Einsiedler-Katok-Lindenstrauss \cite{EKL}). Some of the
ideas came from discussions with Amir Mohammadi.  
But the main
strategy is to replace polynomial divergence by the ``exponential drift''
idea of Benoist-Quint \cite{Benoist:Quint}.

\bold{Stationary measures.}
Let $\mu$ be an $SO(2)$-invariant compactly supported measure on
$SL(2,\reals)$ which is absolutely continuous with respect to Lebesgue
measure. A measure $\nu$ on $\cH_1(\alpha)$ is called 
{\em $\mu$-stationary} if $\mu * \nu = \nu$, where
\begin{displaymath}
\mu * \nu = \int_{SL(2,\reals)} (g_* \nu) \, d\mu(g). 
\end{displaymath}

Recall that by a theorem of Furstenberg \cite{F1}, \cite{F2}, restated
as \cite[Theorem 1.4]{Nevo:Zimmer}, there exists a probability measure
$\rho$ on $SL(2,\reals)$ such that $\nu \to \rho \ast \nu$ is a
bijection between ergodic $P$-invariant measures and ergodic 
$\mu$-stationary measures. Therefore, 
Theorem~\ref{theorem:P:measures} implies the following:
\begin{theorem}
\label{theorem:stationary:invariant}
Any ergodic $\mu$-stationary measure on $\cH_1(\alpha)$ is
$SL(2,\reals)$-invariant and affine. 
\end{theorem}
\medskip

\bold{Counting periodic trajectories in rational billiards.}
Let $Q$ be a rational polygon, and let $N(Q,T)$ denote the number of
cylinders of periodic trajectories of length at most $T$ 
for the billiard flow on $Q$. By a theorem of H. Masur
\cite{Masur:upper} \cite{Masur:lower}, there exist
$c_1$ and $c_2$ depending on $Q$ such that for all  $t >1$, 
\begin{displaymath}
c_1 e^{2t} \le N(Q,e^t) \le c_2 e^{2t}. 
\end{displaymath}
Theorem~\ref{theorem:P:measures}  and
Proposition~\ref{prop:countable} together with some extra work (done
in \cite{Eskin:Mirzakhani:Mohammadi}) 
imply the following ``weak asymptotic
formula'' (cf.\ \cite{Athreya:Eskin:Zorich}):
\begin{theorem}
\label{theorem:weak:asymptotics}
For any rational polygon $Q$, there exists a constant $c = c(Q)$ such
that
\begin{displaymath}
\lim_{t \to \infty} \frac{1}{t} \int_0^{t} N(Q, e^{s}) e^{-2s} \, ds =
c. 
\end{displaymath}
\end{theorem}
The constant $c$ in Theorem~\ref{theorem:weak:asymptotics} 
is the Siegel-Veech constant (see \cite{Veech:Siegel},
\cite{Eskin:Masur:Zorich}) associated to the affine invariant
submanifold $\cM = \overline{SL(2,\reals)S}$ where $S$ is the flat
surface obtained by unfolding $Q$. 

It is natural to conjecture that the extra averaging on
Theorem~\ref{theorem:weak:asymptotics} is not necessary, and one has
$\lim_{t \to \infty} N(Q,e^t) e^{-2t}
= c$. This can be shown if one obtains a classification of the
measures invariant under the subgroup $N$ of $SL(2,\reals)$. Such a
result is in general beyond the reach of the current methods. However
it is known in a few very special cases, 
see \cite{Eskin:Masur:Schmoll}, \cite{Eskin:Marklof:Morris},
\cite{Calta:Wortman} and \cite{Bainbridge:L:shaped}. 
\medskip

\bold{Other applications to rational billiards.}
All the above theorems apply also to the moduli spaces of flat
surfaces with marked points. Thus one should expect applications to
the ``visibility'' and ``finite blocking'' problems in rational
polygons as in \cite{Hubert:Schmoll:Trubetzkoy}. It is likely that
many other applications are possible.

\medskip

\noindent 
{\bf Acknowledgments.}
We thank  Amir Mohammadi for many useful
discussions relating to all aspects of this project. In particular
some of the
ideas for the proof of Theorem~\ref{theorem:main:partone} came during
discussions with Amir Mohammadi. 
We also thank Vadim Kaimanovich and Emmanuel Breuillard 
for their insights into the work of Benoist and Quint and Elon
Lindenstrauss for his helpful comments. We also thank the
anonymous referee for his truly extraordinary effort and his numerous
detailed and helpful comments. The paper has vastly improved as a
result of his contribution.

\section{Outline of the paper}
\label{sec:assumptions}

\subsection{Some notes on the proofs}
\label{sec:notes:on:proofs}
The theorems of \S\ref{sec:theorems} are inspired by the results of
several authors on unipotent flows on homogeneous spaces, and in
particular by Ratner's seminal work. In particular, the 
analogues of Theorem~\ref{theorem:P:measures} and
Theorem~\ref{theorem:closure:submanifold} in homogeneous dynamics 
are due to Ratner \cite{RatnerSolv}, \cite{RatnerSS}, \cite{RatnerMeas},
\cite{RatnerEqui}. (For an introduction to these ideas, and also to
the proof by Margulis and Tomanov \cite{Margulis:Tomanov:Ratner} see the book
\cite{Morris-Ratner}.)  
The homogeneous analogue of the fact that
$P$-invariant measures are $SL(2,\reals)$-invariant is due to Mozes
\cite{Mozes:epimorphic} and is based on Ratner's work. All of these
results are based in part on the ``polynomial divergence''
of the unipotent flow on homogeneous spaces. 

However, in our setting, 
the dynamics of the unipotent flow (i.e.\ the action of $N$) on
$\cH_1(\alpha)$ is poorly understood, and plays no role 
in our proofs.  The main strategy is to
replace the ``polynomial divergence'' of unipotents by the
``exponential drift'' idea in the recent breakthrough paper by Benoist
and Quint \cite{Benoist:Quint}. 

One major difficulty is that we have no apriori control over the
Lyapunov spectrum of the geodesic flow (i.e.\ the action of $A$). 
By \cite{Avila:Viana} the
Lyapunov spectrum is simple for the case of Lebesgue
(i.e.\ Masur-Veech) measure, but for
the case of an arbitrary $P$-invariant measure this is not
always true, see e.g.\ \cite{Forni:handbook}, \cite{Forni:Matheus}. 

In order to use the Benoist-Quint exponential drift argument, we
must show that the Zariski closure (or more precisely the algebraic
hull, as defined by Zimmer \cite{ZimmerBook}) of the Kontsevich-Zorich cocycle
is semisimple. The proof proceeds in
the following steps:

\bold{Step 1.} 
We use an entropy argument inspired by the ``low entropy method'' of 
\cite{EKL} 
(using \cite{Margulis:Tomanov:Ratner} together with some ideas from 
\cite{Benoist:Quint}) to show that any $P$-invariant measure $\nu$ on
$\cH_1(\alpha)$ is in fact $SL(2,\reals)$ invariant. We also prove
Theorem~\ref{theorem:main:partone} which gives control over the
conditional measures of $\nu$. This argument
occupies \S\ref{sec:semi:markov}-\S\ref{sec:proof:theorem:partone}
and is outlined in more detail in \S\ref{sec:outline:step1}.

\bold{Step 2.} By some results of Forni (see Appendix~\ref{sec:app:forni}), 
for an $SL(2,\reals)$-invariant measure $\nu$, the absolute cohomology part
of the Kontsevich-Zorich cocycle
$A:
SL(2,\reals) \times \cH_1(\alpha) \to Sp(2g, \zed)$ is semisimple, i.e.\ has
semisimple algebraic hull. For an exact statement see
Theorem~\ref{theorem:KZ:semisimple}.

\bold{Step 3.} We pick an $SO(2)$-invariant compactly supported measure
$\mu$ on $SL(2,\reals)$ which is absolutely continuous with respect to Lebesgue
measure, and work
in the random walk setting as in \cite{F1} \cite{F2} and
\cite{Benoist:Quint}. Let $B$ denote the space of
infinite sequences $g_0, g_1, \dots, $ where $g_i \in SL(2,\reals)$.
We then have a skew product shift map $T: B \times \cH_1(\alpha) \to B
\times \cH_1(\alpha)$ as in 
\cite{Benoist:Quint}, so that $T(g_0, g_1, \dots ; x) = (g_1, g_2, \dots;
g_0^{-1} x)$. 
Then, we use (in Appendix~\ref{sec:appendixC}) a
modification of the arguments by Guivarc'h and Raugi
\cite{Guivarch:Raugi:Frontiere}, \cite{Guivarch:Raugi:Contraction}, as
presented by Goldsheid and Margulis in
\cite[\S{4-5}]{Goldsheid:Margulis}, and an argument of Zimmer (see
\cite{Zimmer:amenable:reduction} or \cite{ZimmerBook}) to prove
Theorem~\ref{theorem:semisimple:lyapunov} which states  that
the Lyapunov spectrum of $T$ is always ``semisimple'', which means
that for each $SL(2,\reals)$-irreducible component of the cocycle, there is
a $T$-equivariant non-degenerate inner product on the Lyapunov
subspaces of $T$ (or more precisely on the successive quotients of the
Lyapunov flag of $T$).  This statement is trivially true if the
Lyapunov spectrum of $T$ is simple.

\bold{Step 4.} We can now use the Benoist-Quint exponential drift
method to show that the measure $\nu$ is affine. This is done in
\S\ref{sec:random:walks}-\S\ref{sec:martingale}. 
At one point, to avoid a problem with relative
homology,  we need to use a result, Theorem~\ref{theorem:no:yeti} 
about the isometric (Forni) subspace of
the cocycle, which is proved in joint work with A.~Avila and
M.~M\"oller \cite{Avila:Eskin:Moeller:yeti}. 
\medskip

Finally, we note that the proof relies heavily on
various recurrence to compact sets results for the $SL(2,\reals)$
action, such as those of \cite{Eskin:Masur} and \cite{Jayadev:thesis}. 
All of these results originate in the ideas of Margulis and
Dani, \cite{Margulis:return}, \cite{Dani:invariant:1979},
\cite{Eskin:Margulis:Mozes:31}, \cite{Eskin:Margulis:Mozes:22}.

\subsection{Notational conventions}
\label{sec:subsec:notation}
For $t \in \reals$, 
let 
\begin{displaymath}
g_t = \begin{pmatrix} e^t & 0 \\ 0 & e^{-t} \end{pmatrix}, 
\qquad u_t
= \begin{pmatrix} 1 & t \\ 0 & 1 \\ \end{pmatrix}. \index{$g_t$} \index{$u_t$}
\end{displaymath}
Let $A = \{ g_t \st t \in \reals \}$, $N = \{ u_t \st t \in \reals
\}$. Let $P = AN$. \index{$A$} \index{$N$} \index{$P$}

Let $X_0$ \index{$X_0$} denote a finite cover of 
the stratum $\cH_1(\alpha)$ which is a manifold (see
\S\ref{sec:semi:markov}). Let $\tilde{X}_0$ denote the universal cover
of $X_0$. Let $\pi: \tilde{X}_0 \to X_0$ \index{$\pi$} denote the
natural projection map.

We will need at some point to consider a certain measurable finite
cover \index{$X$} $X$ of $X_0$. This cover will be constructed in
\S\ref{sec:subsec:the:cover:X} below. 
Let $\tilde{X}$ 
\index{$X$@$\tilde{X}$} denote the ``universal cover'' of $X$, see
\S\ref{sec:subsec:the:cover:X} for the exact definition. We abuse
notation by denoting the covering map from $\tilde{X}$ to $X$ also by
the letter \index{$\pi$} $\pi$.

If $f$ is a function on $X_0$ or $X$ 
we sometimes abuse notation by denoting $f
\circ \pi$ by $f$ and write $f(x)$ instead of $f(\pi(x))$. 
A point of $\cH(\alpha)$ is a pair $(M,\omega)$, where $M$ \index{$M$}
is a compact Riemann surface,
and $\omega$ \index{$\omega$} is a holomorphic $1$-form on $M$.
 Let $\Sigma$ \index{$\Sigma$} denote the
set of zeroes of $\omega$. The cohomology class of
$\omega$ in the relative cohomology group $H^1(M,\Sigma,\cx) \isom
H^1(M,\Sigma, \reals^2)$ is a local coordinate on $\cH(\alpha)$ (see \cite{Forni:Deviation}). 
For $x  \in \tilde{X}_0$, 
let $V(x)$ denote a subspace of $H^1(M,\Sigma,\reals^2)$. Then we
denote by \index{$V[x]$} the image of $V(x)$ under the affine exponential
  map, i.e.\ 
\begin{displaymath}
V[x] = \{ y \in \tilde{X}_0 \st y - x \in V(x) \}. 
\end{displaymath}
(For some subspaces $V$, we can define $V[x]$ for $x \in \tilde{X}$ as
well. This will be explained in \S\ref{sec:subsec:the:cover:X}. Also,
depending on the context,
we sometimes consider $V[x]$ to be a subset of $X$ or $X_0$.)

Let $p:H^1(M,\Sigma,\reals) \to H^1(M,\reals)$ \index{$p$} \index{$H^1(M,\Sigma,\reals)$}
denote the natural
map. Let \index{$H^1_\perp(x)$}
\begin{equation}
\label{eq:def:H1perp}
H^1_\perp(x) = \{ v \in H^1(M,\Sigma,\reals) \st p(\Re x) \wedge p(v)
  = p(\Im x)  \wedge p(v) = 0 \}.
\end{equation}
where we are considering the ``real part map'' $\Re$ \index{$Re$@$\Re$} 
and the ``imaginary part map'' $\Im$ \index{$Im$@$\Im$} as maps from
$H^1(M,\Sigma, \cx) \isom H^1(M,\Sigma,\reals^2)$ to
$H^1(M,\Sigma,\reals)$. 
Let \index{$W(x)$} 
\begin{displaymath}
W(x) = \reals (\Im x) \oplus H^1_\perp(x) \subset H^1(M,\Sigma,\reals),
\end{displaymath}
so that
\begin{displaymath}
W(x) = \{ v \in H^1(M,\Sigma, \reals) \st p(\Im x) \wedge p(v) = 0 \}.  
\end{displaymath}
Let $\pi^-_x: W(x)  \to
H^1(M,\Sigma,\reals)$ \index{$\pi^-_x$} denote the map (defined
  for a.e. $x \in \tilde{X}_0$)
\begin{equation}
\label{eq:def:pi:minus}
\pi^-_x( c \,\Im x + v) = c \,\Re x + v \qquad c \in \reals, v \in
  H^1_\perp(x),
\end{equation}
so that
\begin{displaymath}
\pi_x^{-} (W(x)) = \{ v \in H^1(M,\Sigma, \reals) \st p(\Re x) \wedge
p(v) = 0 \}.   
\end{displaymath}
We have $H^1(M,\Sigma,\reals^2) =
  \reals^2 \tensor H^1(M,\Sigma,\reals)$. For a subspace $V(x) \subset
  W(x)$, we write 
\begin{displaymath}
V^+(x) = (1,0) \tensor V(x), \qquad V^-(x) = (0,1) \tensor \pi_x^-(V(x)). 
\end{displaymath}
Then $W^+[x]$ and $W^-[x]$ \index{$W^+[x]$} \index{$W^-[x]$}
play the role of the unstable and stable foliations for
the action of $g_t$ on $X_0$ for $t > 0$, see Lemma~\ref{lemma:forni}.

\bold{Starred Subsections.}
Some technical proofs are relegated to subsections marked with a
star. These subsections can be skipped on first reading. The general
rule is that no statement from  a starred subsection is used in
subsequent sections.

\subsection{Outline of the proof of Step 1}
\label{sec:outline:step1}

The general strategy is based on the idea of additional invariance which
was used in the proofs of Ratner
\cite{RatnerSolv}, \cite{RatnerSS}, \cite{RatnerMeas},
\cite{RatnerEqui} and Margulis-Tomanov \cite{Margulis:Tomanov:Ratner}.

The aim of Step 1 is to prove the following:
\begin{theorem}
\label{theorem:main:partone}
Let $\nu$ \index{$\nu$} be an ergodic $P$-invariant measure on $X_0$. Then
$\nu$ is $SL(2,\reals)$-invariant. In addition, there exists an 
$SL(2,\reals)$-equivariant system of 
subspaces $\cL(x) \subset W(x)$ such that for almost
all $x$, the
conditional measures of $\nu$ along $W^+[x]$ are the Lebesgue
measures along $\cL^+[x]$, and the
conditional measures of $\nu$ along $W^-[x]$ are the Lebesgue
measures along $\cL^-[x]$. 
\end{theorem}

In the sequel, we will often refer to a (generalized) subspace $U^+[x] \subset
W^+[x]$ \index{$U^+[x]$} \index{$U^+(x)$} 
on which we already proved that the conditional measure of
$\nu$ is Lebesgue. The proof of Theorem~\ref{theorem:main:partone} 
will be by induction,
and in the beginning of the induction, $U^+[x] = N x$. 
(Note: generalized subspaces are defined in
\S\ref{sec:divergence:subspaces}). 

In this introductory subsection, let $\index{$U^+(x)$}U^+(x) \subset W^+(x)$ 
denote the subspace $\{ y - x \st y \in U^+[x] \}$. (This definition
has to be modified when we are dealing with generalized subspaces, see
\S\ref{sec:divergence:subspaces}).   

\makefig{Outline of the proof of
  Theorem~\ref{theorem:main:partone}}{fig:outline}{\includegraphics{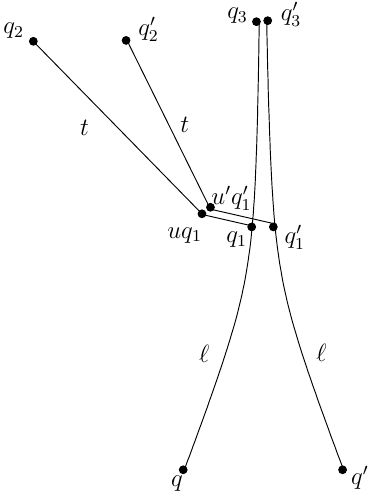}}

\bold{Outline of the proof Theorem~\ref{theorem:main:partone}.} 
Let \index{$\nu$}$\nu$ be an ergodic $P$-invariant
probability measure on $X_0$.  Since $\nu$ is $N$-invariant, the conditional
measure $\nu_{W^+}$ of $\nu$ along $W^+$ is non-trivial. This implies
that the entropy of $A$ is positive, and thus the conditional measure
$\nu_{W^-}$ of $\nu$ along $W^-$ is non-trivial (see e.g.\ \cite{Einsiedler:Lindenstrauss:Pisa}). This implies that
on a set of almost full measure, we can pick points $q$ and $q'$ in
the support of $\nu$ such that $q$ and $q'$ are in the same leaf of
$W^-$ and $d(q,q') \approx 1/100$, see Figure~\ref{fig:outline}.

Let $\ell > 0$ be a large parameter. Let $q_1 = g_\ell q$ and let
$q_1' = g_\ell q'$. Then $q_1$ and $q_1'$ are very close
together. We pick $u \in U^+(q_1)$ 
with $\|u \| \approx 1/100$, and pick (as described below) $u' \in U^+(q_1')$.
Consider the
points $u q_1$ and $u' q_1'$. With our choice of $u'$, the points 
$u q_1$ and $u'
q_1'$ will be close, but they are no longer in the same leaf of
$W^-$, and we expect them to diverge under the action of $g_t$ as $t
\to +\infty$. Let $t$ be chosen so that $q_2 = g_t u q_1$ and $q_2' =
g_t u' q_1'$ be such that $d(q_2,q_2') \approx \epsilon$, where
$\epsilon > 0$ is fixed. 

Consider the bundle (which we will denote for short $H^1$) 
whose fiber above $x \in \cH(\alpha)$ is
$H^1(M,\Sigma, \reals)$. 
The presence of the integer lattice $H^1(M,\Sigma, \zed)$ in
$H^1(M,\Sigma,\reals)$ allows one to identify the fibers at nearby
points. This defines a flat connection, called the Gauss-Manin
connection on this bundle. 

The action  of $SL(2,\reals)$ and in particular the geodesic flow $g_t$ on
$\cH(\alpha)$, extends to an action on the bundle $H^1$, where the
action on the fibers is by parallel transport with respect to the
Gauss-Manin connection. The action on the bundle takes the form
\begin{displaymath}
g_t (x, v) = (g_t x, \index{$A(g_t, v)$}A(g_t, v)),
\end{displaymath}
where $A: SL(2,\reals) \cross \cH_1(\alpha) \to GL(H^1(M,\Sigma,\reals))$
is the Kontsevich-Zorich cocycle. It is continuous (in fact locally
constant) and $\log$-integrable. Thus the multiplicative ergodic
theorem  can be applied. 

Let 
\begin{displaymath}
1 = \lambda_1(H^1) > \lambda_2(H^1) \ge \dots \ge \lambda_{k-1}(H^1) >
\lambda_{k}(H^1) = -1
\end{displaymath}
denote the Lyapunov spectrum of the Kontsevich-Zorich cocycle. (The
fact that $\lambda_2 < 1$ is due to Veech \cite{veech:gauss}
\mccc{check reference?} and Forni \cite{Forni:Deviation}). 
We have
\begin{displaymath}
H^1(M,\Sigma,\reals) = \bigoplus_{i=1}^{k} \cV_i(H^1)(x)
\end{displaymath}
where $\cV_i(H^1)(x)$  
is the Lyapunov subspace corresponding to $\lambda_i(H^1)$
(see \S\ref{sec:general:cocycle}). 
Note that $\cV_1(H^1)(x)$
corresponds to the unipotent direction inside the $SL(2,\reals)$
orbit.  In the first step of the induction, $U^+(x) = \cV_1(H^1)(x)$.

In general, for $y \in U^+[x]$, if we identify $H^1$ at $x$ and $y$
using the Gauss-Manin connection, we have (see
Lemma~\ref{lemma:properties:lyapunov:flag}), 
\begin{equation}
\label{eq:tmp:subspaces:general}
\cV_i(H^1)(y) \subset \bigoplus_{j \le i} \cV_j(H^1)(x).
\end{equation}
We say that the Lyapunov exponent $\lambda_i(H^1)$ is $U^+$-inert if
for a.e.\! $x$, $\cV_i(H^1)(x) \not\subset U^+(x)$ and also, for
a.e.\ $y \in U^+[x]$,
\begin{displaymath}
\cV_i(H^1)(y) \subset U^+(x) + \cV_i(H^1)(x).
\end{displaymath}
(In other words, $\cV_i(H^1)(x)$ is constant (modulo $U^+$) along $U^+[x]$.)
Note that in view of (\ref{eq:tmp:subspaces:general}), $\lambda_2(H^1)$ is
always $U^+$-inert. We now assume for simplicity that $\lambda_2(H^1)$ is
the only $U^+$-inert exponent. 

We may write
\begin{displaymath}
u' q_1' - u q_1 = w_+ + g_s (u q_1) + w_- 
\end{displaymath}
where $w_+ \in W^+(uq_1)$, $w_- \in W^-(uq_1)$, and $s \in \reals$.  
Furthermore, due to the assumption that $\lambda_2$ is the only inert
exponent, after possibly making a small change to $u$ and $u'$ (see
\S\ref{sec:divergence:subspaces}), 
we may write 
\begin{displaymath}
w_+ = \sum_{i=2}^n v_i
\end{displaymath}
where $v_i \in \cV_i(H^1)(u q_1)$, and furthermore, 
$\|v_2\|/\|u'q_1' - u q_1\|$ is bounded from below. Then,  
$q_2' - q_2$ will be approximately in the direction
of $\cV_2(H^1)(q_2)$, see \S\ref{sec:first:divergence} for the details.

Let $f_2(x)$
denote the conditional measure of $\nu$ along $(\cV_1+\cV_2)(H^1)[x]$. (This
conditional measure can be defined since $\nu$ is $U^+$-invariant). 
Let $q_3 = g_s q_1$ and $q_3' = g_s q_1'$ where $s > 0$  is such that
the amount of expansion along $\cV_2(H^1)$ from $q_1$ to $q_3$ is equal to
the amount of expansion along $\cV_2(H^1)$ from $uq_1$ to $q_2$. Then, as in
\cite{Benoist:Quint}, 
\begin{equation}
\label{eq:intro:A:Aprime}
f_2(q_2) = A_* f_2(q_3), \text{ and } f_2(q_2') =
A'_* f_2(q_3'), 
\end{equation}
where $A$ and $A'$ are essentially the same bounded linear map.  
But $q_3$ and $q_3'$ approach each other, so that
\begin{displaymath}
f_2(q_3) \approx f_2(q_3'). 
\end{displaymath}
Hence 
\begin{equation}
\label{eq:f1:approx:same}
f_2(q_2) \approx f_2(q_2').
\end{equation}
Taking a limit as $\ell \to \infty$ 
of the points $q_2$ and $q_2'$ we obtain points $\tilde{q}_2$
and $\tilde{q}_2'$ in the same leaf of $(\cV_1+\cV_2)(H^1)$ and distance 
$\epsilon$ apart such that
\begin{equation}
\label{eq:intro:f1:q2:f1:q2prime}
f_2(\tilde{q}_2) = f_2(\tilde{q}_2'). 
\end{equation}
This means that the
conditional measure $f_2(\tilde{q}_2)$ is invariant under a
shift of size approximately $\epsilon$. 
Repeating this argument with $\epsilon \to 0$ we obtain a point
$p$ such that $f_2(p)$ is invariant under arbitrarily small
shifts. This implies that the conditional measure $f_2(p)$ restricts
to Lebesgue measure 
on some subspace $U_{new}$ of $(\cV_1+\cV_2)(H^1)$, which is distinct from the orbit of
$N$. Thus, we can enlarge $U^+$ to be $U^+ \oplus U_{new}$. 
 
\bold{Technical Problem \#1.}
The
argument requires that all eight points $q$, $q'$, $q_1$, $q_1'$,
$q_2$, $q_2'$, $q_3$, $q_3'$ belong to some ``nice'' set $K$ of almost
full measure. We will give a very rough outline of the solution to
this problem here; a more detailed outline is given at the beginning
of \S\ref{sec:conditional}. 

We have the following elementary statement:
\begin{lemma}
\label{lemma:non:trivial}
If $\nu_{W^-}$ is non-trivial, then for any $\delta > 0$ there exist
constants $c(\delta) > 0$ and $\rho(\delta) > 0$ such that
for any compact $K \subset X_0$ with $\nu(K) > 1-\delta$ there exists
a compact subset $K' \subset K$ with $\nu(K') > 1-c(\delta)$ so that 
for any $q \in K'$ there exists $q' \in K \cap
W^-[q]$ with 
\begin{displaymath}
\rho(\delta) < d(q,q') < 1/100. 
\end{displaymath}
Furthermore, $c(\delta) \to 0$ as $\delta \to 0$. 
\end{lemma}

In other words, there is a set $K' \subset K$ of almost full measure such that
every point $q \in K'$ has a ``friend'' $q' \in W^-[q]$, with $q'$
also in the ``nice'' set $K$, such that 
\begin{displaymath}
d(q,q') \approx 1/100. 
\end{displaymath}
Thus, $q$ can be chosen essentially anywhere in $X_0$. 
(In fact we use a variant of
Lemma~\ref{lemma:non:trivial}, namely
Proposition~\ref{prop:can:avoid:most:Mu}
in \S\ref{sec:conditional}.)

We also note the following trivial statement:
\begin{lemma}
\label{lemma:trivial:bilipshitz}
Suppose $\nu$ is a measure on $X_0$ invariant under the flow $g_t$. Let
$\hat{\tau}: 
X_0 \times \reals \to \reals$ be a function such that there exists
$\kappa > 1$ so that for all $x \in X_0$ and for $t > s$,
\begin{equation}
\label{eq:trivial:bilipshitz}
\kappa^{-1}(t-s) \le \hat{\tau}(x,t) - \hat{\tau}(x,s) \le \kappa(t - s).
\end{equation}
Let $\psi_t: X_0 \to X_0$ be given by $\psi_t(x) = g_{\hat{\tau}(x,t)}
x$. Then, for any $K^c \subset X_0$ and any $\delta >
  0$, there exists a subset $E \subset \reals$ of density at least
  $(1-\delta)$ such that for $t \in E$, 
\begin{displaymath}
\nu(\psi_t^{-1}(K^c)) \le (\kappa^2/\delta) \nu(K^c).
\end{displaymath}
\end{lemma}
(We remark that the maps $\psi_t$ are not a flow, since $\psi_{t+s}$
is not in general $\psi_t \circ \psi_s$. However,
Lemma~\ref{lemma:trivial:bilipshitz} still holds.)

In \S\ref{sec:new:bilipshitz} we show that roughly, 
$q_2 = \psi_t(q)$, where $\psi_t$ is as in
Lemma~\ref{lemma:trivial:bilipshitz}. (A more precise
statement, and the strategy for
dealing with this problem is given at the beginning of
\S\ref{sec:conditional}). Then, to make sure that $q_2$ avoids a ``bad
set'' $K^c$ of small measure, we make sure that $q \in
\psi_t^{-1}(K)$ which by Lemma~\ref{lemma:trivial:bilipshitz} has
almost full measure. Combining this with
Lemma~\ref{lemma:non:trivial}, we can see that we can choose $q$, $q'$
and $q_2$ all in an a priori prescribed subset $K$ of almost full
measure.  A similar argument can be done for all eight points, 
see \S\ref{sec:inductive:step}, where the precise arguments are assembled.

\bold{Technical Problem \#2.} Beyond the first step of 
the induction, the subspace
$U^+(x)$ may not be locally constant as $x$ varies along
$W^+(x)$. This complication has a ripple effect on the proof. In
particular, instead of dealing with the divergence of the points $g_t
u q_1$
and $g_t u' q_1'$ we need to deal with the divergence of the 
affine subspaces
$U^+[g_t u q_1]$ and $U^+[g_t u' q_1']$. As a first step, we project
$U^+[g_t u' q_1']$ to the leaf of $W^+$ containing $U^+[g_t u q_1]$,
to get a new affine subspace $\cU'$. One way to keep track of the
relative location of $U^+ = U^+[g_t u' q_1']$ and $\cU'$ is (besides keeping
track of the linear parts of $U^+$ and $\cU'$) 
to pick a
transversal $Z(x)$ to $U^+[x]$, and to keep track of the
intersection of $\cU'$ and $Z(x)$, see Figure~\ref{fig:divergence}.

\makefignocenter{\begin{minipage}{5.5in}
\begin{itemize}
\item[(a)] We keep track of the relative position of the subspaces
  $U^+[x]$ and $\cU'$ in part by picking a transversal $Z(x)$ to $U^+[x]$,
and noting the distance between $U^+[x]$ and $\cU'$ along $Z[x]$. 
\item[(b)] If we apply the
flow $g_t$ to the entire picture in (a), we see that the transversal
$g_t Z[x]$ can get almost parallel to $g_t U^+[x]$. Then, the distance
between $g_t U^+[x]$ and $g_t \cU'$ along $g_t Z[x]$ may be much
larger then the distance between $g_t x \in g_t U^+[x]$ and the closest point
in $g_t \cU'$.
\end{itemize}
\end{minipage}
}{fig:divergence}{\includegraphics{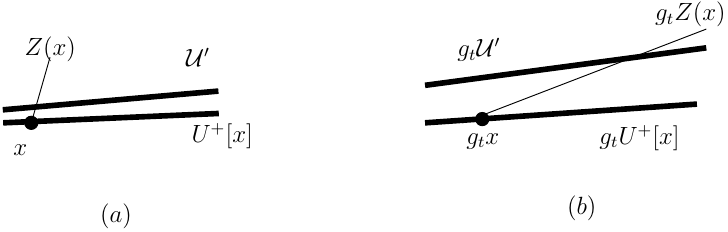}}

However, since we do not know at this point that the cocycle is
semisimple, we cannot pick $Z$ in a way which is invariant under the
flow. Thus, we have no choice except to pick some transversal $Z(x)$
to $U^+(x)$ at $\nu$-almost every point $x \in X_0$, and then deal with
the need to change transversal. 

It turns out that the formula for computing how $\cU' \cap Z$ changes
when $Z$ changes is non-linear (it involves inverting a certain
matrix). However, we would really like to work with linear maps. This
is done in two steps: first we show that we can choose the
approximation $\cU'$ and the transversals $Z(x)$ 
in such a way that changing transversals involves inverting a
unipotent matrix. This makes the formula for changing transversals
polynomial. In the second step, we embed the space of 
parameters of affine subspaces near $U^+[x]$ 
into a certain tensor power space $\bH(x)$ so that on the level
of $\bH(x)$ the change of transversal map becomes linear. The
details of this construction are in \S\ref{sec:divergence:subspaces}.

\bold{Technical Problem \#3.}
There may be more than one $U^+$-inert Lyapunov exponent. In that
case, we do not have precise control over how $q_2$ and $q_2'$ diverge. In
particular the assumption that $q_2 - q_2'$ is nearly in the direction
of $\cV_2(H^1)(q_2)$ is not justified. Also we really need to work with
$U^+[q_2]$ and $U^+[q_2']$. So let $\bfv \in \bH(q_2)$ denote the vector
corresponding to (the projection to $W^+(q_2)$ of) the affine subspace
$U^+[q_2']$. (This vector $\bfv$ takes on the role of $q_2 -
q_2'$). We have no a-priori control over the direction of $\bfv$ (even
though we know that $\|\bfv \| \approx \epsilon$, and we know that
$\bfv$ is almost contained in $\bE(q_2) \subset \bH(q_2)$, where
$\bE(x)$ is defined in \S\ref{sec:first:divergence} as the union of
the Lyapunov subspaces of $\bH(x)$ corresponding to the $U^+$-inert
Lyapunov exponents.) 

The idea is to vary $u$ (while keeping $q_1$, $q_1'$, $\ell$
fixed). To make this work, we need to define a {\em finite}
collection of subspaces $\bE_{[ij],bdd}(x)$ of $\bH(x)$ (which
actually only make sense on a certain finite measurable cover $X$ of
$X_0$) such that
\begin{itemize}
\item[{\rm (a)}] By varying $u$ (while keeping $q_1$, $q_1'$, $\ell$ fixed) we
  can make sure that the vector $\bfv$ becomes close to one of the
  subspaces $\bE_{[ij],bdd}$,
and
\item[{\rm (b)}] For a suitable choice of point $q_3 = q_{3,ij} = g_{s_{ij}}
  q_1$, the map 
\begin{displaymath}
(g_t u g_{-s_{ij}})_* \bE_{[ij],bdd}(q_3) \to
  \bE_{[ij],bdd}(q_2)
\end{displaymath}
is a linear map whose norm is bounded independently of the parameters.  

\item[{\rm (c)}] Also, for a suitable choice of point $q_3' =
  q'_{3,ij} = g_{s'_{ij}}   q_1$, the map 
\begin{displaymath}
(g_t u g_{-s'_{ij}})_* \bE_{[ij],bdd}(q_3') \to
  \bE_{[ij],bdd}(q_2') 
\end{displaymath}
is a linear map whose norm is bounded independently of the parameters.  
\end{itemize}
For the precise conditions see
Proposition~\ref{prop:some:fraction:bounded} and
Proposition~\ref{prop:ej:bdd:transport:bounded}. 
This construction is done in detail in
\S\ref{sec:bounded:synchornized}. 
The general idea is as follows: Suppose  $\bfv \in \bE_i(x) \oplus \bE_j(x)$
where $\bE_i(x)$ and $\bE_j(x)$ are the Lyapunov subspaces
corresponding to the $U^+$-inert (simple) Lyapunov exponents $\lambda_i$ and
$\lambda_j$. Then, if while varying $u$, the vector $\bfv$ does not
swing towards either $\bE_i$ or $\bE_j$, we say that $\lambda_i$ and
$\lambda_j$ are ``synchronized''. In that case, we consider the subspace
$\bE_{[i]}(x) = \bE_i(x) \oplus \bE_j(x)$ and show that (b) and (c) hold.

The conditions (b) and (c) allow us to define in
\S\ref{sec:equivalence} conditional measures
$f_{ij}$ on $W^+(x)$ which are associated to each subspace 
$\bE_{[ij],bdd}$. In fact the measures are supported on the points $y
\in W^+[x]$ such that the affine subspace $U^+[y]$ maps to a vector in
$\bE_{[ij],bdd}(x) \subset \bH(x)$. 

\bold{Technical Problem \#4.} More careful analysis (see the
discussion following the statement of
Proposition~\ref{prop:nearby:linear:maps})
shows that the
maps $A$ and $A'$ of (\ref{eq:intro:A:Aprime}) are not exactly the
same. Then, when 
one passes to the limit $\ell \to \infty$ one gets, instead of
(\ref{eq:intro:f1:q2:f1:q2prime}), 
\begin{displaymath}
f_{ij}(\tilde{q}_2) = P^+(\tilde{q}_2,\tilde{q}_2')_* f_{ij}(\tilde{q}_2')
\end{displaymath}
where $P^+: W^+(\tilde{q}_2) \to W^+(\tilde{q}_2')$ is a certain
unipotent map (defined in \S\ref{sec:subsec:connection}). Thus the
conditional measure $f_{ij}(\tilde{q}_2)$ is invariant under the
composition of a translation  of size $\epsilon$ and a unipotent
map. Repeating the argument with  $\epsilon \to 0$ we obtain a point
$p$ such that the conditional measure at $p$ is invariant under
arbitrarily small combinations of (translation + unipotent map). This
does {\em not} imply that the conditional measure $f_{ij}(p)$ restricts to
Lebesgue measure on some subspace of $W^+$, but it does imply that it
is in the Lebesgue measure class along some polynomial curve in
$W^+$. More precisely, for $\nu$-a.e $x \in X$ 
there is a subgroup $U_{new} = U_{new}(x)$ of the affine
group of $W^+(x)$ such that the conditional measure of $f_{ij}(x)$ on
the polynomial curve $U_{new}[x] \subset W^+[x]$ is induced from the
Haar measure on $U_{new}$. (We call such a set a ``generalized
subspace''). The exact definition is given in
\S\ref{sec:divergence:subspaces}. 

Thus, during the induction steps, we need to deal with generalized
subspaces. This is not a very serious complication since the general
machinery developed in \S\ref{sec:divergence:subspaces} can deal
with generalized subspaces as well as with ordinary affine
subspaces. 
\medskip

\bold{Completion of the proof of Theorem~\ref{theorem:main:partone}.}
Let \index{$L$@$\cL(x)$} $\cL(x) \subset H^1(M, \Sigma, \reals)$
be the smallest subspace such that $\nu_{W^-(x)}$ is
supported on \index{$L$@$\cL^-(x)$} $\cL^-(x)$. 
Roughly, the above argument can be iterated until we
know the
conditional measure $\nu_{W^+(x)}$ is Lebesgue on a subspace
$\cU^+[x]$,  where $\cU(x) \subset H^1(M, \Sigma, \reals)$
contains $\cL(x)$. (The precise condition for when the induction stops
is given by Lemma~\ref{lemma:Lplus:Splus:equiv} and
Proposition~\ref{prop:stop:induction:condition}). 
Then a Margulis-Tomanov style entropy
comparison argument (see \S\ref{sec:proof:theorem:partone})
shows that $\cU(x) = \cL(x)$, and the conditional
measures along $\cL^-(x)$ are Lebesgue. Since $\cU^+(x)$ contains the
orbit of the unipotent direction $N$, this implies that $\cL^-(x)$ contains the
orbit of the opposite unipotent direction $\bar{N} \subset
SL(2,\reals)$. Thus, the conditional measure along the orbit of
$\bar{N}$ is Lebesgue, which means that $\nu$ is
$\bar{N}$-invariant. This, together with the assumption that $\nu$ is
$P = AN$-invariant implies that $\nu$ is $SL(2,\reals)$-invariant,
completing the proof of Theorem~\ref{theorem:main:partone}.

\section{Hyperbolic properties of the geodesic flow}
\label{sec:semi:markov}

\bold{The spaces $X_0$ and $\tilde{X}_0$.} 
Let \index{$X_0$}$X_0$ 
be a finite cover of the stratum $\cH_1(\alpha)$ which is a
manifold.  (Such a cover may be obtained by choosing a level $3$
structure, i.e.\ a basis for the mod $3$ homology of the surface). Let \index{$X$@$\tilde{X}_0$}$\tilde{X}_0$ 
be the universal cover of $X_0$. Then the fundamental group 
\index{$\pi_1(X_0)$}$\pi_1(X_0)$
acts properly discontinuously on $\tilde{X}_0$. Let $\nu$ be a
$P$-invariant ergodic probability measure on $X_0$.




We recall the following standard fact:
\begin{lemma}[Mautner Phenomenon]
\label{lemma:mautner}
Let $\nu$ be an ergodic $P$-invariant measure on a space $Z$. Then
$\nu$ is $A$-ergodic. 
\end{lemma}

\bold{Proof.} See e.g.\ \cite{Mozes:epimorphic}.  
\qed\medskip

\begin{lemma}
\label{lemma:embedded:leaves}
For almost all $x \in X_0$, the affine exponential map from
$W^+(x)$ to $W^+[x]$ is globally defined and is bijective, endowing
$W^+[x]$ with a global affine structure. The same holds for $W^-[x]$. 
\end{lemma}
\bold{Proof.} Since $W^-$ and $W^+$ play the role of the stable and unstable
foliations for the action of $g_t \in A$ (cf.\ Lemma~\ref{lemma:forni}),
this follows from the Poincar\'e recurrence theorem. 
\qed\medskip

\bold{The bundle $H^1$.} 
Let \index{$H^1$}$H^1$ denote the bundle whose fiber above $x \in X_0$ is
$H^1(M,\Sigma,\reals)$. We denote the fiber above the point $x \in X_0$
by \index{$H^1(x)$}$H^1(x)$. 

The geodesic flow acts on $H^1$ by parallel transport using the
Gauss-Manin connection (see \S\ref{sec:outline:step1}).

\bold{The bundles $H^1_+$ and $H^1_-$.} 
Let \index{$H^1_+$}$H^1_+$ denote the same bundle as $H^1$ except that the action of
$g_t$ on $H^1_+$ includes an extra multiplication by $e^t$ on the
fiber. (In other words, if $h_t(x,v) = (x, e^{t} v)$ and 
$i: H^1 \to H^1_+$ is the identity map, 
then $g_t \circ i(x,v) = h_t \circ i \circ g_t(x,v)$). Similarly, let
\index{$H^1_-$}$H^1_-$ denote the same bundle as $H^1$ except that the
action of $g_t$ includes an extra multiplication by $e^{-t}$ on the fiber. 

We use the notation \index{$H^1_+(x)$}$H^1_+(x)$ and 
\index{$H^1_-(x)$}$H^1_-(x)$ to refer to the fiber of
the corresponding bundle above the point $x \in X_0$.

\bold{The bundles $H_{big}$, $H_{big}^{(+)}$, $H_{big}^{(-)}$, $H_{big}^{(++)}$ 
and $H_{big}^{(--)}$.}  
In this paper, 
we will need to deal with several bundles derived from the Hodge
bundle $H^1$. It is convenient to introduce a bundle
\index{$H_{big}$}$H_{big}$ 
so that every
bundle we will need will be a subbundle of $H_{big}$.  Let $d \in
\natls$ be a large integer chosen later (it will be chosen in
\S\ref{sec:divergence:subspaces} and will depend only on the Lyapunov spectrum
of the Kontsevich-Zorich cocycle). Let 

\begin{displaymath}
\hat{H}_{big}(x) = \bigoplus_{k=1}^d \bigoplus_{j=1}^k 
\left(\bigotimes_{i=1}^j H^1(x) \otimes \bigotimes_{l=1}^{k-j} (H^1(x))^*
\right), 
\end{displaymath}
\begin{displaymath}
\hat{H}_{big}^{(+)}(x) = \bigoplus_{k=1}^d \bigoplus_{j=1}^k 
\left(\bigotimes_{i=1}^j H^1_+(x) \otimes \bigotimes_{l=1}^{k-j} (H^1_+(x))^*
\right), 
\end{displaymath}
\begin{displaymath}
\hat{H}_{big}^{(-)}(x) = \bigoplus_{k=1}^d \bigoplus_{j=1}^k 
\left(\bigotimes_{i=1}^j H^1_-(x) \otimes \bigotimes_{l=1}^{k-j} (H^1_-(x))^*
\right), 
\end{displaymath}
and let
\begin{displaymath}
\tilde{H}_{big}(x) = \hat{H}_{big}(x) \oplus \hat{H}_{big}^{(+)}(x) \oplus
\hat{H}_{big}^{(-)}(x). 
\end{displaymath}
Suppose $L_1 \subset L_2 \subset \tilde{H}_{big}$ are
 $g_t$-invariant subbundles. We say that $L_2/L_1$ is an admissible
 quotient if the cocycle on $L_2/L_1$ is measurably conjugate to a
 conformal cocyle (see Lemma~\ref{lemma:jordan:canonical:form}), and
 also $L_2/L_1$ is maximal in the sense that if $L_2' \supset L_2$ and
 $L_1' \subset L_1$ are $g_t$-invariant subbundles with the cocycle
 $L_2'/L_1'$ measurably conjugate to a conformal cocycle, then
 $L_2'=L_2$ and $L_1'=L_1$.
We then let $\Delta_{big}$ denote the set of all
admissible quotients of $\tilde{H}_{big}$ and let 
\begin{displaymath}
\index{$H_{big}(x)$}H_{big}(x) = \bigoplus_{Q \in \Delta_{big}} Q(x). 
\end{displaymath}
(We apply a similar operation to the bundles $\hat{H}_{big}^{(+)}$ and 
$\hat{H}_{big}^{(-)}$ to get bundles
\index{$H_{big}^{(+)}(x)$}$H_{big}^{(+)}$ and \index{$H_{big}^{(-)}(x)$} 
$H_{big}^{(-)}$. )

The flow $g_t$ acts on the bundle $H_{big}$ in the natural way. We
denote the action on the fibers by $(g_t)_*$. \index{$g_t$@$(g_t)_*$} 
Let $H_{big}^{(++)}(x)$ \index{$H_{big}^{(++)}(x)$} 
denote the direct sum of the positive Lyapunov subspaces of
$H_{big}(x)$. Similarly, let \index{$H_{big}^{(--)}(x)$}
$H_{big}^{(--)}(x)$ denote the direct sum of the negative Lyapunov
subspaces of $H_{big}(x)$. 

\begin{lemma}
The subspaces $H_{big}^{(++)}(x)$ are locally constant along
  $W^+[x]$, i.e.\ for almost all $x \in \tilde{X}_0$ and almost all 
$y \in W^+[x]$ close to $x$ we have
  $H_{big}^{(++)}(y) = H_{big}^{(++)}(x)$. Similarly, the subspaces
$H_{big}^{(--)}(x)$ are locally constant along $W^-[x]$.
\end{lemma}

\bold{Proof.} Note that
\begin{displaymath}
H_{big}^{(++)}(x) = \left\{ v \in H_{big}(x) \st \lim_{t \to \infty}
  \frac{1}{t} 
\log \frac{ \|(g_{-t})_* v\|}{\|v\|} < 0 \right\}
\end{displaymath}
Therefore, the subspace $H_{big}^{(++)}(x)$ depends only on the trajectory
$g_{-t} x$ as $t \to 
\infty$. However, if $y \in W^+[x]$ then $g_{-t} y$ will for large $t$
be close to $g_{-t} x$, and so in view of the affine
structure, $(g_{-t})_*$ will be the same linear map on $H_{big}(x)$ and
$H_{big}(y)$. This implies that $H_{big}^{(++)}(x)=
H_{big}^{(++)}(y)$.  
\qed\medskip

\bold{The Avila-Gou\"ezel-Yoccoz norm.} The Avila-Gou\"ezel-Yoccoz
 norm on the relative cohomology group
$H^1(M,\Sigma,\reals)$ is described in Appendix~\ref{sec:app:forni}. This then
induces a norm which we will denote by $\|\cdot\|_Y$ 
\index{abs@$\norm{\cdot}_Y$}
and then, as the projective cross norm, also on $H_{big}$.
We also use the notation  \index{abs@$\norm{\cdot}_{Y,x}$}$\norm{\cdot}_{Y,x}$
to denote the AGY norm at $x \in X_0$. 

\bold{The distance $d^+(x,y)$.} 
Since the tangent
space to $W^+[x]$ is included in $H^1(M,\Sigma,\reals)$, the AGY norm on
$H^1(M,\Sigma,\reals)$ defines a distance on $W^+[x]$. We denote this
distance by \index{$d^+(\cdot,\cdot)$}$d^+(\cdot, \cdot)$. (Thus, for $y,z \in W^+[x]$,
$d^+(y,z)$ is the length of the shortest path in $W^+[x]$ connecting
$y$ and $z$, where lengths of paths are measured using the AGY norm).  

\bold{The ball $B^+(x,r)$.}
Let $\index{$B^+(x,r)$}B^+(x,r) \subset W^+[x]$ denote the ball of
radius $r$ centered at $x$, in the metric $d^+(\cdot, \cdot)$. 
\medskip

The following is a rephrasing of \cite[Proposition~5.3]{Avila:Gouezel}: 
\begin{proposition}
\label{prop:AGY:regularity}
For all $x \in X_0$, $x+v$ is well defined for $v\in W^+(x)$ with
$\|v\|_Y \le 1/2$. Also, for all $y, z \in B^+(x,1/50)$, we have 
\begin{displaymath}
\tfrac{1}{2} \| y - z \|_{Y,y} \le \| y - z \|_{Y,z} \le 2 \|y - z\|_{Y,y},
\end{displaymath}
and
\begin{displaymath}
\tfrac{1}{2} \| y - z \|_{Y,y} \le d^+(y,z) \le 2 \| y - z \|_{Y,y}.  
\end{displaymath}
\end{proposition}

Note that we have a similar distance $d^-(\cdot,\cdot)$ on
$W^-[x]$, and the analogue of Proposition~\ref{prop:AGY:regularity}
holds.

\bold{The ``distance'' \index{$d^{X_0}(\cdot, \cdot)$}$d^{X_0}(\cdot, \cdot)$.}
Suppose $x,y \in \tilde{X}_0$ are not far apart. Then, there exist unique 
$z  \in W^+[x]$ and $t \in \reals$ such that $g_t z \in W^-[y]$. We
then define
\begin{displaymath}
d^{X_0}(x,y) = d^+(x,z) + |t| + d^-(g_t z, y).
\end{displaymath}
Thus, if $y \in W^+[x]$ then $d^{X_0}(x,y) = d^+(x,y)$, and if $y \in
W^-[x]$, then $d^{X_0}(x,y) = d^-(x,y)$. 

We sometimes abuse notation by using the notation $d^{X_0}(x,y)$ where $x$, $y
\in X_0$. By this we mean $d^{X_0}(\tilde{x},\tilde{y})$ where
$\tilde{x}$ and $\tilde{y}$ are appropriate lifts of $x$ and $y$.
\medskip

Choose a compact subset \index{$K_{thick}'$}$K_{thick}' \subset X_0$ with
$\nu(K_{thick}') \ge 5/6$. Let \index{$K_{thick}$}$K_{thick} = \{ x \in X_0 \st
d^{X_0}(x,K_{thick}') \le 1/100 \}$. 

\begin{lemma}
\label{lemma:forni}
There exists $\alpha > 0$  
such that the following holds: 
\begin{itemize}
\item[{\rm (a)}] Suppose $x \in X_0$ and $t > 0$ are such
that the geodesic segment from $x$ to $g_t x$ spends at least half the
time in $K_{thick}$. Then,
for all $v \in W^-(x)$, 
\begin{displaymath}
\|(g_t)_* v\|_Y \le e^{-\alpha t} \|v\|_Y.
\end{displaymath}
\item[{\rm (b)}] Suppose $x \in X_0$ and $t > 0$ are as in (a). Then,
for all $v \in W^+(x)$, 
\begin{displaymath}
\|(g_t)_* v\|_Y \ge e^{\alpha t} \|v\|_Y.
\end{displaymath}
\item[{\rm (c)}] 
For every $\epsilon >
  0$ there exist a compact subset
\index{$K_{thick}''$}$K_{thick}'' \subset X_0$ with $\nu(K_{thick}'') > 1-\epsilon$ 
and $t_0 > 0$ such that for $x \in K_{thick}''$, $t >
t_0$ and all $v \in H_{big}^{(++)}(x)$, 
\begin{displaymath}
\|(g_t)_* v\|_Y  \ge e^{\alpha t} \|v\|_Y. 
\end{displaymath}
\item[{\rm (d)}] For all $v \in W^+(x)$, all $x \in X_0$ and all $t >
  0$,
\begin{displaymath}
\|(g_t)_* v \|_Y \ge \|v\|_Y. 
\end{displaymath}
\end{itemize}
\end{lemma}

\bold{Proof.} Parts (a), (b) and (d) follow from
Theorem~\ref{theorem:decay:relative:norm}. Part (c) follows
immediately from the Osceledets multiplicative ergodic theorem.
\qed\medskip

We also have the following simpler statement:
\begin{lemma}
\label{lemma:forni:upper}
There exists $N > 0$ such that for all $x \in X_0$, all
  $t \in \reals$, and all $v \in 
  H_{big}(x)$, 
\begin{displaymath}
e^{-N |t|} \|v\|_Y \le \|(g_t)_* v\|_Y \le e^{N |t|} \|v\|_Y. 
\end{displaymath}
For $v \in W^+[x]$, we can take $N=2$. 
\end{lemma}

\bold{Proof.} This follows immediately from
Theorem~\ref{theorem:decay:relative:norm}. 
\qed\medskip

\begin{proposition}
\label{prop:semi:markov}
Suppose $\cC \subset X_0$ is a set with $\nu(\cC) > 0$, and $T_0: \cC \to
\reals^+$ is a measurable function which is finite a.e. 
Then we can find $x_0 \in \tilde{X}_0$,
a subset $\cC_1 \subset W^-[x_0] \cap \pi^{-1}(\cC)$ and for each $c
\in \cC_1$ a subset $\index{$E^+[c]$}E^+[c] \subset W^+[c]$ of
  diameter in the AGY metric at most $1/200$ and a number
$\index{$t(c)$} t(c) > 0$ such
that if we let
\begin{displaymath}
\index{$J_c$} J_c = \bigcup_{0 \le t < t(c)} g_{-t} E^+[c],
\end{displaymath}
then the following holds:
\begin{itemize}
\item[{\rm (a)}] $E^+[c]$ is relatively open in $W^+[c]$. 
\item[{\rm (b)}] $\pi(J_{c}) \cap \pi(J_{c'}) = \emptyset$ if $c \ne c'$. 
\item[{\rm (c)}] $\pi(J_{c})$ is embedded in $X_0$, i.e.\ if $\pi(g_{-t}
  x) = \pi(g_{-t'} x')$ where $x, x'\in E^+[c]$ and $0 \le t < t(c)$, $0 \le t' < t(c)$
  then $x = x'$ and $t = t'$. 
\item[{\rm (d)}] $\bigcup_{c \in \cC_1} \pi(J_c)$ is conull in $X_0$.
\item[{\rm (e)}] For every $c \in \cC_1$ there exists $c' \in \cC_1$ such
  that $\pi(g_{-t(c)} E^+[c]) \subset \pi(E^+[c'])$. 
\item[{\rm (f)}] $t(c) > T_0(c)$ for all $c \in \cC_1$. 
\end{itemize}
\end{proposition}
\bold{Remark.} All the construction in \S\ref{sec:semi:markov} will
depend on the choice of $\cC$ and $T_0$, but we will suppress this from
the notation. The set $\cC$ and the function $T_0$ will be finally
chosen in Lemma~\ref{lemma:exists:C:T0}.
\medskip

The proof of Proposition~\ref{prop:semi:markov} relies on the following:
\begin{lemma}
\label{lemma:semi:markov}
Suppose $\cC \subset X_0$ is a set with $\nu(\cC) > 0$, and $T_0: \cC \to
\reals^+$ is a measurable function which is finite a.e. 
Then we can find $x_0 \in \tilde{X}_0$, a subset $\cC_1 \subset W^-[x_0] \cap \pi^{-1}(\cC)$ and for each $c
\in \cC_1$ a subset $E^+[c] \subset W^+[c]$ of
  diameter in the AGY metric at most $1/200$ so that the following hold:
\begin{itemize}
\item[{\rm (0)}] $E^+[c]$ is a relatively open subset of $W^+[c]$. 
\item[{\rm (1)}] The set $E = \pi\left(\bigcup_{c \in \cC_1}
    E^+[c]\right)$ is embedded in $X_0$, i.e.\ if $\pi(x) = 
  \pi(x')$ where $x \in E^+[c]$ and $x' \in E^+[c']$, then $x = x'$
  and $c = c'$. 
\item[{\rm (2)}] For some $\epsilon > 0$, 
        $\nu(\bigcup_{t\in (0,\epsilon)} g_t E) > 0$. 
\item[{\rm (3)}] If $t > 0$ and $c \in \cC_1$ is such that $\pi(g_{-t}E^+[c])
  \cap E \ne \emptyset$, then $\pi(g_{-t}E^+[c]) \subset \pi(E^+[c'])$ for some
  $c' \in \cC_1$.  
\item[{\rm (4)}] Suppose $t$, $c$, $c'$ are as in {\rm (3)}. Then $t >
  T_0(c)$. 
\end{itemize}
\end{lemma}
\bold{Proof.} 
This proof is essentially identical to the proof of
  Lemma~\ref{lemma:mt:91}, except that we need to take care that (4)
  is satisfied. In this proof, for 
$x \in \cC$, we denote by $\nu_{W^\pm[x]}$
the conditional measure of $\nu$ along $W^\pm[x] \cap \cC$.

\mccc{make sure atoms of partition have radius < 1/100}

Choose $T_1 > 0$ so that if we let $\cC_4
= \{ x \in \cC \st T_0(x) < T_1 \}$ then $\nu(\cC_4) > \nu(\cC)/2$.  Let
$X_{per}$ denote the union of the periodic orbits of $g_t$. By the
$P$-invariance of $\nu$ and the 
ergodicity of $g_t$, $\nu(X_{per}) = 0$, and the same is true of the set
$X_{per}' = \bigcup_{x \in X_{per}} W^-[x]$. Therefore there
exists $x_0 \in \pi^{-1}(\cC_4)$ 
and a compact subset $\cC_3 \subset W^-[x_0] \cap \pi^{-1}(\cC_4)$ 
with $\nu_{W^-[x_0]}(\cC_3) > 0$ 
such that for $x \in \cC_3$ and $0 < t < T_1$, $\pi(g_{-t} x) \notin
\pi(\cC_3)$. Then, since $\cC_3$ is compact, we can find a small neighborhood
$V^+ \subset W^+$ of the origin such that the set 
$\pi\left(\bigcup_{c \in \cC_3} V^+[c]\right)$ is embedded in $X_0$ and for $x \in
\bigcup_{c \in \cC_3} V^+[c]$ and $0 < t < T_1$, $\pi(g_{-t} x) \not\in
\pi\left(\bigcup_{c \in \cC_3} V^+[c]\right)$. 

There exists $\cC_2 \subset \cC_3$ with
$\nu_{W^-[x_0]}(\cC_2) > 0$ and $N >
T_1$ such that for all $c \in \cC_2$ and all $T > N$, 
\begin{displaymath}
|\{ t \in [0,T] \st \pi(g_{-t} c) \in K_{thick}' \}| \ge T/2. 
\end{displaymath}
Then, for $c \in \cC_2$, $T > N$ and $x \in V^+[c]$, 
\begin{displaymath}
|\{ t \in [0,T] \st \pi(g_{-t} x) \in K_{thick} \}| \ge T/2. 
\end{displaymath}
Let 
\begin{displaymath}
M = \sup\left\{ \frac{\|v\|_{Y,x}}{\|v\|_{Y,y}} \st x \in V^+[c], y \in
  V^+[c], c \in \cC_2, v \in W^+(x) \right\}
\end{displaymath}
Let $\alpha > 0$ be as in Lemma~\ref{lemma:forni},
and choose $N_1 > N$ such that $M^2 e^{-\alpha N_1} < 1/10$.
Then, for $c \in \cC_2$, 
$x,y \in \pi(V^+[c])$ and $t>N_1$ such that $g_{-t} x \in
\pi\left(\bigcup_{c \in   \cC_2} V^+[c]\right)$,  
in view of Lemma~\ref{lemma:forni} and Proposition~\ref{prop:AGY:regularity},
\begin{displaymath}
d^{X_0}(g_{-t} x, g_{-t} y) \le \frac{1}{10} d^{X_0}(x,y). 
\end{displaymath} 
Now choose $\cC_1 \subset \cC_2$
with $\nu_{W^-[x_0]}(\cC_1) > 0$ so that if we let $Y = \pi\left(\bigcup_{c \in
  \cC_1} V^+[x]\right)$ then $g_{-t} Y \cap Y = \emptyset$ for $0 < t <
\max(T_1,N_1)$, in 
other words, the first return time to $Y$ is at least
$\max(T_1,N_1)$. (This can be done e.g.\ by Rokhlin's Lemma).
Condition (4) now follows since $T_0(c) < T_1$ for all $c \in
\cC_1$. 
The rest of the proof is essentially the same as the proof of
Lemma~\ref{lemma:mt:91}, applied to the first return map of $g_{-t}$
to $Y$.   \mcc{explain better}
\qed\medskip

\bold{Proof of Proposition~\ref{prop:semi:markov}.}
For $x \in E$,  let $t(x) \in \reals^+$ be the smallest such that
$g_{-t(x)} x \in E$. By property (3), the function $t(x)$ is
constant on each set of the form $\pi(E^+[c])$. Let $F_t = \{ x \in E \st
t(x) = t \}$. (We have $F_t = \emptyset$ if $t < N_1$). By property
(2) and the ergodicity of $g_{-t}$, up to a null set, 
\begin{displaymath}
X_0 = \bigsqcup_{t > 0} \bigsqcup_{s < t} g_{-s} F_t. 
\end{displaymath}
Then properties (a)-(f) are easily verified. 
\qed\medskip

\bold{Notation.} For $x \in X_0$, let 
\index{$J[x]$}$J[x]$ denote the set $\pi(J_c)$ containing $x$. 
For $x \in \tilde{X}_0$, let \index{$J[x]$}$J[x]$ 
denote $\gamma J_c$ where $\gamma
\in \pi_1(X_0)$ is such that $\gamma^{-1} x \in J_c$.

\begin{lemma}
\label{lemma:stay:in:same:orbit:nontimechanged}
Suppose $x \in \tilde{X}_0$, $y \in W^+[x] \cap J[x]$. Then for any $t > 0$, 
\begin{displaymath}
g_{-t} y \in J[g_{-t} x] \cap W^+[g_{-t} x]. 
\end{displaymath}
\end{lemma}

\bold{Proof.} This follows immediately from property (e) of
Proposition~\ref{prop:semi:markov}. 
\qed\medskip

\bold{Notation.} For $x \in X_0$, let
\begin{displaymath}
\index{$B$@$\gB_t[x]$}\gB_t[x] = \pi(g_{-t}(J[g_t \tilde{x}] \cap W^+[g_t
\tilde{x}])), \qquad \text{ where $\tilde{x}$ is any element of
  $\pi^{-1}(x)$.}
\end{displaymath}

\begin{lemma} $ $
\label{lemma:gB:properties}
\begin{itemize}
\item[{\rm (a)}] For $t' > t \ge 0$, $\gB_{t'}[x] \subset \gB_t[x]$. 
\item[{\rm (b)}] Suppose $t \ge 0, t' \ge 0$, $x \in X_0$ and $x' \in X_0$ are such that 
$\gB_t[x] \cap \gB_{t'}[x'] \ne \emptyset$. Then either $\gB_t[x]
\supseteq \gB_{t'}[x']$ or $\gB_{t'}[x'] \supseteq \gB_t[x]$ (or
both). 
\end{itemize}
\end{lemma}

\bold{Proof.} Part (a) is a restatement of
Lemma~\ref{lemma:stay:in:same:orbit:nontimechanged}. 
For (b), without loss of generality, we may assume that
$t' \ge t$. Then, by (a), we have $\gB_t[x] \cap \gB_t[x'] \ne
\emptyset$. 

Suppose $y \in \gB_t[x] \cap \gB_t[x']$. Then $g_t y \in \gB_0[g_t x]$ 
and $g_{t} y \in \gB_0[g_t x']$. 
Since the sets $\gB_0[z]$, $z \in X_0$ form a partition, we must have
$\gB_0[g_t x] = \gB_0[g_t x']$. Therefore, $\gB_t[x] = \gB_t[x']$, and
thus, by (a), 
\begin{displaymath}
\gB_{t'}[x'] \subset \gB_t[x'] = \gB_t[x].
\end{displaymath}
\qed\medskip

By construction, the sets $\gB_0[x]$ are the atoms of a measurable partition of
$X_0$ subordinate to $W^+$ (see
Definition~\ref{def:subordinate}). \mcc{check this}
Then, let \index{$\nu_{W^+[x]}$} $\nu_{W^+[x]}$ denote the conditional measure of $\nu$
along the atom of the partition containing $x$. For
  notational simplicity, for $E \subset W^+[x]$, we sometimes write
  \index{$\nu_{W^+}$} $\nu_{W^+}(E)$ instead of $\nu_{W^+[x]}(E)$.

\begin{lemma}
\label{lemma:gB:vitali}
Suppose $\delta > 0$ and $K \subset X_0$ is such that
$\nu(K) > 1-\delta$. Then there exists a subset $K^* \subset
K$ with $\nu(K^*) > 1 - \delta^{1/2}$ such that for any $x \in K^*$,
and any $t > 0$, 
\begin{displaymath}
\nu_{W^+}(K \cap \gB_t[x]) \ge (1-\delta^{1/2}) \nu_{W^+}(\gB_t[x]). 
\end{displaymath}
\end{lemma}

\bold{Proof.} Let $E = K^c$, so $\nu(E) \le \delta$. Let $E^*$
denote the set of $x \in X_0$ such that there exists some $\tau \ge 0$ with
\begin{equation}
\label{eq:gB:lots:of:bad:set}
\nu_{W^+}(E \cap \gB_\tau[x]) \ge \delta^{1/2} \nu_{W^+}(\gB_\tau[x]). 
\end{equation}
It is enough to show that $\nu(E^*) \le \delta^{1/2}$. 
Let $\tau(x)$ be the smallest $\tau>0$ so that (\ref{eq:gB:lots:of:bad:set})
holds for $x$. Then the (distinct) sets $\{\gB_{\tau(x)}[x]\}_{x \in E^*}$
cover $E^*$ and are pairwise disjoint by
Lemma~\ref{lemma:gB:properties} (b). Let 
\begin{displaymath}
F = \bigcup_{x \in E^*} \gB_{\tau(x)}[x].
\end{displaymath}
Then $E^* \subset F$. For every set of the form $\gB_0[y]$, let
$\Delta(y)$ denote the set of distinct sets $\gB_{\tau(x)}[x]$ where
$x$ varies over $\gB_0[y]$.  Then, by (\ref{eq:gB:lots:of:bad:set})
\begin{multline*}
\nu_{W^+}(F \cap \gB_0[y]) = \sum_{\Delta(y)} \nu_{W^+}(\gB_{\tau(x)}) \le \\
\le \delta^{-1/2} \sum_{\Delta(y)} \nu_{W^+}(E \cap \gB_{\tau(x)}[x]) \le
\delta^{-1/2} \nu_{W^+}(E \cap \gB_0[y]). 
\end{multline*}
Integrating over $y$, we get $\nu(F) \le \delta^{-1/2}
\nu(E)$. Hence, 
\begin{displaymath}
\nu(E^*) \le \nu(F) \le \delta^{-1/2} \nu(E) \le \delta^{1/2}.  
\end{displaymath}
\qed\medskip

\section{General cocycle lemmas}
\label{sec:general:cocycle}

\subsection{Lyapunov subspaces and flags}
\label{sec:subsec:lyapunov:flags}
Let \index{$V$@$\cV_i(H^1)(x)$} $\cV_i(H^1)(x)$,  $1 \le i \le k$
denote the Lyapunov subspaces of
the Kontsevich-Zorich cocycle under the action of the geodesic flow
$g_t$, and let \index{$\lambda_i(H^1)$}$\lambda_i(H^1)$, $1 \le i \le k$
denote the (distinct) Lyapunov
exponents. Then we have for almost all $x \in X_0$, 
\begin{displaymath}
H^1(M, \Sigma, \reals) =
\bigoplus_{i=1}^{k}  \cV_i(H^1)(x)  
\end{displaymath}
and for all non-zero $v \in \cV_i(H^1)(x)$, 
\begin{displaymath}
\lim_{t \to \pm\infty} \frac{1}{t} \log \frac{ \|(g_t)_* v
  \|}{\|v\|} = \lambda_i(H^1), 
\end{displaymath}
where $\| \cdot \|$ is any reasonable 
norm on $H^1(M,\Sigma,\reals)$ for example
the Hodge norm or the AGY norm 
defined in \S\ref{sec:subsec:app:hodge}. 
By the
notation \index{$g_t$@$(g_t)_*$} $(g_t)_* v$ we mean the action of the geodesic flow
(i.e.\ parallel transport using the Gauss-Manin connection) on the
Hodge bundle $H^1(M,\Sigma,\reals)$. 
We note that the Lyapunov exponents of the geodesic flow (viewed as a
diffeomorphism of $X_0$) are in fact $1+\lambda_i$,  $1 \le i \le k$ and
$-1+\lambda_i$, $1 < i \le k$. 

We have 
\begin{displaymath}
1 = \lambda_1(H^1) > \lambda_2(H^1) > \dots > \lambda_k(H^1) = -1.
\end{displaymath}
It is a standard fact that $\dim \cV_1(H^1) = \dim
  \cV_k(H^1) = 1$, $\cV_1(H^1)$ 
corresponds to the direction of the unipotent $N$ and
$\cV_k(H^1)$ 
corresponds to the direction of $\bar{N}$. 
Let \index{$p$}$p: H^1(M,\Sigma, \reals) \to H^1(M,\reals)$ denote the natural
map. Recall that if $x \in X_0$ denotes the pair $(M,\omega)$, then
\begin{displaymath}
H^1_\perp(x) = \{ \alpha \in H^1(M,\Sigma, \reals) \st p(\alpha) \wedge
\Re(\omega) = p(\alpha) \wedge \Im(\omega) = 0\}. 
\end{displaymath}
Then
\begin{displaymath}
H^1_\perp(x) = \bigoplus_{i
      =2}^{k-1} \cV_i(H^1)(x).  
\end{displaymath}
We note that the subspaces $H^1_\perp(x)$ are equivariant under the
$SL(2,\reals)$ action on $X_0$ (since so is the subspace spanned by $\Re
\omega$ and $\Im \omega$). Since the cocycle preserves the symplectic
form on $p(H^1_\perp)$, we have 
\begin{displaymath}
\lambda_{k+1-i}(H^1) = -\lambda_i(H^1), \quad 1 \le i \le k. 
\end{displaymath}
Let
\begin{displaymath}
\index{$V$@$\cV_{\le i}(H^1)(x)$}\cV_{\le i}(H^1)(x) = \bigoplus_{j=1}^i \cV_j(H^1)(x),\qquad \index{$V$@$\cV_{\ge i}(H^1)(x)$}\cV_{\ge i}(H^1)(x) =
\bigoplus_{j=i}^k \cV_j(H^1)(x). 
\end{displaymath}
Then we have the Lyapunov flags
\begin{displaymath}
\{0\} = \cV_{\le 0}(H^1)(x) \subset \cV_{\le 1}(H^1)(x) \subset \dots \subset
\cV_{\le k}(H^1)(x) = H^1(M,\Sigma,\reals)
\end{displaymath}
and
\begin{displaymath}
\{0\} = \cV_{>k}(H^1)(x) \subset \cV_{>k-1}(H^1)(x) \subset \dots \subset 
\cV_{> 0}(H^1)(x) = H^1(M,\Sigma,\reals).
\end{displaymath}

We record some simple properties of the Lyapunov flags:
\begin{lemma}
\label{lemma:properties:lyapunov:flag}
$ $
\begin{itemize} 
\item[{\rm (a)}] The subspaces $\cV_{\le i}(H^1)(x)$ are locally constant along
  $W^+[x]$, i.e.\ for almost all $x \in X_0$, for almost all $y \in W^+[x]$ close to $x$ we have $\cV_{\le i}(H^1)(y) =
  \cV_{\le i}(H^1)(x)$ for all $1 \le i \le k$. (Here and in (b) we identify
  $H^1(x)$ with $H^1(y)$ using the Gauss-Manin connection).  
\item[{\rm (b)}] The subspaces $\cV_{\ge i}(H^1)(x)$ are locally constant along
  $W^-[x]$, i.e.\ for almost all $x \in X_0$ and for almost all $y \in W^-[x]$ close to $x$ we have
  $\cV_{\ge i}(H^1)(y) = \cV_{\ge i}(H^1)(x)$ for all $1 \le i \le k$. 
\end{itemize}
\end{lemma}

\bold{Proof.} To prove (a), note that
\begin{displaymath}
\cV_{\le i}(H^1)(x) = \left\{ v \in H^1(M,\Sigma,\reals) \st \lim_{t \to \infty} \frac{1}{t}
\log \frac{ \|(g_{-t})_* v\|}{\|v\|} \le -\lambda_i \right\}.
\end{displaymath}
Therefore, the subspace $\cV_{\le i}(H^1)(x)$ depends only on the trajectory
$g_{-t} x$ as $t \to 
\infty$. However, if $y \in W^+[x]$ then $g_{-t} y$ will for large $t$
be close to $g_{-t} x$, and so in view of the affine
structure, $(g_{-t})_*$ will be the same linear map on
$H^1(M,\Sigma,\reals)$ at $x$ and $y$, 
as in 
\S\ref{sec:semi:markov}. 
This implies that $\cV_{\le i}(H^1)(x)= \cV_{\le i}(H^1)(y)$. The proof of property (b) is
identical. 
\qed\medskip

\bold{The action on $H^1_+$ and $H^1_-$.} 
Recall that the bundles $H^1_+$ and $H^1_-$ were defined in
\S\ref{sec:semi:markov}. All of the results of
\S\ref{sec:subsec:lyapunov:flags} also apply to these bundles. Also, 
\begin{displaymath}
\index{$\lambda_i(H^1_+)$}\lambda_i(H^1_+) = 1 + \lambda_i(H^1),
\qquad \index{$\lambda_i(H^1_-)$}\lambda_i(H^1_-) = -1 + \lambda_i(H^1).  
\end{displaymath}
Furthermore, under the natural identification by the identity map, for
all $x \in X_0$, 
\begin{displaymath}
\index{$V$@$\cV_i(H^1_+)(x)$}\cV_i(H^1_+)(x) = \index{$V$@$\cV_i(H^1_-)(x)$}\cV_i(H^1_-)(x) = \cV_i(H^1)(x). 
\end{displaymath}

\subsection{Equivariant measurable flat connections.}
\label{sec:subsec:connection}

Let $L$ be a subbundle of $H_{big}^{(++)}$. \mcc{(or possibly a quotient
bundle of a subbundle of $H_{big}^{(++)}$.)}  
Recall that by Lemma~\ref{lemma:embedded:leaves}, typical leaves of
$W^+$ are simply connected. 
By an equivariant measurable flat $W^+$-connection on $L$ 
we mean a measurable collection of
linear ``parallel transport'' maps:
\begin{displaymath}
F(x,y): L(x) \to L(y)
\end{displaymath}
defined for $\nu$-almost all $x \in X_0$ and $\nu_{W^+[x]}$ almost all
$y \in W^+[x]$ such that
\begin{equation}
\label{eq:def:connection:a}
F(y,z) F(x,y) = F(x,z),
\end{equation}
and 
\begin{equation}
\label{eq:def:connection:b}
(g_t)_* \circ F(x,y) = F(g_t x, g_t y) \circ (g_t)_*. 
\end{equation}
For example, if $L = W^+(x)$, then the Gauss-Manin connection (which
in period local coordinates is the identity map) is an equivariant measurable
flat $W^+$ connection on $H^1$. 
However, there is another important equivariant measurable flat
$W^+$-connection on $H^1$ which we describe below. 

\bold{The maps $P^+(x,y)$ and $P^-(x,y)$.} \index{$P^+(x,y)$} \index{$P^-(x,y)$}
Recall that $\cV_i(H^1)(x) \subset H^1(x)$ 
are the Lyapunov subspaces for the flow
$g_t$. Recall that the $\cV_i(H^1)(x)$ are not locally constant along leaves
of $W^+$, but by Lemma~\ref{lemma:properties:lyapunov:flag}, 
the subspaces $\cV_{\le i}(H^1)(x) = \sum_{j=1}^i \cV_j(W^+)(x)$ are
locally constant along the leaves of $W^+$.  
Now suppose $y \in W^+[x]$. Any vector $v \in \cV_i(H^1)(x)$ can be
written uniquely as
\begin{displaymath}
v = v' + v'' \qquad v' \in \cV_i(H^1)(y), \quad v'' \in \cV_{< i}(H^1)(y). 
\end{displaymath}
Let $P_i^+(x,y): \cV_i(H^1)(x) \to \cV_i(H^1)(y)$ be the linear map sending $v$
to $v'$. Let $P^+(x,y)$ be the unique linear map which restricts to
$P_i^+(x,y)$ on each of the subspaces $\cV_i(H^1)(x)$. We call $P^+(x,y)$
the ``parallel transport'' from $x$ to $y$. The following is immediate
from the definition:
\begin{lemma}
\label{lemma:properties:Pplus}
Suppose $x,y \in W^+[z]$. Then
\begin{itemize}
\item[{\rm (a)}] $P^+(x,y) \cV_i(H^1)(x) = \cV_i(H^1)(y)$. 
\item[{\rm (b)}] $P^+(g_t x, g_t y) = (g_t)_* \circ P^+(x,y) \circ
  (g_t^{-1})_*$.
\item[{\rm (c)}] $P^+(x,y) \cV_{\le i}(H^1)(x) = \cV_{\le i}(H^1)(y)$. If we identify $H^1(x)$
  with $H^1(y)$ using the Gauss-Manin connection, 
 then the map $P^+(x,y)$ is unipotent. 
\item[{\rm (d)}] $P^+(x,z) = P^+(y,z) \circ P^+(x,y)$.  
\end{itemize}
\end{lemma}

Note that the map $P^+$ on $H^1_+$ is the same as on $H^1$, provided
we identify $H^1_+$ with $H^1$ via the identity map.

The statements (b) and (d) imply that the maps $P^+(x,y)$ define an
equivariant measurable flat $W^+$-connection on $H^1$. This connection is
in general different from the Gauss-Manin connection, and is only
measurable.

If $y \in W^-[x]$,
then we can define a similar map which we denote by $P^-(x,y)$. This
yields an equivariant measurable flat $W^-$-connection on $H^1$. 

Clearly the connection $P^+(x,y)$ induces an equivariant measurable
flat $W^+$-connection on $H_{big}^{(++)}$.  This
connection preserves the Lyapunov subspaces of the $g_t$-action on
$H_{big}^{(++)}$, as in Lemma~\ref{lemma:properties:Pplus} (a). 
In view of Proposition~\ref{prop:sublyapunov:locally:constant}
below, the connection $P^+(x,y)$ also induces an equivariant
measurable flat $W^+$-connection on any  $g_t$-equivariant subbundle of
$H_{big}^{(++)}$. 

\bold{Equivariant measurable flat $U^+$-connections.}  Suppose
$U^+[x] \subset W^+[x]$ is a $g_t$-equivariant family of algebraic
subsets, with $U^+[y] = U^+[x]$ for $y \in U^+[x]$. In fact, we will
only consider families compatible with $\nu$ as defined in
Definition~\ref{def:compatible:family}. We denote the conditional
measure of $\nu$ along $U^+[x]$ by $\nu_{U^+[x]}$. In the cases we
will consider, these measures are well defined a.e.\ and are in the
Lebesgue measure class, see \S\ref{sec:divergence:subspaces}.

By an equivariant measurable flat $U^+$-connection on a
bundle $L \subset H_{big}^{(++)}$ we mean a measurable collection 
of linear maps $F(x,y): L(x) \to L(y)$ 
satisfying (\ref{eq:def:connection:a}) and
(\ref{eq:def:connection:b}),
defined for $\nu$-almost all $x \in X_0$ and $\nu_{U^+[x]}$-almost all
$y \in U^+[x]$.

\subsection{The Jordan Canonical Form of a cocycle}
\label{sec:subsec:jordan:form}
$ $

\bold{Zimmer's Amenable reduction.}
The following is a general fact about linear cocycles over an action
of $\reals$ or $\zed$. It is often called ``Zimmer's amenable
reduction''. We state it only for the cases which will be
used. 

\begin{lemma}
\label{lemma:jordan:canonical:form}
Suppose $L_i$ is a $g_t$-equivariant subbundle of
$H_{big}^{(++)}$. (For example, we
could have $L_i(x) = \cV_i(H^1_+)(x)$). Then, there exists a
measurable finite cover $\sigma_{L_i}: X_{L_i} \to X_0$ such that for $\sigma_{L_i}^{-1}(\nu)$-a.e $x \in X_{L_i}$
there exists an invariant flag
\begin{equation}
\label{eq:jordan:flag}
\{0\} = L_{i,0}(x) \subset L_{i,1}(x) \subset \dots \subset
L_{i,n_i}(x) = L_i(x),
\end{equation}
and on each $L_{ij}(x)/L_{i,j-1}(x)$ there exists a nondegenerate quadratic
form \index{$\langle \cdot, \cdot \rangle_{ij,x}$}$\langle \cdot, \cdot \rangle_{ij,x}$ and a cocycle
$\lambda_{ij}: X_{L_i} \cross \reals \to \reals$ such that for all $u, v \in
L_{ij}(x)/L_{i,j-1}(x)$, 
\begin{displaymath}
\langle (g_t)_* u, (g_t)_* v \rangle_{ij,g_t x} = e^{\lambda_{ij}(x,t)} \langle
u, v  \rangle_{ij,x}.
\end{displaymath}
(Note: For each $i$, 
the pullback measures $\sigma_{L_i}^{-1}(\nu)$ is uniquely defined by
the condition that for almost all $x_0 \in X_0$, the 
conditional of $\sigma_{L_i}^{-1}(\nu)$  on the (finite) set 
$\sigma_{L_i}^{-1}(x_0)$ is the normalized counting measure.)
\end{lemma}

\bold{Remark.} The statement of
Lemma~\ref{lemma:jordan:canonical:form} is the assertion that on the
finite cover $X_{L_i}$ one can
make a change of basis at each $x \in X_{L_i}$ so that in the new basis, the
matrix of the cocycle restricted  to $L_i$ is of the form
\begin{equation}
\label{eq:jordan:block}
\begin{pmatrix} C_{i,1}   & *        & \dots & * \\
                   0      & C_{i,2}  & \dots & * \\
                   \vdots & \vdots   & \ddots & * \\
                   0      & 0        & \dots & C_{i,n_i} \\
\end{pmatrix},
\end{equation}
where each $C_{i,j}$ is a conformal matrix (i.e.\ is the composition of
an orthogonal matrix and a scaling factor $\lambda_{ij}$). 

We call a cocycle {\em block-conformal} if all the off-diagonal
entries labeled $\ast$ in (\ref{eq:jordan:block}) are $0$.

\bold{Proof of Lemma~\ref{lemma:jordan:canonical:form}.} 
See \cite{Arnold} (which uses many of the ideas of Zimmer). 
The statement differs slightly from that of \cite[Theorem 5.6]{Arnold}
in that we want the cocycle in each block to be
conformal (and not just block-conformal). However, our statement is in
fact equivalent because we are willing to replace the
original space $X_0$ by a finite cover $X_{L_i}$. 
\qed\medskip

\subsection{Covariantly constant subspaces}
The main result of this subsection is the following:
\begin{proposition}
\label{prop:sublyapunov:new:locally:constant}
Suppose $L$ is a $g_t$-equivariant subbundle over the base $X_0$. We can
write
\begin{displaymath}
L(x) = \bigoplus_i L_i(x), 
\end{displaymath}
where $L_i(x) \equiv \cV_i(L)(x)$ 
is the Lyapunov subspace corresponding to the Lyapunov
exponent $\lambda_i$. Suppose
there exists an equivariant flat measurable $W^+$-connection $F$ on
$L$, such that 
\begin{equation}
\label{eq:Fxy:preserves:lyapunov:subspaces}
F(x,y)L_i(x) = L_i(y). 
\end{equation}
Suppose that $\cM$ is a finite collection of subspaces of $L$ which is
$g_t$-equivariant. Then, for almost all $x \in X_0$ 
and almost all $y \in \gB_0[x]$, 
\begin{displaymath}
F(x,y) \cM(x) = \cM(y),
\end{displaymath}
i.e.\ the collection of subspaces $\cM$ is locally covariantly constant
with respect to the connection $F$. 
\end{proposition}

\bold{Remark.} The same result holds if
$F$ is only assumed to be a measurable $U^+$-connection, and
$\gB_0[x]$ is replaced by  $\cB[x]$.
\medskip

The following is a generalization of
Lemma~\ref{lemma:properties:lyapunov:flag}: 
\begin{corollary}
\label{cor:sublyapunov:locally:constant}
Suppose $M \subset H^1(M,\Sigma, \reals)$ is a $g_t$-equivariant
subbundle over the base $X_0$. Suppose also for a.e $x \in X_0$, 
$\cV_{< i}(x) \subset 
M(x) \subset \cV_{\le i}(x)$.  Then, (up to a set of measure $0$), 
$M(x)$ is locally constant along
$W^+(x)$. 
\end{corollary}

\bold{Proof of Corollary~\ref{cor:sublyapunov:locally:constant}.} 
By Lemma~\ref{lemma:properties:lyapunov:flag}, $L(x)
\equiv \cV_{\le i}(x)/\cV_{< i}(x)$ is locally constant along $W^+[x]$. 
Let $F(x,y)$ denote the Gauss-Manin connection (i.e.\ the identity map)
on $L(x)$. Note that the action of $g_t$ on $L(x)$ has only one
Lyapunov exponent, namely $\lambda_i$. Thus,
(\ref{eq:Fxy:preserves:lyapunov:subspaces}) is trivially satisfied. 
Then, by
Proposition~\ref{prop:sublyapunov:new:locally:constant}, 
$M(x)/\cV_{< i}(x) \subset L(x)$ is locally constant along $W^+[x]$. Since
$\cV_{< i}(x)$ is also locally constant (by
Lemma~\ref{lemma:properties:lyapunov:flag}), this implies that $M(x)$
is locally constant. 
\qed\medskip

\bold{Remark.} Our proof of
Proposition~\ref{prop:sublyapunov:new:locally:constant} is essentially by
reference to \cite[Theorem~1]{Ledrappier:Positivity}.  It 
is given in
\S\ref{sec:starredsubsec:proof:of:prop:sublyapunov:locally:constant}
and can be skipped on first reading.
For similar results in a partially hyperbolic
setting see \cite{Avila:Viana:extremal},
\cite{Avila:Santamaria:Viana}, \cite{Kalinin:Sadovskaya}. 
\medskip

\subsection{Some estimates on Lyapunov subspaces.}
\label{sec:subsec:estimates:Lyapunov:subspaces}
Let $(V, \| \cdot \|_Y)$ be a normed vector space. By a splitting 
$E=(E_1,\dots, E_n)$ of $V$ we mean a direct sum decomposition
\begin{displaymath}
V = E_1 \oplus \dots \oplus E_n
\end{displaymath}
Suppose $E=(E_1, \dots, E_n)$ and $E'=(E_1', \dots, E_n')$ are two
splittings of $V$, with $\dim E_i = \dim E_i'$ for $1 \le i
\le n$. 

We define
\begin{displaymath}
\index{$D^+(E,E')$}D^+(E,E') = \max_{1 \le i \le n} \sup_{ v \in
  \bigoplus\limits_{j \le i} E_j\setminus \{0\}} \inf \left\{ \frac{\| w \|_Y}{\|v\|_Y} \st
v+w \in \bigoplus_{j \le i} E_j', \text{ and } w \in \bigoplus_{j > i} E_j \right\}, 
\end{displaymath}
and
\begin{displaymath}
\index{$D^-(E,E')$}D^-(E,E') = \max_{1 \le i \le n} \sup_{ v \in
  \bigoplus\limits_{j \ge i} E_j\setminus \{0\}} \inf \left\{ \frac{\| w \|_Y}{\|v\|_Y} \st
v+w \in \bigoplus_{j \ge i} E_j', \text{ and } w \in \bigoplus_{j < i} E_j \right\}. 
\end{displaymath}
Note that $D^+(E,E')$ depends on $E'$ only via the flag $\bigoplus_{j
  \le i} E_j'$, $1 \le i \le n$. Similarly, $D^-(E,E')$ 
depends on $E'$ only via the flag $\bigoplus_{j
  \ge i} E_j'$, $1 \le i \le n$.
Also $D^+(E,E') = D^-(E,E') = 0$ 
if $E = E'$, and $D^+(E,E') = \infty$ if some $\bigoplus_{j \le i} E_j'$
has non-trivial intersection with $\bigoplus_{j > i} E_j$.

In this subsection, we write $\cV_i(x)$ for $\cV_i(H^1)(x)$, etc. 
For almost all $x$ in  $\tilde{X}_0$, we have the splitting
\begin{displaymath}
H^1(x) = \cV_1(x) \oplus \dots \oplus \cV_n(x).
\end{displaymath}

For $x, y \in \tilde{X}_0$, we have the Gauss-Manin connection 
$P^{GM}(x,y)$, which is a linear map from $H^1(x)$ to $H^1(y)$ (see
\S\ref{sec:outline:step1}). Let
\begin{displaymath}
\index{$D^+(x,y)$}D^+(x,y) = D^+((\cV_1(x), \dots, \cV_n(x)),
(P^{GM}(y,x)\cV_1(y), \dots, P^{GM}(y,x) \cV_n(y)). 
\end{displaymath}
\begin{displaymath}
\index{$D^-(x,y)$}D^-(x,y) = D^-((\cV_1(x), \dots, \cV_n(x)),
(P^{GM}(y,x)\cV_1(y), \dots, P^{GM}(y,x) \cV_n(y)). 
\end{displaymath}

\bold{Distance between subspaces.}
For a subspace $V$ of $H^1(x)$, let \index{$SV$}$SV$ denote the intersection of
$V$ with the unit ball in the AGY norm. 

For subspaces $V_1, V_2$ of $H^1(x)$, we define
\begin{equation}
\label{eq:def:distance:between:subspaces}
\index{$d_Y(V_1,V_2)$}d_Y(V_1,V_2) 
= \text{The Hausdorff distance between $SV_1$ and $SV_2$}
\end{equation}
measured with respect to the AGY norm at $x$.

\begin{lemma}
\label{lemma:one:sided:distance}
There exists a continuous function $C_0: X_0 \to \reals^+$ such that for
subspaces $V_1, V_2$ of $H^1(x)$ of the same dimension,
\begin{displaymath}
C_0(x)^{-1} d_Y(V_1,V_2) \le 
\delta_Y(V_1,V_2) \le d_Y(V_1,V_2),
\end{displaymath}
where
\begin{displaymath}
\index{$\delta_Y(\cdot,\cdot)$}\delta_Y(V_1,V_2) = \max_{ v_1 \in SV_1} \min_{v_2 \in SV_2} \|v_1- v_2\|_Y.
\end{displaymath}
\end{lemma}

\bold{Proof.} Since$d_Y(V_1,V_2) = \max(\delta_Y(V_1,V_2), \delta_Y(V_2,V_1))$,
the inequality on the left follows immediately from the
definition of the Hausdorff distance. To prove the inequality on
the right it is enough to show that for some continuous function $C_0:
X_0 \to \reals^+$, 
\begin{equation}
\label{eq:delta:Y:asymp:lower}
C_0(x)^{-1} \delta_Y(V_2,V_1) \le \delta_Y(V_1,V_2). 
\end{equation}
To prove (\ref{eq:delta:Y:asymp:lower}),
pick some arbitrary inner product $\langle\cdot,\cdot\rangle_0$
on $H^1(M,\Sigma,\reals)$, and let $\| \cdot \|_0$ be the associated
norm. Then, there exists a continuous function $C_1: X_0 \to \reals^+$
such that for all $v \in H^1(x)$, 
\begin{displaymath}
C_1(x)^{-1} \|v\|_0 \le \|v\|_Y \le C_1(x) \|v\|_0.
\end{displaymath}
Let $\delta_0(\cdot, \cdot)$ and $d_0(\cdot, \cdot)$ be the analogues
of $\delta_Y(\cdot, \cdot)$ and $d_Y(\cdot, \cdot)$ for the norm $\|
\cdot \|_0$. Then, it is enough to prove that there exists a constant
$c_2 > 0$ depending only on the dimension such that for subspaces $V_1$, $V_2$
of equal dimension,  
\begin{equation}
\label{eq:goal:subspace:comparasion}
c_2 \, \delta_0(V_2,V_1) \le \delta_0(V_1,V_2). 
\end{equation}
For subspaces $U,V$ of equal dimension $n$, let $u_1, \dots, u_n$ and
$v_1, \dots, v_n$ be orthonormal bases for $U$ and $V$
respectively. Then, we have
\begin{equation}
\label{eq:subspace:basis:identity}
\left(\sum_{i=1}^n \inf_{v \in V} \| u_i - v\|_0^2
\right)^{1/2} = \left(n - \sum_{i=1}^n \sum_{j=1}^n \langle u_i,
  v_j\rangle^2_0 \right)^{1/2}
\end{equation}
Note that the expression on the left in
(\ref{eq:subspace:basis:identity}) is independent of the basis for
$V$, and the expression on the right of
(\ref{eq:subspace:basis:identity}) is symmetric in $U$ and $V$. Thus,
the expression in (\ref{eq:subspace:basis:identity}) is independent of
the basis for $U$ as well, and thus defines a function
$d_H(U,V)$. (This function is called the Frobenius or chordal distance
between subspaces, see e.g. \cite{Encyclopaedia:distances},
\cite{WWF}).

From the expression on the left of (\ref{eq:subspace:basis:identity})
it is clear that there exists a constant $c_3$ depending only on the
dimension so that
\begin{displaymath}
c_3 \, d_H(V_1,V_2) \le d_0(V_1,V_2) \le c_3^{-1} d_H(V_1,V_2). 
\end{displaymath}
Since $d_H(V_1,V_2) = d_H(V_2,V_1)$,
(\ref{eq:goal:subspace:comparasion}) follows. 
\qed\medskip

\begin{lemma}
\label{lemma:subspaces:stay:close}
There exists $\alpha > 0$ depending only on the Lyapunov spectrum, and 
a function $C: X_0 \to \reals^+$ finite almost
everywhere such that the following holds:
\begin{itemize}
\item[{\rm (a)}] For all $t > 0$, and all 
$x \in \tilde{X}_0$, and all $y \in
\tilde{X}_0$ such that $d^{X_0}(g_s x, g_s y) \le 1/100$ for $0 \le s \le t$, we have,
for all $1 \le i \le n$, 
\begin{displaymath}
d_Y(\cV_{\le i}(g_t x), P^{GM}(g_t y,g_t x)\cV_{\le i}(g_t y)) \le \min_{0
  \le s \le t} C(g_s x) (1+D^+(x,y)) e^{-\alpha t}. 
\end{displaymath}
\item[{\rm (b)}] For all $t > 0$, and all 
$x \in \tilde{X}_0$, and all $y \in
\tilde{X}_0$ such that $d^{X_0}(g_{-s} x, g_{-s} y) \le 1/100$ 
for $0 \le s \le t$, we have,
for all $1 \le i \le n$, 
\begin{displaymath}
d_Y(\cV_{\ge i}(g_{-t} x), P^{GM}(g_{-t} y,g_{-t} x)\cV_{\ge i}(g_{-t} y)) 
\le \min_{0 \le s \le t} C(g_{-s} x) (1+D^-(x,y)) e^{-\alpha t}. 
\end{displaymath}
\end{itemize}
\end{lemma}

The proof of Lemma~\ref{lemma:subspaces:stay:close} is a
straightforward but tedious argument using the Osceledets
multiplicative ergodic theorem. It is done in 
\S\ref{sec:starredsubsec:proof:of:lemma:subspaces:stay:close}.

\begin{lemma}
\label{lemma:can:bound:D}
There exists a function $C_3: X_0 \to \reals^+$ finite almost
everywhere, such that for all $x \in \tilde{X}_0$, all $y \in W^-[x]$
with $d^{X_0}(x,y) < 1/100$ we have $D^+(x,y) \le C_3(x) C_3(y)$. Similarly, 
for all $x \in \tilde{X}_0$, all $y \in W^+[x]$
with $d^{X_0}(x,y) < 1/100$ we have $D^-(x,y) \le C_3(x) C_3(y)$.
\end{lemma}

\bold{Proof of Lemma~\ref{lemma:can:bound:D}.}
For $\epsilon > 0$, let $K_\epsilon \subset X_0$ be a compact set with
measure  at
least $1-\epsilon$ on which the functions $x \to \cV_i(x)$ are
continuous. Then there exists $\rho = \rho(\epsilon)$ such that if $x'
\in \pi^{-1}(K_\epsilon)$, $y' \in W^-[x] \cap \pi^{-1}(K_\epsilon)$ and $d^{X_0}(x'.y') < \rho$
then $D^+(x',y') < 1$. Then, by the
Birkhoff ergodic theorem and Lemma~\ref{lemma:forni}, 
there exists a compact $K_\epsilon' \subset
X_0$ with $\nu(K_\epsilon') > 1-2\epsilon$ and $C_2 = C_2(\epsilon)$
such that for all $x \in \pi^{-1}(K_\epsilon')$, all $y \in
W^-[x] \cap \pi^{-1}(K_\epsilon')$ with $d^{X_0}(x,y) < 1/100$ there exists $C_2(\epsilon) <
t' < 2 C_2(\epsilon)$ with $g_{t'} x \in K_\epsilon$, $g_{t'} y \in K_\epsilon$
and $d^{X_0}(x,y) < \rho(\epsilon)$. Thus, $D^+(g_{t'} x, g_{t'} y) < 1$, which
implies that $D^+(x,y) < C_2' = C_2'(\epsilon)$. Without loss of
generality, we may assume that $C_2' \ge 1$ and that
$K'_\epsilon$ and $C_2'(\epsilon)$ both
decrease as functions of $\epsilon$. 
Now for $x \in X_0$, let $\Upsilon(x) = \{ \epsilon \st x \in
K_\epsilon' \}$, and let
\begin{displaymath}
C_3(x) = \inf \{ C_2'(\epsilon) \st \epsilon \in \Upsilon(x) \}. 
\end{displaymath}
The proof of the second assertion is identical. 
\qed\medskip

\begin{corollary}
\label{cor:staying:close}
There exists a measurable function $C_1: X_0 \to \reals^+$ finite a.e such
that if $x \in X_0$, $y \in W^-[x]$ with $d^{X_0}(x,y) < 1/100$, we
have for all $t > 0$, 
\begin{equation}
\label{eq:Pminus:fast:to:zero}
\| P^-(g_t x, g_t y) P^{GM}(g_t y, g_t x) - I \|_Y \le C_1(x) C_1(y)
e^{-\alpha t},
\end{equation}
where $\alpha > 0$ depends only on the Lyapunov spectrum. Consequently, 
for almost all $x \in X_0$, and almost all $y \in W^-[x]$, 
\begin{equation}
\label{eq:Pminus:to:zero}
\lim_{t \to \infty} \| P^-(g_t x, g_t y) P^{GM}(g_t y, g_t x) - I \|_Y
= 0. 
\end{equation}
The same assertions hold if $W^-$ is replaced by $W^+$, 
$g_t$ by $g_{-t}$ and $P^-$ by $P^+$. 
\end{corollary}

\bold{Proof of Corollary~\ref{cor:staying:close}.}
Let $C_1(x) = C(x) C_3(x)$, where $C(\cdot)$ is as in
Lemma~\ref{lemma:subspaces:stay:close} and $C_3(\cdot)$ is as in
Lemma~\ref{lemma:can:bound:D}. Then, by
Lemma~\ref{lemma:subspaces:stay:close} and
Lemma~\ref{lemma:can:bound:D},
\begin{displaymath}
d_Y(\cV_{\le i}(g_t x), P^{GM}(g_t y, g_t x) \cV_{\le i}(g_t y)) \le
C_1(x) C_1(y) e^{-\alpha t}. 
\end{displaymath}
Since by
Lemma~\ref{lemma:properties:lyapunov:flag}, $\cV_{\ge i}(x) = P^{GM}(y,x) \cV_{\ge i}(y)$, 
we get, for $t > 0$, 
\begin{displaymath}
d_Y(\cV_i(g_t x), P^{GM}(g_t y,g_t x) \cV_i(g_t y)) \le C_1(x) C_1(y)
e^{-\alpha t}. 
\end{displaymath}
This, by the definition of $P^-(x,y)$, implies that
(\ref{eq:Pminus:fast:to:zero}) holds
as required. Even if 
we do not assume that $d^{X_0}(x,y) < 1/100$, then for almost
all $x$ and almost all $y \in W^-[x]$, for $t$ large enough $d^{X_0}(g_t x,
g_t y) < 1/100$, and thus, in view of (\ref{eq:Pminus:fast:to:zero}),
(\ref{eq:Pminus:to:zero}) holds.  
\qed\medskip

\subsection{The cover $X$.}
\label{sec:subsec:the:cover:X}

Let $L = H_{big}$ viewed as a bundle over $X_0$. Let $L_i = \cV_i(L)$.
By Lemma~\ref{lemma:jordan:canonical:form}, there exists a 
measurable finite cover \index{$X$}$X$ of $X_0$ such that 
Lemma~\ref{lemma:jordan:canonical:form} holds on $X$ for all the
$L_i$.  We always assume that the degree of
the covering map $\index{$\sigma_0$}\sigma_0: X \to X_0$ is as small as
possible.  

\bold{The set $\Delta(x_0)$.}
For $x_0 \in X_0$, let \index{$\Delta_i(x_0)$}$\Delta_i(x_0)$ 
denote the set of flags
\begin{displaymath}
\Delta_i(x_0) = \left\{ \{0\} = L_{i,0}(x) \subset L_{i,1}(x) \subset
  \dots \subset L_{i,n_i}(x) = L_i(x) \st x \in \sigma_0^{-1}(x_0) \right\}.
\end{displaymath}
Let \index{$\Delta(x_0)$}$\Delta(x_0)$ 
denote the Cartesian product of the $\Delta_i(x_0)$. 
Then, we can think of a point $x \in X$ as a pair $(x_0, \gF)$
where $\index{$F$@$\gF$}\gF \in \Delta(x_0)$. 

\bold{The measure $\nu$ on $X$.}
We can use $\sigma_0$ to define a pullback of the invariant measure $\nu$ on
$X_0$ to $X$, by requiring that the pushforward of the pullback
measure by $\sigma_0$ is $\nu$, and that the conditionals of the pullback
measure on the fibers of $\sigma_0$ are the (normalized) counting measure. 
We abuse notation by denoting the pullback measure also
by \index{$\nu$}$\nu$. 

\begin{lemma}
\label{lemma:cover:X:ergodic}
The measure $\nu$ is ergodic for the action of $g_t$ on $X$. 
\end{lemma}

\bold{Proof.} Suppose $E$ is a $g_t$-invariant set of $X$ with $\nu(E)
> 0$. Then by the
ergodicity of the action of $g_t$ on $X_0$, $\sigma(E)$ is
conull. Let $N(x_0)$ denote the cardinality of $\sigma_0^{-1}(x_0)
\cap E$. Then, again by the ergodicity of $g_t$, $N(x_0)$ is constant
almost everywhere. If $E$ does not have full measure, then we have
that $N(x_0)$ is smaller than the degree of the cover
$\sigma_0$. Then, we could replace $X$ by $E$, contradicting the
assumption that the degree of the covering map $\sigma_0$ is as small
as possible. 
\qed\medskip

\bold{The space $\tilde{X}$.} Recall that $\tilde{X}_0$ is the
universal cover of $X_0$. Let \index{$X$@$\tilde{X}$}$\tilde{X}$
denote the cover of $\tilde{X}_0$ corresponding to the cover
$\sigma_0: X \to X_0$. More precisely, 
\begin{displaymath}
\tilde{X} = \{ (x_0, \gF) \st x_0 \in \tilde{X}_0, \ \gF \in
\Delta(x_0) \}. 
\end{displaymath}
We denote the covering map from $\tilde{X}$ to $\tilde{X}_0$
again by \index{$\sigma_0$}$\sigma_0$.

\bold{Stable and Unstable manifolds for $X$ and $\tilde{X}$.}
Suppose $x = (x_0, \gF) \in \tilde{X}$. We define
\begin{gather}
\label{eq:def:Wplus:X}
\index{$W^+[x]$}W^+[x] = \{ (y_0, \gF' ) \in \tilde{X} \st y_0 \in W^+[x_0], \text{ and }
\gF' = P^+(x_0,y_0) \gF \}. \\
\label{eq:def:Wminus:X}
\index{$W^-[x]$}W^-[x] = \{ (y_0, \gF' ) \in \tilde{X} \st y_0 \in W^-[x_0], \text{ and }
\gF' = P^-(x_0,y_0) \gF \}.
\end{gather}
This definitions make sense, since by
Proposition~\ref{prop:sublyapunov:new:locally:constant}, 
\begin{displaymath}
P^+(x_0,y_0) \Delta(x_0) = \Delta(y_0) \qquad\text{ for $y_0 \in W^+[x_0]$,} 
\end{displaymath}
\begin{displaymath}
P^-(x_0,y_0) \Delta(x_0) = \Delta(y_0) \qquad\text{ for $y_0 \in W^-[x_0]$.}
\end{displaymath}

\bold{Remark.} Even though $\tilde{X}$ itself does not have a manifold
structure, for almost all $x \in \tilde{X}$, the sets $W^+[x]$ and $W^-[x]$
have the structure of an affine manifold (intersected with a
  set of full measure in $\tilde{X}$), see
Lemma~\ref{lemma:embedded:leaves}.
Lemma~\ref{lemma:X:stables:unstables} below asserts that
these can be interpreted as the strong stable and strong unstable
manifolds for the action of $g_t$ on $\tilde{X}$. 

\bold{Notation.} If $x \in \tilde{X}$ and 
$V$ is a subspace of $W^+(x)$ or $W^-(x)$  we
write
\begin{displaymath}
\index{$V[x]$}V[x] = \{  y \in W^\pm[x] \st y - x \in V(x) \}.  
\end{displaymath}

\bold{The ``distance'' $d^X( \cdot, \cdot)$.} 
For $x = (x_0, \gF) \in \tilde{X}$, and $y = (y_0, \gF') \in
\tilde{X}$ and $y \in W^+[x]$ or $W^-[x]$ define
\begin{equation}
\label{eq:def:dX}
\index{$d^X(x,y)$}d^X(x,y) = d^{X_0}(x_0,y_0) + d_Y(\gF, P^{GM}(y_0,x_0)\gF'),
\end{equation}
where we extend the distance $d_Y$ between subspaces defined in 
(\ref{eq:def:distance:between:subspaces}) to a
distance between flags. 

\begin{lemma}
\label{lemma:X:stables:unstables}
For almost all $x \in \tilde{X}$ and almost all $y \in W^+[x]$, $d^X(g_t x, g_t
y) \to 0$ as $t \to -\infty$. 
Similarly, for almost all $x \in \tilde{X}$ 
and almost all $y \in W^-[x]$, we have
$d^X(g_t x, g_t y) \to 0$ as $t \to \infty$.  
\end{lemma}

\bold{Proof.} This follows immediately from
Corollary~\ref{cor:staying:close}.  
\qed\medskip

\bold{Notational Convention.}
If $f$ is an object on $X_0$, and $x \in X$, we write $f(x)$ instead of
$f(\sigma_0(x))$. Thus, we can define $\cV_i(H_{big})(x)$ for $x \in
X$, $P^+(x,y)$ for $x \in X$ and $y \in W^+[x]$, etc.
Also, if $x \in \tilde{X}$, we write $f(x)$ instead of $f(\pi \circ
\sigma_0(x))$ etc. 
 

\bold{The partitions $\gB_t$ of $X$.}
Suppose $x = (x_0, \gF) \in X$. 
We define \index{$\gB_0[x]$}
\begin{displaymath}
\index{$\gB_t[x]$}\gB_t[x] = \{ (y_0,\gF') \st y_0 \in \gB_t[x_0], \quad
\gF' = P^+(x_0,y_0) \gF \}. 
\end{displaymath}
Then $\gB_t$ is a measurable partition of $X$ subordinate to $W^+$. In
a similar way, we can define sets \index{$J[x]$}$J[x]$ for $x \in X$ 
and \index{$E^+[c]$}$E^+[c]$ for $c \in \sigma_0^{-1}(\cC_1)$, where
$\cC_1$ is as in Proposition~\ref{prop:semi:markov}. 
Proposition~\ref{prop:semi:markov} and all subsequent 
results of \S\ref{sec:semi:markov} apply to $X$ as well as $X_0$. 
\medskip

The following is an alternative version of
Proposition~\ref{prop:sublyapunov:new:locally:constant} adapted to the
cover $X$. 
\begin{prop}
\label{prop:sublyapunov:locally:constant}
Suppose $L$ is a $g_t$-equivariant subbundle of $H_{big}$.
For almost all $x \in X$, we can write
\begin{displaymath}
L(x) = \bigoplus_i L_i(x), 
\end{displaymath}
where $L_i(x)$ is the Lyapunov subspace corresponding to the Lyapunov
exponent $\lambda_i$. Suppose
there exists an equivariant flat measurable $W^+$-connection $F$ on
$L$, such that 
\begin{displaymath}
F(x,y)L_i(x) = L_i(y), 
\end{displaymath}
and that $M \subset L$ is a
$g_t$-equivariant subbundle. Then, 
\begin{itemize}
\item[{\rm (a)}] For almost all $y \in \gB_0[x]$, 
\begin{displaymath}
F(x,y) M(x) = M(y),
\end{displaymath}
i.e.\ the subbundle $M$ is locally covariantly constant with respect to the
connection $F$. 
\item[{\rm (b)}] For all $i$, the
decomposition 
(\ref{eq:jordan:flag}) of $L_i$ is 
locally covariantly constant along $W^+$,
i.e.\ for $\nu_{W^+[x]}$-almost all $y \in \gB_0[x]$, 
for all $i \in I$ and
for all $1 \le j \le n_i$, 
\begin{equation}
\label{eq:Lijy:Lijx}
L_{ij}(y) = F(x,y) L_{ij}(x).
\end{equation}
Also, up to a scaling factor, the quadratic forms $\langle \cdot,
\cdot \rangle_{i,j}$ are locally covariantly constant along $W^+$,
i.e.\ for almost all $y \in
\gB_0[x]$, and for $v, w
\in L_{ij}(x)/L_{i,j-1}(x)$, 
\begin{equation}
\label{eq:locally:constant:quadratic:form}
\langle F(x,y)v, F(x,y)w \rangle_{ij,y} = c(x,y) \langle v, w 
\rangle_{ij,x}.
\end{equation}
\end{itemize}
\end{prop}

Proposition~\ref{prop:sublyapunov:locally:constant} will be proved in
\S\ref{sec:starredsubsec:proof:of:prop:sublyapunov:locally:constant}.
The proof also shows the following:
\begin{remark}
\label{remark:Uplus:connection}
Proposition~\ref{prop:sublyapunov:locally:constant} applies also to
$U^+$-connections, provided the measure along $U^+[x]$ is in the Lebesgue
measure class, and provided that in the statement, 
the set $\gB_0[x]$ is replaced by $\index{$B$@$\cB[x]$}\cB[x] = \gB_0[x] \cap U^+[x]$. 
\end{remark}

\subsection{Dynamically defined norms}
\label{sec:subsec:dynamical:norms}
In this subsection we work on the cover $X$. We
define a norm on \index{abs@$\norm{\cdot}$}
$\| \cdot \|$ on
$H_{big}^{(++)}$, which has some advantages over the AGY norm $\|
\cdot \|_Y$. 

\bold{Notation.} In \S\ref{sec:subsec:dynamical:norms}
we let $L$ denote the entire bundle
$H_{big}^{++}$, write $L_i$ for $\cV_i(L)$, and for each $i$, 
consider the decomposition (\ref{eq:jordan:flag}). 

\bold{The function $\Xi(x)$.} For $x \in X$, let
\begin{displaymath}
\Xi^+(x) = \sup_{ij} \sup \left\{ \langle v, v \rangle_{ij,x}^{1/2},
  \st v \in
L_{ij}(x)/L_{i,j-1}(x), \ \| v \|_{Y,x} = 1 \right\}, 
\end{displaymath}
and let 
\begin{displaymath}
\Xi^-(x) = \inf_{ij} \inf \left\{ \langle v, v \rangle_{ij,x}^{1/2},
  \st v \in
L_{ij}(x)/L_{i,j-1}(x), \ \| v \|_{Y,x} = 1 \right\}.
\end{displaymath}
Let
\begin{displaymath}
\index{$\Xi(x)$}
\Xi(x) = \Xi^+(x)/\Xi^-(x).
\end{displaymath}
We have $\Xi(x) \ge 1$ for all $x \in X$. For $x_0 \in X_0$, we define
$\Xi(x_0)$ to be $\max_{x \in \sigma_0^{-1}(x_0)} \Xi(x)$.  
\medskip

Let $d_Y(\cdot, \cdot)$ be the distance between subspaces defined in
(\ref{eq:def:distance:between:subspaces}).
Let $\cC_0 \subset X_0$ with $\nu(\cC_0)> 0$ and 
  $M_0 \ge 1$ be chosen later. (We will choose them immediately before
  Lemma~\ref{lemma:Lyapunov:exponents:on:bH} in
  \S\ref{sec:divergence:subspaces}.)
\begin{lemma}
\label{lemma:exists:C:T0}
Fix $\epsilon > 0$ smaller than $\min_i |\lambda_i|$, and smaller than
$\min_{i\ne j} |\lambda_i - \lambda_j|$, where the $\lambda_i$ are the
Lyapunov exponents of $H_{big}^{(++)}$. 
There exists a compact subset $\cC \subset \cC_0 \subset X_0$
with $\nu(\cC) > 0$ and a function $T_0: \cC \to \reals^+$ with
$T_0(x) < \infty$ for $\nu$ a.e.\ $x \in \cC$ such that the
following hold:
\begin{itemize}
\item[{\rm (a)}] There exists $\sigma > 0$ such that for all $c \in
  \cC$, and any subset $S$ of the Lyapunov exponents, 
\begin{displaymath}
d_Y(\bigoplus_{i \in S} L_i(c), \bigoplus_{j \not\in S} L_j(c)) \ge \sigma. 
\end{displaymath}
\item[{\rm (b)}] There exists $M' > 1$ such that for all $c \in \cC$,
  $\Xi(c) \le M'$. 
\item[{\rm (b')}] There exists a constant $M'' < \infty$ such that
for all $x \in \pi^{-1}(\cC)$, for all $y \in \pi^{-1}(\cC) \cap
W^+[x]$ with $d^{X_0}(x,y) < 1/100$, \mccc{say this better} the
Gauss-Manin connection $P^{GM}$ satisfies the estimate: 
\begin{displaymath}
\|P^{GM}(x, y)\|_Y \equiv \sup_{v
  \ne 0} \frac{\|P^{GM}(x, 
  y)v\|_{Y,y}}{\|v\|_{Y,x}} \le M''.
\end{displaymath}
\item[{\rm (c)}] For 
all $c \in \cC$, for all $t > T_0(c)$ and for any subset $S$ of the
Lyapunov spectrum, 
\begin{displaymath}
d_Y(\bigoplus_{i \in S} L_i(g_{-t} c), \bigoplus_{j \not\in S} L_j(g_{-t}
c)) \ge e^{-\epsilon t}. 
\end{displaymath}
Hence, for all $c \in \cC$ and
all $t > T_0(c)$ and all $c' \in \cC \cap W^+[g_{-t} c]$ with
$d^{X_0}(g_{-t} c,c') < 1/100$, 
\begin{equation}
\label{eq:bound:norm:Pplus}
M_0^{-2}\rho_1  e^{-\epsilon t} \le \|P^+(g_{-t} c, c')\|_Y \equiv \sup_{v
  \ne 0} \frac{\|P^+(g_{-t} c, 
  c')v\|_{Y,c'}}{\|v\|_{Y, g_{-t} c}} 
 \le M_0\rho_1^{-1} e^{\epsilon t},
\end{equation}
where $\rho_1 = \rho_1(M', \sigma, M'',M_0) > 0$. 
\item[{\rm (d)}] There exists $\rho > 0$ such that for all $c \in \cC$,
for all $t > T_0(c)$,   for all $i$ and all $v \in L_i(c)$, 
\begin{displaymath}
e^{-(\lambda_i+\epsilon) t} \rho_1 \rho^2 \|v\|_{Y,c} \le \|g_{-t} v\|_{Y,g_{-t}c} \le \rho_1^{-1} \rho^{-2} e^{-(\lambda_i-\epsilon) t} \|v\|_{Y,c}.
\end{displaymath}
\end{itemize}
\end{lemma}

\bold{Proof.} Parts (a) and
(b) hold since the inverse of the angle between Lyapunov subspaces and
the ratio of the norms are finite a.e., therefore
bounded on a set of almost full measure.  To see (c), note that by the
Osceledets multiplicative ergodic theorem, \cite[Theorem S.2.9
(2)]{Katok:Hasselblatt} for $\nu$-a.e.\ $x \in X_0$, 
\begin{displaymath}
\lim_{t \to \infty} \frac{1}{t} \log | \sin \angle(\bigoplus_{i \in S}
L_i(g_{-t}x),\bigoplus_{j \not\in S} L_j(g_{-t}x))| = 0. 
\end{displaymath}
Also, (d) follows immediately from the multiplicative ergodic
theorem. 
\qed\medskip

We now choose the set $\cC$ and the function $T_0$ of
Proposition~\ref{prop:semi:markov} and Lemma~\ref{lemma:semi:markov}
to be as in Lemma~\ref{lemma:exists:C:T0}.

The main result of this subsection is the following:
\begin{proposition} 
\label{prop:properties:dynamical:norm:Hbig}
For almost all $x \in X$ there
exists an inner product \index{$\langle \cdot, \cdot \rangle_x$}$\langle \cdot, \cdot \rangle_x$ on
$H_{big}^{(++)}(x)$ (or on any bundle for which the
  conclusions of Lemma~\ref{lemma:exists:C:T0} hold)
with the following properties: 
\begin{itemize}
\item[{\rm (a)}] For a.e.\ $x \in X$, the distinct eigenspaces $L_i(x)$
  are orthogonal. 
\item[{\rm (b)}] Let $L'_{ij}(x)$ denote the orthogonal complement, 
relative to the inner product $\langle \cdot, \cdot \rangle_x$  
of $L_{i,j-1}(x)$ in $L_{ij}(x)$. Then, for a.e.\ $x \in X$, all $t \in
\reals$  and all $v \in L_{ij}'(x) \subset H_{big}^{(++)}(x)$,
\begin{displaymath}
(g_t)_* v = e^{\index{$\lambda_{ij}(x,t)$}\lambda_{ij}(x,t)} v' + v'',
\end{displaymath}
where $\lambda_{ij}(x,t) \in \reals$, $v' \in L'_{ij}(g_t x)$, $v'' \in L_{i,j-1}(g_t x)$, 
and $\| v' \| = \| v\|$. Hence (since $v'$ and $v''$ are orthogonal), 
\begin{displaymath}
\|(g_t)_* v\| \ge e^{\lambda_{ij}(x,t)} \|v\|. 
\end{displaymath}
\item[{\rm (c)}] There exists a constant $\kappa > 1$ such that for
  a.e.\ $x \in X$ and for all $t > 0$, 
\begin{displaymath}
\kappa^{-1} t \le \lambda_{ij}(x,t) \le \kappa t.
\end{displaymath}
\item[{\rm (d)}] There exists a constant $\kappa > 1$ such that
for a.e $x \in X$ and for all $v \in H_{big}^{(++)}(x)$, and all $t
\ge 0$, 
\begin{displaymath}
e^{\kappa^{-1} t } \|v\| \le \| (g_t)_* v \| \le e^{\kappa t} \|
v\|. 
\end{displaymath}
\item[{\rm (e)}] For a.e.\ $x \in X$, and a.e.\ $y \in \gB_0[x]$ 
and all $t \le 0$, 
\begin{displaymath}
\lambda_{ij}(x,t) = \lambda_{ij}(y,t). 
\end{displaymath}
\item[{\rm (f)}]   For a.e.\ $x \in X$, a.e.\ $y \in \gB_0[x]$, and any
$v,w \in H_{big}^{(++)}(x)$,
\begin{displaymath}
\langle P^+(x,y) v, P^+(x,y) w \rangle_y = \langle v,w \rangle_x.
\end{displaymath}
\end{itemize}
\end{proposition}
We often omit the subscript from \index{$\langle \cdot, \cdot \rangle$}
$\langle \cdot, \cdot \rangle_x$ and
from the associated norm \index{abs@$\norm{\cdot}_x$}$\| \cdot \|_x$.  
\medskip

The inner product $\langle \cdot, \cdot \rangle_x$ is first defined
for  $x \in E^+[c]$ for $c \in \sigma_0^{-1}(\cC_1)$ (in the notation of
\S\ref{sec:semi:markov}, see also \S\ref{sec:subsec:the:cover:X}).  
We then interpolate between $x \in E^+[c]$ and $g_{-t(c)} x$ (again in
the notation of \S\ref{sec:semi:markov}). The details of the proof of
Proposition~\ref{prop:properties:dynamical:norm:Hbig}, which can be
skipped on first reading, are given in 
\S\ref{sec:starredsubsec:proof:of:prop:properties:dynamical:norm:Hbig}.

\bold{The dynamical norm $\| \cdot \|$ on $X_0$.}
The dynamical inner product $\langle \cdot, \cdot \rangle_x$ and the
dynamical norm $\| \cdot \|_x$ of
Proposition~\ref{prop:properties:dynamical:norm:Hbig} are defined for
$x \in X$. For $x_0 \in X_0$, and $v,w \in H_{big}(x_0)$ we define
\begin{equation}
\label{eq:def:dynamical:norm:X0}
\index{$\langle\cdot,\cdot\rangle_{x_0}$}\langle v,w \rangle_{x_0} = \frac{1}{|\sigma_0^{-1}(x_0)|} \sum_{x \in \sigma_0^{-1}(x_0)} \langle
  v, w \rangle_x, \qquad 
\index{abs@$\norm{\cdot}_{x_0}$}\|v\|_{x_0} = \langle v,v\rangle_{x_0}^{1/2}.
\end{equation}

\begin{remark}
\label{remark:dynamical:norm:X0}
The inner product and norm on $X_0$ satisfy properties (a) and
(d) of Proposition~\ref{prop:properties:dynamical:norm:Hbig}. 
\end{remark}

\mcc{MAKE SURE THIS LEMMA IS QUOTED IN THE RIGHT PLACES}

\begin{lemma}
\label{lemma:hodge:norm:vs:dynamical:norm}
For every $\delta > 0$ there exists a compact subset $K(\delta)
\subset X_0$ with
$\nu(K(\delta)) > 1-\delta$ and a number $C_1(\delta) < \infty$ such
that for all $x \in K(\delta)$ and all $v$ on $H_{big}^{(++)}(x)$ or
$H_{big}^{(--)}(x)$, 
\begin{displaymath}
C_1(\delta)^{-1} \le \frac{\|v\|_{x}}{\|v\|_{Y,x}} \le C_1(\delta),
\end{displaymath}
where $\|\cdot \|_x$ is the dynamical norm defined in this subsection
and $\|\cdot\|_{Y,x}$ is the AGY norm. 
\end{lemma}
\bold{Proof.} Since any two norms on a finite dimensional vector space
are equivalent, there exists a function $\index{$\Xi_0(x)$}\Xi_0: X \to \reals^+$ finite
a.e.\ such that for all $x \in X$ and all $v \in H_{big}^{(++)}(x)$,
\begin{displaymath}
\Xi_0(x)^{-1} \|v\|_{Y,x} \le \|v\|_x \le \Xi_0(x) \|v\|_{Y,x}. 
\end{displaymath}
Since $\bigcup_{N \in \natls} \{ x \st \Xi_0(x) < N \}$ is conull in $X$, we
can choose $K(\delta) \subset X$ and $C_1 = C_1(\delta)$ so that
$\Xi_0(x) < C_1(\delta)$ for $x \in K(\delta)$ and $\nu(K(\delta)) \ge
(1-\delta)$. 
\qed\medskip

\starredsubsection{Proof of Lemma~\ref{lemma:subspaces:stay:close}.}
\label{sec:starredsubsec:proof:of:lemma:subspaces:stay:close}
We first prove (a). 
Note that the action of $g_t$ commutes with $P^{GM}$, i.e.\ 
\begin{displaymath}
P^{GM}(g_t x, g_t y) \circ g_t = g_t \circ P^{GM}(x,y). 
\end{displaymath}
Let $\alpha_0 = \min_{i \ne j} |\lambda_i
-\lambda_j|$, where the $\lambda_i = \lambda_i(H^1)$. 
We will choose $0 < \epsilon < \alpha_0/100$. For every $\epsilon > 0$ there
exists a compact set $K_0 = K_0(\epsilon) 
\subset X_0$ with $\nu(K_0) > 1 - \epsilon/4$ and $\sigma = \sigma(\epsilon)
> 0$ such that for any subset $S$ of the Lyapunov exponents, 
\begin{equation}
\label{eq:transverse:Lyapunov}
d_Y(\bigoplus_{i \in S} \cV_i(x),\bigoplus_{j \not\in S} \cV_j(x)) > \sigma \qquad \text{for all $x \in \pi^{-1}(K_0)$.}
\end{equation}
By the multiplicative ergodic theorem and the Birkhoff ergodic
theorem, there exists a set $K = K(\epsilon) \subset K_0$ with $\nu(K) >
1-\epsilon/2$ and a constant $C = C(\epsilon)$ such that
such that for all $z \in \pi^{-1}(K)$, all $s \in \reals$ and all $v \in
\cV_i(z)$, 
\begin{equation}
\label{eq:stay:osceledets}
C(\epsilon)^{-1/2} \|v\|_Y e^{\lambda_i s - (\epsilon/6) |s|} \le \|g_s
v \|_Y \le C(\epsilon)^{1/2} \|v\|_Y e^{\lambda_i s + (\epsilon/6) |s|},
\end{equation}
and also for any interval $I \subset \reals$ containing the origin of
length at least $4 \log C(\epsilon)/\alpha_0$, and any $z \in \pi^{-1}(K)$, 
\begin{equation}
\label{eq:stay:birkhoff}
|\{ s \in I \st g_s z \in K_0 \}| \ge (1-\epsilon)|I|.
\end{equation}
Suppose the set $\{ g_s x \st 0 \le s \le t\}$ intersects $K$. We will
show that for all $y \in \tilde{X}_0$ such that $d^{X_0}(g_s x, g_s y) \le 1/100$ for $0 \le s \le t$,
\begin{equation}
\label{eq:goal:stay:together}
d_Y(\cV_{\le i}(g_t x), P^{GM}(g_t y,g_t x)\cV_{\le i}(g_t y)) \le
C_0(x) C(\sigma) C(\epsilon)^2 (1+D^+(x,y)) e^{-\alpha t},
\end{equation}
where $C_0(x)$ is as in Lemma~\ref{lemma:one:sided:distance}.
Let $\Upsilon(x) = \{ \epsilon \st x \in K(\epsilon) \}$ and let
\begin{displaymath}
C(x) = C_0(x) \inf \{ C(\sigma) C(\epsilon)^2 \st \epsilon \in
\Upsilon(x) \}.  
\end{displaymath}
Since the union as $\epsilon \to 0$ of the sets $K = K(\epsilon)$ is conull,
(\ref{eq:goal:stay:together}) implies part (a) of the lemma. 

We now prove (\ref{eq:goal:stay:together}). 
We may assume that $t > 4 \log C(\epsilon)/\alpha_0$, otherwise (\ref{eq:goal:stay:together})
trivially holds. Then, by (\ref{eq:stay:birkhoff}), 
there exists $(1-\epsilon)t < t' \le t$ with $g_{t'} x \in K_0$. 
In view of Lemma~\ref{lemma:forni:upper} the inequality
(\ref{eq:goal:stay:together}) for $t'$ implies the inequality
(\ref{eq:goal:stay:together}) for $t$
(after replacing $\alpha$ by $\alpha - 4\epsilon$). Thus, we may
assume without loss of generality that $g_t x \in K_0$. 

By assumption, there exists $0 < s < t$ such that $g_s x \in
K$. Let $z = g_s x$. 
Then, applying (\ref{eq:stay:osceledets}) twice at $z$, we get, 
for all $v \in \cV_i(x)$, 
\begin{equation}
\label{eq:stay:newosceledets}
C(\epsilon)^{-1} \|v\|_Y e^{\lambda_i t - (\epsilon/3) t} \le \|g_t v \|_Y \le
C(\epsilon) \|v\|_Y e^{\lambda_i t + (\epsilon/3) t}.
\end{equation}
Let $v' \in P^{GM}(g_t y,g_t x)\cV_{\le i}(g_t y)$ be such that
$\|v'\|_Y = 1$ and 
\begin{displaymath}
d_Y(v', \cV_{\le i}(g_t x)) = \delta_Y(P^{GM}(g_t y,g_t x)\cV_{\le i}(g_t y), \cV_{\le i}(g_t x)),
\end{displaymath}
where $\delta_Y(\cdot, \cdot)$ is as in
Lemma~\ref{lemma:one:sided:distance}.
Then, $v' = g_t v$ for some $v \in P^{GM}(y,x) \cV_{\le i}(y)$. 
We may write
\begin{displaymath}
v = v_0 + w, \qquad\text{ $v_0 \in \cV_{\le i}(x)$, $w \in
  \cV_{>i}(x)$. }
\end{displaymath}
We have, by the definition of $D^+(\cdot, \cdot)$, 
\begin{displaymath}
\|w\|_Y \le D^+(x,y) \|v_0\|_Y.
\end{displaymath}
Then, we have 
\begin{displaymath}
v' = g_t v = g_t v_0 + g_t w,
\end{displaymath}
and by (\ref{eq:stay:newosceledets}), 
\begin{displaymath}
\|g_t v_0 \|_Y \ge C(\epsilon)^{-1} e^{(\lambda_i - \epsilon/3)t} \|v_0\|_Y,
\end{displaymath}
and
\begin{displaymath}
\|g_t w \|_Y \le C(\epsilon) e^{(\lambda_{i+1} + \epsilon/3)t} \|w\|_Y.
\end{displaymath}
Thus, 
\begin{displaymath}
\|g_t w \|_Y \le C(\epsilon)^2 D^+(x,y) e^{-(\alpha_0 - 2 \epsilon/3) t} \|g_t v_0\|_Y.
\end{displaymath}
Since $g_t v_0 \in \cV_{\le i}(g_t x)$ and $g_t w \in \cV_{>i}(g_t
x)$, this, together with (\ref{eq:transverse:Lyapunov}) implies
\begin{displaymath}
d_Y(v', \cV_{\le i}(g_t x)) \le C(\sigma) C(\epsilon)^2 (1+D^+(x,y)) e^{-(\alpha_0 - 2 \epsilon/3) t}. 
\end{displaymath}
This, together with Lemma~\ref{lemma:one:sided:distance}, completes
the proof of (\ref{eq:goal:stay:together}).  

The proof of (b) is identical. 
\qed\medskip

\starredsubsection{Proof of
  Proposition~\ref{prop:sublyapunov:new:locally:constant} and 
Proposition~\ref{prop:sublyapunov:locally:constant}.}
\label{sec:starredsubsec:proof:of:prop:sublyapunov:locally:constant}

The proof of Proposition~\ref{prop:sublyapunov:new:locally:constant}
will essentially be by reference to 
\cite[Theorem 1]{Ledrappier:Positivity}. We recall the setup (in our
notation):

Let $(X,\nu)$ be a measure space, and let $T: X \to X$ be a measure
preserving transformation. Let $\gB$ be a $\sigma$-subalgebra of the
$\sigma$-algebra of Borel sets on $X$, such that $\gB$ is
$T$-decreasing (i.e.\ $T^{-1} \gB \subset \gB$).
Let $\gB_{-\infty}$ denote the $\sigma$-algebra generated by all the
$\sigma$-algebras $T^n \gB$, $n \in \zed$. 

Let $V$ be a vector space, and let $A: X \to GL(V)$ be a log-integrable
$\gB$-measurable function. Let
\begin{displaymath}
A^{(n)}(x) = A(T^{n-1} x) \dots A(x) \qquad \text{ for $n > 0$ }
\end{displaymath}
\begin{displaymath}
A^{(0)}(x) = Id
\end{displaymath}
and
\begin{displaymath}
A^{(n)}(x) = A^{-1}(T^{n} x) \dots A^{-1}(T^{-1} x) \qquad \text{ for $n < 0$ }
\end{displaymath}
We have a skew-product map $\hat{T}: X \cross V \to X \cross V$ given by
\begin{displaymath}
\hat{T}(x,v) = (Tx, A(x) v),
\end{displaymath}
and then,
\begin{displaymath}
\hat{T}^n(x,v) = (T^n x, A^{(n)}(x) v). 
\end{displaymath}
Let
\begin{displaymath}
\gamma_+ = \lim_{n \to \infty} \frac{1}{n} \int_X \log \|A^{(n)}(x) \|
\, d\nu(x),
\end{displaymath}
\begin{displaymath}
\gamma_- = \lim_{n \to \infty} -\frac{1}{n} \int_X \log \|(A^{(n)}(x))^{-1} \|
\, d\nu(x).
\end{displaymath}
where $\| \cdot \|$ is the operator norm. The limits exist by the
subadditive ergodic theorem. 

The matrix $A$ also naturally acts on the projective space
$\proj(V)$. We use the notation $\hat{T}$ to denote also the
associated skew-product map $X \cross \proj(V) \to X \cross \proj(V)$.  

We have the following:
\begin{theorem}[Ledrappier, \protect{\cite[Theorem 1]{Ledrappier:Positivity}}]
\label{theorem:ledrappier}
Suppose 
\begin{itemize}
\item[{\rm (a)}] $\gamma_+ = \gamma_-$.
\item[{\rm (b)}] $x \to \nu_x$ is a family of measures on
  $\proj(V)$ defined for almost every $x$
  such that $A(x) \nu_x = \nu_{Tx}$ and such that the map $x \to
  \nu_x$ is $\gB_{-\infty}$-measurable.
\end{itemize}
Then, $x \to \nu_x$ is $\gB$-measurable. 
\end{theorem}

\bold{Proof of Proposition~\ref{prop:sublyapunov:new:locally:constant}.}
We first make some preliminary reductions. 
For $x \in X_0$, write 
$\cM(x) = \{ M^1(x), \dots, M^k(x) \}$. 
Since $\cM(x)$ is
$g_t$-equivariant, for $1 \le j \le k$,
\begin{displaymath}
M^j(x) = \bigoplus M_i^j(x), \qquad M_i^j(x) \subset L_i(x). 
\end{displaymath}
Let $\cM_i(x) = \{ M_i^1(x), \dots, M_i^k(x) \}$. 
Thus, it is enough to show that
\begin{displaymath}
F(x,y) \cM_i(x) = \cM_i(y). 
\end{displaymath}
Without loss of generality,
we may assume that for a fixed $i$, all the $M_i^j$ have the same dimension.
Suppose $x \in J_c$, where $J_c$ is as in
Proposition~\ref{prop:semi:markov}. Then the sets $\{ g_{-t} c \st 0
\le t \le t(c) \}$ and $\gB_0[x] = J_c \cap W^+[x]$ intersect at a unique point
$x_0 \in X_0$. Then,  we can replace the bundle $L(x)$ by $\tilde{L}(x)
\equiv F(x,x_0) L(x)$. Then, for $y \in \gB_0[x]$, 
\begin{displaymath}
\tilde{L}(y) = F(y,x_0) L(y) = F(y,x_0) F(x,y) L(x) = F(x,x_0) L(x) =
\tilde{L}(x), 
\end{displaymath}
i.e.\ $\tilde{L}(x)$ is locally constant along $W^+(x)$.
Also, by
(\ref{eq:def:connection:b}), the action of $(g_t)_*$ on $\tilde{L}$ is
locally constant. Thus, without loss of generality, we may assume that
$F$ is locally constant (or else we replace $L$ by $\tilde{L}$). 
Thus, it is enough to show that
assuming the subspaces $L_i(x)$ are almost everywhere
  locally constant along $W^+$,
the set of subspaces $\cM_i(x)$ is also almost everywhere 
locally constant along $W^+$. In other words, we assume that the functions
$x \to L_i(x)$ are $\gB_0$-measurable, and would like to show that the 
functions $x \to \cM_i(x)$ are $\gB_0$-measurable. 

Let $T=g_1$ denote the time $1$ map of the geodesic flow. Fix $i$ and
$j$, and let $d_i = \dim M_i^1 = \dots = \dim M_i^k$. 
Let $V(x) = \bigwedge^{d_i} (L_i(x)/L_{i-1}(x))$. Note that $V(x)$ is
$\gB_0$-measurable and $g_t$-equivariant. 

We can write the action of $(g_t)_*$ (for $t=1$) on the bundle $V$ as
\begin{displaymath}
(g_1)_* (x,v) = (g_1 x, A(x) v). 
\end{displaymath}
Then, $A(x)$ is $\gB_1$-measurable (where $\gB_t$ is as in
\S\ref{sec:semi:markov}). Also, the condition $\gamma_+ = \gamma_-$
follows from the multiplicative ergodic theorem. (In fact, 
$\gamma_+ = \gamma_- = d_i \lambda_i$, where $\lambda_i$ is the Lyapunov
exponent corresponding to $L_i$). 

Let $\nu^j_x$ denote the Dirac measure on (the line through ) 
$v_1 \wedge \dots \wedge v_d$,
where $\{ v_1, \dots, v_d \}$ is any  basis for $M_i^j(x)$, and let
\begin{displaymath}
\nu_x = \frac{1}{k} \sum_{j =1}^k \nu^j_x.
\end{displaymath}
Then, since the $\cM_i(x)$ are $g_t$-equivariant, the measures $\nu_x$ are
$\hat{T}$-invariant. Also note that $\gB_{-\infty}$ is the partition
into points. Thus, we can apply
Theorem~\ref{theorem:ledrappier} (with $\gB = \gB_1$). We conclude
that the function $x \to \nu_x$ is $\gB_1$-measurable, which implies
that the $\cM_i(x)$ are locally constant on atoms of $\gB_1$. Since
the $\cM_i(x)$ are $g_t$-equivariant, this implies that the
$\cM_i(x)$ are also locally constant (in particular the function $x
\to \cM_i(x)$ is $\gB_0$-measurable).
\qed\medskip

\bold{Proof of Proposition~\ref{prop:sublyapunov:locally:constant}.}
Note that (a) and also (\ref{eq:Lijy:Lijx}) follow immediately from
Proposition~\ref{prop:sublyapunov:new:locally:constant}.

We now prove (\ref{eq:locally:constant:quadratic:form}). After making
the same reductions as in the proof of
Proposition~\ref{prop:sublyapunov:new:locally:constant}, we may assume
that the $L_{ij}$ and $F$ are locally constant. 
Let $K \subset K_\epsilon$ denote a compact subset with $\nu(K) > 0.9$ 
where $\langle \cdot,  \cdot
\rangle_{ij}$ is uniformly continuous. Consider the points $g_t x$ and
$g_t y$, as $t \to -\infty$. Then $d^{X_0}(g_t x, g_t y) \to 0$. 
Let 
\begin{displaymath}
v_t = e^{-\lambda_{ij}(x,t)} (g_t)_* v, \quad 
w_t = e^{-\lambda_{ij}(x,t)} (g_t)_* w,  
\end{displaymath}
where $\lambda_{ij}(x,t)$ is as in
Lemma~\ref{lemma:jordan:canonical:form}. 
Then, by Lemma~\ref{lemma:jordan:canonical:form}, 
we have
\begin{equation}
\label{eq:rescaled:un}
\langle v_t, w_t \rangle_{ij,g_t x} = \langle v, w \rangle_{ij,x}, \quad 
\langle v_t, w_t \rangle_{ij,g_t y} = c(x,y,t) \langle v, w \rangle_{ij,y}.
\end{equation}
where $c(x,y,t) = e^{\lambda_{ij}(x,t)-\lambda_{ij}(y,t)}$. 

Now take a sequence $t_k \to \infty$ with $g_{t_k} x \in K$, 
$g_{t_k} y \in K$ (such a sequence exists for
$\nu$-a.e.\ $x$ and  $y$ with $y \in \gB_0[x]$). Then, since the
$L_{ij}(x)$ and the connection $F$ are assumed to be locally constant,
$c(x,y,t_k)$ is bounded between two constants. Also, 
\begin{displaymath}
\langle v_{t_k}, w_{t_k} \rangle_{ij,g_{t_k} x} - \langle
v_{t_k}, w_{t_k} 
\rangle_{ij,g_{t_k} y} \to 0.
\end{displaymath}
Now the equation (\ref{eq:locally:constant:quadratic:form})
 follows from (\ref{eq:rescaled:un}). 
\qed\medskip

\starredsubsection{Proof of
Proposition~\ref{prop:properties:dynamical:norm:Hbig}}
\label{sec:starredsubsec:proof:of:prop:properties:dynamical:norm:Hbig}
To simplify notation, we assume that $M_0=1$ (where
  $M_0$ is as in Lemma~\ref{lemma:exists:C:T0}).

\bold{The inner products $\langle \cdot, \cdot \rangle_{ij}$ on $E^+[c]$.}
Note that the inner products $\langle \cdot, \cdot\rangle_{ij}$ and
the $\reals$-valued cocycles $\lambda_{ij}$ of
Lemma~\ref{lemma:jordan:canonical:form} are not unique, since we can
always multiply $\langle \cdot, \cdot\rangle_{ij,x}$ by a scalar factor
$c(x)$, and then replace $\lambda_{ij}(x,t)$ by $\lambda_{ij}(x,t) +
\log c(g_t x) - \log c(x)$. 
In view of (\ref{eq:locally:constant:quadratic:form}) in
Proposition~\ref{prop:sublyapunov:locally:constant} (b), 
we may (and will) use this freedom to make
$\langle \cdot, \cdot \rangle_{ij,x}$ 
constant on each set $E^+[c]$, where $c \in \sigma_0^{-1}(\cC_1)$ and
$E^+[c]$ is as in \S\ref{sec:semi:markov} (see also
\S\ref{sec:subsec:the:cover:X}).

\bold{The inner product $\langle \cdot, \cdot \rangle_x$ 
on $E^+[c]$.}
Let 
\begin{equation}
\label{eq:Lyapunov:flag:Hbig:plusplus}
\{0\} = \cV_{\le 0} \subset \cV_{\le 1}  \subset \dots
\end{equation}
be the Lyapunov flag for $H_{big}^{(++)}$, and for each $i$, let
\begin{equation}
\label{eq:Hbig:maximal:invariant:refinement}
\cV_{\le i-1} = \cV_{\le i,0} \subset \cV_{i,1} \subset \dots \cV_{\le
  i,n_i} = \cV_{\le i}
\end{equation}
be a maximal invariant refinement.  

Let $L_i = \cV_i(H_{big}^{(++)})$ denote the Lyapunov subspaces for
$H_{big}^{(++)}$. Then we have a maximal invariant flag
\begin{displaymath}
\{0\} = L_{i,0} \subset  L_{i,1} \subset \dots \subset L_{i,n_i} = L_i,
\end{displaymath}
where $L_{ij} = L_i \cap \cV_{\le i,j}$. 

Let $c \in \sigma_0^{-1}(\cC_1)$, $E^+[c]$ be as in
\S\ref{sec:semi:markov} and \S\ref{sec:subsec:the:cover:X}. By
Lemma~\ref{lemma:exists:C:T0} (b), we can (and do) rescale the inner
products $\langle \cdot, \cdot \rangle_{ij,c}$ so that after the
rescaling, for all $v \in L_{ij}(c)/L_{i,j-1}(c)$,
\begin{displaymath}
(M')^{-1} \| v \|_{Y,c} \le \langle v, v \rangle_{ij,c}^{1/2} \le M'
\|v\|_{Y,c},
\end{displaymath}
where $\| \cdot \|_{Y,c}$ is the AGY norm at $\sigma_0(c)$ and $M' >
1$ is as in
Lemma~\ref{lemma:exists:C:T0}. We then choose $L'_{ij}(c) \subset
L_{ij}(c)$ to be a complementary subspace to $L_{i,j-1}(c)$ in
$L_{ij}(c)$, so that for all $v \in L_{i,j-1}(c)$ and all $v' \in
L'_{ij}(c)$, 
\begin{displaymath}
\| v + v' \|_{Y,c} \ge \rho'' \max( \|v\|_{Y,c}, \|v'\|_{Y,c} ),
\end{displaymath}
and $\rho'' > 0$ depends only on the dimension. 

Then, 
\begin{displaymath}
L'_{ij}(c) \isom L_{ij}(c)/L_{i,j-1}(c) \isom \cV_{\le i,j}(c)/\cV_{\le i,j-1}(c).
\end{displaymath}
Let $\pi_{ij}: \cV_{\le i,j} \to \cV_{\le i,j}/\cV_{\le i,j-1}$ be the natural quotient
map. Then the restriction of $\pi_{ij}$ to $L_{ij}'(c)$ is an
isomorphism onto $\cV_{\le i,j}(c)/\cV_{\le i,j-1}(c)$. 

We can now define for $u,v \in H_{big}^{(++)}(c)$, 
\begin{multline*}
\langle u, v \rangle_c \equiv \sum_{ij} \langle \pi_{ij}(u_{ij}),
\pi_{ij}(v_{ij}) 
\rangle_{ij,c},  \\
\text{ where $u = \sum_{ij} u_{ij}$, $v = \sum_{ij} v_{ij}$,
  $u_{ij} \in L'_{ij}(c)$, $v_{ij} \in L'_{ij}(c)$. }
\end{multline*}
In other words, the distinct $L'_{ij}(c)$ are orthogonal, and the
inner product on each $L'_{ij}(c)$ coincides with $\langle \cdot, \cdot
\rangle_{ij,c}$ under the identification $\pi_{ij}$ of $L'_{ij}(c)$ with
$\cV_{\le i,j}(c)/\cV_{\le i,j-1}(c)$. 

We now define, for $x \in E^+[c]$, and $u,v \in H_{big}^{(++)}(x)$
\begin{displaymath}
\langle u , v \rangle_x \equiv \langle P^+(x,c) u, P^+(x,c)
v\rangle_c,
\end{displaymath}
where $P^+(\cdot, \cdot)$ is the connection defined in
\S\ref{sec:subsec:connection}. 
Then for $x \in E^+[c]$, 
the inner product $\langle \cdot, \cdot
  \rangle_x$ induces the inner product $\langle  \cdot, \cdot
  \rangle_{ij,x}$ on $\cV_{\le i,j}(x)/\cV_{\le i,j-1}(x)$. 
\mcc{(check this induces the correct scaling)}

\bold{Symmetric space interpretation.}
We want to define the inner product $\langle
\cdot, \cdot \rangle_x$ for any $x
\in J[c]$ by interpolating between $\langle \cdot, \cdot \rangle_c$
and $\langle \cdot, \cdot\rangle_{c'}$, where $c'$ is such that
$g_{-t(c)} c \in E^+[c']$. To define this interpolation, we recall
that the set of inner products on a vector space $V$ is canonically
isomorphic to  $SO(V) \bs GL(V)$, where $GL(V)$ is the general linear
group of $V$ and $SO(V)$ is the subgroup preserving the inner product
on $V$. In our case, $V = H_{big}^{(++)}(c)$ with the inner product
$\langle \cdot, \cdot \rangle_c$. 

Let $K_c$ denote the subgroup of $GL(H_{big}^{(++)}(c))$ which
preserves the inner product $\langle \cdot, \cdot \rangle_c$. 
Let $\cQ$ denote the parabolic subgroup of $GL(H_{big}^{(++)}(c))$ which
preserves the flags  (\ref{eq:Lyapunov:flag:Hbig:plusplus}) and
(\ref{eq:Hbig:maximal:invariant:refinement}),  and on each successive
quotient $\cV_{\le i,j}(c)/\cV_{\le i,j-1}(c)$ preserves $\langle \cdot, \cdot
\rangle_{ij,c}$. Let $K_c A'$ denote
the point in $K_c \bs GL(H_{big}^{(++)}(c))$ which represents the inner
product $\langle \cdot, \cdot \rangle_{c'}$, i.e.\
\begin{displaymath}
\langle u,v \rangle_{c'} = \langle A' u, A' v\rangle_{c}.
\end{displaymath}
Then, since 
$\langle \cdot, \cdot \rangle_{c'}$ induces the inner products
$\langle \cdot, \cdot \rangle_{ij,c'}$ on the space $\cV_{\le
  i,j}(c')/\cV_{\le i,j-1}(c')$
which is the same as $\cV_{\le i,j}(g_{-t(c) c})/\cV_{\le
  i,j-1}(g_{-t(c) c})$), we may assume that the matrix
product $A' (g_{-t(c)})_*$ is in  $\cQ$. 

Let $N_{\cQ}$ be the normal
subgroup of $\cQ$ in which all diagonal blocks are the identity, and let
$\cQ' = \cQ/N_\cQ$.  (We may consider $\cQ'$ to be the subgroup of $\cQ$ in
which all off-diagonal blocks are $0$). 
Let $\pi'$ denote the natural map $\cQ \to \cQ'$. 

\begin{claim}
\label{claim:estimate:A:prime}
We may write
\begin{displaymath}
 A' (g_{-t(c)})_*  = \Lambda A'',
\end{displaymath}
where $\Lambda \in \cQ'$ is the diagonal matrix which is scaling by
$e^{-\lambda_i t(c)}$ on $L_i(c)$,  $A'' \in \cQ$ and $\| A'' \| = 
O(e^{\epsilon t(c)})$.
\mcc{say which norm and what is epsilon.}
\end{claim}

\bold{Proof of claim.} 
Suppose $x \in E^+[c]$ and $t = -t(c) < 0$ where $c \in \cC_1$ and $t(c)$ is as in
Proposition~\ref{prop:semi:markov}. 
By construction, $t(c) > T_0(c)$, where $T_0(c)$ is as in
Lemma~\ref{lemma:exists:C:T0}. Then, the claim follows from
(\ref{eq:bound:norm:Pplus}) and Lemma~\ref{lemma:exists:C:T0} (d). 
\qed\medskip

\bold{Interpolation.} We may write $A'' = D A_1$, where $D$ is
diagonal, and $\det A_1 = 1$. In view of Claim~\ref{claim:estimate:A:prime}, 
$\|D \| = O(e^{\epsilon t})$ and $\|A_1\| = O(e^{\epsilon
  t})$.

We now connect $K_c \bs A_1$ to the identity by the shortest possible path
$\Gamma: [-t(c),0] \to K_c \bs K_c \cQ$, which stays 
in the subset $K_c \bs K_c \cQ$ of the symmetric space $K_c \bs SL(V)$. (We
parametrize the path so it has constant speed). 
This path
has length $O(\epsilon t)$ where the implied constant depends only on the
symmetric space. \mcc{this needs more explanation}

Now for $-t(c) \le t \le 0$, let 
\begin{equation}
\label{eq:def:A(t)}
A(t) = (\Lambda D)^{-t/t(c)} \Gamma(t). 
\end{equation}
Then $A(0)$ is the identity map, and $A(-t(c)) = A' (g_{-t(c)})_*$. Then, we define,
for $x \in E^+[c]$ and $-t(c) \le t \le 0$, 
\begin{displaymath}
\langle (g_t)_* u, (g_t)_* v \rangle_{g_t x} = \langle A(t) u, A(t) v
\rangle_{x}.  
\end{displaymath}

\bold{Proof of Proposition~\ref{prop:properties:dynamical:norm:Hbig}.}
Suppose first that $x = c$, where $c$ and $E^+[c]$ are
as in \S\ref{sec:semi:markov} and \S\ref{sec:subsec:the:cover:X}. 
Then, by construction, (a) and (b) hold. Also,
from the construction, it is clear that the inner 
product  $\langle \cdot, \cdot \rangle_c$ induces the inner product
$\langle \cdot, \cdot \rangle_{ij,c}$ on
$L_{ij}(c)/L_{i,j-1}(c)$. 

Now by Proposition~\ref{prop:sublyapunov:locally:constant}, for $x \in
E^+[c]$, $P^+(x,c) L_{ij}(x) = L_{ij}(c)$, and for $\bar{u}, \bar{v}
\in L_{ij}(x)/L_{i,j-1}(x)$, 
$\langle u , v \rangle_{ij,x} = \langle P^+(x,c) u, P^+(x,c) v
\rangle_{ij,c}$.  Therefore, (a), (b), (e) and (f) hold for $x \in
E^+[c]$, and also 
for $x \in E^+[c]$, the inner product $\langle \cdot, \cdot \rangle_x$
induces \mcc{(a multiple of ?)} the inner product 
$\langle \cdot, \cdot \rangle_{ij,x}$ on $L_{ij}(x)/L_{i,j-1}(x)$. 
Now, (a),(b),(e) and (f) hold for arbitrary $x \in J[c]$ since $A(t) \in
\cQ$. 

Let $\psi_{ij}: \cQ' \to \reals_+$ denote the homomorphism taking the
block-conformal matrix $\cQ'$ to the scaling part of block corresponding
to $L_{ij}/L_{i,j-1}$.
Let $\varphi_{ij} = \psi_{ij} \circ \pi'$; then $\varphi_{ij}: \cQ
\to \reals_+$ is a homomorphism. 

From (\ref{eq:def:A(t)}), we have, for $x \in E^+[c]$ and $-t(c)
\le t \le 0$, 
\begin{displaymath}
\lambda_{ij}(x,t) = \log \varphi_{ij}(A(t)) = t \lambda_i + \gamma_{ij}(x,t),
\end{displaymath}
where $t \lambda_i$ is the contribution of $\Lambda^{t/t(c)}$ and
$\gamma_{ij}(x,t)$ is the contribution of $D^{t/t(c)} \Gamma(t)$. 
By Claim~\ref{claim:estimate:A:prime}, for all $-t(c) \le t \le 0$,
\begin{equation}
\label{eq:ttemp:gamma:prime:epsilon}
|\frac{\partial}{\partial t} \gamma_{ij}(x,t)| = O( \epsilon) 
\end{equation}
where $\epsilon > 0$ is as in Claim~\ref{claim:estimate:A:prime}, and
the implied constant depends only on the symmetric
space. Without loss of generality, the function
  $T_0(x)$ in Lemma~\ref{lemma:exists:C:T0} can be chosen large enough
  so that since $t(c) > T_0(c)$, (c) holds.
\mcc{fix epsilon vs. multiple of epsilon}

The lower bound in (d) now follows immediately from (b) and (c). The
upper bound in (d) follows from (\ref{eq:ttemp:gamma:prime:epsilon}). 
\qed\medskip

\section{Conditional measure lemmas}
\label{sec:conditional}

In \S\ref{sec:conditional}-\S\ref{sec:first:divergence} we work on
$X_0$ (and not on $X$).

\bold{Motivation.} We use notation from \S\ref{sec:outline:step1}.
Recall that $\cL^-(q)$ is the smallest linear
subspace of $W^-(q)$ containing the support of the conditional measure
$\nu_{{W^-}(q)}$. 
For two (generalized) subspaces $\cU'$ and $\cU''$ and $x \in \tilde{X}_0$ let
\index{$hd_x^{X_0}(\cdot,\cdot)$}$hd_x^{X_0}(\cU',\cU'')$ denote the Hausdorff distance between $\cU' \cap
B^{X_0}(x,1/100)$ and $\cU'' \cap B^{X_0}(x,1/100)$, 
where \index{$B^{X_0}(x,r)$}$B^{X_0}(x,r)$ denotes $\{ y \in \tilde{X}_0 \st
d^{X_0}(x,y) < r\}$. For $x \in X_0$, we will sometimes
write
\index{$hd_x^{X_0}(\cdot,\cdot)$}$hd_x^{X_0}(\cU', \cU'')$ instead of
$hd_{\tilde{x}}^{X_0}(\cU', \cU'')$ as long as the proper lift
  $\tilde{x} \in \tilde{X}_0$ of $x$ is clear from the context.

We can write
\begin{displaymath}
hd_{q_2}^{X_0}(U^+[q_2'], U^+[q_2]) = Q_t(q' - q), 
\end{displaymath}
where $Q_t: \cL^-(q) \to \reals$ is a map depending on $q$, $u$,
$\ell$, and $t$. 
The map $Q_t$ is essentially the composition of
flowing forward for time $\ell$, shifting by $u \in U^+$ and then flowing
forward again for time $t$. We then adjust $t$ so that
$hd_{q_2}^{X_0}(U^+[q_2'], U^+[q_2])\approx \epsilon$, 
where $\epsilon > 0$ is a priori fixed. 

In order to solve ``technical difficulty \#1'' of
\S\ref{sec:outline:step1}, it is crucial to ensure that $t$ does
not depend on the precise choice of $q'$ (it can depend on $q$, $u$,
$\ell$). The idea is to use the following trivial:
\begin{lemma}
\label{lemma:bad:subspace}
For any $\rho > 0$ there is a constant $c(\rho)$ with the following
property:  
Let $A: \cV \to \cW$ be a linear map between Euclidean spaces.
Then there exists a proper subspace
$\cM \subset \cV$ such that for any $v$ with
$\|v\|=1$ and $d(v,\cM) > \rho$, we have 
\begin{displaymath}
\|A\| \ge \|A v\| \ge c (\rho) \|A\|. 
\end{displaymath}
\end{lemma}
\bold{Proof of Lemma~\ref{lemma:bad:subspace}.} 
The matrix $A^t A$ is symmetric, so it has a complete
orthogonal set of eigenspaces $W_1, \dots, W_m$ corresponding to
eigenvalues $\mu_1 > \mu_2 > \dots \mu_m$. Let $\cM = W_1^{\perp}$. 
\qed\medskip

Now suppose the map $Q_t: \cL^-(q) \to \reals$ is of the form
$Q_t(v) = \| \cQ_t(v) \|$ where 
$\cQ_t: \cL^-(q) \to \bH(q_2)$ is a linear map, and $\bH(q_2)$ a
vector space. This
in fact happens in the first step of the induction where $U^+$ is the
unipotent $N$ (and we can take $\bH(q_2) = W^+(q_2)/N$). 
We can then choose $t$, depending only on $q$, $u$
and $\ell$, such that the operator norm
\begin{displaymath}
\|\cQ_t\| \equiv \sup_{ v \in \cL^-(q)}
\frac{\|\cQ_t(v)\|}{\|v\|} = \epsilon. 
\end{displaymath}
Then, we need to prove
that we can choose $q' \in \cL^-[q]$ such that $\|q' - q\| \approx 1/100$,
$q'$ avoids an a priori given set of small measure, 
and also $q' - q$ is at least $\rho$ away from the ``bad subspace'' $\cM =
\cM_u(q,\ell)$ of Lemma~\ref{lemma:bad:subspace}. (Actually, since we
do not want the choice of $q'$ to depend on the choice of $u$, we want
to choose $q'$ such that $q'-q$ avoids most of the subspaces $\cM_u$
as $u \in U^+$ varies over a unit box). Then, for most $u$, 
\begin{displaymath}
c(\rho) \epsilon \le \|\cQ_t(q_2' - q_2) \| \le \epsilon,
\end{displaymath}
and thus
\begin{equation}
\label{eq:wishful:thinking:q2prime:q2}
c(\rho) \epsilon \le hd_{q_2}(U^+[q_2],U^+[q_2']) \le \epsilon,
\end{equation}
as desired. 
In general we do not know that the map $\cQ_t$ is linear, because we do
not know the dependence of the subspace $U^+(q)$ on $q$. To handle
this problem, we can write
\begin{displaymath}
\cQ_t(q'-q) = \cA_t(F(q') - F(q))
\end{displaymath}
where the map $\cA_t: \cL_{ext}[q]^{(r)} 
\to W^+(q_2)$ is linear (and can
depend on $q$, $u$, $\ell$), and the measurable map 
$F: \cL^-[q] \to \cL_{ext}[q]^{(r)}$ depends only on $q$. (See
Proposition~\ref{prop:reason:cA:F} below for a precise statement). 
The map $F$ and the space $\cL_{ext}[q]^{(r)}$ 
are defined in this section, and
the linear map $\cA_t = \cA(q,u,\ell,t)$ 
is defined in \S\ref{sec:subsec:themap:cA}. 

We then proceed in the same way. We choose $t =
\hat{\tau}(q,u,\ell,\epsilon)$ so that $\|\cA_t\| = \epsilon$. (A
crucial bilipshitz type property of the function $\hat{\tau}$ similar
to (\ref{eq:trivial:bilipshitz}) is
proved in \S\ref{sec:new:bilipshitz}). In this section we prove
Proposition~\ref{prop:can:avoid:most:Mu}, 
which roughly states that (for most
$q$) we can choose $q' \in \cL^-[q]$ while avoiding an a priori given
set of small measure, so that $\|F(q') - F(q)\| \approx 1/100$ and also
$F(q') - F(q)$ avoids most of a family of linear subspaces of
$\cL_{ext}[q]^{(r)}$ (which will be the ``bad subspaces'' of the linear maps
$\cA_t$ as $u$ varies over $U^+$). Then as above, for most $u$,
(\ref{eq:wishful:thinking:q2prime:q2}) holds. 
We can then proceed using (a variant of)
  Lemma~\ref{lemma:trivial:bilipshitz} as outlined in
  \S\ref{sec:outline:step1}. 

In view of the above discussion, we need to keep track of the way
$U^+[y]$ varies as $y$ varies over $W^-[x]$. In view of
Proposition~\ref{prop:sublyapunov:locally:constant}(a), all bundles
equivariant with respect to the geodesic flow are, when restricted to
$W^-$, equivariant with respect to the connection $P^-(x,y)$ defined
in \S\ref{sec:subsec:connection}. Thus, it
will be enough for us to keep track of the maps $P^-(x,y)$. However,
this is a bit awkward, since $P^-(x,y)$ depends on two points $x$ and
$y$. Thus, it is convenient to prove the following:
\begin{lemma}
\label{lemma:there:exists:gP}
There exists a subbundle $\index{$Y$@$\cY$}\cY \subset H_{big}^{(-)}$, 
locally constant under the Gauss-Manin connection along $W^-$, and for
almost all $x \in X_0$ an invertible linear map $\index{$P$@$\gP$}\gP(x): X_0 \to
Hom(\cY(x), H^1(M,\Sigma,\reals))$,  such that for almost all $x$, $y$, 
\begin{equation}
\label{eq:Pminus:gP}
P^-(x,y) = \gP(y) \circ \gP(x)^{-1}.
\end{equation}
\end{lemma}
The proof of Lemma~\ref{lemma:there:exists:gP} is simple, but
notationally heavy, and is relegated to
\S\ref{sec:starredsubsec:there:exists:gP}. It may be skipped on first
reading. 

\bold{The spaces $\cL^-(x)$ and $\cL_{ext}(x)$.}
Let the subspace \index{$L$@$\cL^-(x)$}$\cL^-(x) \subset W^-(x)$ be the smallest such that
the conditional measure $\nu_{W^-[x]}$ is supported on
$\cL^-[x]$. 
Since $\nu$ is invariant under $N$, the entropy of any $g_t \in A$ is
positive. Therefore for $\nu$-almost all $x \in X_0$, $\cL^-(x) \ne
\{0\}$ (see Proposition~\ref{prop:mt:93}). 

In the same spirit, let
\begin{displaymath}
\index{$L$@$\cL_{ext}(x)$} \cL_{ext}[x] \subset Hom(\cY(x), H^1(M,\Sigma,\reals))
\end{displaymath}
denote the smallest affine subspace which for almost every $y \in
W^-[x]$ contains the vector $\gP(y)$.  
(This makes sense since $\cY(x)$ is locally constant along
$W^-[x]$.) We also set $\cL_{ext}(x)$ to be the vector space spanned
by all vectors of the form $\gP(y) - \gP(x)$ as $y$ varies over
$W^+[x]$. Then, 
\begin{displaymath}
\cL_{ext}(x) = \cL_{ext}[x] -\gP(x).
\end{displaymath}
Note that for almost all $x$ and almost all $y \in W^-[x]$, 
$\cL_{ext}[y] = \cL_{ext}[x]$.

\bold{The space $\cL_{ext}(x)^{(r)}$ and the function $F$.}
For a vector space $V$ we use
the notation \index{$V^{\tensor m}$}$V^{\tensor m}$ to denote the $m$-fold tensor product of
$V$ with itself.  If $f: V \to W$ is a linear map, we write 
\index{$f^{\tensor m}$}$f^{\tensor m}$ for the induced linear map from $V^{\tensor m}$ to
$W^{\tensor m}$. Let $\index{$j^{\tensor m}$}j^{\tensor m}: V \to V^{\tensor m}$ denote the map $v \to
v \tensor \dots \tensor v$ ($m$-times). 

Let \index{$V^{\uplus m}$}$V^{\uplus m}$ denote $\bigoplus_{k=1}^m V^{\tensor k}$. 
If $f: V \to W$ is a linear map, we write 
\index{$f^{\uplus m}$}$f^{\uplus m}$ for the induced linear map from $V^{\uplus m}$ to
$W^{\uplus m}$ given by 
\begin{displaymath}
f^{\uplus m}(v) = (f^{\tensor 1}(v), f^{\tensor {2}}(v), \dots
f^{\tensor {m}}(v)). 
\end{displaymath}
Now if $V$ and $W$ are affine spaces, then we can still canonically define
$V^{\uplus m}$ and $W^{\uplus m}$, and an affine map $f: V \to W$
induces an affine map $f^{\uplus m}: V^{\uplus m} \to W^{\uplus m}$.

Let \index{$r$}$r$ be an integer to be chosen later. 
Let $\index{$F$}F: X_0 \to \cL_{ext}[x]^{\uplus r}$ denote the diagonal embedding
\begin{displaymath}
F(x) = \gP(x)^{\uplus r}.
\end{displaymath}
Let
\begin{displaymath}
\index{$L$@$\cL_{ext}[x]^{(r)}$}\cL_{ext}[x]^{(r)} \subset \cL_{ext}(x)^{\uplus r}
\end{displaymath}
denote the smallest affine subspace which 
contains the vectors $F(y)$ for almost all $y \in W^-[x]$. 
We also set
\begin{displaymath}
\index{$L$@$\cL_{ext}(x)^{(r)}$}\cL_{ext}(x)^{(r)} = \cL_{ext}[x]^{(r)} - F(x).
\end{displaymath}
Note that for $y \in W^-[x]$, $\cL_{ext}[y]^{(r)} =
\cL_{ext}[x]^{(r)}$.
\medskip

In this section, let $(\cB, | \cdot |)$ be a finite measure space. (We
will use the following proposition with $\cB \subset U^+$ is a ``unit
box''. The precise setup will be given in
\S\ref{sec:divergence:subspaces}).  
\medskip

To carry out the program outlined at the beginning of
\S\ref{sec:conditional}, we need the following:

\begin{proposition}
\label{prop:can:avoid:most:Mu}
For every $\delta > 0$ there exist constants 
$c_1(\delta) > 0$,  $\epsilon_1(\delta) > 0$
with $c_1(\delta) \to 0$ and $\epsilon_1(\delta) \to 0$ 
as $\delta \to 0$, and also constants
$\rho(\delta) > 0$, $\rho'(\delta) > 0$, and
$C(\delta) < \infty$ such that the following holds: 

For any subset $K' \subset X_0$ with $\nu(K') > 1-\delta$, 
there exists a subset $K \subset K'$ 
with $\nu(K) > 1-c_1(\delta)$
such that the following holds: suppose for each $x\in X_0$ we have a
measurable map from $\cB$ to proper subspaces of
$\cL_{ext}(x)^{(r)}$, written as 
$u \to \cM_u(x)$, where $\cM_u(x)$ is
a proper subspace of $\cL_{ext}(x)^{(r)}$. Then, for 
any $q \in K$ there exists $q' \in K'$ with
\begin{equation}
\label{eq:rho:prime:delta:le:d:q:qprime}
\rho'(\delta) \le d^{X_0}(q,q') \le 1/100
\end{equation}
and
\begin{equation}
\label{eq:rho:delta:le:Fq:minus:Fqprime}
\rho(\delta) \le \|F(q') - F(q) \|_Y \le C(\delta)
\end{equation}
and so that
\begin{equation}
\label{eq:qprime:avoids:Mu}
d_Y(F(q')-F(q), 
\cM_u(q)) > \rho(\delta) \qquad \text{\rm  for at least
  $(1-\epsilon_1(\delta))$-fraction of $u \in \cB$. } 
\end{equation}
\end{proposition}

This proposition is proved in 
\S\ref{sec:starredsubsec:proof:of:prop:can:avoid:most:Mu}.
The proof uses almost
nothing about the maps $F$ or the measure $\nu$, other than the
definition of $\cL_{ext}(x)$. It may be skipped on first reading. 

\mccc{Check that $K$ does not depend on the $\cM_u$. The
  referee seems to think that it does.}

\starredsubsection{Proof of Lemma~\ref{lemma:there:exists:gP}.}
\label{sec:starredsubsec:there:exists:gP}
As in \S\ref{sec:subsec:lyapunov:flags},  
let $\cV_i(x) \equiv \cV_i(H^1)(x) \subset H^1(M , \Sigma, \reals)$ 
denote the subspace corresponding to
the (cocycle) Lyapunov exponent $\lambda_i$.  Let
\begin{displaymath}
\cY(x) = \bigoplus_{i=1}^k \cV_{\ge i}(x)/\cV_{> i}(x),
\end{displaymath}
where $\cV_{\ge j}$ and $\cV_{> j}$ are as in \S\ref{sec:subsec:lyapunov:flags}.
Let $\pi_i: \cV_{\ge i}(x)
\to \cV_{\ge i}(x)/\cV_{> i}(x)$ denote the natural
projection. 

For $x \in X_0$, let $\index{$P_{i,x}$}P_{i,x} \in Hom(
\cV_{\ge i}(x)/\cV_{> i}(x), H^1(M, \Sigma, \reals))$ 
denote the unique linear 
map such that for $\bar{x} \in \cV_{\ge i}(x)/\cV_{> i}(x)$, 
$P_{i,x}(\bar{x}) \in \cV_i(H^1)(x)$ and 
$\pi_i(P_{i,x}(\bar{x})) = \bar{x}$. Note that the
$P_{i,x}$ satisfy the following:
\begin{equation}
\label{eq:Pix:equivariant}
P_{i,g_t x} = g_t \circ P_{i,x} \circ g_t^{-1}, 
\end{equation}
and
\begin{equation}
\label{eq:difference:P:lower:flag}
P_{i,x}(\bar{u}) - P_{i,y}(\bar{u}) \in \cV_{> i}(x). 
\end{equation}

\bold{Example.} The space $\cV_{\ge 1}/\cV_{>1}$  
is one dimensional, and corresponds to the Lyapunov
exponent $\lambda_1 = 1$. If we identify it with $\reals$ in the
natural way then $P_{1,x}: \reals \to H^1(M, \Sigma,
\reals)$ is given by the formula
\begin{equation}
\label{eq:P1x:example}
P_{1,x}(\xi) = (\Im x) \xi
\end{equation}
where for $x = (M,\omega)$, we write $\Im x$ for the imaginary part of
$\omega$.  
\medskip

Let
\begin{displaymath}
\index{$P$@$\gP$}\gP: X_0 \to \bigoplus_{i=1}^k Hom(\cV_{\ge i}(x)/\cV_{> i}(x),
H^1(M, \Sigma, \reals)) 
\end{displaymath}
be given by 
\begin{displaymath}
\gP(x) = (P_{1,x}, \dots P_{k,x}).
\end{displaymath}
Then, we can think of $\gP(x)$ as a map from $\cY(x)$ to $H^1(M,
\Sigma, \reals)$ and (\ref{eq:Pminus:gP}) holds, 
where $P^-(x,y)$ is as in \S\ref{sec:subsec:connection}. 
\qed\medskip

\starredsubsection{Proof of Proposition~\ref{prop:can:avoid:most:Mu}}
\label{sec:starredsubsec:proof:of:prop:can:avoid:most:Mu}

$ $

\bold{The measure $\tilde{\nu}_x$.}
Let $\tilde{\nu}_x = F_* \left(\nu_{W^-[x]}\right)$ denote the
pushforward of $\nu_{W^-}$ under $F$. 
Then $\tilde{\nu}_x$ is a measure supported on $\cL_{ext}[x]^{(r)}$. 
(Note that for $y \in W^-[x]$,
$\tilde{\nu}_x = \tilde{\nu}_y$).

\begin{lemma}
\label{lemma:not:supported:finite:union}
For $\nu$-almost all $x \in X_0$, for any $\epsilon >
0$ (which is allowed to depend on $x$), the restriction of the 
measure $\tilde{\nu}_x$ to the ball $B(F(x),\epsilon) \subset
\cL_{ext}[x]^{(r)}$  is not supported on a finite union of proper affine 
subspaces of $\cL_{ext}[x]^{(r)}$. 
\end{lemma}

\bold{Outline of proof.} Suppose not. Let $N(x)$ be the minimal
integer $N$ such that for some $\epsilon = \epsilon(x) > 0$, 
the restriction of $\tilde{\nu}_x$ to $B(F(x),\epsilon)$ is
supported on $N$ affine subspaces.  
Note that in view of (\ref{eq:Pix:equivariant}) and
(\ref{eq:difference:P:lower:flag}), the induced action on $\cL_{ext}$
(and hence on $\cL_{ext}^{(r)}$)  
of $g_{-t}$ for $t \ge 0$ is expanding. 
Then $N(x)$ is invariant under $g_{-t}$, $t \ge 0$. This
implies that $N(x)$ is constant for $\nu$-almost all $x$, and also
that the only affine subspaces of $\cL_{ext}[x]^{(r)}$ which contribute to
$N(\cdot)$ pass through $F(x)$. Then, $N(x) >1$ almost everywhere is
impossible. Indeed, suppose $N(x) = k$ a.e., then pick $y$ near $x$ 
such that $F(y)$ is in one of the affine subspaces through $F(x)$; then
there must be exactly $k$ affine subspaces of non-zero measure passing
though $F(y)$, but then at most one of them passes through $F(x)$. Thus,
the measure restricted to a neighborhood of $F(x)$ gives positive
weight to at least $k+1$ subspaces, contradicting our
assumption. Thus, we must have $N(x) = 1$ almost everywhere; but then
(after flowing by $g_{-t}$ for sufficiently large $t > 0$) we see that
for almost all $x$, $\tilde{\nu}_x$ is supported on a proper subspace of
$\cL_{ext}[x]^{(r)}$ passing through $x$, which 
contradicts the definition of $\cL_{ext}(x)^{(r)}$. 
\qed\medskip

\bold{Remark.} Besides Lemma~\ref{lemma:not:supported:finite:union},
the rest of the proof of
Proposition~\ref{prop:can:avoid:most:Mu} 
uses only the measurability of the map $F$.

\bold{The measure $\hat{\nu}_x$.} Let \index{$B$@$\gB_0^-$}$\gB_0^-$ 
be the analogue of the partition $\gB_0$ constructed in \S\ref{sec:semi:markov}
but along the stable leaves $W^-$. (The only
properties we use here is that $\gB_0^-$ is a measurable partition
subordinate to $W^-$ with atoms of diameter at most $1/100$). Let
$\gB_0^-[x] \subset W^-[x]$ 
denote the atom of the partition $\gB_0^-$ containing $x$. 

Let $\hat{\nu}_x = F_*(\nu_{W^-[x]}|_{\gB_0^-[x]})$, 
i.e.\ $\hat{\nu}_x$ is 
the pushforward
under $F$ of the restriction of $\nu_{W^-[x]}$ to $\gB_0^-[x]$. 
Then, for $y \in \gB_0^-[x]$, $\hat{\nu}_x =
\hat{\nu}_y$. Suppose $\delta > 0$ is given. Since
\begin{displaymath}
\lim_{C \to \infty} \hat{\nu}_x(B(F(x),C)) = \hat{\nu}_x(\cL_{ext}[x]^{(r)}),
\end{displaymath}
there exists a function $c(x)> 0$ finite almost everywhere such that for
almost all $x$, 
\begin{displaymath}
\hat{\nu}_x(B(F(x),c(x))) > (1-\delta^{1/2}) \hat{\nu}_x(\cL_{ext}[x]^{(r)}).
\end{displaymath}
Therefore, we can find $C = C(\delta)> 0$ and a compact set $K_\delta'$
with $\nu(K_\delta') > 1-\delta^{1/2}$ such that for each $x \in
K_\delta'$, 
\begin{equation}
\label{eq:def:Kprime:delta}
\hat{\nu}_x(B(F(x),C)) > (1-\delta^{1/2}) \hat{\nu}_x(\cL_{ext}[x]^{(r)})
\quad\text{ for all $x \in K_\delta'$.}
\end{equation}
In the rest of
\S\ref{sec:starredsubsec:proof:of:prop:can:avoid:most:Mu},
$C$ will refer to the constant of (\ref{eq:def:Kprime:delta}).

\begin{lemma}
\label{lemma:cond:proper:subspace:not:full:measure}
For every $\eta > 0$ and every $N > 0$ there exists $\beta_1 =
\beta_1(\eta,N)  > 0$, $\rho_1 = \rho_1(\eta,N) > 0$ 
and a compact subset 
$K_{\eta,N}$ of measure at least $1-\eta$ such that
for all $x \in K_{\eta,N}$, and any proper subspaces $\cM_1(x), \dots,
\cM_N(x)\subset \cL_{ext}(x)^{(r)}$,
\begin{equation}
\label{eq:cond:proper:subspace:not:full:measure}
\hat{\nu}_x( B(F(x),C) \setminus \bigcup_{k=1}^N\Nbhd(\cM_k(x), \rho_1) )
\ge \beta_1 \hat{\nu}_x(B(F(x),C)). 
\end{equation}
\end{lemma}

\bold{Outline of Proof.} By
Lemma~\ref{lemma:not:supported:finite:union}, 
there exist $\beta_x = \beta_x(N) > 0$ and $\rho_x = \rho_x(N) > 0$ such
that for any subspaces $\cM_1(x), \dots \cM_N(x) \subset \cL_{ext}(x)^{(r)}$, 
\begin{equation}
\label{eq:temp}
\hat{\nu}_x( B(F(x),C) \setminus \bigcup_{k=1}^N \Nbhd(\cM(x), \rho_x) )
\ge \beta_x \hat{\nu}_x(B(F(x),C)). 
\end{equation}
Let $E(\rho_1,\beta_1)$ be the set of $x$ such that
(\ref{eq:cond:proper:subspace:not:full:measure}) holds. By (\ref{eq:temp}),
\begin{displaymath}
\nu\left(\bigcup_{\stackrel{\rho_1>0}{\beta_1>0}} E(\rho_1,\beta_1)
\right) = 1. 
\end{displaymath}
Therefore, we can choose $\rho_1 > 0$ and $\beta_1 > 0$ such that
$\nu(E(\rho_1, \beta_1)) > 1-\eta$. 
\qed\medskip

\begin{lemma}
\label{lemma:avoid:most:Mu}
For every $\eta > 0$ and every $\epsilon_1 > 0$ there exists
$\beta = \beta(\eta,\epsilon_1) > 0$, a compact set $K_\eta =
K_\eta(\epsilon_1)$ of measure at least
$1-\eta$, and $\rho = \rho(\eta,\epsilon_1) > 0$ 
such that the following holds: 
Suppose for each $u \in \cB$ let $\cM_u(x)$ be a proper subspace of
$\cL_{ext}(x)^{(r)}$. Let 
\begin{multline*}
E_{good}(x) = \{ v \in B(F(x),C) \st \text{ for at least
  $(1-\epsilon_1)$-fraction of $u$ in $\cB$, } \\
 d_Y(v-F(x),\cM_u(x)) > \rho/2 \}.
\end{multline*}
Then, for $x \in K_\eta$, 
\begin{equation}
\label{eq:Egood:large:measure}
\hat{\nu}_x(E_{good}(x)) \ge \beta \hat{\nu_x}(B(F(x),C)).  
\end{equation}
\end{lemma}
\bold{Proof.} Let $n = \dim \cL_{ext}[x]^{(r)}$. By
  considering determinants, it is easy to show that for any $C > 0$ there
  exists a constant $c_n = c_n(C) > 0$ depending on $n$ and $C$ such
  that for any $\eta > 0$ and any points
  $v_1, \dots, v_n$ in a ball of radius $C$ with the property  
  that for all $1 < i \le n$, $v_i$ 
  is not within $\eta$ of the subspace spanned by $v_1, \dots,
  v_{i-1}$, then $v_1, \dots, v_n$ are not within $c_n \eta$ of any
  $n-1$ dimensional subspace. 
Let $k_{max} \in \natls$
denote the smallest integer greater then $1+n/\epsilon_1$, and let
$N = N(\epsilon_1) = \begin{pmatrix} k_{max} \\
  n-1 \end{pmatrix}$. 
Let $\beta_1$, $\rho_1$ and $K_{\eta,N}$ be as in
Lemma~\ref{lemma:cond:proper:subspace:not:full:measure}. 
Let $\beta = \beta(\eta,\epsilon_1) = \beta_1(\eta,N(\epsilon_1))$,
$\rho = \rho(\eta,\epsilon_1) =
\rho_1(\eta,N(\epsilon_1))/c_n$, 
$K_\eta(\epsilon_1) = K_{\eta,N(\epsilon_1)}$. 
Let $E_{bad}(x) = B(F(x),C)\setminus E_{good}(x)$. To
  simplify notation, we choose coordinates so that $F(x) = 0$.
We claim that $E_{bad}(x)$ is contained in the union of the
$\rho_1$-neighborhoods of at most $N$
subspaces. Suppose this is not 
true. Then, for $1 \le k \le k_{max}$
we can inductively pick points $v_1, \dots, v_k \in E_{bad}(x)$ such
that $v_j$ is not within $\rho_1$  
of any of the subspaces spanned by
$v_{i_1}, \dots, v_{i_{n-1}}$ where $i_1 \le \dots \le i_{n-1} < j$. 
Then, any $n$-tuple of points $v_{i_1}, \dots, v_{i_{n}}$ is not
contained within $\rho = c_n \rho_1$ of a single subspace. 
Now, since $v_i \in E_{bad}(x)$, there exists $U_i \subset \cB$ with
$|U_i| \ge \epsilon_1 |\cB|$ such that for all $u \in U_i$, $d_Y(v_i, \cM_u) <
\rho/2$. We now claim that for any $1 \le i_1 < i_2 < \dots < i_{n}
\le k$, 
\begin{equation}
\label{eq:intersection:empty}
U_{i_1} \cap \dots \cap U_{i_{n}} = \emptyset. 
\end{equation}
Indeed, suppose $u$ belongs to the intersection. Then each of the
$v_{i_1}, \dots v_{i_{n}}$ is within $\rho/2$ of the single subspace
$\cM_u$, but this
contradicts the choice of the $v_i$. This proves
(\ref{eq:intersection:empty}). Now, 
\begin{displaymath}
\epsilon_1 k_{max} |\cB| \le \sum_{i=1}^{k_{max}} 
|U_i| \le n \left|\bigcup_{i=1}^{k_{max}} U_i\right|
  \le n |\cB|. 
\end{displaymath}
This is a contradiction, since $k_{max} > 1+ n/\epsilon_1$. 
This proves the claim. 
Now (\ref{eq:cond:proper:subspace:not:full:measure}) implies that
\begin{displaymath}
\hat{\nu}_x(E_{good}(x)) \ge \hat{\nu}_x( B(F(x),C) \setminus \bigcup_{k=1}^N\Nbhd(\cM_k(x), \rho_1) )
\ge \beta \hat{\nu}_x(B(F(x),C)). 
\end{displaymath}
\qed\medskip

\bold{Proof of Proposition~\ref{prop:can:avoid:most:Mu}.}
Let 
\begin{displaymath}
K'' = \{ x \in X_0 \st \nu_{W^-[x]}( K' \cap \gB_0^-[x]) \ge (1-\delta^{1/2})
\nu_{W^-[x]}(\gB_0^-[x]) \}. 
\end{displaymath}
By Lemma~\ref{lemma:gB:vitali}, we have $\nu(K'') \ge
1-\delta^{1/2}$.

We have, for $x \in K''$, 
\begin{equation}
\label{eq:tmp:hat:nu:x:F:Ktwoprime}
\hat{\nu}_x (F(K' \cap \gB_0^-[x])) \ge (1-\delta^{1/2}) \hat{\nu}_x(\cL_{ext}[x]^{(r)}).
\end{equation}

Let $\beta(\eta,\epsilon_1)$ be as in
Lemma~\ref{lemma:avoid:most:Mu}. 
Let
\begin{displaymath}
c(\delta) = \delta + \inf\{ (\eta^2 + \epsilon_1^2)^{1/2} \st
\beta(\eta,\epsilon_1) \ge 8 \delta^{1/2} \}. 
\end{displaymath}
We have $c(\delta) \to 0$ as $\delta \to 0$. By the definition of
$c(\delta)$ we can choose $\eta =
\eta(\delta) < c(\delta)$ and $\epsilon_1 = \epsilon_1(\delta) <
c(\delta)$ so that $\beta(\eta,\epsilon_1) \ge 8\delta^{1/2}$. 

Now suppose $x \in K'' \cap K'_\delta$. Then, by
(\ref{eq:def:Kprime:delta}) and (\ref{eq:tmp:hat:nu:x:F:Ktwoprime}),
\begin{equation}
\label{eq:tmp:intersection:FKprime}
\hat{\nu}_x(F(K' \cap \gB_0^-[x]) \cap B(F(x),C)) \ge 
(1-2 \delta^{1/2}) \hat{\nu}_x(B(F(x),C)). 
\end{equation}
\mcc{check last formula}
By (\ref{eq:Egood:large:measure}), for $x \in K_\eta$, 
\begin{equation}
\label{eq:tmp:nu:x:Egood:eight:delta}
\hat{\nu}_x(E_{good}(x)) \ge 8 \delta^{1/2} \hat{\nu}_x(B(F(x),C)). 
\end{equation}
Let $K = K' \cap K'' \cap K'_\delta \cap K_\eta$. We have $\nu(K) \ge
1-\delta - 2\delta^{1/2} - c(\delta)$, so $\nu(K) \to 1$ as $\delta \to 0$. 
Also, if $q \in K$, by
(\ref{eq:tmp:intersection:FKprime}) and
(\ref{eq:tmp:nu:x:Egood:eight:delta}), 
\begin{displaymath}
F(K' \cap \gB_0^-[q]) \cap E_{good}(q) \cap
B(F(x),C) \ne \emptyset.
\end{displaymath}
Thus, we can  choose $q' \in K' \cap \gB_0^-[q]$ such that $F(q') \in
E_{good}(q) \cap B(F(q),C)$. Then (\ref{eq:qprime:avoids:Mu}) holds
with $\rho = \rho(\eta(\delta),\epsilon_1(\delta)) > 0$. Also the
upper bound in (\ref{eq:rho:prime:delta:le:d:q:qprime}) holds since
$\gB_0^-[q]$ has diameter at most $1/100$, and the upper bound in 
(\ref{eq:rho:delta:le:Fq:minus:Fqprime}) holds since $F(q') \in
B(F(q),C)$. Since all $\cM_{u}(q)$ contain the origin $q$, the lower
bound in (\ref{eq:rho:delta:le:Fq:minus:Fqprime}) follows from
(\ref{eq:qprime:avoids:Mu}). Finally, the lower bound in
(\ref{eq:rho:prime:delta:le:d:q:qprime}) follows from 
lower bound in (\ref{eq:rho:delta:le:Fq:minus:Fqprime}) since
  in view of (\ref{eq:P1x:example}), $q-q'$ is essentially a component
  of $F(q)-F(q')$.
\qed\medskip

\section{Divergence of generalized subspaces} 
\label{sec:divergence:subspaces}

\bold{The groups $\cG$, $\cG_+$ and $\cG_{++}$.}
Recall that \index{$H^1(x)$}$H^1(x)$ 
denotes $H^1(M,\Sigma,\reals)$. (In fact the dependence on
$x$ is superfluous, but we find it useful to consider $H^1(x)$ as the
fiber over $X_0$ of a flat bundle.) 
Let $\index{$G$@$\cG(x)$}\cG(x)=(SL(H^1) \ltimes H^1)(x)$ which is isomorphic to the group
of affine maps of $H^1(x)$ to itself. We can write $g \in \cG(x)$ as a
pair $(L,v)$ where $L \in SL(H^1(x))$ and $v \in H^1(x)$. We call $L$
the linear part of $g$, and $v$ the translational part. 
 
Let \index{$Q_+(x)$}$Q_+(x)$ denote the group of
linear maps from $H^1(x)$ to itself which preserve the
flag $\{  0 \}
\subset \cV_{\le 1}(H^1)(x) \subset \dots \subset \cV_{\le k}(H^1)(x)
= H^1(x)$, and let $\index{$Q_{++}(x)$}Q_{++}(x)
\subset Q_+(x)$ denote the unipotent subgroup of maps which are the
identity on $\cV_{\le i}(H^1)(x)/\cV_{< i}(H^1)(x)$ for all $i$.
Let \index{$G$@$\cG_+(x)$}$\cG_+(x)$ 
denote the subgroup of $\cG(x)$ in which the linear
part lies in $Q_+(x)$, and let \index{$G$@$\cG_{++}(x)$}$\cG_{++}(x)$
denote the subgroup of $\cG_+(x)$ in which the linear
part lies in $Q_{++}(x)$. Note that $\cG_{++}(x)$ is unipotent. Also, 
since $W^+(x) = \cV_{\le k-1}(H^1)(x)$, $\cG_{++}(x)$ preserves
$W^+(x)$. 

For $y$ near $x$, we have the Gauss-Manin connection $\index{$P^{GM}(x,y)$}P^{GM}(x,y):
H^1(x) \to H^1(y)$. This induces a map $\index{$P^{GM}_*$}P^{GM}_*(x,y): \cG(x) \to
\cG(y)$. In view of
Lemma~\ref{lemma:properties:lyapunov:flag}, for $y \in W^+[x]$,
\begin{multline*}
P^{GM}_*(y,x) \cG_+(y) = \cG_+(x), \qquad 
P^{GM}_*(y,x) Q_+(y) = Q_+(x), \\
P^{GM}_*(y,x) Q_{++}(y) = Q_{++}(x) \qquad
\text{and} \qquad P^{GM}_*(y,x)\cG_{++}(y) = \cG_{++}(x). 
\end{multline*}

We may consider elements of $\cG_+(x)$ and
$\cG_{++}(x)$ as affine 
maps from $W^+[x]$ to $W^+[x]$. More precisely, $g = (L,v) \in \cG(x)$
corresponds to the affine map $W^+[x] \to W^+[x]$ given by:
\begin{equation}
\label{eq:def:affine:action}
z \to x+L(z-x) + v. 
\end{equation}
Then, $Q_{++}(x)$ is the stabilizer of $x$ in $\cG_{++}(x)$. 
We denote by \index{$L$@$\Lie(\cG_{++})(x)$}$\Lie(\cG_{++})(x)$ 
the Lie algebra of $\cG_{++}(x)$, etc.

We will often identify $W^+(x)$ with the translational part of $\Lie(\cG_{++})(x)$. Then, we have an exponential map $\index{$\exp$}\exp:
W^+(x) \to \cG_{++}(x)$, taking $v \in W^+(x)$ to $\exp v \in
\cG_{++}(x)$. Then, $\exp v: W^+[x] \to W^+[x]$ is translation by $v$. 

\bold{The maps $Tr(x,y)$ and $tr(x,y)$.}
For $h \in \cG(x)$, let $\index{$Conj$}Conj(h)$ to be the conjugation
map $g \to h g h^{-1}$, and let $\index{$Ad$}Ad(h): \Lie(\cG)(x) \to
\Lie(\cG)(x)$ be the adjoint map. 
Suppose $y \in W^+[x]$. Let $\index{$Tr(x,y)$}Tr(x,y): \cG(x) \to \cG(y)$ 
and $\index{$tr$@$tr(x,y)$}tr(x,y): \Lie(\cG)(x) \to \Lie(\cG)(y)$ be 
defined as 
\begin{displaymath}
Tr(x,y) = P^{GM}_*(x,y) \circ Conj(\exp(x-y)),  
\end{displaymath}
\begin{displaymath}
tr(x,y) = P^{GM}_*(x,y) \circ Ad(\exp(x-y)).
\end{displaymath}
The following lemma is clear from the definitions:
\begin{lemma}
\label{lemma:trivial:affine:maps}
Suppose $y \in W^+[x]$. Then the elements $g_x \in \cG(x)$ and $g_y
\in \cG(y)$ correspond to the same affine map of $W^+[x] = W^+[y]$ 
(in the sense of (\ref{eq:def:affine:action})) 
if and only if $g_y = Tr(x,y) g_x$. 
\end{lemma}

\bold{Admissible Partitions.}
By an admissible measurable partition we mean any partition $\gB_0$
as constructed in
\S\ref{sec:semi:markov} (with some choice of $\cC$ and $T_0(x)$).

\bold{Generalized subspaces.}
Let $\index{$U'(x)$}U'(x) \subset \cG_{++}(x)$ be a connected Lie subgroup. We write
\begin{displaymath}
\index{$U'[x]$}U'[x] = \{ u x  \st u \in U'(x) \}
\end{displaymath}
and call $U'[x]$ a generalized subspace. We have $U'[x] \subset
W^+[x]$. 

\begin{definition}
\label{def:compatible:family}
Suppose that for almost all $x \in X_0$ we have a distinguished subgroup
\index{$U^+(x)$}$U^+(x)$ of $\cG_{++}(x)$. We say that the family of subgroups
$U^+(x)$ is \textit{compatible with $\nu$} if the following hold:
\begin{itemize}
\item[{\rm (i)}] The assignment $x \to U^+(x)$ is measurable and
  $g_t$-equivariant.  
\item[{\rm (ii)}] For any admissible measurable partition $\gB'$
  of $X_0$, the sets of the form
$\index{$U^+[x]$}U^+[x] \cap 
  \gB'[x]$ are a measurable partition of $X_0$. 
\item[{\rm (iii)}] For any admissible measurable
  partition $\gB'$ of $X_0$, for almost every $x \in X_0$, the conditional
  measure of $\nu$
  along $U^+[x] \cap \gB'[x]$ is a multiple of the unique $U^+(x)$
  invariant measure on 
  $U^+[x] \isom U^+(x)/(U^+(x) \cap Q_{++}(x))$. (Note that both $U^+(x)$ and $U^+(x) \cap
  Q_{++}(x)$ are unimodular, since they are unipotent. Hence there is
  a well-defined Haar measure on the quotient $U^+(x)/(U^+(x) \cap
  Q_{++}(x))$). 
\item[{\rm (iv)}] We have, for almost all $x \in X_0$ and almost all $u \in U^+(x)$, 
\begin{equation}
\label{eq:Uplus:ux:Ux:weird}
U^+(u x) = Tr(x,ux) U^+(x). 
\end{equation}
(This is motivated by Lemma~\ref{lemma:trivial:affine:maps} and
the fact that we want $U^+[ux] = U^+[x]$). Thus, 
\begin{equation}
\label{eq:LieUplus:x:ux:transform}
\index{$L$@$\Lie(U^+)$}\Lie(U^+)(ux) = tr(x,ux) \Lie(U^+)(x).
\end{equation}
\item[{\rm (v)}] $U^+(x) \supset \exp N(x)$ where $\index{$N(x)$}N(x) \subset W^+(x)$ is the
  direction of the orbit of the unipotent $N \subset SL(2,\reals)$. 
\end{itemize}
\end{definition}

\bold{Standing Assumption.} 
We are assuming that for almost every $x \in X_0$ there is a distinguished 
subgroup $U^+(x)$ of $\cG_{++}(x)$ so that the family of subgroups
$U^+(x)$ is compatible with $\nu$ in the sense of
Definition~\ref{def:compatible:family}. This will be used as an
inductive assumption in \S\ref{sec:inductive:step}.

We emphasize that $U^+(x)$ is defined for $x \in X_0$. Using our
notational conventions, for $x \in X$, we write $U^+(x)$ for
$U^+(\sigma_0(x))$ etc.

\bold{The unipotent $N$ as a compatible system of measures.}
At the
start of the induction we have
$U^+(x) = \exp N(x) \subset \cG_{++}(x)$. We now verify that
$U^+(x) = \exp N(x)$ is a family of subgroups compatible with $\nu$ in
the sense of Definition~\ref{def:compatible:family}. Note that
$N(x) = \cV_{\le 1}(H^1)(x) = \cV_1(H^1)(x)$. In particular, by
Lemma~\ref{lemma:properties:lyapunov:flag}, for $y \in W^+[x]$,
\begin{equation}
\label{eq:NyNx}
N(y) = P_{GM}(x,y) N(x).
\end{equation}
This implies (i) and (ii) of Definition~\ref{def:compatible:family}. 

The subgroup $U^+(x) = \exp N(x) \subset \cG_{++}(x)$ consists of pure
translations (i.e. $U^+(x) \cap \cQ_{++}(x)$ is only the identity
map). In particular, $U^+[x] = N[x]$. This, together with the
$N$-invariance of $\nu$ implies (iii) of
Definition~\ref{def:compatible:family}. 

Note that since $U^+(x)$ consists of pure translations,
for any $y \in W^+[x]$, $Conj(\exp(y-x))(U^+(x)) = U^+(x)$. This,
together with (\ref{eq:NyNx}) implies (iv) of
Definition~\ref{def:compatible:family}.


\bold{The sets $\cB[x]$, $\cB_t[x]$ and $\cB(x)$.}
Recall the partitions $\gB_t[x]$ from \S\ref{sec:semi:markov}. 
Let \index{$B$@$\cB_t[x]$} $\cB_t[x] = U^+[x] \cap \gB_t[x]$. We will also use the notation
\index{$B$@$\cB[x]$} $\cB[x]$ for $\cB_0[x]$.

For notational reasons, we will make the following construction: 
let 
\begin{displaymath}
\index{$B$@$\cB_t(x)$} \cB_t(x) = \{ u \in U^+(x)/(U^+(x) \cap Q_{++}(x)) 
\st u x \in \cB_t[x] \}.  
\end{displaymath}
We also write \index{$B$@$\cB(x)$}$\cB(x)$ for $\cB_0(x)$. 

\bold{The Haar measure.}
Let $|\cdot|$ denote the conditional measure of $\nu$ on $\cB[x]$. (By
our assumptions, this measure is $U^+(x)$-invariant where it makes sense.)  
We also denote the Haar measure (with some normalization) on $\cB(x)$
by \index{abs@$\abs{\cdot}$}$|\cdot|$. Unless otherwise specified,
  all statements will be independent of the choice of normalization.

The same argument as Lemma~\ref{lemma:gB:vitali} 
also proves the following:
\begin{lemma}
\label{lemma:cB:vitali:substitute}
Suppose $\delta > 0$, $\theta' > 0$ and $K \subset X$, with 
$\nu(K) > 1-\delta$. 
Then there exists a subset $K^* \subset
K$ with $\nu(K^*) > 1 - \delta/\theta'$ 
such that for any $x \in K^*$, and any $t > 0$, 
\begin{displaymath}
|K \cap \cB_t[x]| \ge (1-\theta') |\cB_t[x]|,  
\end{displaymath}
and thus
\begin{displaymath}
|\{ u \in \cB_t(x) \st u x \in K \}| \ge (1-\theta') 
|\cB_t(x)|. 
\end{displaymath}
\end{lemma}

\bold{The ``ball'' $\cB(x,r)$.}
For notational reasons, for 
$0 < r \le 1/50$, and $x \in X_0$ we define
\begin{displaymath}
\index{$B$@$\cB(x,r)$} \cB(x,r) = \{ u \in U^+(x)/(U^+(x) \cap Q_{++}(x)) 
\st d^+(u x,x) <r \},  
\end{displaymath}
where $d^+(\cdot, \cdot)$ is as in \S\ref{sec:semi:markov}. 
In view of Proposition~\ref{prop:AGY:regularity}, we will normally use
the ball \index{$B$@$\cB(x,1/100)$}$\cB(x,1/100) \subset U^+(x)/(U^+(x) \cap Q_{++}(x))$.


\bold{Lyapunov subspaces.} Suppose $W$ is a subbundle of $H_{big}$. 
Let $\index{$\lambda_i(W)$}\lambda_1(W) > \lambda_2(W) > \dots >
\lambda_n(W)$ denote the Lyapunov 
exponents of the action of $g_t$ on $W$, and for $x \in X_0$ let
\index{$V$@$\cV_i(W)(x)$}$\cV_i(W)(x)$ denote the corresponding subspaces. 
Let $\index{$V$@$\cV_{\le i}(W)$}\cV_{\le i}(W) = \bigoplus_{j=1}^i \cV_i(W)$.

\bold{Notational convention.} In this subsection, 
we write \index{$V$@$\cV_i(x)$}$\cV_i(x)$, 
\index{$V$@$\cV_{\le i}(x)$}$\cV_{\le i}(x)$ and 
\index{$\lambda_i$}$\lambda_i$ instead of
$\cV_i(\Lie(\cG_{++}))(x)$, $\cV_{\le i}(\Lie(\cG_{++}))(x)$ and $\lambda_i(\Lie(\cG_{++}))$. 
\medskip


Since $\Lie(U^+)(x)$ and $\Lie(Q_{++})(x)$ are equivariant under the $g_t$
action, we have
\begin{displaymath}
\Lie(U^+)(x) = \bigoplus_i \Lie(U^+)(x) \cap \cV_i(x), \qquad \Lie(Q_{++})(x)
= \bigoplus_i \Lie(Q_{++})(x) \cap \cV_i(x). 
\end{displaymath}

\bold{The spaces $\cH_+(x)$ and $\cH_{++}(x)$.}
Let $\index{$H$@$\cH_+(x)$}\cH_+(x) = \Hom(\Lie(U^+)(x),\Lie(\cG_{++})(x))$. (Here, $\Hom$
  means linear maps between vector spaces, not Lie algebra homomorphisms).

For every $M \in
\cH_+(x)$, we can write
\begin{equation}
\label{eq:def:Mij}
M = \sum_{ij} M_{ij} \quad \text{ where $M_{ij} \in \Hom(
    \Lie(U^+)(x) \cap \cV_j(x), \Lie(\cG_{++})(x) \cap \cV_i(x))$. }
\end{equation}
Let
\begin{displaymath}
\index{$H$@$\cH_{++}(x)$}\cH_{++}(x) = \{ M \in \cH_+(x) \st M_{ij} = 0 \text{ if } \lambda_i \le
\lambda_j \}.
\end{displaymath}
Then, $\cH_{++}$ is the direct sum of all the positive Lyapunov subspaces
of the action of $g_t$ on $\cH_+$.

\bold{Parametrization of generalized subspaces.}
Suppose $M \in \cH_+(x)$ is such that $(I+M)\Lie(U^+)(x)$ is a
subalgebra of $\Lie(\cG_{++})(x)$. 
We say that the
pair $(M,v) \in \cH_+(x) \cross W^+(x)$ parametrizes the generalized
subspace $\cU$ if
\begin{displaymath}
\cU = \{ \exp[(I + M)u] \, (x+v) \st u \in \Lie(U^+)(x) \}.
\end{displaymath}
(Thus, $\cU$ is the orbit of the subgroup $\exp[(I+M)\Lie(U^+)(x)]$
through the point $x+v \in W^+[x]$.) 
In this case we write $\cU = \index{$U$@$\cU(M,v)$}\cU(M,v)$. 

\bold{Remark.} In this discussion, $\cU$ is a generalized
  subspace which passes near the point $x \in X_0$. However, $\cU$ need not be
  $U^+[x]$, or even $U^+[y]$ for any $y \in X_0$.

\bold{Remark.} 
From the definitions, it is clear that any
generalized subspace $\cU \subset W^+[x]$ can be parametrized by a
pair $(M,v) \in \cH_+(x) \cross W^+(x)$. Also, if $v = v'$ and
\begin{equation}
\label{eq:alternative:parametrization}
I + M = (I + M') \circ J, 
\end{equation}
where $J: \Lie(U^+)(x) \to \Lie(U^+)(x)$ is a linear map,
then $(M,v) \in \cH_{+}(x) \cross W^+(x)$ and
$(M',v') \in \cH_{+}(x) \cross W^+(x)$ are two parametrizations of the
same generalized subspace $\cU$. 


\bold{Example 1.} We give an example of a non-linear generalized
subspace. (The example does not satisfy condition (v) of
Definition~\ref{def:compatible:family} but this is not relevant for
the discussion). 
Suppose for simplicity that $W^+$ has two Lyapunov exponents
$\lambda_1(W^+)$ and $\lambda_2(W^+)$ with $\lambda_1(W^+) = 2 \lambda_2(W^+)$. Let
$e_1(x)$ and $e_2(x)$ be unit vectors so that $\cV_1(W^+)(x) = \reals e_1(x)$,
and $\cV_2(W^+)(x) = \reals e_2(x)$. 

Let $i: W^+(x) \to \reals^3$ be the map sending $a e_1(x) + b e_2(x)
\to (a,b,1) \in \reals^3$. We identify $W^+(x)$ with its image in
$\reals^3$ under $i$. Then, we can identify 
\begin{displaymath}
\cG_{++}(x) = \begin{pmatrix} 1 & * & * \\ 0 & 1 & * \\ 0 & 0 & 1
\end{pmatrix}, \qquad
\Lie(\cG_{++}(x)) = \begin{pmatrix} 0 & * & * \\ 0 & 0 & * \\ 0 & 0 & 0
\end{pmatrix}.
\end{displaymath}
Suppose 
\begin{displaymath}
U^+(x) = \left\{\begin{pmatrix} 1 & t & \tfrac{t^2}{2} \\ 0 & 1 & t \\ 0 & 0 & 1
\end{pmatrix} \st t \in \reals \right\}, \qquad
\Lie(U^+(x)) = \left\{\begin{pmatrix} 0 & t & 0 \\ 0 & 0 & t \\ 0 & 0 & 0
\end{pmatrix} \st t \in \reals \right\}.
\end{displaymath}
Then, $U^+[x]$ is the parabola $\{ x+ t e_2(x) + \frac{t^2}{2}
  e_1(x) \st t \in \reals \} \subset W^+[x]$.

\bold{Transversals.} 
Note that we have, as vector spaces, 
\begin{displaymath}
\Lie(\cG_{++})(x) = \Lie(Q_{++})(x) \oplus W^+(x)
\end{displaymath}
where we
identify $W^+(x)$ with the subspace of $\Lie(\cG_{++})(x)$
corresponding to pure translations.

For each $i$, and each $x \in X_0$, let $\index{$Z_{i1}(x)$}Z_{i1}(x) \subset W^+(x) \cap
\cV_i(x) \subset \Lie(\cG_{++})(x) \cap \cV_i(x)$ be a linear subspace so that
\begin{displaymath}
\Lie(\cG_{++})(x) \cap \cV_i(x) = Z_{i1}(x) \oplus ((\Lie(U^+) +
\Lie(Q_{++}))(x) \cap \cV_i(x)). 
\end{displaymath}
Let $\index{$Z_{i2}(x)$}Z_{i2}(x) \subset \Lie(Q_{++})(x) \cap \cV_i(x)$ be such that
\begin{displaymath}
(\Lie(U^+)+\Lie(Q_{++}))(x) \cap \cV_i(x) = (\Lie(U^+)(x) \cap \cV_i(x))
\oplus Z_{i2}(x). 
\end{displaymath}
Let $\index{$Z_i(x)$}Z_i(x) = Z_{i1}(x) \oplus Z_{i2}(x)$, and let $Z(x) = \bigoplus_i
Z_i(x)$. We always assume that the function $x \to Z(x)$ is measurable.
We say that $\index{$Z(x)$}Z(x) \subset \Lie(\cG_{++})(x)$ is a {\em
  Lyapunov-admissible transversal} to $\Lie(U^+)(x)$. All of our
transversals will be of this type, so we will sometimes simply use the
word ``transversal''. 

Note that $Z_{i1}(x) = Z(x) \cap W^+(x) \cap \cV_i(x)$.

\bold{Example 2.} Suppose $U^+(x)$ is as in Example~1. Then, (since
$\lambda_1(W^+) - \lambda_2(W^+) = \lambda_2(W^+)$),
\begin{displaymath}
\lambda_1 \equiv \lambda_1(\Lie(\cG_{++}) = \lambda_1(W^+) \qquad 
\lambda_2 \equiv \lambda_2(\Lie(\cG_{++}) = \lambda_2(W^+),
\end{displaymath}
\begin{displaymath}
\cV_1 \equiv \cV_1(\Lie(\cG_{++})) = \begin{pmatrix} 0 & 0 & *
  \\ 0 & 0 & 0 \\ 0 & 0 & 0 \end{pmatrix},
\cV_2 \equiv \cV_2(\Lie(\cG_{++}))(x) = \begin{pmatrix} 0 & * & 0 \\ 0 & 0 & * \\ 0 & 0 & 0
\end{pmatrix}, 
\end{displaymath}
\begin{displaymath}
(\Lie(Q_{++}) \cap \cV_2)(x) = \begin{pmatrix} 0 & * & 0 \\ 0 & 0 & 0 \\ 0 & 0 & 0
\end{pmatrix}, \quad
(\Lie(U^+) \cap \cV_2)(x) = \left\{\begin{pmatrix} 0 & t & 0 \\ 0 & 0 & t \\ 0 & 0 & 0
\end{pmatrix} \st t \in \reals \right\},
\end{displaymath}
and $(\Lie(U^+) \cap \cV_1)(x) = (\Lie(Q_{++}) \cap \cV_1)(x) = \{0\}$. Therefore,
$Z_{12}(x) = \{0\}$, and 
\begin{displaymath}
Z_{22}(x) = 
\begin{pmatrix} 0 & * & 0 \\ 0 & 0 & 0 \\ 0 & 0 & 0
\end{pmatrix} \qquad Z_{11}(x) = \begin{pmatrix} 0 & 0 & * \\ 0 & 0
  & 0 \\ 0 & 0 & 0 \end{pmatrix}, \qquad Z_{21}(x) = \{0 \}. 
\end{displaymath}

We note that in this example, the transversal $Z$ was uniquely
determined (and is in fact invariant under the flow $g_t$). 
This is a consequence of the fact that we chose an example
with simple Lyapunov spectrum, and would not be true in general.

\bold{Parametrization adapted to a transversal.}
We say that the parametrization $(M,v) \in \cH_+(x) \cross W^+(x)$ of a
generalized subspace $\cU = \cU(M,v)$ is
adapted to the transversal $Z(x)$ if
\begin{displaymath}
v \in Z(x) \cap W^+(x)
\end{displaymath}
and
\begin{displaymath}
Mu \in Z(x) \quad \text{ for all $u \in \Lie(U^+)(x)$.}
\end{displaymath}
\smallskip

The following lemma implies that adapting a parametrization to a
transversal is similar to inverting a nilpotent matrix.
\begin{lemma}
\label{lemma:reparametrization}
Suppose the pair $(M',v') \in \cH_{++}(x) \cross W^+(x)$ parametrizes a
generalized subspace $\cU$.  Let $Z(x)$ be a Lyapunov-admissible
transversal. Then, there exists a unique pair $(M,v) \in \cH_{++}(x)
\cross W^+(x)$ which parametrizes
$\cU$ and is adapted to $Z(x)$. If we write
\begin{displaymath}
M' = \sum_{ij} M'_{ij}
\end{displaymath}
as in (\ref{eq:def:Mij}), 
and 
\begin{displaymath}
v' = \sum_j v'_j,
\end{displaymath}
where $v'_j \in W^+(x) \cap \cV_j(x)$, 
then $M = \sum_{ij} M_{ij}$ and $v = \sum_{i} v_i$ 
are given by formulas of the form
\begin{equation}
\label{eq:reparam:v}
v_i = L_i v_i' + p_i(v', M')
\end{equation}
\begin{equation}
\label{eq:reparam:M}
M_{ij} = L_{ij} M'_{ij} + p_{ij}(M')
\end{equation}
where $L_i$ is a linear map and
$p_i$ is a polynomial in the $v'_j$ and $M'_{jk}$ which depends
only on the $v'_j$ with $\lambda_j < \lambda_i$ and the $M'_{jk}$ with
$\lambda_j - \lambda_k < \lambda_i$. Similarly, $L_{ij}$
is a linear map, and $p_{ij}$ is a polynomial which depends on the
$M_{kl}'$ with $\lambda_k - \lambda_l <
  \lambda_i - \lambda_j$. 

If we assume in addition that $(M',v')$ is adapted to 
another Lyapunov-admissible transversal $Z'(x)$, then 
$L_i$ and $L_{ij}$ can be taken to be invertible linear maps
(depending only on $Z(x)$ and $Z'(x)$). 
\end{lemma}

The proof of Lemma~\ref{lemma:reparametrization} is a straightforward
but tedious calculation. It is done in 
\S\ref{sec:starredsubsec:proof:reparametrization}.

\bold{The map $S_{x}^Z$.} 
Suppose $Z$ is a Lyapunov-admissible transversal to $U^+(x)$. Then, let
$\index{$S_{x}^Z$}S_{x}^Z: \cH_{++}(x)\cross W^+(x) \to \cH_{++}(x) \cross W^+(x)$ be given by 
\begin{displaymath}
S_{x}^Z(M',v') = (M,v)
\end{displaymath}
where $M$ and $v$ are given by (\ref{eq:reparam:M}) and
(\ref{eq:reparam:v}) respectively. Note that $S_{x}^Z$ is a
polynomial, but is {\em
  not} a linear map in the entries of $M'$ and $v'$.  To deal with the
non-linearity, we work with certain tensor product spaces defined below. 

\bold{Tensor Products: the spaces $\hat{\bH}$, $\tilde{\bH}$ 
and the maps $\bfj$.}
As in \S\ref{sec:conditional}, for a vector space $V$ and a map $f: V
\to W$ we use
the notations $V^{\tensor m}$, $V^{\uplus m}$, $f^{\tensor m}$,
$f^{\uplus m}$, $j^{\tensor m}$, $j^{\uplus m}$.

Let $m$ be the number of distinct Lyapunov exponents on $\cH_{++}$, and
let $n$ be the number of distinct Lyapunov exponents on $W^+$.  
Let $(\alpha;\beta) = (\alpha_1, \dots, \alpha_m; \beta_1, \dots,
\beta_n)$ be a multi-index,  and let 
\begin{displaymath}
\tilde{\bH}^{(\alpha;\beta)}(x) = \bigotimes_{i=1}^m (\cV_i(\cH_{++})(x))^{\tensor \alpha_i} \tensor \bigotimes_{j=1}^n (\cV_j(W^+)(x))^{\tensor \beta_j}
\end{displaymath}
and let
\begin{displaymath}
\hat{\bH}^{(\alpha;\beta)}(x) = \bigotimes_{i=1}^m \cH_{++}(x)^{\tensor
    \alpha_i} \tensor \bigotimes_{j=1}^n W^+(x)^{\tensor \beta_j}. 
\end{displaymath}
We have a natural map $\hat{\pi}^{(\alpha;\beta)}:
\hat{\bH}^{(\alpha;\beta)}(x) 
\to \tilde{\bH}^{(\alpha;\beta)}(x)$ given by
\begin{multline*}
\hat{\pi}^{(\alpha;\beta)}( Y_1 \tensor \dots
\tensor Y_m \tensor (Y_1') 
\tensor \dots \tensor (Y_n') ) = \\
= \pi_1^{\tensor \alpha_1}(Y_1)
\tensor \dots \tensor \pi_m^{\tensor \alpha_m}(Y_m) \tensor
(\pi_1')^{\tensor \beta_1}(Y_1') \tensor \dots
\tensor (\pi_n')^{\tensor \beta_n}(Y_n'),
\end{multline*}
where $\pi_i: \cH_{++}(x) \to \cV_i(\cH_{++})(x)$ and $\pi_j': W^+(x)
\to \cV_j(W^+)(x)$ are the natural projections associated to the direct sum
decompositions $\cH_{++}(x) = \bigoplus_{i=1}^m \cV_i(\cH_{++})(x)$ and
$W^+(x) = \bigoplus_{j=1}^n \cV_j(W^+)(x)$.

Let $\cS$ be a finite collection of multi-indices (chosen in
Lemma~\ref{lemma:tensor:product} below). 
Then, let
\begin{equation}
\label{eq:def:tilde:bH}
\index{$H$@$\tilde{\bH}_0(x)$}\tilde{\bH}_0(x) = \bigoplus_{(\alpha;\beta) \in \cS}
\tilde{\bH}^{(\alpha;\beta)}, 
\qquad 
\index{$H$@$\hat{\bH}_0(x)$}\hat{\bH}_0(x) = \bigoplus_{(\alpha;\beta) \in \cS} \hat{\bH}^{(\alpha;\beta)}
\end{equation}
Let $\index{$\pi$@$\hat{\pi}$}\hat{\pi}: 
\hat{\bH}_0(x) \to \tilde{\bH}_0(x)$ be the linear map
with coincides with $\hat{\pi}^{(\alpha;\beta)}$ on each
$\hat{\bH}^{(\alpha;\beta)}$.

Let $\hat{\bfj}^{(\alpha;\beta)}: \cH_{++}(x) \cross W^+(x) \to
\hat{\bH}^{(\alpha;\beta)}(x)$ be the ``diagonal embedding''
\begin{displaymath}
\hat{\bfj}^{(\alpha;\beta)}(M,v) = M \tensor M \dots \tensor M \tensor v
\tensor \dots \tensor v,
\end{displaymath}
and let $\index{$j$@$\hat{\bfj}$}\hat{\bfj}: \cH_{++}(x) \cross W^+(x) \to \hat{\bH}_0(x)$ be the
linear map 
$\bigoplus_{(\alpha;\beta) \in \cS} \hat{\bfj}^{(\alpha;\beta)}$. Let
\begin{equation}
\label{eq:tmp:def:bfj}
\index{$j$@$\bfj$}\bfj: \cH_{++}(x) \cross W^+(x) \to \tilde{\bH}_0(x)
\end{equation}
denote $\hat{\pi} \circ \hat{\bfj}$.
Let \index{$H$@$\hat{\bH}(x)$}$\hat{\bH}(x)$ denote the linear span of the image of
$\hat{\bfj}$, and let \index{$H$@$\tilde{\bH}(x)$}$\tilde{\bH}(x)$ 
denote the linear span of the image of $\bfj$.

\bold{Induced linear maps on $\hat{\bH}(x)$ and $\tilde{\bH}(x)$.}
Suppose $F_t: \cH_{++}(x) \to \cH_{++}(y)$ 
and $F_t': W^+(x) \to W^+(y)$ are
linear maps. Let $f_t = (F_t, F_t')$. Then, $f_t$ induces a linear map
$\hat{\bff}_t: \hat{\bH}(x) \to \hat{\bH}(y)$. 
If $F_t$ sends each $\cV_i(\cH_{++})(x)$ to each $\cV_i(\cH_{++})(y)$ and $F_t'$ sends each
$\cV_j(W^+)(x)$ to $\cV_j(W^+)(y)$, then $f_t$ also induces a
linear map $\tilde{\bff}_t: \tilde{\bH}(x) \to \tilde{\bH}(y)$. 

\medskip

Note that $\tilde{\bH}(x) \subset \hat{\bH}(x) \subset H_{big}^{(++)}(x)$ where
$H_{big}^{(++)}(x)$ is as in \S\ref{sec:semi:markov}. 


\bold{Notation.}
For an invertible linear map $A: W^+(x) \to W^+(y)$, let $\index{$A_*$}A_*: \Lie(\cG_{++})(x) \to
\Lie(\cG_{++})(y)$ denote the map
\begin{equation}
\label{eq:def:star:notation}
A_*(Y) = A \circ Y_1 \circ A^{-1} + A \circ Y_2
\end{equation}
where for $Y \in \Lie(\cG_{++})(x)$, $Y_1$ is the linear part of $Y$ and
$Y_2$ is the pure translation part.

\begin{lemma}
\label{lemma:def:u:lower:star}
Suppose $x \in X_0$, $u \in U^+(x)$. Then, there exists a linear map
$\index{$u_*$}u_*: \cH_{++}(x) \cross W^+(x) \to \cH_{++}(ux) \cross
W^+(ux)$ with the following properties:
\begin{itemize}
\item[{\rm (a)}] If $(M',v') \in \cH_{++}(x) \cross W^+(x)$ parametrizes
  a generalized subspace $\cU$, then $(M,v) = u_*(M',v')$ parametrizes the
  same generalized subspace $\cU$.
\item[{\rm (b)}] If $(M,v) = u_*(M',v')$, then $M$ and $v$ are given
  by formulas of the form (\ref{eq:reparam:v}) and (\ref{eq:reparam:M}).
\end{itemize}
\end{lemma}

\bold{Proof.} In fact we claim that
\begin{equation}
\label{eq:def:u:lower:star}
u_*(M',v') = (tr(x,ux) \circ M' \circ tr( u x, x), \exp((I+M')Y)(x+v') - \exp(Y)x), 
\end{equation}
where $Y = \log u$. 

This can be verified as follows. 
Let $\cU = \cU(M',v')$ denote the generalized subspace
parametrized by $(M',v')$, 
 and let $U' = \exp((I + M') \Lie(U^+)(x))$, so that $U'$ is a subgroup of
$\cG_{++}(x)$. Then, for any $w \in \cU$, $\cU = U' w$. Then, in view
of Lemma~\ref{lemma:trivial:affine:maps} and 
(\ref{eq:def:affine:action}), 
\begin{displaymath}
\cU = Tr(x,ux) U' (ux + (w -ux)).
\end{displaymath}
Thus, $(M,v) \in \cH_{++}(ux) \cross W^+(ux)$ parametrizes $\cU$ if
\begin{equation}
\label{eq:Tr:first:compatible}
\exp((I+M)\Lie(U^+)(ux)) = Tr(x,ux) U'
\end{equation}
and
\begin{equation}
\label{eq:Tr:second:compatible}
v = w - ux  \qquad\text{ for some $w \in \cU$.}
\end{equation}
Now let $(M,v)$ be the right-hand-side of
(\ref{eq:def:u:lower:star}). 
We claim that (\ref{eq:Tr:first:compatible}) and (\ref{eq:Tr:second:compatible}) hold. 

Indeed, by (\ref{eq:LieUplus:x:ux:transform}), 
\begin{displaymath}
tr(ux,x) \Lie(U^+)(ux) = \Lie(U^+)(x), 
\end{displaymath}
and furthermore, $tr(ux,x) (\Lie(U^+)\cap \cV_{\le i})(ux) =
(\Lie(U^+)\cap \cV_{\le i})(x)$. Now, 
\begin{multline*}
Tr(x,ux) U' = \exp(tr(x,ux)\Lie(U')) = 
\exp(tr(x,ux)(I+M')\Lie(U^+)(x))= \\
\exp(tr(x,ux)(I+M')tr(ux,x)\Lie(U^+(ux))) = \exp((I+M)\Lie(U^+)(ux)), 
\end{multline*}
verifying (\ref{eq:Tr:first:compatible}). 
Also,  let
\begin{displaymath}
w=\exp((I+M')Y)(x+v') \in \cU = \cU(M',v'). 
\end{displaymath}
Therefore, since $\exp(Y) x = u x$, 
\begin{displaymath}
w -ux =  (\exp((I+M')Y)(x+v') - \exp(Y)x) =v,
\end{displaymath}
and hence (\ref{eq:Tr:second:compatible}) holds. 
Thus, $u_*(M',v') \in \cH_{++}(ux) \cross W^+(u x)$ as defined in 
(\ref{eq:def:u:lower:star}) parametrizes the same generalized subspace
$\cU$ as $(M',v') \in \cH_{++}(x) \cross W^+(x)$.  This completes the proof
of part (a). 

It is clear from (\ref{eq:def:u:lower:star}) that part (b) of the lemma
holds.
\qed\medskip

\begin{lemma}
\label{lemma:tensor:product}
For an appropriate choice of $\cS$, the following hold:
\begin{itemize}
\item[{\rm (a)}] Let $Z(x)$ be a Lyapunov-admissible transversal to $U^+(x)$. 
There exists a linear map  $\index{$S$@$\bS_{x}^{Z(x)}$}\bS_{x}^{Z(x)} : \tilde{\bH}(x) \to \tilde{\bH}(x)$
such that for all $(M,v) \in \cH_{++}(x) \cross W^+(x)$, 
\begin{displaymath}
(\bS_{x}^{Z(x)} \circ \bfj)(M,v) = (\bfj \circ S_{x}^{Z(x)})(M,v). 
\end{displaymath}
\item[{\rm (b)}] Suppose $u \in U^+(x)$, and 
let $Z(ux)$ be a Lyapunov-admissible transversal to $U^+(ux)$. Then, there
exists a linear map $\index{$u$@$(u)_*$}(u)_*: 
\tilde{\bH}(x) \to \tilde{\bH}(ux)$ such
that for all $(M,v) \in \cH_{++}(x) \cross W^+(x)$, 
\begin{displaymath}
((u)_* \circ \bfj )(M,v) = (\bfj \circ S_{ux}^{Z(ux)} \circ u_*)(M,v), 
\end{displaymath}
where $u_*: \cH_{++}(x) \cross W^+(x) \to \cH_{++}( u x)
\cross W^+(u x)$ is as in (\ref{eq:def:u:lower:star}). 
\end{itemize}
\end{lemma}

\bold{Proof.} Part (a) formally follows from the universal property of the
tensor product and the partial ordering in  (\ref{eq:reparam:v}) and
 (\ref{eq:reparam:M}). 
We now make a brief outline: see also Example~3 below. 

Let $\tilde{\bH}_{\cS}(x)$ and $\bfj^{\cS}$ 
be as in (\ref{eq:def:tilde:bH}) and (\ref{eq:tmp:def:bfj}) with the
dependence on $\cS$ explicit. Let $\cS_0$ denote the set of
multi-indices of the form $(0, \dots, 0,1, 0,\dots, 0; 0, \dots, 0)$ or $(0,
\dots, 0; 0, \dots, 0,1, 0, \dots, 0)$. Then $\bfj^{\cS_0}$ is an
isomorphism between $\cH_{++}(x) \cross W^+(x)$ and $\tilde{\bH}_{\cS_0}(x)$.

Let $(M,v) = S_x^{Z(x)}(M',v')$.  
By (\ref{eq:reparam:v}), (\ref{eq:reparam:M}) and 
the universal property of the tensor product, there
exists $\cS_1 \supset \cS_0$ 
and a linear map $\bS_1: \tilde{\bH}_{\cS_1}(x) \to
\tilde{\bH}_{\cS_0}(x)$ such that 
\begin{displaymath}
\bfj^{\cS_0}(M,v) = \bS_1 \circ \bfj^{\cS_1}(M',v').
\end{displaymath}
We now repeat this procedure to get a sequence $\cS_j$ of
multi-indices. More precisely, 
at each stage, for each $(\alpha;\beta) \in \cS_j$, 
we may write, by (\ref{eq:reparam:v}), (\ref{eq:reparam:M}) and the
universal property of the tensor product,
\begin{displaymath}
\bfj^{(\alpha;\beta)}(M,v) = \bL^{(\alpha;\beta)}
\left(\bfj^{(\alpha;\beta)}(M,v)\right) + \bS_{j+1}^{(\alpha;\beta)}\left(
\bigoplus_{(\alpha';\beta') \in \cS(\alpha;\beta)}
\bfj^{(\alpha';\beta')}(M',v')\right),
\end{displaymath}
where $\bL^{(\alpha;\beta)}$ and 
$\bS_{j+1}^{(\alpha;\beta)}$ are linear maps; we then define
$\cS_{j+1} = \cS_j \cup \bigcup_{(\alpha;\beta) \in \cS_j}
\cS(\alpha;\beta)$. Putting these maps together, we then get a linear
map $\bS_j$ such that
\begin{displaymath}
\bfj^{\cS_j}(M,v) = \bS_j \circ \bfj^{\cS_{j+1}}(M',v').
\end{displaymath}
Because of the partial order in (\ref{eq:reparam:v})
and  (\ref{eq:reparam:M}), we may assume that $\cS(\alpha;\beta)$ consists of
multi-indices $(\alpha';\beta')$ where either $\alpha'$ has more zero
entries than $\alpha$ or $\beta'$ has more zero entries than
$\beta$. 
Therefore, this procedure eventually terminates, so that
$\cS_{j+1} = \cS_j$ for large enough $j$. 
We then define $\cS$ to be the eventual common
value of the $\cS_j$; then part (a) of
Lemma~\ref{lemma:tensor:product} holds.

To prove part (b) of Lemma~\ref{lemma:tensor:product}, 
note that part (b) of Lemma~\ref{lemma:def:u:lower:star} and the proof
of part (a) of Lemma~\ref{lemma:tensor:product} show that there exists a
map $\tilde{u}_*: \tilde{\bH}(x) \to \tilde{\bH}(ux)$ 
such that $\tilde{u}_* \circ
\bfj = \bfj \circ u_*$, where $u_*$ is as in
(\ref{eq:def:u:lower:star}). Now, we can define $(u)_* : \tilde{\bH}(x) \to
\tilde{\bH}(ux)$ to be $\bS_{ux}^{Z(ux)} \circ \tilde{u}_*$, where
$\bS_{ux}^{Z(ux)}$ is as in (a). Thus $(u)_*$ denotes the
  induced action of $u$ on $\bH(x)$.
\qed\medskip

\bold{Example 3.}  Suppose $U^+$ is as in Example~1 and Example~2. Let 
\begin{displaymath}
F = \begin{pmatrix} 0 & 1 & 0 \\ 0 & 0 & 1 \\ 0 & 0 & 0
\end{pmatrix}, E_1 = \begin{pmatrix} 0 & 0 & 1
  \\ 0 & 0 & 0 \\ 0 & 0 & 0 \end{pmatrix}, E_2 = \begin{pmatrix} 0 & 0 & 0
  \\ 0 & 0 & 1 \\ 0 & 0 & 0 \end{pmatrix}.
\end{displaymath}
Then, $(\Lie(U^+)\cap \cV_2)(x) = \reals F$, $(\Lie(\cG_{++})\cap \cV_1)(x) =
\reals E_1$. Then, 
for $M \in \cH_{++}(x)$, the only non-zero component is $M_{12} \in 
\Hom((\Lie(U^+)\cap \cV_2)(x), (\Lie(\cG_{++})\cap \cV_1)(x))$, which is
$1$-dimensional. Let 
\begin{displaymath}
\Psi \in
\Hom((\Lie(U^+)\cap \cV_2)(x), (\Lie(\cG_{++})\cap \cV_1)(x)) 
\end{displaymath}
denote the element such that $\Psi F = E_1$, so that $\cH_{++} = \reals
\Psi$. 

With the choice of transversal $Z$ given in Example~2, the equations
(\ref{eq:reparam:v}) and (\ref{eq:reparam:M}) become:
\begin{equation}
\label{eq:explicit:example}
v_1 = -M_{12}' v_2' + v_1' -(v_2')^2, \quad v_2 = 0, \quad M_{12} = M_{12}'.
\end{equation}
Then we can choose $\cS = \{ (1; 0,0), (0; 1,0), (0;0,1) ,(1; 0, 1), (0;0,2)\}$, so that (dropping the $(x)$),
\begin{displaymath}
\tilde{\bH}_0 = \cH_{++} \oplus \cV_1(W^+) \oplus \cV_2(W^+) \oplus
(\cH_{++} \tensor \cV_2(W^+)) \oplus (\cV_2(W^+)\tensor \cV_2(W^+)).
\end{displaymath}
(Since for any vector space $V$,  $V^{\tensor 0} = \reals$, we have
omitted such factors in the above formula). 
Let $\bS = \bS_{x}^{Z(x)}$. 
Then, the linear map $\bS: \tilde{\bH}(x) \to \tilde{\bH}(x)$ is given by
\begin{displaymath}
\bS(\Psi) = \Psi, \qquad \bS(E_1) = E_1, \qquad \bS(E_2) = 0, \qquad
\bS(\Psi \tensor E_2) = -E_1, \quad \bS(E_2 \tensor E_2) = -E_1.
\end{displaymath}

\bold{Example 4.} We keep all notation from Examples 1-3. 
Suppose $u = \exp Y$, where $Y = t F$. We now compute the map
$(u)_*$. 

Note that by Lemma~\ref{lemma:properties:lyapunov:flag}, 
we have $e_1(u x) = e_1(x)$. 
Also note that by Example~1, 
at $x$, the tangent vector to $U^+[x]$ coincides with
$e_2(x)$. 
Recall that we are assuming that the foliation whose leaves are $U^+[x]$ is
invariant under the geodesic flow. This implies that at the point
$ux$, the tangent vector to the parabola $U^+[x]$ is $e_2(u
x)$. Therefore, 
\begin{displaymath}
e_1(u x) = e_1(x), \qquad e_2(u x) = t e_1(x) + e_2(x). 
\end{displaymath}
\mcc{check above signs}
Therefore, 
\begin{displaymath}
P^+(x,ux) e_1(x) = e_1(ux), \qquad P^+(x,ux) e_2(x) = e_2(ux) = t
e_1(x) + e_2(x). 
\end{displaymath}

Suppose $\cU$ is
parametrized by $(M',v')$, where $M' = M_{12}' \Psi$, $v' = v_1'
e_1(x) + v_2' e_2(x)$. 
Then 
\begin{displaymath}
\exp[(I+M')Y] = 
\begin{pmatrix} 
1 & t & \tfrac{1}{2} t^2 + M_{12}' t \\ 
0 & 1 & t \\
0  & 0 & 1
\end{pmatrix}, \qquad
\exp(Y) = 
\begin{pmatrix} 
1 & t & \tfrac{1}{2}t^2\\ 
0 & 1 & t \\
0  & 0 & 1
\end{pmatrix}.
\end{displaymath}
Therefore, 
\begin{displaymath}
\exp[(I+M')Y] (x+v') - \exp(Y) x = \begin{pmatrix} v_1' + t v_2' + t
  M_{12}' \\
v_2' \\
0
\end{pmatrix}.
\end{displaymath}
Let
$\Psi' \in
\Hom((\Lie(U^+)\cap \cV_2)(ux), (\Lie(\cG_{++})\cap \cV_1)(ux))$
be the analogue of $\Psi$, but at the point $ux$. 
Then, 
\begin{multline*}
u_*(M',v') = u_*(M'_{12} \Psi, v_1' e_1(x) + v_2' e_2(x)) =
(M_{12}'\Psi', (v_1' + 
t v_2' + t M_{12}') e_1(x) + v_2' e_2(x)) \\
= (M_{12}'\Psi', (v_1' + t M_{12}') e_1(u x) + v_2' e_2(u x))
\end{multline*}
Then, in view of  (\ref{eq:explicit:example}), 
$(S_{ux}^{Z(ux)} \circ u_*)(M',v') =
(M_{12} \Psi', v_1 e_1(ux)+ v_2 e_2(ux))$, where
\begin{displaymath}
v_1 = -M_{12}' v_2' + v_1' + t M_{12}' -(v_2')^2, \qquad v_2
= 0, \qquad M_{12} = M_{12}'.
\end{displaymath}
Then, $(u)_*: \tilde{\bH}(x) \to \tilde{\bH}(ux)$ is given by 
\begin{multline*}
(u)_* (\Psi) = \Psi' + t E_1, \quad (u)_*(E_1) = E_1, \quad (u)_*(E_2)
= 0, \\ 
(u)_*(\Psi \tensor E_2) = -E_1, \quad (u)_*(E_2 \tensor E_2) = -E_1.  
\end{multline*}
\mcc{CHECK EXAMPLE 4 CAREFULLY}


\bold{The dynamical system $G_t$.}
Suppose we fix some Lyapunov-admissible transversal $Z(x)$ for every $x \in
X_0$. Suppose $(M,v) \in \cH_{++}(x) \cross W^+(x)$ 
is adapted to $Z(x)$. Let 
\begin{displaymath}
\index{$G_t$}G_t(M,v) = S_{g_t x}^{Z(g_t x)}( g_t \circ M \circ g_t^{-1}, (g_t)_*
v) \in \cH_{++}(g_t x) \cross W^+(g_t x), 
\end{displaymath}
where $(g_t)_*$ on the right-hand side is $g_t$ acting on $W^+(x)$, and 
$g_t$ on the right-hand side is the natural map $\Lie(Q_{++})(x) \to
\Lie(Q_{++})(g_t x)$, which maps $\Lie(U^+)(x)$ to $\Lie(U^+)(g_t
x)$. 
Then, if $\cU'$ is the generalized subspace
parametrized by $(M,v)$ then 
$(M'',v'') = G_t (M,v) \in \cH_{++}(g_t x) \cross W^+(g_t x)$  
parametrizes $g_t\cU'$ and is adapted to $Z(g_t x)$. 
From the definition, we see that 
\begin{displaymath}
G_{t+s} = G_t \circ G_s. 
\end{displaymath}
Also, it is easy to see that for $(M,v) \in \cH_{++}(x) \cross W^+(x)$,
\begin{displaymath}
G_t(M,v) = (g_t \circ M' \circ g_t^{-1}, (g_t)_* v'), \qquad \text{
  where $(M',v') = S_x^{g_t^{-1} Z(g_t x)}(M,v)$. }
\end{displaymath}

\bold{The bundle $\bH(x)$.} Suppose we are given a Lyapunov adapted
transversal $Z(x)$ at each $x \in X_0$. 
Let 
\begin{displaymath}
\index{$H$@$\bH(x)$}\bH(x) = \bS_x^{Z(x)} \tilde{\bH}(x)
\end{displaymath}
denote the image of $\tilde{\bH}(x)$ under $\bS_x^{Z(x)}$. Then, if
$(M,v) \in \cH_{++}(x) \cross W^+(x)$ is adapted to $Z(x)$, then
$\bfj(M,v) \in \bH(x)$. We can also consider \index{$u$@$(u)_*$}$(u)_*$ 
as defined in Lemma~\ref{lemma:tensor:product} (b) to be a map 
\begin{displaymath}
(u)_*: \bH(x) \to \bH(u x).
\end{displaymath}

\bold{The bundle $\bH$ and the flow $g_t$.} Let $Z(x)$ be an
admissible transversal to $U^+(x)$ for every $x \in X_0$. 
Let $\index{$g_t$@$(g_t)_*$}(g_t)_*: \bH(x) \to \bH(g_t x)$ be given by
\begin{equation}
\label{eq:def:action:on:bHplus}
(g_t)_* = \bS_{g_t x}^{Z(g_t x)} \circ \tilde{\bff}_t
\qquad\text{
  where $f_t(M,v) = (g_t \circ M \circ g_t^{-1}, (g_t)_* v)$},
\end{equation}
$\tilde{\bff}_t$ is the map induced by $f_t$ on $\tilde{\bH} \supset \bH$, 
$(g_t)_*$ on the right-hand side is $g_t$ acting on $W^+(x)$, 
$g_t$ on the right-hand side is the natural map $\Lie(U^+)(x) \to
\Lie(U^+)(g_t x)$, and $\bS_{x}^Z$ is as in
Lemma~\ref{lemma:tensor:product}.  
Then $(g_t)_*$ is a linear map, 
and for $(M,v) \in \cH_{++}(x) \cross W^+(x)$, 
\begin{equation}
\label{eq:gt:coupling:Gt}
(g_t)_*(\bfj(M,v)) = \bfj(G_t(M,v)). 
\end{equation}
Since $G_t \circ G_s = G_{t+s}$, and the linear span of 
$\bfj(\cH_{++}(x) \cross W^+(x))$ 
is $\tilde{\bH}(x) \supset \bH(x)$, it follows from (\ref{eq:gt:coupling:Gt})
that $(g_t)_* \circ (g_s)_* = (g_{t+s})_*$.

\begin{lemma}
\label{lemma:independent:repn:u}
$ $
\begin{itemize} 
\item[{\rm (a)}] Suppose $u'x = u x \in U^+[x]$ and $\bfv \in
  \bH(x)$. Then $(u)_* \bfv = (u')_* \bfv$. 
\item[{\rm (b)}] Suppose $u \in U^+(g_t x)$. Then there exists $u' \in
  U^+(x)$ such that $g_t u' x = u g_t x$. Furthermore, for any choice
  of $u'$ satisfying $g_t u' x = u g_t x$ and any $\bfv \in \bH(x)$,
  we have $(u)_* (g_t)_* \bfv = (g_t)_* (u')_* \bfv$. 
\end{itemize}
\end{lemma}
\bold{Proof.} It is enough to prove (a) for $\bfv = \bfj(M,v)$
where $(M,v) \in \cH_{++}(x) \cross W^+(x)$. Let $\cU$ be the generalized
subspace parametrized by $(M,v)$. Then, $(u)_* \bfv = \bfj(M',v')$
where 
$(M',v') \in \cH_{++}(ux) \cross W^+(ux)$ is the (unique) parametrization
of $\cU$ adapted to $Z(ux)$. But then $(u')_* \bfv$ is also a
parametrization of $\cU$ adapted to $Z(ux)$. Therefore $(u')_* \bfv =
(u)_* \bfv$. 

The proof of (b) is essentially the same. 
\qed\medskip

\bold{Choosing $M_0$ and $\cC_0$.} 
For a.e.\ $x \in X$, let $M^+(x) = \|\bS_x^{Z(x)}\|$, and let
\begin{displaymath}
M^-(x) = \sup_{\bfw \in \bS_x^{Z(x)}(\tilde{\bH}(x))} \frac{1}{\|\bfw\|} \inf\{
  \|\bfv\| \st  \bfv \in \tilde{\bH}(x), \ \bS_x^{Z(x)}(\bfv) = \bfw \}.  
\end{displaymath}
Choose $M_0 > 1$ sufficiently large so that $\cC_0 \equiv \{ x \in X_0 \st
\max(M^+(x),M^-(x)) < M_0 \}$ has positive measure. 
Let $\cC \subset \cC_0$ and $T_0: \cC \to \reals$ be as in
Lemma~\ref{lemma:exists:C:T0} (with this choice of $M_0$, $\cC_0$).

\bold{Adjusting the transversal $Z(x)$.} 
For $c \in \cC$, let $E^+[c]$, $t(c)$ and $J_c$ be as in
Proposition~\ref{prop:semi:markov}. For $x \in E^+[c]$ we define $Z(x)
= P^+(c,x)_* Z(c)$, and for $0 \le t < t(c)$, we define $Z(g_{-t} x) =
g_{-t} Z(x)$. This defines $Z(y)$  for  $y \in J_c$. From now on, we assume that
the transversal $Z$ is obtained via this construction.

\begin{lemma}
\label{lemma:Lyapunov:exponents:on:bH}
Let $(g_t)_*: \bH(x) \to \bH(g_t x)$ and $\tilde{\bff}_t: \tilde{\bH}(x) \to
\tilde{\bH}(g_t x)$ be as in (\ref{eq:def:action:on:bHplus}). 
Then the Lyapunov subspaces for $(g_t)_*$ at $x$ are the image under
$\bS_x^{Z(x)}$ of the Lyapunov subspaces of $\tilde{\bff}_t$ at $x$, and 
the Lyapunov exponents of $g_t$ are those Lyapunov exponents of
$\tilde{\bff}_t$ whose Lyapunov subspace at a generic point $x$ is not
contained in the kernel of $\bS_x^{Z(x)}$. 
\end{lemma}

\bold{Proof.} 
Let $\cV_i(\tilde{\bH})(x)$ and \index{$V$@$\cV_i(\bH)(x)$}$\cV_i(\bH)(x)$ 
denote the Lyapunov subspaces of the flow $\tilde{\bff}_t$ and $g_t$
respectively, and let
\index{$\lambda_i(\tilde{\bH})$}$\lambda_i(\tilde{\bH})$  and
$\lambda_i(\bH)$ denote the corresponding Lyapunov exponents.  Then, for $\bfv
\in \cV_i(\tilde{\bH})$, by the
multiplicative ergodic theorem, for every $\epsilon > 0$, 
\begin{displaymath}
\|g_t \bS_x^{Z(x)} \bfv \| = \| \bS_{g_t x}^{Z(g_t x)} \tilde{\bff}_t \bfv \|_Y \le
\|\bS_{g_t x}^{Z(g_t x)}\|  \| \tilde{\bff}_t \bfv \| \le C_\epsilon(x) C_1(g_t
x) e^{\lambda_i(\tilde{\bH}) t  + \epsilon |t|}.
\end{displaymath}
Taking $t \to \infty$ and $t \to
-\infty$ we see that $\lambda_i(\bH) = \lambda_i(\tilde{\bH})$ and $\bS_x^{Z(x)} \bfv \in \cV_i(\bH)(x)$. 
\qed\medskip

\bold{The measurable flat connection $\bP^+(x,y)$.} 
Recall that the measurable flat $g_t$-equivariant 
$W^+$-connection map
$P^+$ on $H^1$ induces a measurable flat $g_t$-equivariant connection on
$H_{big}^{(++)}$, and thus on $\tilde\bH$. We will call this
connection $\tilde{\bP}^+(x,y)$. Then, we can define a measurable flat
$W^+$-connection $\index{$P$@$\bP^+(x,y)$}\bP^+(x,y): \bH(x) \to \bH(y)$ by
\begin{equation}
\label{eq:def:bold:Pplus}
\bP^+(x,y) = \bS_y^{Z(y)} \circ \tilde{\bP}^+(x,y), \qquad y \in W^+[x].
\end{equation}
Without loss of generality, we may assume that
  Lemma~\ref{lemma:jordan:canonical:form} applies to subbundles of
  $\bH$ as well as subbundles of $H_{big}^{(++)}$ (or else we can replace $X$
  by a measurable finite cover). 
Then, Proposition~\ref{prop:sublyapunov:locally:constant}
applies to $\bP^+$.

\bold{The dynamical inner product $\langle\cdot,\cdot\rangle_x$ and
  the dynamical norm $\| \cdot \|_x$ on $\bH$.}
Even though $\bH$ is not formally a subbundle of
$H_{big}^{(++)}$, $\bH \subset \tilde{\bH} \subset
H_{big}^{(++)}$. Thus, the AGY norm makes sense in $\bH$.
Note that by our choices of $\cC_0$ and $M_0$, (\ref{eq:bound:norm:Pplus}) holds for $\bP^+$ in place of $P^+$ (and $1$ in
  place of $M_0$). Then, the proof of 
Proposition~\ref{prop:properties:dynamical:norm:Hbig} goes
through. Thus, Proposition~\ref{prop:properties:dynamical:norm:Hbig}
also applies to $\bH$, with a norm which may be different from the
norm obtained from thinking of $\bH$ as a subset of $H_{big}^{(++)}$.

\subsection{Approximation of generalized subspaces and the map
  $\cA(\cdot, \cdot, \cdot, \cdot)$} 
\label{sec:subsec:themap:cA}

$ $


\bold{Hausdorff distance between generalized subspaces.}
For $x \in \tilde{X}_0$, and two generalized subspaces $\cU'$ and $\cU''$, let 
\index{$hd_x^{X_0}(\cdot,\cdot)$}$hd_x^{X_0}(\cU',\cU'')$ 
denote the Hausdorff distance using
  the metric $d^{X_0}(\cdot, \cdot)$ defined in
  \S\ref{sec:semi:markov} between $\cU' \cap
B^{X_0}(x,1/100)$ and $\cU'' \cap B^{X_0}(x,1/100)$. (The balls
$B^{X_0}(\cdot,\cdot)$ are defined in \S\ref{sec:conditional}). 

\begin{lemma}
\label{lemma:hausdorff:distance:to:norm}
Suppose $x \in \tilde{X}_0$, $(M,v) \in \cH_{++}(x) \cross W^+(x)$, and $$hd_x^{X_0}(U^+[x], \cU(M,v)) \le 1/100.$$ 
\begin{itemize}
\item[{\rm (a)}] We have for some absolute
constant $C > 0$,
\begin{displaymath}
hd_x^{X_0}(U^+[x], \cU(M,v)) \le C \max( \|v\|_Y, \|M\|_Y). 
\end{displaymath}
Also if $(M,v) \in \cH_{++}(x) \cross W^+(x)$ is adapted to $Z(x)$, then there exists $c(x) > 0$ such that
\begin{displaymath}
hd_x^{X_0}(U^+[x],\cU(M,v)) \ge c(x) \max( \|v\|_Y, \|M\|_Y). 
\end{displaymath}
\item[{\rm (b)}] For some $c_1(x) >0$, 
  we have, for $(M,v) \in \cH_{++}(x) \cross W^+(x)$ adapted to $Z(x)$,
\begin{displaymath}
c_1(x) \|\bfj(M,v)\|_Y \le hd_x^{X_0}(U^+[x], \cU(M,v)) \le c_1(x)^{-1}
\|\bfj(M,v)\|_Y. 
\end{displaymath}
\end{itemize}
\end{lemma}
\bold{Proof.} Part (a) is immediate from the definitions and
Proposition~\ref{prop:AGY:regularity}. To see (b)
note that part (a) implies that $\max(\|M\|_Y,\|v\|_Y) = O(1)$, and thus
all the higher order terms in $\bfj(M,v)$ which are
  polynomials in $M_{ij}$ and $v_j$, have size bounded by a constant
  multiple of the size of the first order terms, i.e.\  by  $\max(\|M\|_Y, \|v\|_Y)$. 
\qed\medskip

We will be dealing with Hausdorff distances of particularly
well-behaved sets (i.e. generalized subspaces parametrized by
elements of $\cH_{++}(x) \cross W^+(x)$.) For such subspaces, the
following holds:
\begin{lemma}
\label{lemma:crude:divergence:subspaces}
Suppose $x \in \tilde{X}_0$, and $\cU' \subset W^+[x]$ 
is a generalized subspace.  Then,
\begin{itemize}
\item[{\rm (a)}] We have, for $t \in \reals$, 
\begin{displaymath}
e^{-2|t|} hd^{X_0}_{x}( U^+[x], \cU')
\le hd^{X_0}_{g_t x}(U^+[g_t x], (g_t)_* \cU') \le
e^{2|t|} hd^{X_0}_x( U^+[x], \cU'),
\end{displaymath}
provided the quantity on the right is at most $1/100$. (The first
inequality in the above line holds as long as the quantity in the
middle is at most $1/100$).
\item[{\rm (b)}] Suppose that $\cU'$ is 
parametrized by an element of $\cH_{++}(x)
\cross W^+(x)$.
There exists a function $C: X_0 \to \reals^+$ 
finite almost everywhere and 
$\beta > 0$ depending only on the Lyapunov spectrum, such that, for $t
\ge 0$, 
\begin{displaymath}
C(x)^{-1} e^{\beta t}
hd_x^{X_0}(U^+[x],\cU') \le hd_{g_t x}^{X_0}( U^+[g_t x], (g_t)_*
\cU'),
\end{displaymath}
provided the quantity on the right is at most $1/100$.
Also, for $t < 0$,
\begin{displaymath}
hd_{g_t x}^{X_0}( U^+[g_t x], (g_t)_*
\cU') \le C(x) e^{-\beta|t|} hd_x^{X_0}(U^+[x],\cU'),
\end{displaymath}
provided the quantity on the right is at most $1/100$. 
\end{itemize}
\end{lemma}

\bold{Proof.} Recall that $B^+(x,r) = B^{X_0}(x,r) \cap W^+[x]$
denotes the ball of radius $r$ in the metric $d^+(\cdot,
\cdot)$. Suppose $t \ge 0$.  Note that, by Lemma~\ref{lemma:forni}(d),
for $t > 0$,
\begin{displaymath}
B^+_t[x] \equiv g_t^{-1} B^+(g_t x, 1/100) \subset B^+(x,1/100). 
\end{displaymath}
Note that the action of $g_t$ can expand in any direction by at most
$e^{2t}$, see also  
Lemma~\ref{lemma:forni:upper}. 
Therefore, 
\begin{displaymath}
hd^{X_0}_{g_t x}( (g_t)_* U^+[x], (g_t)_* \cU') \le e^{2t}
hd^{X_0}_x( U^+[x] \cap 
B^+_t[x], \cU' \cap B^+_t[x]) \le e^{2t} hd^{X_0}_x(U^+[x], \cU').
\end{displaymath}
This completes the proof of the second inequality in (a). The first
inequality in (a) follows after renaming $x$ to $g_{t} x$. 

We now begin the proof of (b). We assume $t \ge 0$ (the proof for the
case $t < 0$ is identical). It is enough to show that for any 
$\delta > 0$ there exists $C=C(\delta) < \infty$ and a set $K(\delta)$
with measure at least $1-\delta$ such that for $x \in K(\delta)$ and
$t > 0$,
\begin{equation}
\label{eq:goal:crude:divergence}
C(\delta)^{-1} e^{\beta t}
hd_x^{X_0}(U^+[x],\cU') \le hd_{g_t x}^{X_0}( U^+[g_t x], (g_t)_* \cU').
\end{equation}
For any $\eta> 0$ let $K_\eta$ be
the set where $c_1(x) > \eta$, where $c_1(x)$ is as in
Lemma~\ref{lemma:hausdorff:distance:to:norm}. Choose $\eta$ so that
$K_\eta$ has measure at least $1-\delta/4$. By the Birkhoff ergodic
theorem we may find a set $K'$ of measure at least $1-\delta/2$ and
$t_1 > 0$ such that for $x \in K'$ and
$t > t_1$, there exists $t' \in \reals$ with
$|t-t'| < \epsilon t$, and $g_{t'} x \in K_\eta$.

Let $\alpha > 0$ be as
in Lemma~\ref{lemma:forni}. Choose $\epsilon < \alpha/2$.
By Lemma~\ref{lemma:forni} (c), we may 
find a set $K'' \subset K_\eta$ of measure at least $1-\delta/2$,
and a constant $t_2 = t_2(\delta)$ such that for all $x \in K''$
all $t > t_2$ and all $\bfv \in \bH(x)$,
\begin{equation}
\label{eq:tmp:tmp:bfv}
\|(g_t)_* \bfv \|_Y \ge e^{\alpha t} \|\bfv\|_Y. 
\end{equation}
Let $K(\delta) = K' \cap K''$, and let $t_0 = \max(t_1,t_2)$. 
If $0 \le t < (1+\epsilon)t_0$, then (\ref{eq:goal:crude:divergence})
holds in view of
Lemma~\ref{lemma:crude:divergence:subspaces} (a). Suppose $t >
(1+\epsilon)t_0$, and let $t'$ be as in the definition of $K'$. Since $x \in
K_\eta$ and $g_{t'} x \in K_\eta$, 
by Lemma~\ref{lemma:hausdorff:distance:to:norm} and (\ref{eq:tmp:tmp:bfv}),
\begin{displaymath}
hd_{g_{t'} x}^{X_0}( U^+[g_{t'} x], (g_{t'})_* \cU') \ge \eta^2
e^{\alpha t} hd_x^{X_0}(U^+[x],\cU').
\end{displaymath}
Then, again using Lemma~\ref{lemma:crude:divergence:subspaces} (a), we
get
\begin{displaymath}
hd_{g_t x}^{X_0}( U^+[g_t x], (g_t)_* \cU') \ge e^{-\epsilon t}
hd_{g_{t'} x}^{X_0}( U^+[g_{t'} x], (g_{t'})_* \cU').
\end{displaymath}
Now, (\ref{eq:goal:crude:divergence}) follows, with $\beta =
(\alpha-\epsilon)$.
\qed\medskip

\bold{Motivation.} We work in the universal cover $\tilde{X}_0$. Let
$q_1$, $q_1'$ be as in \S\ref{sec:outline:step1}, so in particular,
$q_1' \in W^-[q_1]$.
Suppose $u \in \cB(q_1,1/100)$ and $t > 0$. Note that the generalized
subspace $U^+[g_t q_1] = U^+[g_t u q_1]$ passes through the point $g_t
u q_1$. If $t$ is not too large, the generalized subspace $U^+[g_t
q_1']$ will pass near $g_t u q_1$. These subspaces are not on the same leaf of
$W^+$ (even though the leaf $W^+[g_t q_1']$ 
containing $U^+[g_t q_1']$ gets closer
to the leaf $W^+[g_t q_1] = W^+[g_t u q_1]$ 
containing $U^+[g_t u q_1]$ as $t \to
\infty$). It is convenient to find a way to ``project'' the part of 
$U^+[g_t q_1']$ near $g_t u q_1$ 
to $W^+[g_t u q_1]$. In particular, we want the projection 
to be again a generalized subspace (i.e.\ an orbit of a subgroup of
$\cG_{++}(g_t u q_1)$). We also want the projection to be exponentially
close, in a ball of radius $1/100$ about $g_t u q_1$, 
to the original generalized subspace $U^+[g_t q_1']$. 
Furthermore, in order to carry out the program outlined in the
beginning of \S\ref{sec:conditional}, we want the pair $(M'',v'')$ 
parametrizing the projection to be such that $\bfj(M'',v'') \in
\bH(g_t u q_1)$ depends polynomially on $P^-(q_1,q_1')$. Then it
will depend linearly on $F(q) - F(q')$ since any fixed degree
polynomial in $P^-(q_1,q_1')$ can be expressed as a linear function of
$F(q) - F(q')$ as long as $r$ in the definition of $\cL_{ext}(q)^{(r)}$
is chosen large enough.

More precisely, we need the following:
\begin{proposition}
\label{prop:reason:cA:F}
Suppose $\alpha_3 > 0$ is a constant. 
We can choose $r$ sufficiently large (depending only on
$\alpha_3$ and the Lyapunov spectrum) so that there exists a linear map
$\index{$A$@$\cA(q_1,u,\ell,t)$}\cA(q_1,u,\ell,t): \cL_{ext}(g_{-\ell}q_1)^{(r)} \to \bH(g_t
u q_1)$, defined for almost all $q_1 \in \tilde{X}_0$, almost all $u \in
U^+[x]$, all $\ell \ge 0$ and all $t \ge 0$, and  
a constant $\alpha_1 > 0$ depending only on $\alpha_3$ and the Lyapunov
spectrum  such that the following hold:
\begin{itemize}
\item[{\rm (i)}] We have
\begin{equation}
\label{eq:reason:cA:F:equivariance}
\cA(q_1, u, \ell+\ell', t+t') = g_{t'} \circ \cA(q_1,u,\ell,t) \circ g_{\ell'}.
\end{equation}
\item[{\rm (ii)}] Suppose $\delta > 0$, and $\ell$ is sufficiently
large depending on $\delta$. There exists a set $K = K(\delta)$ with
$\nu(K) > 1-\delta$ and constants $C_1(\delta)$ and $C_2(\delta)$
such that the following holds: Suppose $q_1 \in \pi^{-1}(K)$.  Let
$q = g_{-\ell} q_1$ (see Figure~\ref{fig:outline}).  Suppose
$q' \in \pi^{-1}(K) \cap W^-[q]$ satisfies the upper bounds in
(\ref{eq:rho:prime:delta:le:d:q:qprime}) and
(\ref{eq:rho:delta:le:Fq:minus:Fqprime}) with the same constant
$\delta$, and write 
$q_1' = g_\ell q'$.  For all $u \in \cB(q_1,1/100)$ such that
$u q_1 \in \pi^{-1}(K)$, and any $t> 0$ such that
\begin{equation}
\label{eq:reason:cA:F:tbound}
t \le \alpha_3 \ell,
\end{equation}
\begin{equation}
\label{eq:reason:cA:F:point:bound}
d^{X_0}(g_t uq_1, U^+[g_t q_1'])   \le 1/100,
\end{equation}
and also 
\begin{equation}
\label{eq:reason:cA:F:alt:newcondition}
C_1(\delta)  e^{-\alpha_1 \ell} \le hd^{X_0}_{g_t u q_1}(U^+[g_t u q_1],
U^+[g_t q_1']), 
\end{equation}
we have
\begin{multline}
\label{eq:reason:cA:norm}
C(g_t u q_1)^{-1} \|\cA(q_1,u,\ell,t)(F(q')-F(q))\|_Y \le \\
\le hd_{g_t u q_1}^{X_0}(U^+[g_t u q_1], U^+[g_t q_1']) \le C(g_t u q_1)
\|\cA(q_1,u,\ell,t)(F(q')-F(q))\|_Y,
\end{multline}
where $C: X_0 \to \reals^+$ is a measurable function finite almost
everywhere.  
\item[{\rm (iii)}] Suppose $\delta$, $\ell$, $q$, $u$, $q'$, $q_1'$, 
  are as in {\rm (ii)}, and $t$ satisfies
  (\ref{eq:reason:cA:F:tbound}) and
  (\ref{eq:reason:cA:F:point:bound}). 
  Then, we have
\begin{equation}
\label{eq:reason:cA:F:unravel}
  \cA(q_1,u,\ell,t)(F(q')-F(q)) = \bfj(M'',v''),
\end{equation}
where the pair $(M'',v'') \in \cH_{++}(g_t u q_1) \cross W^+(g_t u
q_1)$ (which will be chosen in the proof) is adapted to $Z(g_t u q_1)$ and 
parametrizes a generalized subspace $\cU(M'',v'') \subset W^+(g_t u
q_1)$ satisfying 
\begin{equation}
\label{eq:reason:cA:F:estimate}
hd_{g_t u q_1}^{X_0}(U^+[g_t q_1'], \cU(M'',v''))  \le
C_3(\delta) e^{-\alpha_1 \ell}.
\end{equation}
\end{itemize}
\end{proposition}

Part (ii) of Proposition~\ref{prop:reason:cA:F} is key to
resolving ``Technical Problem \#1'' of \S\ref{sec:outline:step1}
(see the discussion at the beginning of \S\ref{sec:conditional}). 
We claim part (ii) of Proposition~\ref{prop:reason:cA:F} follows 
easily from part (iii) of Proposition~\ref{prop:reason:cA:F} and
Lemma~\ref{lemma:hausdorff:distance:to:norm}(b). Indeed, by the
triangle inequality,
\begin{multline}
\label{eq:tmp:simple:triangle}
hd_{g_t u q_1}^{X_0}(U^+[g_t u q_1], U^+[g_t q_1']) = hd_{g_t u
  q_1}^{X_0}(U^+[g_t u q_1], \cU(M'',v'')) + \\ + O(hd_{g_t u q_1}^{X_0}(\cU(M'',v''),
U^+[g_t q_1'])),
\end{multline}
The $O(\cdot)$ term on the
right-hand-side of (\ref{eq:tmp:simple:triangle}) is bounded by
(\ref{eq:reason:cA:F:estimate}), and the size of the first term on the
right-hand-side of (\ref{eq:tmp:simple:triangle}) is comparable to
$\|\bfj(M'',v'')\|_Y$ by Lemma~\ref{lemma:hausdorff:distance:to:norm}
(b). Thus, (\ref{eq:reason:cA:norm}) follows from
(\ref{eq:reason:cA:F:unravel}). 

\begin{lemma}
\label{lemma:reason:cA:F:altcondition}
For any $\delta > 0$, there exists $K' = K'(\delta)
\subset X_0$ with $\nu(K') > 1-c(\delta)$ where $c(\delta) \to 0$ as
$\delta  \to 0$, 
and constants
$C_1'(\delta) > 0$, $C_2'(\delta) > 0$ and $C_4'(\delta)
> 0$ such that in
Proposition~\ref{prop:reason:cA:F}(ii) and (iii), the conditions
(\ref{eq:reason:cA:F:point:bound}) and 
(\ref{eq:reason:cA:F:alt:newcondition}) can be replaced by either
\begin{itemize}
\item[{\rm (a)}]$g_t u q_1 \in K'$ and
\begin{equation}
\label{eq:reason:cA:F:newcondition}
C_1'(\delta)  e^{-\alpha_1 \ell} \le \|\cA(q_1,u,\ell,t)(F(q')-F(q)) \|_Y
\le C_2'(\delta).
\end{equation}
or by  
\item[{\rm (b)}] $ $ 
\begin{equation}
\label{eq:reason:cA:F:hd:Uplus:cU:bound}
C_4'(\delta) e^{-\alpha \ell} \le hd^{X_0}_{g_t u q_1}(U^+[g_t u
q_1], \cU(M'',v'')) \le 1/400,
\end{equation}
where $\cU(M'',v'')$ is as in (\ref{eq:reason:cA:F:unravel}).
\end{itemize}
\end{lemma}

\bold{Proof of Lemma~\ref{lemma:reason:cA:F:altcondition}.}
Let $c_1(x)$ be as in
Lemma~\ref{lemma:hausdorff:distance:to:norm}(b). 
There exists a
compact $K' \subset X_0$ with $\nu(K') > 1-c(\delta)$, with $c(\delta)
\to 0$ as $\delta \to 0$, and a constant
$1 < C'(\delta') < \infty$ with $C'(\delta') \to \infty$ as $\delta
\to 0$ 
such that $c_1(x)^{-1} < C'(\delta')$ for all $x \in
K'$. Then, in view of Lemma~\ref{lemma:hausdorff:distance:to:norm}(b),
there exist $0 < C_1'(\delta) < C_2'(\delta)$ and $C_4'(\delta) > 0$
such that for $t$ such that $g_t u q_1 \in K'$ and
(\ref{eq:reason:cA:F:newcondition}) holds,
(\ref{eq:reason:cA:F:hd:Uplus:cU:bound}) also holds. Thus, it is
enough to show that if for some $t > 0$ (\ref{eq:reason:cA:F:tbound}) and
(\ref{eq:reason:cA:F:hd:Uplus:cU:bound}) hold, then
(\ref{eq:reason:cA:F:point:bound}) and
(\ref{eq:reason:cA:F:alt:newcondition}) also hold. 

Let $t_{max}=\min\{s \in \reals^+ \st d^{X_0}(g_s uq_1, U^+[g_s q_1'])
\ge 1/100\}$, so that for $0 \le t \le t_{max}$
(\ref{eq:reason:cA:F:point:bound}) holds. If $t_{max} \ge \alpha_3 \ell$,
then for $t \in [0,\alpha_3 \ell)$, (\ref{eq:reason:cA:F:point:bound})
is automatically satified. Now assume $t_{max} < \alpha_3 \ell$. 
Then, by the definition of $t_{max}$ and 
Proposition~\ref{prop:reason:cA:F} (iii),
(i.e. (\ref{eq:reason:cA:F:unravel}) and (\ref{eq:reason:cA:F:estimate})),
and assuming $\ell$ is suffciently large (depending on
$\delta$) we have
\begin{displaymath}
d^{X_0}(g_{t_{max}} uq_1, \cU(M'',v'')) \ge 1/200.
\end{displaymath}
Let $\cU_0 = g_{-t_{max}} \cU(M'',v'') \subset W^+[uq_1]$. By
Proposition~\ref{prop:reason:cA:F}(iii), for $0 \le t \le t_{max}$,
$g_t \cU_0$ is parametrized by $(M_t,v_t)$ satisfying
(\ref{eq:reason:cA:F:unravel}). 

Suppose $t>0$ satisfies (\ref{eq:reason:cA:F:tbound}) and
(\ref{eq:reason:cA:F:hd:Uplus:cU:bound}).   
Let 
\begin{displaymath}
t_1 = \max\{ s \in \reals^+ \st d^{X_0}(g_s uq_1, g_s \cU_0) \le
1/200 \}. 
\end{displaymath}
Since by Lemma~\ref{lemma:forni}(iv) the function $s \to d^{X_0}(g_s
uq_1, g_s \cU_0)$ is monotone increasing, we have $t < t_1 \le
t_{max}$. Thus, since $t < t_{max}$,
(\ref{eq:reason:cA:F:point:bound}) holds. In particular,
Proposition~\ref{prop:reason:cA:F}(iii) applies and then, 
(\ref{eq:reason:cA:F:estimate}) and
(\ref{eq:reason:cA:F:hd:Uplus:cU:bound}) (with a proper choice of
$C_4(\delta)$) imply (\ref{eq:reason:cA:F:alt:newcondition}). 
\qed\medskip



\begin{corollary}
\label{cor:reason:cA:F:divergence}
Suppose $\delta$, $\ell$, $q$, $u$, $q'$, $q_1'$, 
are as in Proposition~\ref{prop:reason:cA:F}{\rm (ii)}, and $s\ge 0$ is such that
(\ref{eq:reason:cA:F:tbound}),   (\ref{eq:reason:cA:F:point:bound}),
and (\ref{eq:reason:cA:F:alt:newcondition}) hold for $s$ in place of
$t$. Suppose $t \in \reals$ is such that $0 < t+s < \alpha_3 \ell$. 
Then, there exists $C_4(\delta) > 0$ such that
\begin{itemize}
\item[{\rm (a)}] We have, for $t \in \reals$ such that $0<t+s<
\alpha_3 \ell$, 
\begin{multline*}
e^{-2|t|} hd^{X_0}_{g_s u q_1}( U^+[g_s u q_1], U^+[g_s q_1']) -
C_4(\delta) e^{-\alpha \ell} 
\le hd^{X_0}_{g_{s+t} u q_1}(U^+[g_{s+t} u q_1], U^+[g_{s+t}
q_1']) \le \\ \le
e^{2|t|} hd^{X_0}_{g_s u q_1}( U^+[g_s u q_1], U^+[g_s q_1']) +
C_4(\delta) e^{-\alpha \ell}.
\end{multline*}
provided the quantity on the right is at most $1/800$. (The first
inequality in the above line holds as long as the quantity in the
middle is at most $1/800$).
\item[{\rm (b)}]
There exists a function $C: X_0 \to \reals^+$ 
finite almost everywhere and 
$\beta > 0$ depending only on the Lyapunov spectrum, such that, for $t
\ge 0$, 
\begin{multline*}
C(g_s u q_1)^{-1} e^{\beta t} hd_{g_s u q_1}^{X_0}(U^+[g_s u
q_1],U^+[g_s q_1']) - C_4(\delta) e^{-\alpha \ell}
 \le  \\ \le hd_{g_{s+t} u q_1}^{X_0}( U^+[g_{s+t} u q_1],
U^+[g_{s+t} q_1']),
\end{multline*}
provided the quantity on the right is at most $1/800$.
Also, for $t < 0$,
\begin{multline*}
hd_{g_{s+t} u q_1}^{X_0}( U^+[g_{s+t} u q_1], U^+[g_{s+t} q_1']) \le \\
\le C(g_s u q_1) e^{-\beta|t|} hd_{g_s u q_1}^{X_0}(U^+[g_s u
q_1],U^+[g_s q_1']) + C_4(\delta) e^{-\alpha \ell},
\end{multline*}
provided the quantity on the right is at most $1/800$.
\end{itemize}
\end{corollary}

\bold{Proof.} Suppose $0 \le t \le \alpha_3 \ell$, and $\ell$ is
sufficiently large depending on $\delta$. 
Let $\cU_t$ denote the generalized subspace of
Proposition~\ref{prop:reason:cA:F}(iii). Then, by
Proposition~\ref{prop:reason:cA:F} if
$d^{X_0}(g_t u q_1, \cU^+[g_t q_1']) < 1/200$, then
$d^{X_0}(g_t u q_1, \cU_t) < 1/100$. 
Conversely, by (the proof of)
Lemma~\ref{lemma:reason:cA:F:altcondition}(b), if
$d^{X_0}(g_t u q_1, \cU_t) < 1/400$, then
$d^{X_0}(g_t u q_1, \cU^+[g_t q_1']) < 1/200$. 
Also, by Proposition~\ref{prop:reason:cA:F}(iii) and
Lemma~\ref{lemma:reason:cA:F:altcondition}(b), if 
either of these conditions holds, then
(\ref{eq:reason:cA:F:estimate}) holds. Thus, the corollary follows
from Lemma~\ref{lemma:crude:divergence:subspaces}.
\qed\medskip

Proposition~\ref{prop:reason:cA:F} is proved by constructing a linear map
$\tilde{P}_s(uq_1,q_1'): W^+(u q_1) \to W^+(q_1')$ with nice
properties; then the approximating subspace $\cU(M'',v'')$ is given by
$g_t \tilde{P}_s(uq_1,q_1')^{-1} U^+[q_1']$. 
The construction is technical, and is
postponed to \S\ref{sec:starredsubsec:construction:cA}. Then,
Proposition~\ref{prop:reason:cA:F} is proved in
\S\ref{sec:starredsubsec:proofs:reason:cA:F}.  
From the proof, we will also deduce
the following lemma (which will be used in \S\ref{sec:inductive:step}):

\begin{lemma}
\label{lemma:parts:staying:close}
For every $\delta > 0$ there exists $\epsilon > 0$ and a compact set
$K \subset X_0$ with $\nu(K) >
1-\delta$ so that the following holds: Suppose $\epsilon_0 <
1/100$. Suppose $q \in \pi^{-1}(K)$,
$\ell > 0$ is sufficiently large depending on $\delta$,
and suppose $q' \in W^-[q] \cap \pi^{-1}(K)$ is such that
(\ref{eq:rho:prime:delta:le:d:q:qprime}) and
(\ref{eq:rho:delta:le:Fq:minus:Fqprime}) hold.  
Let $q_1 = g_{\ell} q$,
$q_1' = g_{\ell} q'$ (see Figure~\ref{fig:outline}). Fix $u \in \cB(q_1,1/100)$,
and suppose $t> 0$ is such that 
\begin{displaymath}
hd_{g_t u q_1}^{X_0}(U^+[g_t u q_1], U^+[g_t q_1']) \le \epsilon \ll \epsilon_0. 
\end{displaymath}
Furthermore, suppose $q_1$, $q_1'$, $u q_1$, $q_1'$, and $g_t u q_1$ 
all belong to $\pi^{-1}(K)$. Suppose $x \in U^+[g_t u q_1] \cap
B^{X_0}(g_t u q_1, 1/100)$. Let
\begin{displaymath}
A_t = U^+[ g_t u q_1] \cap B^{X_0}(x,\epsilon_0),
\end{displaymath}
\begin{displaymath}
A_t' = U^+[ g_t q_1'] \cap B^{X_0}(x,\epsilon_0).
\end{displaymath}
Then, 
\begin{displaymath}
\kappa^{-1} \frac{|g_{-t}A_t |}{|U^+[q_1] \cap B^+(q_1,1/100)|} \le
\frac{|g_{-t}A_t'|}{|U^+[q_1'] \cap B^+(q_1',1/100)|} \le \kappa
\frac{|g_{-t}A_t|}{|U^+[q_1] \cap B^+(q_1,1/100)|},
\end{displaymath}
where $\kappa$ depends only on $\delta$ and the Lyapunov spectrum,
the ``Haar measure'' $| \cdot |$ is defined at the
beginning of \S\ref{sec:divergence:subspaces},  and the ball
$B^+(x, r)$ is defined in \S\ref{sec:semi:markov}. Also, 
\begin{displaymath}
hd^{X_0}(g_{-t}A_t, g_{-t}A_t') \le e^{-\alpha \ell},
\end{displaymath}
where $hd^{X_0}(\cdot, \cdot)$ denotes the Hausdorff distance,
and $\alpha$ depends only on the Lyapunov
spectrum. \mccc{make sure proof takes care of $\epsilon_0$}
\end{lemma}

This lemma will also be proved in
\S\ref{sec:starredsubsec:proofs:reason:cA:F}. 

\subsection{The stopping condition.}
\label{sec:subsec:stopping}

We now state and prove Lemma~\ref{lemma:Lplus:Splus:equiv} and
Proposition~\ref{prop:stop:induction:condition} 
which tell us when the
inductive procedure outlined in \S\ref{sec:outline:step1} stops.


Recall the notational conventions \S\ref{eq:def:pi:minus}.

\bold{The sets $L^-(q)$ and $L^-[q]$.}
For a.e $q \in X_0$, let \index{$L^-[q]$}$L^-[q] 
\subset W^-[q]$ denote the smallest real-algebraic subset
containing, for some $\epsilon > 0$, the intersection of the 
ball of radius $\epsilon$ with
the support of the measure $\nu_{W^-[q]}$, which is the conditional measure
of $\nu$ along $W^-[q]$. Then, $L^-[q]$ is $g_t$-equivariant. Since
the action of $g_{-t}$ is expanding along $W^-[q]$, we see that for
almost all $q$ and any $\epsilon > 0$, $L^-[q]$ is the smallest
real-algebraic subset of $W^-(q)$ such that $L^-[q]$ contains
$\operatorname{support} (\nu_{W^-[q]}) \cap B^{X_0}(q,\epsilon)$. 
Let $\index{$L^-(q)$}L^-(q) = L^-[q] - q$.

\bold{The sets $L^+(q)$ and $L^+[q]$.}
Let $\index{$\hat{\pi}^+$}\hat{\pi}^+: W(x) \to W^+(x)$ 
and $\index{$\hat{\pi}^-$}\hat{\pi}^-: W(x) \to W^-(x)$
denote the maps
\begin{displaymath}
\hat{\pi}_{q_1}^+(v) = (1,0) \tensor v, \qquad \hat{\pi}_{q_1}^-(v) =
(0,1) \tensor \pi^-_{q_1}(v), 
\end{displaymath}
where $\pi_{q_1}^-$ is as in (\ref{eq:def:pi:minus}). 
Let $\index{$L^+(q)$}L^+(q) = \hat{\pi}_q^+ \circ (\hat{\pi}_q^{-})^{-1} L^-(q)$, and let
$\index{$L^+[q]$}L^+[q] = q + L^+(q)$.

\bold{The automorphism $h_t$ and the set $S^+[x]$.}
Let $\index{$h_t$}h_t$ denote the automorphism of 
the affine group $\cG_{++}(x)$ which is the
identity on the linear part and multiplication by $e^{2t}$ on the
translational part. For $x \in X_0$, let
\begin{displaymath}
\index{$S^+[x]$}S^+[x] = \bigcap_{t \in \reals} h_t (U^+)[x].
\end{displaymath}
It is clear from the definition that $S^+[x]$ is relatively closed in
$W^+[x]$, $S^+[x]
\subset U^+[x]$, and also $S^+[x]$ is
star-shaped relative to $x$ (so that if $x+v \in S^+[x]$, so is $x +
tv$ for all $t > 0$). 

\begin{lemma}
\label{lemma:Lplus:Splus:equiv}
The following are equivalent:
\begin{itemize}
\item[{\rm (a)}] $L^+[x] \subset S^+[x]$ for almost all $x \in X_0$.
\item[{\rm (b)}] There exists $E \subset X_0$ with $\nu(E) > 0$ such that
    $L^+[x] \subset S^+[x]$ for $x \in E$. 
\item[{\rm (c)}] There exists $E \subset X_0$ with $\nu(E) > 0$ such that
  $L^+[x] \subset U^+[x]$ for $x \in E$. 
\end{itemize}
\end{lemma}

\bold{Proof.} It is immediately clear that (a) implies (b). Also,
since $S^+[x] \subset U^+[x]$, (b) immediately implies (c). It remains
to prove that (c) implies (a). 

Now suppose (c) holds. Let $\Omega \subset X_0$ be the
set such that for $q_1 \in \Omega$, $g_t q_1$ spends a positive
proportion of the time in $E$. 
Then, by the ergodicity of $g_t$, $\Omega$ is conull. For $q_1 \in
\Omega$, we have, for a positive fraction of $t$,
\begin{displaymath}
L^+[g_t q_1] \subset U^+[g_t q_1].
\end{displaymath}
Let $A(x,t)$ denote the Kontsevich-Zorich cocycle. Then $g_t$ acts on
$W^+$ by $e^t A(q_1,t)$ and acts on $W^-$ by $e^{-t}A(q_1,t)$. Therefore,
$L^-(g_t q_1) = e^{-t} A(q_1,t) L^-(q_1)$, and thus $L^+(g_t q_1) = e^{-t}
A(q_1,t) L^+(q_1)$. Also, we have $U^+(g_t q_1) = e^{t} A(q_1,t)
U^+(q_1)$. Thus, for a positive measure set of $t$, we have
\begin{equation}
\label{eq:tmp:Lplus}
L^+(q_1) \subset e^{2t} U^+(q_1) = h_t(U^+)(q_1),
\end{equation}
where $h_t$ is as in the statement of
Proposition~\ref{prop:stop:induction:condition}. 
Since both sides of (\ref{eq:tmp:Lplus})
depend analytically on $t$, we see that (\ref{eq:tmp:Lplus}) holds for all
$t$. Then, $L^+[q_1] \subset S^+[q_1]$. 
\qed\medskip

\begin{proposition}
\label{prop:stop:induction:condition}
Suppose the equivalent conditions of
Lemma~\ref{lemma:Lplus:Splus:equiv} do not hold. 
Then, there exist constants $\alpha'_1 > 0$, $\alpha'_2>0$ and
$\alpha_1'' > 0$ depending only
on the Lyapunov spectrum, such that  for any $\delta > 0$ and any
sufficiently small (depending on $\delta$) $\epsilon > 0$, 
there exist $\ell_0(\delta,\epsilon) > 0$ and  
a compact $K \subset X_0$ with $\nu(K) > 1-\delta$ such
that for $q_1 \in K$ there exists a subset $Q(q_1) \subset \cB(q_1,1/100)$ 
with $|Q(q_1)| > (1-\delta) |\cB(q_1,1/100)|$, 
such that for $\ell > \ell_0(\delta,\epsilon)$,  for $u \in
Q(q_1)$, and for $t> 0$ such that
\begin{equation}
\label{eq:stop:condition:range}
- \alpha_1'' \ell \le \alpha_2' t - \alpha_1' \ell \le 0,
\end{equation}
we have
\begin{equation}
\label{eq:stop:condition:estimate}
\|\cA(q_1,u, \ell, t)\| \ge e^{-\alpha_1' \ell} e^{\alpha_2' t}.
\end{equation}
Consequently, if $\epsilon >0$ is sufficiently small depending on
$\delta$, $\ell > \ell_0(\delta,\epsilon)$, 
$q_1 \in K$, $u \in Q(q_1)$,
and $t > 0$ is chosen to be as small as possible so that
\begin{displaymath}
\|\cA(q_1,u, \ell, t)\| = \epsilon,
\end{displaymath}
then $t < \frac{1}{2}\alpha_3 \ell$, where
$\alpha_3=\alpha_1'/\alpha_2'$  depends only on the
Lyapunov spectrum. 
\end{proposition}

\bold{Remark.} The constant $\alpha_3$ constructed during the proof of
Proposition~\ref{prop:stop:induction:condition} depends only on the
Lyapunov spectrum. This value of $\alpha_3$ is then used in
Proposition~\ref{prop:reason:cA:F} to construct the function
$\cA(\cdot, \cdot, \cdot, \cdot)$, which is referred to in
(\ref{eq:stop:condition:estimate}).

\starredsubsection{Proof of \protect{Proposition~\ref{prop:stop:induction:condition}}}

\begin{lemma}
\label{lemma:Lminus:avoid:algebraic:subset}
Suppose $k \in \natls$, and $\epsilon > 0$. 
For every sufficiently small $\delta > 0$, and every compact $K'$ with
$\nu(K') > 1-\delta$, there exists a constant $\beta(\epsilon,k,\delta) > 0$ 
and compact set $K''= K''(\epsilon,K',k,\delta) \subset K'$ with $\nu(K'') >
1-c_1(\delta)$ where $c_1(\delta) \to 0$ as $\delta \to 0$ such that
the following holds:

Suppose $q \in \pi^{-1}(K'')$ and $H \subset L^-[q]$ is a connected, degree at most $k$,
$\reals$-algebraic set which is a proper subset of $L^-[q]$. Then there
exists $q' \in \pi^{-1}(K') \cap L^-[q]$ with $d^{X_0}(q',q) <
\epsilon$ and 
\begin{displaymath}
d^{X_0}(q', H) > \beta.
\end{displaymath}
\end{lemma}
\bold{Proof.} This argument is virtually identical to the proof of 
Lemma~\ref{lemma:not:supported:finite:union} and of
Lemma~\ref{lemma:cond:proper:subspace:not:full:measure}. 
\qed\medskip

\begin{lemma}
\label{lemma:algebraic:divergence}
Suppose $k \in \natls$, $m \in \natls$, 
$q_1 \in \tilde{X}_0$, and $\cU' \subset W^+[q_1]$ is the image of a
polynomial map of degree at most $k$ from $\reals^m$ to $W^+[q_1]$. 
Suppose furthermore
that $U^+[q_1]$ is also the image of a polynomial map of degree at
most $k$ from $\reals^m$ to $W^+[q_1]$, 
and $\epsilon > 0$ is such that there exists $u \in
\cB(q_1,1/100)$ with 
\begin{displaymath}
d^{X_0}(u q_1, \cU') = \epsilon.
\end{displaymath}
Suppose $\delta > 0$. Then, for at least $(1-\delta)$-fraction of $u
\in \cB(q_1,1/100)$, 
\begin{displaymath}
d^{X_0}(uq_1, \cU' ) > \beta \epsilon,
\end{displaymath}
where $\beta > 0$ depends only on $k$, $m$, $\delta$ and the dimension. 
\end{lemma}

\bold{Proof.} 
This is a compactness argument. If the lemma was false, we would
(after passing to a limit) obtain polynomial maps whose images
are Hausdorff distance $\epsilon > 0$ apart, yet coincide on a
set of measure at least $\delta$. This leads to a contradiction. 
\qed\medskip

The following lemma is stated in terms of the distance $d^{X_0}(\cdot,
\cdot)$. However, in view of Proposition~\ref{prop:AGY:regularity}, it
is equivalent to the analogous statement for the Euclidean distance on
$W^+[x]$. 
\begin{lemma}
\label{lemma:lojasiewicz}
There exists $C: X_0 \to \reals^+$ 
finite a.e and $\alpha > 0$ depending only on
the Lyapunov spectrum such that for all $q_1 \in \tilde{X}_0$ and all $z \in
L^+[x]$ with $d^{X_0}(z,q_1) < 1/100$, 
\begin{displaymath}
d^{X_0}(z, U^+[x]) \ge C(x) d^{X_0}(z, U^+[x] \cap L^+[x])^{\alpha}. 
\end{displaymath}
\end{lemma}

\bold{Proof.} By the \L{}ojasiewicz inequality
\cite[Theorem~2]{Lojasiewicz} for any $x \in \tilde{X}_0$
and any $k$-algebraic sets $U \subset W^+[x]$, $L \subset W^+[x]$, 
and any $z$ with $d^{X_0}(z,x) < 1/100$,
\begin{displaymath}
d^{X_0}(z, U) + d^{X_0}(z,L) \ge c(U,L) d^{X_0}(z, U \cap L)^\alpha,
\end{displaymath}
where $c(U,L) > 0$ and $\alpha > 0$ depends only on $k$ and the
dimension.

In our case, $U = U^+[x]$. $L = L^+[x]$, and $z \in L^+[x]$. The lemma
follows. 
\qed\medskip

Recall that for $x$ near $q_1$,
$\pi_{W^+(q_1)}(x)$ is the unique point in $W^+[q_1] \cap A W^-[x]$.
Let $n_\tau = \begin{pmatrix} 1 & \tau \\ 0 & 1 \end{pmatrix} \subset N
\subset SL(2,\reals)$.
\begin{lemma}
\label{lemma:pi:x:q1prime}
Suppose $q_1 \in \tilde{X}_0$, $q_1' \in W^-[q_1]$. 
Then, we have
\begin{displaymath}
\pi_{W^+(q_1)}(n_\tau q_1') = n_{\tau'} (q_1 + (1,0) \tensor \tau(1+c\tau)^{-1}
(\hat{\pi}_{q_1}^-)^{-1}(q_1'-q_1)), 
\end{displaymath}
where $c = p(v) \wedge p(\Im q_1)$, $q_1' - q_1 = (0,1) \tensor
v$, and $\tau' = (1-c)\tau(1+c\tau)^{-1}$. 
\end{lemma}
\bold{Proof.} Abusing notation, we work in period
coordinates. 
Since $q_1' \in W^-[q_1]$, we can write $q_1' = q_1 + (0,1) \tensor v$, where
$p(v) \wedge p(\Re q_1) = 0$. 
Then, 
\begin{displaymath}
n_\tau q_1' = (1,0) \tensor (\Re q_1 + \tau(\Im q_1 + v)) + (0,1) \tensor
(\Im q_1 + v).
\end{displaymath}
Let 
\begin{displaymath}
w = v + c\tau(1+c\tau)^{-1} \Im q_1. 
\end{displaymath}
Then, $p(w) \wedge p(\Re (n_\tau q_1')) = 0$, and thus, $(0,1) \tensor w
\in W^-(n_\tau q_1')$. 
Therefore,
\begin{displaymath}
n_\tau q_1' - (0,1) \tensor w = (1,0) \tensor ( \Re q_1 + \tau(\Im q_1 + v)) + (0,1)
\tensor (1+c\tau)^{-1} 
\Im q_1 \in W^-[n_\tau q_1'].
\end{displaymath}
We have $\begin{pmatrix} (1+c\tau)^{-1} & 0 \\ 0 & 1+c\tau \end{pmatrix} \in
A$. Therefore, 
\begin{equation}
\label{eq:tmp:10:tensor}
(1,0) \tensor (1+c\tau)^{-1} (\Re q_1 + \tau(\Im q_1 + v)) + (0,1) \tensor
\Im q_1 \in A W^-[n_\tau q_1'].
\end{equation}
It is easy to check that (\ref{eq:tmp:10:tensor}) is in
$W^+[q_1]$. Therefore, 
\begin{multline*}
\pi_{W^+(q_1)}(n_\tau q_1') = (1,0) \tensor (1+c\tau)^{-1} (\Re q_1 + \tau(\Im q_1 + v)) + (0,1) \tensor
\Im q_1 = \\
= q_1 + (1,0) \tensor \tau (1+c\tau)^{-1} (\Im q_1 + v'),
\end{multline*}
where $v' \in H^1_\perp$ is such that $v = c \, \Re q_1 + v'$. 
Also
\begin{multline*}
n_{\tau'}^{-1} \pi_{W^+(q_1)}(n_\tau q_1') = \pi_{W^+(q_1)}(n_\tau q_1') - (1,0)
\tensor \tau' \, \Im q_1 = \\ = \pi_{W^+(q_1)}(n_\tau q_1') - (1,0) \tensor
(1-c)\tau(1+c\tau)^{-1}\, \Im q_1. 
\end{multline*}
Therefore,
\begin{displaymath}
n_{\tau'}^{-1} \pi_{W^+(q_1)}(n_\tau q_1') = q_1 + (1,0) \tensor \tau (1+c\tau)^{-1} (c\, \Im q_1 + v').
\end{displaymath}
Also,
\begin{displaymath}
c \, \Im q_1 + v' = (\pi_{q_1}^-)^{-1}(v) = (\hat{\pi}_{q_1}^-)^{-1}(q_1' - q_1).
\end{displaymath}
This completes the proof of the lemma.
\qed\medskip

\bold{Proof of Proposition~\ref{prop:stop:induction:condition}.}
Suppose the equivalent conditions of
Lemma~\ref{lemma:Lplus:Splus:equiv} do not hold. For $x \in X_0$, let $U^-(x) =
\hat{\pi}_x^{-} \circ (\hat{\pi}_x^+)^{-1} U^+(x)$, and let $U^-[x] =
x + U^-(x)$. 
Then, for a.e $x \in
X_0$, $L^-[x] \not\subset U^-[x]$, and hence $U^-[x] \cap L^-[x]$ is a
proper algebraic subset of $L^-[x]$. 


By Lemma~\ref{lemma:Lminus:avoid:algebraic:subset}, 
there exists a $K' \subset X_0$ with
$\nu(K') > 1-\delta/4$ and $K'' \subset X_0$ with $\nu(K'') > 1-\delta/2$ such
that for any $q \in \pi^{-1}(K'')$ and any degree $k$ proper real algebraic
subset $H$ of $L^-[q]$, there exists $q' \in L^-[q]$ satisfying the
upper bounds in 
(\ref{eq:rho:prime:delta:le:d:q:qprime}) and
(\ref{eq:rho:delta:le:Fq:minus:Fqprime})
such that $d^{X_0}(q', H) > \beta'(\delta)$. 

Now assume that $q \equiv g_{-\ell} q_1 \in \pi^{-1}(K'')$. (We will later
remove this assumption). Then, we apply
Lemma~\ref{lemma:Lminus:avoid:algebraic:subset}  with $H =
g_{-\ell} (U^-[q_1] \cap L^-[q_1])$ to get $q' \in L^-[q] \cap \pi^{-1}(K'')$
satisfying the upper bounds in (\ref{eq:rho:prime:delta:le:d:q:qprime}) and
(\ref{eq:rho:delta:le:Fq:minus:Fqprime})
and so that 
\begin{displaymath}
d^{X_0}(q', g_{-\ell} (U^-[q_1] \cap L^-[q_1])) \ge \beta'(\delta). 
\end{displaymath}

In view of Lemma~\ref{lemma:forni:upper} and
Proposition~\ref{prop:AGY:regularity}, there exists $N > 0$
such that for all $x \in \tilde{X}_0$ and all $y \in W^-[x]$ with
$d^{X_0}(x,y) < 1/100$ and all $t > 1$,  
\begin{displaymath}
d^{X_0}(g_t x, g_t y) > e^{-N t} d^{X_0}(x,y).
\end{displaymath}
Let $q_1' = g_\ell q'$. Then, $q_1' \in L^-[q_1]$, and
\begin{displaymath}
d^{X_0}(q_1', U^-[q_1] \cap L^-[q_1]) \ge \beta'(\delta) e^{-N \ell}. 
\end{displaymath}
Let $z \in L^+[q_1]$ be such that $\hat{\pi}^+_{q_1} \circ
(\hat{\pi}^{-}_{q_1})^{-1} (q_1') = z$. Then, we have
\begin{displaymath}
d^{X_0}(z, U^+[q_1] \cap L^+[q_1]) \ge \beta'(\delta) e^{-N \ell},
\end{displaymath}
and thus by Lemma~\ref{lemma:lojasiewicz}, 
\begin{equation}
\label{eq:d:z:U:plus:q1}
d^{X_0}(z, U^+[q_1]) \ge \beta(\delta) \beta'(\delta) e^{-\alpha N \ell}. 
\end{equation}

Let $\cU = U^+[q_1']$. Then, $\cU$ is a generalized
subspace, and $q_1' \in \cU$. Furthermore,  both 
$\cU$ and $U^+[q_1]$ are invariant under the action of $N \subset SL(2,\reals)$. 

Without loss of generality, we may assume that $\ell$ is large enough
so that the constant $c$ in Lemma~\ref{lemma:pi:x:q1prime} satisfies
$c < 1/2$.  Now choose $\tau$ so that $\tau(1+c\tau)^{-1} =
1$, and let $\tau'$ be as in Lemma~\ref{lemma:pi:x:q1prime}. 

Let $\cU' = \pi_{W^+(q_1)}(\cU)$. Then, since $n_\tau q_1' \in \cU$, we
have, by Lemma~\ref{lemma:pi:x:q1prime}, 
\begin{displaymath}
n_{\tau'} z = \pi_{W^+(q_1)}( n_\tau q_1') \in \cU'.
\end{displaymath}
But, since $U^+[q_1]$ is $N$-invariant and (\ref{eq:d:z:U:plus:q1}) holds,
we have 
\begin{displaymath}
d^{X_0}(n_{\tau'} z, U^+[q_1]) >  \beta''(\delta) e^{-\alpha N \ell}.
\end{displaymath}
Thus, (because of $n_{\tau'} z$ and Lemma~\ref{lemma:algebraic:divergence}),
\begin{displaymath}
hd_{q_1}^{X_0}(U^+[q_1], \cU') > \beta''(\delta)
e^{-\alpha N \ell}.
\end{displaymath}
Then, by Lemma~\ref{lemma:algebraic:divergence}, 
for $(1-\delta)$-fraction of $u \in \cB(q_1,1/100)$, 
\begin{equation}
\label{eq:d:uq1:cUprime:lower}
d^{X_0}(uq_1, \cU') > \beta'''(\delta) e^{-\alpha N \ell}. 
\end{equation}
By Lemma~\ref{lemma:forni}, and Proposition~\ref{prop:AGY:regularity}.
there exists a compact set $K_2$ of measure at least $(1-\delta)$ and
$\lambda_{min}$ depending only on the Lyapunov spectrum 
such that for $x \in \pi^{-1}(K_2)$ and $y \in W^+[x]$, 
\begin{displaymath}
d^{X_0}(g_t x, g_t y) > c(\delta) e^{\lambda_{min} t} d^{X_0}(x,y),
\end{displaymath}
as long as $t > 0$ and $d^{X_0}(g_t x, g_t y) < 1/100$. 
Let $t_0>0$ be the smallest such that $d^{X_0}(g_{t_0} x, g_{t_0}
  \cU') = 1/100$.
Therefore, assuming $u
q_1 \in \pi^{-1}(K_2)$ in addition to (\ref{eq:d:uq1:cUprime:lower})
we have, for $0 < t < t_0$, 
\begin{displaymath}
d^{X_0}(g_t u q_1, g_t \cU') > c(\delta) \beta'''(\delta) e^{\lambda_{min} t -
  \alpha N \ell}.
\end{displaymath}
Hence, for $0 < t < t_0$,
\begin{displaymath}
hd_{g_t u q_1}^{X_0}(U^+[g_t u q_1], g_t \cU') > c_1(\delta) e^{\lambda_{min} t
  - \alpha N \ell},
\end{displaymath}
and thus, in view of Proposition~\ref{prop:reason:cA:F}(iii) and
Lemma~\ref{lemma:reason:cA:F:altcondition}(b),
\begin{displaymath}
hd_{g_t u q_1}^{X_0}(U^+[g_t u q_1], g_t \cU) > c_2(\delta) e^{\lambda_{min} t -
  \alpha  N \ell}.
\end{displaymath}
Let $\alpha_2' = \lambda_{min}/2$, $\alpha_1' = 2 \alpha N$, and
let $\alpha_3=\alpha_1'/\alpha_2'$. Let $\alpha_1 > 0$ be as in
Proposition~\ref{prop:reason:cA:F} for this choice of $\alpha_3$. Then
we can choose $\alpha_1'' > 0$ to be smaller than $\alpha_1$, so that
if (\ref{eq:stop:condition:range}) holds and $\ell$ is
sufficiently large then (\ref{eq:reason:cA:F:alt:newcondition}) holds. 
Hence, by Proposition~\ref{prop:reason:cA:F} (ii) 
if (\ref{eq:stop:condition:range}) holds, $0 < t < t_0$,
 (and assuming that $g_t u q_1 \in
\pi^{-1}(K''')$ where $K'''$ is a compact set of measure at least $1-\delta$),
\begin{displaymath}
\|\cA(q_1, u, \ell, t) (F(q')-F(q))\| \ge c_3(\delta) e^{\lambda_{min}
  t - \alpha N \ell} 
\end{displaymath}
Then, for $0 < t < t_0$ satisfying (\ref{eq:stop:condition:range}),
\begin{displaymath}
\|\cA(q_1, u, \ell, t) \| \ge c_4(\delta) e^{\lambda_{min}
  t - \alpha N \ell}
\end{displaymath}
If $t \ge t_0$ satisfies (\ref{eq:stop:condition:range}), then
\begin{displaymath}
\|\cA(q_1, u, \ell, t) \| \ge \|\cA(q_1, u, \ell, t_0) \| \ge
c_5(\delta) \ge c_5(\delta) e^{\lambda_{min} t - \alpha N \ell}
\end{displaymath}
Thus, for all $t$
such that (\ref{eq:stop:condition:range}) holds, 
\begin{displaymath}
\|\cA(q_1, u, \ell, t) \| \ge c_6(\delta) e^{\lambda_{min}
  t - \alpha N \ell}.
\end{displaymath}
This implies (\ref{eq:stop:condition:estimate}), assuming that $\ell$
is sufficiently large (depending on $\delta$), $q \in
\pi^{-1}(K'')$ and $g_t u q_1 \in \pi^{-1}(K''')$. 

For the general case (i.e.\ without the assumptions that $q \in
\pi^{-1}(K'')$ and $g_t u q_1 \in \pi^{-1}(K''')$), 
note that we can assume that $g_{-\ell} q_1 \in
\pi^{-1}(K'')$ for a set of 
$\ell$ of density at least $(1-2\delta)$, and also 
$g_t u q_1 \in \pi^{-1}(K''')$ for a set of $t$ of density at least
$(1-2\delta)$. 
Now the general case of (\ref{eq:stop:condition:estimate})
follows from the special case, Proposition~\ref{prop:reason:cA:F} (i) and 
Lemma~\ref{lemma:forni:upper}.  
\qed\medskip

\starredsubsection{Proof of Lemma~\ref{lemma:reparametrization}.}
\label{sec:starredsubsec:proof:reparametrization}

We can choose a
subspace $T(x) \subset \Lie(U^+)(x)$,  
so that
\begin{displaymath}
\Lie(U^+)(x) + \Lie(Q_{++})(x) = T(x) \oplus \Lie(Q_{++})(x).
\end{displaymath}
(In particular, if $\Lie(U^+)(x) \cap \Lie(Q_{++}) = \{0 \}$, $T(x) =
\Lie(U^+)(x)$.)
Then, 
\begin{displaymath}
\Lie(\cG_{++})(x) = (Z(x) \cap W^+(x)) \oplus T(x) \oplus \Lie(Q_{++})(x).
\end{displaymath}
Thus, for any vector $Y \in \Lie(\cG_{++})(x)$, we can write
\begin{equation}
\label{eq:three:projections}
Y = \pi_Q(Y) + \pi_Z(Y) + \pi_T(Y), 
\end{equation}
where $\pi_Q(Y) \in
  \Lie(Q_{++})(x)$, $\pi_Z(Y) \in Z(x) \cap W^+(x)$, $\pi_T(Y) \in T(x)$.

Suppose there exists $\tilde{u} \in T(x)$ such that (in $W^+(x)$)
\begin{equation}
\label{eq:tmp:v:equiv}
x+ v \equiv \exp[(I+M')\tilde{u}] (x+v') \in x+ Z(x) \cap W^+(x). 
\end{equation}
Then there exists $q \in \Lie(Q_{++})(x)$, $z \in Z(x) \cap W^+(x)$ such
that in $\cG_{++}(x)$, 
\begin{equation}
\label{eq:tmp:exp:IplusM:u}
\exp[(I+M')\tilde{u}] \exp(v') = \exp(z) \exp (q).
\end{equation}
In this subsection, we write $\cV_i(x)$ for $\cV_i(\Lie(\cG_{++})(x)$,
and $\lambda_i$ for $\lambda_i(\Lie(\cG_{++}))$. We also write 
$\cV_{< i}(x) = \bigoplus_{j=1}^{i-1} \cV_j(x)$. 

Write $\tilde{u} = \sum_i \tilde{u}_i$, where $\tilde{u}_i \in
(\Lie(U^+)\cap \cV_i)(x)$.  
Also, write $q = \sum_i q_i$, where $q_i \in (\Lie(Q_{++})\cap \cV_i)(x)$, 
$v = \sum_i v_i$, where $v_i \in (W^+ \cap \cV_i)(x)$, and
$z = \sum_i z_i$ where $z_i \in Z_{i1}(x) = Z(x) \cap W^+(x)\cap \cV_i(x)$. 

For $h \in \cG_{++}(x)$ we may write $h = h_1 h_2$ where $h_1 \in
\cQ_{++}(x)$, and $h_2 \in W^+(x)$ is a pure translation. Let $\hat{i}(h)$
denote the element of $\Lie(\cG_{++})(x)$ whose linear part is $h_1 - I$
and whose pure translation part is $h_2$. Then, $\hat{i}: \cG_{++}(x) \to
\Lie(\cG_{++})(x)$ is a bijective $g_t$-equivariant map. 

Recall
that our Lyapunov exponents are numbered so that $\lambda_i > \lambda_j$ for
$i < j$. Then, we claim that 
\begin{multline}
\label{eq:tmp:left:exponential}
\hat{i}\bigg(\exp[(I+M')\tilde{u}] \exp(v')\bigg) 
+ \cV_{< i}(x) = \\
= \tilde{u}_i + v_i' +
\hat{i}\left( \exp\left[(I+M') \sum_{j > i} \tilde{u}_j\right]
  \exp\left[\sum_{j > i} v'_j\right] \right) +
\cV_{< i}(x). 
\end{multline}
Indeed, any term involving $\tilde{u}_j$ or $v'_j$ for $j < i$ would
belong to $\cV_{< i}(x)$ (since it  would lie in a subspace with
Lyapunov exponent bigger
than $\lambda_i$). Also, for the same reason, 
any terms involving $\tilde{u}_i$ or $v'_i$
other than those written on the left-hand-side of
(\ref{eq:tmp:left:exponential}) would belong to $\cV_{< i}(x)$. 
Similarly, 
\begin{multline}
\label{eq:tmp:right:exponential}
\hat{i}\bigg(\exp(z) \exp (q)\bigg) + \cV_{< i}(x) = \\ = z_i + q_i +
\hat{i}\left(\exp\left(\sum_{j>i} z_j\right) \exp\left(\sum_{j>i}
  q_j\right)\right) +\cV_{< i}(x). 
\end{multline}
We now apply $\hat{i}$ to both sides of (\ref{eq:tmp:exp:IplusM:u}),
plug in (\ref{eq:tmp:left:exponential}) and
(\ref{eq:tmp:right:exponential}), 
and compare terms in
$\cV_i(x)$. We get equations of the form
\begin{displaymath}
\tilde{u}_i + v_i' + p_i = z_i + q_i, 
\end{displaymath}
where $p_i$ is a polynomial in the $\tilde{u}_j$ and $q_j$ for $\lambda_j <
\lambda_i$, and in the $M'_{jk}$ for $\lambda_j - \lambda_k <
\lambda_i$. Then, the equation can be solved inductively, starting
with the equation with $i$ maximal (and thus $\lambda_i$
minimal). Thus, the equation (\ref{eq:tmp:v:equiv}) can indeed
  be solved for $\tilde{u}$ and we get,
\begin{displaymath}
\tilde{u}_i = -\pi_T (v_i' + p_i), \quad z_i = \pi_Z( v_i' + p_i), \quad
q_i = \pi_Q(v_i' + p_i), 
\end{displaymath}
where $\pi_Q$, $\pi_T$ and $\pi_Z$ as in
(\ref{eq:three:projections}).  This shows that $v = \exp(z) v'$ has
the form given in (\ref{eq:reparam:v}). \mcc{give more details}

Let $U' = \exp((I+M') \Lie(U^+)(x))$. By our assumptions, $U'$ is a
subgroup of $\cG_{++}$. Therefore, for $\tilde{u}$ as in (\ref{eq:tmp:v:equiv}), 
\begin{displaymath}
\cU = U' \cdot(x+v') = U' \exp(-(I+M')\tilde{u}) \cdot(x+v) = U' \cdot(x+v). 
\end{displaymath}
Then, $(M',v)$ is also a parametrization of $\cU$. To make $M'$
adapted to $Z(x)$ we proceed as follows:

For $u \in \Lie(\cG_{++})(x)$, we can write $u = u'' + z''$, where $u''
\in \Lie(U^+)(x)$ and $z'' \in Z(x)$. Let $\pi^Z_{U^+}: \Lie(\cG_{++}) \to
\Lie(U^+)$ be the linear map sending $u$ to $u''$.

In view of (\ref{eq:alternative:parametrization}), we need to find a
linear map $J: \Lie(U^+)(x) \to \Lie(U^+)(x)$, so that if we define
$M$ via the formula (\ref{eq:alternative:parametrization}), then
$M$ is adapted to $Z(x)$. Write $u' = J u$. Then, 
$u'\in \Lie(U^+)(x)$ must be such that $u' + M' u' = u + z$, where $z \in
Z$. Then, 
\begin{displaymath}
u' + \pi^Z_{U^+} (M' u')) = u,
\end{displaymath}
hence $u' = J u$ must be given by the formula
\begin{displaymath}
u' = (I + \pi^Z_{U^+} \circ M')^{-1} u. 
\end{displaymath}
Thus, in view of (\ref{eq:alternative:parametrization}), we define $M$
by
\begin{equation}
\label{eq:tmp:reparam:M}
M = (I+M')(I + \pi^Z_{U^+} \circ M')^{-1} - I. 
\end{equation}
Then for all $u \in \Lie(U^+)(x)$, 
$Mu = (I+M) u - u = (I+M') u' - u  \in Z(x)$. 
Thus $(M,v)$ is adapted to $Z(x)$. 
Since $M' \in \cH_{++}(x)$,
\begin{displaymath}
\pi^{Z}_{U^+} \circ M' = \sum_{i < j} \pi^Z_{U^+} \circ M'_{ij},
\end{displaymath}
where $M_{ij}' \in \Hom(\Lie(U^+)\cap \cV_j, \Lie(\cG_{++})\cap \cV_i)$. 
Since $Z(x)$ is a Lyapunov-admissible transversal, 
$\pi^{Z}_{U^+}$ takes $\Lie(\cG_{++})\cap \cV_j$ to
$\Lie(\cU^+)\cap \cV_i$. Therefore, 
\begin{displaymath}
\pi^Z_{U^+} \circ M_{ij}' \in \Hom(\Lie(U^+)\cap \cV_j,
\Lie(U^+)\cap \cV_i).
\end{displaymath}
Thus, $\pi^{Z}_{U^+} \circ M'$ is nilpotent. 
Then (\ref{eq:reparam:M}) follows from (\ref{eq:tmp:reparam:M}).  

This argument shows the existence of a pair $(M,v)$ which parametrizes
$\cU$ and is adapted to $Z(x)$. The uniqueness follows from the same
argument. Essentially one shows that any $(M,v)$ which parametrizes
$\cU$ and is adapted to $Z(x)$ must satisfy equations whose
unique solution is given by (\ref{eq:reparam:v}) and
(\ref{eq:reparam:M}). 
\qed\medskip

\starredsubsection{Construction of  the map $\cA(q_1, u, \ell, t)$.}
\label{sec:starredsubsec:construction:cA}
$ $

\bold{Motivation.} Suppose $q_1 \in \tilde{X}_0$, $q_1' \in W^-[q_1]$,
$u \in U^+(q_1)$, so $u q_1 \in W^+[q_1]$. 
To construct the generalized subspace $\cU =
\cU(M'',v'')$ of Proposition~\ref{prop:reason:cA:F}, we first let $\cU
= g_t \cU_0$ and construct the generalized subspace $\cU_0 \subset
W^+[uq_1]$. Let $z =
\index{$\pi_{W^+(q_1')}$}\pi_{W^+(q_1')}(uq_1)$, so that $z$ is the unique point in $W^+[q_1']
\cap A W^-[u q_1]$. In particular, $W^+[q_1'] = W^+[z]$. (Note that we
are not assuming any ergodic properties of $z$; in particular the
Lyapunov subspaces at $z$ may not be defined). 

We will construct a $\pi_1(X_0)$-equivariant 
linear map $\index{$P$@$\tilde{P}_s(u q_1,q_1')$}\tilde{P}_s(u q_1,q_1'): W^+(u q_1) \to W^+(z)$, and let
$\cU_0 =  \tilde{P}_s(u q_1,q_1')^{-1} U^+[q_1']$. (This makes sense
since $U^+[q_1'] \subset W^+[q_1'] = W^+[z]$). 
We 
want $\tilde{P}_s(u q_1,q_1')$ to have the following properties:
\begin{itemize}
\item[{\rm (P1)}] 
$\tilde{P}_s(u q_1,q_1')$ depends only on $W^+[q_1']$,
  i.e.\ for $z' \in W^+[q_1']$, we have $\tilde{P}_s(u q_1,z') =
  \tilde{P}_s(uq_1,q_1')$. In particular, for any $u' \in U^+(q_1')$,
  $\tilde{P}_s(uq_1,q_1') = \tilde{P}_s(uq_1,u'q_1')$. 
\item[{\rm (P2)}]  For nearby $x,y \in \tilde{X}_0$, let $P^{GM}(x,y): H^1(x)
  \to H^1(y)$ denote the Gauss-Manin connection. 
For $u \in \cB(q_1,1/100)$,  
  $u' \in \cB(q_1',1/50)$ and $t\ge 0$ with 
\begin{equation}
\label{eq:gt:uq1:gtprime:uprime:close}
d^{X_0}(g_t u q_1, g_t u' q_1') < 1/100, 
\end{equation}
let $z' = \pi_{W^+(g_t u' q_1')}(g_t u q_1)$.
Then,  there exists $\alpha_1 > 0$ depending only on the
  Lyapunov spectrum such that 
$\| \tilde{P}_s( g_t u q_1, g_t u' q_1')^{-1} P^{GM}(g_t u q_1, z') - I \|_Y =
O(e^{-\alpha_1 \ell})$, for all $t>0$ such that
(\ref{eq:gt:uq1:gtprime:uprime:close}) holds. 
(Also note that the points $u q_1$ and $u' q_1'$ satisfy
$d^{X_0}(g_{-\tau} u q_1, g_{-\tau} u' q_1') = O(1)$ for all $0 \le \tau \le \ell$).

Note that as long as (\ref{eq:gt:uq1:gtprime:uprime:close}) 
 holds, $d^{X_0}(g_t u q_1, z')=O(1)$ and $d^{X_0}(z',
g_t u' q_1') = O(1)$ so that $P^{GM}(g_t u q_1, z')$ connects nearby
points. This 
would not be the case if we defined $\tilde{P}_s(u q_1, q_1')$ to be a
linear map from $W^+(uq_1)$ to $W^+(q_1')$, since 
$g_t uq_1$ and $g_t q_1'$ would quickly become far apart. 
\item[{\rm (P3)}] The (entries of the matrix) $\tilde{P}_s( u q_1,
  q_1')^{-1}$ are polynomials of degree at most $s$ in (the entries of the
  matrix) $P^-(q_1,q_1')$. \mccc{say which basis -- refer} 
\item[{\rm (P4)}] The generalized subspace $\cU =  \tilde{P}_s(u
  q_1,q_1')^{-1} U^+[q_1']$ can be parametrized by $(M'',v'') \in \cH_{++}(u
  q_1) \cross W^+(u q_1)$ (and not by an arbitrary element of
  $\cH_+(u q_1) \cross   W^+(u q_1)$). 
\end{itemize}
The construction will take place in several steps.

\bold{Notation.} In this subsection, $\cV_i(x)$ refers to
$\cV_i(H^1)(x)$. 

\bold{The map $\hat{P}(x,y)$.}
There exists a set $K$ of full measure such that each point $x$ in $K$
is Lyapunov-regular with respect to the bundle $W^+$, i.e.\
\begin{displaymath}
H^1(x) = \bigoplus_i \cV_i(x),
\end{displaymath}
where $\cV_i(x) = \cV_i(H^1)(x)$ 
are the Lyapunov subspaces, and the multiplicative
ergodic theorem holds. We have the flag
\begin{equation}
\label{eq:basic:flag}
\{ 0 \} \subset \cV_{\le 1}(x) \subset \dots \subset \cV_{\le n}(x)  = H^1(x),
\end{equation}
where $\cV_{\le j}(x) = \bigoplus_{i=1}^j \cV_i(x)$. Note that
$\cV_{\le n-1}(x) = W^+(x)$. 
If $y \in W^+[x]$ is also Lyapunov-regular, then the flag
(\ref{eq:basic:flag}) at $y$ agrees with the flag at $x$, provided we
identify $H^1(y)$ with $H^1(x)$ using the Gauss-Manin connection. 
Thus, we
may define (\ref{eq:basic:flag}) at any point $x$ such that
$W^+[x]$ contains a regular point. 

Now suppose $x$ and $y$ are restricted to a subset where the $\cV_i$
vary continuously. Then, for nearby $x$ and $y$,
we have, for each $i$, 
\begin{equation}
\label{eq:decomp:two}
H^1(x) = \cV_{\le i}(y) \oplus \bigoplus_{j=i+1}^n \cV_j(x).
\end{equation}
Let $z = \pi_{W^+(y)}(x)$, and let $\hat{P}_i: \cV_i(x) \to H^1(z)$ be
the map taking $v \in \cV_i(x)$ to its $\cV_{\le i}(y)$ component
under the decomposition (\ref{eq:decomp:two}). Let
$\index{$P$@$\hat{P}(x,y)$}\hat{P}(x,y): H^1(x) \to H^1(z)$ be the
linear map which agrees with $\hat{P}_i$ on each $\cV_i(x)$. Note that
$\hat{P}(x,y)$ is defined for all nearby $x$, $y$ such that
(\ref{eq:decomp:two}) holds for all $i$. Let $\hat{P}[x,y]$ be the
affine map from $W^+[x]$ to $W^+[y]$ whose linear part is
$\hat{P}(x,y)$ and such that $x$ maps to $z = \pi_{W^+(y)}(x)$. To
simplify notation, we will denote $\hat{P}[x,y]$ also by
$\hat{P}(x,y)$.

We have
\begin{displaymath}
\hat{P}(g_t x,g_t y) = g_t \circ \hat{P}(x,y) \circ g_{-t},
\end{displaymath}
and
\begin{equation}
\label{eq:hat:P:flag:preserving}
\hat{P}(x,y) \cV_{\le i}(x) = P^{GM}(y,z)\cV_{\le i}(y) =  \cV_{\le i}(z).
\end{equation}
(Since $z \in W^+[y]$, we can define $\cV_{\le i}(z)$ to be
$P^{GM}(y,z)\cV_{\le i}(y)$ even if $\cV_i(z)$ were not originally defined).



The following lemma essentially states that the map
$\hat{P}(uq_1,q_1')$ has properties (P1) and (P2). 
\begin{lemma}
\label{lemma:reduce:to:P:hat}
Suppose $\delta > 0$, $\alpha_3 > 0$ and $\ell$ is sufficiently large
  depending on $\delta$ and $\alpha_3$. 
Suppose $q \in \tilde{X}_0$ and $q' \in W^-[q]$ 
satisfy the upper bounds in
(\ref{eq:rho:prime:delta:le:d:q:qprime}) and
(\ref{eq:rho:delta:le:Fq:minus:Fqprime}). 
Let $q_1 = g_\ell q$ (see Figure~\ref{fig:outline}), and write
$q_1' = g_\ell q'$. Then, 
for almost all $u \in \cB(q_1,1/100)$ 
and $t$ with $0 < t < \alpha_3 \ell$ such that 
\begin{displaymath}
d^{X_0}( g_t u q_1, U^+[g_t q_1']) < 1/100,
\end{displaymath}
the following holds:

Let $\widehat{\cU} = \hat{P}(u q_1, q_1')^{-1}(U^+[q_1'])$. Then
$\widehat{\cU} \subset W^+[q_1]$ is a generalized subspace, and
\begin{displaymath}
hd_{g_t u q_1}^{X_0}( g_t \widehat{\cU}, U^+[g_t q_1']) \le C(q_1) C(uq_1) 
e^{-\alpha (t+\ell)},
\end{displaymath}
where $\alpha > 0$ depends only on $\alpha_3$ and the Lyapunov
spectrum, and $C: X_0 \to \reals^+$
is finite almost everywhere. 
\end{lemma}

\bold{Proof.} In this proof, we write $\cV_i(x)$ for $\cV_i(H^1)(x)$
and $\cV_{\le i}(x)$ for $\cV_{\le
  i}(H^1)(x)$. For convenience, we also choose $u' \in \cB(q_1',1/50)$ with
\begin{displaymath}
d^{X_0}(g_t u q_1, g_t u' q_1') = d^{X_0}(g_t u q_1, U^+[g_t q_1']) \le 1/100.
\end{displaymath}
(Nothing in the proof will depend on the choice of $u'$). 

Let $q_2 = g_t u q_1$, $q_2' = g_t u' q_1'$. 
We claim that 
\begin{equation}
\label{eq:tmp:d:Viq2:Viq2prime}
d_Y(\cV_{\le i}(q_2),P^{GM}(q_2',q_2)\cV_{\le i}(q_2')) \le C(q_1)C(uq_1) e^{-\alpha (t+\ell)},
\end{equation}
where $\alpha > 0$ depends only on the Lyapunov spectrum, and $C: X_0
\to \reals^+$ (which depends on $\delta$) is finite a.e. 

We will apply Lemma~\ref{lemma:subspaces:stay:close} (with $t+\ell$ in
place of $t$) to the points
$x = g_{-(t+\ell)} q_2$ and $y = g_{-(t+\ell)} q_2'$. Thus, we need to
bound $D^+(x,y)$. 
In the following argument, we identify $H^1(x)$, $H^1(y)$, $H^1(q)$
and $H^1(q')$ using the Gauss-Manin connection, while suppressing
$P^{GM}$ from the notation. 

Suppose $v' \in \cV_{\le i}(y)$ realizes the supremum in the
definition of $D^+(x,y)$, i.e.\ $v' = v + w$ where $v \in \cV_{\le
  i}(x)$, $w \in \cV_{>i}(x)$, and $D^+(x,y) = \|w\|_Y/\|v\|_Y$. 
 
Note that $\cV_{\le i}(x) = \cV_{\le i}(q)$ and
$\cV_{\le i}(y) = \cV_{\le i}(q')$.  Thus, $v' \in \cV_{\le
  i}(q')$. Also note that $\cV_{>i}(q') = \cV_{>i}(q)$ for all $i$,
$P^-(q',q) \cV_i(q') = \cV_i(q)$, and by
Lemma~\ref{lemma:properties:Pplus} (c), $P^-(q',q)$ is lower
triangular and unipotent.  By the upper bound in
(\ref{eq:rho:delta:le:Fq:minus:Fqprime}),
$\|P^-(q',q)\|_Y \le C'(\delta)$. (In particular, we have a
  lower bound, depending on $\delta$, on the angles between the
  Lyapunov subpaces $\cV_i(q')$). Hence we can write
\begin{displaymath}
v' = v'' + w'' \qquad \text{$v'' \in \cV_{\le i}(q)$, $w'' \in
  \cV_{>i}(q)$, $\|w''\|_Y \le C(\delta) \|v''\|_Y$.}
\end{displaymath}
Since $\cV_{\le i}(x) = \cV_{\le i}(q)$, we have $v'' \in \cV_{\le
  i}(x)$. By Corollary~\ref{cor:staying:close} (applied with $x =
q_1$, $y = u q_1$ and $t = \ell$) we can write
\begin{displaymath}
w'' = v_2 + w_2 \qquad \text{$v_2 \in \cV_{\le i}(x)$, $w_2 \in
  \cV_{>i}(x)$,  and $\|v_2\|_Y \le C_1(q_1) C_1(u q_1) e^{-\alpha \ell } \|w''\|_Y$. } 
\end{displaymath}
Thus, 
\begin{displaymath}
  v = v'' + v_2, \qquad w = w_2.
\end{displaymath}
If $\ell$ is bounded depending on $C_1(q_1) C_1(uq_1)$ and $\delta$,
then (in view of the condition $t < \alpha_3 \ell$), the desired
estimate (\ref{eq:tmp:d:Viq2:Viq2prime}) is trivially true. Thus, we may assume that $\ell$ is sufficiently large so that
\begin{displaymath}
C_1(q_1) C_1(u q_1) e^{-\alpha \ell} \le 1. 
\end{displaymath}
Then, 
\begin{displaymath}
\|w_2\|_Y \le \|w''\|_Y + \|v_2\|_Y \le 2\|w''\|_Y \le 2 C(\delta) \|v''\|_Y.
\end{displaymath}
But, 
\begin{displaymath}
\|v_2 \|_Y \le C_1(q_1) C_1(u q_1) e^{-\alpha \ell} \|w_2\|_Y \le
2 C(\delta) C_1(q_1) C_1( u q_1) e^{-\alpha \ell} \|v''\|_Y.
\end{displaymath}
Arguing as above, we may assume, without loss of generality,  
that $\ell$ is sufficiently large so that
\begin{displaymath}
\|v \|_Y \ge \|v''\|_Y - \|v_2\|_Y \ge (1/2) \|v''\|_Y. 
\end{displaymath}
Then, 
\begin{displaymath}
D^+(x,y) = \frac{\| w_2 \|_Y}{\|v\|_Y} \le 4 C(\delta).
\end{displaymath}
Hence, by Lemma~\ref{lemma:subspaces:stay:close},
(\ref{eq:tmp:d:Viq2:Viq2prime}) follows. 

By Lemma~\ref{lemma:exists:C:T0} (c), for any $\epsilon > 0$ and any
subset $S$ of the Lyapunov exponents, 
\begin{equation}
\label{eq:transverse:Lyapunov1}
d_Y(\bigoplus_{i \in S} \cV_i(q_2),\bigoplus_{j \not\in S} \cV_j(q_2))
> C_\epsilon(u q_1) e^{-\epsilon t} >
C_\epsilon(u q_1) e^{-\epsilon(t + \ell)}.
\end{equation}
Choose $\epsilon < \alpha/2$, where $\alpha$ is as in
(\ref{eq:tmp:d:Viq2:Viq2prime}). 
Then, by (\ref{eq:transverse:Lyapunov1}), (\ref{eq:tmp:d:Viq2:Viq2prime}),
and the definition of $\hat{P}(q_2,q_2') = \hat{P}(g_t u q_1, g_t u' q_1')$, \mcc{maybe give more details}
\begin{equation}
\label{eq:hatP:minus:I:small}
\|\hat{P}(g_t u q_1, g_t u' q_1')^{-1} P^{GM}(g_t u q_1, g_t' u' q_1')
- I\|_Y \le C_\epsilon'(u q_1) C'(q_1) e^{-\alpha'(\ell + t)},
\end{equation}
where $\alpha' = \alpha - \epsilon$ depends only on the Lyapunov
spectrum, and $C'(\cdot)$, $C_\epsilon'(\cdot)$ are finite a.e. 
Also note that by the upper bound in
(\ref{eq:rho:prime:delta:le:d:q:qprime}) and Lemma~\ref{lemma:forni}, we have
\begin{displaymath}
d_Y(u q_1,z) \le C_\epsilon(q_1) e^{-\alpha' \ell},
\end{displaymath}
and again by Lemma~\ref{lemma:forni},
\begin{equation}
\label{eq:image:z:is:close}
d_Y(g_t u q_1, g_t z) < C_\epsilon(u q_1) e^{-\alpha' t} d_Y(u
  q_1, z) \le C_\epsilon(q_1) C_\epsilon(u q_1) e^{-\alpha'(t+\ell)}.
\end{equation}

Note that $\widehat{\cU}$ is the orbit of a subgroup $\hat{U}$ of
$\cG(u q_1)$ whose Lie algebra is 
\begin{displaymath}
\hat{P}(u q_1,q_1')^{-1}_*
\Lie(U^+)(q_1')
\end{displaymath}
(and we are using the notation
(\ref{eq:def:star:notation})). 
By (\ref{eq:hat:P:flag:preserving}) and the fact that
$\Lie(U^+)(q_1') \in \cG_{++}(q_1')$ we have $\Lie(\hat{U}) \in \cG_{++}(u
q_1)$. Thus, $\widehat{\cU}$ is a generalized subspace. 

Since $U^+[q_1']$ is a generalized subspace, for all $u' \in
U^+(q_1')$, $U^+[q_1'] = U^+[u' q_1']$. We have 
\begin{displaymath}
g_t \widehat{\cU} = g_t \hat{P}(u q_1, u' q_1')^{-1} U^+[u' q_1'] =
\hat{P}(g_t u q_1, g_t u' q_1')^{-1} U^+[g_t u' q_1'].  
\end{displaymath}
Therefore, the lemma follows from (\ref{eq:hatP:minus:I:small})
and (\ref{eq:image:z:is:close}).
\qed\medskip

\bold{Motivation.} Suppose $q_1 \in \tilde{X}_0$,
$u \in U^+(q_1)$, $q_1' \in W^-[q_1]$.
In view of Lemma~\ref{lemma:reduce:to:P:hat},
$\hat{P}(uq_1, q_1')$ has properties (P1) and (P2). We claim that it
does not in general have the properties (P3) and (P4). 

Let $z = \pi_{W^+(q_1')}(uq_1)$ so in particular $\hat{P}(uq_1, q_1')
= \hat{P}(uq_1, z)$ and let
\begin{equation}
\label{eq:def:hat:Q}
\index{$Q$@$\hat{Q}(u q_1; q_1')$}\hat{Q}(u q_1; q_1') = \hat{P}(u q_1, z)^{-1} P^{GM}(q_1', z) P^-(q_1,q_1') \circ
P^+( u q_1, q_1), 
\end{equation}
so that
\begin{equation}
\label{eq:hatP:hatQ}
\hat{P}(uq_1, z) \hat{Q}(u q_1; q_1') = P^{GM}(q_1',z) P^-(q_1, q_1') P^+(u q_1,
q_1). 
\end{equation}
Then, $\hat{Q}( u q_1; q_1'): H^1(u q_1) \to H^1( u
q_1)$ and $\hat{Q}(u q_1; q_1') \cV_{\le i}( u q_1) = \cV_{\le i}( u q_1)$,
hence $\hat{Q}(u q_1; q_1') \in Q_+(u q_1)$. In particular
$\hat{Q}(uq_1; q_1') W^+(uq_1) = W^+(uq_1)$. 

We now show how to compute $\hat{P}(uq_1, q_1')$ and $\hat{Q}(uq_1;
q_1')$ in terms of $P^+ = P^+(u q_1, q_1)$ and $P^- = P^-(q_1,
q_1')$. In view of
Lemma~\ref{lemma:properties:Pplus}, $P^+$ is upper triangular with $1$'s along the diagonal 
in terms of a basis
adapted to $\cV_i(u q_1)$. Also by Lemma~\ref{lemma:properties:Pplus}
applied to $P^-$ instead of $P^+$, 
$P^-$ is lower triangular with $1$'s along the diagonal 
in terms of a basis adapted to
$\cV_i(q_1)$. Therefore, since $P^+$ takes $\cV_i(u q_1)$ to
$\cV_i(q_1)$, $(P^+)^{-1} P^-P^+$ is lower triangular with $1$'s along
the diagonal in terms of a basis adapted to $\cV_i(u q_1)$. 

Let $\hat{P} = \hat{P}(uq_1, q_1')$, $\hat{Q} = \hat{Q}(uq_1; q_1')$.
Then, in view of the definition of $\hat{P}$, $\hat{P}$ is lower
triangular with $1$'s along the diagonal in terms of a basis adapted
to $\cV_i(uq_1)$ (and we identify $H^1(q_1')$ with $H^1(u q_1)$ using
the Gauss-Manin connection). Also, since $\hat{Q}$ preserves the flag
$\cV_{\le i}(u q_1)$, $\hat{Q}$ is upper triangular in terms of the basis
adapted to $\cV_i(u q_1)$. Thus, (\ref{eq:hatP:hatQ}) can we written
as 
\begin{equation}
\label{eq:LU}
\hat{P} \hat{Q} = P^- P^+ = P^+ ((P^+)^{-1} P^- P^+)
\end{equation}
Recall that the Gaussian elimination algorithm shows that any matrix
$A$ in neighborhood of the identity $I$ can be written uniquely as $A = L U$
where $L$ is lower triangular with $1$'s along the diagonal and $U$ is
upper triangular. Thus, $\hat{P} = \hat{P}(uq_1, q_1')$ and $\hat{Q} =
\hat{Q}(uq_1; q_1')$ are the $L$ and $U$ parts of the $LU$
decomposition of the matrix $A = P^-(q_1,q_1') P^+(uq_1, q_1)$. (Note
that we are given $A = U' L'$ where $U' = P^+$ is upper triangular and $L' = 
(P^+)^{-1} P^- P^+$ is lower triangular, so we are really solving the
equation $LU = U' L'$ for $L$ and $U$). 

Since the Gaussian elimination algorithm involves division, the
entries of $\hat{P}(uq_1, q_1')^{-1}$ are rational functions of the
entries of $P^+(uq_1,q_1)$ and $P^-(q_1,q_1')$, but not in general
polynomials. This means that $\hat{P}(uq_1, q_1')$ does not in general
have property (P3). Also, the diagonal entries of $\hat{Q}(uq_1;
q_1')$ are not $1$. This eventually translates to the failure of the
property (P4). Both problems are addressed below.

\bold{The maps $\hat{P}_s( uq_1, q_1')$ and $\tilde{P}_s( u q_1,
  q_1')$.}
 For $s > 1$, 
let \index{$Q$@$\hat{Q}_s(u q_1; q_1')$}$\hat{Q}_s(u q_1; q_1')$ 
be the order $s$ Taylor approximation to $\hat{Q}(u q_1; q_1')$,
where the variables are the entries of $P^-(q_1,q_1')$ (and $u$, $q_1$
and the entries of $P^+(uq_1, q_1)$ are considered constants). 
Then, $\hat{Q}_s = \hat{Q}_s(u q_1; q_1') \in Q_+(u q_1)$. We may write
\begin{displaymath}
\hat{Q}_s = D_s + \tilde{Q}_s,
\end{displaymath}
where $D_s$ preserves all the subspaces
$\cV_i(u q_1)$ and $\tilde{Q}_s = \tilde{Q}_s(uq_1; q_1') \in Q_{++}(u q_1)$.
Let $\tilde{P}_s( uq_1,
q_1')=\tilde{P}_s(uq_1,z)$ be defined by the relation:
\begin{equation}
\label{eq:tilde:P:s}
\tilde{P}_s(u q_1, q_1')^{-1} = \tilde{Q}_s(u q_1; q_1') P^+(q_1, u
q_1) P^-(q_1', q_1) P^{GM}(z,q_1').  
\end{equation}

\bold{Motivation.} We will effectively show that for $s$ sufficiently
large, (chosen at the end of the proof of
Proposition~\ref{prop:reason:cA:F}) the map
$\tilde{P}_s(uq_1,q_1')$ has the properties (P1),(P2), (P3) and (P4). 

We have, by (\ref{eq:tilde:P:s}),  
\begin{displaymath}
\tilde{P}_s(u q_1,q_1')^{-1} \cV_{\le i}(q_1') = 
\tilde{P}_s(u q_1,q_1')^{-1} \cV_{\le i}(z) = \cV_{\le i}(u q_1). 
\end{displaymath}
As a consequence, 
\begin{displaymath}
\tilde{P}_s(u q_1,q_1')^{-1} \circ Y \circ \tilde{P}_s(u q_1,q_1') 
\in \cG_{++}( u q_1)  \qquad
\text{ for all $Y \in   \cG_{++}(q_1')$. }
\end{displaymath}
Thus, for any subalgebra $L$ of $\Lie(\cG_{++})(q_1')$, it follows that $\tilde{P}_s(uq_1,q_1')_*^{-1}(L)$ is a
subalgebra of $\Lie(\cG_{++})(uq_1)$, where $\tilde{P}_s(uq_1, q_1')^{-1}_*:
\Lie(\cG_{++})(q_1') \to \Lie(\cG_{++})(u q_1)$ is as in
(\ref{eq:def:star:notation}).

\bold{The map $i_{u,q_1,s}$.}

\bold{Motivation.} For $q_1 \in X_0$ and $u \in \cB(q_1,1/100)$, 
we want $\index{$i_{u,q_1,s}$}i_{u,q_1,s}: \cL_{ext}(q_1) \to \cH_{++}(u
q_1) \cross W^+( u q_1)$ to be such that 
\begin{displaymath}
i_{u,q_1,s}(\gP(q_1') - \gP(q_1)) = (M_s,v_s),
\end{displaymath}
where the pair $(M_s,v_s) \in \cH_{++}(u q_1) \cross W^+(u q_1)$
parametrizes the approximation $\tilde{P}_s( uq_1, q_1')^{-1}
U^+[q_1']$ to $U^+[q_1']$ constructed above. Furthermore, 
we want $i_{u,q_1,s}$ to be a
polynomial map of degree at most $s$ in the entries of $\gP(q_1') -
\gP(q_1)$. 
\medskip

By Proposition~\ref{prop:sublyapunov:locally:constant} (a), we have
\begin{equation}
\label{eq:Lie:U:plus:q1:prime}
\Lie(U^+)(q_1') = P^-(q_1,q_1')_* \circ P^+(uq_1,q_1)_*(\Lie(U^+)(u q_1)), 
\end{equation}
where we used the notation (\ref{eq:def:star:notation}). 
Let $\cU'_s = \tilde{P}_s( u q_1,q_1')^{-1} U^+[q_1']$. We first find
$(M_s',v'_s) \in \cH_+(q_1) \cross W^+(q_1)$ which parametrizes $\cU'_s$. 
Let 
\begin{displaymath}
v_s = \tilde{P}_s( u q_1,q_1')^{-1} q_1' \in \cU'_s \subset
W^+[q_1] = W^+[u q_1].
\end{displaymath}
By (\ref{eq:Lie:U:plus:q1:prime}), $\cU'_s = U_s \cdot v_s$ where the
subgroup $U_s$ of $\cG_{++}(u q_1)$ is such that
\begin{displaymath}
\Lie(U_s) = \tilde{P}_s(u q_1, z)^{-1}_* \circ P^{GM}(q_1',z)_* \circ
P^-(q_1, q_1')_* \circ P^+(uq_1, q_1)_* (\Lie(U^+)(u q_1)). 
\end{displaymath}
By (\ref{eq:tilde:P:s}), 
\begin{equation}
\label{eq:Lie:Us}
\Lie(U_s) = \tilde{Q}_s(uq_1; q_1')_* \Lie(U^+)(u q_1). 
\end{equation}
Let 
\begin{displaymath}
M_s = \tilde{Q}_s(uq_1; q_1')_* - I.
\end{displaymath}
Then $(M_s,v_s)$ parametrizes $\cU_s'$. Since $\tilde{Q}_s(uq_1; q_1')
\in Q_{++}(uq_1)$, $M_s \in \cH_{++}(q_1)$. 

Note that by (\ref{eq:P1x:example}), we can recover $\Im q_1$ from
$\gP(q_1)$. Also, since $q_1$ is considered known and fixed here, 
knowing $\Im q_1'$ is equivalent to knowing $q_1'$ since $\Re q_1 =
\Re q_1'$.

Also, since by
Proposition~\ref{prop:sublyapunov:locally:constant} (a), for $q_1' \in
W^-[q_1]$, 
\begin{equation}
\label{eq:Lie:Uplus:q1prime:recover}
\Lie(U^+)(q_1') = P^-(q_1,q_1')_* \Lie(U^+)(q_1) = (\gP(q_1') \circ
\gP(q_1)^{-1})_* 
\Lie(U^+)(q_1), 
\end{equation}
we can reconstruct $U^+(q_1')$ if we know $\gP(q_1)$, $U^+(q_1)$ and
$\gP(q_1')$. Now 
let $\index{$i_{u,q_1,s}$}i_{u,q_1,s}: \cL_{ext}(q_1) \to \cH_{++}(u q_1) \cross W^+( u q_1)$
be the map taking $\gP(q_1') - \gP(q_1)$ to $(M_s, v_s)$. In view of 
(\ref{eq:Lie:Us}), this is a
polynomial map, since $\tilde{Q}_s$ is a polynomial, and 
both $\Im q_1'$ and $\Lie(U^+)(q_1')$ can be
recovered from $\gP(q_1')$ using (\ref{eq:P1x:example}) and 
(\ref{eq:Lie:Uplus:q1prime:recover}). (Note that $q_1$
  is considered fixed here, so knowing $\gP(q_1') - \gP(q_1)$ is
  equivalent to knowing $\gP(q_1')$).

\bold{The maps $(i_{u,q_1,s})_*$ and $\bfi_{u,q_1,s}$.}
For $a \in \natls$, let $\hat{\bfj}^{\tensor a}: \cL_{ext}(x) \to
\cL_{ext}(x)^{\tensor a}$ be the ``diagonal embedding''
\begin{displaymath}
\hat{\bfj}^{\tensor a}(v) = v
\tensor \dots \tensor v, \qquad\text{($a$ times)}
\end{displaymath}
and let $\hat{\bfj}^{\uplus a}$ denote the corresponding map 
$\cL_{ext}(x) \to \cL_{ext}(x)^{\uplus a}$.

Since
$i_{u,q_1,s}: \cL_{ext}(q_1) \to \cH_{++}(uq_1) \cross W^+(u q_1)$ is a
polynomial map, 
by the universal property of the tensor product, there exists $a
> 0$ and a {\em linear} map $\index{$i$@$(i_{u,q_1,s})_*$}(i_{u,q_1,s})_*: \cL_{ext}(q_1)^{\uplus a} \to
\cH_{++}(u q_1) \cross W^+(u q_1)$ such that
\begin{displaymath}
i_{u,q_1,s} = (i_{u,q_1,s})_* \circ \hat{\bfj}^{\uplus a}.
\end{displaymath}
Furthermore, there exists $r > a$ and a linear map $\index{$i$@$\bfi_{u,q_1,s}$}\bfi_{u,q_1,s}:
\cL_{ext}(q_1)^{\uplus r} \to \tilde{\bH}(u q_1)$ such that
\begin{equation}
\label{eq:sufficiently:large:r}
\bfj \circ (i_{u,q_1,s})_* = \bfi_{u,q_1,s} \circ \hat{\bfj}^{\uplus
  r},
\end{equation}
where $\bfj$ is as in (\ref{eq:tmp:def:bfj}). 
Then $\bfi_{u,q_1,s}$ takes $F(q_1') - F(q_1) \in
\cL_{ext}(q_1)^{\uplus r}$ to $\bfj(M_s,v_s) \in \tilde{\bH}(u
q_1)$, where $(M_s,v_s)$
is a parametrization of the approximation $\tilde{P}_s(uq_1;q_1')^{-1}
U^+[q_1']$ to $U^+[q_1']$. 

\bold{Construction of the map $\cA(q_1,u,\ell,t)$.}
Let $s \in \natls$ be a sufficiently large integer to be chosen
later. (It will be chosen near the end of the proof of 
Proposition~\ref{prop:reason:cA:F}, depending
only on the Lyapunov spectrum). 
Let $r \in \natls$ be such that
(\ref{eq:sufficiently:large:r}) holds. 
Suppose $q_1 \in X_0$ and $u \in \cB(q_1,1/100)$. For $\ell > 0$ and
$t > 0$, let
\begin{displaymath}
\cA(q_1,u,\ell,t): \cL_{ext}(g_{-\ell} q_1)^{(r)} \to \bH(g_t u q_1),
\end{displaymath}
be given by 
\begin{displaymath}
\cA(q_1,u,\ell,t) = (g_t)_* \circ \bS_{uq_1}^{Z(u q_1)} \circ \hat{\pi}  \circ \bfi_{u,q_1,s} \circ (g_\ell)_*^{\uplus r}
\end{displaymath}
where $(g_\ell)_* : \cL_{ext}(q) \to \cL_{ext}(g_\ell q)$ is given by 
\begin{displaymath}
(g_\ell)_*(P) = g_\ell \circ P \circ g_\ell^{-1}.
\end{displaymath}
Then $\cA(q_1,u,\ell,t)$ is a {\em linear} map. Unraveling the
definitions, we have, for $P \in \cL_{ext}(g_{-\ell}q_1)$,
\begin{displaymath}
\cA(q_1,u,\ell,t)(\hat{\bfj}^{\uplus r}(P)) = \bfj(G_t^+ \circ
S^{Z(uq_1)}_{uq_1}  \circ (i_{u,q_1,s}) \circ (g_\ell)_* (P))
\end{displaymath}
\mcc{(maybe more details here?)} 
Thus, for $q'$ satisfying the upper bounds in
(\ref{eq:rho:prime:delta:le:d:q:qprime})
and (\ref{eq:rho:delta:le:Fq:minus:Fqprime}), 
\begin{equation}
\label{eq:cA:better:unravel}
\cA(q_1,u,\ell,t)(F(q)-F(q')) = \bfj(M'',v''),
\end{equation}
where $(M'',v'') \in \cH_{++}(g_t u q_1) \cross W^+(uq_1)$ is a
parametrization of the approximation $$g_t \tilde{P}_s(uq_1,
u'q_1')^{-1} U^+[u' q_1']$$ to $U^+[g_t u' q_1']$, 
where $u' q_1' \in U^+[q_1']$ is such that 
$d^{X_0}(g_t uq_1, g_t u' q_1') < 1/100$.


\starredsubsection{Proofs of Proposition~\ref{prop:reason:cA:F} and
  Lemma~\ref{lemma:parts:staying:close}.}
\label{sec:starredsubsec:proofs:reason:cA:F}

$ $

\bold{Proof of Proposition~\ref{prop:reason:cA:F}.} Note that
Proposition~\ref{prop:reason:cA:F} (i) follows immediately from the
definition of $\cA(\cdot, \cdot, \cdot, \cdot)$. We now begin the
proof of Proposition~\ref{prop:reason:cA:F} (iii). 
Let $P = \gP(q')-\gP(q) \in \cL_{ext}(q)$. 
Let
\begin{displaymath}
P_1 = (g_\ell)_*(P) = g_\ell \circ P \circ g_{\ell}^{-1} \in
\cL_{ext}(q_1). 
\end{displaymath}
Let
\begin{displaymath}
(M_s,v_s) = i_{u,q_1,s}(P_1). 
\end{displaymath}
Let $\widetilde{\cU}_s = \widetilde{\cU}_s(M_s,v_s)$ 
be the generalized subspace parametrized by $(M_s,v_s)$. Then 
\begin{equation}
\label{eq:def:tilde:cU:s}
\widetilde{\cU}_s = \tilde{P}_s(uq_1,q_1')^{-1} U^+[q_1'].
\end{equation}
Let 
\begin{equation}
\label{eq:def:hat:cU:hat:cU:s}
\widehat{\cU} = \hat{P}(uq_1, q_1')^{-1} U^+[q_1'], \qquad
\widehat{\cU}_s = \hat{P}_s(u q_1, q_1')^{-1} U^+[q_1]. 
\end{equation}
Suppose (\ref{eq:reason:cA:F:point:bound}) holds. By Lemma~\ref{lemma:reduce:to:P:hat}, 
\begin{equation}
\label{eq:hd:widehat:cU:Uplus:q1prime}
hd_{g_t u q_1}^{X_0}( \widehat{\cU}, U^+[g_t u' q_1']) = O_{u q_1}(e^{-\alpha_1 t}),
\end{equation}
where $\alpha_1$ depends only on the Lyapunov spectrum. We have, in
view of (\ref{eq:rho:prime:delta:le:d:q:qprime}) and 
(\ref{eq:rho:delta:le:Fq:minus:Fqprime}), for $\ell$ sufficiently
large depending on $\delta$, 
\begin{equation}
\label{eq:estimate:Pminus:q1:q1prime}
\|P^-(q_1,q_1')P^{GM}(q_1',q_1) - I \|_Y = O_{q_1}(e^{-\alpha_2 \ell})
\end{equation}
where $\alpha_2$ depends only on the Lyapunov spectrum. Therefore,
\begin{displaymath}
hd_{uq_1}^{X_0}(U^+[u q_1], U^+[q_1']) =  O_{q_1}(e^{-\alpha_2 \ell})
\end{displaymath}

To go from $\hat{Q}$ to $\hat{Q}_s$ we are doing order $s$ Taylor expansion of
the solution to (\ref{eq:LU})
in the entries of $P^-(q_1,q_1') P^{GM}(q_1',q_1) - I$. Thus, 
by (\ref{eq:estimate:Pminus:q1:q1prime}),  
\begin{displaymath} 
\|\hat{Q}_s(u q_1;q_1') - \hat{Q}(u q_1; q_1')\|_Y =
O_{q_1,uq_1}(e^{-\alpha_2 (s+1) \ell}) 
\end{displaymath}
and thus, by  (\ref{eq:tilde:P:s}),
\begin{equation}
\label{eq:hat:Ps:close:to:hatP}
\|\hat{P}_s(uq_1, q_1')^{-1} - \hat{P}(uq_1,q_1')^{-1}\|_Y = O_{q_1,uq_1}(e^{-\alpha_2 (s+1) \ell}) 
\end{equation}
Then, by (\ref{eq:def:hat:cU:hat:cU:s}), 
\begin{displaymath}
hd_{u q_1}^{X_0} ( \widehat{\cU}, \widehat{\cU}_s) = O_{q_1,uq_1}(e^{-\alpha_2(s+1)
  \ell}). 
\end{displaymath}
Then, by Lemma~\ref{lemma:crude:divergence:subspaces}(a), 
\begin{equation}
\label{eq:hd:hatU:hatUs}
hd_{g_t u q_1}^{X_0}( g_t \widehat{\cU}, g_t \widehat{\cU}_s) =
O_{q_1,uq_1}(e^{-\alpha_2(s+1)\ell + 2t}).  
\end{equation}
Also, by (\ref{eq:estimate:Pminus:q1:q1prime}), (\ref{eq:def:hat:Q})
and (\ref{eq:hatP:minus:I:small}),  
we have
\begin{displaymath}
\|\hat{Q}(uq_1; q_1') - I \|_Y = O_{q_1, uq_1}(e^{-\alpha_2 \ell}), 
\end{displaymath}
and therefore 
\begin{displaymath}
\|\hat{Q}_s(uq_1; q_1') - I \|_Y = O_{q_1, uq_1}(e^{-\alpha_2 \ell}), 
\end{displaymath}
Thus, 
\begin{displaymath}
\|D_s\|_Y = \|\tilde{Q}_s(u q_1;q_1') - \hat{Q}_s(u q_1; q_1')\|_Y = 
O_{q_1}(e^{-\alpha_2 \ell})
\end{displaymath}
Therefore, since $D_s$ preserves all the eigenspaces $\cV_i$, 
and the Osceledets multiplicative ergodic theorem, for 
sufficiently small $\epsilon > 0$ (depending on the Lyapunov
spectrum), 
\begin{displaymath}
\| g_t \circ D_s \circ g_t^{-1} \|_Y \le C_1(q_1) C_2(u q_1, \epsilon)e^{-\alpha_2 \ell +
  \epsilon t} \le C_1(q_1) C_2'(u q_1)e^{-(\alpha_2/2) \ell}. 
\end{displaymath}
Thus, 
\begin{equation}
\label{eq:tilde:Ps:close:to:hat:Ps}
\|\tilde{P}_s(g_t u q_1, g_t u' q_1')^{-1} - \hat{P}_s(g_t u q_1, g_t u'
  q_1')^{-1} \| = O_{u q_1}(e^{-(\alpha_2/2)\ell})
\end{equation}
and hence by (\ref{eq:def:tilde:cU:s}) and (\ref{eq:def:hat:cU:hat:cU:s}),
\begin{equation}
\label{eq:hd:widehat:cUs:widetilde:cUs}
hd_{g_t u q_1}^{X_0} (g_t \widehat{\cU}_s, g_t \widetilde{\cU}_s) = O_{u q_1}(\| g_t
\circ D_s \circ g_t^{-1}\|_Y) = O_{u q_1}(e^{-(\alpha_2/2)\ell}). 
\end{equation}
We now choose $s$ so that $\alpha_2 \alpha_3 (s+1)
-3 > \alpha_2$. Then, by (\ref{eq:reason:cA:F:tbound}), 
(\ref{eq:hd:widehat:cU:Uplus:q1prime}), 
(\ref{eq:hd:hatU:hatUs}), and
(\ref{eq:hd:widehat:cUs:widetilde:cUs}),
\begin{equation}
\label{eq:final:in:proof:reason:cA:F}
hd_{g_t u q_1}^{X_0}( g_t \widetilde{\cU}_s, U^+[g_t q_1']) \le C(q_1)
C(u q_1) e^{-\alpha \ell}, 
\end{equation}
where $\alpha$ depends only on the Lyapunov spectrum. 
In view of (\ref{eq:cA:better:unravel}), the pair $(M'', v'')$
parametrizes $g_t 
\widetilde{\cU}_s$. Therefore, (\ref{eq:reason:cA:F:unravel})
holds.  Finally, (\ref{eq:reason:cA:F:estimate}) 
is an immediate consequence of (\ref{eq:final:in:proof:reason:cA:F}). 
This completes the proof of
Proposition~\ref{prop:reason:cA:F} (iii). (Note that is was shown
immediately after the statement of Proposition~\ref{prop:reason:cA:F}
that Proposition~\ref{prop:reason:cA:F} (iii) implies
Proposition~\ref{prop:reason:cA:F} (ii).) 
\qed\medskip


\bold{Proof of Lemma~\ref{lemma:parts:staying:close}.}In the proof of
this lemma we normalize the measure $|\cdot|$ on $U^+[q_1]$ so
that $|U^+[q_1] \cap B^+(q_1,1/100)| = 1$ and similarly we normalize the measure
$|\cdot|$ on $U^+[q_1']$ so that $|U^+[q_1'] \cap B^+(q_1',1/100)|=1$.
As in the proof of Lemma~\ref{lemma:reduce:to:P:hat}, we 
choose $u' \in \cB(q_1',1/50)$ with $\cV_i(g_t u' q_1')$ and
$U^+[g_t u' q_1'] = U^+[g_t q_1']$ defined and 
\begin{displaymath}
d^{X_0}(g_t u q_1, g_t u' q_1') \le hd_{g_t u q_1}^{X_0}(U^+[g_t u
q_1], U^+[g_t q_1']) \le \epsilon.
\end{displaymath}
(Nothing in the proof will depend on the choice of $u'$). 

Let $A_0 = g_{-t}A_t$, $A_0' = g_{-t}A_t'$. 
Let $\tilde{P}_s$ be as in (\ref{eq:tilde:P:s}). 
Let $\tilde{A}_t = \tilde{P}_s(g_t u q_1, g_t u' q_1')^{-1} A'_t$. Then,
\begin{displaymath}
\tilde{A}_0 \equiv g_{-t}\tilde{A}_t = \tilde{P}_s(u q_1, u' q_1')^{-1} A'_0.
\end{displaymath}
As in the proof of Proposition~\ref{prop:reason:cA:F}, (i.e.\
by combining (\ref{eq:hatP:minus:I:small}),
(\ref{eq:hat:Ps:close:to:hatP}) and (\ref{eq:tilde:Ps:close:to:hat:Ps})),
we have
\begin{displaymath}
\|\tilde{P}_s(u q_1, u' q_1')^{-1} P^{GM}(u q_1, u' q_1')
- I\|_Y = O(e^{-\alpha \ell}).
\end{displaymath}
\begin{displaymath}
\|\tilde{P}_s(g_t u q_1, g_t u' q_1')^{-1} P^{GM}(g_t u q_1, g_t' u' q_1')
- I\|_Y = O(e^{-\alpha \ell}).
\end{displaymath}
Hence, $|\tilde{A}_t|$ is comparable to $|A'_t|$ and $|\tilde{A}_0|$
is comparable to $|A_0'|$. 
Thus, it is enough to show that $|\tilde{A}_0|$ is comparable to $|A_0|$.

As in the proof of Proposition~\ref{prop:reason:cA:F}, let $(M'',v'')$ be
the pair parametrizing 
$g_t\tilde{\cU}_s = \tilde{P}_s( g_t u q_1, g_t u' q_1)^{-1} U^+[g_t
u'q_1']$. Let $\tilde{f}_t: \Lie(U^+)(g_t u q_1) \to g_t\tilde{\cU}_s$ be the
``parametrization'' map
\begin{displaymath}
\tilde{f}_t(Y) = \exp[(I+M'')Y](g_t u q_1)(g_t u q_1 + v''). 
\end{displaymath}
Similarly, let $f_t: \Lie(U^+)(g_t u q_1) \to U^+[g_t u q_1]$ be the
exponential map
\begin{displaymath}
f_t(Y) = \exp(Y)g_t u q_1. 
\end{displaymath}
Then, provided that $\epsilon$ is sufficiently small, we have
\begin{equation}
\label{eq:triple:inclusion}
0.5 f^{-1}(A_t) \subset \tilde{f}_t^{-1}(\tilde{A}_t) \subset 2 f^{-1}(A_t)
\end{equation}
Let $M_0 = g_t^{-1} \circ M'' \circ g_{t}$, $v_0 = g_t^{-1}
v''$. Then, $g_t^{-1} \circ \tilde{f}_t \circ g_t = \tilde{f}_0$,
where $\tilde{f}_0: \Lie(U^+)(uq_1) \to \cU_s$ is given by
\begin{displaymath}
\tilde{f}_0(Y) = \exp[(I+M_0)Y](g_t u q_1)(g_t u q_1 + v_0). 
\end{displaymath}
Similarly, $g_t^{-1} \circ f_t \circ g_t = f_0$, where $f_0:
\Lie(U^+)(u q_1) \to U^+[u q_1]$ is given by the exponential map
\begin{displaymath}
f_0(Y) = \exp(Y) u q_1. 
\end{displaymath}
Then, it follows from applying $g_t^{-1}$ to (\ref{eq:triple:inclusion}) that
\begin{equation}
\label{eq:new:triple:inclusion}
0.5 f_0^{-1}(A_0) \subset \tilde{f}_0^{-1}(\tilde{A}_0) \subset 2 f_0^{-1}(A_0)
\end{equation}
Thus, $|\tilde{f}_0^{-1}(\tilde{A}_0)|$ is comparable to
$|f_0^{-1}(A_0)| = |A_0|$. But, since $M'' \in \cH_{++}(g_t u q_1)$ and
$v'' \in W^+(g_t u q_1)$ are $O(\epsilon)$, $M_0$ and $v_0$ are
exponentially small. Therefore, the map $\tilde{f}_0$ is close to
$f_0$ (and since $Y$ is small, it is close to the
identity). Therefore, $|\tilde{f}_0^{-1}(\tilde{A}_0)|$ is comparable
to $|\tilde{A_0}|$. The second assertion of the Lemma also follows
from (\ref{eq:new:triple:inclusion}) and the fact that $M_0$ and $v_0$
are exponentially small. \mcc{this proof needs to be more formal}
\qed\medskip

\section{Bilipshitz estimates}
\label{sec:new:bilipshitz}

In this section, we continue working on $X_0$ (and not $X$). 
Let $\|\cdot\|$ be the norm on $H_{big}^{(++)}$ defined in
(\ref{eq:def:dynamical:norm:X0}). Since $\bH \subset
H_{big}^{(++)}$, $\| \cdot \|$ is also a norm on $\bH$. 
We can also define a norm on $H_{big}^{(--)}$ in an analogous way. 
Since
$\cL_{ext}(x)^{(r)} \subset H_{big}^{(--)}(x)$, the norm $\| \cdot \|_x$
is also a norm on $\cL_{ext}(x)^{(r)}$. 
Let $\index{$A(q_1,u,\ell,t)$}A(q_1,u,\ell,t) =
\|\cA(q_1,u,\ell,t)\|$ where the operator norm is with respect to the
dynamical norms $\| \cdot \|$ at $g_{-\ell} q_1$ and $g_t u q_1$. 
In the rest of this section we assume that the equivalent conditions
of Lemma~\ref{lemma:Lplus:Splus:equiv} do not hold, and then by
Proposition~\ref{prop:stop:induction:condition},
(\ref{eq:stop:condition:estimate}) holds.

For $1/100 > \epsilon > 0$, almost all $q_1 \in X_0$, almost all $u \in \cB(q_1,1/100)$ and $\ell > 0$, let 
\begin{displaymath}
\index{$\tau$@$\hat{\tau}_{(\epsilon)}(q_1,u, \ell)$}\hat{\tau}_{(\epsilon)}(q_1,u, \ell) = \sup \{t \st t > 0 \text{ and }
A(q_1,u,\ell,t) \le \epsilon  \}.
\end{displaymath}
Note that $\hat{\tau}_{(\epsilon)}(q_1,u,0)$ need not be $0$.

For $x \in X_0$, 
let $\index{$A$@$\cA_+(x,t)$}\cA_+(x,t): \bH(x) \to \bH(g_t x)$ denote the
action of $g_t$ on $\bH$ as in (\ref{eq:def:action:on:bHplus}).  
Let $\index{$A$@$\cA_-(x,s)$}\cA_-(x,s): \cL_{ext}^{(r)}(x) \to \cL_{ext}^{(r)}(g_s x)$ denote
the action of $g_s$ on $\cL_{ext}^{(r)}(x)$. 
\begin{lemma}
\label{lemma:forni2}
There exist absolute constants $N > 0$, $\alpha> 0$ such that
for almost all $x$, and $t > 0$, 
\begin{displaymath}
e^{-\alpha t} \ge \|\cA_-(x,t)\| \ge e^{-Nt}, \qquad  e^{\alpha t} \le
\|\cA_+(x,t)\| \le e^{Nt}. 
\end{displaymath}
and,
\begin{displaymath}
e^{ N t} \ge \|\cA_-(x,-t)\| \ge e^{ \alpha t}, \qquad  e^{-N t} \le
\|\cA_+(x,-t)\| \le e^{-\alpha t}. 
\end{displaymath}
\end{lemma}
\medskip
\bold{Proof.} This follows immediately from
Proposition~\ref{prop:properties:dynamical:norm:Hbig}.
\qed\medskip

\begin{lemma}
\label{lemma:monotonicity1}
Suppose $0 < \epsilon  < 1/100$. 
There exists $\kappa_1 > 1$ (depending only on the
Lyapunov spectrum) with the following property: for almost all $q_1
\in X_0$, $u \in \cB(q_1,1/100)$, for all $\ell > 0$ and $s > 0$, 
\begin{displaymath}
\hat{\tau}_{(\epsilon)}(q_1,u, \ell + s) > \hat{\tau}_{(\epsilon)}(q_1,u,\ell) + \kappa_1^{-1} s. 
\end{displaymath}
\end{lemma}

\bold{Proof.} Note that by (\ref{eq:reason:cA:F:equivariance}), 
\begin{displaymath}
\cA(q_1,u,\ell+s,t+\tau) = \cA_+(g_t u q_1,\tau) \cA(q_1,u,\ell,t)
\cA_-(g_{-(\ell+s)} q_1, s). 
\end{displaymath}
Let $t = \hat{\tau}_{(\epsilon)}(q_1,u,\ell)$,  so that
$A(q_1,u,\ell,t) = \epsilon$.  
Therefore, 
\begin{multline*}
A(q_1,u,\ell+s,t+\tau) \le \|\cA_+(g_t u q_1,\tau)\|
A(q_1,u,\ell,t) \|\cA_-(g_{-(\ell+s)}q_1,s)\| \le \\
\epsilon \|\cA_+(q_t u q_1,\tau)\| \|\cA_-(g_{-(\ell+s)}q_1,s)\| \le
\epsilon e^{N\tau-\alpha s}, 
\end{multline*}
where we have used the fact that $A(q_1,u,\ell,t) = \epsilon$ and
Lemma~\ref{lemma:forni2}. 
If $t+\tau = \hat{\tau}_{(\epsilon)}(q_1,u,\ell+s)$
then $A(q_1,u,\ell+s,t+\tau) = \epsilon$. It follows that $N\tau - \alpha s \ge
0$, i.e.\ $\tau \ge (\alpha/N) s$. Hence,
\begin{displaymath}
\hat{\tau}_{(\epsilon)}(q_1,u,\ell+s) \ge \hat{\tau}_{(\epsilon)}(q_1, u, \ell) + (\alpha/N) s. 
\end{displaymath}
\qed\medskip

\begin{lemma}
\label{lemma:monotonicity2}
Suppose $0 < \epsilon < 1/100$. 
There exists $\kappa_2 > 1$ (depending only on the Lyapunov spectrum) 
such that for almost all $q_1 \in X_0$, almost all $u \in \cB(q_1,1/100)$,
all $\ell > 0$ and all $s > 0$, 
\begin{displaymath}
\hat{\tau}_{(\epsilon)}(q_1,u, \ell + s) <
\hat{\tau}_{(\epsilon)}(q_1,u,\ell) + \kappa_2 s. 
\end{displaymath}
\end{lemma}

\bold{Proof.}
We have
\begin{displaymath}
\cA(q_1,u,\ell,t) = \cA_+(g_{t+\tau} u q_1,-\tau) \cA(q_1,u,\ell+s,t+\tau)
\cA_-(g_{-\ell} q_1, -s). 
\end{displaymath}
Let $t+\tau =
\hat{\tau}_{(\epsilon)}(q_1,u,\ell+s)$. 
Then, by Lemma~\ref{lemma:forni2},
\begin{multline*}
A(q_1,u,\ell,t) \le \|\cA_+(q_{t+\tau} u q_1,-\tau)\|
A(q_1,u,\ell+s,t+\tau) \|\cA_-(g_{-\ell}q_1,-s)\| \le \\
\epsilon \|\cA_+(q_{t+\tau} u q_1,-\tau)\| \|\cA_-(g_{-\ell}q_1,-s)\| \le 
\epsilon e^{-\alpha\tau+Ns}, 
\end{multline*}
where we have used the fact that $A(q_1,u,\ell+s,t+\tau) = \epsilon$. 
Since $A(q_1,u,\ell,t)=\epsilon$, it follows that $-\alpha\tau + N s >
0$, i.e.\ $\tau < (N/\alpha) s$. It follows that
\begin{displaymath}
\hat{\tau}_{(\epsilon)}(q_1,u, \ell + s) <
\hat{\tau}_{(\epsilon)}(q_1,u,\ell) + (N/\alpha) s
\end{displaymath}
\qed\medskip

\begin{proposition}
\label{prop:bilip:hattau:epsilon}
There exists $\kappa > 1$ depending only on the Lyapunov spectrum, and
such that for almost all $q_1  \in X_0$, almost all $u \in
\cB(q_1,1/100)$, any $\ell > 0$ and any measurable 
subset $E_{bad} \subset \reals^+$, 
\begin{align*}
|\hat{\tau}_{(\epsilon)}(q_1,u,E_{bad}) \cap [\hat{\tau}_{(\epsilon)}(q_1,u,0),\hat{\tau}_{(\epsilon)}(q_1,u,\ell)]| & \le \kappa
|E_{bad} \cap [0,\ell]| \\
|\{ t \in [0,\ell]| \st \hat{\tau}_{(\epsilon)}(q_1,u,t) \in E_{bad}
\}| & \le \kappa |E_{bad} \cap
      [\hat{\tau}_{(\epsilon)}(q_1,u,0),\hat{\tau}_{(\epsilon)}(q_1,u,\ell)]|.
\end{align*}
\end{proposition}

\bold{Proof.} 
Let $\kappa = \max(\kappa_1^{-1},\kappa_2)$, where $\kappa_1$,
$\kappa_2$ are as in Lemma~\ref{lemma:monotonicity1} and
Lemma~\ref{lemma:monotonicity2}. Then, for fixed $q_1$, $u$,
$\hat{\tau}_{(\epsilon)}(q_1,u,\ell)$ is $\kappa$-bilipshitz as a
function of $\ell$. The proposition follows immediately.  
\qed\medskip

\section{Preliminary divergence estimates}
\label{sec:first:divergence}

In this section, we continue working on $X_0$ (and not $X$). 

\bold{Motivation.} Suppose in the notation of
\S\ref{sec:outline:step1}, $q_1$ and $q_1'$ are fixed, but $u \in \cB(q_1,1/100)$
and $u' \in \cB(q_1',1/100)$ 
vary. Then, as $u$ and $u'$ vary, so do the points $q_2$ and $q_2'$, and
thus the subspaces $U^+[q_2]$ and $U^+[q_2']$. Let $\cU =
\cU(M''(u),v''(u))$ be the 
approximation to $U^+[q_2']$ given by Proposition~\ref{prop:reason:cA:F},
and as in Proposition~\ref{prop:reason:cA:F}, 
let $\bfv(u) = \bfj(M''(u),v''(u)) \in \bH(q_2)$ be the associated vector
in $\bH(q_2)$. 

In this section we define a certain $g_t$-equivariant and
$(u)_*$-equivariant subbundle $\index{$E$@$\bE$}\bE \subset \bH$
such that, for fixed $q_1$, $q_1'$, for most $u \in U^+[q_1]$, 
$\bfv = \bfv(u)$ is near $\bE(q_2)$ (see
Proposition~\ref{prop:most:inert} (a) below for the precise
statement).  We call $\bE$ the $U^+$-inert subbundle of $\bH$.
The subbundle $\bE$ is the direct sum of subbundles
$\index{$E$@$\bE_j(x)$}\bE_i$, where $\bE_i$ is contained in the $i$-th Lyapunov subspace of
$\bH$, and also each $\bE_i$ is both $g_t$-equivariant and
$(u)_*$-equivariant.

\subsection{The $U^+$-inert subspaces $\bE(x)$}


We apply the Osceledets multiplicative
ergodic theorem to the action on $\bH(x)$ (see
(\ref{eq:def:action:on:bHplus})). We often drop the $*$ and denote the
action simply by $g_t$. In this section, $\lambda_i$ denotes the $i$-th
Lyapunov exponent of the flow $g_t$ on the bundle $\bH$. 

Let
\begin{displaymath}
\index{$V$@$\bV_{\le i}(x)$}\bV_{\le i}(x) = \bigoplus_{j \le i}
\cV_j(\bH)(x), \qquad 
\index{$V$@$\bV_{< i}(x)$}\bV_{< i}(x) = \bigoplus_{j < i} \cV_j(\bH)(x),
\end{displaymath}
\begin{displaymath}
\index{$V$@$\bV_{\geq i}(x)$}\bV_{\geq i}(x) = \bigoplus_{j \geq i} \cV_j(\bH)(x),  \qquad \index{$V$@$\bV_{> i}(x)$}\bV_{> i}(x) = \bigoplus_{j > i} \cV_j(\bH)(x).  
\end{displaymath}
This means that for almost all $x \in X_0$ and 
for $\bfv \in \bV_{\le i}(x)$ such that $\bfv \not\in \bV_{< i}(x)$, 
\begin{equation}
\label{eq:bVi:growth:exactly:lambda:i}
\lim_{t \to -\infty} \frac{1}{t} \log \frac{\|g_t \bfv\|}{\|\bfv\|}
= \lambda_i,
\end{equation}
and for $\bfv \in \bV_{\geq i}(x)$ such that $\bfv \not\in \bV_{> i}(x)$, 
\begin{equation}
\label{eq:growth:hat:V:n:plus:one:minus:i}
\lim_{t \to \infty} \frac{1}{t} \log \frac{\|g_t \bfv\|}{\|\bfv\|} =
\lambda_{i}.
\end{equation}
By e.g.\ \cite[Lemma 1.5]{Goldsheid:Margulis}, we have for a.e.\ $x\in X_0$, 
\begin{equation}
\label{eq:Lyapunov:transversality}
\bH(x) = \bV_{\le i}(x) \dirsum \bV_{> i}(x).
\end{equation}

Let 
\begin{equation}
\label{eq:def:bFj}
\index{$F$@$\bF_{\ge j}(x)$}\bF_{\ge j}(x) = \{ \bfv \in \bH(x) \st \text{ for almost 
all $u \in \cB(x)$, $(u)_* \bfv \in \bV_{\geq j}(u x)$} \}, 
\end{equation}
where $(u)_*$ is as in Lemma~\ref{lemma:tensor:product}. 
In other words, if $\bfv \in \bF_{\ge j}(x)$, then  
for almost all $u \in \cB(x)$,
\begin{equation}
\label{eq:meaning:bFj}
\limsup_{t \to \infty} \frac{1}{t} \log \|(g_t)_* (u)_* \bfv\| \le \lambda_j.
\end{equation}
From the definition of $\bF_{\ge j}(x)$, we have
\begin{equation}
\label{eq:flag:bFj}
\{0\} = \bF_{\ge n+1}(x) \subset \bF_n(x) \subset \bF_{\ge n-1}(x) \subset
\dots \bF_2(x) \subset \bF_1(x) = \bH(x).
\end{equation}
Let 
\begin{displaymath}
\index{$E$@$\bE_j(x)$}\bE_j(x) = \bF_{\ge j}(x) \cap \bV_{\le j}(x). 
\end{displaymath}
In particular, $\bE_1(x) = \bV_{\le 1}(x) = \cV_1(\bH)(x)$.  We may have $\bE_j(x) = \{0\}$
if $j \ne 1$. 
\begin{lemma}
\label{lemma:grow:exactly:lambdaj}
For almost all $x \in X_0$ the following holds: suppose $\bfv \in
\bE_j(x)\setminus \{0\}$. Then for almost all $u \in \cB(x)$, 
\begin{equation}
\label{eq:Ej:grows:exactly:lambda:j}
\lim_{t \to \infty} \frac{1}{t} \log \|(g_t)_* (u)_* (\bfv)\| = \lambda_j.
\end{equation}
Thus (recalling that $\cV_j(\bH)$ denotes the
  subspace of $\bH$ 
corresponding to the Lyapunov exponent $\lambda_j$), we have
for almost all $x$, using Fubini's theorem,
\begin{displaymath}
\bE_j(x) \subset \cV_j(\bH)(x). 
\end{displaymath}
\mccc{(For this, we need to apply
  (\ref{eq:Ej:grows:exactly:lambda:j}) for $u = 0$, which needs extra Fubini).} 
In particular, if $i \ne j$, $\bE_i(x) \cap \bE_j(x) = \{0\}$ for
almost all $x \in X_0$. 
\end{lemma}
\bold{Proof.} Suppose $\bfv \in \bE_j(x)$. Then $\bfv \in \bV_{\le j}(x)$. Since
in view of (\ref{eq:bVi:growth:exactly:lambda:i}),
$\bV_{\le j}(u x) = (u)_* \bV_{\le j}(x)$ for all $u \in U^+(x)$, 
we have for almost all $u \in \cB(x)$, $(u)_* \bfv \in \bV_{\le j}(ux)$. It
follows from (\ref{eq:Lyapunov:transversality}) that (outside of a set
of measure $0$), $(u)_* \bfv \not\in \bV_{> j}(ux)$. Now
(\ref{eq:Ej:grows:exactly:lambda:j}) follows from
(\ref{eq:growth:hat:V:n:plus:one:minus:i}). 
\qed\medskip

\begin{lemma}
\label{lemma:Ej:equivariant}
After possibly modifying $\bE_j(x)$ and $\bF_{\ge j}(x)$ 
on a subset of measure $0$ of $X$,
the following hold:
\begin{itemize}
\item[{\rm (a)}] $\bE_j(x)$ and $\bF_{\ge j}(x)$ are $g_t$-equivariant, i.e.\ $(g_t)_*
  \bE_j(x) = \bE_j(g_t x)$, and $(g_t)_* \bF_{\ge j}(x) = \bF_{\ge j}(g_t x)$.
\item[{\rm (b)}] For almost all $u \in U^+(x)$, $\bE_j(ux) =
  (u)_* \bE_j(x)$, and $\bF_{\ge j}(ux) = (u)_* \bF_{\ge j}(x)$. 
\end{itemize}
\end{lemma}
\bold{Proof.} Note that for $t > 0$, 
$g_t \cB[x] \supset \cB[g_t x]$. Therefore, (a)
for the case $t > 0$ follows immediately from the definitions of
$\bE_j(x)$ and $\bF_{\ge j}(x)$.  
Since the flow $\{g_t\}_{t > 0}$ is ergodic, it follows that almost
everywhere (\ref{eq:def:bFj}) holds with $\cB[x]$ replaced by arbitrary
large balls in $U^+[x]$. This implies that almost everywhere, 
\begin{displaymath}
\bF_{\ge j}(x) = \{ \bfv \in \bH(x) \st \text{ for almost 
all $u \in U^+$, $(u)_* \bfv \in \bV_{\geq j}(u x)$} \},
\end{displaymath}
where $(u)_* \bfv$ is as in Lemma~\ref{lemma:tensor:product}. 
Therefore (b) holds. Then, (a) for $t < 0$ also holds, as long as both
$x$ and $g_t x$ belong to a subset of full measure. By considering a
transversal for the flow $g_t$, it is easy to
check that it is possible to modify  $\bE_j(x)$ and $\bF_{\ge j}(x)$ 
on a subset of measure $0$ of $X_0$ in such a way that (a) holds for $x$
in a subset of full measure and all $t \in \reals$.
\mcc{more details?}
\qed\medskip

\begin{lemma}
\label{lemma:zero:one:substitute}
For $x \in X_0$, let 
\begin{displaymath}
Q(\bfv) = \{ u \in \cB(x) \st (u)_* \bfv \in  \bV_{\geq j}(u x) \}. 
\end{displaymath}
Then for almost all $x$, 
either $|Q(\bfv)| = 0$, or $|Q(\bfv)| = |\cB(x)|$ (and thus
$\bfv \in \bF_{\ge j}(x)$). 
\end{lemma}
\bold{Proof.} For a subspace $\bV \subset \bH(x)$, let 
\begin{displaymath}
Q(\bV) = \{ u \in \cB(x) \st (u)_* \bV  \subset  \bV_{\geq j}(u x) \}. 
\end{displaymath}
Let $d$ be the maximal number such that there exists 
$E'\subset X_0$ with $\nu(E') > 0$ such that for $x \in E'$ there exists a
subspace $\bV \subset \bH(x)$ 
of dimension $d$ with $|Q(\bV)| > 0$. 
For a fixed $x \in E'$, let
$\cW(x)$ denote the set of subspaces $\bV$ of dimension $d$ 
for which $|Q(\bV)| > 0$. Then, by
the maximality of $d$, if $\bV$ and $\bV'$ are distinct elements of
$\cW(x)$ then $Q(\bV) \cap Q(\bV')$ has measure $0$. Let $\bV_x \in \cW(x)$ be
such that $|Q(\bV_x)|$ is maximal (among elements of $\cW(x)$).

Let $\epsilon > 0$ be arbitrary, and
suppose $x \in E'$. By the same Vitali-type argument as in the proof
of Lemma~\ref{lemma:gB:vitali}, \mcc{more details?}
there exists $t_0 >0$ and 
a subset $Q(\bV_x)^* \subset Q(\bV_x) \subset \cB(x)$ 
such that for all $u \in Q(\bV_x)^*$ and
all $t > t_0$, 
\begin{equation}
\label{eq:ux:point:of:density}
|\cB_t( u x ) \cap Q(\bV_x)| \ge (1-\epsilon) |\cB_t(ux)|.   
\end{equation}
(In other words, $Q(\bV_x)^*$ are ``points of density'' for $Q(\bV_x)$,
relative to the ``balls'' $\cB_t$.) 
Let 
\begin{displaymath}
E^* = \{ u x \st x \in E', \quad u \in Q(\bV_x)^* \}. 
\end{displaymath}
Then, $\nu(E^*) > 0$. Let $\Omega = \{ x \in X_0 \st g_{-t} x \in E^*
\text{ for an unbounded set of $t > 0$ }\}$. Then $\nu(\Omega) = 1$. Suppose $x \in
\Omega$. We can choose $t >t_0$ such that $g_{-t} x \in
E^*$. Note that
\begin{equation}
\label{eq:relation:balls}
\cB[x] = g_t \cB_t[g_{-t}x].
\end{equation}
\mcc{Better to work with $\cB[x]$
  instead of $\cB(x)$}
Let $x' = g_{-t}x$, and let $\bV_{t,x} = (g_t)_*
\bV_{x'}$. Then in view of (\ref{eq:ux:point:of:density}) and
(\ref{eq:relation:balls}),  
\begin{displaymath}
|Q(\bV_{t,x})| \ge (1-\epsilon)|\cB(x)|. 
\end{displaymath}
By the maximality of $d$ (and assuming $\epsilon <
  1/2$), $\bV_{t,x}$ does not depend on $t$. 
Hence, for every $x \in \Omega$, there exists $\bV \subset \bH(x)$ such that
$\dim \bV = d$ and $|Q(\bV)| \ge (1-\epsilon) |\cB(x)|$. Since $\epsilon
> 0$ is arbitrary, for each $x \in \Omega$, there exists $\bV \subset
\bH(x)$ with $\dim \bV =d$, and $|Q(\bV)| = |\cB(x)|$. 
Now the maximality of
$d$ implies that if $\bfv \not\in \bV$ then $|Q(\bfv)| = 0$. 
\qed\medskip

By Lemma~\ref{lemma:grow:exactly:lambdaj}, $\bE_j(x) \cap \bE_k(x)
= \{0\}$ if $j \ne k$. Let 
\begin{displaymath}
\Lambda' = \{ i \st \bE_i(x) \ne \{0\} \text{ for a.e.\ $x$} \}.
\end{displaymath}
Let the $U^+$-inert subbundle $\bE$ be defined by
\begin{displaymath}
\index{$E$@$\bE(x)$}\bE(x) = \bigoplus_{i \in \Lambda'} \bE_i(x). 
\end{displaymath}
Then $\bE(x) \subset \bH(x)$.

In view of (\ref{eq:meaning:bFj}), (\ref{eq:flag:bFj})  and
Lemma~\ref{lemma:grow:exactly:lambdaj}, we have $\bF_{\ge j}(x) =
\bF_{\ge j+1}(x)$ 
unless $j \in \Lambda'$. 
Therefore if we write the elements of $\Lambda'$ in decreasing order
as $i_1, \dots, i_m$ we have the flag (consisting of distinct subspaces)
\begin{equation}
\label{eq:distinct:Fj:flag}
\{0\} = \bF_{\ge i_{m+1}} \subset \bF_{\ge i_m}(x) \subset \bF_{\ge i_{m-1}}(x) \subset
\dots \bF_{\ge i_2}(x) \subset \bF_{\ge i_1}(x) = \bH(x).
\end{equation}
For a.e.\ $x \in X_0$, and $1 \le r \le m$, 
let \index{$F$@$\bF_j'(x)$}$\bF_{i_r}'(x)$ be the 
orthogonal complement (using the inner
product $\langle \cdot, \cdot \rangle_x$ defined in
\S\ref{sec:subsec:dynamical:norms})
to $\bF_{\ge i_{r+1}}(x)$ in $\bF_{\ge i_r}(x)$. 
\begin{lemma}

\label{lemma:effective:zero:one:substitute}
Given $\delta > 0$ there exists a
compact $K_{01} \subset X_0$ with $\nu(K_{01}) > 1-\delta$,
$\beta(\delta) > 0$, $\beta'(\delta) > 0$, 
and for every $x \in K_{01}$ any $j \in \Lambda'$ any $\bfv' \in
\proj(\bF_j')(x)$ a subset $Q_{01} = Q_{01}(x,\bfv')
\subset \cB(x)$ with $|Q_{01}| > (1-\delta)|\cB(x)|$ 
such that for any $j \in \Lambda'$ any $\bfv' \in
\bF_j'(x)$ and any $u \in Q_{01}$, we can write 
\begin{displaymath}
(u)_* \, \bfv' = \bfv_u + \bfw_u,  \quad \bfv_u \in \bE_j(ux),
\quad \bfw_u \in \bV_{>j}(ux),  
\end{displaymath}
with $\|\bfv_u\| \ge \beta(\delta)\|\bfv'\|$, and $\|\bfv_u\| >
\beta'(\delta) \| \bfw_u\|$.  
\end{lemma}

\bold{Proof.} This is a corollary of
Lemma~\ref{lemma:zero:one:substitute}. 
Let $\Phi \subset X_0$ be the conull set where
(\ref{eq:Lyapunov:transversality}) holds and where $\bF_{\ge i}(x) =
\bF_{\ge i+1}(x)$ for all $i \not\in \Lambda'$. 
Suppose $x \in \Phi$. 

Let $\bF_{\ge k}(x) \subset \bF_{\ge j}(x)$ be the next subspace in the flag
(\ref{eq:distinct:Fj:flag}), 
(i.e.\ $\bF_{\ge k} = \{0\}$ if $j$ is the maximal index in $\Lambda'$ and
otherwise we have $k > j$ be minimal such that $k \in \Lambda'$.) 
Then $\bF_{\ge j+1}(x) = \bF_{\ge k}(x)$. Since $\bF_j'(x)$ is complementary to
$\bF_{\ge k}(x)$ we have that $\bF_j'(x)$ is complementary to $\bF_{\ge j+1}(x)$.

By Lemma~\ref{lemma:Ej:equivariant}, $\bF_{\ge j}$ is $g_t$-equivariant, and
therefore, by the multiplicative ergodic theorem applied to $\bF_{\ge j}$,
$\bF_{\ge j}$ is the direct sum of its Lyapunov subspaces. Therefore, in view of
(\ref{eq:Lyapunov:transversality}), for almost all $y \in X_0$,  
\begin{equation}
\label{eq:Fj:Lyapunov:transversality}
\bF_{\ge j}(y) = (\bF_{\ge j}(y) \cap \bV_{\le j}(y)) \dirsum (\bF_{\ge j}(y) \cap
\bV_{>j}(y)). 
\end{equation}
Since $\bF_j'(x) \subset \bF_{\ge j}(x)$, we have by
Lemma~\ref{lemma:Ej:equivariant}, $(u)_* \bfv' \in \bF_{\ge j}(ux)$ for almost
all $u \in \cB(x)$. By the definition of
$\bF_{\ge j+1}(x)$, since $\bfv' \not\in \bF_{\ge j+1}(x)$, 
for almost all $u$ if we decompose using
(\ref{eq:Fj:Lyapunov:transversality}), 
\begin{displaymath}
(u)_* \, \bfv' = \bfv_{u} + \bfw_{u},  \quad \bfv_{u} \in \bF_{\ge j}(ux)
\cap \bV_{\le j}(u x),
\quad \bfw_{u} \in  \bF_{\ge j}(ux) \cap \bV_{>j}(u x),  
\end{displaymath}
then $\bfv_{u} \ne 0$. Since by definition $\bF_{\ge j}(u x) \cap
\bV_{\le j}(u x) =
\bE_j(u x)$ we have $\bfv_u \in \bE_j(u x)$. Let
\begin{displaymath}
E_n(x) = \{ \bfv' \in \proj(\bF'(x)) \st |\{ u \in \cB(x) \st \|
\bfv_u \| \ge \tfrac{1}{n} \|\bfv'\| \}| > (1-\delta/2)|\cB(x)| \}. 
\end{displaymath}
Then the $E_n(x)$ are an increasing family of 
open sets, and $\bigcup_{n=1}^\infty E_n(x) =
\proj(\bF_j'(x))$. 
Since $\proj(\bF_j'(x))$ is compact, there exists
$n(x)$ such that $E_{n(x)}(x) = \proj(\bF_j'(x))$. We can now choose
$K_{01}' \subset \Phi$ with $\nu(K_{01}') > 1-\delta/2$ such that for $x
\in K_{01}'$, $n(x) < 1/\beta(\delta)$. This shows that for $x \in
K_{01}'$, for any $\bfv' \in \proj(\bF_j'(x))$, 
for $(1-\delta/2)$-fraction of $u \in \cB(x)$ we have
$\|\bfv_u\| >
\beta(\delta)\|\bfv'\|$. 

To prove the final estimate note that there exists a set $K_{01}''$
with $\nu(K_{01}'') > 1-\delta/2$ and a constant $C(\delta)$ such that 
for all $x \in K_{01}''$ and at least $(1-\delta/2)$-fraction of $u
\in \cB(x)$, we have $\|(u)_* \bfv' \| \le C(\delta)
\|\bfv'\|$. Let $K_{01} = K_{01}' \cap K_{01}''$. Then, for at least
$(1-\delta)$-fraction of $u \in \cB(x)$, we have
\begin{displaymath}
\|\bfw_u\| \le \|(u)_* \bfv' \| \le C(\delta) \|\bfv'\| \le C(\delta)
\beta(\delta)^{-1} \|\bfv_u\|. 
\end{displaymath}
\qed\medskip



\begin{proposition}
\label{prop:most:inert}
$ $
\begin{itemize}
\item[{\rm (a)}] For every $\delta > 0$ there exists $K \subset X_0$ of measure at least
$1-\delta$ and a number $L_2(\delta) > 0$ 
such that the following holds: Suppose $x \in K$, 
$\bfv \in \bH(x)$. Then, 
for any $L' > L_2(\delta)$ there exists $L' <
t < 2L'$ such that for at least $(1-\delta)$-fraction of $u \in
\cB(g_{-t} x)$, 
\begin{displaymath}
d\left( \frac{(g_s)_* (u)_* (g_{-t})_* \bfv}{\|(g_s)_* (u)_*
    (g_{-t})_* \bfv\|} , 
     \bE(g_s u g_{-t} x) \right) \le C(\delta) e^{-\alpha t}, 
\end{displaymath}
where $s> 0$ is such that 
\begin{equation}
\label{eq:tmp:most:inert:def:s}
\|(g_s)_* (u)_*
    (g_{-t})_* \bfv\| = \|\bfv\|, 
\end{equation}
and 
$\alpha$ depends only on the Lyapunov spectrum. 
\item[{\rm (b)}] There exists $\epsilon' > 0$ 
(depending only on the
  Lyapunov spectrum) and for every $\delta >
  0$ a compact set $K''$
  with $\nu(K'') > 1 - c(\delta)$ where $c(\delta) \to 0$ as $\delta
  \to 0$ such that the following holds: 
Suppose there exist arbitrarily large $t > 0$ with $g_{-t} x \in K''$
 so that for at least $(1-\delta)$-fraction of $u \in \cB(x)$, the
 number $s > 0$ satisfying (\ref{eq:tmp:most:inert:def:s}), also satisfies
\begin{equation}
\label{eq:tmp:most:inert:b}
s \ge (1-\epsilon') t.
\end{equation}
Then $\bfv \in \bE(x)$. 
\end{itemize}
\end{proposition}
\bold{Proof.} Let $\epsilon > 0$ be smaller than one third of the
difference between any two Lyapunov exponents for the action on
$\bH$. 
By the Osceledets multiplicative ergodic theorem, there
exists a compact subset $K_1 \subset X_0$ 
with $\nu(K_1) > 1- \delta^2$ and $L > 0$ 
such that for $x \in K_1$ and all $j$ and all $t > L$, 
\begin{displaymath}
\| (g_t)_* \bfv \| \le e^{(\lambda_j+\epsilon)t} \|\bfv\|, \qquad \bfv
\in \bV_{\geq j}(x)
\end{displaymath}
and
\begin{displaymath}
\| (g_t)_* \bfv \| \ge e^{(\lambda_j-\epsilon)t} \|\bfv\|, \qquad \bfv
\in {\bV}_{\le j}(x).
\end{displaymath}
By Fubini's theorem there exists $K_1^* \subset X_0$ with $\nu(K_1^*) >
1-2\delta$ such that for $x \in K_1^*$, 
\begin{displaymath}
|\{ u \in \cB(x) \st u x \in K_1 \}| \ge (1-\delta/2) |\cB(x)|. 
\end{displaymath}
Let $K'' = K_{01} \cap K_1^*$, where $K_{01}$ is as in
Lemma~\ref{lemma:effective:zero:one:substitute} (with $\delta$
replaced by $\delta/2$).  
Let $K$, $L_2(\delta)$ be such that for all $x \in K$
and all $L' > L_2$, there exists $t$ with $L' < t < 2L'$ and $g_{-t} x
\in K''$. 
Write 
\begin{equation}
\label{eq:tmp:decomp:g:minus:t:v}
(g_{-t})_* \, \bfv = \sum_{j \in \Lambda'} \bfv_j', \qquad 
  \bfv_j' \in \bF_j'(g_{-t}x).  
\end{equation}
We have $g_{-t} x \in K_{01} \cap K_1^*$. 
Suppose $u \in Q_{01}(g_{-t x})$ and $u g_{-t} x \in
K_1$. Then, by
Lemma~\ref{lemma:effective:zero:one:substitute}, we have
\begin{equation}
\label{eq:tmp:decomp:ustar:gminus:t:bfv}
(u)_* (g_{-t})_* \, \bfv = \sum_{j \in \Lambda'} (\bfv_j + \bfw_j), 
\end{equation}
where $\bfv_j \in \bE_j(u g_{-t} x)$, $\bfw_j \in \bV_{>j}(u g_{-t}
  x)$, and for all $j \in \Lambda'$,
\begin{equation}
\label{eq:tmp:bfvj:ge:beta:delta}
\|\bfv_j\| \ge \beta'(\delta) \|\bfw_j\|. 
\end{equation}
Then, 
\begin{displaymath}
\|(g_s)_* \bfw_j \| \le e^{(\lambda_{j+1}+ \epsilon)s} \|\bfw_j\|,
\end{displaymath}
and,
\begin{equation}
\label{eq:tmp:most:inert:almost:done}
\|(g_s)_* \bfv_j \| \ge e^{(\lambda_j - \epsilon)s} \|\bfv_j\| \ge
e^{(\lambda_j - \epsilon)s} \beta'(\delta) \|\bfw_j\|.
\end{equation}
Thus, for all $j \in \Lambda'$, 
\begin{displaymath}
\|(g_s)_* \bfw_j \| \le
e^{-(\lambda_j-\lambda_{j+1}+2\epsilon)s} 
\beta'(\delta)^{-1} \|(g_s)_* \bfv_j \|.
\end{displaymath}
Since $(g_s)_* \bfv_j \in \bE$ and using part (a) of
Proposition~\ref{prop:properties:dynamical:norm:Hbig}, we get (a) of 
Proposition~\ref{prop:most:inert}.

To prove (b), suppose $\bfv \not \in \bE(x)$. 
We may write
\begin{displaymath}
\bfv = \sum_{i \in \Lambda'} \hat{\bfv}_i, \qquad  \hat{\bfv}_i
  \in \bF_i'(x)
\end{displaymath}
Let $j$ be minimal such that $\hat{\bfv}_j \not\in \bE_j(x)$. Let
$k>j$ be such that $\bF_{\ge k}(x) \subset \bF_{\ge j}(x)$ is the subspace
preceding $\bF_{\ge j}(x)$ in (\ref{eq:distinct:Fj:flag}). Then,
$\bF_{\ge i}(x) = \bF_{\ge j}(x)$ for $k+1 \le i \le j$.  

Since $\hat{\bfv}_j \not\in \bE_j(x)$, $\hat{\bfv}_j$ must have a
component in $\cV_i(\bH)(x)$ for some $i \ge j+1$. Therefore, by
looking only at the component in $\cV_i(\bH)$, we get
\begin{displaymath}
\|(g_{-t})_* \bfv \| \ge C(\bfv) e^{-(\lambda_{j+1}+\epsilon) t},
\end{displaymath}
Also since $\bF_{\ge k}$ is $g_t$-equivariant we have $\bF_{\ge k}(x)
= \bigoplus_{m}
\bF_{\ge k}(x) \cap \cV_m(\bH)$. Note that by the multiplicative
ergodic theorem, the restriction of $g_{-t}$ to $\cV_i(\bH)$ is of the
form $e^{-\lambda_i t} h_t$, where $\|h_t\| = O(e^{\epsilon t})$.
Therefore, (again by looking only at the
component in $\cV_i(\bH)$ and using
Proposition~\ref{prop:properties:dynamical:norm:Hbig} (a)), we get
\begin{displaymath}
d((g_{-t})_* \bfv, \bF_{\ge k}(g_{-t} x)) \ge C(\bfv)
e^{-(\lambda_{j+1}+2\epsilon) t}.
\end{displaymath}
(Here and below, \index{$d(\cdot,\cdot)$}$d(\cdot,\cdot)$ 
denotes the distance on $\bH(x)$
given by the dynamical norm $\| \cdot \|_x$.)
Therefore, (since $(g_{-t})_* \bfv \in \bF_{\ge j}(g_{-t}x)$),
we see that if we decompose $(g_{-t})_* \bfv$ as in
(\ref{eq:tmp:decomp:g:minus:t:v}), we get
\begin{displaymath}
\|\bfv_j' \| \ge C(\bfv) e^{-(\lambda_{j+1}+2\epsilon)t},
\end{displaymath}
We now decompose $(u)_* (g_{-t})_* \bfv$ as in
(\ref{eq:tmp:decomp:ustar:gminus:t:bfv}).
Then, from (\ref{eq:tmp:bfvj:ge:beta:delta}) and
(\ref{eq:tmp:most:inert:almost:done}), 
\begin{equation}
\label{eq:tmp:bfvj:lower:bound}
\|(g_s)_* \bfv_j \| \ge e^{(\lambda_j - \epsilon)s} \|\bfv_j\| \ge
e^{(\lambda_j - \epsilon)s} \beta(\delta) \|\bfv_j'\| \ge
e^{(\lambda_j - \epsilon)s} \beta(\delta) C(\bfv) e^{-(\lambda_{j+1}+2\epsilon)t}.
\end{equation}
If $s$ satisfies (\ref{eq:tmp:most:inert:def:s}), then $\|(g_s)_*
\bfv_j \| = O(1)$. 
Therefore, in view of (\ref{eq:tmp:bfvj:lower:bound}), 
\begin{displaymath}
e^{(\lambda_j-\epsilon) s} e^{-(\lambda_{j+1}+2\epsilon)t}
\le c = c(\bfv,\delta). 
\end{displaymath}
Therefore, 
\begin{displaymath}
s \le \frac{(\lambda_{j+1} + 2\epsilon)t + \log c(\bfv,
  \delta)}{(\lambda_j - \epsilon)}. 
\end{displaymath}
Since $\lambda_j > \lambda_{j+1}$,
this contradicts (\ref{eq:tmp:most:inert:b}) 
if $\epsilon$ is sufficiently small and 
$t$ is sufficiently large. \mcc{(more details here?)}
\qed\medskip

\section{The action of the cocycle on $\bE$}
\label{sec:action:Eplus}

In this section, we work on the finite cover $X$ defined in
\S\ref{sec:subsec:the:cover:X}. Recall that if $f(\cdot)$ is an object
defined on $X_0$, then for $x \in X$ we write $f(x)$ instead of
$f(\sigma_0(x))$ (where $\sigma_0: X \to X_0$ is the covering map). 

In this section and in \S\ref{sec:bounded:synchornized}, 
assertions will hold at best for a.e $x \in X$, and never for all $x \in X$.
This will be sometimes suppressed from the statements of the lemmas. 

\subsection{The Jordan canonical form of the cocycle on $\bE(x)$}



We consider the action of the cocycle on $\bE$. The Lyapunov
exponents are $\lambda_i$, $i \in \index{$\Lambda'$}\Lambda'$. 
We note that by
Lemma~\ref{lemma:Ej:equivariant}, 
the bundle $\bE$ admits the equivariant measurable flat
$U^+$-connection given by the maps $(u)_*: \bE(x) \to \bE(y)$, where
$(u)_*$ is as in Lemma~\ref{lemma:tensor:product}. This connection
satisfies the condition (\ref{eq:Fxy:preserves:lyapunov:subspaces}), since by 
Lemma~\ref{lemma:Ej:equivariant}, $(u)_* \bE_j(x) = \bE_j(y)$. 
For each $i \in \Lambda'$, we have the maximal flag as in
Lemma~\ref{lemma:jordan:canonical:form}, 
\begin{equation}
\label{eq:flag:bEij}
\{ 0\} \subset \index{$E$@$\bE_{ij}(x)$}\bE_{i1}(x) \subset \dots \subset \bE_{i,n_i}(x) =
\bE_i(x).
\end{equation}
Let \index{$\Lambda''$}$\Lambda''$ denote the set of pairs $ij$ which appear in
(\ref{eq:flag:bEij}). 
By Proposition~\ref{prop:sublyapunov:locally:constant} and
Remark~\ref{remark:Uplus:connection}, we have for a.e.\ $u \in \cB(x)$, 
\begin{displaymath}
(u)_* \bE_{ij}(x) = \bE_{ij}(u x). 
\end{displaymath}
Let \index{abs@$\norm{\cdot}_x$}$\| \cdot \|_x$ and $\index{$\langle \cdot, \cdot \rangle_x$}\langle \cdot, \cdot \rangle_x$ denote the
restriction to $\bE(x)$ of the norm and inner product on
$\bH(x)$ 
defined in \S\ref{sec:subsec:dynamical:norms} and
\S\ref{sec:divergence:subspaces}.
(We will often omit the subscript from $\langle \cdot, \cdot \rangle_x$ and
$\| \cdot \|_x$.) \index{abs@$\norm{\cdot}$} \index{$\langle \cdot,
  \cdot \rangle$}
Then, the distinct $\bE_i(x)$ are orthogonal. For each $ij \in \Lambda''$ let
\index{$E$@$\bE'_{ij}(x)$}$\bE'_{ij}(x)$
be the orthogonal complement (relative to the inner
product $\langle \cdot, \cdot \rangle_x$) to
$\bE_{i,j-1}(x)$ in $\bE_{ij}(x)$. 

Then, by Proposition~\ref{prop:properties:dynamical:norm:Hbig}, 
we can write, for $\bfv \in \bE'_{ij}(x)$,
\begin{equation}
\label{eq:gt:star:bfv}
(g_t)_* \bfv = e^{\index{$\lambda_{ij}(x,t)$}\lambda_{ij}(x,t)} \bfv' + \bfv'',
\end{equation}
where $\bfv' \in \bE'_{ij}(g_t x)$, $\bfv'' \in \bE_{i,j-1}(g_t x)$, 
and $\|\bfv' \| = \|\bfv\|$. Hence (since $\bfv'$ and $\bfv''$ are
orthogonal), 
\begin{displaymath}
\|(g_t)_* \bfv\| \ge e^{\lambda_{ij}(x,t)} \|\bfv\|. 
\end{displaymath}
In view of Proposition~\ref{prop:properties:dynamical:norm:Hbig}
there exists a constant $\kappa > 1$ such that
for a.e $x \in X$ and for all $\bfv \in \bE(x)$ and all $t \ge 0$, 
\begin{equation}
\label{eq:hyperbolic:dynamical:norm:E}
e^{\kappa^{-1} t } \|\bfv\| \le \| (g_t)_* \bfv \| \le e^{\kappa t} \|
\bfv\|. 
\end{equation}

\begin{lemma}
\label{lemma:u:star:agrees:with:Pplus}
For a.e.\ $x \in X$ and for a.e.\ $y = u x \in \cB[x]$, the 
connection $(u)_*: \bE(x) \to \bE(y)$ 
agrees with the restriction to $\bE$ of the 
connection \index{$P$@$\bP^+(x,y)$}$\bP^+(x,y)$ induced from
the map $P^+(x,y)$ defined in \S\ref{sec:subsec:connection}.
\end{lemma}

\bold{Proof.} Let $\cV_{\le i}(x) = \cV_{\le i}(H^1)(x)$ and $\cV_i(x)
= \cV_i(H^1)(x)$, where $\cV_{\le i}(H^1)(x)$ and $\cV_i(H^1)(x)$ are as in
  \S\ref{sec:subsec:lyapunov:flags}.
Consider the definition
(\ref{eq:def:u:lower:star}) of $u_*$ in \S\ref{sec:divergence:subspaces}. 
For a fixed $Y = \log u \in \Lie(U^+)(x)$ and $M \in \cH_{++}(x)$, let
$h: W^+(x) \to W^+(ux)$ be given by
\begin{displaymath}
h(v) = \exp((I+M)Y)(x+v) - \exp(Y)x.
\end{displaymath}
From the form of $h$, we see that $h(\cV_{\le i}(x)) = \cV_{\le i}(ux)$, 
and also, $h$ induces the identity map on $\cV_{\le i}(x)/\cV_{< i}(x) =
\cV_{\le i}(ux)/\cV_{< i}(ux)$. Thus, for $v \in \cV_i(x)$, 
\begin{displaymath}
h(v) \in P^+(x,ux) v + \cV_{< i}(ux). 
\end{displaymath}
Similarly, $M'' \equiv tr(x,ux) \circ M \circ tr(ux,x)$ agrees with
$M$ up to higher Lyapunov exponents. Then, in view of
(\ref{eq:def:u:lower:star}), (\ref{eq:def:bold:Pplus})
  and Lemma~\ref{lemma:Lyapunov:exponents:on:bH}, for $\bfv \in \bE_i(x)$,
\begin{displaymath}
(u)_* \bfv \in \bP^+(x,ux) \bfv + \bV_{< i}(ux).
\end{displaymath}
But, for $\bfv \in \bE_i(x)$, $(u)_* \bfv \in \bE_i(ux)$ (and thus has
no component in $\bV_{< i}(ux)$). Hence, for all $\bfv \in \bE_i(x)$,
we have $(u)_* \bfv = \bP^+(x,ux) \bfv.$ 
\qed\medskip

\subsection{Time changes}
\label{sec:subsec:timechanges}
$ $

\bold{The flows $g^{ij}_t$ and the time changes $\hat{\tau}_{ij}(x,t)$.}
We define the time changed flow \index{$g^{ij}_t$}$g^{ij}_t$ so that (after the time
change) the cocycle $\lambda_{ij}(x,t)$ of (\ref{eq:gt:star:bfv})
becomes $\lambda_i t$. 
We write $g^{ij}_t x = g_{\hat{\tau}_{ij}(x,t)} x$. Then, by construction,
$\lambda_{ij}(x,\index{$\tau$@$\hat{\tau}_{ij}(x,t)$}\hat{\tau}_{ij}(x,t)) = \lambda_i t$.
We note the following:
\begin{lemma}
\label{lemma:stay:in:same:orbit}
Suppose $y \in \gB_0[x]$. Then for any $ij \in \Lambda''$
and any $t > 0$, 
\begin{displaymath}
g_{-t}^{ij} y \in \gB_0[g_{-t}^{ij} x]. 
\end{displaymath}
\end{lemma}

\bold{Proof.} This follows immediately from property (e) of
Proposition~\ref{prop:properties:dynamical:norm:Hbig}, and the
definition of the flow $g_{-t}^{ij}$. 
\qed\medskip

In view of Proposition~\ref{prop:properties:dynamical:norm:Hbig}, we have
\begin{equation}
\label{eq:bilip:hat:tau:ij}
\frac{1}{\kappa} |t-t'| \le |\hat{\tau}_{ij}(x,t) -  \hat{\tau}_{ij}(x,t')|
\le \kappa |t-t'| 
\end{equation}
where $\kappa$ depends only on the Lyapunov spectrum.

\subsection{The foliations $\cF_{ij}$, $\cF_\bfv$ and the parallel transport
  $R(x,y)$}

\label{sec:foliations}


For $x \in \tilde{X}$, let
\begin{displaymath}
\index{$G[x]$}G[x] = \{ g_s u g_{-t} x \st t \ge 0, s \ge 0, u \in \cB(g_{-t} x) \}
\subset \tilde{X}. 
\end{displaymath}
For $y =g_s u g_{-t} x \in G[x]$, let 
\begin{displaymath}
\index{$R(x,y)$}R(x,y) = (g_s)_* (u)_* (g_{-t})_*.
\end{displaymath}
Here $(g_s)_*$ is as in (\ref{eq:def:action:on:bHplus})
 and $(u)_*:
\bH(g_{-t} x) \to \bH(u g_{-t} x)$ is as in
Lemma~\ref{lemma:tensor:product}. 
It is easy to see using Lemma~\ref{lemma:independent:repn:u}
that $R(x,y): \bH(x) \to \bH(y)$ depends only on $x,y$ and not on the
choices of $t$, $u$, $s$. We will usually consider $R(x,y)$ as
a map from $\bE(x) \to \bE(y)$.  

In view of (\ref{eq:gt:star:bfv}),
Lemma~\ref{lemma:u:star:agrees:with:Pplus} and
Proposition~\ref{prop:properties:dynamical:norm:Hbig} (e) and (f), 
we have, for $\bfv \in \bE'_{ij}(x)$, and any
$y = g_s u g_{-t} x \in G[x]$, 
\begin{equation}
\label{eq:Rxy:bfv}
R(x,y) \bfv = e^{\lambda_{ij}(x,y)} \bfv' + \bfv''
\end{equation}
where $\bfv' \in \bE'_{ij}(y)$, $\bfv'' \in \bE_{i,j-1}(y)$, 
and $\|\bfv' \| = \|\bfv\|$. In (\ref{eq:Rxy:bfv}), we have
\begin{equation}
\label{eq:def:lambdaij:xy}
\index{$\lambda_{ij}(x,y)$}\lambda_{ij}(x,y) = \lambda_{ij}(x,-t) + \lambda_{ij}(u g_{-t}x, s).
\end{equation}

\bold{Notational convention.} 
We sometimes use the notation $R(x,y)$ when $x \in X$ (instead of
$\tilde{X}$) and $y \in G[x]$. 
\medskip

For $x \in \tilde{X}$ and $ij \in \Lambda''$, let 
\index{$F$@$\cF_{ij}[x]$}$\cF_{ij}[x]$ denote the set of $y \in G[x]$ such that there exists
$\ell \ge 0$ so that
\begin{equation}
\label{eq:def:Fv}
g^{ij}_{-\ell}y \in \cB[g^{ij}_{-\ell} x]. 
\end{equation}
By Lemma~\ref{lemma:stay:in:same:orbit}, if (\ref{eq:def:Fv}) holds
for some $\ell$, it also holds for any bigger $\ell$. 
Alternatively, 
\begin{displaymath}
\cF_{ij}[x] = \{ g_{\ell}^{ij} u g_{-\ell}^{ij} x \st \ell \ge 0, \quad u \in
\cB(g_{-\ell}^{ij} x) \} \subset \tilde{X}. 
\end{displaymath}
As above, when $x \in X$, we can think of the leaf of the
foliation $\cF_{ij}[x]$ as a subset of $X$ (not $\tilde{X}$). 

In view of (\ref{eq:def:lambdaij:xy}), it follows that 
\begin{equation}
\label{eq:lambdaij:zero:on:Fij}
\lambda_{ij}(x,y) = 0 \qquad \text{ if $y \in \cF_{ij}[x]$.}
\end{equation}
We refer to the sets $\cF_{ij}[x]$ as {\em leaves}. Locally, the leaf
$\cF_{ij}[x]$ through $x$ is a piece of $U^+[x]$. More precisely, for $y
\in \cF_{ij}[x]$, 
\begin{displaymath}
\cF_{ij}[x] \cap \gB_0[y] \subset U^+[y]. 
\end{displaymath}
Then, for any compact subset $A \subset \cF_{ij}[x]$ there exists $\ell$
large enough so that $g^{ij}_{-\ell}(A)$ is contained in a set of
the form $\cB[z] \subset U^+[z]$. Then the same holds for
$g^{ij}_{-t}(A)$, for any $t > \ell$. 

Recall (from the start of \S\ref{sec:divergence:subspaces}) 
that the sets $\cB[x]$ support a ``Lebesgue measure'' $| \cdot |$, 
namely the
pushforward of the Haar measure on $U^+(x)/(U^+(x) \cap Q_{++}(x))(x)$ 
to $\cB[x]$ under the
map $u \to u x$. (Recall that $Q_{++}(x)$ is the stabilizer of $x$ in the
affine group $\cG_{++}(x)$). 
As a consequence, the leaves $\cF_{ij}[x]$ also support a Lebesgue
measure (defined up to normalization), which we also denote by 
\index{abs@$\abs{\cdot}$}$| \cdot |$. More precisely, if $A \subset \cF_{ij}[x]$ and $B \subset \cF_{ij}[x]$
are compact subsets, we define
\begin{equation}
\label{eq:def:lebesgue:leaf}
\frac{|A|}{|B|} \equiv \frac{|g^{ij}_{-\ell}(A)|}{|g^{ij}_{-\ell}(B)|}, 
\end{equation}
where $\ell$ is chosen large enough so that both
$g^{ij}_{-\ell}(A)$ and $g^{ij}_{-\ell}(B)$ are contained in a
set of the form $\cB[z]$, 
$z \in X$. It is clear that if we replace
$\ell$ by a larger number, the right-hand-side of
(\ref{eq:def:lebesgue:leaf}) remains the same. \mcc{(explain this?)}

We define the ``balls'' $\cF_{ij}[x,\ell] \subset \cF_{ij}[x]$ by 
\begin{equation}
\label{eq:def:Fij:ball}
\index{$F$@$\cF_{ij}[x,\ell]$}\cF_{ij}[x,\ell] = \{ y \in \cF_{ij}[x] \st g^{ij}_{-\ell} y \in \cB[
g^{ij}_{-\ell} x] \}.
\end{equation}

\begin{lemma}
\label{lemma:symmetric:difference:balls}
Suppose $x \in \tilde{X}$ and $y \in \cF_{ij}[x]$. Then, for $\ell$ large
enough, 
\begin{displaymath}
\cF_{ij}[x,\ell] = \cF_{ij}[y,\ell]. 
\end{displaymath}
\end{lemma}

\bold{Proof.} Suppose $y \in \cF_{ij}[x]$. Then, for $\ell$ large
enough, $g^{ij}_{-\ell} y \in \cB[g^{ij}_{-\ell} x]$, 
and then $\cB[g^{ij}_{-\ell} y] = \cB[g^{ij}_{-\ell} x]$. 
\qed\medskip

\bold{The ``flows'' $g^\bfv_t$.}
Suppose $x \in \tilde{X}$  
and $\bfv \in \bE(x)$. Let $\index{$g^\bfv_t$}g^\bfv_t x = g_{\hat{\tau}_\bfv(x,t)} x$,
  where the time change \index{$\tau$@$\hat{\tau}_\bfv(x,t)$}$\hat{\tau}_\bfv(x,t)$ is chosen so that 
\begin{displaymath}
\|(g^\bfv_t)_* \bfv\|_{g^{\bfv}_t x} = e^t \|\bfv\|_x.
\end{displaymath}
(Note that we are not defining $g_t^{\bfv} y$ for $y \ne x$). We have,
for $x \in \tilde{X}$, 
\begin{displaymath}
g_{t+s}^{\bfv} x = g_s^{(g_t)_* \bfv} g_t^{\bfv} x.
\end{displaymath}
By (\ref{eq:hyperbolic:dynamical:norm:E}), 
(\ref{eq:bilip:hat:tau:ij}) holds for $\hat{\tau}_\bfv$ instead of
$\hat{\tau}_{ij}$.

For $y \in G[x]$ and $\ell \in \reals$, let 
\begin{equation}
\label{eq:def:tilde:g}
\index{$g$@$\tilde{g}^{\bfv,x}_{-\ell}$}\tilde{g}^{\bfv,x}_{-\ell}
= g_{-\ell}^{\bfw} y, \quad \text{ where $\bfw =
  R(x,y)\bfv$. }
\end{equation}
(When there is no potential for confusion about the point $x$ and the
vector $\bfv$ used, we
denote $\tilde{g}^{\bfv,x}_{-\ell}$ by $\index{$g$@$\tilde{g}_{-\ell}$}\tilde{g}_{-\ell}$.)
Note that Lemma~\ref{lemma:stay:in:same:orbit} still holds if
$g^{ij}_{-t}$ is replaced by $\tilde{g}^{\bfv,x}_{-t}$.

\bold{The foliations $\cF_{\bfv}$.} For $\bfv \in \bE(x)$ we
can define the foliations  \index{$F$@$\cF_{\bfv}[x]$}$\cF_{\bfv}[x]$
and the ``balls'' 
\index{$F$@$\cF_{\bfv}[x,\ell]$}$\cF_{\bfv}[x,\ell]$ as in
(\ref{eq:def:Fv})  and
(\ref{eq:def:Fij:ball}), with $\tilde{g}^{\bfv,x}_{-t}$
  replacing the role of $g^{ij}_{-t}$.

For $y \in \cF_{\bfv}[x]$, we have
\begin{displaymath}
\cF_{\bfv}[x] = \cF_{\bfw}[y], \qquad \text{ where } \bfw = R(x,y) \bfv. 
\end{displaymath}
\mcc{(give proof of above formula?)}

We can 
define the measure (up to normalization) \index{abs@$\abs{\cdot}$}$| \cdot |$ on
$\cF_{\bfv}[x,\ell]$ as in (\ref{eq:def:lebesgue:leaf}).
Lemma~\ref{lemma:symmetric:difference:balls} holds for $\cF_{\bfv}[x]$
without modifications.

The following follows immediately from the
construction: 
\begin{lemma}
\label{lemma:Fv:norm}
For a.e.\ $x \in \tilde{X}$, any $\bfv \in \bE(x)$, and a.e.\ $y \in
\cF_\bfv[x]$, we have
\begin{displaymath}
\|R(x,y) \bfv \|_y = \|\bfv\|_x. 
\end{displaymath}
\end{lemma}

\subsection{A maximal inequality}
\begin{lemma}
\label{lemma:fake:ergodicity:Fij}
Suppose $K \subset X$ with $\nu(K) > 1-\delta$. Then, for any
$\theta' > 0$ there exists a
subset $K^* \subset X$ with $\nu(K^*) > 1-2\kappa^2
\delta/\theta'$ such that 
for any $x \in K^*$ and any $\ell > 0$, 
\begin{equation}
\label{eq:fake:ergodicity:Fij}
|\cF_{ij}[x,\ell] \cap K| > (1-\theta') |\cF_{ij}[x,\ell]|. 
\end{equation}
\end{lemma}


\bold{Proof.} For $t > 0$ let 
\begin{displaymath}
  \cB_t^{ij}[x] = g_{-t}^{ij}
  (\gB_0[g_t^{ij} x] \cap U^+[g_t^{ij} x]) = \cB_\tau[x], 
\end{displaymath}
where $\tau$ is such that $g_\tau x = g_t^{ij} x$. 
Let $s > 0$ be arbitrary. Let $K_s = g_{-s}^{ij} K$. Then
$\nu(K_s) > 1-\kappa \delta$. Then, by 
Lemma~\ref{lemma:cB:vitali:substitute}, there exists a subset $K_s'$ with
$\nu(K_s') \ge (1-2 \kappa \delta/\theta')$ such that
for $x \in K_s'$ and all $t > 0$, 
\begin{displaymath}
| K_s \cap \cB_t^{ij}[x] | \ge (1-\theta'/2) |K_s|. 
\end{displaymath}
Let $K_s^* = g_s^{ij} K_s'$, and 
note that $g_s^{ij}\cB_t^{ij}[x] = \cF_{ij}[g_s^{ij} x, s-t]$. Then, for all 
$x \in K_s^*$ and all $0 < s - t < s$, 
\begin{displaymath}
|\cF_{ij}[x,s-t] \cap K| \ge (1-\theta'/2) |\cF_{ij}[x,s-t]|. 
\end{displaymath}
We have $\nu(K_s^*) \ge (1-2\kappa^2 \delta/\theta')$. Now take a
sequence $s_n \to\infty$, and let $K^*$
be the set of points which are in infinitely many $K_{s_n}^*$. 
\qed\medskip

\section{Bounded subspaces and synchronized exponents}
\label{sec:bounded:synchornized}

Recall that $\Lambda''$ indexes the ``fine Lyapunov spectrum'' on $\bE$.
In this section we define an equivalence relation called
``synchronization'' on $\Lambda''$; the equivalence class of $ij \in
\Lambda''$ is denoted by $[ij]$ and the set of equivalence classes is
denoted by \index{$\Lambda$@$\tilde{\Lambda}$}$\tilde{\Lambda}$. For each $ij \in \Lambda''$ we define
a $g_t$-equivariant and locally $(u)_*$-equivariant (in the sense of
Lemma~\ref{lemma:tensor:product} (b)) 
subbundle  \index{$E$@$\bE_{ij,bdd}(x)$}$\bE_{ij,bdd}$ of the bundle
$\bE_{i} \equiv \cV_i(\bE)$ and we define
\begin{displaymath}
\index{$E$@$\bE_{[ij],bdd}(x)$}\bE_{[ij],bdd}(x) = \sum_{kr \in [ij]} \bE_{kr,bdd}(x).
\end{displaymath}
In fact we will show that there exists a subset $[ij]' \subset
  [ij]$ such that 
\begin{equation}
\label{eq:def:ij:prime}
\index{$E$@$\bE_{[ij],bdd}(x)$}\bE_{[ij],bdd}(x) = \bigoplus_{kr \in [ij]'} \bE_{kr,bdd}(x).
\end{equation}
Then, we claim that the following three propositions hold:
\begin{proposition}
\label{prop:some:fraction:bounded}
There exists $\theta > 0$ depending only on $\nu$ and $n \in \natls$
depending only on the dimension of $X$ such that the
following holds: 
for every $\delta > 0$ and every $\eta > 0$, 
there exists a subset $K=K(\delta,\eta)$ of measure at least
$1-\delta$ and $L_0 = L_0(\delta,\eta) > 0$ 
such that the following holds: Suppose $x \in X$, 
$\bfv \in \bE(x)$, $L \ge L_0$, 
and 
\begin{displaymath}
| g_{[-1,1]} K \cap \cF_\bfv[x,L]| \ge (1-(\theta/2)^{n+1})
|\cF_\bfv[x,L]|. 
\end{displaymath}
Then, for at least $(\theta/2)^n$-fraction of $y \in \cF_\bfv[x,L]$,
\begin{displaymath}
d\left(\frac{R(x,y) \bfv}{\|R(x,y) \bfv\|},
    \bigcup_{ij \in \tilde{\Lambda}} \bE_{[ij],bdd}(y) \right) < \eta. 
\end{displaymath}
\end{proposition}


\begin{proposition}
\label{prop:ej:bdd:transport:bounded}
There exists a function $C_3: X \to \reals^+$ finite almost
everywhere so that for all $x \in \tilde{X}$, for all $y \in
\cF_{ij}[x]$, for all $\bfv \in
\bE_{[ij],bdd}(x)$, 
\begin{displaymath}
C_3(x)^{-1}C_3(y)^{-1} \|\bfv\| \le \|R(x,y) \bfv \| \le 
C_3(x)C_3(y) \|\bfv\|.  
\end{displaymath}
(Recall from \S\ref{sec:subsec:notation} that by $C_3(x)$ we mean $C_3(\pi(x))$.) 
\end{proposition}

\begin{proposition}
\label{prop:Rxy:v:small:implied:sync:bounded}
There exists $\theta > 0$ (depending only on $\nu$) and a subset
$\Psi \subset X$ with $\nu(\Psi) = 1$ such that the
following holds: 

Suppose $x \in \Psi$, $\bfv \in \bH(x)$, 
and there exists $C > 0$ such that for all $\ell
> 0$, and at least
$(1-\theta)$-fraction of $y \in \cF_{ij}[x,\ell]$, 
\begin{displaymath}
\| R(x,y) \bfv \| \le C \|\bfv \|.  
\end{displaymath}
Then, $\bfv \in \bE_{[ij],bdd}(x)$. 
\end{proposition}

Proposition~\ref{prop:some:fraction:bounded} is what allows us to
choose $u$ so that there exists $u'$ such that the vector in $\bH$
associated to the difference between the generalized subspaces
$U^+[g_t u' q_1']$  and $U^+[g_t u q_1]$ points close to a controlled
direction, i.e.\ close to $\bE_{[ij],bdd}(g_t u q_1)$. 
This allows us to address
``Technical Problem \#3'' from \S\ref{sec:outline:step1}. 
Then, 
Proposition~\ref{prop:ej:bdd:transport:bounded} and
Proposition~\ref{prop:Rxy:v:small:implied:sync:bounded} 
are used in
\S\ref{sec:equivalence}  to define and
control conditional measures $f_{ij}$ associated to each $[ij] \in
\tilde{\Lambda}$, so we can implement the outline in
\S\ref{sec:outline:step1}. We note that it is important for us to
define a family of subspaces so that all three propositions hold.

The number $\theta > 0$, the synchronization relation and the
subspaces $\bE_{ij,bdd}$ are defined in
\S\ref{sec:starredsubsec:bounded:synchronized}. Also
Proposition~\ref{prop:some:fraction:bounded} is proved in
\S\ref{sec:starredsubsec:bounded:synchronized}. 
Proposition~\ref{prop:ej:bdd:transport:bounded} and
Proposition~\ref{prop:Rxy:v:small:implied:sync:bounded} 
are proved in
\S\ref{sec:starredsubsec:invariant:measures:proj}. Both subsections
may be skipped on first reading.

\bold{Example.} To completely understand the example below, it
necessary to read at least
\S\ref{sec:starredsubsec:bounded:synchronized}. However, we include it
here to give some flavor of the construction. 

Suppose we have a basis $\{\bfe_1(x), \bfe_2(x), \bfe_3(x), \bfe_4(x)\}$
for $\bE(x)$, relative to which the cocycle has the form (for $y
\in G[x]$): 
\begin{displaymath}
R(x,y) = \begin{pmatrix} 
e^{\lambda_{11}(x,y)}& u_{12}(x,y)           & 0                 &0\\
0                    & e^{\lambda_{12}(x,y)} & 0                 &0\\
0                    &  0                    & e^{\lambda_{31}(x,y)}&0\\
0                    &  0                    & 0 &
e^{\lambda_{41}(x,y)}\\
\end{pmatrix}.
\end{displaymath}
Suppose $\bE_1(x) = \reals \bfe_1(x) \oplus \reals \bfe_2(x)$ (so
$\bfe_1$ and $\bfe_2$ correspond to the Lyapunov exponent
$\lambda_1$), $\bE_3(x) = \reals \bfe_3(x)$, $\bE_4(x) = \reals
\bfe_4(x)$ (so that $\bfe_3$ and $\bfe_4$ correspond to the Lyapunov
exponents $\lambda_3$ and $\lambda_4$ respectively). Therefore the
Lyapunov exponents $\lambda_3$ and $\lambda_4$ have multiplicity $1$,
while $\lambda_1$ has multiplicity $2$. 

Then, we have
\begin{displaymath}
\bE_{31,bdd}(x) = \reals \bfe_3(x), \qquad \bE_{41,bdd}(x) = \reals
\bfe_4(x), \qquad \bE_{11,bdd}(x) = \reals \bfe_1(x).  
\end{displaymath}
(For example, if $y \in \cF_{31}[x]$ then $\lambda_{31}(x,y) = 0$, so
that by (\ref{eq:Rxy:bfv}), $\|R(x,y) \bfe_3 \| = \| \bfe_3 \|$.)

Now suppose that $31$ and $41$ are
synchronized, but all other pairs are not synchronized. 
(See Definition~\ref{def:synchronized} for the exact definition of
synchronization, 
but roughly this means that $|\lambda_{41}(x,y)|$ is bounded as $y$
varies over $\cF_{31}[x]$, but for all other distinct pairs $ij$ and $kl$,
$|\lambda_{ij}(x,y)|$ is essentially
unbounded as $y$ varies over $\cF_{kl}[x]$). Then, 
\begin{displaymath}
\bE_{[31],bdd}(x) = \reals \bfe_3(x) \oplus \reals \bfe_4(x), 
\end{displaymath}
Depending on the boundedness behavior  of $u_{12}(x,y)$ as  $y$ varies
over $\cF_{12}[x]$ 
we would have either
\begin{displaymath}
\bE_{12,bdd}(x) = \{0\} \qquad \text{ or } \qquad \bE_{12,bdd}(x) = \reals
\bfe_2(x). 
\end{displaymath}
Since $[11]' = \{11\}$ and $[12]' = \{12\}$, we have $\bE_{[11],bdd}(x) =
\bE_{11,bdd}(x)$ and $\bE_{[12],bdd}(x) = \bE_{12,bdd}(x)$.

\starredsubsection{Bounded subspaces and synchronized exponents.}
\label{sec:starredsubsec:bounded:synchronized}


For $x \in \tilde{X}$, $y \in \tilde{X}$, let
\begin{displaymath}
\index{$\rho(x,y)$}\rho(x,y) = \begin{cases} |t| & \text{ if }  y = g_t x, \\
 \infty & \text{otherwise.}
\end{cases}
\end{displaymath}
If $x \in \tilde{X}$ and $E \subset \tilde{X}$, we let $\index{$\rho(x,E)$}\rho(x,E) =
\inf_{y \in E} \rho(x,y)$.

\begin{lemma}
\label{lemma:bounded:implies:lower:rank}
For every $\eta > 0$ and $\eta' > 0$ there exists $h = h(\eta',\eta)$
such that the following holds: 
Suppose $\bfv \in \bE_{ij}(x)$ and 
\begin{displaymath}
d\left(\frac{\bfv}{\|\bfv\|}, \bE_{i,j-1}(x)\right) > \eta'. 
\end{displaymath}
Then if
$y \in \cF_{\bfv}[x]$ and 
\begin{displaymath}
\rho(y, \cF_{ij}[x]) > h
\end{displaymath}
then
\begin{displaymath}
d(R(x,y) \bfv, \bE_{i,j-1}(y)) \le \eta \|\bfv\|.
\end{displaymath}
\end{lemma}

\bold{Proof.} There exists $t \in \reals$ such that
$y' = g_t y \in \cF_{ij}[x]$. Then
\begin{displaymath}
\rho(y, \cF_{ij}[x]) = \rho(y,y') = |t| > h.
\end{displaymath}
We have the orthogonal decomposition 
$\bfv = \hat{\bfv} + \bfw$,
where $\hat{\bfv} \in \bE_{ij}'(x)$ and 
$\bfw \in \bE_{i,j-1}(x)$. 
Then by (\ref{eq:Rxy:bfv}) we have the orthogonal decomposition.
\begin{displaymath}
R(x,y') \hat{\bfv} = e^{\lambda_{ij}(x,y')} \bfv' + \bfw', 
\quad\text{ where $\bfv' \in \bE'_{ij}(y')$, $\bfw' \in
  \bE_{i,j-1}(y')$, $\|\hat{\bfv}\| = \|\bfv'\|$. }  
\end{displaymath}
Since $R(x,y') \bfw \in \bE_{i,j-1}(y')$, we have
\begin{displaymath}
\|R(x,y') \bfv \|^2 = e^{2\lambda_{ij}(x,y')} \|\hat{\bfv}\|^2 
+ \|\bfw' + R(x,y')\bfw\|^2 \ge e^{2\lambda_{ij}(x,y')} \|\hat{\bfv}\|^2. 
\end{displaymath}
By (\ref{eq:lambdaij:zero:on:Fij}), we have $\lambda_{ij}(x,y') = 0$. 
Hence, 
\begin{displaymath}
\|R(x,y') \bfv \| \ge \| \hat{\bfv} \| \ge \eta' \|\bfv\|.
\end{displaymath}
Since $y \in \cF_{\bfv}[x]$, $\|R(x,y) \bfv \| = \| \bfv \|$.
Since $|t| > h$, we have either $t > h$ or $t < -h$. If $t < -h$, then
by (\ref{eq:hyperbolic:dynamical:norm:E}) and Lemma~\ref{lemma:Fv:norm},
\begin{displaymath}
\|\bfv\| = \|R(x,y) \bfv\| = \| (g_{-t})_* R(x,y')
\bfv \| \ge e^{\kappa^{-1} h} \|R(x,y') \bfv\| \ge e^{\kappa^{-1}h} \eta' \|\bfv\|, 
\end{displaymath}
which is a contradiction if $h > \kappa \log(1/\eta')$. Hence
we may assume that $t > h$. We have,
\begin{displaymath}
R(x,y) \bfv = e^{\lambda_{ij}(x,y)} \bfv'' + \bfw'' 
\end{displaymath}
where $\bfv'' \in \bE'_{ij}(y)$ with $\|\bfv''\| = \|\hat{\bfv}\|$,
and $\bfw'' \in \bE_{i,j-1}(y)$. Hence, 
\begin{displaymath}
d(R(x,y) \bfv, \bE_{i,j-1}(y)) = e^{\lambda_{ij}(x,y)} \|\hat{\bfv}\|
\le e^{\lambda_{ij}(x,y)} \|\bfv\|.
\end{displaymath}
But, 
\begin{displaymath}
\lambda_{ij}(x,y) = \lambda_{ij}(x,y') + \lambda_{ij}(y',-t) \le -\kappa^{-1}t
\end{displaymath}
by (\ref{eq:lambdaij:zero:on:Fij}) and
Proposition~\ref{prop:properties:dynamical:norm:Hbig}. Therefore, 
\begin{displaymath}
d(R(x,y) \bfv, \bE_{i,j-1}(y)) \le e^{-\kappa^{-1} t} \|\bfv\| \le
e^{-\kappa^{-1} h} \|\bfv\|. 
\end{displaymath}
\qed\medskip

\bold{The bounded subspace.}
Fix $\index{$\theta$}\theta > 0$. (We will eventually choose $\theta$ sufficiently
small depending only on the dimension). 
\begin{definition}
\label{def:bounded:subspace}
Suppose $x \in \tilde{X}$. 
A vector $\bfv\in \bE_{ij}(x)$ is called {\em $(\theta,ij)$-bounded} if there
exists $C < \infty$ such that 
for all $\ell > 0$ and for ($1-\theta$)-fraction
of $y \in \cF_{ij}[x,\ell]$,
\begin{equation}
\label{eq:def:thetabounded}
\|R(x,y)\bfv \| \le C\|\bfv\|.
\end{equation}
\end{definition}

\bold{Remark.} From the definition and (\ref{eq:Rxy:bfv}), 
it is clear that every
vector in $\bE_{i1}(x)$ is $(\theta,i1)$-bounded for every
$\theta$. Indeed, we have $\bE_{i1}' = \bE_{i1}$, and $\lambda_{i1}(x,y)
= 0$ for $y \in \cF_{i1}[x]$, 
thus for $y \in \cF_{i1}[x]$ and $\bfv
\in E_{i1}(x)$, $\|R(x,y) \bfv \| = \| \bfv \|$.

\begin{lemma}
\label{lemma:bounded:subspace}
Let $n = \dim \bE_{ij}(x)$ (for a.e $x$). If there exists no non-zero
$\theta/n$-bounded vector in $\bE_{ij}(x) \setminus \bE_{i,j-1}(x)$,
we set $\bE_{ij,bdd} = \{0\}$. Otherwise, we define 
$\index{$E$@$\bE_{ij,bdd}(x)$}\bE_{ij,bdd}(x) \subset \bE_{ij}(x)$ to
be the linear span of the $\theta/n$-bounded vectors in $\bE_{ij}(x)$.
This is a subspace of $\bE_{ij}(x)$, and any vector in this
subspace is $\theta$-bounded. Also, 
\begin{itemize}
\item[{\rm (a)}] $\bE_{ij,bdd}(x)$ is $g_t$-equivariant, i.e.\ $(g_t)_*
  \bE_{ij,bdd}(x) = \bE_{ij,bdd}(g_t x)$. 
\item[{\rm (b)}] For almost all $u \in \cB(x)$, $\bE_{ij,bdd}(ux) =
  (u)_* \bE_{ij,bdd}(x)$. 
\end{itemize}
\end{lemma}

\bold{Proof.} 
Let $\bE_{ij,bdd}(x) \subset \bE_{ij}(x)$ denote the linear span of 
all $(\theta/n, ij)$-bounded vectors. 
If $\bfv_1, \dots, \bfv_n$ 
are any $n$ $(\theta/n,ij)$-bounded vectors, then there
exists $C > 1$ such that for $1-\theta$ fraction of $y$ in $\cF_{ij}[x,L]$,
(\ref{eq:def:thetabounded}) holds. But then
(\ref{eq:def:thetabounded}) holds (with a different $C$) for any
linear combination of the $\bfv_i$. This shows that any vector in
$\bE_{ij,bdd}(x)$ is $(\theta,ij)$-bounded. 
To show that (a) holds, suppose
that $\bfv \in \bE_{ij}(x)$ is $(\theta/n,ij)$-bounded, and $t < 0$. In
view of Lemma~\ref{lemma:Ej:equivariant}, it is enough to show that
$\bfv'\equiv (g_t^{ij})_* \bfv \in \bE_{ij}(g_t^{ij} x)$ is
$(\theta/n,ij)$-bounded. (This would show that for $t < 0$, $(g_t^{ij})_*
\bE_{ij,bdd}(x) \subset \bE_{ij,bdd}(g_t^{ij} x)$ which, in view of the
ergodicity of the action of $g_t$, would imply (a).)

\makefig{Proof of Lemma~\ref{lemma:bounded:subspace}
  (a).}{fig:foliations}{\includegraphics{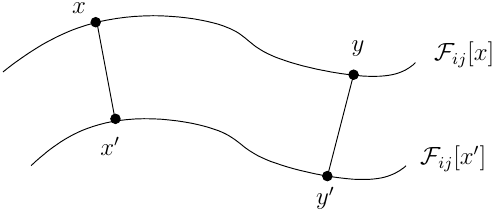}} 

Let $x' = g_t^{ij} x$. 
By (\ref{eq:hyperbolic:dynamical:norm:E}), 
there exists $C_1= C_1(t)$ such that for all $z \in
  X$ and all
$\bfw \in \bE(z)$,
\begin{equation}
\label{eq:temp:w}
C_1^{-1} \|\bfw\| \le \|(g^{ij}_t)_* \bfw\| \le C_1 \|\bfw\|. 
\end{equation}
Suppose $y \in \cF_{ij}[x,L]$ satisfies (\ref{eq:def:thetabounded}). Let
$y' = g_t^{ij} y$. 
Then $y' \in \cF_{ij}[x']$. 
Let  $\bfv' = (g^{ij}_t)_* \bfv$. (See Figure~\ref{fig:foliations}). 
Note that 
\begin{displaymath}
R(x',y')\bfv' = R(y,y')R(x,y)R(x',x)\bfv' = R(y,y')R(x,y)\bfv 
\end{displaymath}
hence by (\ref{eq:temp:w}), (\ref{eq:def:thetabounded}), and again
(\ref{eq:temp:w}), 
\begin{displaymath}
\|R(x',y')\bfv'\| \le  C_1 \|R(x,y) \bfv\| \le C_1 C \|\bfv\| \le
C_1^2 C \|\bfv'\|. 
\end{displaymath}
Hence, for $y \in \cF_{ij}[x,L]$ satisfying (\ref{eq:def:thetabounded}),
$y' = g^{ij}_t y \in \cF_{ij}[x']$ satisfies
\begin{equation}
\label{eq:yprime:satisfies:thetabounded}
\|R(x',y') \bfv'\| < C C_1^2 \|\bfv'\|.
\end{equation}
Therefore, since $\cF_{ij}[g^{ij}_t x, L+t] = g^{ij}_t \cF_{ij}[x,L]$,
we have that for $1-\theta/n$ fraction of  $y' \in \cF_{ij}[x', L+t]$, 
(\ref{eq:yprime:satisfies:thetabounded}) holds. Therefore, 
$\bfv'$ is $(\theta/n,ij)$-bounded. Thus, $\bE_{ij,bdd}(x)$ is
$g_t$-equivariant.  This completes the proof of (a). Then (b) follows
immediately from (a) since 
Lemma~\ref{lemma:symmetric:difference:balls} implies that $\cF_{ij}[
ux, L] = \cF_{ij}[x, L]$ for $L$ large enough. 
\qed\medskip


\begin{remark}
\label{remark:def:theta}
Formally, from its definition, the subspace
$\bE_{ij,bdd}(x)$ depends on the choice of $\theta$. It is clear that
as we decrease $\theta$, the subspace $\bE_{ij,bdd}(x)$ decreases. 
In view of Lemma~\ref{lemma:bounded:subspace}, there
exists $\theta_0 > 0$ and $m \ge 0$ 
such that for all $\theta < \theta_0$ and almost
all $x \in X$, the dimension of $\bE_{ij,bdd}(x)$ is $m$. We will
always choose $\index{$\theta$}\theta \ll \index{$\theta_0$}\theta_0$. 
\end{remark}
\medskip

\bold{Synchronized Exponents.} 
\begin{definition}
\label{def:synchronized}
Suppose $\theta > 0$. We say that $ij \in \Lambda''$ 
and $kr \in \Lambda''$  are
{\em $\theta$-synchronized} if there exists $E \subset X$ with $\nu(E) > 0$,
and $C < \infty$,  
such that for all $x \in \pi^{-1}(E)$, for all $\ell > 0$, 
for at least $(1-\theta)$-fraction of $y \in \cF_{ij}[x,\ell]$, we have
\begin{displaymath}
\rho(y,\cF_{kr}[x]) < C. 
\end{displaymath}
\end{definition}

\begin{remark}
\label{remark:sync:one}
By the same argument as in the proof of
Lemma~\ref{lemma:bounded:subspace} (a), if $ij$ and $kr$ are
$\theta$-synchronized then we can replace the set $E$ in
Definition~\ref{def:synchronized} by $\bigcup_{|s| < t} g_s E$. 
Therefore, we can take $E$ in Definition~\ref{def:synchronized}
to have measure arbitrarily close to $1$.
\end{remark}

\begin{remark}
\label{remark:sync:two}
Clearly if $ij$ and $kr$ are not $\theta$-synchronized,
then they are also not $\theta'$-synchronized for any $\theta' <
\theta$. Therefore there exists $\index{$\theta_0'$}\theta_0' > 0$ such that if any pairs
$ij$ and $kr$ are not $\theta$-synchronized for some $\theta > 0$ then they are
also not $\theta_0'$-synchronized. We will always consider $\index{$\theta$}\theta \ll
\theta_0'$, and will sometimes use the term ``synchronized'' with no
modifier to mean $\theta$-synchronized for $\theta \ll \theta_0'$. Then
in view of Remark~\ref{remark:sync:one}, synchronization is an
equivalence relation.  
\end{remark}
We now fix $\theta \ll \min(\theta_0,\theta_0')$.

\medskip

If $\bfv \in \bE(x)$, we can write
\begin{equation}
\label{eq:v:decomp:vij}
\bfv = \sum_{ij \in I_\bfv} \bfv_{ij}, \qquad\text{ where $\bfv_{ij} \in
  \bE_{ij}(x)$, but $\bfv_{ij} \not\in \bE_{i,j-1}(x)$}.
\end{equation}
In the sum, $I_\bfv$ is a finite set of pairs $ij$ where $i \in
\Lambda'$ and $1 \le j \le n_i$. (Recall that $\Lambda'$ denotes the
Lyapunov spectrum of $\bE$). Since for a fixed $i$ the
$\bE_{ij}(x)$ form a flag, without loss of generality we
may (and always will) assume that $I_\bfv$ contains at most one pair $ij$
for each $i \in \Lambda'$.

For $\bfv \in \bE(x)$, and $y \in \cF_\bfv[x]$, let 
\begin{displaymath}
\index{$H_\bfv(x,y)$}H_\bfv(x,y) = \sup_{ij \in I_\bfv} \rho(y,\cF_{ij}[x]).
\end{displaymath}

\begin{lemma}
\label{lemma:H:bounded}
There exists a set $\Psi \subset X$ with $\nu(\Psi) = 1$ 
such that the following holds: Suppose 
$x \in \Psi$, $C < \infty$,
and there exists $\bfv \in \bE(x)$ so that for each $L > 0$, 
for at least
$(1-\theta)$-fraction of $y \in \cF_\bfv[x,L]$ 
\begin{displaymath}
H_\bfv(x,y) < C. 
\end{displaymath}
Then, if we write $\bfv =
\sum_{ij \in I_\bfv} \bfv_{ij}$ as in (\ref{eq:v:decomp:vij}), 
then all $\{ ij \}_{ij \in I_\bfv}$ 
are synchronized, and also for all $ij \in I_\bfv$, $\bfv_{ij} \in
\bE_{ij,bdd}(x)$. 
\end{lemma}
\bold{Proof.} Let $\Psi = \bigcup_{t \in \reals} g_t E$, where $E$ is
as in Definition~\ref{def:synchronized}. (In view of
Remark~\ref{remark:sync:one}, we may assume that the same $E$ works
for all synchronized pairs).
Suppose $ij \in I_\bfv$ and $kr \in I_\bfv$. 
We have for at least
$(1-\theta)$-fraction of $y \in \cF_\bfv[x,L]$, 
\begin{displaymath}
\rho(y,\cF_{ij}[x]) < C, \qquad \rho(y,\cF_{kr}[x]) < C. 
\end{displaymath}
Let $y_{ij} \in \cF_{ij}[x]$ be such that $\rho(y,\cF_{ij}[x]) =
\rho(y,y_{ij})$. Similarly, let $y_{kr} \in \cF_{kr}[x]$ be such that
$\rho(y,\cF_{kr}[x]) = \rho(y,y_{kr})$. We have 
\begin{equation}
\label{eq:tmp:rho:yij:ykr}
\rho(y_{ij},y_{kr}) \le \rho(y_{ij},y) + \rho(y,y_{kr}) \le 2 C. 
\end{equation}
Note that $\tilde{g}^{\bfv,x}_{-L}(\cF_\bfv[x,L]) =
g^{ij}_{-L'}(\cF_{ij}[x,L'])$, where $L'$ is chosen so that
$g^{\bfv}_{-L} x = g^{ij}_{-L'} x$, where the notation $\tilde{g}$ is
as in (\ref{eq:def:tilde:g}). Hence, in view of
(\ref{eq:tmp:rho:yij:ykr}) and
(\ref{eq:def:lebesgue:leaf}), for any $L' > 0$, for
$(1-\theta)$-fraction of $y_{ij} \in \cF_{ij}[x,L']$, $\rho(y_{ij},
\cF_{kr}[x]) \le 2C$. 
Then, for any $t \in \reals$, for any $L'' > 0$, for
$(1-\theta)$-fraction of $y_{ij} \in \cF_{ij}[g_t x,L'']$,
$\rho(y_{ij}, \cF_{kr}[g_t x]) \le C(t)$. Since $x \in \Psi$, we can
choose $t$ so that $g_t x \in E$ where $E$ is as in
Definition~\ref{def:synchronized}.
This implies that $ij$ and $kr$ are synchronized. 

Recall that $I_\bfv$ contains at most one $j$ for each $i \in
\Lambda'$. Since $R(x,y)$ preserves each $\bE_i$, and the
distinct $\bE_i$ are orthogonal, for all $y'' \in G[x]$, 
\begin{displaymath}
\|R(x,y'') \bfv \|^2 = \sum_{ij \in I_\bfv} \|R(x,y'') \bfv_{ij} \|^2. 
\end{displaymath}
Therefore, for each $ij \in I_\bfv$,  and all $y'' \in G[x]$, 
\begin{displaymath}
\|R(x,y'') \bfv_{ij}\| \le  \|R(x,y'') \bfv\|. 
\end{displaymath}
In particular, 
\begin{displaymath}
\|R(x,y_{ij}) \bfv_{ij}\| \le  \|R(x,y_{ij}) \bfv\|. 
\end{displaymath}
We have for $(1-\theta)$-fraction of $y_{ij} \in \cF_{ij}[x,L']$,
$\rho(y_{ij},y) < C$, where $y \in \cF_{\bfv}(x)$.  
We have, by Lemma~\ref{lemma:Fv:norm},
$\|R(x,y) \bfv \| = \|\bfv\|$, and hence, by
(\ref{eq:hyperbolic:dynamical:norm:E}),  for
$(1-\theta)$-fraction of $y_{ij} \in \cF_{ij}[x,L]$, 
\begin{displaymath}
\| R(x,y_{ij}) \bfv \| \le C_2 \|\bfv\|. 
\end{displaymath}
Hence, for $(1-\theta)$-fraction of $y_{ij} \in \cF_{ij}[x,L']$, 
\begin{displaymath}
\| R(x,y_{ij}) \bfv_{ij} \| \le C_2 \|\bfv\|. 
\end{displaymath}
This implies that $\bfv_{ij} \in \bE_{ij,bdd}(x)$.
\qed\medskip

We write \index{$ij \sim kr$}$ij \sim kr$ if $ij$ and $kr$ are synchronized. With our
choice of $\theta > 0$, synchronization is an equivalence
relation, see Remark~\ref{remark:sync:two}.
We write $\index{$ij$@$[ij]$}[ij] = \{kr \st kr \sim ij \}$. Let
\begin{displaymath}
\bE_{[ij],bdd}(x) = \sum_{kr \in [ij]} \bE_{kr,bdd}(x). 
\end{displaymath}

For $\bfv \in \bE(x)$, write $\bfv = \sum_{ij \in I_\bfv} \bfv_{ij}$, as in
(\ref{eq:v:decomp:vij}). 
Define
\begin{displaymath}
\index{$height$@$\operatorname{height}(\bfv)$}\operatorname{height}(\bfv) = \sum_{ij \in I_\bfv} (\dim \bE) + j
\end{displaymath}

The height is defined so it would have the following properties:
\begin{itemize}
\item If $\bfv \in \bE_{ij}(x) \setminus \bE_{i,j-1}(x)$
  and $\bfw \in \bE_{i,j-1}(x)$ then $\operatorname{height}(\bfw) <
  \operatorname{height}(\bfv)$. 
\item If $\bfv = \sum_{i \in I_\bfv} \bfv_i$, $\bfv_i \in \bE_i$, $\bfv_i
  \ne 0$, and $\bfw = \sum_{j \in J} \bfw_j$, $\bfw_j \in \bE_j$, $\bfw_j
  \ne 0$, and also the cardinality of $J$ is smaller then the
  cardinality of $I_\bfv$, then $\operatorname{height}(\bfw) <
  \operatorname{height}(\bfv)$. 
\end{itemize}


Let $\index{$P$@$\cP_k(x)$}\cP_k(x) \subset \bE(x)$ denote the set of 
vectors of height at most $k$. This is a closed subset of $\bE(x)$. 

\begin{lemma}
\label{lemma:push:to:next}
For every $\delta > 0$ and every $\eta  > 0$ 
there exists a subset $K \subset X$
of measure at least $1-\delta$ and 
$L'' > 0$ such that for any $x \in K$ 
and any unit vector $\bfv \in \cP_k(x)$ with 
$d(\bfv, \bigcup_{ij} \bE_{[ij],bdd}) > \eta$ and $d(\bfv,\cP_{k-1}(x)) >
\eta$, there exists $0 < L' < L''$ 
so that for at least $\theta$-fraction of $y \in \cF_\bfv[x,L']$,
\begin{displaymath}
d\left(\frac{R(x,y) \bfv}{\|R(x,y) \bfv\|}, \cP_{k-1}(y)\right) < \eta.
\end{displaymath}
\end{lemma}

\bold{Proof.} 
Suppose $C>1$ (we will later choose $C$ depending on $\eta$). 
We first claim that we can choose $K$ with $\nu(K) > 1-\delta$ 
and $L'' > 0$ so that for every $x \in g_{[-1,1]}K$ 
and every $\bfv \in
\cP_k(x)$ such that $d(\bfv, \bigcup_{ij} \bE_{[ij],bdd}) > \eta$
there exists 
$0 < L' < L''$ so that 
for $\theta$-fraction of $y \in \cF_\bfv[x,L']$,
\begin{equation}
\label{eq:temp:Ru:geC}
H_\bfv(x,y) \ge C.
\end{equation}
(Essentially, this follows from Lemma~\ref{lemma:H:bounded}, but the
argument given below is a bit more elaborate since we want to choose $L''$
uniformly over all $\bfv \in \cP_{k}(x)$ satisfying $d(\bfv,
\bigcup_{ij} \bE_{[ij],bdd}) > \eta$). 
Indeed, let $E_L \subset \cP_k(x)$ 
denote the set of unit vectors $\bfv \in \cP_k(x)$ such
that for all $0 < L' < L$, 
for at least $(1-\theta)$-fraction of  
$y \in \cF_\bfv[x,L']$, 
$H_\bfv(x,y) \le C$. Then, the $E_L$ are closed sets which are
decreasing as $L$ increases, and by Lemma~\ref{lemma:H:bounded}, 
\begin{displaymath}
\bigcap_{L =1}^\infty E_L \subset \left(\bigcup_{ij \in
      \tilde{\Lambda}} \bE_{[ij],bdd}(x) \right) \cap \cP_k(x).  
\end{displaymath}
Let $F$ denote the subset of the unit sphere in $\cP_k(x)$ which is the
complement of the $\eta$-neighborhood of $\bigcup_{ij}
\bE_{[ij],bdd}(x)$. 
Then the $E_L^c$ are an open cover of $F$, and since $F$ is compact, there
exists $L= L_x$ such that $F \subset E_L^c$. Now for any $\delta > 0$
we can choose $L''$ so that $L'' > L_x$ for all $x$ in a set $K$ of
measure at least $(1-\delta)$. 

Now suppose $\bfv \in F$. Since $F \subset E_{L''}^c$, $\bfv \not\in
E_{L''}$, hence there exists $0 < L' < L''$ 
(possibly depending on $\bfv$) such that 
the fraction of $y \in \cF_\bfv[x,L']$ which satisfies
$H_\bfv(x,y) \ge C$ is greater than $\theta$. Then, 
(\ref{eq:temp:Ru:geC}) holds. 

Now suppose (\ref{eq:temp:Ru:geC}) holds (with a yet to be chosen $C =
C(\eta)$). Write
\begin{displaymath}
\bfv = \sum_{ij \in I_\bfv} \bfv_{ij}
\end{displaymath}
as in (\ref{eq:v:decomp:vij}). Let
\begin{displaymath}
\bfw  = R(x,y) \bfv, \qquad \bfw_{ij} = R(x,y) \bfv_{ij}. 
\end{displaymath}
Since $y \in \cF_\bfv[x]$, by Lemma~\ref{lemma:Fv:norm},
$\|\bfw\| = \|\bfv\| = 1$. 
Let $ij \in I_\bfv$ be such that the supremum in the definition of $H_\bfv(x,y)$
is achieved for $ij$. 
If $\|\bfw_{ij}\| < \eta/2$ we are done, since $\bfw' = \sum_{kr\ne
  ij} \bfw_{kr}$ has smaller height than $\bfv$, 
and $d(\bfw,\frac{\bfw'}{\|\bfw'\|}) < \eta$. Hence we may
assume that $1 \ge \|\bfw_{ij}\| \ge \eta/2$.  

Since $d(\bfv, \cP_{k-1}(x)) \ge \eta$, we have
\begin{displaymath}
d(\bfv_{ij},\bE_{i,j-1}(x)) \ge \eta \ge \eta \|\bfv_{ij}\|.
\end{displaymath}
where the last inequality follows from the fact that $\|\bfv_{ij}\|
\le 1$.  In particular, we have $1 \ge \|\bfv_{ij}\| \ge \eta$.  

Let $y' = g_t y$ be such that $y' \in \cF_{\bfv_{ij}}[x]$. Note that
\begin{displaymath}
1 = \|R(x,y') \bfv_{ij} \| = \|R(y,y') \bfw_{ij} \| = \|(g_t)_*
\bfw_{ij} \|
\qquad\text{and}\qquad 1 \ge \|\bfw_{ij} \| \ge \eta/2. 
\end{displaymath}
Then, in
view of (\ref{eq:hyperbolic:dynamical:norm:E}), $|t| \le
C_0(\eta)$, and hence $\|R(y',y)\| \le C_0'(\eta)$. 

Let $C_1 = C_0(\eta)+ h(\eta,\frac{1}{2}\eta/C_0'(\eta))$, where $h(\cdot,\cdot)$ is as in
Lemma~\ref{lemma:bounded:implies:lower:rank}. We now choose the
constant $C$ in (\ref{eq:temp:Ru:geC}) to be $C_1$. 
If $H_\bfv(x,y) > C_1$ then, by the choice of $ij$, $\rho(y,\cF_{ij}[x])>
C_1$. Since $y' = g_t y$ and $|t| \le C_0(\eta)$, we
have 
\begin{displaymath}
\rho(y',\cF_{ij}[x]) > C_1 - C_0(\eta) = h(\eta,\tfrac{1}{2}\eta/C_0'(\eta)).
\end{displaymath}
Then,
by Lemma~\ref{lemma:bounded:implies:lower:rank} applied to $\bfv_{ij}$
and $y' \in \cF_{\bfv_{ij}}[x]$, 
\begin{displaymath}
d(R(x,y') \bfv_{ij}, \bE_{i,j-1}(y')) \le \tfrac{1}{2}(\eta/C_0'(\eta))
\|\bfv_{ij}\| \le \tfrac{1}{2}\eta/C_0'(\eta). 
\end{displaymath}
Then, since $\bfw_{ij} = R(y',y) R(x,y') \bfv_{ij}$, 
\begin{displaymath}
\|d(\bfw_{ij}, \bE_{i,j-1}(y)) \| \le \|R(y',y)\| d(R(x,y') \bfv_{ij},
\bE_{i,j-1}(y')) \le \|R(y',y)\| (\eta/C_0'(\eta)) \le \frac{\eta}{2}. 
\end{displaymath}
Let $\bfw'_{ij}$ be the closest vector to $\bfw_{ij}$ in
$\bE_{i,j-1}(y)$, and let $\bfw' = \bfw'_{ij} + \sum_{kr\ne ij}
\bfw_{ij}$. Then $d(\bfw,\frac{\bfw'}{\|\bfw'\|}) < \eta$ and $\bfw' \in \cP_{k-1}$. 
\qed\medskip

\bold{Proof of Proposition~\ref{prop:some:fraction:bounded}.} 
Let $n$ denote the maximal possible height of a
vector. Let $\delta' = \delta/n$. 
Let $\eta_n = \eta$. Let $L_{n-1}= L_{n-1}(\delta',\eta_n)$ and
$K_{n-1} = K_{n-1}(\delta',\eta_n)$ be chosen so that
Lemma~\ref{lemma:push:to:next} holds for $k=n-1$, $K = K_{n-1}$, 
$L'' = L_{n-1}$ and $\eta =
\eta_n$. Let $\eta_{n-1}$ be chosen so that 
$\exp(N (L_{n-1}+1))\eta_{n-1} \le \eta_n$, where $N$ is as in
Lemma~\ref{lemma:forni2}. We repeat this
process until we choose $L_1$, $\eta_0$. Let $L_0 = L_1+1$. 
Let $K = K_0 \cap \dots \cap K_{n-1}$. Then $\nu(K) >
1-\delta$. 

Let 
\begin{displaymath}
E_k' = \left\{ y \in \cF_\bfv[x,L] \st d\left(\frac{R(x,y)\bfv}
{\|R(x,y) \bfv\|}, \cP_k(y) \cup \bigcup_{ij \in \tilde{\Lambda}}
  \bE_{[ij],bdd}(y) \right) < \eta_k \right\}.
\end{displaymath}
and let 
\begin{displaymath}
E_k = \tilde{g}_{-L}(E_k'),
\end{displaymath}
so $E_k \subset \cB[z]$, where $z =  \tilde{g}_{-L} x$. Since $E_n' =
\cF_\bfv[x,L]$, we have $E_n = \cB[z]$. 
Let $Q = \tilde{g}_{-L}(g_{[-1,1]}K \cap \cF_\bfv[x,L])$. Then, by assumption, 
\begin{equation}
\label{eq:Q:lower:bound}
|Q| \ge (1 - (\theta/2)^{n+1})|\cB[z]|.
\end{equation}
By Lemma~\ref{lemma:push:to:next}, for every
point $uz \in (E_k \cap Q) \setminus E_{k-1}$ 
there exists a ``ball''
$\cB_t[uz]$ (where $t = L -L'$ and $L'$ is as in
Lemma~\ref{lemma:push:to:next}) such that 
\begin{equation}
\label{eq:balls:cover}
|E_{k-1} \cap \cB_t[uz]| \ge \theta |\cB_t[uz]|.
\end{equation}
(When we are applying Lemma~\ref{lemma:push:to:next} we do not have $\bfv
\in \cP_k$ but rather $d(\bfv/\|\bfv\|, \cP_k) < \eta_k$; however by the
choice of the $\eta$'s and the $L$'s this does not
matter). \mcc{explain better}
The collection of balls $\{\cB_t[uz]\}_{uz \in (E_k \cap Q) 
\setminus E_{k-1}}$ as
in (\ref{eq:balls:cover})
are a cover of $(E_k \cap Q)\setminus E_{k-1}$. These balls satisfy
the condition of Lemma~\ref{lemma:gB:properties} (b); hence we may
choose a pairwise disjoint subcollection which still covers 
$(E_k \cap Q)\setminus E_{k-1}$. 
We get $|E_{k-1}| \ge \theta |E_k \cap Q|$. Hence, by (\ref{eq:Q:lower:bound})
and induction over $k$, we have
\begin{displaymath}
|E_{k}| \ge (\theta/2)^{n-k} |\cB[z]|.
\end{displaymath}
Hence, $|E_0| \ge (\theta/2)^n |\cB[z]|$. Therefore $|E_0'| \ge (\theta/2)^n
|\cF_\bfv[x,L]|$. Since $\cP_0 = \emptyset$, the Proposition follows
from the definition of $E_0'$. 
\qed\medskip

\starredsubsection{Invariant measures on $X \cross \proj(\bL)$.}
\label{sec:starredsubsec:invariant:measures:proj}

In this subsection we prove
Proposition~\ref{prop:ej:bdd:transport:bounded}. 

Recall that any bundle is measurably trivial. 
\begin{lemma}
\label{lemma:proj:general}
Suppose $\bL(x)$   is an invariant subbundle or quotient
bundle of $\bH(x)$. 
(In fact the arguments in this subsection apply to arbitrary
vector bundles).  
Let \index{$\mu$@$\tilde{\mu}_\ell$}$\tilde{\mu}_\ell$ be the measure on $X \cross \proj(\bL)$ defined by
\begin{equation}
\label{eq:def:tilde:mu}
\tilde{\mu}_\ell(f) = \int_X \int_{\proj(\bL)}
\frac{1}{|\cF_{ij}[x,\ell]|} \int_{\cF_{ij}[x,\ell]}
f(x,R(y,x)\bfv) \, dy \, d\rho_0(\bfv) \, d\nu(x)
\end{equation}
where $\rho_0$ is the ``round'' measure on $\proj(\bL)$. 
(In fact, $\rho_0$ can be any measure on $\proj(\bL)$ in the measure class of
Lebesgue measure, independent of x and fixed once and for all). 
Let \index{$\mu$@$\hat{\mu}_\ell$}$\hat{\mu}_\ell$ be
the measure on $X \cross \proj(\bL)$ defined by
\begin{equation}
\label{eq:def:hat:mu}
\hat{\mu}_\ell(f) = \int_X \int_{\proj(\bL)}
\frac{1}{|\cF_{ij}[x,\ell]|} \int_{\cF_{ij}[x,\ell]}
f(y,R(x,y)\bfv) \, dy \, d\rho_0(\bfv) \, d\nu(x).
\end{equation}
Then $\hat{\mu}_\ell$ is in the same measure class as $\tilde{\mu}_\ell$, and
\begin{equation}
\label{eq:hat:nu:tilde:nu}
\kappa^{-2} \le \frac{d\hat{\mu}_\ell}{d\tilde{\mu}_\ell} \le \kappa^2,
\end{equation}
where $\kappa$ is as in
Proposition~\ref{prop:properties:dynamical:norm:Hbig}. 
\end{lemma}
\bold{Proof.} Let
\begin{displaymath}
F(x,y) = \int_{\proj(\bL)} f(x,R(y,x) \bfv) \, d\rho_0(\bfv). 
\end{displaymath}
Then, 
\begin{equation}
\label{eq:tilde:mu:ell:f}
\tilde{\mu}_\ell(f) = \int_{X}  \frac{1}{|\cF_{ij}[x,\ell]|}
\int_{\cF_{ij}[x,\ell]} F(x,y) \, dy \, d\nu(x)
\end{equation}
\begin{equation}
\label{eq:hat:mu:ell:f}
\hat{\mu}_\ell(f) = \int_{X} \frac{1}{|\cF_{ij}[x,\ell]|}
\int_{\cF_{ij}[x,\ell]} F(y,x) \, dy \, d\nu(x)
\end{equation}
Let $x' = g^{ij}_{-\ell} x$. Then, in view of
Proposition~\ref{prop:properties:dynamical:norm:Hbig}, 
$\kappa^{-1} \, d\nu(x)
\le d\nu(x') \le \kappa \, d\nu(x)$. Then, 
\begin{displaymath}
\frac{1}{\kappa} \tilde{\mu}_\ell(f) \le \int_{X}  \frac{1}{|\cB[x']|}
\int_{\cB[x']} F(g^{ij}_\ell x', g^{ij}_\ell z) \, dz \, d\nu(x') \le \kappa \tilde{\mu}_\ell(f),
\end{displaymath}
and
\begin{displaymath}
\frac{1}{\kappa} \hat{\mu}_\ell(f) \le \int_{X}  \frac{1}{|\cB[x']|}
\int_{\cB[x']} F(g^{ij}_\ell z, g^{ij}_\ell x') \, dz \, d\nu(x') \le \kappa \hat{\mu}_\ell(f)
\end{displaymath}
Let $X''$ consist of one point from each $\cB[x]$. In view of
Definition~\ref{def:compatible:family} (iii),
we now disintegrate $d\nu(x') = d\beta(x'') dz'$ where $x'' \in
X''$, $z' \in \cB[x']$. 
\begin{align*}
\int_{X}  \frac{1}{|\cB[x']|}
\int_{\cB[x']} F(g^{ij}_\ell x', g^{ij}_\ell z) \, dz \, d\nu(x') & = \int_{X''}
\int_{\cB[x''] \cross \cB[x'']} F(g^{ij}_\ell z', g^{ij}_\ell z) \, dz' \, dz \,
d\beta(x'') \\
& = \int_{X''}
\int_{\cB[x''] \cross \cB[x'']} F(g^{ij}_\ell z, g^{ij}_\ell z') \, dz' \, dz \,
d\beta(x'') \\
& = \int_{X}  \frac{1}{|\cB[x']|}
\int_{\cB[x']} F(g^{ij}_\ell z, g^{ij}_\ell x') \, dz \, d\nu(x'). 
\end{align*}
Now (\ref{eq:hat:nu:tilde:nu}) follows from (\ref{eq:tilde:mu:ell:f}) and
(\ref{eq:hat:mu:ell:f}). 
\qed\medskip

\begin{lemma}
\label{lemma:lambda:x:invariant}
Let $\tilde{\mu}_\infty$ be any weak-star limit of the measures
$\tilde{\mu}_\ell$. Then, 
\begin{itemize}
\item[{\rm (a)}] We may disintegrate $d\tilde{\mu}_\infty(x,\bfv) =
  d\nu(x) \, d\lambda_x(\bfv)$, where for each $x \in X$,
  $\lambda_x$ is   a measure on $\proj(\bL)$.  
\item[{\rm (b)}] For $x \in \tilde{X}$ and $y \in \cF_{ij}[x]$, 
\begin{displaymath}
\lambda_y = R(x,y)_* \lambda_x,
\end{displaymath}
(where to simplify notation, we write $\lambda_x$ and $\lambda_y$
instead of $\lambda_{\pi(x)}$ and $\lambda_{\pi(y)}$). 
\item[{\rm (c)}] Let $\bfw \in \proj(\bL)$ be a point. For $\eta > 0$ let  
\begin{displaymath}
B(\bfw,\eta) = \{ \bfv \in \proj(\bL) \st d(\bfv, \bfw) \le \eta\}.
\end{displaymath}
Then, for any $t < 0$ there exists
  $c_1 = c_1(t,\bfw) > 0$ and $c_2 = c_2(t, \bfw) > 0$ such that for
  $x \in X$, 
\begin{displaymath}
\lambda_{g_t x}(B(g_t \bfw,c_1\eta)) \ge c_2
\lambda_x(B(\bfw,\eta)).
\end{displaymath}
Consequently, for $t < 0$, the support of $\lambda_{g_t x}$
contains the support of $(g_t)_*\lambda_x$. 
\item[{\rm (d)}] For almost all $x \in X$ there exist a measure
  $\psi_x$ on $\proj(\bL)$ such that 
\begin{displaymath}
\lambda_x = h(x) \psi_x 
\end{displaymath}
for some $h(x) \in SL(\bL)$, and also for almost all $y \in \cF_{ij}[x]$,
$\psi_y = \psi_x$ (so that $\psi$ is constant on the leaves 
$\cF_{ij}$). The maps $x \to \psi_x$ and $x \to h(x)$ are both
$\nu$-measurable. 
\end{itemize}
\end{lemma}

\bold{Proof.} If $f(x,\bfv)$ is independent of the second variable, then
it is clear from the definition of $\tilde{\mu}_\ell$ that
$\tilde{\mu}_\ell(f) = \int_X f \, d\nu$. This implies (a). 
To prove (b), note that $R(y',y) = R(x,y) R(y',x)$. Then, 
\begin{align*}
\lambda_y & = \lim_{k \to \infty}
\frac{1}{|\cF_{ij}[y,\ell_k]|} \int_{\cF_{ij}[y,\ell_k]} \left(R(y',y)_*
  \rho_0 \right) \, dy'  \\
&= R(x,y)_* \lim_{k \to \infty}
\frac{1}{|\cF_{ij}[y,\ell_k]|} \int_{\cF_{ij}[y,\ell_k]} \left(R(y',x)_*
  \rho_0 \right) \, dy'  \\
&= R(x,y)_* \lim_{k \to \infty} 
\frac{1}{|\cF_{ij}[x,\ell_k]|} \int_{\cF_{ij}[x,\ell_k]} \left(R(y',x)_*
  \rho_0 \right) \, dy' \\
&=R(x,y)_* \lambda_x
\end{align*}
where to pass from the second line to the third we used the fact that
$\cF_{ij}[x,\ell] = \cF_{ij}[y,\ell]$ for $\ell$ large enough. This
completes the proof of (b). 

We now begin the proof of (c). Let $\bfw(x) = \bfw$. Working in the
universal cover, we define for $y \in G[x]$, 
$\bfw(y) = R(x,y) \bfw(x)$. We define
\begin{displaymath}
\bfw_\eta(x) = \{ \bfv \in \proj(\bL(x)) \st d(\bfv, \bfw(x)) \le \eta\}.
\end{displaymath}
(Here we are thinking of
the space as $X \cross \proj(\bL)$ and using the same metric on all
the $\proj(\bL)$ fibers).

Let $x' = g_t^{ij} x$, $y' = g_t^{ij} y$.  We have
\begin{displaymath}
R(y',x') = R(x,x') R(y,x) R(y',y). 
\end{displaymath}
Since $\|R(x,x')^{-1}\| \le c^{-1}$, where $c$ depends on $t$, we have
$R(x,x')^{-1} \bfw_{c\eta}(x') \subset \bfw_\eta(x)$.  Then,  
\begin{align*}
\rho_0 \{ \bfv \st R(y',x') \bfv \in \bfw_{c\eta}(x') \} & = \rho_0 \{ \bfv \st
R(y,x) R(y',y) \bfv \in R(x,x')^{-1} \bfw_{c\eta}(x') \} \\
& \ge \rho_0 \{ \bfv \st R(y,x) R(y',y) \bfv \in \bfw_\eta(x) \} \\
& =  \rho_0 \{ R(y,y')^{-1} \bfu \st R(y,x) \bfu  \in \bfw_\eta(x) \} \\
& =  R(y,y')^{-1}_* \rho_0 \{ \bfu \st R(y,x) \bfu  \in \bfw_\eta(x) \} \\
& \ge c' \rho_0 \{ \bfu \st R(y,x) \bfu  \in \bfw_\eta(x) \}. 
\end{align*}
Note that for $t < 0$, $g^{ij}_t\cF_{ij}[x, \ell] \subset
\cF_{ij}[g^{ij}_t x,\ell]$ and $|g^{ij}_t \cF_{ij}[x, \ell]| \ge
c(t)|\cF_{ij}[g^{ij}_t x,\ell]|$.
Substituting into (\ref{eq:def:tilde:mu}) completes the
proof of (c). \mcc{more details?}

To prove part (d), let $\cM$ denote the space of measures on
$\proj(\bL)$. Recall that by \cite[Theorem 3.2.6]{ZimmerBook} 
the orbits of the
special linear group $SL(\bL)$ on
$\cM$ are locally closed. 
Then, by \cite[Theorem 2.9 (13), Theorem 2.6(5)]{Effros} \footnote{The
  ``condition C'' of \cite{Effros} is satisfied since 
$SL(\bL)$ is locally compact and $\cM$ is Hausdorff.}
there exists a
Borel cross section $\phi: \cM/SL(\bL) \to \cM$. Then, let $\psi_x =
\phi(\pi(\lambda_x))$ where $\pi: \cM \to \cM/SL(\bL)$ is the quotient
map.  
\qed\medskip

We also recall the following well known Lemma of Furstenberg (see
e.g.\ \cite[Lemma 3.2.1]{ZimmerBook}):
\begin{lemma}
\label{lemma:furstenberg:two:subspace}
Let $\bL$ be a vector space, and suppose $\mu$ and $\nu$ are two
probability measures on $\proj(\bL)$. Suppose $g_i \in SL(\bL)$ are such
that $g_i \to \infty$ and $g_i \mu \to \nu$. Then the support of $\nu$
is contained in a union of two proper subspaces of $\bL$. 

In particular, if the support of a measure $\nu$ on $\proj(\bL)$ 
is not contained in a union of
two proper subspaces, then the stabilizer of $\nu$ in $SL(\bL)$ is
bounded. 
\end{lemma}

\begin{lemma}
\label{lemma:switching:Rxy:Ryx}
Suppose $\bL$ is either a subbundle or a quotient bundle of
$\bH$. Suppose that $\theta > 0$, and 
suppose that for all $\delta > 0$ there exists a set $K \subset X$
with $\nu(K) > 1-\delta$ and a constant $C_1 < \infty$, such that
such that for all $x \in K$, 
all $\ell > 0$ and at least $(1-\theta)$-fraction of $y \in
\cF_{ij}[x,\ell]$, 
\begin{equation}
\label{eq:tmp:Rxy:v:small}
\|R(x,y) \bfv\| \le C_1 \|\bfv \| \quad\text{ for all $\bfv \in \bL$.}
\end{equation}
Then for all $\delta > 0$ and
for all $\ell > 0$ there exists a subset $K''(\ell) \subset X$
with $\nu(K''(\ell)) > 1-c(\delta)$ where $c(\delta) \to 0$ as
$\delta \to 0$, and there exists $\theta'' = \theta''(\theta,\delta)$ with
$\theta'' \to 0$ as $\theta \to 0$ and $\delta \to 0$
such that for all $x \in
K''(\ell)$, for at least $(1-\theta'')$-fraction of $y \in
\cF_{ij}[x,\ell]$, 
\begin{equation}
\label{eq:tmp:Rxy:v:twosided:small}
C_1^{-1} \|\bfv\| \le \|R(x,y) \bfv\| \le C_1 \|\bfv \| \quad\text{
  for all $\bfv \in \bL$}.
\end{equation}
\end{lemma}

\bold{Proof.} Let $f$ be the characteristic
function of $K \cross \proj(\bL)$. By (\ref{eq:def:tilde:mu}),
$\tilde{\mu}_\ell(f) \ge (1-\delta)$. 
By Lemma~\ref{lemma:proj:general} we have $\hat{\mu}_\ell(f) \ge
(1-\kappa^2 \delta)$. Therefore, by (\ref{eq:def:hat:mu}) 
there exists a subset $K'(\ell) \subset
X$ with $\nu(K'(\ell)) \ge 1-(\kappa^2 \delta)^{1/2}$ such that
such that for all $x \in K'(\ell)$, 
\begin{displaymath}
|\cF_{ij}[x,\ell] \cap K| \ge (1-(\kappa^2 \delta)^{1/2})
|\cF_{ij}[x,\ell]|. 
\end{displaymath}
For $x_0 \in X$, let
\begin{displaymath}
Z_\ell[x_0] = \{ (x,y) \in \cF_{ij}[x_0,\ell] \cross \cF_{ij}[x_0,\ell] \st x
\in K, \quad y \in K, \quad \text{ and (\ref{eq:tmp:Rxy:v:small})
  holds } \}.
\end{displaymath}
Then, if $x_0 \in K'(\ell)$ and $\theta' = \theta + (\kappa^2
\delta)^{1/2}$ then, by Fubini's theorem, 
\begin{displaymath}
|Z_\ell[x_0]| \ge (1-\theta') |\cF_{ij}[x_0,\ell] \cross
\cF_{ij}[x_0,\ell]|. 
\end{displaymath}
Let 
\begin{displaymath}
Z_\ell[x_0]^t = \{ (x,y) \in \cF_{ij}[x_0,\ell] \cross \cF_{ij}[x_0,\ell]
\st (y,x) \in Z_\ell[x_0] \}.
\end{displaymath}
Then, for $x_0 \in K'(\ell)$, 
\begin{displaymath}
|Z_\ell[x_0] \cap Z_\ell[x_0]^t| \ge (1-2 \theta') |\cF_{ij}[x_0,\ell]
\cross \cF_{ij}[x_0,\ell]|. 
\end{displaymath}
For $x \in \cF_{ij}[x_0,\ell]$, let 
\begin{displaymath}
Y_\ell'(x) = \{ y \in \cF_{ij}[x,\ell] \st (x,y) \in Z_\ell[x] \cap
Z_\ell[x]^t \}. 
\end{displaymath}
Therefore, by Fubini's theorem, for all $x_0 \in K'(\ell)$ and
$\theta'' = (2\theta')^{1/2}$,  
\begin{equation}
\label{eq:tmp:fubini:cFij:x0:ell}
|\{ x \in \cF_{ij}[x_0,\ell] \st | Y_\ell'(x)| \ge (1-\theta'')
|\cF_{ij}[x_0,\ell]| \}| \ge (1-\theta'') |\cF_{ij}[x_0,\ell]|.
\end{equation}
(Note that $\cF_{ij}[x_0,\ell] = \cF_{ij}[x,\ell]$.) Let 
\begin{displaymath}
K''(\ell) = \{ x \in X \st | Y_\ell'(x)| \ge (1-\theta'')
|\cF_{ij}[x,\ell]| \}. 
\end{displaymath}
Therefore, by (\ref{eq:tmp:fubini:cFij:x0:ell}), 
for all $x_0 \in K'(\ell)$, 
\begin{displaymath}
|\cF_{ij}[x_0,\ell] \cap K''(\ell)| \ge (1-\theta'') |\cF_{ij}[x_0,\ell]|.
\end{displaymath}
Then, by the definition of $\hat{\mu}_\ell$, 
\begin{displaymath}
\hat{\mu}_\ell( K''(\ell) \cross \proj(\bL)) \ge (1-\theta'')
\nu(K'(\ell)) \ge (1-2\theta''), 
\end{displaymath}
and therefore, by Lemma~\ref{lemma:proj:general}, 
\begin{displaymath}
\nu(K''(\ell)) = \tilde{\mu}_\ell(K''(\ell) \cross \proj(\bL)) \ge
(1-2\kappa^2 \theta''). 
\end{displaymath}
Now, for $x \in K''(\ell)$, and $y \in Y_\ell'(x)$,
(\ref{eq:tmp:Rxy:v:twosided:small}) holds. 
\qed\medskip

\begin{lemma}
\label{lemma:Vbdd:undependent:u}
Suppose $\bL(x) = \bE_{ij,bdd}(x)$. Then
there exists a $\Gamma$-invariant function $C: \tilde{X} 
\to \reals^+$ finite almost
everywhere such that for all $x \in \tilde{X}$, 
all $\bfv \in \bL(x)$, 
and all $y \in \cF_{ij}[x]$, 
\begin{displaymath}
C(x)^{-1}C(y)^{-1}\|\bfv\| \le \|R(x,y) \bfv\| \le C(x)C(y)\|\bfv\|,
\end{displaymath}
\end{lemma}

\bold{Proof.} Let $\tilde{\mu}_\ell$ and $\hat{\mu}_\ell$ be as in
Lemma~\ref{lemma:proj:general}. Take a sequence $\ell_k \to \infty$ 
such that $\tilde{\mu}_{\ell_k} \to \tilde{\mu}_\infty$, and
$\hat{\mu}_{\ell_k} 
\to \hat{\mu}_\infty$. Then by Lemma~\ref{lemma:lambda:x:invariant} (a),
we have $d\tilde{\mu}_\infty(x,\bfv) = d\nu(x) \, d\lambda_x(\bfv)$ where
$\lambda_x$ is a measure on $\proj(\bL)$. Let $E \subset X$ be such
that for $x \in E$, $\lambda_x$ is supported on at most two subspaces.
We will show that $\nu(E) = 0$. 

Suppose not;  
then $\nu(E) > 0$, and for $x \in E$, $\lambda_x$ is
supported on $\bF_1(x) \cup \bF_2(x)$, where $\bF_1(x)$ and $\bF_2(x)$ are
subspaces of $\bL(x)$. We always choose $\bF_1(x)$ and $\bF_2(x)$ to be of
minimal dimension, and if $\lambda_x$ is supported on a single
subspace $\bF(x)$ (of minimal dimension), we let $\bF_1(x) = \bF_2(x)
= \bF(x)$. Then, for $x\in E$, $\bF_1(x) \cup \bF_2(x)$ is uniquely
determined by $x$. After possibly replacing $E$ by a smaller subset of
positive measure, we may assume that $\dim \bF_1(x)$ and $\dim
\bF_2(x)$ are independent of $x \in E$.

Let 
\begin{displaymath}
\Psi = \{ x \in X \st g_t x  \in E \text{ and } g_{-s} x \in E \text{
  for some $t > 0$ and $s > 0$.}\}
\end{displaymath}
Then, $\nu(\Psi) = 1$. If $x \in \Psi$, then, by
Lemma~\ref{lemma:lambda:x:invariant} (c), 
\begin{equation}
\label{eq:tmp:gt:union:subspaces}
(g_{s})_* \bF_1(g_{-s} x) \cup (g_s)_* \bF_2(g_{-s} x) \subset \operatorname{supp} \lambda_x \subset (g_{-t})_* \bF_1(g_t x) \cup (g_{-t})_* \bF_2(g_t x),
\end{equation}
Since $\bF_i(g_t x)$ and $\bF_i(g_{-s} x)$ have the same dimension,
the sets on the right and on the left of (\ref{eq:tmp:gt:union:subspaces})
coincide. Therefore, $E \supset \Psi$ (and so $E$ has
full measure) and 
the set $\bF_1(x) \cup \bF_2(x)$ is $g_t$-invariant. By
Proposition~\ref{prop:sublyapunov:new:locally:constant} (see also the
remark immediately following the Proposition)
the set $\bF_1(x) \cup \bF_2(x)$ is also
$U^+$-invariant.

Fix $\delta > 0$ (which will be chosen sufficiently small later). 
Suppose $\ell > 0$ is arbitrary. 
Since $\bL = \bE_{ij,bdd}$, 
there exists a constant $C_1$ independent of $\ell$
and a compact subset $K
\subset X$ with $\nu(K) > 1-\delta$ 
and for each $x \in K$ a subset $Y_\ell(x)$
of $\cF_{ij}[x,\ell]$ with $|Y_\ell(x)| \ge (1-\theta)
|\cF_{ij}[x,\ell]|$, such that for $x \in K$
and $y \in Y_\ell(x) \cap K$ we have
\begin{displaymath}
\|R(x,y) \bfv\| \le C_1 \|\bfv \| \quad\text{ for all $\bfv \in \bL$.}
\end{displaymath}
Therefore by Lemma~\ref{lemma:switching:Rxy:Ryx}, there exists
$0 < \theta'' <1/2$, $K''(\ell) \subset X$ and for each $x \in K''(\ell)$ a
subset $Y_\ell'(x) \subset \cF_{ij}[x,\ell]$ with $|Y_\ell'(x)| \ge
(1-\theta'') |\cF_{ij}[x,\ell]|$ such that for $x \in K''(\ell)$ and
$y \in Y_\ell'(x)$, (\ref{eq:tmp:Rxy:v:twosided:small}) holds.

Let
\begin{displaymath}
\bZ(x,\eta) = \{ \bfv \in \proj(\bL) \st d(\bfv, \bF_1(x) \cup \bF_2(x))
\ge \eta \}.
\end{displaymath}
We may choose $\eta > 0$ small enough so that there exists $K' \subset
X$ with $\nu(K''(\ell) \cap K') > 0$ such that for all $x \in K'$, 
\begin{displaymath}
\rho_0( \bZ(x,C_1 \eta)) > 1/2. 
\end{displaymath}
Let 
\begin{displaymath}
S(\eta) = \{ (x,\bfv) \st x \in X, \quad \bfv \in \bZ(x,\eta) \}
\end{displaymath}
Let $f$ denote the characteristic function of the set 
\begin{displaymath}
\{ (x, \bfv ) \st  x \in K''(\ell) \cap K', \quad \bfv \in \bZ(x,\eta)\}
\subset S(\eta). 
\end{displaymath}
We now claim that for any $\ell$, 
\begin{equation}
\label{eq:tilde:nu:ell}
\hat{\mu}_\ell(f) \ge \nu(K''(\ell) \cap K') (1-\theta'') (1/2). 
\end{equation}
Indeed, if we restrict in (\ref{eq:def:hat:mu})  
to $x \in K''(\ell) \cap K'$, $y \in Y_\ell'(x)$, and  $\bfv \in
\bZ(x,C_{1} \eta)$, then by (\ref{eq:tmp:Rxy:v:twosided:small}),
$f(x,R(x,y) \bfv) = 1$. This implies (\ref{eq:tilde:nu:ell}). Thus,
(provided $\delta > 0$ and $\theta > 0$ in
Definition~\ref{def:bounded:subspace} are sufficiently small), 
there exists $c_0 > 0$ such that for all $\ell$, 
$\hat{\mu}_\ell(S(\eta)) \ge c_0 > 0$. Therefore, by
Lemma~\ref{lemma:proj:general}, 
$\tilde{\mu}_\ell(S(\eta)) \ge c_0/\kappa^2$.

There exists compact $K_0
\subset X$ with $\nu(K_0) > 1 - c_0/(2 \kappa^2)$ such that the map $x
\to \bF_1(x) \cap \bF_2(x)$ is continuous on $K_0$. Let $K_0' = \{ (x,
\bfv) \st x \in K_0\}$. Then $S(\eta) \cap K_0'$
is a closed set with $\tilde{\mu}_\ell(S(\eta) \cap K_0') \ge c_0/(2
\kappa^2)$. Therefore, 
$\tilde{\mu}_\infty(S(\eta) \cap K_0') > c_0/(2\kappa^2) > 0$,
which is a contradiction to the
fact that $\lambda_x$ is supported on $\bF_1(x) \cup \bF_2(x)$.  

Thus, for almost all $x$, 
$\lambda_x$ is not supported on a union of two subspaces. Thus the
same holds for the measure $\psi_x$ of
Lemma~\ref{lemma:lambda:x:invariant} (d). By
combining (b) and (d) of Lemma~\ref{lemma:lambda:x:invariant} we see
that for almost all $x$ and almost all $y \in \cF_{ij}[x]$, 
\begin{displaymath}
R(x,y) h(x) \psi_x = h(y) \psi_x,
\end{displaymath}
hence $h(y)^{-1} R(x,y) h(x)$ stabilizes $\psi_x$. Hence by
Lemma~\ref{lemma:furstenberg:two:subspace}, 
\begin{displaymath}
h(y)^{-1} \bar{R}(x,y) h(x) \in K(x)
\end{displaymath}
where $K(x)$ is a compact subset of $SL(\bL)$, and $\bar{R}(x,y)$ is
the image of $R(x,y)$ under the natural map $GL(\bL) \to SL(\bL)$.
Thus, $\bar{R}(x,y) \in h(y) K(x) h(x)^{-1}$, and thus
\begin{equation}
\label{eq:bar:R:upper}
\|\bar{R}(x,y)\| \le C(x) C(y). 
\end{equation}
Since $\bar{R}(x,y)^{-1} = \bar{R}(y,x)$, we get, by exchanging
  $x$ and $y$,
\begin{equation}
\label{eq:bar:R:lower}   
\|\bar{R}(x,y)^{-1}\| \le C(x) C(y). 
\end{equation}
Note that by Lemma~\ref{lemma:bounded:subspace},
there exists $\bfv \in \bL(x)=
\bE_{ij,bdd}(x) \subset \bE_{ij}(x)$ such that $\bfv \not\in
\bE_{i,j-1}(x)$. Then, (\ref{eq:Rxy:bfv}) and the fact that
$\lambda_{ij}(x,y) = 0$ for $y \in \cF_{ij}[x]$ shows that 
(\ref{eq:bar:R:upper}) and (\ref{eq:bar:R:lower}) must
hold for $R(x,y)$ in place of $\bar{R}(x,y)$.
This implies the statement of the lemma.
\qed\medskip

\begin{lemma}
\label{lemma:synchronized:always}
Suppose that for all $\delta > 0$ there exists a constant $C > 0$
and  a compact subset $K \subset X$ 
with $\nu(K) > 1-\delta$ and for each $\ell > 0$ and 
$x \in K$ a subset $Y_\ell(x)$
of $\cF_{ij}[x,\ell]$ with $|Y_\ell(x)| \ge (1-\theta)
|\cF_{ij}[x,\ell]|$, such that for $x \in K$
and $y \in Y_\ell(x)$ we have
\begin{equation}
\label{eq:lambda:kr:upper}
\lambda_{kr}(x,y) \le C. 
\end{equation}
Then, $ij$ and $kr$ are synchronized, and there exists a function $C:
X \to \reals^+$ finite $\nu$-almost everywhere such that
for all $x \in X$, and all $y \in \cF_{ij}[x]$, 
\begin{equation}
\label{eq:rho:y:Fkr:x}
\rho(y,\cF_{kr}[x]) \le C(x)C(y). 
\end{equation}
\end{lemma}

\bold{Proof.} 
The proof is a simplified version of the proof of
Lemma~\ref{lemma:Vbdd:undependent:u}. Let $\bL_1 =
\bE_{ij}/\bE_{i,j-1}$, $\bL_2 = \bE_{kr}/\bE_{k,r-1}$, and $\bL =
\bL_1 \cross \bL_2$.

We have, for $y \in G[x]$, and
$(\bar{\bfv}, \bar{\bfw}) \in \bL$, \mccc{($R$ is not defined in this context)}
\begin{multline}
\label{eq:Rxy:bar:bfu:bar:bfv}
R(x,y) (\bar{\bfv}, \bar{\bfw}) = (e^{\lambda_{ij}(x,y)} \bar{\bfv}',
e^{\lambda_{kr}(x,y)} \bar{\bfw}'), \\ \text{ where $\|\bar{\bfv}' \|=
  \|\bar{\bfv}\|$ and $\|\bar{\bfw}' \|=
  \|\bar{\bfw}\|$.}
\end{multline}
Recall that $\lambda_{ij}(x,y) = 0$ for all $y \in
\cF_{ij}[x]$. Therefore, (\ref{eq:lambda:kr:upper}) implies that for
all $x \in K$, all $\ell > 0$ and all $y \in Y_\ell(x)$, 
\begin{displaymath}
\|R(x,y) (\bar{\bfv}, \bar{\bfw}) \| \le C_1
\|(\bar{\bfv},\bar{\bfw})\|. 
\end{displaymath}
Therefore, by Lemma~\ref{lemma:switching:Rxy:Ryx}, there exists a
subset $K''(\ell) \subset X$ with $\nu(K''(\ell)) > 1-c(\delta)$ where
$c(\delta) \to 0$ as $\delta \to 0$, and for each $x \in K''(\ell)$ a
subset $Y_\ell' \subset \cF_{ij}[x,\ell]$ with $|Y_\ell'| > (1-\theta'')
|\cF_{ij}[x,\ell]|$ such that for all $y \in Y_\ell'$,
\begin{displaymath}
C_1^{-1} \|(\bar{\bfv},\bar{\bfw})\| \le \|R(x,y) (\bar{\bfv},
\bar{\bfw}) \| \le C_1 \|(\bar{\bfv},\bar{\bfw})\|. 
\end{displaymath}
This implies that for $x \in K''(\ell)$, $y \in Y_\ell'(x)$, 
\begin{equation}
\label{eq:lambda:kr:xy:abs}
|\lambda_{kr}(x,y)| = |\lambda_{ij}(x,y) - \lambda_{kr}(x,y)| \le C_1.
\end{equation}

Let $\tilde{\mu}_\ell$ and $\hat{\mu}_\ell$ be as in
Lemma~\ref{lemma:proj:general}. Take a sequence $\ell_m \to \infty$ 
such that $\tilde{\mu}_{\ell_m} \to \tilde{\mu}_\infty$, and
$\hat{\mu}_{\ell_m} \to \hat{\nu}_\infty$. 
Then by Lemma~\ref{lemma:lambda:x:invariant} (a),
we have $d\tilde{\mu}_\infty(x,\bfv) = d\nu(x) \, d\lambda_x(\bfv)$ where
$\lambda_x$ is a measure on $\proj(\bL)$. We will show that for almost
all $x \in X$, $\lambda_x$ is not supported on $\bL_1 \cross \{0\}
\cup \{0 \} \cross \bL_2$. 

Suppose that for a set of positive measure $\lambda_x$ is
supported on $(\bL_1 \cross \{0\})
\cup (\{0 \} \cross \bL_2)$.  Then, in view of the ergodicity of $g_t$ and
Lemma~\ref{lemma:lambda:x:invariant} (c), $\lambda_x$ is supported on
$(\bL_1 \cross \{0\}) \cup (\{0 \} \cross \bL_2)$ for almost all $x \in
X$. Let 
\begin{displaymath}
\bZ(x,\eta) = \{ (\bar{\bfv},\bar{\bfw}) \in \bL(x), \quad
\|(\bar{\bfv},\bar{\bfw})\| = 1, 
\quad d(\bar{\bfv},\bL_1) \ge \eta, \quad d(\bar{\bfw}, \bL_2) \ge \eta
\}. 
\end{displaymath}
and let 
\begin{displaymath}
S(\eta) = \{ (x, (\bar{\bfv},\bar{\bfw})) \st x \in X,
(\bar{\bfv},\bar{\bfw}) \in \bZ(x,\eta) \}.
\end{displaymath}
Then we have $\tilde{\mu}_\infty(S(\eta)) = 0$. Therefore, by
Lemma~\ref{lemma:proj:general}, $\hat{\mu}_\infty(S(\eta)) = 0$. 

By (\ref{eq:Rxy:bar:bfu:bar:bfv}) and (\ref{eq:lambda:kr:xy:abs}), 
for $x \in K''(\ell_m)$ and $y \in Y'_{\ell_m}(x)$,
\begin{equation}
\label{eq:Rxy:bZxC1eta}
R(x,y) \, \bZ(x,C_1 \eta) \subset \bZ(y,\eta). 
\end{equation}
Choose $\eta >0$ so that there exists $K'=K'(\ell_m) \subset
X$ with $\nu(K''(\ell_m) \cap 
K') > 0$ such that for $x \in K'$,  
$\rho_0(\bZ(x,C_1 \eta)) > (1/2)$. 
Let $f$ be the characteristic function of $S(\eta)$. 
Then, if we restrict in (\ref{eq:def:hat:mu})
to $x \in K''(\ell_m) \cap K'$, $y \in
Y_{\ell_m}'(x)$, and  $\bfv \in 
\bZ(x,C_1 \eta)$, then by (\ref{eq:Rxy:bZxC1eta}), $f(x,R(x,y) \bfv) =
1$. This implies that for all $m$, 
\begin{displaymath}
  \hat{\mu}_{\ell_m}(S(\eta)) \ge \nu(K''(\ell_m)
  \cap K') (1-\theta'') (1/2). 
\end{displaymath}
Hence $\hat{\mu}_{\infty}(S(\eta)) > 0$ which is a contradiction. 
Therefore, for almost all $x$, 
$\lambda_x$ is not supported on $\bL_1 \cross \{0\}
\cup \{0 \} \cross \bL_2$. 
Thus the same holds for the measure $\psi_x$ of
Lemma~\ref{lemma:lambda:x:invariant} (d). By
combining (b) and (d) of Lemma~\ref{lemma:lambda:x:invariant} we see
that for almost all $x\in X$ and almost all $y \in \cF_{ij}[x]$, 
\begin{displaymath}
R(x,y) h(x) \psi_x = h(y) \psi_x,
\end{displaymath}
hence $h(y)^{-1} R(x,y) h(x)$ stabilizes $\psi_x$. Note that in view
of (\ref{eq:Rxy:bar:bfu:bar:bfv}), 
\begin{multline*}
h(y)^{-1} R(x,y) h(x) (\bar{\bfv}, \bar{\bfw}) = (e^{\alpha(x,y)} \bar{\bfv}',
e^{\alpha'(x,y)} \bar{\bfw}'), \\ \text{ where $\alpha(x,y) \in
  \reals$, $\alpha'(x,y) \in \reals$, $\|\bar{\bfv}' \|=
  \|\bar{\bfv}\|$ and $\|\bar{\bfw}' \|=
  \|\bar{\bfw}\|$.} 
\end{multline*}
For $i=1,2$ let $\Conf_x(\bL_i)$ denote the subgroup of $GL(\bL_i)$
which preserves the inner product $\langle\cdot,\cdot\rangle_x$ up to
a scaling factor. Let $\Conf_x(\bL) = \Conf_x(\bL_1) \cross
\Conf_x(\bL_2)$.
Then, by an elementary variant of Lemma~\ref{lemma:furstenberg:two:subspace}, 
since $\psi_x$ is not supported on $\bL_1 \cross \{0\}
\cup \{0 \} \cross \bL_2$, we get
\begin{displaymath}
h(y)^{-1} R(x,y) h(x) \in K(x)
\end{displaymath}
where $K(x)$ is a compact subset of $\Conf_x(\bL)$.
Thus, $R(x,y) \in h(y)
K(x) h(x)^{-1}$, and thus
\begin{displaymath}
\|R(x,y)\| \le C(x) C(y). 
\end{displaymath}
Note that by reversing $x$ and $y$ we get $\|R(x,y)^{-1}\| \le
C(x) C(y)$. Therefore, by (\ref{eq:Rxy:bar:bfu:bar:bfv}), 
\begin{displaymath}
|\lambda_{ij}(x,y) - \lambda_{kr}(x,y)| \le C(x) C(y). 
\end{displaymath}
This completes the proof of (\ref{eq:rho:y:Fkr:x}).

For any $\delta > 0$ we can choose a compact $K \subset X$
with $\nu(K) > 1-\delta$ and $N < \infty$ such that $C(x) < N$ for $x \in K$. 
Now, the fact that $ij$ and $kr$ are synchronized follows from
applying Lemma~\ref{lemma:fake:ergodicity:Fij} to $K$.
\qed\medskip

\bold{Proof of Proposition~\ref{prop:ej:bdd:transport:bounded}.}
This follows immediately from
Lemma~\ref{lemma:synchronized:always} and 
Lemma~\ref{lemma:Vbdd:undependent:u}.
\qed\medskip

\bold{Proof of
  Proposition~\ref{prop:Rxy:v:small:implied:sync:bounded}.}
Choose $\epsilon < \epsilon'/10$, where $\epsilon'$ is as in Proposition~\ref{prop:most:inert} (b). 
By the multiplicative ergodic theorem, there exists a set $K''_1
\subset X$ with
$\nu(K''_1) > 1-\theta$ and $T > 0$, such that for $x \in K''_1$ and
$t > T$, 
\begin{equation}
\label{eq:tmp:lambdaij:lambdai}
|\lambda_{ij}(x,t) - \lambda_i t| < \epsilon t,
\end{equation}
where $\lambda_{ij}(x,t)$ is as in (\ref{eq:gt:star:bfv}). Then, by
Fubini's theorem there exists a set $K''_2 \subset K''_1$ with $\nu(K''_2)
> 1-3\theta$ such that for $x \in K''_2$, for $(1-\theta)$-fraction of $u \in \cB(x)$, $u x \in K''_1$. 

Let $K''$ be as in Proposition~\ref{prop:most:inert} (b) with $\delta
= \theta$. We may 
assume that the conull set $\Psi$ in
Proposition~\ref{prop:Rxy:v:small:implied:sync:bounded}  
is such so that for $x \in \Psi$, $g_{-t} x \in
K'' \cap K''_2$ for arbitrarily large $t > 0$. 

Suppose $g_{-t} x \in K'' \cap K''_2$ and $y \in \cF_{ij}[x]$. We may write 
\begin{displaymath}
y = g^{ij}_{t'} u g^{ij}_{-t'} x  = g_{s'} u g_{-t} x.
\end{displaymath}
By the definition of $\cF_{ij}[x,t']$, and since $g_{-t} x \in K''_2$,
we have $g_{-t} x \in K''_1$ and 
for at least $(1-\theta)$-fraction of $y \in \cF_{ij}[x,t']$, we
have $u g_{-t} x \in K''_1$, and thus, in view of (\ref{eq:tmp:lambdaij:lambdai}), 
\begin{displaymath}
|s' - \lambda_i t'| \le \epsilon t \quad\text{and}\quad |t - \lambda_i
t'| \le \epsilon t.  
\end{displaymath}
Therefore for $(1-\theta)$-fraction of $y \in \cF_{ij}[x,t']$ or
equivalently for $(1-\theta)$-fraction of $u \in \cB(g_{-t} x)$, 
\begin{equation}
\label{eq:tmp:time:changes:within:epsilon}
|s' - t| \le 2 \epsilon t.
\end{equation}
Now suppose $\bfv \in \bH(x)$. Note that if
$\|R(x,y) \bfv \| \le C \|\bfv\|$, and $s$ is as in
Proposition~\ref{prop:most:inert}, then $s > s' - O(1)$ (where the
implied constant depends on $C$.) Therefore, in view of
(\ref{eq:tmp:time:changes:within:epsilon}), for $(1-\theta)$-fraction
of $u \in \cB(g_{-t} x)$, (\ref{eq:tmp:most:inert:b})
holds. Thus, 
by Proposition~\ref{prop:most:inert}(b), we have $\bfv \in
\bE(x)$. 
Thus, we can write
\begin{displaymath}
\bfv = \sum_{kr \in I_\bfv} \bfv_{kr}
\end{displaymath}
where the indexing set $I_\bfv$ contains at most one $r$ for each $k \in
\Lambda'$.
Without loss
of generality,  $\Psi$
is such that for $x \in \Psi$, $g_{-t} x$ satisfies the conclusions of
Proposition~\ref{prop:properties:dynamical:norm:Hbig} infinitely
often. Note that for $y \in \cF_{ij}[x]$, 
\begin{displaymath}
\|R(x,y) \bfv\| \ge \|R(x,y) \bfv_{kr} \| \ge e^{\lambda_{kr}(x,y)}
\|\bfv_{kr} \|. 
\end{displaymath}
By assumption, for 
all $\ell > 0$ and for at least $1-\theta$ fraction of $y \in
\cF_{ij}[x,\ell]$, $\|R(x,y) \bfv\| \le C$. Therefore, for all $\ell >
0$ and for at least $(1-\theta)$ fraction of $y \in \cF_{ij}[x,\ell]$,
(\ref{eq:lambda:kr:upper}) holds. 
Then, by Lemma~\ref{lemma:synchronized:always}, for all $kr \in
I_\bfv$, $kr$ and $ij$ are synchronized, i.e.\ $kr \in [ij]$. Therefore, for
at least $(1-2\theta)$-fraction of $y' \in \cF_{kr}[x,\ell]$, 
\begin{displaymath}
\|R(x,y') \bfv_{kr}\| \le \|R(x,y') \bfv \| \le C'.
\end{displaymath}
Now, by
Definition~\ref{def:bounded:subspace}, 
$\bfv_{kr}(x) \in
\bE_{kr,bdd}(x)$. Therefore, $\bfv \in \bE_{[ij],bdd}(x)$. \mcc{(more
  details here)}
\qed\medskip 

It follows from the proof of
Proposition~\ref{prop:Rxy:v:small:implied:sync:bounded} that
(\ref{eq:def:ij:prime}) holds.

\section{Equivalence relations on $W^+$}
\label{sec:equivalence}

Let \index{$G$@$\GSpc$}$\GSpc$ denote the space of generalized
subspaces of $W^+$.  Let
\index{$H$@$\bar{\cH}_{++}(x)$}$\bar{\cH}_{++}(x)$ denote the set of
$M \in \cH_{++}(x)$ such that $(I+M)\Lie(U^+)(x)$ is a subalgebra of
$\Lie(\cG_{++})(x)$.  We have a map $\index{$U$@$\cU_x$}\cU_x: \bar{\cH}_{++}(x) \cross
W^+(x) \to \GSpc$ taking the pair $(M,v)$ to the generalized subspace
it parametrizes. Let $\cU_x^{-1}$ denote the inverse of this map
(given a Lyapunov-adapted transversal $Z(x)$).

For $ij \in \tilde{\Lambda}$, let 
\begin{displaymath}
\index{$E$@$\cE_{ij}[x]$}\cE_{ij}[x] = \{ \cQ \in \GSpc \st \bfj(\cU_x^{-1}(\cQ)) \in
\bE_{[ij],bdd}(x) \}.
\end{displaymath}

\bold{Motivation.} In view of
Proposition~\ref{prop:ej:bdd:transport:bounded} and
Lemma~\ref{lemma:hausdorff:distance:to:norm} (b),  for any
sufficiently small $\epsilon > 0$, 
the conditions that $\cQ \in \cE_{ij}[x]$ and $hd_{x}^{X_0}(\cQ,
U^+[x]) = O(\epsilon)$ 
imply the following: for ``most'' $y \in \cF_{ij}[x]$, 
\begin{displaymath}
hd_y^{X_0}(R(x,y) \cQ, U^+[y]) = O(\epsilon). 
\end{displaymath}


\bold{A partition of $W^+[x]$.} Let $\gB_0$ denote the measurable partition
constructed in \S\ref{sec:semi:markov}, (see also
\S\ref{sec:subsec:the:cover:X}). We denote
the atom containing $x$ by $\gB_0[x]$, and let $\gB_0(x) = \{ v \in W^+(x) \st
v+x \in \gB_0[x] \}$. In this section, the only properties of $\gB_0$ we
will use is that it is subordinate to $W^+$, and that the atoms
$\gB_0[x]$ are relatively open in $W^+[x]$. 

\bold{Equivalence relations.} 
Fix $x_0 \in X$. 
For $x$, $x' \in W^+[x_0]$ we say that
\begin{displaymath}
x' \index{$\sim_{ij}$}\sim_{ij} x \text{ if $x' \in \gB_0[x]$ and $U^+[x'] \in \cE_{ij}[x]$. }
\end{displaymath}
\begin{proposition}
\label{prop:sim:ij:is:equivalence:relation}
The relation $\sim_{ij}$ is a (measurable) equivalence relation. 
\end{proposition}
The main part of the proof of
Proposition~\ref{prop:sim:ij:is:equivalence:relation} is the following:
\begin{lemma}
\label{lemma:Eij:bdd:locally:constant}
There exists a subset $\index{$\Psi$}\Psi \subset X$ 
with $\nu(\Psi) = 1$ such
that for any $ij \in \tilde{\Lambda}$, 
if $x_0 \in \Psi$, $x_1 \in \Psi$, $x_1 \in \gB_0[x_0]$
(so in particular $d^{X_0}(x_0,x_1) < 1/100$), and $U^+[x_1] \in \cE_{ij}[x_0]$, then 
$\cE_{ij}[x_1] = \cE_{ij}[x_0]$. 
\end{lemma}

\bold{Warning.}
We will consider the condition $x' \sim_{ij} x$ to be undefined unless
$x$ and $x'$ both belong to the set $\Psi$ of
Lemma~\ref{lemma:Eij:bdd:locally:constant}. 

\bold{Motivation.}
In view of Proposition~\ref{prop:some:fraction:bounded}, we can
ensure, in the notation of \S\ref{sec:outline:step1} that for some $ij
\in \tilde{\Lambda}$, $U^+[q_2']$ is
close to $\cE_{ij}[q_2]$; then in the limit we would have
$U^+[\tilde{q}_2'] \in \cE_{ij}[\tilde{q}_2]$, and thus $\tilde{q}_2'
\sim_{ij} \tilde{q}_2$.

\bold{Proof of Proposition~\ref{prop:sim:ij:is:equivalence:relation},
  assuming Lemma~\ref{lemma:Eij:bdd:locally:constant}.}
We have $0 \in \bE_{[ij],bdd}(x)$,
therefore, 
\begin{equation}
\label{eq:Uplus:x:in:cEij:x}
U^+[x] \in \cE_{ij}[x]. 
\end{equation}
Thus $x \sim_{ij} x$. 

Suppose $x' \sim_{ij} x$. Then, $x' \in \gB_0[x]$, and so $x \in \gB_0[x']$. 
By (\ref{eq:Uplus:x:in:cEij:x}), $U^+[x] \in \cE_{ij}[x]$, and by
Lemma~\ref{lemma:Eij:bdd:locally:constant}, $\cE_{ij}[x'] =
\cE_{ij}[x]$. Therefore, $U^+[x] \in \cE_{ij}[x']$, and thus $x
\sim_{ij} x'$. 

Now suppose $x' \sim_{ij} x$ and $x'' \sim_{ij} x'$. Then, $x'' \in
\gB_0[x]$. Also, $U^+[x''] \in \cE_{ij}[x'] = \cE_{ij}[x]$, therefore $x''
\sim_{ij} x$. 
\qed\medskip

\bold{Remark.} By Lemma~\ref{lemma:Eij:bdd:locally:constant}, for
$x,x' \in \Psi$, $x' \sim_{ij} x$ if and only if $x' \in \gB_0[x]$ and
$\cE_{ij}[x'] = \cE_{ij}[x]$.

\bold{Outline of the proof of
  Lemma~\ref{lemma:Eij:bdd:locally:constant}.} Intuitively, the
condition $U^+[x_1] \in \cE_{ij}[x_0]$ is the same as ``$\cF_{ij}[x_1]$
and $\cF_{ij}[x_0]$ stay close'', and 
``$U^+[x_1]$ and $U^+[x_0]$ stay close as we travel along
$\cF_{ij}[x_0]$ or $\cF_{ij}[x_1]$'', 
which is clearly an equivalence relation. We give some more detail below. 
Throughout the proof we will be using
Lemma~\ref{lemma:stay:in:same:orbit}, without mentioning it
explicitly. 

Fix $\epsilon \ll 1/100$.
Suppose $x_1 \in \gB_0[x_0]$, so in particular 
$d^{X_0}(x_0,x_1) < 1/100$,  and suppose
\begin{displaymath}
hd_{x_0}^{X_0}(U^+[x_1], U^+[x_0]) = \epsilon. 
\end{displaymath}
Then, by Lemma~\ref{lemma:hausdorff:distance:to:norm} (b), 
\begin{displaymath}
\bfj(\cU_{x_0}^{-1}(U^+[x_1])) = O(\epsilon). 
\end{displaymath}
We are given that
$U^+[x_1] \in \cE_{ij}[x_0]$, thus $\bfj(\cU_{x_0}^{-1}(U^+[x_1])) \in
\bE_{[ij],bdd}(x_0)$. Then, by
Proposition~\ref{prop:ej:bdd:transport:bounded}, for most $y_0 \in
\cF_{ij}[x_0]$, 
\begin{displaymath}
\|R(x_0,y_0) \bfj(\cU_{x_0}^{-1}(U^+[x_1]))\| = O(\epsilon). 
\end{displaymath}
We have
\begin{displaymath}
R(x_0,y_0) \bfj(\cU_{x_0}^{-1}(U^+[x_1])) = \bfj(\cU_{y_0}^{-1}(U^+[y_1'])),
\end{displaymath}
for some $y_1' \in G[x_1]$. 
Then, by
Lemma~\ref{lemma:hausdorff:distance:to:norm} (b), for most $y_0 \in
\cF_{ij}[x_0]$, 
\begin{displaymath}
hd_{y_0}^{X_0}(U^+[y_1'],U^+[y_0]) = O(\epsilon) \quad \text{ for some $y_1'
  \in G[x_1]$. }
\end{displaymath}
It is not difficult to show that $y_1'$ is
near a point $y_1 \in \cF_{ij}[x_1]$. Thus, for most $y_0 \in
\cF_{ij}[x_0]$, 
\begin{equation}
\label{eq:tmp:hd:y0:Uplus:y1:Uplus:y0}
hd_{y_0}^{X_0}(U^+[y_1],U^+[y_0]) = O(\epsilon) \quad \text{ for some $y_1
  \in \cF_{ij}[x_1]$. }
\end{equation}
Thus, most of the time
$\cF_{ij}[x_0]$ and $\cF_{ij}[x_1]$ remain close, 
and also that for most $y_0 \in \cF_{ij}[x_0]$, 
$U^+[y_1]$ and $U^+[y_0]$ remain close, for some $y_1 \in
\cF_{ij}[x_1]$. 

Now suppose $\cQ_1 \in \cE_{ij}[x_1]$, and 
\begin{displaymath}
hd_{x_1}^{X_0}(\cQ_1, U^+[x_1]) = O(\epsilon). 
\end{displaymath}
Then,
$\bfj(\cU_{x_1}^{-1}(\cQ_1)) \in \bE_{[ij],bdd}(x_1)$, and thus, for
most $y_1 \in \cF_{ij}[x_1]$, using
Proposition~\ref{prop:ej:bdd:transport:bounded} and
Lemma~\ref{lemma:hausdorff:distance:to:norm} (b) twice as above, we
get that for most $y_1 \in \cF_{ij}[x_1]$, 
\begin{equation}
\label{eq:tmp:hd:y1:Rx1:y1:cQ1:Uplus:y1}
hd_{y_1}^{X_0}(R(x_1,y_1) \cQ_1, U^+[y_1]) = O(\epsilon). 
\end{equation}
In our notation, $R(x_1,y_1) \cQ_1$ is the same generalized subspace
(i.e.\ the same subset of $W^+$) 
as $R(x_0,y_0) \cQ_1$ for $y_0 \in \cF_{ij}[x_0]$ close to
$y_1$. 
Then, from (\ref{eq:tmp:hd:y0:Uplus:y1:Uplus:y0}) and
(\ref{eq:tmp:hd:y1:Rx1:y1:cQ1:Uplus:y1}), for most $y_0 \in
\cF_{ij}[x_0]$, 
\begin{displaymath}
hd_{y_0}^{X_0}(R(x_0,y_0) \cQ_1, U^+[y_0]) = O(\epsilon).
\end{displaymath}
Thus, using Lemma~\ref{lemma:hausdorff:distance:to:norm} (b) again, we
get that for most $y_0 \in \cF_{ij}[x_0]$, 
\begin{displaymath}
\|R(x_0,y_0) \bfj(\cU_{x_0}^{-1}(\cQ_1)) \| = O(\epsilon). 
\end{displaymath}
By Proposition~\ref{prop:Rxy:v:small:implied:sync:bounded}, this
implies that $\bfj(\cU_{x_0}^{-1}(\cQ_1)) \in \bE_{[ij],bdd}(x_0)$,
and thus $\cQ_1 \in \cE_{ij}[x_0]$. Thus, $\cE_{ij}[x_1] \subset
\cE_{ij}[x_0]$.


Conversely, if $\cQ_0 \in \cE_{ij}[x_0]$, then  the same argument
shows that $\cQ_0 \in \cE_{ij}[x_1]$. Therefore, $\cE_{ij}[x_0] =
\cE_{ij}[x_1]$. 
\qed\medskip

\mcc{handle $\epsilon$ vs $1$ in above proof.}


The (tedious) formal verification of
Lemma~\ref{lemma:Eij:bdd:locally:constant} is given in
\S\ref{sec:starredsubsec:proof:lemma:Eij:bdd:locally:constant} below. 

\bold{The equivalence classes $\cC_{ij}[x]$.} 
For $x \in \Psi$ we define the equivalence class
\begin{displaymath}
\index{$C$@$\cC_{ij}[x]$}\cC_{ij}[x] = \{ x' \in \gB_0[x] \st x' \sim_{ij} x \}.
\end{displaymath}
Let $\cC_{ij}$ denote the $\sigma$-algebra
of $\nu$-measurable sets which are unions of the equivalence classes
$\cC_{ij}[x]$. We do not distinguish between $\sigma$-algebras which
are equivalent mod sets of $\nu$-measure $0$, so we can assume that 
$\cC_{ij}$ is countably generated (see
\cite[\S{1.2}]{Climenhaga:Katok}). 
We now want to show that (away from a set of measure
$0$), the atoms of the $\sigma$-algebra $\cC_{ij}$ are the sets
$\cC_{ij}[x]$. More precisely, we want to show that the partition
$\cC_{ij}$ whose atoms are the sets $\cC_{ij}[x]$ is
a measurable partition in the sense of
\cite[Definition 1.10]{Climenhaga:Katok}. 

To see this, note that each set $\cE_{ij}[x]$ is an algebraic subset
of $\GSpc$, and is thus parametrized by a finite dimensional space
$Y$. Let $\psi_{ij}: X \to Y$ be the map taking $x$ to the
parametrization of $\cE_{ij}[x]$. We note that the functions
$\psi_{ij}$ are measurable. \mcc{check} 
Also, in view of Lemma~\ref{lemma:Eij:bdd:locally:constant}, we have
\begin{displaymath}
x \sim_{ij} y \qquad \text{if and only if $y \in \gB_0[x]$ and
$\psi_{ij}(y) = \psi_{ij}(x)$.}
\end{displaymath}
By Lusin's theorem, for each $ij$, 
there exists a Borel function $\tilde{\psi}_{ij}$ such that
$\nu$-almost everywhere, 
$\tilde{\psi}_{ij} = \psi_{ij}$.
Now the measurability of $\cC_{ij}$ follows from
\cite[Theorem~1.14]{Climenhaga:Katok}.

\begin{lemma}
\label{lemma:cCij:equivariant}
Suppose $t \in \reals$, $u \in U^+(x)$.
\begin{itemize}
\item[{\rm (a)}] 
$g_t \cC_{ij}[x] \cap \gB_0[g_t x] \cap g_t \gB_0[x] = \cC_{ij}[g_t x] \cap
\gB_0[g_t x] \cap g_t \gB_0[x]$.
\item[{\rm (b)}]
$u\cC_{ij}[x] \cap \gB_0[ux] \cap u\gB_0[x] = \cC_{ij}[ux] \cap \gB_0[ux] \cap
u \gB_0[x]$. 
\end{itemize}
\end{lemma}

\bold{Proof.} Note that the sets $U^+[x]$ and $\bE_{[ij],bdd}(x)$ are
$g_t$-equivariant. Therefore, so are the $\cE_{ij}[x]$, which implies
(a). 
Part (b) is also clear, since 
locally, by Lemma~\ref{lemma:Ej:equivariant}, $(u)_* \bE_{ij}(x) =
\bE_{ij}(ux)$. 
\qed\medskip

\bold{The measures $f_{ij}[x]$.}
We now define \index{$f_{ij}[x]$}$f_{ij}[x]$ 
to be the conditional measure
of $\nu$ along the $\cC_{ij}[x]$. In other words, $f_{ij}[x]$
is defined so that for any measurable $\phi: X \to \reals$,
\begin{displaymath}
\mathbb{E}(\phi \mid \cC_{ij})(x) = \int_X \phi \, df_{ij}[x]. 
\end{displaymath}
We view $f_{ij}[x]$ as a measure on $W^+[x]$ 
which is supported on $\cC_{ij}[x]$. 

\bold{The measures $f_{ij}(x)$.}
We can identify $W^+[x]$ with the vector space $W^+(x)$,
where $x$ corresponds to the origin. Let \index{$f_{ij}(x)$}$f_{ij}(x)$ be
the pullback to $W^+(x)$ of $f_{ij}[x]$ under this
identification. We will also call the $f_{ij}(x)$ conditional
measures. (The term ``leaf-wise'' measures is used in
\cite{Einsiedler:Lindenstrauss:Pisa} in a related context). We abuse
notation slightly and write formulas such as  
\begin{displaymath}
\mathbb{E}(\phi \mid \cC_{ij})(x) = \int_X \phi \, df_{ij}(x). 
\end{displaymath}

\bold{The ``distance'' $d_*(\cdot, \cdot)$.} 
  Suppose $E_1, E_2$ are open subsets of a normed
  vector space $V$, with $E_1 \cap
E_2 \ne \emptyset$. Suppose that for $i=1,2$, $\mu_i$ is a finite measure on
$E_i$, with $\mu_i(E_1 \cap E_2) > 0$.  Then, let $d_*(\mu_1, \mu_2)$ denote the
Kontorovich-Rubinstein distance between (the normalized versions of)
$\bar{\mu}_1$ and
$\bar{\mu_2}$, i.e.  \index{$d_*(\cdot,\cdot)$}
\begin{displaymath}
d_*(\mu_1,\mu_2) = \sup_{f} \left|\frac{1}{\mu_1(E_1 \cap
    E_2)}\int_{E_1 \cap E_2} f \,  
  d\mu_1 - \frac{1}{\mu_2(E_1 \cap E_2)}
  \int_{E_1 \cap E_2} f \, d\mu_2
\right|,  
\end{displaymath}
where the sup is taken over all $1$-Lipshitz functions $f: E_1 \cap E_2\to
\reals$ with $\sup |f(x)| \le  1$.

The only property of $d_*(\cdot, \cdot)$ we will use is that it
induces the topology of weak-* convergence on the
domain of common definition of the measures, up to
normalization. \mcc{say this better}


\makefig{Proposition~\ref{prop:nearby:linear:maps}}{fig:outline2}{\includegraphics{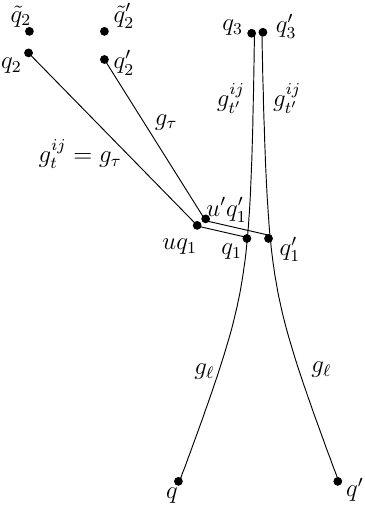}}

\medskip
The following Proposition 
is the rigourous version of (\ref{eq:f1:approx:same}) in
\S\ref{sec:outline:step1}: 

\begin{proposition}
\label{prop:nearby:linear:maps}
There exists $0 < \alpha_0 < 1$ depending only on the Lyapunov
spectrum, and 
for every $\delta > 0$ there exists a compact set $K_0 \subset X$ with
$\nu(K_0) > 1-\delta$ such that the following holds:
Suppose $ij \in \tilde{\Lambda}$,
$1 < C_1 < \infty$, $0 < \epsilon < C_1^{-1}/100$, 
$C < \infty$, $t > 0$, $t' > 0$, 
and $|t'-t| < C$. Furthermore  
suppose $q \in \pi^{-1}(K_0)$ and $q' \in W^-[q] \cap \pi^{-1}(K_0)$
are such that $d^{X}(q,q') < 1/100$. \mcc{check that this is enough}
Let
$q_1 = g_\ell q$, $q_1' = g_{\ell} q'$. 
Also let $q_3 = g^{ij}_{t'} q_1$, $q_3' = g^{ij}_{t'} q_1'$. Suppose $q_1,
q_1'$, $q_3$, $q_3'$ all belong to $\pi^{-1}(K_0)$. 

Suppose $u \in \cB(q_1,1/100)$, $u' \in \cB(q_1',1/100)$. Let
  $q_2 = g^{ij}_t u q_1'$. We write $q_2 = g_\tau u q_1$ for some $\tau >
0$, and let $q_2' = g_\tau u' q_1'$ (see
Figure~\ref{fig:outline2}). Also
suppose $u q_1 \in \pi^{-1}(K_0)$, $u' q_1' \in \pi^{-1}(K_0)$, 
$q_2 \in \pi^{-1}(K_0)$, $q_2' \in \pi^{-1}(K_0)$ and
\begin{displaymath}
C_1^{-1} \epsilon \le hd_{q_2}^{X_0}(U^+[q_2], U^+[q_2']) \le
C_1 \epsilon  \text{ and } 
\ell > \alpha_0 \tau. 
\end{displaymath}
In addition, suppose
there exist $\tilde{q}_2 \in \pi^{-1}(K_0)$ and $\tilde{q}_2' \in
\pi^{-1}(K_0)$ such that $\sigma_0(\tilde{q}_2') \in
W^+[\sigma_0(\tilde{q}_2)]$, and also $d^X(\tilde{q}_2, q_2) < \xi$ and
$d^X(\tilde{q}_2', q_2') < \xi$. 
Then, provided $\xi$ is small enough and $t$ is large enough (depending on $K_0$), 
\begin{equation}
\label{eq:goal:same:sheet}
\tilde{q}_2' \in W^+[\tilde{q}_2]. 
\end{equation}
Also, there exists $\xi''' > 0$ (depending on $\xi$, $K_0$ and $C$
and $t$)
with $\xi''' \to 0$ as $\xi \to 0$ and $t \to \infty$
such that 
\begin{equation}
\label{eq:nearby:maps:final}
d_*( P^+(\tilde{q}_2,\tilde{q}_2') f_{ij}(\tilde{q}_2),
f_{ij}(\tilde{q}_2')) \le \xi'''. 
\end{equation}
(In (\ref{eq:nearby:maps:final})  we
think of $f_{ij}(\tilde{q}_2')$ as a measure on $\gB_0[\tilde{q}_2']$,
$P^+(\tilde{q}_2,\tilde{q}_2') f_{ij}(\tilde{q}_2)$ as a measure on
$P^+(\tilde{q}_2,\tilde{q}_2') \gB_0[\tilde{q}_2]$, and we use the
AGY norm $\| \cdot \|_Y$ on $W^+(\tilde{q}_2')$
for the norm in the definition of $d_*(\cdot, \cdot)$.)
\end{proposition}

Proposition~\ref{prop:nearby:linear:maps} is proved in
\S\ref{sec:starredsubsec:proof:nearby:linear:maps}. We give an outline
of the argument below. 

\bold{Outline of the proof of
  Proposition~\ref{prop:nearby:linear:maps}.} 
The initial intuition behind the proof of
Proposition~\ref{prop:nearby:linear:maps} is that ``one
goes from $q_3'$ to $q_2'$ by nearly the same linear map as from $q_3$
to $q_2$; since this map is bounded on the relevant subspaces,
$f_{ij}(q_2)$ should be related to $f_{ij}(q_3)$ and $f_{ij}(q_2')$
should be related to $f_{ij}(q_2)$; since $f_{ij}(q_3)$ and
$f_{ij}(q_3')$ are close, $f_{ij}(q_2')$ should be related to
$f_{ij}(q_2)$.''

There are several problems with this argument. First, because of the
need to change transversals, there is no linear map from the space
$\GSpc(q_3)$ of generalized
subspaces near $q_3$ to the space $\GSpc(q_2)$ of generalized
subspaces near $q_2$. This difficulty is easily handled by working
instead with the linear maps $R(q_3,q_2): \bH(q_3) \to \bH(q_2)$ and
$R(q_3', q_2'): \bH(q_3') \to \bH(q_2')$. 

The second difficulty is connected to the first. We would like to say
that the two maps $R(q_3,q_2)$ and $R(q_3',q_2')$ 
are close, but the domains and ranges of the maps
are different. Thus we need ``connecting'' linear maps from
$\bH(q_3)$ to $\bH(q_3')$, and also from $\bH(q_2)$ to
$\bH(q_2')$. The first map is easy to construct: since $q_3$ and
$q_3'$ are in the same leaf of $W^-$, we can just use the linear map 
\index{$P$@$\bP^-(x,y)$}$\bP^-(q_3,q_3')$ induced by the ``$W^-$-connection map''
$P^-(q_3,q_3')$ defined in \S\ref{sec:subsec:connection}. 

Instead of constructing directly a map from $\bH(q_2)$ to
$\bH(q_2')$ we construct, using the choice of transversal $Z(\cdot)$,  
linear maps $\index{$P$@$\bP^{Z(x)}(x,y)$}\bP^{Z(q_2)}(q_2,\tilde{q}_2) : \bH(q_2) \to
\bH(\tilde{q}_2)$ and $\bP^{Z(q_2')}(q_2', \tilde{q}_2'):
\bH(q_2') \to \bH(\tilde{q}_2')$. 
Since $q_2$ and $\tilde{q}_2$ are close, and also since $q_2'$ and
$\tilde{q}_2'$ are close, these maps are in a suitable sense close to the
identity.  Then, since $\tilde{q}_2$ and
$\tilde{q}_2'$ are on the same leaf of $W^+$, we have the map
$\bP^+(\tilde{q}_2, \tilde{q}_2')$ induced by the $W^+$-connection
map $P^+(\tilde{q}_2,\tilde{q}_2')$ of \S\ref{sec:subsec:connection}. 

Thus, finally we have two maps from $\bH(q_3)$ to
$\bH(\tilde{q}_2')$:
\begin{displaymath}
\bA = \bP^+(\tilde{q}_2,\tilde{q}_2') \circ \bP^{Z(q_2)}(q_2,
\tilde{q}_2) \circ R(q_3,q_2)
\end{displaymath}
and
\begin{displaymath}
\bA' = \bP^{Z(q_2')}(q_2',\tilde{q}_2') \circ R(q_3',q_2') \circ
\bP^-(q_3,q_3').
\end{displaymath}
Even though $\bA$ and $\bA'$ are defined on $\bH(q_3)$, in what follows we
only need to consider their restrictions to $\bE_{[ij],bdd}(q_3) \subset
\bH(q_3)$; we will denote the restrictions by $\bB$ and $\bB'$
respectively. 

We would like to show that $\bB$ and $\bB'$ 
are close. By linearity, it is enough to show
that the restrictions of $\bB$ and $\bB'$ to each $\bE_{ij,bdd}(q_3)
\subset \bE_{[ij],bdd}(q_3)$ are close. Note that by
Proposition~\ref{prop:sublyapunov:locally:constant} (a), $\bP^-(q_3,q_3')
\bE_{ij,bdd}(q_3) = \bE_{ij,bdd}(q_3')$. Continuing this argument, we
see that the two subspaces  $\bB \,\bE_{ij,bdd}(q_3)$ and $\bB'\,
\bE_{ij,bdd}(q_3)$ are 
close to $\bE_{ij,bdd}(\tilde{q}_2')$ (and thus are close to each
other). Also, from the construction and
Proposition~\ref{prop:ej:bdd:transport:bounded}, we see that both
$\bB$ and $\bB'$ are uniformly bounded linear maps. However, this is
still not enough to conclude that $\bB$ and $\bB'$ are close. In fact
we also check that $\bB$ and $\bB'$ are close modulo
$\bV_{< i}(\tilde{q}_2)$. (This part of the argument uses the
assumptions on $q$, $q'$, $q_1$, $q_1'$, etc).  
Then we apply the elementary Lemma~\ref{lemma:lin:alg:B:Bprime}
  below with $E = \bE_{ij,bdd}(q_3)$, $L = \bH(\tilde{q}_2')$, 
$F = \bE_{ij,bdd}(\tilde{q}_2')$, $V = \bV_{< i}(\tilde{q}_2')$ to get
\begin{equation}
\label{eq:norm:bA:minus:bAprime}
\|\bB - \bB'\| \to 0 \text{ as } \xi \to 0. 
\end{equation}
The final part of the proof of
Proposition~\ref{prop:nearby:linear:maps} consists of deducing
(\ref{eq:nearby:maps:final}) from
(\ref{eq:norm:bA:minus:bAprime}) and the fact that $\bB$ and $\bB'$
are uniformly bounded
(Proposition~\ref{prop:ej:bdd:transport:bounded}). 

\begin{lemma}
\label{lemma:lin:alg:B:Bprime}
Suppose $L$ is a finite-dimensional 
normed vector space, $F$ and $V$ are subspaces of $L$, with
$F \cap V = \{0\}$. Let $S$ denote the unit sphere in $L$, and let
$hd(\cdot, \cdot)$ denote the Hausdorff distance induced by the norm
on $L$. 
Suppose $E$ is another finite-dimensional normed vector space, and
$\bB: E \to L$ and $\bB': E \to L$ are two linear maps each of norm at
most $C$. Let $\pi_V$ denote the projection $L \to L/V$.
Suppose $\xi > 0$ is such that
\begin{itemize}
\item[{\rm (i)}] $\| \pi_V \circ \bB - \pi_V \circ \bB' \| \le \xi$.
\item[{\rm (ii)}] $hd(\bB(E) \cap S, F \cap S) \le \xi$. 
\item[{\rm (iii)}] $hd(\bB'(E) \cap S,F \cap S) \le \xi$.
\end{itemize}
Then, $\|\bB - \bB' \| \le \xi'$,
where $\xi'$ depends on $\xi$, $C$ and the angle between $V$ and
$F$. Furthermore, $\xi' \to 0$ as $\xi \to 0$ (and the other
parameters remain fixed). 
\end{lemma}

In the course of the proof, we will prove the following lemma, which
will be used in \S\ref{sec:inductive:step}:

\begin{lemma}
\label{lemma:tau:ij:nearby:close}
For every $\delta > 0$ there exists a compact set $K_0 \subset X$ with
$\nu(K_0) > 1-\delta$ such that the following holds:
Suppose $x,x',y,y' \in \pi^{-1}(K_0)$, $y \in W^+[x]$, $y'
  \in W^+[x']$ and $x' \in W^-[x]$. Suppose further that $d^{X_0}(x,y)
  \le 1/100$, $d^{X_0}(y,y') \le 1/100$, 
and that there exists $s > 0$ such that 
for all $|\tau| \le s$, $d^{X_0}(g_\tau x , g_\tau x') \le 1/100$ and
$d^{X_0}(g_\tau y, g_\tau y') \le 1/100$. 
Furthermore, suppose $0 < \alpha_0 < 1$ and 
that $0 < t < \alpha_0^{-1} s$ is such
that $d^{X_0}(g_t y, g_t y') < 1/100$, $g_t y \in K_0$
and $g_t y' \in K_0$. Then, for
all $ij \in \Lambda''$, 
\begin{equation}
\label{eq:tau:ij:nearby:close}
|\hat{\tau}_{ij}(y,t) - \hat{\tau}_{ij}(y',t)| \le C,
\end{equation}
where $C$ depends only on $\delta$, $\alpha_0$ and the Lyapunov spectrum.
\end{lemma}

\starredsubsection{Proof of Lemma~\ref{lemma:Eij:bdd:locally:constant}.}
\label{sec:starredsubsec:proof:lemma:Eij:bdd:locally:constant}

\mcc{make sure to convert correctly between hodge and dynamical norms here}

Let $\theta_1 > 0$ and $\delta > 0$ be small
  constants to be chosen later.
Let $K \subset X$ and $C > 0$ be such
that $\nu(K) > 1-\delta$, for $x \in K$ the 
Lemma~\ref{lemma:hausdorff:distance:to:norm} (b) holds with $c_1(x) >
C^{-1}$, and for all $x \in K$, all $\bfv \in \bE_{[ij],bdd}(x)$ and all
$\ell > 0$, for at least $(1-\theta_1)$ fraction of $y \in \cF_{ij}[x,\ell]$, 
\begin{equation}
\label{eq:tmp:redef:bounded}
\|R(x,y) \bfv\| < C \|\bfv\|. 
\end{equation}
By Lemma~\ref{lemma:fake:ergodicity:Fij} there exists a subset $K^*
\subset K$ with $\nu(K^*) \ge (1-2\kappa^2 \delta^{1/2})$ such that for
$x \in K^*$, (\ref{eq:fake:ergodicity:Fij}) holds with $\theta' =
\delta^{1/2}$. Furthermore, we may ensure that for $x \in K^*$, 
$K^* \cap \cF_{ij}[x]$ is relatively open in $\cF_{ij}[x]$. (Indeed,
suppose $z \in \cF_{ij}[x]$ is near $x \in K^*$. Then, there exists
$\ell_0$ such that 
for $\ell > \ell_0$, $\cF_{ij}[x,\ell] = \cF_{ij}[z,\ell]$
and thus (\ref{eq:fake:ergodicity:Fij}) holds for $z$. For $\ell <
\ell_0$, (\ref{eq:fake:ergodicity:Fij}) holds for $z$ sufficiently
close to $x$ by continuity.) Let
\begin{displaymath}
\Psi  = \{ x \in X \st \lim_{T \to \infty} |\{ t \in [0,T] \st
  g_{-t} x \in K^* \}| \ge (1-2 \kappa^2 \delta^{1/2}). 
\end{displaymath}
Then $\nu(\Psi) = 1$. From its definition, $\Psi$ is invariant
under $g_t$. Since $K^* \cap \cF_{ij}[x]$ is
relatively open in $\cF_{ij}[x]$, $\Psi$ is saturated by the leaves
of $\cF_{ij}$. This implies that $\Psi$ is (locally) invariant under
$U^+$.  Now, let 
\begin{displaymath}
K_N = \{ x \in \Psi \st \text{ for all $T > N$, } \quad |\{ t \in [0,T] \st
  g_{-t} x \in K^* \}| \ge (1-4\kappa^2 \delta^{1/2})T \}. 
\end{displaymath}
(We may assume that $4\kappa^2 \delta^{1/2} \ll 1$.) 
We have $\bigcup_{N} K_N = \Psi$. 

\makefig{Proof of Lemma~\ref{lemma:Eij:bdd:locally:constant}}{fig:equivalence1}{\includegraphics{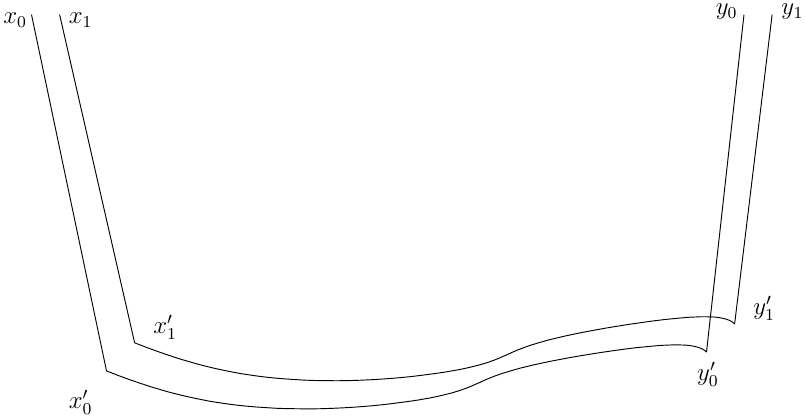}}

\makefignocenter{Proof of Lemma~\ref{lemma:Eij:bdd:locally:constant}
\noindent\smallskip

In (b), the subspaces $U^+[y_0']$ and $U^+[y_1']$ stay close since
$x_1' \in \cE_{ij}(x_0')$, and also for $k\in \{0,1\}$, the subspaces
$R(x_k',y_k') \cQ_k'$ and $U^+[y_k']$ stay close since $\cQ_k' \in
\cE_{ij,bdd}(x_k')$. }
{fig:equivalence2}{\includegraphics{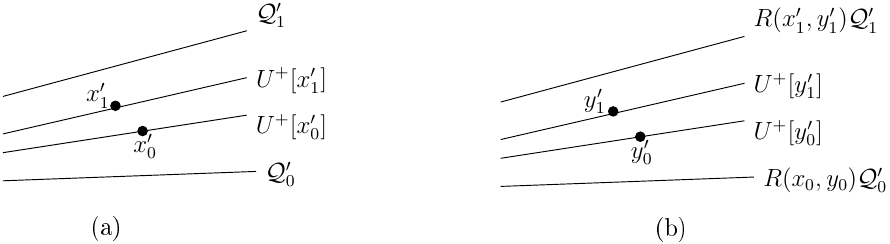}}

Suppose $x_0 \in K_N$, $x_1 \in \gB_0[x_0] \cap K_N$, so $d^{X_0}(x_0,x_1) < 1/100$. 
For $k=0,1$, let $\cQ_k \subset \cE_{ij}[x_k]$ be such that
\begin{displaymath}
hd_{x_k}^{X_0}(\cQ_k,U^+[x_k]) \le 1/100, 
\end{displaymath}
and the vector
\begin{displaymath}
\bfv_k = \bfj(\cU_{x_k}^{-1}(\cQ_{1-k})) 
\end{displaymath}
satisfies $\|\bfv_k\| \le 1/100$.

We claim that $\bfv_k \in \bH(x_k)$. Indeed, we may write 
$\cU_{x_{1-k}}^{-1}(\cQ_{1-k}) = (M_{1-k}, v_{1-k})$. Also we may write 
$\cU_{x_k}^{-1}(U^+[x_{1-k}]) = (M'_k, v_k')$. 
Then, $\cQ_{1-k}$ is
parametrized (from $x_k$) 
by a pair $(M_k'', w_k)$ where $w_k \in W^+(x_k)$, and 
\begin{displaymath}
M_k'' = (I+M_{1-k}) \circ (I+ M_k')-I
\end{displaymath}
(This parametrization is not
necessarily adapted to $Z(x_k)$.) 
Since $M_{1-k}$ and $M'_k$ are
both in $\cH_{++}$, $M_k'' \in \cH_{++}(x_k)$. Thus, $\bfv_k =
\bS_{x_k}(\bfj(M_k'',w_k)) \in \bH(x_k)$. 

For $C_1(N)$ sufficiently large, we can find
$C_1(N) < t < 2C_1(N)$ such that
$x_0' \equiv g_{-t}^{ij} x_0 \in K^*$, 
$x_1' \equiv g^{ij}_{-t} x_1 \in K^*$.  By
Lemma~\ref{lemma:stay:in:same:orbit}, $x_1' \in \gB_0[x_0']$. 
Let
$\bfv_k' = g_{-t}^{ij} \bfv_k$, $\cQ_k' = g_{-t}^{ij} \cQ_k$. 
By choosing $C_1(N)$ sufficiently large (depending on $N$), we can
ensure that
\begin{displaymath}
hd_{x_k'}^{X_0}(U^+[x_k'],U^+[x_{1-k}']) \le C^{-3}, \qquad
hd_{x_k'}^{X_0}(\cQ_k',U^+[x_k']) \le C^{-3}. 
\end{displaymath}
By Lemma~\ref{lemma:hausdorff:distance:to:norm}, since $x_k' \in K$, 
\begin{equation}
\label{eq:estimates:jplus:xk:prime}
\|\bfj(\cU_{x_k'}^{-1}(U^+[x_{1-k}']))\| \le C^{-2}, \qquad
\|\bfj(\cU_{x_k'}^{-1}(\cQ_k'))\| \le C^{-2}. 
\end{equation}
Let $\ell > 0$ be arbitrary, and let $\ell'$ be such that $g_t^{ij}
\cF_{ij}[x, \ell'] = \cF_{ij}[x,\ell]$. 
Then, for $k=0,1$, since $x_k' \in K^*$, 
\begin{displaymath}
|\{ y_k' \in \cF_{ij}[x_k',\ell'] \st y_k' \in K\}| \ge (1-\delta^{1/2}) 
|\cF_{ij}[x',\ell']|,
\end{displaymath}
Since $U^+[x_1] \in \cE_{ij}[x_0]$, we have $U^+[x_1'] \in
\cE_{ij}[x_0']$, and thus $\bfj(\cU_{x_0'}^{-1}(U^+[x_1'])) \in
\bE_{[ij],bdd}(x_0')$. Since $x_0' \in K$, we have by
(\ref{eq:tmp:redef:bounded}), 
for at least $(1-\theta_1)$-fraction of $y_0' \in
\cF_{ij}[x_0', \ell']$,  
\begin{equation}
\label{eq:j:Uplus:y0:close:j:Uplus:y1}
\|R(x_0',y_0') \bfj(\cU_{x_0'}^{-1}(U^+[x_1']))\| \le C
\|\bfj(\cU_{x_0'}^{-1}(U^+[x_1']))\| \le C^{-1}, 
\end{equation}
where we have used (\ref{eq:estimates:jplus:xk:prime}) 
for the last estimate. 
Let $\theta'' = 2\theta_1 + 2\delta^{1/2}$. Then, for
at least $1-\theta''/2$ 
fraction of $y_0' \in \cF_{ij}[x_0',\ell']$, $y_0' \in K$ and
(\ref{eq:j:Uplus:y0:close:j:Uplus:y1}) holds. 
Therefore, by Lemma~\ref{lemma:hausdorff:distance:to:norm}, 
for at least $(1-\theta''/2)$-fraction of $y_0' \in \cF_{ij}[x_0',
\ell']$, for a suitable $y_1' \in \cF_{ij}[x_1',\ell']$,
\begin{equation}
\label{eq:hd:U:y0:prime:U:y1:prime:bounded}
hd_{y_0'}^{X_0}(U^+[y_0'],U^+[y_1']) \le 1/100. 
\end{equation}
Also, since $\cQ_k \in \cE_{ij}[x_k]$, $\cQ_k' \in \cE_{ij}[x_k']$,
and thus $\bfj(\cU_{x_k'}^{-1}(\cQ_k')) \in \bE_{[ij],bdd}(x_k')$. 
Hence, by (\ref{eq:tmp:redef:bounded}), 
for at least $(1-\theta)$-fraction of $y_k' \in \cF_{ij}[x_k', \ell']$, 
\begin{equation}
\label{eq:j:Uplus:y0:close:image:Qk:prime}
\|R(x_k',y_k') \bfj(\cU_{x_k'}^{-1}(\cQ_k'))\| \le C
\| \bfj(\cU_{x_k'}^{-1}(\cQ_k'))\| \le C^{-1}. 
\end{equation}
where we used (\ref{eq:estimates:jplus:xk:prime}) for the last
estimate. Then, for at least $(1-\theta''/2)$-fraction of $y_k' \in
\cF_{ij}[x_k',\ell']$, $y_k' \in K$ and  
(\ref{eq:j:Uplus:y0:close:image:Qk:prime}) holds. Therefore, by
Lemma~\ref{lemma:hausdorff:distance:to:norm}, 
for at least $(1-\theta''/2)$-fraction of $y_k' \in \cF_{ij}[x_k', \ell']$, 
\begin{displaymath}
hd_{y_k'}^{X_0}(U^+[y_k'],R(x_k',y_k')\cQ_k') \le 1/100. 
\end{displaymath}
Therefore, by (\ref{eq:hd:U:y0:prime:U:y1:prime:bounded}), 
 for at least $(1-\theta'')$-fraction of
$y_k' \in \cF_{ij}[x_k',\ell']$, for a suitable $y_{1-k}' \in
\cF_{ij}[x_{1-k}',\ell']$, 
\begin{equation}
\label{eq:final:hd:bound}
hd_{y_k'}^{X_0}(U^+[y_k'],R(x_{1-k}',y_{1-k}') \cQ_{1-k}') \le 1/50.
\end{equation}
Let 
\begin{displaymath}
\bfw_k' = \bfj(\cU_{y_k'}^{-1}(R(x_{1-k}',y_{1-k}') \cQ_{1-k}')) =
R(x_k',y_k') \bfv_k'.  
\end{displaymath}
Then, assuming $y_0' \in K$ and (\ref{eq:final:hd:bound}) holds, 
by Lemma~\ref{lemma:hausdorff:distance:to:norm}, 
\begin{displaymath}
\|\bfw_k'\| \le C. 
\end{displaymath}
Let $y_k = g_t^{ij} y_k'$, and let 
\begin{displaymath}
\bfw_k = R(y_k',y_k) \bfw_k' = R(x_k,y_k) \bfv_k. 
\end{displaymath}
Then, for at least $(1-\theta'')$-fraction of $y_k \in
\cF_{ij}[x_k,\ell]$, $\|R(x_k,y_k)\bfv_k\| \le C_2(N)$. This
implies, by Proposition~\ref{prop:Rxy:v:small:implied:sync:bounded}, that
$\bfv_k \in \bE_{[ij],bdd}(x_k)$. (By making
$\theta_1 > 0$ and $\delta > 0$ sufficiently
small, we can make sure that $\theta'' < \theta$ where $\theta > 0$ is
as in Proposition~\ref{prop:Rxy:v:small:implied:sync:bounded}.) 

Thus, for all $\cQ_k \in \cE_{ij}[x_k]$ such that
$\bfj(\cU_{x_{1-k}}^{-1}(\cQ_k)) \le 1/100$, 
we have $\bfj(\cU_{x_{1-k}}^{-1}(\cQ_k)) \in
\bE_{[ij],bdd}(x_{1-k})$. Since both $\cU_{x_{1-k}}^{-1}$ and $\bfj$ are
analytic, this implies that 
$\bfj(\cU_{x_{1-k}}^{-1}(\cQ_k)) \in
\bE_{[ij],bdd}(x_{1-k})$ for all $\cQ_k \in \cE_{ij}[x_k]$. Thus, for $k=0,1$,
$\cE_{ij}[x_k] \subset \cE_{ij}[x_{1-k}]$. This implies that
$\cE_{ij}[x_0] = \cE_{ij}[x_1]$.  
\qed\medskip

\starredsubsection{Proof of
  Proposition~\ref{prop:nearby:linear:maps}.}
\label{sec:starredsubsec:proof:nearby:linear:maps}

Let $\cO \subset X$ be an open set contained in the fundamental
domain, and let $x \to u_x \in U^+(x)$ 
be a function which is  constant
on each set of the form $U^+[x] \cap \cO$. Let $T_u:\cO \to X$ be the
map which takes $x \to u_x x$. 

\begin{lemma}
\label{lemma:Tu:measure:preserving}
Suppose $E \subset \cO$. Then $\nu(T_u(E)) = \nu(E)$. 
\end{lemma}

\mcc{rewrite the following proof}

\bold{Proof.} Without loss of generality, we may assume that $T_u(\cO)
\cap \cO = \emptyset$. 
For each $x \in \cO$, let $\tilde{U}[x]$ be a finite piece
of $U^+[x]$ which contains both $U[x] \cap \cO$ and $T_u(U[x] \cap
\cO)$. We may assume that $\tilde{U}[x]$ is the same for all $x \in
U[x] \cap \cO$. 
Let $\tilde{\cU}$ be the $\sigma$-algebra of functions which are
constant along each $\tilde{U}[x]$. 
Then, for any measurable $\phi: X \to \reals$,
\begin{displaymath}
\int_X \phi \, d\nu = \int_X \mathbb{E}(\phi \mid \tilde{\cU}) \, d\nu
\end{displaymath}
Now suppose $\phi$ is supported on $\cO$. 
We have $\mathbb{E}(\phi \circ T_u \mid \tilde{\cU}) = \mathbb{E}(\phi \mid
\tilde{\cU})$ since the
conditional measures along $U^+$ are Haar, and $T_u$ restricted to
$\cO \cap U^+[x]$ is a translation. \mcc{explain} 
Thus
\begin{displaymath}
\int_X \phi \circ T_u \, d\nu = \int_X \mathbb{E}(\phi \circ T_u \mid \cU) \,
d\nu = \int_X \mathbb{E}(\phi \mid \tilde{\cU}) \, d\nu = \int_X \phi \, d\nu.
\end{displaymath}
\qed\medskip

We also recall the following standard fact:
\begin{lemma}
\label{lemma:change:of:var:conditional:measure}
Suppose $\Psi: X \to X$ preserves $\nu$, and also for
almost all $x$, $\cC_{ij}[\Psi(x)] \cap \gB_0[\Psi(x)] \cap  \Psi(\gB_0[x])
= \Psi(\cC_{ij}[x]) \cap \gB_0[\Psi(x)] \cap \Psi(\gB_0[x])$. Then, 
\begin{displaymath}
f_{ij}(\Psi(x)) \propto \Psi_* f_{ij}(x),
\end{displaymath}
in the sense that the restriction of both measures to the set $\gB_0[\Psi(x)]
\cap \Psi(\gB_0[x])$ where both make sense is the same up to
normalization. 
\mcc{write this more formally} 
\end{lemma}

\bold{Proof.} See \cite[Lemma 4.2(iv)]{Einsiedler:Lindenstrauss:Pisa}.
\qed\medskip

\begin{lemma}
\label{lemma:change:cond:meas}
We have (on the set where both are defined):
\begin{displaymath}
f_{ij}(g_t T_u g_{-s} x) \propto (g_t T_u g_{-s})_* f_{ij}(x). 
\end{displaymath}
\end{lemma}

\bold{Proof.} This follows immediately from
Lemma~\ref{lemma:Tu:measure:preserving} and
Lemma~\ref{lemma:change:of:var:conditional:measure}. 
\qed\medskip

\bold{The maps $\phi_x$.}
We have the map $\phi_x: W^+(x) \to \cH_{++}(x) \cross
W^+(x)$ given by
\begin{equation}
\label{eq:def:phi:x}
\index{$\phi_x(z)$}\phi_x(z) = \cU_x^{-1}(U^+[z]).
\end{equation}
(Here $\cU_x^{-1}$ is defined using the transversal $Z(x)$.)

Suppose $Z(x)$ is an admissible transversal to $U^+(x)$. 
Since $f_{ij}(x)$ is Haar
along $U^+$, we can recover $f_{ij}(x)$ from its restriction to
$Z(x)$. More precisely, the following holds:

Let $\pi_2: \cH_{++}(x) \cross W^+(x) \to W^+(x)$ be projection onto the second
factor. Then, for $z \in Z(x)$, $\pi_2(\phi(z)) = z$. 
Now, suppose $Z'$ is another transversal to $U^+(x)$. Then, 
\begin{displaymath}
(f_{ij}\mid_{Z'})(x) = (\pi_2 \circ S_x^{Z'} \circ \phi)_*
(f_{ij}\mid_{Z(x)}). 
\end{displaymath}

\bold{The measures $\bff_{ij}(x)$.}
Let 
\begin{displaymath}
\index{$f$@$\bff_{ij}(x)$}\bff_{ij}(x) = (\bfj \circ \phi_x)_* f_{ij}(x). 
\end{displaymath}
Then, $\bff_{ij}(x)$ is a measure on $\bH(x)$. 

\begin{lemma}
\label{lemma:transformation:rule:bff}
For $y \in \cF_{ij}[x]$, we have (on the set where both
  are defined), 
\begin{displaymath}
\bff_{ij}(y) \propto R(x,y)_* \bff_{ij}(x). 
\end{displaymath}
\end{lemma}

\bold{Proof.} 
Suppose $t > 0$ is such that $x' = g_{-t}^{ij} x$ and $y' =
g_{-t}^{ij} y$ satisfy $y' \in \cB[x']$. Working in the universal
cover, let $Z[x] = \{ z \st z - x \in Z(x) \}$. 
Let $Z[x'] = g_{-t}^{ij}
Z[x]$, and let $Z[y'] = g_{-t}^{ij} Z[y]$. For $z \in Z[x']$ near $x'$,
let $u_z$ be such that $u_z z \in Z[y']$. We extend the function $z
\to u_z$ to be locally constant along $U^+$ in a neighborhood of
$Z[x']$. Then, let  
\begin{displaymath}
\Psi = g_{t}^{ij} \circ T_u \circ g_{-t}^{ij}. 
\end{displaymath}
\mcc{worry about whether time change is locally constant along small
  neighborhood, so measure is preserved. Seems ok.}

Note that $\Psi$ takes $Z[x]$ into $Z[y]$, and 
by Lemma~\ref{lemma:change:cond:meas}, 
\begin{equation}
\label{eq:Psi:star:fij}
\Psi_* f_{ij}(x) \propto f_{ij}(y). 
\end{equation}
By the definition of $u_*$ in \S\ref{sec:divergence:subspaces}, 
for $z \in Z[x]$, 
\begin{displaymath}
(R(x,y) \circ \bfj \circ \cU_x^{-1}) U^+[z] = (\bfj \circ
\cU_y^{-1}) U^+[\Psi(z)]. 
\end{displaymath}
Hence, by (\ref{eq:def:phi:x}),
\begin{equation}
\label{eq:Rxy:circ:bfjplus:circ:phi:x}
(R(x,y) \circ \bfj \circ \phi_x) (z) = (\bfj \circ
\phi_y  \circ \Psi)(z),
\end{equation}
where $\phi_y$ is relative to the transversal $Z(y)$ and $\phi_x$
is relative to the transversal $Z(x)$.  
(Here we have used the fact that $\Psi (U^+[z]) = U^+[\Psi(z)]$ which
follows from the equivariance of $U^+$. 
Also, in (\ref{eq:Rxy:circ:bfjplus:circ:phi:x}), $R(x,y)$ is as in
\S\ref{sec:foliations}.)
\mcc{Maybe this needs extra proof; change transversals from $Z[x']$ to
  $Z[y']??$. Seems ok.}
Now the lemma follows from (\ref{eq:Psi:star:fij}) and
(\ref{eq:Rxy:circ:bfjplus:circ:phi:x}).  
\qed\medskip


Let $P^+(x,y)$ and $P^-(x,y)$ be as in \S\ref{sec:subsec:connection}. 
The maps $P^+(x,y)_*: \Lie(\cG_{++})(x) \to \Lie(\cG_{++})(y)$ 
(where we use the notation (\ref{eq:def:star:notation})) are an
equivariant measurable flat $W^+$-connection on the bundle
$\Lie(\cG_{++})$ satisfying (\ref{eq:Fxy:preserves:lyapunov:subspaces}). 
Then, by Proposition~\ref{prop:sublyapunov:locally:constant}(a),
\begin{equation}
\label{eq:Pplus:preserves:U}
P^+(x,y)_* \Lie(U^+)(x) = \Lie(U^+)(y). 
\end{equation}

\bold{The maps $\bP^+(x,y)$ and $\bP^-(x,y)$.}
In view of (\ref{eq:Pplus:preserves:U}), 
the maps $P^+(x,y)$ naturally induce a linear map (which we denote by
$\tilde{\bP}^+(x,y)$) from $\tilde{\bH}(x)$ to $\tilde{\bH}(y)$, so
that for $(M,v) \in \cH_{++}(x)$, 
\begin{displaymath}
\index{$P$@$\tilde{\bP}^+(x,y)$}\tilde{\bP}^+(x,y) \circ \bfj(M,v) = \bfj( P^+(x,y) \circ M \circ
P^+(x,y)^{-1}, P^+(x,y) v). 
\end{displaymath}
Let $\bP^+(x,y) = \bS_y^{Z(y)} \circ \tilde{\bP}^+(x,y)$. Then
the maps $\bP^+(x,y): \bH(x) \to \bH(y)$ are an
equivariant measurable flat $W^+$-connection on the bundle
$\bH$ satisfying (\ref{eq:Fxy:preserves:lyapunov:subspaces}). 
Then, by Proposition~\ref{prop:sublyapunov:locally:constant}(a),
we have 
\begin{equation}
\label{eq:bPplus:preserves:bEijbdd}
\bP^+(x,y) \bE_{ij,bdd}(x) = \bE_{ij,bdd}(y).  
\end{equation}
For $y \in W^-[x]$, we have a map $\bP^-(x,y)$ with analogous
properties. 

\bold{The maps $P^Z(x,y)$ and $\bP^Z(x,y)$.}
We also need to define a map between $\bH(x)$ and $\bH(y)$ even if
$x$ and $y$ are not on the same leaf of $W^+$ or $W^-$. For every $v_i
\in \cV_i(x) \equiv \cV_i(H^1)(x)$, and $i \in \Lambda$ (where
  $\Lambda$ is the  
Lyapunov spectrum) \mcc{(we also include the orbit direction here)}
we can write
\begin{displaymath}
v_i = v_i' + v_i'' \quad v_i'  \in \cV_i(H^1)(y), \quad v_i'' \in \bigoplus_{j
  \ne i} \cV_i(H^1)(y).
\end{displaymath}
Let
\index{$P$@$P^\sharp(x,y)$}$P^\sharp(x,y): H^1(x) \to H^1(y)$ be the linear map whose restriction
to $\cV_i(H^1)(x)$ sends $v_i$ to $v_i'$. By definition, $P^\sharp(x,y)$
sends $\cV_i(H^1)(x)$ to $\cV_i(H^1)(y)$, but it is not clear that
$P^\sharp(x,y)_* \Lie(U^+)(x) = \Lie(U^+)(y)$. 
To correct this, given a Lyapunov-adapted transversal
$Z(x)$, note that (for $y$ near $x$),
\begin{displaymath}
\Lie(\cG_{++})(x) = P^\sharp(x,y)^{-1}_* \Lie(U^+)(y) \oplus Z(x).
\end{displaymath}
Then, given $v \in \Lie(U^+)(x) \subset \Lie(\cG_{++})(x)$, we can decompose
\begin{equation}
\label{eq:def:Mxy:one}
v = v' + v'', \qquad v' \in  P^\sharp(x,y)^{-1}_* \Lie(U^+)(y), \qquad
v'' \in Z(x). 
\end{equation}
Define $\index{$M(x;y)$}M(x;y): \Lie(U^+)(x) \to 
\Lie(\cG_{++})(x)$ by 
\begin{equation}
\label{eq:def:Mxy:two}
M v = - v''. 
\end{equation}
Then, since $Z(x)$ is Lyapunov
adapted, $M(x;y): \Lie(U^+)(x) \to \Lie(\cG_{++})(x)$ 
is the linear map such that
\begin{equation}
\label{eq:def:Mxy}
(I + M(x;y)) \Lie(U^+)(x) = P^\sharp(x,y)^{-1}_* \Lie(U^+)(y), 
\end{equation}
and $M(x;y)\cV_i(\Lie(U^+))(x) \subset Z_i(x)$,
where $Z_i(x) = Z(x)
\cap \cV_i(\Lie(\cG_{++}))(x)$ is as in
\S\ref{sec:divergence:subspaces}.  
Then, let $P^{Z(x)}(x,y): \cH_{++}(x) \to \cH_{++}(y)$ be the map
taking $f \in \cH_{++}(x)$ to 
\begin{displaymath}
\index{$P^{Z(x)}(x,y)$}P^{Z(x)}(x,y) f \equiv P^\sharp(x,y)_* \circ f \circ (I + M(x;y))^{-1} \circ P^\sharp(x,y)_*^{-1} 
\in \cH_{++}(y).
\end{displaymath}
Then, since $M(x;y) \cV_i(\Lie(U^+))(x) \subset
\cV_i(\Lie(\cG_{++}))(x)$ we have for a.e.\ $x$, $y$,
\begin{displaymath}
P^{Z(x)}(x,y) \cV_i(\cH_{++})(x) = \cV_i(\cH_{++})(y). 
\end{displaymath}
Then $P^{Z(x)}$ gives a map $\index{$\tilde{P}^{Z(x)}(x,y)$}\tilde{P}^{Z(x)}(x,y): \cH_{++}(x) \times W^+(x) \to
\cH_{++}(y) \times W^+(y)$ given by \mcc{check}
\begin{displaymath}
\tilde{P}^{Z(x)}(x,y)(f,v) = (P^{Z(x)}(x,y)f, P^\sharp(x,y) v).
\end{displaymath}
Therefore, (after possibly composing with a change in transversal map
$\bS$) $\tilde{P}^{Z(x)}(x,y)$ induces a map we will call \index{$P$@$\bP^{Z(x)}(x,y)$}$\bP^{Z(x)}(x,y)$ between
$\bH(x)$ and $\bH(y)$. This map satisfies
\begin{equation}
\label{eq:bPZ:preserves:eigenspaces}
\bP^{Z(x)}(x,y) \cV_i(\bH)(x) = \cV_i(\bH)(y),
\end{equation}
and has the equivariance property
\begin{displaymath}
\bP^{g_{-t} Z(x)}(g_{-t} x, g_{-t} y) = g_{-t} \circ \bP^{Z(x)}(x,y)
\circ g_t.   
\end{displaymath}
\begin{lemma}
\label{lemma:Psharp:is:Pplus}
For $y \in W^+[x]$, and any choice of $Z(x)$, we have
\begin{equation}
\label{eq:PZ:is:Pplus}
\bP^{Z(x)}(x,y) = \bP^+(x,y).
\end{equation}
\end{lemma}
\bold{Proof.} Suppose $y \in W^+[x]$. Then by
Lemma~\ref{lemma:properties:lyapunov:flag}, $P^\sharp(x,y) =
  P^+(x,y)$, thus 
\begin{displaymath}
P^\sharp(x,y)^{-1}_* \Lie(U^+)(y) = P^+(x,y)^{-1}_*\Lie(U^+)(x) = \Lie(U^+)(x)
\end{displaymath}
where for the last equality we used
Proposition~\ref{prop:sublyapunov:locally:constant}(a). Hence,
$M(x;y) = 0$ and (\ref{eq:PZ:is:Pplus}) follows.
\qed\medskip

\mcc{try to merge the following argument with the one in
  \S\ref{sec:divergence:subspaces}}

\begin{lemma}
\label{lemma:P:sharp}
For any $\delta > 0$ there exists a compact subset $K \subset X_0$ with
$\nu(K) > 1-\delta/2$ such that the following holds:
Suppose $x$ and $y \in \pi^{-1}(K)$, and $s > 0$ are such that for all
$|t| < s$, 
$d^{X_0}(g_t x, g_t y) < 1/100$.  
Then, there exists $\alpha > 0$ depending only
on the Lyapunov spectrum, and $C = C(\delta)$ such that for all $i$, 
\begin{displaymath}
d_Y(P^{GM}(x,y) \cV_i(H^1)(x), \cV_i(H^1)(y)) \le C(\delta) e^{-\alpha s}. 
\end{displaymath}
\end{lemma}

\bold{Proof.} There exists a compact subset $K_1 \subset X_0$ such that the
functions $x \to \cV_i(H^1)(x)$ are uniformly continuous. (Here we are
using the Gauss-Manin connection to identify $H_1(x)$ with $H^1(y)$
for $y$ near $x$). Then, there exists $\sigma > 0$ such that if $x \in
\pi^{-1}(K_1)$, $y \in \pi^{-1}(K_1)$ and $d^{X_0}(x,y) < \sigma$ then  
$D^(x,y) < 1$ and $D^-(x,y) < 1$.
(See
\S\ref{sec:subsec:estimates:Lyapunov:subspaces} for the definition of
$D^\pm(\cdot, \cdot)$). We also may assume that there exists a
constant $C_0(\delta)$ such that $C(x) < C_0(\delta)$ for all $x \in
K_1$, where $C(\cdot)$ is as in
Lemma~\ref{lemma:subspaces:stay:close}. 
Then there exists a compact subset $K \subset X$ with 
$\nu(K) > 1-\delta$, and $t_0 > 0$ such that for $x \in K$, for $t >
t_0$, for $(1-\delta)$-fraction of $t \in [0,s]$, $g_t x \in K_1$,
$g_{-t} x \in K_1$ also for at least half the fraction of $t \in
[0,s]$, $g_t x$ and $g_{-t} x$ belong to $K_{thick}$ where $K_{thick}$
is as in Lemma~\ref{lemma:forni}. 

Suppose $x \in \pi^{-1}(K)$, and $y \in \pi^{-1}(K)$. 
Then, by Lemma~\ref{lemma:forni}, there exists $\alpha_1 > 0$
depending only on the Lyapunov spectrum such that there exists 
$t \in [\alpha_1 s, s]$ with $g_t x \in K_1$, $g_t y \in K_1$ and
$d^{X_0}(g_t x, g_t y) < \sigma$. Then, $D^-(g_t x, g_t y) < 1$. Then, 
by Lemma~\ref{lemma:subspaces:stay:close}, applied to the points $g_t
x$, $g_t y$, we get
\begin{displaymath}
d_Y(\cV_{\ge i}(H^1)(x), \cV_{\ge i}(H^1)(y)) \le C(\delta) e^{-\alpha
  t} = C(\delta) e^{-\alpha \alpha_1 s}.
\end{displaymath}
Similarly, there exists $t \in [\alpha_1 s, s]$ with $g_{-t} x \in
K_1$ and $g_{-t} y \in K_1$. Then, we get
\begin{displaymath}
d_Y(\cV_{\le i}(H^1)(x), \cV_{\le i}(H^1)(y)) \le C(\delta) e^{-\alpha
  t} = C(\delta) e^{-\alpha \alpha_1 s}.
\end{displaymath}
The lemma follows. 
\qed\medskip

For every $\delta > 0$ and every $0 < \alpha < 1$
there exist
compact sets $K_0 \subset K^\sharp \subset X$ with $\nu(K_0) > 1-\delta$ such
that the following hold: 
\begin{itemize}
\item[($K^\sharp$1)]  
The functions $U^+(x)$, $\cV_i(H^1)(x)$ and more generally,
$\cV_i(H_{big})(x))$ for all $i$,
are uniformly continuous on $K^\sharp$.
\item[($K^\sharp$2)]  
The functions $Z(x)$ are uniformly continuous on $K^\sharp$. 
\item[($K^\sharp$3)] 
The functions $\bE_{ij,bdd}(x)$ are uniformly continuous on $K^\sharp$. 
\item [($K^\sharp$4)] 
The functions $f_{ij}(x)$ and $\bff_{ij}(x)$ are uniformly
  continuous on $K^\sharp$ (in the weak-* convergence topology).   
\item [($K^\sharp$5)] 
There exists $t_0 > 0$ and $\epsilon' < 0.25 \alpha \min_{i \ne j}|\lambda_i -
  \lambda_j|$ such that for $t > t_0$, $x \in K^\sharp$, all $i$, and
  any $v \in \cV_i(H^1)(x)$, 
\begin{displaymath}
e^{(\lambda_i - \epsilon')t} \|v\| \le \|(g_t)_* v\| \le e^{(\lambda_i + \epsilon')t} \|v\|
\end{displaymath}
\item[($K^\sharp$6)] The function $C_3(\cdot)$ of
  Proposition~\ref{prop:ej:bdd:transport:bounded} is uniformly bounded
  on $K^\sharp$. 
\item[($K^\sharp$7)] $\bE_{ij,bdd}(x)$ and $\bV_{< i}(x)$ are
  transverse for $x \in K^\sharp$.
\item[($K^\sharp$8)] $K^\sharp \subset K_{thick}''$ where
  $K_{thick}''$ is as in Lemma~\ref{lemma:forni} (c).  Also $K^\sharp
  \subset K$ where $K$ is as in Lemma~\ref{lemma:P:sharp}. 
\item[($K^\sharp$9)] There exists $c_0(\delta) > 0$ with $c_0(\delta) \to
  0$ as $\delta \to 0$ such that for all $x \in K^\sharp$, 
$d^{X_0}(x,\partial \gB_0[x]) >
  c_0(\delta)$ where $\gB_0[x]$ is as in \S\ref{sec:semi:markov}.
\item[($K^\sharp$10)] There exists a constant $C_4(\delta)$ such that
  for all $x \in K^\sharp$ and all $v \in H_{big}(x)$,
  $C_4(\delta)^{-1} \| v\| \le \|v\|_Y \le C_4(\delta) \|v\|$. 
\item[($K^\sharp$11)] There exists a constant $C_1 = C_1(\delta) < \infty$
  such that for $x  \in K_0$ and all $T > C_1(\delta)$ and all $ij$ we have 
\begin{displaymath}
|\{ t \in [C_1, T] \st g^{ij}_{-t} x \in K^\sharp \}| \ge 0.99 (T-C_1).
\end{displaymath}
\end{itemize}

\begin{lemma}
\label{lemma:Psharp:easybounds}
Suppose $x,x',y,y' \in \pi^{-1}(K_0)$, $y \in W^+[x]$, $y'
  \in W^+[x']$ and $x' \in W^-[x]$. Suppose further that $d^{X_0}(x,y) <1/100$, $d^{X_0}(y,y') <1/100$, and that there exists $s > 0$ such that 
for all $|t| \le s$, $d^{X_0}(g_t x , g_t x') <1/100$ and $d^{X_0}(g_t
y, g_t y') < 1/100$. Then, 
\begin{itemize}
\item[{\rm (a)}] There exists $\alpha_2$ depending only on the
  Lyapunov spectrum, such that
\begin{equation}
\label{eq:Psharp:bound}
\|P^\sharp(y,y') P^{GM}(y',y) - I \|_Y = O(e^{-\alpha_2  s}).
\end{equation}
\item[{\rm (b)}] There exists $\alpha_6$ depending only on the
  Lyapunov spectrum such that
\begin{equation}
\label{eq:triple:loop}
\|P^+(x',y') \circ P^-(x,x') - P^{GM}(y,y') \circ
P^+(x,y)  \|_Y =
O(e^{-\alpha_6 s}).
\end{equation}
\end{itemize}
\end{lemma}

\bold{Proof.} Note that part (a) follows immediately from 
Lemma~\ref{lemma:P:sharp},  since we are assuming that $d^{X_0}(g_t y, g_t y') \le 1/100$ for all $t$ with $|t| \le s$. 

To prove (b) we abuse notation by identifying $H^1_+$ at all
four points $x$, $y$, $x'$, $y'$ using the Gauss-Manin connection. We
write $\cV_i(x)$ for $\cV_i(H^1_+)(x)$.  
Since
\begin{displaymath}
P^+(x',y') \circ P^-(x,x') \circ P^+(x,y)^{-1} \cV_i(y) = \cV_i(y'),
\end{displaymath}
and by Lemma~\ref{lemma:P:sharp},
\begin{displaymath}
d_Y(\cV_i(y),\cV_i(y')) = O(e^{-\alpha_2 s} ),
\end{displaymath}
it is enough to check that for $v \in \cV_i(y)$,
\begin{equation}
\label{eq:triple:mod:loop}
\|(P^+(x',y') \circ P^-(x,x') \circ P^+(x,y)^{-1}-I) v + \cV_{< i}(y)\|_Y =
O(e^{-\alpha_6 s} \|v\|_Y). 
\end{equation}
But (\ref{eq:triple:mod:loop}) follows from the following:
\begin{itemize}
\item $P^+(x,y)^{-1}$ is the identity map on $\cV_{\le i}(y)/\cV_{< i}(y) =
  \cV_{\le i}(x)/\cV_{< i}(x)$. 
\item $P^-(x,x') \cV_{\le i}(x) = \cV_{\le i}(x')$ and 
by Lemma~\ref{lemma:P:sharp}), $\|P^-(x,x') - I\|_Y =
O(e^{-\alpha_2 s})$.
\item $P^+(x',y')$ is the identity on $\cV_{\le i}(x')/\cV_{< i}(x') =
  \cV_{\le i}(y')/\cV_{< i}(y')$.  
\item $d_Y(\cV_{\le i}(y), \cV_{\le i}(y')) = O(e^{-\alpha_2 s})$. 
\end{itemize}
This completes the proof of (\ref{eq:triple:mod:loop}) and thus
(\ref{eq:triple:loop}).
\qed\medskip

\begin{lemma}
\label{lemma:tworoutes}
$ $
\begin{itemize} 
\item[{\rm (a)}] Suppose $x$, $\tilde{x}$, $y$, $\tilde{y}$ all belong to
  $\pi^{-1}(K^\sharp)$, $d^{X_0}(x,y) < 1/100$, 
$\tilde{y} \in W^+[\tilde{x}]$, $d^X(x,\tilde{x}) \le \xi$
  and $d^X(y,\tilde{y}) \le \xi$. 
Then
\begin{displaymath}
\| \bP^+(\tilde{x}, \tilde{y}) \circ \bP^{Z(x)}(x,\tilde{x}) -
\bP^{Z(y)} (y,\tilde{y}) \circ \bP^{Z(x)}(x,y) \| \le \xi',
\end{displaymath}
where $\xi' \to 0$ as $\xi \to 0$.
\item[{\rm (b)}] Suppose $x,x',y,y' \in \pi^{-1}(K_0)$, $y \in W^+[x]$, $y'
  \in W^+[x']$ and $x' \in W^-[x]$. Suppose further that $d^{X_0}(x,y)
  \le 1/100$, $d^{X_0}(y,y') \le 1/100$, 
and that there exists $s > 0$ such that 
for all $|t| \le s$, $d^{X_0}(g_t x , g_t x') \le 1/100$ and
$d^{X_0}(g_t y, g_t y') \le 1/100$. 
Furthermore, suppose $0 < \alpha_0 < 1$ and 
that $0 < \tau < \alpha_0^{-1} s$ is such
that $d^{X_0}(g_\tau y, g_\tau y') < 1/100$ and $g_\tau y \in K^\sharp$. 
Then,
\begin{displaymath}
\|\bP^+(x',y') \circ \bP^-(x,x') - \bP^{g_{-\tau}Z(g_\tau y)}(y,y') \circ \bP^+(x,y)
\| = O(e^{-\alpha s}),
\end{displaymath}
where $\alpha$ depends only on the Lyapunov spectrum and $\alpha_0$.
\end{itemize}
\end{lemma}

\bold{Proof of (a).} 
Since $y \in W^+[x]$, by
Lemma~\ref{lemma:Psharp:is:Pplus} we have $\bP^{Z(x)}(x,y ) =
\bP^+(x,y)$. Since $\bP^{Z(x)}(x,y)$ depends continuously on $x \in
K^\sharp$ and $y \in K^\sharp$, part (a) 
follows from a compactness agrement. 

\bold{Proof of (b).} We first claim that 
\begin{equation}
\label{eq:P:super:g:minus:t:bound}
\|P^{g_{-\tau} Z(g_\tau y)}(y, y') P^{GM}(y',y)_* - I \|_Y = O(e^{-\alpha' s}), 
\end{equation}
where $\alpha'$ depends only on $\alpha_0$ and the Lyapunov
spectrum. 

By $(K^\sharp 1)$ there exists
$\epsilon_0 > 0$ such that for $x_1 \in \pi^{-1}(K^\sharp)$, $y_1 \in \pi^{-1}(K^\sharp)$ with
$d^{X_0}(x_1,y_1) < \epsilon_0$, $hd_{x_1}^{X_0}(U^+[x_1],U^+[y_1]) < 0.01$. 
By $(K^\sharp 10)$ there exists $t > s/2$ with $g_t y \in \pi^{-1}(K^\sharp)$, $g_t y'
\in \pi^{-1}(K^\sharp)$ and $d^{X_0}(g_t y, g_t y') < 1/100$. Therefore, by
Lemma~\ref{lemma:forni} (c) and Proposition~\ref{prop:AGY:regularity}
we have
\begin{displaymath}
hd_{x'}^{X_0}(U^+[y], U^+[y']) = O(e^{-\alpha_3 s}). 
\end{displaymath}
where $\alpha_3$ depends only on the Lyapunov spectrum. Therefore, we get
\begin{displaymath}
d_Y(P^{GM}(y,y')^{-1}_* \Lie(U^+)(y),\Lie(U^+)(y')) = O(e^{-\alpha_3 s}). 
\end{displaymath}
Then, by (\ref{eq:Psharp:bound}), 
\begin{equation}
\label{eq:P:sharp:pullback:close}
d_Y(P^\sharp(y,y')^{-1}_* \Lie(U^+)(y'), \Lie(U^+)(y)) = O(e^{-\alpha_4
s})
\end{equation}
where $\alpha_4$ depends only on the Lyapunov spectrum.

Since $g_\tau y \in \pi^{-1}(K^\sharp)$, by
$(K^\sharp 1)$ and $(K^\sharp 2)$,
\begin{displaymath}
d_Y(Z(g_\tau y) \cap \cV_i(\Lie(\cG_{++}))(g_\tau y), \Lie(U^+)(g_\tau y) \cap
\cV_i(\Lie(\cG_{++}))(g_\tau y)) \ge c(K^\sharp).
\end{displaymath}
By $(K^\sharp 5)$ (i.e.\ the multiplicative ergodic theorem), 
the restriction of $g_\tau$ to
$\cV_i(\Lie(\cG_{++}))$ is $e^{\lambda_i \tau} h_\tau$, where $\|h_\tau \| =
O(e^{\epsilon'  \tau})$.  Therefore, 
\begin{equation}
\label{eq:transversal:not:so:bad}
d_Y(g_{-\tau} Z(g_\tau y) \cap \cV_i(\Lie(\cG_{++}))(y), \Lie(U^+)(y) \cap
\cV_i(\Lie(\cG_{++}))(y)) \ge c 
e^{-\epsilon' s} 
\end{equation}
We may assume (since $\alpha > 0$ in the choice of $K^\sharp$ is
arbitrary), that $\epsilon' < \alpha_4/2$. Then, 
it follows from (\ref{eq:P:sharp:pullback:close}),
(\ref{eq:transversal:not:so:bad}), (\ref{eq:def:Mxy:one}) and
(\ref{eq:def:Mxy:two})  that 
\begin{equation}
\label{eq:M:uprime:q1prime:bound}
\|M(y; y')\|_Y = O(e^{-\alpha_5 s})
\end{equation}
where $M(\cdot; \cdot)$ is as in (\ref{eq:def:Mxy}), and
$\alpha_5$ depends only on $\alpha_0$ and the Lyapunov spectrum.
Now, (\ref{eq:P:super:g:minus:t:bound}) follows from
(\ref{eq:Psharp:bound}) and (\ref{eq:M:uprime:q1prime:bound}).

Combining (\ref{eq:P:super:g:minus:t:bound}),
and (\ref{eq:triple:loop}) we get
\begin{displaymath}
\|P^+(x',y') \circ P^-(x,x') - P^{g_{-\tau}Z(g_\tau y)}(y,y')
\circ P^+(x,y) \|_Y = O(e^{-\alpha_6 s}).
\end{displaymath}
Now (b) of Lemma~\ref{lemma:tworoutes} follows
immediately, see also ($K^\sharp$10). 
\qed\medskip

\begin{lemma}
\label{lemma:on:correct:sheet}
Suppose $q_1 \in K^\sharp$ and $q_1' \in W^-[q] \cap K^\sharp$,
are such that $d^{X_0}(q_1,q_1') < 1/100$. \mcc{check that this is enough}
Suppose $u \in \cB(q_1,1/100)$, $u' \in \cB(q_1',1/100)$, 
with $u q_1 \in K^\sharp$,
$u' q_1' \in K^\sharp$.  
We write $q_2 = g_\tau u q_1$ for some $\tau >
0$, and let $q_2' = g_\tau u' q_1'$ (see
Figure~\ref{fig:outline2} on page~\pageref{fig:outline2}).
Suppose $d^{X_0}(q_2, q_2') < 1/100$, and also 
there exists $\alpha_0 > 0$ depending only on the Lyapunov spectrum
such that for $|t| < \alpha_0 \tau$, $d^{X_0}(g_t u q_1, g_t u' q_1') < 1/100$. 

In addition, suppose
there exist $\tilde{q}_2 \in X$ and $\tilde{q}_2'\in X$ with
$\sigma_0(\tilde{q}_2')  \in W^+[\sigma_0(\tilde{q}_2)]$
such that $d^X(\tilde{q}_2, q_2) < \xi$ and
$d^X(\tilde{q}_2', q_2') < \xi$. 
Suppose further that $q_2$, $q_2'$,
$\tilde{q}_2$ and $\tilde{q}_2'$ all belong to $K^\sharp$.

Then, (assuming $\epsilon'$
in ($K^\sharp$5) is sufficiently small depending on $\alpha_0$ and the
Lyapunov spectrum), $\tau$ is sufficiently large and $\xi$ is
sufficiently small, (both depending only on $K^\sharp$), we have
\begin{displaymath}
\tilde{q}_2' \in W^+[\tilde{q}_2].
\end{displaymath}
\end{lemma}

\bold{Proof.} In this proof, $\alpha$ is a generic constant depending only on
$\alpha_0$ and the Lyapunov spectrum, with its value changing from
line to line. 

By Lemma~\ref{lemma:Psharp:easybounds} (a), 
\begin{displaymath}
\|P^\sharp(uq_1, u'q_1') \circ P^{GM}(u q_1, u' q_1')^{-1} - I \|_Y =
O(e^{-\alpha \tau}).
\end{displaymath}
By Lemma~\ref{lemma:Psharp:easybounds} (b),
\begin{displaymath}
\|P^{GM}(u q_1, u' q_1') \circ P^+(q_1, u q_1) - P^+(q_1',
u' q_1')\circ P^-(q_1,q_1') \|_Y = O(e^{-\alpha \tau}).
\end{displaymath}
Thus, 
\begin{equation}
\label{eq:tmp:triple:triple}
\|P^\sharp(u q_1, u' q_1') \circ P^+(q_1, uq_1) - P^+(q_1',
u' q_1') \circ P^-(q_1,q_1') \|_Y = O(e^{-\alpha \tau}).
\end{equation}
Write $u' q_1' = (\sigma_0(u' q_1'), \gF')$, $u q_1 = (\sigma_0(u
q_1), \gF)$ where $\gF$ and $\gF'$ are as is in \S\ref{sec:subsec:the:cover:X}. 

By Proposition~\ref{prop:sublyapunov:locally:constant}, (see also 
(\ref{eq:def:Wplus:X}) and (\ref{eq:def:Wminus:X})), 
\begin{displaymath}
\gF' =  P^+(q_1', u'q_1') \circ P^-(q_1,q_1')  \circ P^+(uq_1, q_1) \gF.
\end{displaymath}
Therefore, by (\ref{eq:tmp:triple:triple}),
\begin{displaymath}
d_Y(\gF', P^\sharp(u q_1, u' q_1') \gF) = O(e^{-\alpha \tau}).
\end{displaymath}
where the distance $d_Y( \cdot, \cdot)$ between flags is as in
\S\ref{sec:subsec:the:cover:X}. 

We now claim that
\begin{equation}
\label{eq:tmp:gt:gF}
d_Y(g_\tau \gF', g_\tau P^\sharp(u q_1, u' q_1') \gF) = O(e^{-\alpha \tau}).
\end{equation}
Indeed to prove (\ref{eq:tmp:gt:gF}) it is enough to show that for
each $i$, 
\begin{equation}
\label{eq:tmp:gt:gFi}
d_Y(g_\tau \gF_i', g_\tau P^\sharp(u q_1, u' q_1') \gF_i) = O(e^{-\alpha \tau}).
\end{equation}
But $\gF_i' \subset \cV_i(H_{big})(u'q_1')$, 
$\gF_i \subset \cV_i(H_{big})(u q_1)$, and
\begin{displaymath}
P^\sharp(u q_1, u' q_1') \cV_i(H_{big})(u
q_1) = \cV_i(H_{big})(u'q_1').
\end{displaymath}
Thus, we have
\begin{displaymath}
\gF_i' \subset \cV_i(H_{big})(u'q_1'), \qquad 
P^\sharp(uq_1, u'q_1') \gF_i \subset \cV_i(H_{big})(u'q_1')
\end{displaymath}
The geodesic flow $g_\tau$ restricted to $\cV_i(H_{big})(u'q_1')$ is of
the form $e^{\lambda_i \tau} h_\tau$, where $\|h_\tau\|_Y = O(e^{\epsilon'
  \tau})$. Thus, (\ref{eq:tmp:gt:gFi}) and hence (\ref{eq:tmp:gt:gF})
follows. The equivariance property of $P^\sharp$ then implies that
\begin{equation}
\label{eq:tmp2:gt:gF}
d_Y(g_\tau \gF', P^\sharp(q_2, q_2') g_\tau \gF) = O(e^{-\alpha \tau}).
\end{equation}

We have since the $\cV_i$ are continuous on $K^\sharp$  and
Lemma~\ref{lemma:properties:lyapunov:flag}, 
\begin{displaymath}
\|P^{GM}(q_2', \tilde{q}_2') \circ P^\sharp(q_2, q_2') 
- P^+(\tilde{q}_2, \tilde{q}_2') \circ P^{GM}(q_2,\tilde{q}_2) \|_Y \to 0,
\end{displaymath}
as $\xi \to 0$. Combining this with (\ref{eq:tmp2:gt:gF}), we get
\begin{equation}
\label{eq:dPGM:gt:gFprime}
d_Y(P^{GM}(q_2', \tilde{q}_2') g_\tau \gF', P^+(\tilde{q}_2, \tilde{q}_2')
\circ P^{GM}(q_2,\tilde{q}_2) g_\tau \gF) \to 0, 
\end{equation}
as $\xi \to 0$ and $\tau \to \infty$.

Note that $q_2 = (\sigma_0(q_2), g_\tau \gF)$, $q_2' = (\sigma_0(q_2'),
g_\tau \gF')$. 
Write $\tilde{q}_2 = (\sigma_0(\tilde{q}_2), \tilde{\gF})$, 
$\tilde{q}_2' = (\sigma_0(\tilde{q}_2'), \tilde{\gF}')$. Then, since
$d^X(q_2, \tilde{q}_2) \to 0$, in view of (\ref{eq:def:dX}), 
\begin{displaymath}
d_Y(P^{GM}(q_2,\tilde{q}_2) g_\tau \gF, \tilde{\gF}) \le \xi'
\end{displaymath}
\begin{displaymath}
d_Y(P^{GM}(q_2',\tilde{q}'_2) g_\tau \gF', \tilde{\gF}') \le \xi'
\end{displaymath}
where $\xi' \to 0$ as $\xi \to 0$. 
Hence, by (\ref{eq:dPGM:gt:gFprime}), 
\begin{displaymath}
d_Y(\tilde{\gF}', P^+(\tilde{q}_2, \tilde{q}_2') \tilde{\gF}) \to 0
\qquad\text{as $\xi \to 0$ and $\tau \to \infty$.}
\end{displaymath}
This implies that $\tilde{q}_2' \in W^+[\tilde{q}_2]$ by
(\ref{eq:def:Wplus:X}). 
\qed\medskip

\bold{Proof of Lemma~\ref{lemma:tau:ij:nearby:close}.} 
Note that, by the construction of $\bP^{Z(\cdot)}(\cdot, \cdot)$,
for all $i \in \Lambda'$, 
\begin{equation}
\label{eq:PZ:preserves:lyapunov}
\bP^{Z(g_t y)}(g_t y, g_t y') \cV_i(\bH)(g_t y) = \cV_i(\bH)(g_t y').
\end{equation}
However, even though for all $ij \in
\Lambda''$, $\bE_{ij}(x) \subset
\cV_i(\bH)(x)$, we may have 
\begin{displaymath}
\bP^{Z(g_t y)}(g_t y, g_t y')
\bE_{ij}(g_t y) \ne \bE_{ij}(g_t y').
\end{displaymath}


Suppose $\bfv \in \bE_{ij}(y)$, and that $\bfv$ is orthogonal to
$\bE_{i,j-1}(y) \subset
\bE_{ij}(y)$.  Let
\begin{displaymath}
\bfv' = \bP^+(x',y') \circ \bP^-(x,x') \circ \bP^+(y,x) \bfv.
\end{displaymath}
Then, by
Proposition~\ref{prop:sublyapunov:locally:constant}
(a), 
$\bfv' \in \bE_{ij}(y')$. By
$(K^\sharp 1)$, and the fact that
\begin{displaymath}
\bP^+(x',y') \circ \bP^-(x,x') \circ
\bP^+(y,x) \bE_{i,j-1}(y) =
\bE_{i,j-1}(y'), 
\end{displaymath}
we have
\begin{equation}
\label{eq:tmp:tau:ij:start}
C_1^{-1} \| \bfv \| \le \|\bfv' +
\bE_{i,j-1}(y') \|  
\le C_1 \| \bfv \|, 
\end{equation}
where $C_1$ depends only on
$K_0$. By Lemma~\ref{lemma:tworoutes} (b), 
\begin{displaymath}
\| \bP^{g_{-t}Z(g_t y)}(y,y') \bfv -   \bfv' \| = O(e^{-\alpha_1' t}\|\bfv\|),
\end{displaymath}
where $\alpha_1'$ depends only on
$\alpha_0$ and the Lyapunov
spectrum. By
(\ref{eq:PZ:preserves:lyapunov}), 
\begin{displaymath}
\bP^{g_{-t}Z(g_t y)}(y,y') \bfv \in \cV_i(\bH)(y').
\end{displaymath}
Then by the multiplicative ergodic theorem (see also $(K^\sharp 5)$), 
\begin{equation}
\label{eq:tmp:tau:ij:Pv:v:prime}
\|\bP^{Z(g_t y)}(g_t y,g_t y') (g_t \bfv)-    
g_t \bfv'\| =
O(e^{-(\alpha_1'-\epsilon') t}\|g_t
\bfv \|).
\end{equation}
Since $\bfv$ is arbitrary, this
implies that for all $ij \in
\Lambda''$, 
\begin{equation}
\label{eq:tau:ij:subspaces:close}
d\left(\bP^{Z(g_t y)}(g_t y, g_t y') \bE_{ij}(g_t y),
\bE_{ij}(g_t y')\right) = O( e^{-\alpha_1 t}),
\end{equation}
where $\alpha_1$ depends only on $\alpha_0$ and the Lyapunov
spectrum. 

By $(K^\sharp 1)$ and $(K^\sharp 2)$, 
\begin{displaymath}
\| \bP^{Z(g_t y)}(g_t y, g_t y') \| \le C_1'
\end{displaymath}
where $C_1'$ depends only on
$K_0$. Therefore, by
(\ref{eq:tmp:tau:ij:Pv:v:prime}) and (\ref{eq:tau:ij:subspaces:close}), 
\begin{equation}
\label{eq:tmp:tau:ij:almost:done}
C_2^{-1} \| g_t \bfv + \bE_{i,j-1}(g_t y) \| \le
\| g_t \bfv' + \bE_{i,j-1}( g_t y') \| \le  C_2 \| g_t \bfv' +
\bE_{i,j-1}( g_t y) \|,
\end{equation}
where $C_2$ depends only on $K_0$,
$\alpha_0$ and the Lyapunov
spectrum. 

Note that
\begin{displaymath}
\hat{\tau}_{ij}(y,t) = \frac{\| g_t
  \bfv + \bE_{i,j-1}(g_t y) \| } { \| \bfv
  \|}, \qquad 
\hat{\tau}_{ij}(y',t) = \frac{\| g_t
  \bfv' + \bE_{i,j-1}(g_t y') \| } { \| \bfv'
  + \bE_{i,j-1}(y') \|}. 
\end{displaymath}
Now (\ref{eq:tau:ij:nearby:close})
follows from
(\ref{eq:tmp:tau:ij:start}) and
(\ref{eq:tmp:tau:ij:almost:done}). 
\qed\medskip

\begin{proposition}
\label{prop:tmp:nearby:linear:maps}
Suppose $\alpha$, $\epsilon$, $s$, $\ell$, $t$, $t'$, $q$, $q'$, $\tau$,
$q_1$, $q_1'$, $q_3$, $q_3'$, $u$, $u'$, $q_2$, $q_2'$, $\tilde{q}_2$,
$\tilde{q}_2'$, $C$, $C_1$, $\xi$ are as in 
Proposition~\ref{prop:nearby:linear:maps}. Suppose also $\tilde{q}_2'
\in W^+[\tilde{q}_2]$. Then, (assuming $\epsilon'$
in ($K^\sharp$5) is sufficiently small depending on $\alpha_0$ and the
Lyapunov spectrum),
\begin{itemize}
\item[{\rm (a)}] There exists $\xi' > 0$ (depending on
    $\xi$,  $K_0$ and $C$ and $t$)
with $\xi' \to 0$ as $\xi \to 0$ and $t \to \infty$
such that for $\bfv \in
\bE_{[ij],bdd}(q_3)$, 
\begin{multline}
\label{eq:nearby:linear:maps:a}
\| \bP^{Z(q_2')}(q_2',\tilde{q}_2') \circ R(q_3',q_2') \circ
\bP^-(q_3,q_3') \bfv - \\
\bP^+(\tilde{q}_2,\tilde{q}_2') \circ \bP^{Z(q_2)}(q_2,\tilde{q}_2)
\circ R(q_3,q_2) \bfv \| \le \xi' \|\bfv \|.  
\end{multline}
\item[{\rm (b)}] There exists $\xi'' > 0$ (depending on $\xi$, 
  $K_0$, $C$ and $t$) 
with $\xi'' \to 0$ as $\xi \to 0$ and $t \to \infty$ such that 
\begin{displaymath}
d_*(\bP^+(\tilde{q}_2,\tilde{q}_2') 
\bff_{ij}(\tilde{q}_2), \bff_{ij}(\tilde{q}_2')) \le \xi''.
\end{displaymath}
Here $d_*(\cdot, \cdot)$ is any metric which induces the weak-*
convergence topology on the
  domain of common definition of the measures, up to
  normalization. \mcc{say this better} 
\end{itemize}
\end{proposition}

\bold{Proof of (a).}  
Following the outline given after the statement
of Proposition~\ref{prop:nearby:linear:maps}, 
the proof will consist of verifying conditions
(i), (ii) and (iii) of Lemma~\ref{lemma:lin:alg:B:Bprime}, with
$E = \bE_{ij,bdd}(q_3)$, $L = \bH(\tilde{q}_2')$, 
$F = \bE_{ij,bdd}(\tilde{q}_2')$, $V = \bV_{< i}(\tilde{q}_2')$, and
$\bB$ and $\bB'$ as the linear maps on the first and second line of
(\ref{eq:nearby:linear:maps:a}). (We note that $\bB$ and $\bB'$ are
bounded by Proposition~\ref{prop:ej:bdd:transport:bounded}.)
We start with (i). 

Note that by
  (\ref{eq:bilip:hat:tau:ij}), we have 
\begin{equation}
\label{eq:kappa:inv:tau:t:kappa:tau}
\kappa^{-1} \tau \le t \le \kappa \tau,
\end{equation}
where $\kappa$ depends only on the Lyapunov spectrum. Also, by assumption we have 
\begin{displaymath}
\ell > \alpha_0 \tau,
\end{displaymath}
where $\alpha_0$ depends only on the Lyapunov
spectrum. 


Suppose $\bfw \in \bE_{ij,bdd}(q_1)$. We now apply 
Lemma~\ref{lemma:tworoutes}(b), with $x = q_1$, $x' = q_1'$, $y = u
q_1$, $y' = u' q_1'$ and $\tau = \tau$ to get
\begin{displaymath}
\|  \bP^+(u'q_1', uq_1) \circ   \bP^-(q_1,q_1') \bfw - \bP^{g_{-\tau}Z(q_2)}(u q_1, 
u' q_1') \circ \bP^+(q_1,u q_1) \bfw \| =   O(e^{-\alpha
  \tau} \|\bfw\|),  
\end{displaymath}
By Proposition~\ref{prop:sublyapunov:locally:constant}(a),
$\bP^-(q_1,q_1') \bfw \in \bE_{ij,bdd}(q_1') \subset \bE(q_1')$. Therefore, by
Lemma~\ref{lemma:u:star:agrees:with:Pplus}, this can be rewritten as
\begin{displaymath}
\|  (u')_* \circ   \bP^-(q_1,q_1') \bfw - \bP^{g_{-\tau}Z(q_2)}(u q_1, 
u' q_1') \circ (u)_* \bfw \| =   O(e^{-\alpha
  \tau} \|\bfw\|),  
\end{displaymath}
Hence,
\begin{equation}
\label{eq:tmp:bP:super}
(u')_* \circ 
  \bP^-(q_1,q_1') \bfw =  
\bP^{g_{-\tau}Z(q_2)}(u q_1, u' q_1') \circ (u)_* \bfw  +
  \bfw'
\end{equation}
where $\bfw' \in \bH(u' q_1')$ satisfies 
\begin{equation}
\label{eq:tmp:estimate:bfw:prime}
\|\bfw'\| = O(e^{-\alpha \tau} \|\bfw\|) =
O_{\epsilon'}(e^{-(\lambda_i + 
  \alpha - \epsilon') \tau}\|\bfv\|), 
\end{equation}
where we wrote $\bfw = g_{-t'}^{ij} \bfv$ for some $\bfv \in
\bE_{ij}(q_3)$, 
and we have used ($K^\sharp$5), (\ref{eq:kappa:inv:tau:t:kappa:tau})
and the assumption $|t - t' | < C$  for the last estimate. 
We now apply $g_\tau=g^{ij}_{t}$
to both sides of (\ref{eq:tmp:bP:super}) and take
the quotient mod $\bV_{< i}(q_2')$. We get
\begin{multline}
\label{eq:tmp:bP:super:mod:Viminus1}
g_\tau \circ (u')_* \circ
\bP^-(q_1,q_1') 
\bfw + \bV_{< i}(q_2') = \\ = 
\bP^{Z(q_2)}(q_2, q_2') \circ [ g_\tau \circ (u)_* \bfw +
g_\tau \bfw'] + \bV_{< i}(q_2'). 
\end{multline}
We may write 
\begin{displaymath}
\bfw' = \sum_{k} \bfw_k, \qquad \bfw_k \in \cV_k(\bH)(u'q_1').
\end{displaymath}
Then, 
\begin{displaymath}
g_\tau \bfw' + \bV_{< i}(q_2') = \sum_{k} g_\tau \bfw'_k +
\bV_{< i}(q_2')  = \sum_{k \ge i} g_\tau \bfw'_k +
\bV_{< i}(q_2'), 
\end{displaymath}
since for $k <i$, $g_\tau \bfw'_k \in \bV_{< i}(q_2')$. By
($K^\sharp$5), for $k \ge i$, 
\begin{displaymath}
\|g_\tau \bfw_k'\| = O(e^{(\lambda_k + \epsilon')\tau} \|\bfw_k'\|) =
O(e^{(\lambda_i + \epsilon')\tau} \|\bfw_k'\|) = O(e^{-\alpha_5 \tau}\|\bfv\|),  
\end{displaymath}
using (\ref{eq:tmp:estimate:bfw:prime}) 
(and choosing $\epsilon'$ sufficiently small depending on $\alpha_0$
and the Lyapunov spectrum). 
Therefore, substituting into (\ref{eq:tmp:bP:super:mod:Viminus1}), we
get, for  $\bfv \in \bE_{ij,bdd}(q_3)$,
\begin{multline*}
R(q_3',q_2') \circ \bP^-(q_3,q_3')
\bfv + \bV_{< i}(q_2') = \\ = 
\bP^{Z(q_2)}(q_2,q_2') 
\circ R(q_3,q_2) \bfv + O(e^{-\alpha_5 \tau}\|\bfv\|) + 
\bV_{< i}(q_2'). 
\end{multline*}
We now apply $\bP^{Z(q_2')}(q_2',\tilde{q}_2')$ to both sides to get
(using (\ref{eq:bPZ:preserves:eigenspaces}))
\begin{multline}
\label{eq:estimate:modViminus1:prefinal}
\bP^{Z(q_2')}(q_2',\tilde{q}_2') \circ R(q_3',q_2') \circ \bP^-(q_3,q_3')
\bfv + \bV_{< i}(\tilde{q}_2') = \\ = 
\bP^{Z(q_2')}(q_2',\tilde{q}_2') \circ
\bP^{Z(q_2)}(q_2,q_2') 
\circ R(q_3,q_2) \bfv + O(e^{-\alpha_5 \tau}\|\bfv\|) + 
\bV_{< i}(\tilde{q}_2'). 
\end{multline}
Since $q_2$, $\tilde{q}_2$, $q_2'$, $\tilde{q}_2'$ all
belong to $K^\sharp$, we have by Lemma~\ref{lemma:tworoutes}(a),
\begin{displaymath}
\|\bP^{Z(q_2')}(q_2',\tilde{q}_2') \circ
\bP^{Z(q_2)}(q_2,q_2') -
\bP^+(\tilde{q}_2,\tilde{q}_2') \circ \bP^{Z(q_2)}(q_2,\tilde{q}_2) \|
\le \xi_3,
\end{displaymath}
where $\xi_3 \to 0$ as $\xi \to 0$. 
Therefore, substituting into
(\ref{eq:estimate:modViminus1:prefinal}), we get
\begin{multline*}
\bP^{Z(q_2')}(q_2',\tilde{q}_2') \circ R(q_3',q_2') \circ \bP^-(q_3,q_3')
\bfv + \bV_{< i}(\tilde{q}_2') = \\ = 
\bP^+(\tilde{q}_2,\tilde{q}_2') \circ \bP^{Z(q_2)}(q_2,\tilde{q}_2) 
\circ R(q_3,q_2) \bfv + O(e^{-\alpha_5 \tau}\|\bfv\|) +
O(\xi_3 \| \bfv \|) + 
\bV_{< i}(\tilde{q}_2'). 
\end{multline*}
This completes the verification of (i) of
Lemma~\ref{lemma:lin:alg:B:Bprime}.

We now verify (ii) of Lemma~\ref{lemma:lin:alg:B:Bprime}. 
For $\bfv \in \bE_{ij,bdd}(q_3)$, 
we have $R(q_3,q_2) \bfv \in \bE_{ij,bdd}(q_2)$, and then
\begin{displaymath}
\bP^+(\tilde{q}_2, \tilde{q}_2') \circ \bP^{Z(q_2)}(q_2,\tilde{q}_2) 
\circ R(q_3,q_2) \bfv \in \bP^+(\tilde{q}_2, \tilde{q}_2') \circ
\bP^{Z(q_2)}(q_2,\tilde{q}_2) \bE_{ij,bdd}(q_2). 
\end{displaymath}
By ($K^\sharp$2) and ($K^\sharp$3), since $d^X(q_2,\tilde{q}_2) < \xi$, 
\begin{displaymath}
d_Y(\bP^{Z(q_2)}(q_2,\tilde{q}_2) \bE_{ij,bdd}(q_2),
\bE_{ij,bdd}(\tilde{q}_2)) < \xi_0,
\end{displaymath}
where $\xi_0 \to 0$ as $\xi \to 0$. 
Then, using
(\ref{eq:bPplus:preserves:bEijbdd}),  
\begin{displaymath}
d_Y(\bP^+(\tilde{q}_2, \tilde{q}_2') \circ
\bP^{Z(q_2)}(q_2,\tilde{q}_2) \bE_{ij,bdd}(q_2),
\bE_{ij,bdd}(\tilde{q}_2')) < \xi_1.  
\end{displaymath}
where $\xi_1 \to 0$ as $\xi \to 0$. This completes the verification of
condition (ii) of Lemma~\ref{lemma:lin:alg:B:Bprime}.

Also, by (\ref{eq:bPplus:preserves:bEijbdd}) 
(applied to $\bP^-$), we have $\bP^-(q_3,q_3') \bfv \in
\bE_{ij,bdd}(q_3')$. Then, $R(q_3', q_2') \circ \bP^-(q_3,q_3') \bfv
\in \bE_{ij,bdd}(q_2')$, and
\begin{displaymath}
\bP^{Z(q_2')}(q_2',\tilde{q}_2') \circ R(q_3', q_2') \circ \bP^-(q_3,q_3') \bfv
\in \bP^{Z(q_2')}(q_2',\tilde{q}_2') \bE_{ij,bdd}(q_2').
\end{displaymath}
By ($K^\sharp$2) and ($K^\sharp$3), 
\begin{displaymath}
d_Y(\bP^{Z(q_2')}(q_2',\tilde{q}_2')
\bE_{ij,bdd}(q_2'), \bE_{ij,bdd}(\tilde{q}_2')) < \xi_2, 
\end{displaymath}
where $\xi_2 \to 0$ as $\xi \to 0$. This completes the
  verification of condition (iii) of Lemma~\ref{lemma:lin:alg:B:Bprime}.
 
Now (\ref{eq:nearby:linear:maps:a}) for arbitrary $\bfv \in
\bE_{ij,bdd}(q_3)$ follows from Lemma~\ref{lemma:lin:alg:B:Bprime}. 
The general case of
(\ref{eq:nearby:linear:maps:a}) (i.e.\ for an arbitrary $\bfv\in
\bE_{[kr],bdd}(q_3)$) follows since $\bE_{[kr],bdd}(q_3) = 
\bigoplus_{ij \in [kr]} \bE_{ij,bdd}(q_3)$ and all the maps on the
left-hand-side of (\ref{eq:nearby:linear:maps:a}) are linear.

\bold{Proof of (b).} By ($K^\sharp$4), 
\begin{displaymath}
d_*(\bP^-(q_3,q_3')_* \bff_{ij}(q_3), \bff_{ij}(q_3')) \le \xi_1,
\end{displaymath}
where $\xi_1 \to 0$ as and $t \to \infty$.
In view of 
in view of condition ($K^\sharp$6), the assumption $|t-t'| < C$ 
 and 
Proposition~\ref{prop:ej:bdd:transport:bounded}, that $R(q_3,q_2)$ is a linear map with norm bounded
depending only on $K^\sharp$ and $C$. 
It then follows from (a) that
$R(q_3',q_2')$ is also a linear map whose norm is bounded depending
only on $K^\sharp$ and $C$. 
Furthermore, by ($K^\sharp$9) and Lemma~\ref{lemma:forni}
there exists a constant $C_2(\delta)$ such that if 
\begin{equation}
\label{eq:tmp:t:minus:tprime}
C > t - t' > C_2(\delta),
\end{equation}
then if we write $q_2 = g^{ij}_{t} u g^{ij}_{-t'} q_3$, then 
$g^{ij}_{t} u g^{ij}_{-t'} \gB_0[q_3] \cap \cC_{ij}[q_3] \supset
\gB_0[q_2] \cap \cC_{ij}[q_2]$. Then, by
Lemma~\ref{lemma:transformation:rule:bff}, 
\begin{displaymath}
\bff_{ij}(q_2) \propto R(q_3,q_2)_* \bff_{ij}(q_3) \quad\text{ and }\quad
\bff_{ij}(q_2') \propto R(q_3',q_2')_* \bff_{ij}(q_3').
\end{displaymath}
In view of ($K^\sharp$11), we can assume that
(\ref{eq:tmp:t:minus:tprime}) holds: 
 otherwise we can replace $q_3$ and
$q_3'$ by $g^{ij}_{-s}q_3 \in K^\sharp$ and $g^{ij}_{-s}q_3' \in
K^\sharp$ where $C_2(\delta) <
s < 2C_2(\delta)$. (Without loss of generality we may assume that $C >
2 C_2(\delta)$.) 
Hence, we have
\begin{equation}
\label{eq:tmp:propto}
d_*( (R(q_3',q_2') \circ \bP^-(q_3,q_3'))_* \bff_{ij}(q_3), 
\bff_{ij}(q_2')) \le \xi_2,
\end{equation}
where $\xi_2 \to 0$ as $t \to \infty$. 
Thus, by ($K^\sharp$1), ($K^\sharp$2), ($K^\sharp$3), 
\begin{displaymath}
d_*(\bP^{Z(q_2')}(q_2',\tilde{q}_2') \bff_{ij}(q_2'), \bff_{ij}(\tilde{q}_2')) \le \xi_3,
\end{displaymath}
where $\xi_3 \to 0$ as $\xi \to 0$ and $t \to \infty$. Hence, 
\begin{equation}
\label{eq:bound:d:lhs}
d_*((\bP^{Z(q_2')}(q_2',\tilde{q}_2') \circ R(q_3',q_2') \circ
\bP^-(q_3,q_3'))_* \bff_{ij}(q_3), \bff_{ij}(\tilde{q}_2')) \le \xi_4,
\end{equation}
where $\xi_4 \to 0$ as $\xi \to 0$ and $t \to
  \infty$.
Also, in view of (\ref{eq:tmp:propto}), 
and since $\bP^+(\tilde{q}_2,\tilde{q}_2')$ is a linear
map whose norm is bounded depending only on $K^\sharp$, 
\mcc{add more detail to the following equation. Worry about domains of
definition}
\begin{equation}
\label{eq:bound:d:rhs}
d_*(\bP^+(\tilde{q}_2,\tilde{q}_2') \circ \bP^{Z(q_2)}(q_2,\tilde{q}_2) 
\circ R(q_3,q_2))_* \bff_{ij}(q_3), 
\bP^+(\tilde{q}_2,\tilde{q}_2')_* \bff_{ij}(\tilde{q}_2)) \le \xi_5, 
\end{equation}
where $\xi_5 \to 0$ as $\xi \to 0$ and $t \to
  \infty$.  
Now part (b) follows from (\ref{eq:bound:d:lhs}),
(\ref{eq:bound:d:rhs}), and 
(\ref{eq:nearby:linear:maps:a}). 

\bold{Proof of Proposition~\ref{prop:nearby:linear:maps}.} Note that
(\ref{eq:goal:same:sheet}) follows from
Lemma~\ref{lemma:on:correct:sheet}. We assume this from now on. 

Without loss of generality, and to simplify the notation, we may assume that 
$Z(\tilde{q}_2') =
P^+(\tilde{q}_2,\tilde{q}_2')
Z(\tilde{q}_2)$. (Otherwise, we can further compose
  with a reparametrization map at $\tilde{q}_2'$ which will not change
  the result).\mcc{explain.}
We have 
\begin{displaymath}
\bff_{ij}(\tilde{q}_2) = (\bfj \circ \phi_{\tilde{q}_2})_* f_{ij}(\tilde{q}_2)
\end{displaymath}
and 
\begin{displaymath}
\bff_{ij}(\tilde{q}_2') = (\bfj \circ \phi_{\tilde{q}_2'})_* f_{ij}(\tilde{q}_2')
\end{displaymath}
As in \S\ref{sec:divergence:subspaces}, 
let $P^+_*: \cH_{++}(\tilde{q}_2) \cross W^+(\tilde{q}_2) \to
\cH_{++}(\tilde{q}_2') \cross 
W^+(\tilde{q}_2')$ be given by
\begin{equation}
\label{eq:def:PPlus:star}
P^+_*(M,v) = (P^+(\tilde{q}_2,\tilde{q}_2')^{-1} \circ M
\circ P^+(\tilde{q}_2,\tilde{q}_2'), P^+(\tilde{q}_2,\tilde{q}_2') v).   
\end{equation}
Then, 
\begin{equation}
\label{eq:bPplus:circ:j}
\bP^+(\tilde{q}_2,\tilde{q}_2') \circ \bfj(M,v) = \bfj(P^+_*(M,v)) 
\end{equation}
We write $A \approx_{\xi,t} B$ if $d(A,B) \to 0$ as $\xi \to 0$ and $t \to \infty$. Then,
we have, by Proposition~\ref{prop:tmp:nearby:linear:maps},
\begin{displaymath}
(\bfj \circ \phi_{\tilde{q}_2'})_* f_{ij}(\tilde{q}_2') =
\bff_{ij}(\tilde{q}_2') \approx_{\xi,t} \bP^+(\tilde{q}_2,\tilde{q}_2')_*
\bff_{ij}(\tilde{q}_2) = (\bP^+(\tilde{q}_2,\tilde{q}_2') \circ \bfj
\circ \phi_{\tilde{q}_2})_* f_{ij}(\tilde{q}_2) 
\end{displaymath}
By (\ref{eq:bPplus:circ:j}),
\begin{displaymath}
(\bfj \circ \phi_{\tilde{q}_2'})_* f_{ij}(\tilde{q}_2') \approx_{\xi,t} (\bfj \circ P^+_*
\circ\phi_{\tilde{q}_2})_* f_{ij}(\tilde{q}_2).  
\end{displaymath}
Therefore, 
\begin{displaymath}
(\phi_{\tilde{q}_2'})_* f_{ij}(\tilde{q}_2') \approx_{\xi,t} (P^+_* \circ
\phi_{\tilde{q}_2})_* f_{ij}(\tilde{q}_2). 
\end{displaymath}
Let $\pi_2: \cH_{++}(x) \cross W^+(x) \to W^+(x)$ be projection onto the second
factor. Then, applying $\pi_2$ to both sides, we get 
\begin{equation}
\label{eq:tmp:pi2:relation}
(\pi_2 \circ \phi_{\tilde{q}_2'})_* f_{ij}(\tilde{q}_2') \approx_{\xi,t} ( \pi_2
\circ P^+_* \circ \phi_{\tilde{q}_2})_* f_{ij}(\tilde{q}_2). 
\end{equation}
For $z \in Z(\tilde{q}_2)$, $\pi_2(\phi_{\tilde{q}_2}(z)) = z$, and
thus in view of (\ref{eq:def:PPlus:star}),
\begin{equation}
\label{eq:tmp:pi2:rhs}
(\pi_2 \circ P^+_* \circ \phi_{\tilde{q}_2})(z) \approx_{\xi,t}
P^+(\tilde{q}_2,\tilde{q}_2') z. 
\end{equation}
By assumption, we have 
$Z(\tilde{q}_2') = P^+(\tilde{q}_2,\tilde{q}_2')Z(\tilde{q}_2)$.
Then, similarly, for $z \in Z(\tilde{q}_2') =
P^+(\tilde{q}_2,\tilde{q}_2')Z(\tilde{q}_2)$, 
\begin{equation}
\label{eq:tmp:pi2:lhs}
(\pi_2 \circ \phi_{\tilde{q}_2'})(z) = z.
\end{equation}
Since $f_{ij}(x)$ is Haar
along $U^+$, we can recover $f_{ij}(\tilde{q}_2)$ from its restrictions to
$Z(\tilde{q}_2)$ and $f_{ij}(\tilde{q}_2')$ from its restriction to $Z(\tilde{q}_2')$. 
It now follows from (\ref{eq:tmp:pi2:relation}), (\ref{eq:tmp:pi2:rhs}) and
(\ref{eq:tmp:pi2:lhs}) that 
\begin{displaymath}
f_{ij}(\tilde{q}_2') \approx_{\xi,t} P^+(\tilde{q}_2,\tilde{q}_2')_*
f_{ij}(\tilde{q}_2). 
\end{displaymath}
\qed\medskip


\section{The inductive step}
\label{sec:inductive:step}

\begin{proposition}
\label{prop:inductive:step}
Suppose $\nu$ is a $P$-invariant measure on $X_0$. 
Suppose $U^+(x)$ is a family of subgroups of $\cG_{++}(x)$ compatible
with $\nu$ in the sense of Definition~\ref{def:compatible:family}. 
Let $L^-[x]$ and $L^+[x]$ be as in \S\ref{sec:subsec:stopping}, and
suppose the equivalent conditions of
Lemma~\ref{lemma:Lplus:Splus:equiv} do not hold. 
Then, there
exists a family of  subgroups 
$U_{new}^+(x)$ of $\cG_{++}(x)$ compatible with $\nu$ in the sense of
Definition~\ref{def:compatible:family} such that for almost all
$x$, $U_{new}^+(x)$ strictly contains $U^+(x)$.
\end{proposition}

The rest of \S\ref{sec:inductive:step} will consist of the proof of
Proposition~\ref{prop:inductive:step}. We assume that $L^-(x)$, 
$L^+(x)$ and
$U^+(x)$ are as in Proposition~\ref{prop:inductive:step}, and the
equivalent conditions of Lemma~\ref{lemma:Lplus:Splus:equiv} do not
hold. The argument has been outlined in
\S\ref{sec:outline:step1}, and we have kept the same notation (in
particular, see Figure~\ref{fig:outline}). 

Let $f_{ij}(x)$ be the measures on $W^+(x)$ introduced in
\S\ref{sec:equivalence}. We think of $f_{ij}$ as a function
from $X$ to a space of measures (which is
metrizable). Let $P^+(x,y)$ be the map introduced
in \S\ref{sec:subsec:connection}. 
Proposition~\ref{prop:inductive:step} will be derived from the
following:
\begin{proposition}
\label{prop:half:inductive:step}
Suppose $U^+$, $L^+$, $L^-$ 
are as in Proposition~\ref{prop:inductive:step},
and the equivalent conditions of Lemma~\ref{lemma:Lplus:Splus:equiv}
do not hold. 
Then there exists $0 < \delta_0 < 0.1$, a subset $K_* \subset X$ with
$\nu(K_*) > 1-\delta_0$ such that
all the functions $f_{ij}$, $ij \in \tilde{\Lambda}$ 
are uniformly continuous on $K_*$, and $C > 1$ (depending on $K_*$) 
such that for  every $0 < \epsilon < C^{-1}/100$ there exists a 
subset $E \subset K_*$ with 
$\nu(E) > \delta_0$, such that for every $x \in \pi^{-1}(E)$ there exists $ij
\in \tilde{\Lambda}$ and $y \in \cC_{ij}[x] \cap \pi^{-1}(K_*)$ with 
\begin{equation}
\label{eq:dxy:within:C:epsilon}
C^{-1} \epsilon \le hd_x^{X_0}(U^+[x],U^+[y]) \le C \epsilon
\end{equation}
and (on the domain where both are defined)
\begin{equation}
\label{eq:half:inductive:step:two}
f_{ij}(y) \propto P^+(x,y)_* f_{ij}(x). 
\end{equation}
\end{proposition}

We now begin the proof of Proposition~\ref{prop:half:inductive:step}. 

\bold{Choice of parameters \#1.}
Fix $\theta > 0$ as in Proposition~\ref{prop:some:fraction:bounded}
and Proposition~\ref{prop:ej:bdd:transport:bounded}.
We then choose
$\delta > 0$ sufficiently small; the exact value of $\delta$ will we
chosen at the end of this section. All subsequent constants
will depend on $\delta$. (In particular, $\delta \ll \theta$; we will
make this more precise below). Let $\epsilon > 0$ be arbitrary and
$\eta > 0$ be arbitrary;  however we will always assume that
$\epsilon$ and $\eta$ are sufficiently small depending on $\delta$.

We will show that Proposition~\ref{prop:half:inductive:step} holds
with $\delta_0 = \delta/10$. Let $K_* \subset X$ be any subset with
$\nu(K_*) > 1 - \delta_0$ on which all the functions $f_{ij}$
are uniformly continuous. 
It is enough to show that there exists $C
= C(\delta)$ such that for any $\epsilon > 0$ and 
for an arbitrary compact set $K_{00} \subset X$ with $\nu(K_{00}) \ge (1-2
\delta_0)$, there exists $x \in K_{00} \cap K_*$, $ij \in
\tilde{\Lambda}$  and $y \in \cC_{ij}[x] \cap K_*$ 
satisfying (\ref{eq:dxy:within:C:epsilon}) and
(\ref{eq:half:inductive:step:two}). Thus, let $K_{00} \subset X$ be an
arbitrary compact set with $\nu(K_{00}) > 1- 2\delta_0$.

We can choose a compact set $K_0 \subset K_{00} \cap K_*$ with
$\nu(K_0) > 1-5 \delta_0 = 1-\delta/2$ so that 
Proposition~\ref{prop:nearby:linear:maps} holds. 
In addition, 
there exists $\epsilon_0'(\delta) > 0$ such that for all $x \in K_0$, 
\begin{equation}
\label{eq:def:epsilon0:prime}
d^+(x,\partial \gB_0[x]) > \epsilon_0'(\delta).
\end{equation}
(Here, $d^+(\cdot, \cdot)$ is as in \S\ref{sec:semi:markov}) and
  by $\partial \gB_0[x]$ we mean the boundary of 
  $\gB_0[x]$ as a subset of  $W^+[x]$.)
  

Let $\kappa > 1$ be as in
Proposition~\ref{prop:bilip:hattau:epsilon}, and so that
(\ref{eq:bilip:hat:tau:ij}) holds.  
Without loss of generality, assume $\delta < 0.01$. We now choose a
subset $K \subset K_0 \subset X$ with $\nu(K) > 1-\delta$ such that the
following hold:
\begin{itemize}
\item There exists a number $T_0(\delta)$ such that for any $x \in K$
  and any $T > T_0(\delta)$, 
\begin{displaymath}
\{ t \in [-T/2,T/2] \st g_t x \in K_0\} \ge 0.9 T. 
\end{displaymath}
(This can be done by the Birkhoff ergodic theorem).
\item Proposition~\ref{prop:most:inert} (a) holds.
\item Proposition~\ref{prop:some:fraction:bounded} holds. 
\item There exists a constant $C = C(\delta)$ such that for $x \in K$,
  $C_3(x)^2 < C(\delta)$ where $C_3$ is as in 
  Proposition~\ref{prop:ej:bdd:transport:bounded}. 

\item There is a constant $C'' = C''(\delta)$ such that for $x \in K$,
  $C(x) < C''(\delta)$ where $C(x)$ is as in
  Lemma~\ref{lemma:crude:divergence:subspaces} or in
  Corollary~\ref{cor:reason:cA:F:divergence}.
 Also for $x \in K$,
  the function $c_1(x)$ of
  Lemma~\ref{lemma:hausdorff:distance:to:norm} is bounded from below by
  $C''(\delta)^{-1}$. 
\item Lemma~\ref{lemma:hodge:norm:vs:dynamical:norm} holds for
  $K=K(\delta)$ and $C_1 = C_1(\delta)$. 
\item There exists a constant $C' = C'(\delta)$ such that for $x \in K$,
  $C_1(x) < C'$, $C_2(x) < C'$ and $C(x) < C'$ where $C_1(x)$,
  $C_2(x)$ and $C(x)$ are as in
  Proposition~\ref{prop:reason:cA:F}. 
Also $K \subset K'$
and also $C_1'(\delta) < C'$, $C_2'(\delta)< C'$, $C_4'(\delta) < C'$
and $C_4(\delta) < C'$ where
$K'$, $C_1'(\delta)$, $C_2'(\delta)$ and $C_4'(\delta)$ are as in
Lemma~\ref{lemma:reason:cA:F:altcondition}, and $C_4(\delta)$ is as in
Corollary~\ref{cor:reason:cA:F:divergence}. 
\item Lemma~\ref{lemma:parts:staying:close} holds for $K$. 
\item Proposition~\ref{prop:nearby:linear:maps} and
Lemma~\ref{lemma:tau:ij:nearby:close} hold for $K$ (in place of $K_0$).
\end{itemize}


Let 
\begin{displaymath}
\tilde{\cD}_{00}(q_1) =
\tilde{\cD}_{00}(q_1,K_{00},\delta,\epsilon,\eta) = \{ t > 0 \st g_t
q_1 \in K\}.   
\end{displaymath}
For $ij \in \tilde{\Lambda}$, let
\begin{displaymath}
\tilde{\cD}_{ij}(q_1) =
\tilde{\cD}_{ij}(q_1,K_{00},\delta,\epsilon,\eta) 
= \{ \hat{\tau}_{ij}(q_1,t) \st g_t q_1
\in \pi^{-1}(K), \quad t > 0 \}.
\end{displaymath}
Then by the ergodic theorem and (\ref{eq:bilip:hat:tau:ij}),
there exists a set $K_{\cD} = K_{\cD}(K_{00},\delta,\epsilon,\eta)$ 
with $\nu(K_{\cD})
\ge 1-\delta$ and 
$\ell_\cD = \ell_\cD(K_{00},\delta,\epsilon,\eta) > 0$ 
such that for $q_1 \in \pi^{-1}(K_{\cD})$ and all $ij
\in \{00\} \cup \tilde{\Lambda}$,  $\tilde{\cD}_{ij}(q_1)$ has
density at least $1-2\kappa\delta$ for $\ell > \ell_{\cD}$. Let
\begin{displaymath}
E_2(q_1,u) = E_2(q_1,u,K_{00},\delta,\epsilon,\eta) = 
\{ \ell \st g_{\hat{\tau}_{(\epsilon)}(q_1,  u, \ell)}u q_1 \in \pi^{-1}(K) \}, 
\end{displaymath}
\begin{multline*}
E_3(q_1,u) = E_3(q_1,u,K_{00},\delta,\epsilon,\eta) = \\ = \{ \ell \in
E_2(q_1,u) \st \forall ij \in \tilde{\Lambda}, 
\quad \hat{\tau}_{ij}(u q_1,\hat{\tau}_{(\epsilon)}(q_1,u,\ell)) \in
\tilde{\cD}_{ij}(q_1) \}. 
\end{multline*}
Note that $\hat{\tau}_{ij}(u q_1,\hat{\tau}_{(\epsilon)}(q_1,u,\ell)) \in
\tilde{\cD}_{ij}(q_1)$
if and only if
\begin{displaymath}
\hat{\tau}_{ij}(uq_1, \hat{\tau}_{(\epsilon)}(q_1,u,\ell)) =
\hat{\tau}_{ij}(q_1, s) \text{ and } g_s q_1 \in \pi^{-1}(K). 
\end{displaymath}

\begin{claim}
\label{claim:12.3}
There exists $\ell_3 = \ell_3(K_{00},\delta,\epsilon,\eta) > 0$, a  
set $K_3 = K_3(K_{00},\delta,\epsilon,\eta)$ of measure at least  
$1-c_3(\delta)$ and for each $q_1 \in \pi^{-1}(K_3)$ a subset $Q_3 =
Q_3(q_1,K_{00},\ell,\delta,\epsilon,\eta) \subset 
\cB(q_1,1/100)$  of
measure at least $(1-c_3'(\delta))|\cB(q_1,1/100)|$ 
such that for all $q_1 \in \pi^{-1}(K_3)$ and
$u \in Q_3$, $u q_1 \in \pi^{-1}(K)$ and
the density of $E_3(q_1,u)$ (for $\ell > \ell_3$) 
is at least 
$1-c_3''(\delta)$, and we have $c_3(\delta)$, $c_3'(\delta)$ and
$c_3''(\delta) \to 0$ as $\delta \to 0$. 
\end{claim}

\bold{Proof of claim.} We choose $K_2 = K \cap K_\cD$, and
\begin{displaymath}
K_3 = K_2 \cap \{ x \in X \st |\{ u \in \cB(x,1/100) \st ux \in K_2
\}| > (1-\delta)|\cB(x,1/100)|\}.
\end{displaymath}
Suppose $q_1 \in \pi^{-1}(K_3)$,
and $u q_1 \in \pi^{-1}(K_2)$. 
Let 
\begin{displaymath}
E_{bad} = \{ t \st g_t u q_1 \in \pi^{-1}(K^c) \}. 
\end{displaymath}
Then, since $u q_1 \in \pi^{-1}(K_\cD)$, for $\ell > \ell_{\cD}$,
 the density of $E_{bad}$ is at most $2 \kappa \delta$. We have
\begin{displaymath}
E_2(q_1,u)^c = \{ \ell \st
\hat{\tau}_{(\epsilon)}(q_1,u,\ell) \in E_{bad} \}. 
\end{displaymath}
Then, by Proposition~\ref{prop:bilip:hattau:epsilon}, for $\ell >
\kappa \ell_{\cD}$,
the density of $E_2(q_1,u)$ is at least $1-4\kappa^2\delta$.

Let
\begin{displaymath}
\hat{\cD}(q_1,u) = \hat{\cD}(q_1,u,K_{00},\delta,\epsilon,\eta) = \{ t
\st \forall ij \in \tilde{\Lambda}, \quad 
\hat{\tau}_{ij}(u q_1, t) \in \tilde{\cD}_{ij}(q_1) \}. 
\end{displaymath}
Since $q_1 \in \pi^{-1}(K_\cD)$, for each $j$, for $\ell > \ell_{\cD}$, the
density of $\tilde{\cD}_{ij}(q_1)$ is at least
$1-2\kappa \delta$. 
Then, by (\ref{eq:bilip:hat:tau:ij}), for $\ell > \kappa\ell_{\cD}$,
the density of 
$\hat{\cD}(q_1,u)$ is at least
$(1-4|\tilde{\Lambda}|\kappa^2 \delta)$.  Now
\begin{displaymath}
E_3(q_1,u) = E_2(q_1,u) \cap \{ \ell \st \hat{\tau}_{(\epsilon)}(q_1,u,\ell)
\in \hat{\cD}(q_1,u) \}.  
\end{displaymath}
Now the claim follows from Proposition~\ref{prop:bilip:hattau:epsilon}.
\mcc{(old comment: drop including $q=1$ from prop:bilipshitz)}
\qed\medskip

\begin{claim}
There exists  a set $\cD_4 = \cD_4(K_{00},\delta,\epsilon,\eta) 
\subset \reals^+$ and 
a number $\ell_4 = \ell_4(K_{00},\delta,\epsilon,\eta) > 0$ so
that $\cD_4$ has density at least $1-c_4(\delta)$ 
for $\ell > \ell_4$, and for $\ell \in \cD_4$ 
a subset $K_4(\ell) = K_4(\ell,K_{00},\delta,\epsilon,\eta)
\subset X$ with $\nu(K_4(\ell)) > 1- c_4'(\delta)$,  
such that for any $q_1 \in \pi^{-1}(K_4(\ell))$ there exists a subset
$Q_4(q_1,\ell) \subset Q_3(q_1, \ell) \subset \cB(q_1,1/100)$ 
with density at least $1-c_4''(\delta)$, so
that for all $\ell \in \cD_4$, for all $q_1 \in
\pi^{-1}(K_4(\ell))$ and all $u \in Q_4(q_1,\ell)$, 
\begin{equation}
\label{eq:D4:subset:E3}
\ell \in E_3(q_1, u) \subset E_2(q_1,u).  
\end{equation}
(We have $c_4(\delta)$, $c_4'(\delta)$ and
$c_4''(\delta) \to 0$ as $\delta \to 0$). 
\end{claim}

\bold{Proof of Claim.} This follows from Claim~\ref{claim:12.3} by
applying Fubini's theorem to
$X_\cB \cross \reals$, where $X_\cB = \{ (x,u) \st x \in X,\ u \in
\cB(x,1/100) \}$.   
\mcc{write out}
\qed\medskip

Suppose $\ell \in \cD_4$. 
We now apply Proposition~\ref{prop:can:avoid:most:Mu} 
with $K'=g_{-\ell}K_4(\ell)$. We denote
the resulting set $K$ by $K_5(\ell)= K_5(\ell,K_{00},\delta,\epsilon,\eta)$.
In view of the choice of $\epsilon_1$,  
we have $\nu(K_5(\ell)) \ge 1 - c_5(\delta)$, where $c_5(\delta) \to
0$ as $\delta \to 0$.


\medskip

Let $\cD_5 = \cD_4$ and let $K_6(\ell) = g_{\ell}K_5(\ell)$.
\mcc{eventually remove the above line}


\bold{Choice of parameters \#2: Choice of $q$, $q'$, $q_1'$ (depending
  on $\delta$, $\epsilon$, $q_1$, $\ell$).} 
Suppose $\ell \in \cD_5$ and $q_1 \in \pi^{-1}(K_6(\ell))$.  Let 
$q = g_{-\ell} q_1$. Then, $q \in \pi^{-1}(K_5(\ell))$. 
Let $\cA(q,u,\ell,t)$ be as in
\S\ref{sec:divergence:subspaces}. (Note that following our
conventions, we use the notation $\cA(q_1,u,\ell,t)$ for $q_1 \in X$, even
though $\cA(q_1,u,\ell,t)$ was originally defined for $q_1 \in X_0$.)
and for $u \in Q_4(q_1,\ell)$ let $\cM_u$ be the subspace of
Lemma~\ref{lemma:bad:subspace} applied to the linear map $\cA(q_1,u, \ell,
\hat{\tau}_{(\epsilon)}(q_1,u,\ell))$. By
Proposition~\ref{prop:can:avoid:most:Mu} and the
definition of $K_5(\ell)$, we
can choose $q' \in \cL^-[q] \cap \pi^{-1}( g_{-\ell} K_4(\ell))$ with
$\rho'(\delta) \le d^{X_0}(q,q') \le 1/100$
and so that (\ref{eq:rho:delta:le:Fq:minus:Fqprime}) and
(\ref{eq:qprime:avoids:Mu}) hold with
  $\epsilon_1(\delta) \to 0$ as
  $\delta \to 0$.  
Let $q_1' = g_\ell q'$. Then $q_1' \in \pi^{-1}(K_4(\ell))$. 


\bold{Standing Assumption.} We assume $\ell \in \cD_5$, $q_1 \in K_6(\ell)$
and $q$, $q'$, $q_1'$ are as in Choice of parameters \#2. 

\bold{Notation.} For $u \in \cB(q_1,1/100)$, $u' \in \cB(q_1',1/100)$, let
\begin{displaymath}
\tau(u) = \hat{\tau}_{(\epsilon)}(q_1, u, \ell), \qquad \tau'(u') =
\hat{\tau}_{(\epsilon)}(q_1', u', \ell).
\end{displaymath}

\bold{The maps $\psi$ and $\psi'$.} 
For $u \in \cB(q_1,1/100)$, and $u' \in \cB(q_1',1/100)$, let 
\begin{displaymath}
\psi(u) = g_{\tau(u)} u q_1,  \qquad
\psi'(u') = g_{\tau'(u')} u' q_1'.
\end{displaymath}

\begin{claim}
We have
\begin{equation}
\label{eq:RQ4:subset:K}
\psi(Q_4(q_1,\ell)) \subset \pi^{-1}(K), \quad \text{
  and } \quad \psi'(Q_4(q_1',\ell)) \subset \pi^{-1}(K). 
\end{equation}
\end{claim}

\bold{Proof of Claim.} Suppose $u \in Q_4(q_1,\ell)$. Since $q_1 \in K_4$ and 
$\ell\in \cD_4$, it follows from
(\ref{eq:D4:subset:E3}) that $\ell \in E_2(q_1, u)$, and then from the
definition of $E_2(q_1,u)$ is follows that $g_{\tau(u)} u q_1 \in \pi^{-1}(K)$. 
Hence $\psi(Q_4(q_1,\ell)) \subset \pi^{-1}(K)$. Similarly, since $q_1'
\in \pi^{-1}(K_4)$, $\psi'(Q_4(q_1',\ell)) \subset \pi^{-1}(K)$, 
proving (\ref{eq:RQ4:subset:K}). 
\qed\medskip


\bold{The numbers $t_{ij}$ and $t_{ij}'$.}
Suppose $u \in Q_4(q_1,\ell)$, and suppose $ij \in \tilde{\Lambda}$. 
Let $t_{ij}$ be defined by the equation
\begin{equation}
\label{eq:def:tj}
\hat{\tau}_{ij}(u q_1, \hat{\tau}_{(\epsilon)}(q_1,u,\ell)) =
\hat{\tau}_{ij}(q_1, t_{ij}).  
\end{equation}
Then, since $\ell \in \cD_4$ and in view of (\ref{eq:D4:subset:E3}),
we have $\ell \in E_3(q_1,u)$. In view of the definition of $E_3$, it
follows that 
\begin{equation}
\label{eq:gtj:q1:inK}
g_{t_{ij}} q_1 \in \pi^{-1}(K).
\end{equation}
Similarly, suppose $u' \in Q_4(q_1',\ell)$ and $ij \in
\tilde{\Lambda}$. 
Let $t_{ij}'$ be defined by the equation
\begin{equation}
\label{eq:def:tj:prime}
\hat{\tau}_{ij}(u' q_1', \hat{\tau}_{(\epsilon)}(q_1',u',\ell)) =
\hat{\tau}_{ij}(q_1', t_{ij}'). 
\end{equation}
Then, by the same argument, 
\begin{equation}
\label{eq:gtjprime:inK}
g_{t_{ij}'} q_1' \in \pi^{-1}(K).
\end{equation}

\bold{The map $\bfv(u)$ and the generalized subspace $\cU(u)$.}
For $u \in \cB(q_1,1/100)$, let 
\begin{equation}
\label{eq:def:bfv:u}
\bfv(u) = \bfv(q,q',u,\ell,t) = \cA(q, u, \ell, t) (F(q)-F(q'))
\end{equation}
where $t = \hat{\tau}_{(\epsilon)}(q_1,u,\ell)$, $F$ is as in
\S\ref{sec:conditional} and $\cA(\cdot, \cdot, \cdot, \cdot)$ is as in
\S\ref{sec:subsec:themap:cA}. By Proposition~\ref{prop:reason:cA:F},
we may write $\bfv(u) = \bfj(M'',v'')$, where $(M'',v'') \in
\cH_{++}(g_{\tau(u)} u q_1) \cross W^+(g_{\tau(u)} u q_1)$. Let 
$\cU(u) \equiv \cU_{g_{\tau(u)} u q_1}(M'',v'')$ denote the generalized affine
subspace corresponding to $\bfv(u)$.  
Thus, $\cU(u)$ is the approximation to $U^+[g_{\tau(u)} q_1']$ near
$g_{\tau(u)} u q_1$ defined in Proposition~\ref{prop:reason:cA:F}. 
\mcc{(motivate this better?)}

\bold{Standing Assumption.} We have $C(\delta) \epsilon < 1/100$ for
any constant $C(\delta)$ arising in the course of the proof. In
particular, this applies to $C_2(\delta)$ and $C_2'(\delta)$ in the
next claim.

\begin{claim} There exists a subset $Q_5 =
  Q_5(q_1,\ell, K_{00},\delta,\epsilon,\eta) \subset Q_4(q_1,\ell)$ 
with $|Q_5| \ge (1-c_5''(\delta))|\cB(q_1,1/100)|$  
(with $c_5''(\delta) \to 0$ as
$\delta \to 0$), and a number
$\ell_5 = \ell_5(\delta,\epsilon)$ such that for all $u \in Q_5$
and $\ell > \ell_5$, 
\begin{equation}
\label{eq:tau:u:alpha3:ell}
\tau(u) < \frac{1}{2} \alpha_3 \ell,
\end{equation}
where $\alpha_3 > 0$ is as in
Proposition~\ref{prop:stop:induction:condition} and
\S\ref{sec:subsec:themap:cA}. In addition,
\begin{equation}
\label{eq:psi:close:to:psiprime}
C_1(\delta) \epsilon \le hd_{g_{\tau(u)} u q_1}^{X_0}(U^+[g_{\tau(u)} u q_1],U^+[g_{\tau(u)} q_1']) \le
C_2(\delta) \epsilon, 
\end{equation}
\begin{equation}
\label{eq:cU:u:close:to:Uplus:psiprime}
hd_{g_{\tau(u)} u q_1}^{X_0}(U^+[g_{\tau(u)} q_1'], \cU(u)) \le C_7(\delta)
e^{-\alpha \ell},
\end{equation}
where $\alpha$ depends only on the Lyapunov spectrum. 
Also, 
\begin{equation}
\label{eq:vu:has:norm:near:epsilon}
C_1'(\delta) \epsilon \le \|\bfv(u)\| \le C_2'(\delta) \epsilon, 
\end{equation}
and if $u' \in U^+[q_1']$ is such that 
\begin{equation}
\label{eq:d:psi:u:g:tau:u:uprime}
d^{X_0}(g_{\tau(u)} u q_1, g_{\tau(u)} u' q_1') < 1/100,
\end{equation}
then $u' \in \cB(q_1',1/100)$. 
\end{claim}

\bold{Proof of claim.} Let $\cM_u$ be the subspace of
Lemma~\ref{lemma:bad:subspace} applied to the linear map $\cA(q_1,u, \ell,
\hat{\tau}_{(\epsilon)}(q_1,u,\ell))$, where $\cA( , , , )$ is as in
\S\ref{sec:divergence:subspaces}. Let $Q(q_1)$ be as in
Proposition~\ref{prop:stop:induction:condition}, so $|Q(q_1)| \ge
(1-\delta)|\cB(q_1,1/100)|$.  
Let $Q_5' \subset Q_4 \cap Q(q_1)$ be such that for all $u \in
Q_5'$, 
\begin{displaymath}
d_Y(F(q)-F(q'),\cM_u) \ge \beta(\delta)
\end{displaymath}
where $F$ is as in \S\ref{sec:conditional}. Then,
(\ref{eq:tau:u:alpha3:ell}) follows from
Proposition~\ref{prop:stop:induction:condition} and the fact that
$Q_5 \subset Q_1$.  Also, by (\ref{eq:qprime:avoids:Mu}),
\begin{displaymath}
|Q_5'| \ge |Q_4| - (\delta + \epsilon_1(\delta))|\cB(q_1,1/100)| \ge (1-
\delta- \epsilon_1(\delta) - c_4''(\delta))|\cB(q_1,1/100)|.
\end{displaymath}
Then, let 
$Q_5 = \{ u \in Q_5' \st d(u, \partial \cB(q_1,1/100)) > \delta \}$, hence 
$$|Q_5| \ge
(1-c_5'(\delta) - c_4'(\delta)-c_n \delta)|\cB(q_1,1/100)|,$$ 
where $c_n$ depends only on the dimension.

We have $C(\delta)^{-1} \epsilon \le \|\cA(q_1,u,\ell,t)\| \le C(\delta)
\epsilon$ by the definition of $t =
\hat{\tau}_{(\epsilon)}(q_1,u,\ell)$. \mcc{check}
We now apply Lemma~\ref{lemma:bad:subspace} to the linear map
$\cA(q_1,u,\ell,t)$. Then, for all $u \in Q_5$, 
\begin{displaymath}
c(\delta) \|\cA(q_1,u,\ell,t)\| \le \|\cA(q_1,u,\ell, t) (F(q) -
F(q'))\| \le \|\cA(q_1,u,\ell,t)\|. 
\end{displaymath}
Therefore, 
\begin{displaymath}
C'(\delta)^{-1} \epsilon \le 
\|\cA(q_1,u,\ell, t) (F(q) - F(q'))\| \le C'(\delta) \epsilon
\end{displaymath}
This immediately implies (\ref{eq:vu:has:norm:near:epsilon}), in view
of the definition of $\bfv(u)$. 
We now apply Proposition~\ref{prop:reason:cA:F} and
Lemma~\ref{lemma:reason:cA:F:altcondition}(a).  (We assume $\epsilon$ is
  sufficiently small so that (\ref{eq:reason:cA:F:newcondition})
  holds. Also the condition (\ref{eq:reason:cA:F:tbound}) in
Proposition~\ref{prop:reason:cA:F} holds in view of
Proposition~\ref{prop:stop:induction:condition}).  
 Now (\ref{eq:psi:close:to:psiprime}) follows from 
(\ref{eq:reason:cA:norm}). Also (\ref{eq:cU:u:close:to:Uplus:psiprime})
follows from (\ref{eq:reason:cA:F:estimate}). 

Finally, suppose $u \in Q_5$, and $u' \in U^+(q_1')$ is such that
(\ref{eq:d:psi:u:g:tau:u:uprime}) holds. 
Then, by Lemma~\ref{lemma:parts:staying:close},
we have $d^{X_0}(u q_1,u' q_1') = O_\delta(e^{-\alpha \ell})$.  
Then, assuming $\ell$ is sufficiently
large (depending on $\delta$) and using
Proposition~\ref{prop:AGY:regularity}, 
we have $u' \in \cB(q_1',1/100)$. \mcc{more details here}
\qed\medskip

\bold{Standing Assumption.} We assume $\ell > \ell_5$.

\begin{claim} Suppose $u \in Q_5(q_1,\ell)$, $u' \in
  Q_4(q_1',\ell)$ and (\ref{eq:d:psi:u:g:tau:u:uprime}) holds. 
 Then, there exists $C_0 = C_0(\delta)$ such that
\begin{equation}
\label{eq:tau:epsilon:same:as:tau:epsilon:prime}
|\hat{\tau}_{(\epsilon)}(q_1, u, \ell) -
\hat{\tau}_{(\epsilon)}(q_1',u',\ell)| \le C_0(\delta). 
\end{equation}
\end{claim}

\bold{Proof of claim.} Let $t = \hat{\tau}_{(\epsilon)}(q_1, u,
\ell)$, $t' = \hat{\tau}_{(\epsilon)}(q_1',u',\ell)$. 

By Proposition~\ref{prop:reason:cA:F} (ii), (with $q'$ and $q$
reversed) and (\ref{eq:rho:delta:le:Fq:minus:Fqprime}),
\begin{multline*}
\epsilon = \|\cA(q_1',\ell, u', t')\| 
\ge \| \cA(q_1',\ell, u', t') (F(q') - F(q)) \| \ge \\
\ge c(\delta) hd_{g_{t'} u' q_1'}^{X_0}(U^+[g_{t'}
u' q_1'], U^+[g_{t'} u q_1])
\end{multline*}
In view of Corollary~\ref{cor:reason:cA:F:divergence}(b),
(\ref{eq:tau:u:alpha3:ell})
and the fact that $g_{t'} u' q_1' \in \pi^{-1}(K)$,  
this contradicts (\ref{eq:psi:close:to:psiprime}), unless $t' <
t+C(\delta)$. 

It remains give a lower bound on $t'$. Let $\cM'$ denote the subspace
as in Lemma~\ref{lemma:bad:subspace} for $\cA(q',u',\ell,t')$. 
Note that by
Proposition~\ref{prop:can:avoid:most:Mu} (with the function $u\to
\cM_u$ the constant function $\cM'$) we can choose $q'' \in W^-[q]$
with $d_Y(F(q'') - F(q'), \cM') > \rho(\delta)$, and also so the 
upper bounds in 
(\ref{eq:rho:prime:delta:le:d:q:qprime}) and
(\ref{eq:rho:delta:le:Fq:minus:Fqprime}) hold
with $q''$ in place of $q'$. 
Then, 
\begin{displaymath}
\epsilon = \|\cA(q',\ell, u', t')\| \le c(\delta) 
\|\cA(q',\ell, u', t') (F(q'') - F(q') \|. 
\end{displaymath}
Write $q_1''= g_\ell q''$. Then, by Proposition~\ref{prop:reason:cA:F}
(ii), and Lemma~\ref{lemma:reason:cA:F:altcondition}(a),
\begin{equation}
\label{eq:tmp:hd:gt:prime:Uplus:twoprime}
hd_{g_{t'} u' q_1'}^{X_0}(U^+[g_{t'}
u' q_1'], U^+[g_{t'} q_1'']) \ge c_2(\delta) \epsilon. 
\end{equation}
By Corollary~\ref{cor:reason:cA:F:divergence}(a),
(\ref{eq:tau:u:alpha3:ell}) and 
(\ref{eq:psi:close:to:psiprime}), since $g_{t'} u' q_1'
\in \pi^{-1}(K)$, 
\begin{equation}
\label{eq:tmp:hd:u:gt:prime:Uplus:noprime}
hd_{g_{t'} u' q_1'}^{X_0}(U^+[g_{t'}
u' q_1'], U^+[g_{t'} u q_1]) \le \epsilon C(\delta) e^{-\beta(t - t')}
+ C_4(\delta) e^{-\alpha \ell},
\end{equation}
where $\alpha$ and $\beta$ depend only on the Lyapunov spectrum. 
Then, by (\ref{eq:tmp:hd:gt:prime:Uplus:twoprime}),
(\ref{eq:tmp:hd:u:gt:prime:Uplus:noprime}), and the reverse triangle
inequality,
\begin{equation}
\label{eq:tmp:gtprime:u:q1}
hd_{g_{t'} u q_1}^{X_0}(U^+[g_{t'}
u q_1], U^+[g_{t'} q_1'']) \ge \epsilon(c_2(\delta) - C(\delta)
e^{-\beta(t - t')}) -C_4(\delta) e^{-\alpha \ell}.
\end{equation}
But, 
\begin{displaymath}
\epsilon = \|\cA(q,\ell, u, t)\| \ge c_3(\delta) \|\cA(q,\ell, u, t) (F(q'') -
F(q)) \|, 
\end{displaymath}
and thus, by Proposition~\ref{prop:reason:cA:F} (ii) and
Lemma~\ref{lemma:reason:cA:F:altcondition}(a), 
\begin{displaymath}
hd_{g_{t} u q_1}^{X_0}(U^+[g_{t}
u q_1], U^+[g_{t} q_1'']) \le c(\delta) \epsilon
\end{displaymath}
In view of Corollary~\ref{cor:reason:cA:F:divergence}(b) (and the fact
that $g_t u q_1 \in \pi^{-1}(K)$) this contradicts (\ref{eq:tmp:gtprime:u:q1})
unless $t' > t - C_1(\delta)$. 
\qed\medskip

We note the following trivial lemma:
\begin{lemma}
\label{lemma:segments:equal:length}
Suppose $P$ and $P'$ are finite measure subsets of $\reals^n$ with
$|P|=|P'|$, and we have
\begin{displaymath}
P = \bigcup_{j=1}^N P_j, \quad P' = \bigcup_{j=1}^N P'_j, 
\end{displaymath}
Suppose there exists $k \in \natls$ so that any point in $P$ is
contained in  at most $k$ sets $P_j$, and also any point in $P'$ is
contained in at most $k$ sets $P_j'$. 
Also suppose 
$Q \subset P$ and $Q' \subset P'$ are subsets with $|Q| > (1-\delta)|P|$,
$|Q'| > (1-\delta)|P'|$. 

Suppose there exists $\kappa > 1$ such that for all $1 \le j \le N$ such that $P_j \cap Q \ne \emptyset$,
$|P_j| \le \kappa|P_j'|$.  Then there exists $\hat{Q} \subset Q$ with $|\hat{Q}|
\ge (1-2 \kappa k\delta)|P|$ such that if $j$ is such that
$\hat{Q} \cap P_j \ne \emptyset$, then $Q' \cap P_j'  \ne \emptyset$. 
\end{lemma}
\bold{Proof.} Let $J = \{ j \st P_j \cap Q \ne \emptyset\}$, and let 
$J' = \{ j \st Q' \cap P_j' \ne \emptyset \}$, and let
\begin{displaymath}
\hat{Q} = \{ x \in Q \st \text{ for all $j$ with $x \in P_j$, we have
  $j \in J'$.} \}
\end{displaymath}
Thus, if $x \in Q \setminus \hat{Q}$, then there exists $j \in J$ with
$x \in Q \cap P_j$ but $j \not\in J'$. 
Then, 
\begin{displaymath}
|Q\setminus \hat{Q}| \le k \sum_{j \in J\setminus J'} |Q \cap P_j| \le
k \sum_{j \in J\setminus J'} |P_j| \le \kappa k \sum_{j \not\in J'}
|P_j'| \le  \kappa k |(Q')^c|,   
\end{displaymath}
since if $j \not\in J'$ then $P_j' \subset (Q')^c$. Thus, $|Q\setminus
\hat{Q}| \le \kappa k \delta |P|$, and so $|\hat{Q}| \ge (1-2
\kappa k \delta) |P|$. 
\qed\medskip

\bold{The constant $\epsilon_0$.} 
Let \index{$\epsilon_0(\delta)$}$\epsilon_0(\delta)$ be a constant to be chosen
later. (We will choose $\epsilon_0(\delta)$ following
(\ref{eq:prelimit:points:close}) of the form
$\epsilon_0(\delta)  =
\epsilon_0'(\delta)/C(\delta)$), where $\epsilon_0'(\delta)$ is as in
(\ref{eq:def:epsilon0:prime}). 
  We will always assume that $\epsilon
< \epsilon_0(\delta)<\epsilon'(\delta)/10$. 

\begin{claim} There exists a subset $Q_6(q_1,\ell) =
  Q_6(q_1,\ell,K_{00},\delta,\epsilon,\eta) \subset Q_5(q_1,\ell)$ with 
$|Q_6(q_1,\ell)| > (1-c_6'(\delta)) |\cB(q_1,1/100)|$ 
and with $c_6'(\delta) \to 0$ as
$\delta \to 0$ such that for all $u \in Q_6(q_1,\ell)$ there exists
$u' \in Q_4(q_1',\ell)$  such that
\begin{equation}
\label{eq:alt:d:psi:u:g:tau:u:uprime}
d^{X_0}(g_{\tau(u)} u q_1, g_{\tau(u)} u' q_1') < C(\delta) \epsilon_0(\delta).
\end{equation}
\end{claim}

\bold{Proof of Claim.}
Note that the sets $\{ \cB_{\tau(u)}[u q_1] \st u \in Q_5(q_1,\ell) \}$
are a cover of $Q_5(q_1,\ell)q_1$.
Then, since these sets satisfy the condition of
Lemma~\ref{lemma:gB:properties} (b), we can
find a pairwise disjoint subcover, i.e. 
find $u_j \in Q_5(q_1,\ell)$, $1 \le j \le N$,
with $Q_5(q_1,\ell) q_1 = \bigcup_{j=1}^N \cB_{\tau(u_j)}[u_j q_1]$ and
so that
$\cB_{\tau(u_j)}[u_j q_1]$ and $\cB_{\tau(u_k)}[u_k q_1]$ are disjoint for
$j \ne k$. Let
\begin{displaymath}
\cB_j \equiv g_{\tau(u_j)} \cB_{\tau(u_j)}[u_j q_1]=
\cB_0[g_{\tau(u_j)} u_j q_1] \subset \tilde{X}_0
\end{displaymath}
In view of (\ref{eq:def:epsilon0:prime}),
Proposition~\ref{prop:AGY:regularity}, and the Besicovich covering
lemma, there exists $k$, depending only on the dimension 
and points $x_{j,1}, \dots, x_{j,m(j)} \subset \cB_j$
such that 
\begin{displaymath}
\pi^{-1}(K) \cap \cB_j \subset \bigcup_{m=1}^{m(j)} B^{X_0}(x_{j,m},\epsilon_0(\delta))
\cap U^+[g_{\tau(u_j)} u q_1],
\end{displaymath}
and also so that for a fixed $j$, each point 
is contained in at most $k$ balls $B^{X_0}(x_{j,m},
\epsilon_0(\delta))$.
Since $\epsilon_0(\delta) < \epsilon_0'(\delta)/10$,
in view of (\ref{eq:def:epsilon0:prime}) and
(\ref{eq:Pjm:Pjm:prime:close}), the same is true without fixing $j$. 

For $1 \le j \le N$ and $1 \le m \le m(j)$, let
\begin{displaymath}
P_{j,m} = \{ u \in \cB(q_1,1/100) \st g_{\tau(u_j)} u q_1 \in
B^{X_0}(x_{j,m}, \epsilon_0(\delta))\},
\end{displaymath}
and let 
\begin{displaymath}
P_{j,m}' = \{ u' \in \cB(q_1',1/100) \st g_{\tau(u_j)} u' q_1' \in
B^{X_0}(x_{j,m}, \epsilon_0(\delta))\}.
\end{displaymath}
By construction, each point is contained in at most $k$ sets
$P_{j,m}$, and at most $k$ sets $P_{j,m}'$. 

By (\ref{eq:psi:close:to:psiprime}) applied to $u_j$, 
\begin{equation}
\label{eq:Pjm:Pjm:prime:close}
hd^{X_0}_{g_{\tau(u_j)} u_j q_1} (U^+[g_{\tau(u_j)} u_j q_1], 
  U^+[g_{\tau(u_j)} q_1']) \le C_2(\delta) \epsilon.
\end{equation}
Suppose $\epsilon > 0$ is sufficiently small (depending on $\delta$)
so that Lemma~\ref{lemma:parts:staying:close} holds with
$C_2(\delta) \epsilon$ in place of $\epsilon$.  
Since for all $x \in X_0$, $\cB_0[x] \subset B^{X_0}(x,1/200)$
we have $d^{X_0}(x_{j,m},g_t u_j q_1) < 1/200$, and
\begin{equation}
\label{eq:tmp:restricted:hd}
B^{X_0}(x_{j,m}, \epsilon_0(\delta)) \subset B^{X_0}(g_t u_j
q_1, 1/100). 
\end{equation}
By Lemma~\ref{lemma:parts:staying:close}, for $1 \le j \le N$, $1 \le
m \le m(j)$, 
provided $\cB_j \cap Q_5(q_1, \ell) \ne \emptyset$, we have
$\kappa^{-1}|P_{j,m}| \le |P_{j,m}'| \le
\kappa|P_j|$, where $\kappa$ depends only on the Lyapunov
spectrum, and we have normalized the measures $| \cdot |$ so that
$|U^+[q_1] \cap B^+(q_1,1/100)| = |U^+[q_1'] \cap
B^+(q_1',1/100)|=1$. Let $m(0) = 1$ and let 
\begin{displaymath}
P_{0,1} = \cB(q_1,1/100) \setminus \bigcup_{j=1}^N
\bigcup_{m=1}^{m(j)} P_{j,m}, \qquad P_{0,1}' = \cB(q_1',1/100) \setminus \bigcup_{j=1}^N
\bigcup_{m=1}^{m(j)} P_{j,m}'. 
\end{displaymath}
Then,
\begin{displaymath}
\cB(q_1,1/100) = \bigcup_{j=0}^N \bigcup_{m=1}^{m(j)} P_{j,m}, \qquad
\cB(q_1',1/100) = \bigcup_{j=0}^N \bigcup_{m=1}^{m(j)} P_{j,m}'. 
\end{displaymath}
Then, applying Lemma~\ref{lemma:segments:equal:length} with $P =
\cB(q_1,1/100)$, $P' = \cB(q_1', 1/100)$, 
$Q =
Q_5(q_1,\ell)$, $Q' = Q_4(q_1',\ell)$, we get a set $\hat{Q} \equiv
Q_6(q_1,\ell)$ with $|Q_6(q_1,\ell)| \ge
(1-c_6'(\delta))|\cB(q_1,1/100)|$ where $c_6'(\delta) \to 0$ as $\delta
\to 0$, 
so that, in view of (\ref{eq:tmp:restricted:hd}) and
the definitions of $P_{j,m}$ and $P_{j,m}'$, for any
$u \in Q_6(q_1,\ell)$ there exists $u_j \in
Q_5(q_1,\ell)$ with $u q_1 \in \cB_{\tau(u_j)}[u_j q_1]$ and
$u' \in Q_4(q_1',\ell)$ with 
\begin{equation}
\label{eq:tmp:u:prime:almost}
d^{X_0}(g_{\tau(u_j)} u q_1, g_{\tau(u_j)} u' q_1') \le \epsilon_0(\delta).
\end{equation}
It remains to replace $\tau(u_j)$ by $\tau(u)$ in
(\ref{eq:tmp:u:prime:almost}). This can be done as follows:
Since $u q_1 \in \cB_{\tau(u_j)}[u_j q_1]$, we have, by
(\ref{eq:psi:close:to:psiprime}) applied to $u_j$ and 
Lemma~\ref{lemma:algebraic:divergence},
\begin{displaymath}
C_2(\delta)^{-1} \epsilon \le
hd_{g_{\tau(u_j)} u q_1}(U^+[g_{\tau(u_j)} u q_1], U^+[g_{\tau(u_j)}
q_1']) \le C_2(\delta) \epsilon
\end{displaymath}
Then, since $g_{\tau(u)} u q_1 \in \pi^{-1}(K)$, 
by (\ref{eq:psi:close:to:psiprime}),
(\ref{eq:cU:u:close:to:Uplus:psiprime}), (\ref{eq:tau:u:alpha3:ell}) and 
Corollary~\ref{cor:reason:cA:F:divergence}, we have
\begin{equation}
\label{eq:tmp:tau:u:tau:uj}
|\tau(u) - \tau(u_j)| \le C_1(\delta).
\end{equation}
Then, provided $\epsilon$ is small enough depending on $\delta$,
(\ref{eq:alt:d:psi:u:g:tau:u:uprime})  follows from
(\ref{eq:tmp:u:prime:almost}), (\ref{eq:tmp:tau:u:tau:uj}),
and Lemma~\ref{lemma:forni:upper}. 
\qed\medskip

\begin{claim}
There exists a constants $c_7(\delta) > 0$ and $c_7'(\delta)$ with
$c_7(\delta) \to 0$ and $c_7'(\delta) \to 0$ as $\delta \to 0$ and a
subset $K_7(\ell) = K_7(\ell,K_{00},\delta,\epsilon,\eta)$ with
$K_7(\ell) \subset K_6(\ell)$ and $\nu(K_7(\ell)) > 1-
c_7(\delta)$ such that for $q_1 \in \pi^{-1}(K_7(\ell))$, 
\begin{displaymath}
|\cB(q_1) \cap Q_6(q_1,\ell)| \ge (1-c_7'(\delta)) |\cB(q_1)|. 
\end{displaymath}
\end{claim}

\bold{Proof of Claim.}Recall that in view of
  Proposition~\ref{prop:semi:markov},  $\cB(q_1)
  \subset \cB(q_1,1/100)$.
Given $\delta > 0$, there exists $c_7''(\delta)
> 0$ with $c_7''(\delta) \to 0$ as $\delta \to 0$ and  a compact set
$K_7' \subset X$ with $\nu(K_7') > 1-c_7''(\delta)$, 
such that for $q_1 \in \pi^{-1}(K_7')$, $|\cB(q_1) \cap \cB(q_1,1/100)| \ge
c_6'(\delta)^{1/2} |\cB(q_1,1/100)|$. Then, for $q_1 \in \pi^{-1}(K_7'
\cap K_6)$,
\begin{displaymath}
|\cB(q_1) \cap Q_6(q_1,\ell)^c| \le |Q_6(q_1,\ell)^c| \le
c_6'(\delta) |\cB(q_1,1/100)| \le c_6'(\delta)^{1/2} |\cB(q_1)|. 
\end{displaymath}
Thus, the claim holds with $c_7(\delta) = c_6(\delta) + c_7''(\delta)$
and $c_7'(\delta) = c_6'(\delta)^{1/2}$. 
\qed\medskip


\bold{Standing Assumption.} We assume that $q_1 \in \pi^{-1}(K_7(\ell))$. 
\medskip

The next few claims will help us choose $u$ (once the other parameters
have been chosen). Let 
\begin{displaymath}
Q_7(q_1,\ell) = \cB(q_1) \cap Q_6(q_1,\ell)
\end{displaymath}
\mcc{rearrange other standing assumptions so this gets merged}

\begin{claim}
There exists a subset $Q_7^*(q_1,\ell) =
Q_7^*(q_1,\ell,K_{00},\delta,\epsilon,\eta) \subset 
Q_7(q_1,\ell)$ with $|Q_7^*| \ge (1-c_7^*(\delta)) |\cB(q_1)|$
such that for $u \in Q_7^*$ and any
$\ell > \ell_7(\delta)$
we have
\begin{displaymath}
  |\cB_\ell(u q_1) \cap Q_7(q_1,\ell)| \ge
  (1-c_7^*(\delta)) |\cB_\ell(u q_1)|, 
\end{displaymath}
where $c_7^*(\delta) \to 0$ as $\delta \to 0$. 
\end{claim}
\bold{Proof.} This follows immediately from
Lemma~\ref{lemma:cB:vitali:substitute}. 
\qed\medskip

\mcc{move the following claim to a ``tracking'' lemma in another subsection}

\begin{claim}
\label{claim:12:12}
There exist a number $\ell_8 =
  \ell_8(K_{00},\delta,\epsilon,\eta)$ and a constant $c_8(\delta)$ with
$c_8(\delta) \to 0$ as $\delta \to 0$ and for every $\ell > \ell_8$ 
a subset $Q_8(q_1,\ell) = Q_8(q_1,\ell,K_{00},\delta,\epsilon,\eta) 
\subset \cB(q_1)$ with $|Q_8(q_1,\ell)| \ge
(1-c_8(\delta)) |\cB(q_1)|$ so that 
for $u \in Q_8(q_1,\ell)$  we have
\begin{equation}
\label{eq:Ru:R2prime:u:close:Eplus}
d\left(\frac{\bfv(u)}{\|\bfv(u)\|}, \bE(g_{\tau(u)} u q_1)\right) \le
C_8(\delta) e^{-\alpha' \ell},  
\end{equation}
where $\bfv(u)$ is defined in (\ref{eq:def:bfv:u})
and $\alpha'$ depends only on the Lyapunov spectrum.  
\end{claim}

\bold{Proof of claim.} 
Let $L' > L_2(\delta)$ be a constant to be chosen later, where
$L_2(\delta)$ is as in Proposition~\ref{prop:most:inert} (a). Also let
$\ell_8 = \ell_8(\delta,\epsilon,K_{00},\eta)$ be a constant to be
chosen later. Suppose $\ell > \ell_8$, 
and suppose $u \in Q_7^*(q_1,\ell)$, so in particular $g_{\tau(u)} u q_1 \in
\pi^{-1}(K)$. Let $t \in [L',2L']$ 
be such that Proposition~\ref{prop:most:inert} (a) holds for $\bfv =
\bfv(u)$ and $x = g_{\tau(u)} u q_1$. 

Let  
$B_u \subset \cB(q_1)$ denote
$\cB_{\hat{\tau}_{(\epsilon)}(q_1,u,\ell)-t}(u q_1)u$,
(where $\cB_t(x)$ is defined in \S\ref{sec:divergence:subspaces}).    
Suppose $u_1 \in B_u \cap Q_7(q_1,\ell)$, and write
\begin{displaymath}
g_{\tau(u_1)} u_1 q_1 = g_s u_2 g_t^{-1} g_{\tau(u)} u q_1. 
\end{displaymath}
Then, $u_2 \in \cB(g_t^{-1} g_{\tau(u)} u q_1)$ and $t \le  2L'$.

We now claim that
\begin{equation}
\label{eq:s:kappa:L:prime:new}
s \le \frac{1}{2} \kappa t + C_0(\delta) \le \kappa L' + C_0(\delta)
\end{equation}
where $\kappa$
depends only on the Lyapunov spectrum. Let
\begin{displaymath}
\cU_t = U^+[g_{-t} g_{\tau(u)} u q_1],\qquad \cU'_t = U^+[g_{-t}
g_{\tau(u)} q_1']. 
\end{displaymath}
By Corollary~\ref{cor:reason:cA:F:divergence}(b) applied at the point
$g_{\tau(u)} u q_1 \in \pi^{-1}(K)$, 
\begin{displaymath}
hd^{X_0}_{g_{-t} g_{\tau(u)} u q_1} (\cU_t, \cU'_t)
\ge C(\delta) \epsilon e^{-\beta t} - c_0(\delta)e^{-\alpha \ell},
\end{displaymath}
where $\beta$ depends only on the Lyapunov spectrum, 
and by Corollary~\ref{cor:reason:cA:F:divergence}(a) applied at the
point $g_{\tau(u_1)} u_1 q_1 \in \pi^{-1}(K)$, 
\begin{displaymath}
hd^{X_0}_{u_2 g_{-t} g_{\tau(u)} u q_1} (\cU_t,\cU'_t)
\le c(\delta) \epsilon e^{-2 s} + c_0(\delta) e^{-\alpha \ell}
\end{displaymath}
where $\beta'$ also depends only on the Lyapunov spectrum. 
Also, by Lemma~\ref{lemma:algebraic:divergence}, 
\begin{displaymath}
hd^{X_0}_{g_{-t} g_{\tau(u)} u q_1} (\cU_t, \cU'_t) \ge c_1 \,
hd^{X_0}_{u_2 g_{-t} g_{\tau(u)} u q_1} (\cU_t, \cU'_t) -c_0(\delta)
e^{-\alpha \ell}
\end{displaymath}
where $c_1$ is an absolute constant. 
Therefore, 
\begin{displaymath}
\epsilon C(\delta) e^{-\beta t} -c_0(\delta) e^{-\alpha \ell} \le c_1 (c(\delta)
\epsilon e^{-2 s} + c_0(\delta)e^{-\alpha \ell}).  
\end{displaymath}
This implies (\ref{eq:s:kappa:L:prime:new}), assuming
that $\ell$ is sufficiently large depending on $\epsilon$. 

Since $u \in
Q_6(q_1,\ell)$, (\ref{eq:psi:close:to:psiprime})
and (\ref{eq:cU:u:close:to:Uplus:psiprime}) hold. 
Therefore, 
\begin{displaymath}
hd_{g_{\tau(u_1)} u_1 q_1}((g_s u_2 g_t^{-1}) \cU(u), U^+[g_{\tau(u_1)} q_1']) =
O(e^{\kappa' L'} e^{-\alpha \ell}),
\end{displaymath}
where $\kappa'$ and $\alpha$ depend only on the Lyapunov spectrum. 
Thus, using (\ref{eq:cU:u:close:to:Uplus:psiprime}) at the point
$g_{\tau(u_1)} u_1 q_1 \in \pi^{-1}(K)$,
\begin{displaymath}
hd_{g_{\tau(u_1)} u_1 q_1}((g_s u_2 g_t^{-1}) \cU(u), \cU(u_1)) = O(e^{\kappa' L'}
e^{-\alpha \ell}). 
\end{displaymath}
Therefore,
\begin{equation}
\label{eq:tracking:one}
\| (g_s u_2 g_t^{-1})_* \bfv(u) - \bfv(u_1) \| = O(e^{\kappa' L'}
e^{-\alpha \ell}). 
\end{equation}
In view of
(\ref{eq:vu:has:norm:near:epsilon}), $\|\bfv(u_1) \| \approx
\epsilon$. Thus, $\|(g_s u_2 g_t^{-1})_* \bfv(u)\| \approx
  \epsilon$, and
\begin{displaymath}
\left\| \frac{(g_s u_2 g_t^{-1})_* \bfv(u)}
{\|(g_s u_2 g_t^{-1})_* \bfv(u)\|} - \frac{\bfv(u_1)}{\|\bfv(u_1)\|}
\right\| = O_\epsilon( e^{\kappa'L'-\alpha\ell}). 
\end{displaymath}
But, by Proposition~\ref{prop:most:inert} (a), for $1-\delta$ fraction
of $u_2 \in \cB(g_t^{-1}g_{\tau(u)} u q_1)$, 
\begin{displaymath}
d\left( \frac{(g_s u_2 g_{-t})_* \bfv(u)}{\|(g_s u_2 
    g_{-t})_* \bfv(u)\|} , 
     \bE(g_{\tau(u_1)} u_1 q_1) \right) \le C(\delta) e^{-\alpha L'}, 
\end{displaymath}
Note that
\begin{displaymath}
\cB(g_t^{-1}g_{\tau(u)} u q_1) = g_{\hat{\tau}_{(\epsilon)}(q_1, u, \ell) - t}
B_u. 
\end{displaymath}
Therefore, for $1-\delta$ fraction of $u_1 \in B_u$, 
\begin{equation}
\label{eq:tmp:claim:cover}
d\left(\frac{\bfv(u_1)}{\|\bfv(u_1)\|}, \bE(g_{\tau(u_1)} u_1 q_1)\right) 
\le C(\epsilon,\delta) [e^{\kappa'L'-\alpha \ell} + e^{-\alpha L'}]
\end{equation}
We can now choose $L' > 0$ to be $\alpha' \ell$ where $\alpha' > 0$ is
a small constant depending only on the Lyapunov  spectrum, 
and $\ell_8 > 0$ so that for $\ell > \ell_8$ 
the right-hand-side of the above equation is at most $e^{-\alpha'
  \ell}$. 

The collection of balls $\{B_u\}_{u \in Q_7^*(q_1,\ell)}$ 
are a cover of $Q_7^*(q_1,\ell)$. These balls satisfy
the condition of Lemma~\ref{lemma:gB:properties} (b); hence we may
choose a pairwise disjoint subcollection which still covers 
$Q_7^*(q_1,\ell)$. Then, by summing (\ref{eq:tmp:claim:cover}), we see
that 
(\ref{eq:Ru:R2prime:u:close:Eplus}) holds for $u$ in a subset
$Q_8 \subset \cB[q_1]$ of measure at least $(1-c_8(\delta))|\cB[q_1]| =
(1-\delta)(1-c_7^*(\delta))|\cB[q_1]| $. 
\qed\medskip

\begin{claim}
There exists a subset $Q_8^*(q_1,\ell) =
Q_8^*(q_1,\ell,K_{00},\delta,\epsilon,\eta) \subset 
Q_8(q_1,\ell)$ with $|Q_8^*| \ge (1-c_8^*(\delta)) |\cB(q_1)|$
such that for $u \in Q_8^*$ and any
$t > \ell_8(\delta)$
we have
\begin{displaymath}
  |\cB_t(u q_1) \cap Q_8(q_1,\ell)| \ge
  (1-c_8^*(\delta)) |\cB_t(u q_1)|, 
\end{displaymath}
where $c_8^*(\delta) \to 0$ as $\delta \to 0$. 
\end{claim}
\bold{Proof.} This follows immediately from
Lemma~\ref{lemma:cB:vitali:substitute}. 
\qed\medskip

\bold{Choice of parameters \#3: Choice of $\delta$.}
Let $\theta' = (\theta/2)^n$, where $\theta$
and $n$ are as in Proposition~\ref{prop:some:fraction:bounded}. We can
choose $\delta > 0$ so that 
\begin{equation}
\label{eq:choice:of:delta}
c_8^*(\delta) < \theta'/2. 
\end{equation}

\mcc{move the following proof to a ``tracking'' lemma in another subsection}

\begin{claim} There exist sets $Q_9(q_1,\ell) = Q_9(q_1,\ell,K_{00},
  \delta,\epsilon, \eta) \subset Q_8^*(q_1,\ell)$ 
with $|Q_9(q_1,\ell)| \ge (\theta'/2)(1-\theta'/2) |\cB(q_1)|$ and $\ell_9 =
\ell_9(K_{00}, \delta,\epsilon, \eta)$, 
such that for $\ell > \ell_9$ and $u \in Q_9(q_1,\ell)$, 
\begin{equation}
\label{eq:psi:psiprime:close:to:bigcup}
d\left(\frac{\bfv(u)}{\|\bfv(u)\|},
  \bigcup_{ij \in \tilde{\Lambda}}   \bE_{[ij],bdd}(g_{\tau(u)} u q_1)\right) < 4 \eta. 
\end{equation}
\end{claim}

\bold{Proof of claim.}  
Suppose $u \in Q_8^*(q_1,\ell)$. Then, by
(\ref{eq:Ru:R2prime:u:close:Eplus}) and (\ref{eq:vu:has:norm:near:epsilon}),
we may write
\begin{displaymath}
\bfv(u) = \bfv'(u) + \bfv''(u),
\end{displaymath}
where $\bfv'(u) \in \bE(g_{\tau(u)} u q_1)$ and $\|\bfv''(u) \| \le
C(\delta,\epsilon) e^{-\alpha' \ell}$. Arguing in the same way as in
the proof of Claim~\ref{claim:12:12}, we see that for
$(1-O(\delta))$-fraction of $y \in  \cF_{\bfv'(u)}[g_{\tau(u)} u q_1,L]$,
we have $y \in g_{[-1,1]}K$. Then, 
by Proposition~\ref{prop:some:fraction:bounded} applied with $L =
L_0(\delta,\eta)$ and $\bfv = \bfv'(u)$, we get that for a
at least $\theta'$-fraction of $y \in
\cF_{\bfv'}[g_{\tau(u)} u q_1,L]$, 
\begin{displaymath}
d\left(\frac{R(g_{\tau(u)} u q_1,y) \bfv'(u)}{\|R(g_{\tau(u)} u q_1,y) \bfv'(u)\|}, 
\bigcup_{ij \in\tilde{\Lambda}} 
  \bE_{[ij],bdd}(y)\right) < 2 \eta. 
\end{displaymath}
Note that by Proposition~\ref{prop:properties:dynamical:norm:Hbig}
(d), for $y \in
\cF_{\bfv'}[g_{\tau(u)} u q_1,L]$, $\|R(g_{\tau(u)} u q_1,y)\| \le
e^{\kappa^2 L}$, where $\kappa$ is as in
Proposition~\ref{prop:properties:dynamical:norm:Hbig}.
Then, for at least $\theta'$-fraction of $y \in
\cF_{\bfv'}[g_{\tau(u)} u q_1,L]$, 
\begin{equation}
\label{eq:frac:R:psi:u:y}
d\left(\frac{R(g_{\tau(u)} u q_1,y) \bfv(u)}{\|R(g_{\tau(u)} u q_1,y) \bfv(u)\|}, 
\bigcup_{ij \in\tilde{\Lambda}} 
  \bE_{[ij],bdd}(y)\right) < 3 \eta + C(\epsilon,\delta)e^{2\kappa^2L}
e^{-\alpha' \ell}. 
\end{equation}
Let $B_u = \cB_{\hat{\tau}_{(\epsilon)}(q_1,u,\ell)-L}(u q_1)u$.  
In view of (\ref{eq:tracking:one}) and (\ref{eq:vu:has:norm:near:epsilon})
\mcc{explain better} 
there exists $C = C(\epsilon,\delta)$ such that 
\begin{displaymath}
\cF_{\bfv'}[g_{\tau(u)} u q_1,L] \cap \pi^{-1}(K) \subset g_{[-C,C]} \psi(B_u) \quad\text{
  and } \quad \psi(B_u) \cap \pi^{-1}(K) \subset g_{[-C,C]}
\cF_{\bfv'}[g_{\tau(u)} u q_1,L].  
\end{displaymath}
Then, by (\ref{eq:frac:R:psi:u:y}) and (\ref{eq:choice:of:delta}), 
for $(\theta'/2)$-fraction of $u_1 \in B_u$, $g_{\tau(u_1)} u_1 q_1 \in \pi^{-1}(K)$ and 
\begin{displaymath}
d\left(\frac{R(g_{\tau(u)} u q_1,g_{\tau(u_1)} u_1 q_1) \bfv(u)}{\|R(g_{\tau(u)} u q_1,g_{\tau(u_1)} u_1 q_1) \bfv(u)\|}, 
\bigcup_{ij \in\tilde{\Lambda}} 
  \bE_{[ij],bdd}(g_{\tau(u_1)} u_1 q_1)\right) < C_1(\epsilon,\delta) [3 \eta +
e^{2\kappa^2L}e^{-\alpha'\ell}).] 
\end{displaymath}
Then, by (\ref{eq:tracking:one}), \mcc{(make that into separate
  lemma. Referring to that statement is confusing)}
for $(\theta'/2)$-fraction of $u_1 \in B_u$, 
\begin{displaymath}
d\left(\frac{\bfv(u_1)}{\|\bfv(u_1)\|}, 
\bigcup_{ij \in\tilde{\Lambda}} 
  \bE_{[ij],bdd}(g_{\tau(u_1)} u_1 q_1)\right) < C_2(\epsilon,\delta) [3 \eta +
e^{2\kappa^2L}e^{-\alpha'\ell} + e^{-\alpha' \ell}].
\end{displaymath}
Hence, we may choose $\ell_9 = \ell_9(K_{00},\epsilon,\delta,\eta)$ so that
for $\ell > \ell_9$ the right-hand side of the above equation is at
most $4\eta$. Thus, (\ref{eq:psi:psiprime:close:to:bigcup}) holds for
$(\theta'/2)$-fraction of $u_1 \in B_u$. 

The collection of balls $\{B_u\}_{u \in Q_8^*(q_1,\ell)}$ 
are a cover of $Q_8^*(q_1,\ell)$. These balls satisfy
the condition of Lemma~\ref{lemma:gB:properties} (b); hence we may
choose a pairwise disjoint subcollection which still covers 
$Q_8^*(q_1,\ell)$. Then, by summing over the disjoint subcollection, we see
that the claim holds on a set $E$ of measure at least
$(\theta'/2)|Q_8^*| \ge (\theta'/2)(1-c_8^*(\delta)) \ge
(\theta'/2)(1-\theta'/2)$. 
\qed\medskip


\bold{Choice of parameters \#4: Choosing $\ell$, $q_1$, $q$, $q'$, $q_1'$.} 
Choose $\ell > \ell_9(K_{00},\epsilon,\delta,\eta)$. 
Now choose $q_1 \in K_7(\ell)$, and let $q$, $q'$, $q_1'$ be as in Choice of
Parameters \#2.

\bold{Choice of parameters \#5: Choosing $u$, $u'$, $q_2$, $q_2'$,
  $ij$, $q_{3,ij}$,
  $q_{3,ij}'$ (depending on $q_1$, $q_1'$, $u$, $\ell$).}
Choose $u \in Q_9(q_1,\ell)$, $u' \in
Q_4(q_1',\ell)$  so that
(\ref{eq:psi:close:to:psiprime}) and
(\ref{eq:cU:u:close:to:Uplus:psiprime})  hold.
We have $\psi(u) = g_{\tau(u)} u q_1 \in \pi^{-1}(K)$ and $\psi'(u')
\in \pi^{-1}(K)$. By (\ref{eq:tau:epsilon:same:as:tau:epsilon:prime}), 
\begin{displaymath}
|\hat{\tau}_{(\epsilon)}(q_1, u, \ell) -
\hat{\tau}_{(\epsilon)}(q_1',u',\ell)| \le C_0(\delta), 
\end{displaymath}
therefore,
\begin{displaymath}
g_{\tau(u)} u' q_1' \in \pi^{-1}(g_{[-C,C]} K),
\end{displaymath}
where $C = C(\delta)$.

By the definition of $K$ we can
find $C_4(\delta)$ and 
$s \in [0, C_4(\delta)]$ such that 
\begin{displaymath}
q_2 \equiv g_s g_{\tau(u)} u q_1 \in \pi^{-1}(K_0), \qquad q_2' \equiv g_s g_{\tau(u)} u' q_1' \in \pi^{-1}(K_0). 
\end{displaymath}
In view of (\ref{eq:psi:close:to:psiprime}),
(\ref{eq:cU:u:close:to:Uplus:psiprime}), the fact that $s \in
[0,C_4(\delta)]$ and Corollary~\ref{cor:reason:cA:F:divergence}(a) we get
\begin{equation}
\label{eq:prelimit:hd:close}
\frac{1}{C(\delta)} \epsilon \le hd^{X_0}_{q_2}(U^+[q_2],
U^+[q_2'])  \le C(\delta) \epsilon. 
\end{equation}
By (\ref{eq:alt:d:psi:u:g:tau:u:uprime}), the fact that $s \in
[0,C_4(\delta)]$ and Lemma~\ref{lemma:forni:upper} we get 
\begin{equation}
\label{eq:prelimit:points:close}
d^{X_0}(q_2,q_2') = d^+(q_2,q_2') \le
C(\delta) \epsilon_0(\delta).  
\end{equation}
We now choose $\epsilon_0(\delta)$ so that $C(\delta)
\epsilon_0(\delta) < \epsilon_0'(\delta)$, where 
$C(\delta)$ is as in
(\ref{eq:prelimit:points:close}), and $\epsilon_0'(\delta)$ is as in
(\ref{eq:def:epsilon0:prime}).  

Let $ij \in \tilde{\Lambda}$ be
such that 
\begin{equation}
\label{eq:bfvu:close:to:bEij:bdd}
d\left(\frac{\bfv(u)}{\|\bfv(u)\|}, \bE_{[ij],bdd}(g_{\tau(u)} u q_1) \right) \le
4 \eta.
\end{equation}
By  Lemma~\ref{lemma:tau:ij:nearby:close}, 
\begin{displaymath}
|\hat{\tau}_{ij}(uq_1, \hat{\tau}_{(\epsilon)}(q_1, u, \ell)) -
\hat{\tau}_{ij}(u'q_1', \hat{\tau}_{(\epsilon)}(q_1, u, \ell))| \le
C_4'(\delta).
\end{displaymath}
Then, by (\ref{eq:tau:epsilon:same:as:tau:epsilon:prime}) and
(\ref{eq:bilip:hat:tau:ij}),  
\begin{displaymath}
|\hat{\tau}_{ij}(uq_1, \hat{\tau}_{(\epsilon)}(q_1, u, \ell)) -
\hat{\tau}_{ij}(u'q_1', \hat{\tau}_{(\epsilon)}(q_1', u', \ell))| \le
C_4''(\delta).
\end{displaymath}
Hence, by Proposition~\ref{prop:properties:dynamical:norm:Hbig} (e)
(cf. Lemma~\ref{lemma:stay:in:same:orbit}), (\ref{eq:def:tj}) and
(\ref{eq:def:tj:prime}), 
\begin{equation}
\label{eq:t:ij:t:ij:prime:close}
|t_{ij} - t_{ij}'| \le C_5(\delta). 
\end{equation}
Therefore, by (\ref{eq:gtj:q1:inK}) and (\ref{eq:gtjprime:inK}), we
have
\begin{displaymath}
g_{t_{ij}} q_1 \in \pi^{-1}(K), \quad \text{ and } \quad g_{t_{ij}} q_1' \in
\pi^{-1}(g_{[-C_5(\delta),C_5(\delta)]} K). 
\end{displaymath}
By the definition of $K$, we can find $s'' \in [0, C_5''(\delta)]$
such that 
\begin{displaymath}
q_{3,ij} \equiv g_{s''+t_{ij}} q_1 \in \pi^{-1}(K_0), \quad \text{ and } \quad
q_{3,ij}' \equiv g_{s''+t_{ij}} q_1' \in \pi^{-1}(K_0).  
\end{displaymath}
Let $\tau = s+\hat{\tau}_{(\epsilon)}(q_1,u,\ell)$,
$\tau' = s''+t_{ij}$. Then,
\begin{displaymath}
q_2 = g_\tau u q_1, \quad q_2' = g_\tau u' q_1', \quad  q_{3,ij} =
g_{\tau'} q_1, \quad q_{3,ij}' = g_{\tau'} q_1'. 
\end{displaymath}
We may write $q_2 = g_t^{ij} u q_1$, $q_{3,ij} = g_{t'}^{ij} q_1$.
Then,  in view of (\ref{eq:t:ij:t:ij:prime:close}) and
(\ref{eq:bilip:hat:tau:ij}), 
\begin{displaymath}
|t-t'| \le C_6(\delta). 
\end{displaymath}
We note that by Proposition~\ref{prop:stop:induction:condition}, 
$\ell > \alpha_0 \tau$, where $\alpha_0$ depends only on
the Lyapunov spectrum.

\bold{Taking the limit as $\eta \to 0$.} For fixed
  $\delta$ and $\epsilon$, 
we now take a sequence of $\eta_k \to 0$ (this forces $\ell_k \to
\infty$) and  pass to
limits along a subsequence. Let $\tilde{q}_2 \in K_0$ be the limit of the
$q_2$, and $\tilde{q}_2'\in K_0$ be the limit of the $q_2'$.
We may also
assume that along the subsequence $ij \in \tilde{\Lambda}$ is fixed, 
where $ij$ is as in (\ref{eq:bfvu:close:to:bEij:bdd}). By passing to
the limit in (\ref{eq:prelimit:hd:close}), we get
\begin{equation}
\label{eq:postlimit:hd}
\frac{1}{C(\delta)} \epsilon \le hd^{X_0}_{\tilde{q}_2}(U^+[\tilde{q}_2],
U^+[\tilde{q}_2'])  \le C(\delta) \epsilon.
\end{equation}
We now apply Proposition~\ref{prop:nearby:linear:maps} (with $\xi \to
0$ as $\eta_k \to 0$). By (\ref{eq:goal:same:sheet}), $\tilde{q}_2'
\in W^+[\tilde{q}_2]$. By applying $g_s$ to
(\ref{eq:bfvu:close:to:bEij:bdd}) and then passing to the limit, we get
$U^+[\tilde{q}_2'] \in \cE_{ij}(\tilde{q}_2)$. Finally, 
it follows from passing to the limit in
(\ref{eq:prelimit:points:close}) that
$d^+(\tilde{q}_2,\tilde{q}_2') \le \epsilon_0'(\delta)$, and thus,
since $\tilde{q}_2 \in K_0$ and $\tilde{q}_2' \in K_0$,
it follows from
(\ref{eq:def:epsilon0:prime}) that
$\tilde{q}_2' \in \gB_0[\tilde{q}_2]$. Hence, 
\begin{displaymath}
\tilde{q}_2' \in \cC_{ij}(\tilde{q}_2).
\end{displaymath}
Now, by (\ref{eq:nearby:maps:final}), we have
\begin{displaymath}
f_{ij}(\tilde{q}_2) \propto P^+(\tilde{q}_2, \tilde{q}_2')_* f_{ij}(\tilde{q}_2'). 
\end{displaymath}
This concludes the proof of Proposition~\ref{prop:half:inductive:step}. 
We have $\tilde{q}_2 \in \pi^{-1}(K_0) \subset \pi^{-1}(K_{00} \cap
K_*)$, and $\tilde{q}_2' \in \pi^{-1}(K_0 \subset K_*$). 
\qed\medskip

\bold{Applying the argument for a sequence of $\epsilon$'s tending to $0$}.
Take a sequence $\epsilon_n \to 0$. We now apply
Proposition~\ref{prop:half:inductive:step} with $\epsilon =
\epsilon_n$. After passing to a subsequence, we may assume $ij$ is
constant. We get, for each $n$ a set $E_n \subset K_*$ 
with $\nu(E_n) > \delta_0$ and with the property that
for every $x \in E_n$ there exists $y \in \cC_{ij}(x) \cap K_*$ 
such that (\ref{eq:dxy:within:C:epsilon}) and
(\ref{eq:half:inductive:step:two}) hold for $\epsilon = \epsilon_n$.  
Let 
\begin{displaymath}
F = \bigcap_{k=1}^\infty \bigcup_{n=k}^\infty E_n \subset K_*,
\end{displaymath}
(so $F$ consists of the points which are in infinitely many
$E_n$). Suppose $x \in F$. Then
there exists a sequence $y_n \to x$ such that $y_n \in \cC_{ij}[x]$,
$y_n \not\in U^+[x]$, and
so that $f_{ij}(y_n) \propto P^+(x,y_n)_* f_{ij}(x)$. Then, (on the set where
both are defined)
\begin{displaymath}
f_{ij}(x) \propto (\gamma_n)_* f_{ij}(x),
\end{displaymath}
where $\gamma_n \in \cG_{++}(x)$ is the affine map whose linear part is
$P^+(x,y_n)$ and whose translational part is $y_n-x$. (Here we have used
the fact that $y_n \in \cC_{ij}[x]$, and thus by the definition of
conditional measure, $f_{ij}(y_n) = (y_n-x)_* f_{ij}(x)$, where
$(y_n-x)_*: W^+(x) \to W^+(x)$ is translation by $y_n-x$.)

Let $\tilde{f}_{ij}(x)$ denote the measure on $\cG_{++}(x)$ given by
\begin{displaymath}
\tilde{f}_{ij}(x)(h) = \int_{W^+[x]} \bar{h} \, df_{ij}(x),
\end{displaymath}
where for a compactly supported real-valued continuous function $h$ on
$\cG_{++}(x)$, $\bar{h}: W^+[x] \to \reals$ is given by
\begin{displaymath}
\bar{h}(gx) = \int_{\cQ_{++}(x)} h(g q) \, dm(q),
\end{displaymath}
where $m$ is the Haar measure on $\cQ_{++}(x)$. (Thus,
$\tilde{f}_{ij}(x)$ is the pullback of $f_{ij}(x)$ from $W^+[x] \isom
\cG_{++}(x)/\cQ_{++}(x)$ to $\cG_{++}(x)$). 
Then, 
\begin{equation}
\label{eq:gamma:n:star:tilde:f}
(\gamma_n)_* \tilde{f}_{ij}(x) \propto \tilde{f}_{ij}(x)
\end{equation}
on the set where both are defined. 

For $x \in X$, let $U_{new}^+(x)$ denote the maximal connected
subgroup of $\cG_{++}(x)$ 
such that for $u \in U_{new}^+(x)$, (on the domain where both
are defined),
\begin{equation}
\label{eq:u:star:tilde:fij}
(u)_* \tilde{f}_{ij}(x) \propto \tilde{f}_{ij}(x).
\end{equation}
By (\ref{eq:gamma:n:star:tilde:f}) and
Proposition~\ref{prop:measure:inv:closed:subgroup}, for $x\in F$, $U_{new}^+(x)$ strictly contains $U^+(x)$.

Suppose $x \in F$, $y \in F$ and $y \in \cC_{ij}[x]$. Then, since
$\tilde{f}_{ij}(y) = Tr(x,y)_* \tilde{f}_{ij}(x)$, we have that
(\ref{eq:u:star:tilde:fij}) holds for $u \in Tr(y,x)
U_{new}^+(y)$ (see Lemma~\ref{lemma:trivial:affine:maps}). 
Therefore, by the maximality of $U_{new}^+(x)$,
for $x \in F$, $y \in F \cap \cC_{ij}[x]$,  
\begin{equation}
\label{eq:Tr:y:x:Unew}
Tr(y,x) U_{new}^+(y) = U_{new}^+(x). 
\end{equation}
Suppose $x \in F$, $t < 0$ and $g_t x \in F$. Then, since the
measurable partition $\cC_{ij}$ is $g_t$-equivariant (see
Lemma~\ref{lemma:cCij:equivariant}) we have that
(\ref{eq:u:star:tilde:fij}) holds for $u \in g_{-t} U^+(g_t
x)$. Therefore, by the maximality of $U_{new}^+(x)$,
for $x \in F$, $t < 0$  with $g_t x \in F$ we have
\begin{equation}
g_{-t} U_{new}^+(g_t x) = U_{new}^+(x),
\label{eq:g:minus:t:U:new}
\end{equation}
and (\ref{eq:u:star:tilde:fij}) and (\ref{eq:Tr:y:x:Unew}) still
hold.

From (\ref{eq:u:star:tilde:fij}), we get that for $x \in F$ and $u \in
U_{new}^+(x)$, 
\begin{equation}
\label{eq:scaling:cocycle}
(u)_* \tilde{f}_{ij}(x) = e^{\beta_x(u)} \tilde{f}_{ij}(x),
\end{equation}
where $\beta_x: U_{new}^+(x) \to \reals$ is a homomorphism. 
Since $\nu(F) > \delta_0 > 0$ and $g_t$ is
ergodic, for almost
all $x \in X$ there exist arbitrarily large $t > 0$ so that $g_{-t} x
\in F$. Then, we define $U_{new}^+(x)$ to be $g_t
U_{new}^+(g_{-t}x)$. (This is consistent in view of
(\ref{eq:g:minus:t:U:new})). Then, 
(\ref{eq:scaling:cocycle}) holds for a.e.\ $x \in X$.
It follows from (\ref{eq:scaling:cocycle}) that for a.e.\ $x \in X$,
$u \in U_{new}^+(x)$ and $t > 0$, 
\begin{equation}
\label{eq:beta:x:equivariant}
\beta_{g_{-t} x}( g_{-t} u g_t) = \beta_x(u).
\end{equation}
We can write
\begin{displaymath}
\beta_x(u) = L_x(\log u),
\end{displaymath}
where $L_x: \Lie(U^+)(x) \to \reals$ is a Lie algebra homomorphism
(which is in particular a linear map). Let $K \subset X$ be a positive
measure set for which there exists a constant $C$ with $\|L_x\| \le C$
for all $x \in K$. Now for almost all $x \in X$ and $u \in
U^+_{new}(x)$ there exists a
sequence $t_j \to \infty$ so that $g_{-t_j} x\in K$ and $g_{-t_j} u
g_{t_j} \to e$, where $e$ is the identity element of
$U^+_{new}$. Then, (\ref{eq:beta:x:equivariant}) applied to the
sequence $t_j$ implies that $\beta_x(u) = 0$
almost everywhere (cf.\ \cite[Proposition
7.4(b)]{Benoist:Quint}). Therefore, for almost all $x \in X$, 
the conditional measure of $\nu$ along the orbit $U_{new}^+[x]$ is the
push-forward of the Haar measure on $U_{new}^+(x)$. 

The partition whose atoms are $U_{new}^+[x]$ is given by the
refinement of the measurable partition $\cC_{ij}$ into orbits of an
algebraic group. (For the atom $\cC_{ij}[x]$ this group is
$U_{new}^+(y)$ for almost any $y \in \cC_{ij}[x]$; in view of
(\ref{eq:Tr:y:x:Unew}) and Lemma~\ref{lemma:trivial:affine:maps}, this
group, viewed as a group of affine maps of $W^+[x]$ 
is independent of the choice
of $y$). Therefore the partition whose atoms are sets of the form
$U_{new}^+[x] \cap \gB_0[x]$ is a measurable partition.

In view of (\ref{eq:Tr:y:x:Unew}), and since for $u$ near the identity, 
$U_{new}^+[x] \subset \cC_{ij}[x]$ 
we have
that (\ref{eq:Uplus:ux:Ux:weird}) holds for $U_{new}^+$. 
Then, it also holds for any
$u$ in view of $g_t$-equivariance. 
Finally, since $U_{new}^+(x) \supset U^+(x)$ and $U^+(x) \supset \exp
N(x)$, we have $U_{new}^+(x) \supset N(x)$. 

Similarly, recall that the measure $\nu$ on $X$ is the pullback of the
measure on $X_0$ such that the conditionals on the fibers of the covering
map $\sigma_0: X \to X_0$ are the counting measure.

By (\ref{eq:def:Wplus:X}) there exists a subset $\Omega_0 \subset X_0$
of full measure such that for any $x_0 \in \Omega_0$, 
for any $x \in \sigma^{-1}(x_0)$ 
we have an (almost-everywhere defined) identification $\sigma_x$ between $W^+[x]
\subset X$ and $W^+[x_0] \subset X_0$ and under this identification, the conditional
measures coincide, i.e.\  $(\sigma_x)_* \nu_{W^+[x]} =
\nu_{W^+[x_0]}$. Suppose $x_0 \in \Omega_0$ and $x \in
\sigma_0^{-1}(x_0)$. After removing from $\Omega_0$ a set of measure
$0$, we may assume that
Definition~\ref{def:compatible:family}(iii) holds for $x$ and
$U^+_{new}(x)$.  Therefore it also holds for $x_0$ and $\sigma_x \circ
U^+_{new}(x) \circ \sigma_x^{-1} \subset \cG_{++}(x_0)$.  Now for $x_0
\in \Omega_0$ define $U_{new}^+(x_0)$ to be the group
generated by all the groups $\sigma_x \circ
U^+_{new}(x) \circ \sigma_x^{-1}$ where $x$ varies over
$\sigma_0^{-1}(x_0)$. Then,
Definition~\ref{def:compatible:family}(iii) holds for $x_0$ and
$U^+_{new}(x_0)$. In the same way, all of the other parts of
Definition~\ref{def:compatible:family} hold for $x_0$ and
$U^+_{new}(x_0)$ since they hold for $x$ and $U^+_{new}(x)$ for
any $x \in \sigma_0^{-1}(x_0)$.

This completes the proof of 
Proposition~\ref{prop:inductive:step}.
\qed\medskip


\mcc{worry about normalizations of and domains of definition of 
conditional measures}

\mcc{double check normalizations in
  Proposition~\ref{prop:measure:inv:closed:subgroup}.}

\section{Proof of Theorem~\ref{theorem:main:partone}}
\label{sec:proof:theorem:partone}


Let $L^-$, $L^+$, $S^+$ be as in \S\ref{sec:subsec:stopping}. 
Apply
Proposition~\ref{prop:inductive:step} to get an equivariant system of
subgroups $U_{new}^+(x) \subset \cG_{++}(x)$ which is compatible with
$\nu$ in the sense of Definition~\ref{def:compatible:family}. 

We have that $L^-[x]$ is smooth at $x$ for almost all $x \in X$, see 
\cite[\S{3}]{Avila:Eskin:Moeller:yeti}.
Let $T_\reals U^+(x) \subset W^+(x)$ denote the tangent subspace at
$x$ to the smooth manifold $U^+[x]$, and let $T_\reals L^-(x)
\subset W^-(x)$ denote the tangent subspace to $L^-[x]$ at $x$. (This
exists for almost all $x$).

If $L^+[x] \not \subset S^+[x]$ we can apply 
Proposition~\ref{prop:inductive:step} again and repeat the
process. When this process stops, the following hold:
\begin{itemize}
\item[{\rm (a)}] $L^+[x] \subset S^+[x] \subset U^+[x]$. In
  particular,
\begin{displaymath}
T_\reals L^+(x) \equiv \hat{\pi}_x^+ \circ (\hat{\pi}_x^{-1}) T_\reals
L^-(x) \subset T_\reals U^+(x).  
\end{displaymath}
\item[{\rm (b)}] The conditional measures $\nu_{U^+[x]}$ are
  induced from the Haar measure on $U^+[x]$.  These measures are
  $g_t$-equivariant.  
\item[{\rm (c)}] The subspaces $T_\reals U^+(x) \subset W^+(x)$ is 
$P = AN$ equivariant. (This follows from the fact that the $N$
direction is contained in $U^+(x)$, (\ref{eq:Uplus:ux:Ux:weird})
and the fact that the $N$ direction is in the center of $\cG_{++}(x)$). 
\mcc{explain more?}
The subspaces $T_\reals L^-(x)$ are $g_t$-equivariant. 
\item[{\rm (d)}] The conditional measures $\nu_{W^-[x]}$ are supported
  on $L^-[x]$. 
\end{itemize}
Let $H^1_\perp$
denote the subspace of $H^1(M,\Sigma,\reals)$ which is
orthogonal to the $SL(2,\reals)$ orbit, see (\ref{eq:def:H1perp}). 
Let $I$ denote
the Lyapunov exponents (with multiplicity) \mcc{fixme}
of the cocycle in $T_\reals U^+(x) \cap H^1_\perp$, 
$J$ denote the
Lyapunov exponents of the cocycle in $T_\reals L^+(x) \cap
H^1_\perp$. By (a), we have $J \subset I$. 

Since $T_\reals U^+(x) \cap H^1_\perp(x)$ is $AN$-invariant, by
Theorem~\ref{theorem:P:invariant:positive:sum} we have,
\begin{equation}
\label{eq:inv:subspace:pos:sum}
\sum_{i \in I} \lambda_i \ge 0. 
\end{equation}
We now compute the entropy of $g_t$. We have, by
Theorem~\ref{theorem:mt:theorem97}(i) (applied to the flow in the
reverse direction), 
\mcc{(use $\Delta$ notation)}
\begin{equation}
\label{eq:entropy:W:plus}
\frac{1}{t} h(g_t, \nu) \ge 2+ \sum_{i \in I} (1+\lambda_i) = 2 + |I| +
\sum_{i\in I} \lambda_i \ge 2 + |I|
\end{equation}
where the $2$ comes from the direction of $N$, 
and for the last estimate we used (\ref{eq:inv:subspace:pos:sum}). 
Also, by Theorem~\ref{theorem:mt:theorem97}(ii),
\begin{align}
\label{eq:entropy:W:minus}
\frac{1}{t} h(g_{-t}, \nu) & \le 2 + \sum_{j \in J}
                                           (1-\lambda_j),   
&& \text{where the $2$ is the potential contribution of $\bar{N}$}  \notag \\ 
& \le 2 + \sum_{i \in I} (1-\lambda_i) &&\text{since $(1-\lambda_i) \ge
  0$ for all $i$} \notag \\
& \le 2+|I| && \text{by (\ref{eq:inv:subspace:pos:sum})}   
\end{align}
However, $h(g_t, \nu) = h(g_{-t}, \nu)$. Therefore,
all the inequalities in (\ref{eq:entropy:W:plus}) and
(\ref{eq:entropy:W:minus}) are in fact equalities. In particular,  $I
= J$, i.e.\  
\begin{equation}
\label{eq:tangent:Lplus:Uplus}
T_\reals L^+(x) = T_\reals U^+(x). 
\end{equation}
Since $L^+[x] \subset S^+[x]$ and $S^+[x]$ is closed and star-shaped
with respect to $x$, it follows that \mcc{explain more?}
\begin{equation}
\label{eq:tangent:Lplus:Splus}
T_\reals L^+[x] \subset S^+[x].
\end{equation}
Since $S^+[x] \subset U^+[x]$, we get, in view of
(\ref{eq:tangent:Lplus:Uplus}) and (\ref{eq:tangent:Lplus:Splus}) that
\begin{displaymath}
T_\reals U^+[x] \subset S^+[x] \subset U^+[x]. 
\end{displaymath}
Thus $U^+[x]$ is an affine subspace of $W^+[x]$.  Then, in view of
(\ref{eq:tangent:Lplus:Uplus}), and the fact that $L^+[x] \subset
U^+[x]$, we get that $L^+[x] = U^+[x]$. Thus, $L^+[x]$ is an affine
subspace, hence $\cL^-(x) = L^-(x)$. 

We have
\begin{displaymath}
\frac{1}{t} h_{\nu}(g_{-t},W^-) = 2 + \sum_{i \in I} (1-\lambda_i). 
\end{displaymath}
By applying Theorem~\ref{theorem:mt:theorem97}(iii) to the affine
subspaces $\cL^-(x)$, this implies that the
conditional measures $\nu_{\cL^-}(x)$ are Lebesgue, and that
$\nu$ is $\bar{N}$-invariant (where $\bar{N}$ is as in
\S\ref{sec:theorems}). 
Hence $\nu$ is $SL(2,\reals)$-invariant. 

By the definition of $\cL^-$, the conditional measures $\nu_{W^-[x]}$
are supported on $\cL^-[x]$. Thus, the conditional measures
$\nu_{W^-[x]}$ are (up to null sets) precisely the Lebesgue measures
on $\cL^-[x]$. 

Let $\cU^+[x]$ denote the smallest linear subspace of $W^+[x]$ which
contains the support of $\nu_{W^+[x]}$. 
Since $\nu$ is $SL(2,\reals)$-invariant, we can argue
by symmetry that the conditional measures $\nu_{W^+[x]}$ are precisely
the Lebesgue measures on $\cU^+[x]$. Since $U^+[x]$ accounts for all
the entropy of the flow, we must have
$\cU^+[x] = U^+[x]$. \mcc{say some more words here}
Since $\cU^+[x] = \cL^+[x]$, 
this completes the proof of
Theorem~\ref{theorem:main:partone}. 
\qed\medskip

\section{Random walks}
\label{sec:random:walks}

In all of \S\ref{sec:random:walks}-\S\ref{sec:martingale}, we work
with the finite cover $X_0$ (which is a manifold), and do not use the
measurable cover $X$. 

We choose a compactly supported absolutely continuous
measure \index{$\mu$}$\mu$ on $SL(2,\reals)$. We also assume that $\mu$ is spherically
symmetric. 
Let $\nu$ be any ergodic $\mu$-stationary probability measure
on $X_0$. By Furstenberg's theorem \cite[Theorem 1.4]{Nevo:Zimmer}, 
\begin{displaymath}
\nu = \frac{1}{2\pi}\int_0^{2\pi} (r_\theta)_* \nu_0 \, d\theta
\end{displaymath}
where $r_\theta$ is as in \S\ref{sec:theorems} 
and $\nu_0$ is a measure invariant under $P = AN \subset
SL(2,\reals)$. Then, by Theorem~\ref{theorem:main:partone}, $\nu_0$ is
$SL(2,\reals)$-invariant. Therefore the stationary measure $\nu$ is
also in fact $SL(2,\reals)$-invariant.


We can think of $x \in X_0$ as a point in $H^1(M,\Sigma, \cx)$. 
For a subspace $U(x) \subset H^1(M,\Sigma,\reals)$ let
$\index{$U_\cx$}U_\cx = \cx 
\tensor U(x)$ denote its complexification, which is a subspace of 
$H^1(M,\Sigma,\cx)$. In all cases we will consider, $U(x)$
will either contain the space spanned by $\Re x$ and $\Im x$ or will
be symplectically orthogonal to that space. 

Let $\index{$area$@$\area(x,1)$} \area(x,1)
\subset H^1(M,\Sigma,\cx)$ denote the set of $y \in
H^1(M,\Sigma,\cx)$ such that $x+y$ has area $1$. 
We often abuse notation by referring to $U_\cx(x) \cap \area(1,x)$
also as $U_\cx(x)$. We also write $U_{\cx}[x]$ for the
corresponding subset of $X_0$.

The map $p: H^1(M,\Sigma,\reals) \to
H^1(M,\reals)$ naturally extends to a map (also denoted by
$\index{$p$}p$) from $H^1(M,\Sigma, \cx) \to H^1(M,\cx)$.

By Theorem~\ref{theorem:main:partone}, there is a
$SL(2,\reals)$-equivariant 
family of subspaces $\index{$U(x)$}U(x) \subset
H^1(M,\Sigma,\reals)$ containing $\Re x$ and $\Im x$ and
such that the conditional measures of $\nu$
along $U_{\cx}[x]$ are Lebesgue. Furthermore, for almost all $x$,
the conditional measure of $\nu$ along $W^+[x]$ is supported on
$W^+[x] \cap U_{\cx}[x]$, and the conditional measure of $\nu$ along
$W^-[x]$ is supported on $W^-[x] \cap U_{\cx}[x]$.

\begin{lemma}
\label{lemma:volume:form:on:U}
There exists a volume form \index{$V$@$\Vol(x)$}$d\Vol(x)$ on $U(x)$ which is invariant
under the $SL(2,\reals)$ action. This form is
  non-degenerate on compact subsets of $X_0$.
\end{lemma}

\bold{Proof.} The subspaces $p(U(x))$ form an invariant subbundle
$p(U)$ of the
Hodge bundle. By Theorem~\ref{theorem:KZ:semisimple} (a), (after
passing to a finite cover) we may assume that 
$p(U)$ is a direct sum of irreducible subbundles. Then, by
Theorem~\ref{theorem:KZ:semisimple} (b), we have a decomposition
\begin{displaymath}
p(U)(x) = U_{symp}(x) \dirsum U_0(x)
\end{displaymath}
where the symplectic form on $U_{symp}$ is non-degenerate, the
decomposition is orthogonal with respect to the Hodge inner product,
and $U_0$ is isotropic. Then, by
Theorem~\ref{theorem:isotropic:zero:lyapunov} and
Theorem~\ref{theorem:zero:lyapunov:isometric} 
the Hodge inner product on $U_0$ is equivariant under the
$SL(2,\reals)$ action. 

Then we can define the volume form on $p(U)$ 
to be the
product of the appropriate power of the symplectic form on $U_{symp}$
and the volume form induced by the Hodge inner product on
$U_0$. The subbundle $U_{symp}$ is clearly
  $SL(2,\reals)$ equivariant. By
\cite[Corollary~5.4]{Filip:semisimplicity}, 
applied to the section $c_1 \wedge \dots
\wedge c_k$
where $\{ c_1, \dots, c_k \}$ is a symplectic basis for $U_{symp}$,  
it follows that 
the symplectic volume form on $U_{symp}$ agrees with the volume
form induced by the Hodge inner product on $U_{symp}$ 
(which is non-degenerate on compact sets). This gives a volume form on
$p(U)$ with the desired properties. 

Since the Kontsevich-Zorich cocycle acts trivially on $\ker p$, the normalized 
Lebesgue measure on $\ker p$ is well defined. Thus, 
the volume form on $p(U)$
naturally induces a volume form on $U$. \mcc{explain last sentence better}
\qed\medskip

\bold{Remark.} In fact it follows from the results of
\cite{Avila:Eskin:Moeller:yeti} 
that $U_0$ is trivial. 

\begin{lemma}
\label{lemma:exists:invariant:complement}
There exists an $SL(2,\reals)$-equivariant subbundle 
 $p(U)^\perp \subset H^1(M,\reals)$ such that 
\begin{displaymath}
p(U)(x) \dirsum p(U)^\perp(x) = H^1(M,\reals). 
\end{displaymath}
\end{lemma}

\bold{Proof.} This follows from the proof of
Theorem~\ref{theorem:KZ:semisimple}. 
\qed\medskip

\bold{The subbundles $\cL_k$.} 
By Theorem~\ref{theorem:KZ:semisimple} we have
\begin{equation}
\label{eq:KZ:Lk:decomp}
p(U)^\perp(x) = \bigoplus_{k\in \hat{\Lambda}} \index{$L$@$\cL_k(x)$}\cL_k(x),
\end{equation}
where \index{$\Lambda$@$\hat{\Lambda}$}$\hat{\Lambda}$ is an indexing set not containing $0$, and
for each $k \in \hat{\Lambda}$, $\cL_k$ is an $SL(2,\reals)$-equivariant
subbundle of the Hodge bundle. (In our notation, the action of the
Kontsevich-Zorich cocycle may permute some of the $\cL_k$.)
\mccc{(fix this and choose better terminology!)}
Note that $\cL_k(x)$ is symplectically
orthogonal to the $SL(2,\reals)$ orbit of $x$. 
Without loss of generality, we may assume that the decomposition
(\ref{eq:KZ:Lk:decomp}) is maximal, in the sense that
on any (measurable) finite cover of $X_0$  each $\cL_k$ does not
contain a non-trivial proper $SL(2,\reals)$-equivariant subbundle. (If
this was not true, we could without passing to a finite cover, write
a version of (\ref{eq:KZ:Lk:decomp}) with a
larger $k$).
If $U$ does not contain 
the kernel of $p$, then we let $\hat{\lambda}_0 = 0$, and let
$\index{$\Lambda$@$\tilde{\Lambda}$}\tilde{\Lambda} 
= \hat{\Lambda} \cup \{ 0 \}$. 

\bold{The Forni subbundle.} Let
\index{$\lambda$@$\tilde{\lambda}_k$}$\tilde{\lambda}_k$ denote the
top Lyapunov exponent of the geodesic flow $g_t$ restricted to $\cL_k$. Let 
\begin{displaymath}
\index{$F$}F(x) = \bigoplus_{\{k \st \tilde{\lambda}_k = 0\} } \cL_k(x). 
\end{displaymath}
We call $F(x)$ the {\em Forni} subspace of $\nu$. The subspaces $F(x)$
form a subbundle of the Hodge bundle which we call the Forni
subbundle. It is an $SL(2,\reals)$-invariant subbundle, on which the
Kontsevich-Zorich cocycle acts by Hodge isometries. In particular, all
the Lyapunov exponents of $F(x)$ are $0$. Let $F^\perp(x)$ denote the
orthogonal complement to $F(x)$ in the Hodge norm. 
By Theorem~\ref{theorem:properties:forni} (b), 
\begin{displaymath}
F^\perp(x) = \bigoplus_{\{k \st \hat{\lambda}_k \ne 0 \}} \cL_k(x). 
\end{displaymath}
The following is proved in \cite{Avila:Eskin:Moeller:yeti}:
\begin{theorem}
\label{theorem:no:yeti}
There exists a subset $\Phi$ of the stratum with $\nu(\Phi) = 1$ such
that for all $x \in \Phi$ there exists a neighborhood $\cU(x)$ such that
for all $y \in \cU(x) \cap \Phi$ we have $p(y-x) \in F^\perp_{\cx}(x)$. 
\end{theorem}

\bold{The backwards shift map.} \mcc{discuss past and future}
Let \index{$B$}$B$ be the space of (one-sided) infinite sequences of
elements of $SL(2,\reals)$. (We think of $B$ as giving the ``past''
trajectory of the random walk.) 
Let $T: B \to B$ be the shift map. (In our interpretation, \index{$T$}$T$ takes
us one step into the past). 
We define the skew-product map $T: B \cross X_0 \to B \cross X_0$ by 
\begin{displaymath}
T(b,x) = (Tb, b_0^{-1} x), \qquad \text{ where $b = (b_0, b_1, \dots )$}
\end{displaymath}
(Thus the shift map and the skew-product map are
  denoted by the same letter.)
We define the measure \index{$\beta$}$\beta$ on $B$ to be $\mu \cross
\mu \cdots$. The skew product map $T$ naturally acts on the bundle
$H^1(M,\reals)$, and thus on each $\cL_k$ for $k \in \hat{\Lambda}$.  

For each $k\in \hat{\Lambda}$, by the multiplicative ergodic theorem
we have the Lyapunov flag for this action (with respect to the
invariant measure $\beta$): 
\begin{displaymath}
\{ 0 \} = \index{$V$@$\cV_{\le j}^{(k)}$}\cV_{\le 0}^{(k)} \subset \cV_{\le 1}^{(k)}(b,x) \subset \dots
\cV_{\le n_k}^{(k)}(b,x) = \cL_k(x). 
\end{displaymath}

By the multiplicative ergodic theorem applied to the action of
$SL(2,\reals)$ on $\reals^2$, for $\beta$-almost all $b \in B$, 
\begin{displaymath}
\index{$\sigma_0$}\sigma_0 =  \lim_{n \to \infty} \frac{1}{n} \log \| b_0 \dots b_n \|
\end{displaymath}
where $\sigma_0 > 0$ is the Lyapunov exponent for the measure $\mu$  on
$SL(2,\reals)$. Then, the Lyapunov exponents of the flow $g_t$ and the
Lyapunov exponents of the skew-product map $T$ differ by a factor of
$\sigma_0$. Let \index{$\lambda$@$\hat{\lambda}_k$}$\hat{\lambda}_k$
denote the top Lyapunov exponent of $T$ restricted to $\cL_k$. 

\bold{The two-sided shift space.}
Let \index{$B$@$\tilde{B}$}$\tilde{B}$ denote the two-sided shift space. We 
denote the measure $\cdots \cross \mu \cross \mu
\cross \cdots$ on $\tilde{B}$ also by $\beta$.

\bold{Notation.} For
$a, b \in B$ let
\begin{equation}
\label{eq:def:vee}
\index{$a \vee b$}a \vee b = ( \dots, a_2, a_1, b_0, b_1, \dots ) \in \tilde{B}.
\end{equation}
(Note that the indexing for $a \in B$ starts at $1$ not at $0$.)
If $\omega = a \vee b \in \tilde{B}$, we think of the sequence 
\begin{displaymath}
\dots, \omega_{-2}, \omega_{-1} = \dots a_2, a_1  
\end{displaymath}
as the ``future'' of the random walk trajectory. (In general,
following \cite{Benoist:Quint}, we use the symbols $b$, $b'$ etc. to
refer to the ``past'' and the symbols $a$, $a'$ etc. to refer to the
``future'').

\bold{The opposite Lyapunov flag.} Note that on the two-sided shift
space $\tilde{B} \cross X_0$, the map $T$ is invertible.
Thus, for each $a \vee b \in B$, we have the Lyapunov flag for
$T^{-1}$: 
\begin{displaymath}
\{ 0 \} = \cV_{\ge n_k}^{(k)} \subset \cV_{\ge n_k-1}^{(k)}(a,x) \subset \dots
\cV_{\ge 0}^{(k)}(a,x) = \cL_k(x).
\end{displaymath}
(As reflected in the above notation, this flag depends only on the
``future'' i.e.\ ``$a$'' part of $a \vee b$). 

\bold{The top Lyapunov exponent $\hat{\lambda}_k$.}
Recall that $\hat{\lambda}_k \ge 0$ denotes the top Lyapunov exponent in
$\cL_k$. Then, (since $T$ steps into the past), for $v \in \cV_{\le 1}^{(k)}(b,x)$, 
\begin{equation}
\label{eq:tmp:top:lyapunov}
\lim_{n \to \infty} \frac{1}{n} \log \frac{\|T^n(b,x)_* v\|}{\|v\|} =
-\hat{\lambda}_k. 
\end{equation}
In the above equation we used the notation $T^n(b,x)_*$ to denote the
action of $T^n(b,x)$ on $H^1(M,\reals)$. 

Also, for $v \in \cV_{> 1}^{(k)}(a,x)$, for some $\alpha > 0$, 
\begin{displaymath}
\lim_{n \to \infty} \frac{1}{n} \log \frac{\|T^{-n}(a \vee b,x)_* v\|}{\|v\|} <
\hat{\lambda}_k - \alpha.
\end{displaymath}
Here, $\alpha$ is the minimum over $k$ of the difference between the
top Lyapunov exponent in $\cL_k$ and the next Lyapunov exponent.

\bigskip
The following lemma is a consequence of the zero-one law
Lemma~\ref{lemma:non:atomic:lyapunov}(i): 
\begin{lemma}
\label{lemma:zero:one:law:star}
For every $\delta > 0$ and every $\delta' > 0$ 
there exists $E_{good} \subset X_0$ with 
$\nu(E_{good}) > 1 - \delta$ and $\sigma =
\sigma(\delta,\delta') > 0$, 
such that for any $x \in E_{good}$, any $k$  
and any vector $w \in \proj(\cL_k(x))$, 
\begin{equation}
\label{eq:random:subspace}
\beta\left( \{ a \in B \st d_Y(w,\cV^{(k)}_{> 1}(a,x)) > \sigma \}
\right) > 1-\delta' 
\end{equation}
\end{lemma}
(In (\ref{eq:random:subspace}), 
$d_Y(\cdot,\cdot)$ is the distance on the projective space
$\proj(H^1(M,\reals))$ derived from the AGY norm.)

\bold{Proof.} It is enough to prove the lemma for a fixed $k$. 
For $F \subset Gr_{n_k-1}(\cL_k(x))$ (the Grassmannian of
$n_k-1$ dimensional subspaces of $\cL_k(x)$) let
\begin{displaymath}
\hat{\nu}^{(k)}_x(F) = \beta\left( \{ a \in B \st
  \cV^{(k)}_{> 1}(a,x) \in F \} \right),
\end{displaymath}
and let $\hat{\nu}^{(k)}$ denote the measure on the bundle $X_0 \cross
Gr_{n_k-1}(\cL_k)$ given by  
\begin{displaymath}
d\hat{\nu}^{(k)}(x,L) = d\nu(x) \, d\hat{\nu}^{(k)}_x(L). 
\end{displaymath}
Then, $\hat{\nu}^{(k)}$ is a stationary measure for
the (forward) random walk. 
For $w \in \proj(\cL_k(x))$ let $I(w) = \{ L \in
Gr_{n_k-1}(\cL_k(x)) \st w \in L \}$. 
Let 
\begin{displaymath}
Z = \{ x \in X_0 \st  \hat{\nu}^{(k)}_x(I(w)) > 0 \text{ for some  $w
  \in \proj(\cL_k(x)$)} \},
\end{displaymath}
Suppose $\nu(Z) > 0$. Then, for each $x \in Z$ we can choose $w_x \in
\proj(\cL_k(x))$ such that $\hat{\nu}^{(k)}_x(I(w_x)) > 0$. \mcc{check
  measurability?} 
Then,
\begin{equation}
\label{eq:tmp:hat:nu:k}
\hat{\nu}^{(k)}\left( \bigcup_{x\in Z} \{x\} \cross I(w_x) \right) > 0.
\end{equation}
Therefore, (\ref{eq:tmp:hat:nu:k}) holds for some ergodic component of
$\hat{\nu}^{(k)}$. However, this contradicts
Lemma~\ref{lemma:non:atomic:lyapunov} (i), since 
by the definition of 
$\cL_k$, the action of the cocycle on $\cL_k$ is strongly
irreducible. Thus, $\nu(Z) = 0$ and $\nu(Z^c) = 1$. By definition, 
for all $x \in Z^c$ and all $w \in \cL_k(x)$, 
\begin{displaymath}
\beta\left( \{ a \in B \st w \in \cV^{(k)}_{> 1}(a,x) \} \right) = 0. 
\end{displaymath}
Fix $x \in Z^c$. Then, for every $w \in \proj(\cL_k(x))$ there
exists $\sigma_0(x,w,\delta') > 0$ such that 
\begin{displaymath}
\beta\left( \{ a \in B \st  d_Y(\cV^{(k)}_{> 1}(a,x),w) >
  2\sigma_0(x,w,\delta') \} 
\right) > 1 - \delta'.  
\end{displaymath}
Let $\cU(x,w) = \{ z \in \proj(\cL_k(x)) \st d_Y(z,w) < \sigma_0(x,w,\delta')
\}$. Then the $\{\cU(x,w)\}_{w \in \proj(\cL_k(x))}$ form an open
cover of the compact space $\proj(\cL_k(x))$, and therefore there
exist $w_1, \dots w_n$ with
$\proj(\cL_k(x)) = \bigcup_{i=1}^n \cU(x,w_i)$.
Let $\sigma_1(x,\delta') = \min_i \sigma_0(x,w_i,\delta')$. Then, for
all $x \in Z^c$ and all $w \in \proj(\cL_k(x))$,
\begin{displaymath}
\beta\left( \{ a \in B \st  d_Y(\cV^{(k)}_{> 1}(a,x),w) >
  \sigma_1(x,\delta') \} \right) > 1 - \delta'.  
\end{displaymath}
Let $E_N(\delta') = \{ x \in Z^c \st \sigma_1(x,\delta') > \frac{1}{N}\}$. 
Since $\bigcup_{N=1}^\infty E_N(\delta') = Z^c$ and 
$\nu(Z^c) = 1$, there exists $N = N(\delta,\delta')$ such that
$\nu(E_N(\delta')) > 1-\delta$. Let $\sigma = 1/N$ and let $E_{good} =
E_N$. 
\qed\medskip

\bold{Lyapunov subspaces and Relative Homology.}
The following Lemma is well known:
\begin{lemma}
\label{lemma:Lyapunov:exponents:relative}
The Lyapunov spectrum of the Kontsevich-Zorich cocycle acting on
relative homology is the Lyapunov spectrum of the Kontsevich-Zorich
cocycle acting on absolute homology, union $n$ zeroes, where $n = \dim
\ker
p$. 
\end{lemma}
Let $\index{$L$@$\bar{\cL_k}$}\bar{\cL_k} = p^{-1}(\cL_k) \subset H^1(M, \Sigma,\reals)$. We
have the Lyapunov flag
\begin{displaymath}
\{ 0 \} = \index{$V$@$\bar{\cV}_{\le j}^{(k)}$}\bar{\cV}_{\le 0}^{(k)} \subset \bar{\cV}_{\le 1}^{(k)}(b,x) \subset \dots
\bar{\cV}_{\le \bar{n}_k}^{(k)}(b,x) = \bar{\cL}_k(x), 
\end{displaymath}
corresponding to the action on the invariant subspace $\bar{\cL}_k
\subset H^1(M, \Sigma, \reals)$. 
Also for each $a \in B$, we have the opposite Lyapunov flag
\begin{displaymath}
\{ 0 \} = \bar{\cV}_{\ge \bar{n}_k}^{(k)} \subset \bar{\cV}_{\ge \bar{n}_k
  - 1}^{(k)}(a,x)
\subset \dots 
\bar{\cV}_{\ge 0}^{(k)}(a,x) = \bar{\cL}_k(x), 
\end{displaymath}

\begin{lemma}
\label{lemma:Lyapunov:relative:absolute}
Suppose $\hat{\lambda}_k \ne 0$. Then for almost all $(b,x)$, 
\begin{displaymath}
p(\bar{\cV}_{\le 1}^{(k)}(b,x)) = \cV_{\le 1}^{(k)}(b,x), 
\end{displaymath}
and $p$ is an isomorphism between these two subspaces. Similarly, for
almost all $(a,x)$, 
\begin{displaymath}
\bar{\cV}_{>1}^{(k)}(a,x) = p^{-1}(\cV_{>1}^{(k)}(a,x)). 
\end{displaymath}

\end{lemma}

\bold{Proof.} In view of Lemma~\ref{lemma:Lyapunov:exponents:relative}
and the assumption that $\hat{\lambda}_k \ne 0$, $\hat{\lambda}_k$ is
the top Lyapunov exponent on both $\cL_k$ and $\bar{\cL}_k$. 
Note that 
\begin{equation}
\label{eq:top:exponent:relative}
\bar{\cV}_{\le 1}^{(k)} = \{ \bar{v} \in \bar{\cL}_k \st \limsup_{t \to
  \infty} \frac{1}{t} \log \frac{\| T^n
  \bar{v}\|}{\|\bar{v}\|} \le -\hat{\lambda}_k \}.
\end{equation}
Also, 
\begin{equation}
\label{eq:top:exponent:absolute}
\cV_{\le 1}^{(k)} = \{ v \in {\cL}_k \st \limsup_{t \to
  \infty} \frac{1}{t} \log \frac{\| T^n
  v\|}{\|v\|} \le -\hat{\lambda}_k \}.
\end{equation}
It is clear from the definition of the Hodge norm on relative cohomology
(\ref{eq:df:relative:hodge})
that $\|p(v)\| \le C \|v\|$ for some absolute constant
$C$. Therefore, it follows from (\ref{eq:top:exponent:absolute}) and
(\ref{eq:top:exponent:relative}) that $p(\bar{\cV}_{\le 1}^{(k)}) \subset
\cV_{\le 1}^{(k)}$. But by Lemma~\ref{lemma:Lyapunov:exponents:relative},
$\dim(\bar{\cV}_{\le 1}^{(k)}) = \dim(\cV_{\le 1}^{(k)})$. Therefore,
$p(\bar{\cV}_{\le 1}^{(k)})
= \cV_{\le 1}^{(k)}$. 
\qed\medskip

\bold{Remark.} Even though we will not use this, a version of
Lemma~\ref{lemma:Lyapunov:relative:absolute} holds for all 
Lyapunov subspaces for non-zero exponents,
and not just the subspace corresponding to the top
Lyapunov exponent $\hat{\lambda}_k$. 

\bold{The action on $H^1(M,\Sigma,\cx)$.}
By the multiplicative ergodic theorem applied to the action of
$SL(2,\reals)$ on $\reals^2$,  for
$\beta$-almost all $b \in B$ there exists a one-dimensional subspace  
$\index{$W_+(b)$}W_+(b) \subset \reals^2$ such that $v \in W_+(b)$, 
\begin{displaymath}
\lim_{n \to \infty} \frac{1}{n} \log \| b_n^{-1} \dots b_0^{-1} v \| =
-\sigma_0. 
\end{displaymath}
Let
\begin{displaymath}
\index{$W^+(b,x)$}W^+(b,x) = (W_+(b) \tensor H^1(M,\Sigma,\reals))
\cap \area(x,1).
\end{displaymath}
Since we identify $\reals^2
\tensor H^1(M,\Sigma,\reals)$ with  $H^1(M,\Sigma,\cx)$, we may
consider $W^+(b,x)$ as a subspace of $H^1(M,\Sigma,\cx)$. This is the
``stable'' subspace for $T$. (Recall that $T$ moves into the past).

For a ``future trajectory'' $a \in B$, we can similarly define a
$1$-dimensional 
subspace $\index{$W_-(a)$}W_-(a) \subset \reals^2$ such that
\begin{displaymath}
\lim_{n \to \infty} \frac{1}{n} \log \| a_n \dots a_1 v \| =
-\sigma_0 \quad \text{ for $v \in W_-(a)$.} 
\end{displaymath}
Let $A: SL(2,\reals) \cross X_0 \to \operatorname{Hom}(H^1(M, \Sigma,\reals),
H^1(M,\Sigma,\reals))$ denote the Kontsevich-Zorich cocycle. We then
have the cocycle 
\begin{displaymath}
\index{$A$@$\hat{A}$}\hat{A}: SL(2,\reals) \cross X_0 \to \operatorname{Hom}(H^1(M, \Sigma,\cx),
H^1(M,\Sigma,\cx))
\end{displaymath}
given by 
\begin{displaymath}
\hat{A}(g,x)( v \tensor w) = g v \tensor A(g,x) w
\end{displaymath}
and we have made the identification $H^1(M, \Sigma,\cx) = \reals^2
\tensor H^1(M,\Sigma,\reals)$. This cocycle can be
  thought of  as the derivative cocycle for the action of $SL(2,\reals)$.
From the definition we see that the Lyapunov exponents of $\hat{A}$
are of the form $\pm \sigma_0 + \lambda_i$, where
the $\lambda_i$ are the Lyapunov exponents of $A$.

\section{Time changes and suspensions}

There is a 
natural ``forgetful'' map $f: \tilde{B} \to B$.
We extend functions on $B \cross X_0$ to $\tilde{B} \cross X_0$ by 
making them constant along the fibers of $f$. \mcc{(explain better)} 
The measure $\beta \cross \nu$ is a $T$-invariant measure on
  $\tilde{B} \cross X_0$.

\bold{The cocycles $\theta_j$.} By
Theorem~\ref{theorem:KZ:semisimple}, the restriction of the
Kontsevich-Zorich cocycle to each $\cL_j$ is semisimple.
Then by Theorem~\ref{theorem:semisimple:lyapunov}, the Lyapunov
spectrum of $T$ on each $\cL_j$ is semisimple, and the
  restriction of $T$ to the top Lyapunov subspace of each $\cL_j$
  consists of a single conformal block.
This means that 
there is a inner product \index{$\langle \,,\rangle_{j,b,x}$}$\langle \,,\rangle_{j,b,x}$ defined on
  $W_+(b) \tensor \cV_{\le 1}^{(j)}(b,x)$ 
and a function $\index{$\theta_j$}\theta_j: B \cross X_0 \to \reals$ 
such that for all $u,v \in W_+(b) \tensor \cV_{\le 1}^{(j)}(b,x)$, 
\begin{equation}
\label{eq:theta:j:def}
\langle \hat{A}(b_0^{-1},x) u, \hat{A}(b_0^{-1},x) v \rangle_{j,Tb,
  b_0^{-1} x} = e^{-\theta_j(b,x)} \langle u, v
\rangle_{j,b,x}. 
\end{equation}
\mcc{check signs}

To handle relative homology, we need to also consider the case in
which the action of $A(\cdot, \cdot)$ on a subbundle 
is trivial. We thus define an
inner product $\langle\,,\rangle_{0,b}$ on $\reals^2$, and a cocycle
$\theta_0: B \to \reals$ so that for $u, v \in W_+(b)$, 
\begin{equation}
\label{eq:theta:0:def}
\langle b_0^{-1} u, b_0^{-1} v \rangle_{0,Tb} = 
e^{-\theta_0(b)} \langle u, v
\rangle_{0,b}. 
\end{equation}
\mcc{check signs}
For notational simplicity, we let $\theta_0(b,x) = \theta_0(b)$. 


\bold{Switch to positive cocycles.} The cocycle
$\theta_j$ corresponds to the $\hat{A}(\cdot,
  \cdot)$-Lyapunov exponent $\sigma_0 + 
\hat{\lambda}_j$, where $\hat{\lambda}_j$ is the top Lyapunov
exponent of $A( \cdot, \cdot)$ in
$\cL_j$. Since $\sigma_0 > 0$ and $\hat{\lambda}_j \ge 0$, 
\begin{displaymath}
\sigma_0+\hat{\lambda}_j = \int_{B \cross X_0} \theta_j(b,x)
\, d\beta(b)\, 
d\nu(x) > 0. 
\end{displaymath}
Thus, the cocycle $\theta_j$ has positive average on $B \cross
X_0$. However, we do not know that $\theta_j$ is positive, i.e.\ that
for all $(b,x) \in B \cross X_0$, 
$\theta_j(b,x) > 0$. This makes it awkward to use
$\theta_j(b,x)$ to 
define a time change. Following \cite{Benoist:Quint} we use a positive
cocycle $\tau_j$ equivalent to $\theta_j$. 

By \cite[Lemma~2.1]{Benoist:Quint}, we can find a {\em positive} cocycle 
$\index{$\tau_j$}\tau_j: B \cross X_0 \to \reals$ and a measurable function
$\phi_j: B \cross X_0 \to \reals$ such that
\begin{equation}
\label{eq:def:phi:i}
\theta_j - \phi_j \circ T + \phi_j = \tau_j
\end{equation}
and
\begin{displaymath}
\int_{B \cross X_0} \tau_j(b,x) \, d\beta(b)\, d\nu(x) < \infty. 
\end{displaymath}
For $v \in W_+(b) \tensor \cV_{\le 1}^{(j)}(b,x)$ we define
\begin{equation}
\label{eq:def:tmp:prime:norm}
\|v\|'_{j,b,x} = e^{\phi_j(b,x)} \|v\|_{j,b,x},
\end{equation}
where the norm $\langle \cdot, \cdot \rangle_j$ is as in
(\ref{eq:theta:j:def}) and (\ref{eq:theta:0:def}). 
Then
\begin{equation}
\label{eq:hat:A:prime:norm}
\|\hat{A}(b_0^{-1},x) v \|'_{j, T(b,x)} = e^{-\tau_j(b,x)} \|
v\|'_{j,b,x}. 
\end{equation}

\bold{Suspension.}
Let \index{$B^X$}$B^X  = B \cross X_0 \cross (0,1]$. Recall that $\beta$ denotes the
measure on $B$ which is given by $\mu \cross \mu \cdots$.   
Let \index{$\beta^X$}$\beta^X$ denote the measure
on $B^X$ given by $\beta \cross \nu \cross dt$, where $dt$ is the
Lebesgue measure on 
$(0,1]$. In $B^X$ we identify $(b,x,0)$ with $(T(b,x),1)$, so that
$B^X$ is a suspension of $T$. We can then define a
suspension flow $\index{$T_t$}T_t: B^X \to 
B^X$ in the natural way. (Our suspensions are going downwards and not
upwards, since we think of $T$ as going into the past).
Then $T_t$ preserves the measure $\beta^X$. 

Let \index{$\tilde{B}^X$}$\tilde{B}^X  = \tilde{B} \cross X_0 \cross
(0,1]$. The suspension construction, the flow $T_t$, and the invariant
measure $\beta^X$ extend naturally from $B^X$ to $\tilde{B}^X$. 

Let \index{$T_t(b,x,s)_*$}$T_t(b,x,s)_*$ denote the action of $T_t(b,x,s)$
  on $H^1(M,\Sigma,\cx)$ (i.e.\ the derivative cocycle on the tangent
  space). Then, for $t \in \zed$ and $v \in W_+(b) \tensor
  \cV_{\le 1}^{(j)}(b,x)$ and $0 < s \le 1$ we have, in view of
  (\ref{eq:hat:A:prime:norm}), 
\begin{equation}
\label{eq:Tt:star:prime:norm:j}
\|T_t(b,x,s)_* v\|'_{j,T_t(b,x)} = e^{-\tau_j(t,b,x)} \|v\|'_{j,b,x},
\end{equation}
where $\tau_j(t,b,x) = \sum_{n=0}^{t-1} \tau_j(T^n(b,x))$. We can
extend the norm $\| \cdot \|_j'$ from $B \cross X_0$ to $B^X$ by \mcc{check}
\begin{displaymath}
\| v \|'_{j,b,x,s} = \| v \|'_{j,b,x} e^{-(1-s)\tau_j(b,x)}. 
\end{displaymath}
Then (\ref{eq:Tt:star:prime:norm:j}) 
holds for all $t \in \reals$ provided we set for $n \in
\zed$ and $0 \le s < 1$, \mcc{check}
\begin{displaymath}
\tau_j(n+s,b,x) = \tau_j(n,b,x) + s \tau_j(T^n(b,x)).
\end{displaymath}

\bold{The time change.} Here we differ slightly from
\cite{Benoist:Quint} since we would like to have several different
time changes of the flow $T_t$ on the same space. Hence, instead of changing
the roof function, we keep the roof function constant, but
change the speed in which one moves on the $[0,1]$
fibers. 

Let \index{$T^{\tau_j}_t$}$T^{\tau_j}_t: B^X \to B^X$ be the time change of $T_t$ where on $(b,x) \cross
[0,1]$ one moves at the speed $1/\tau_j(b,x)$. More precisely, we set
\begin{equation}
\label{eq:def:T:tau:j}
T^{\tau_j}_t(b,x,s) = (b,x, s  - t/\tau_j(b,x)), \quad \text{ if $0 < s
  - t/\tau_j(b,x) \le 1$,} 
\end{equation}
and extend using the identification $((b,x),0) = (T(b,x),1)$. \mcc{check}

Then $T^{\tau_k}_\ell$ is the operation of moving backwards in
time far enough so that the cocycle multiplies the 
direction of the top Lyapunov exponent in $\cL_k$ 
by $e^{-\ell}$. 
In fact, by   (\ref{eq:Tt:star:prime:norm:j}) and
  (\ref{eq:def:T:tau:j}),  we have, for $v \in W_+(b) \tensor
  \cV^{(k)}_{\le 1}(b,x)$, \mcc{check}
\begin{equation}
\label{eq:time:change:prime:norm}
\|T_\ell^{\tau_k}(b,x,s)_* v\|'_{j,T_\ell^{\tau_k}(b,x,s)} = e^{-\ell}
\|v\|'_{j,b,x,s}.
\end{equation}

\bold{The map $T^{\tau_k}$ and the two-sided shift space.}
On the space $\tilde{B}^{X}$, $T^{\tau_k}$ is invertible, and we denote
the inverse of $T^{\tau_k}_\ell$ by $T^{\tau_k}_{-\ell}$. We write
\begin{equation}
\label{eq:def:T:star}
T^{\tau_k}_{-\ell}(a \vee b,x,s)_* 
\end{equation}
for the linear map on the tangent space $H^1(M, \Sigma, \cx)$ induced by
$T^{\tau_k}_{-\ell}(a \vee b,x,s)$. 
In view of (\ref{eq:def:tmp:prime:norm}) and
(\ref{eq:time:change:prime:norm}),
we have for $v \in W_+(b) \tensor \cV_{\le 1}^{(k)}(b,x)$, 
\begin{equation}
\label{eq:T:tau:k:growth:general}
\|T^{\tau_k}_{-\ell}(a \vee b,x,s)_* v\| = \exp(\ell +
\phi_k(b,x,s) - \phi_k(T^{\tau_k}_{-\ell}(a \vee b,x,s))) \| v\|.
\end{equation}
Here we have omitted the subscripts on the norm $\| \cdot \|_{k,b,x}$
and also extended the function $\phi_k(b,x,s)$ so that for all $(b,x,s)
\in B^X$ and all $v \in W_+(b) \tensor \cV_{\le 1}^{(k)}(b,x)$,
\begin{displaymath}
\| v \|_{k,b,x} = e^{\phi_k(b,x,s)} \| v \|'_{k,b,x,s}.
\end{displaymath}
\mcc{more details here?}

\bold{Invariant measures for the time changed flows.}
Let \index{$\beta^{\tau_j,X}$}$\beta^{\tau_j,X}$ denote the measure on $B^X$ given by 
\begin{displaymath}
d\beta^{\tau_j,X}(b,x,t) = c_j \tau_j(b,x) \, d\beta(b)\, d\nu(x) \, dt,
\end{displaymath}
where the $c_j \in \reals$ is chosen so that $\beta^{\tau_j,X}(B^X) = 1$. 
Then the measures $\beta^{\tau_j,X}$ are invariant under the flows
$T^{\tau_j}_t$. We note the following trivial:
\begin{lemma}
\label{lemma:measures:time:changes}
The measures $\beta^{\tau_j,X}$ are all absolutely continuous with
respect to $\beta^X$. 
For every $\delta > 0$ there exists a compact subset $\cK = \cK(\delta)
\subset B^X$ and $L = L(\delta) < \infty$ such that for all $j$, 
\begin{displaymath}
\beta^{\tau_j,X}(\cK) > 1-\delta,
\end{displaymath}
and also for $(b,x,t) \in \cK$,  
\begin{displaymath}
\frac{d \beta^{\tau_j,X}}{d \beta^X}(b,x,t) \le L, \quad \frac{d
  \beta^X}{d \beta^{\tau_j,X}}(b,x,t) \le L. 
\end{displaymath}
\end{lemma}


\section{The martingale convergence argument}
\label{sec:martingale}


\bold{Standing Assumptions.}
Let
\begin{displaymath}
\index{$W^+[b,x]$}W^+[b,x] = \{ y \st y - x \in W^+(b,x). \}
\end{displaymath}
Then, $W^+[b,x]$ is the stable subspace for $T$.
From the definition, for almost all $b$, (locally) the sets $\{
W^+[b,x] \st x \in X\}$ form a 
measurable partition of $X_0$. Let
\begin{displaymath}
\index{$U^+(b,x)$}U^+(b,x) = W^+(b,x) \cap U_{\cx}(x),
\qquad \index{$U^+[b,x]$}U^+[b,x] = W^+[b,x] \cap U_{\cx}[x]. 
\end{displaymath}
We make the corresponding definitions for
$\index{$W^-(b,x)$}W^-(b,x)$,
\index{$W^-[b,x]$}$W^-[b,x]$,
\index{$U^+[b,x]$}$U^+[b,x]$ and \index{$U^-[b,x]$}$U^-[b,x]$.

It follows from Theorem~\ref{theorem:main:partone} applied to the
flow $r_\theta g_t r_{-\theta}$, using the fact that
$U_{\cx}[r_\theta x] = U_{\cx}[x]$, that for a.e.\ $x$, 
the conditional measures of $\nu$ along
$W^\pm[b,x]$ are supported on $U^{\pm}[b,x]$, and also that the
conditional measures of $\nu$ along $U^{\pm}[b,x]$ are Lebesgue.

\begin{lemma}
\label{lemma:exists:psi}
There exists a subset $\Psi \subset B^X$ with 
$\beta^X(\Psi) = 1$ such that for all $(b,x) \in \Psi$, 
\begin{displaymath}
\Psi \cap W^+[b,x] \cap \text{ ball of radius $1$ } \subset \Psi \cap
  U^+[b,x]. 
\end{displaymath}
\end{lemma}

\bold{Proof.} See \cite{Margulis:Tomanov:Ratner} or
\cite[6.23]{Einsiedler:Lindenstrauss:Pisa}. \mcc{find better reference?}
\qed\medskip

\bold{The parameter $\delta$.}
Let $\delta > 0$ be a parameter which will eventually be chosen
sufficiently small. We use the notation $c_i(\delta)$ and
$c_i'(\delta)$ for functions which tend to $0$ as $\delta \to 0$. In
this section we use the notation $A \approx B$ to mean that the ratio
$A/B$ is bounded between two positive constants depending on $\delta$. 

We first choose a compact subset $K_0 \subset \Psi \cap \Phi$ with
$\beta^X(K_0) 
> 1-\delta> 0.999$, the conull set $\Psi$ is as in
Lemma~\ref{lemma:exists:psi}, and the conull set $\Phi$ is as in
Theorem~\ref{theorem:no:yeti}.  By the multiplicative
  ergodic theorem and (\ref{eq:tmp:top:lyapunov}), we may also assume
  that there exists 
  $\ell_1(\delta) > 0$ such that for all $(b,x,s) \in K_0$
  all $k$ and all $v \in \cV_{\le 1}^{(k)}(b,x)$ and all $\ell > \ell_1(\delta)$,
\begin{equation}
\label{eq:tmp:K0:uniform:osceledets}
\|T_\ell(b,x,s)_* v\| \le e^{-(\lambda_k/2)\ell} \|v\|. 
\end{equation}
(Here, as in (\ref{eq:tmp:top:lyapunov}) 
the notation $T_\ell(b,x,s)_*$ denotes the action on
$H^1(M,\Sigma,\reals)$.) By the norm \index{abs@$"\norm{\cdot}$}$\| \cdot \|$ in
  this section, we mean the AGY norm (see
  \S\ref{sec:subsec:app:hodge}).
\begin{lemma}
\label{lemma:exists:K}
For every $\delta > 0$ there exists $K \subset B^X$ and $C = C(\delta)
<  \infty$, $\beta = \beta(\delta) > 0$  
and $C' = C'(\delta) < \infty$ such that 
\begin{itemize}
\item[(K1)] For all $L > C'(\delta)$, and all $(b,x,s) \in K$, 
\begin{displaymath}
\frac{1}{L}\int_{0}^L \chi_{K_0}(T_t(b,x,s)) \, dt \ge 0.99.
\end{displaymath}
\item[(K2)] $\beta^X(K) > 1 - c_1(\delta)$. Also, for 
all $j$, $\beta^{\tau_j,X}(K) > 1-c_1(\delta)$.
\item[(K3)] For all $j$ and all $(b,x,t) \in K$, $|\phi_j(b,x,t)| <
  C$, where $\phi_j$ is as in (\ref{eq:def:phi:i}). 
\item[(K4)] For all $j$, all $(b,x,t) \in K$ all $k\ne
    0$ and all $v \in
  \bar{V}_{\le 1}^{(k)}(b,x)$, 
\begin{equation}
\label{eq:tmp:pv:ge:beta:delta:v}
\|p(v)\| \ge \beta(\delta) \|v\|. 
\end{equation}
\item[(K5)] There exists $C_0 = C_0(\delta)$ such that for all $(b,x,s) \in
  K$ all $j$ and all $v \in W_+(b) \tensor \cV_{\le 1}^{(j)}(b,x)$, we have
$C_0^{-1} \| v \| \le \|v\|_{j,b,x} \le C_0 \|v\|$.
\end{itemize}
\end{lemma}

\bold{Proof.}  By the Birkhoff ergodic theorem, there exists $\cK''
\subset B^X$ such that $\beta^X(\cK'') > 1-\delta/5$ and (K1) holds for
$\cK''$ instead of $\cK$.  
We can choose $\cK' \subset B^X$ and $C = C(\delta) <
\infty$ such that $\beta^X(\cK') > 1-\delta/5$ and (K3) holds for
$\cK'$ instead of $K$. Let $\cK = \cK(\delta/5)$ and $L = L(\delta/5)$ 
be as in Lemma~\ref{lemma:measures:time:changes} with $\delta/5$
instead of $\delta$. Then choose $K_j \subset \Psi$ with
$\beta^{\tau_j,X}(K_j) > 1-\delta/(5 d L)$, where $d$ is the number of
Lyapunov exponents. In view of
  Lemma~\ref{lemma:Lyapunov:relative:absolute} there exists $\cK'''
  \subset X_0$ with $\beta^X(\cK''') > 1-\delta/5$ so that
  (\ref{eq:tmp:pv:ge:beta:delta:v}) holds. Similarly,
  there exists a set $\cK''''$ with $\cK'''' > 1-\delta/5$ where (K5) holds.
Then, let $K = \cK'''' \cap \cK''' \cap \cK'' \cap \cK' \cap \cK
\cap \bigcap_{j} 
K_j$. The properties (K1), (K2), (K3) and (K4) are easily verified.
\qed\medskip

\bold{Warning.} In the rest of this section, we will often identify
$K$ and $K_0$ with their pullbacks $f^{-1}(K) \subset \tilde{B}^X$ and
$f^{-1}(K_0) \subset \tilde{B}^X$ where $f:
\tilde{B}^X \to B^X$ is the forgetful map. \mcc{fix this later?}

\bold{The Martingale Convergence Theorem.}
Let $\cB^{\tau_j,X}$ denote the $\sigma$-algebra
of $\beta^{\tau_j,X}$ measurable functions on $B^X$.
As in \cite{Benoist:Quint}, let 
\begin{displaymath}
\index{$Q_\ell^{\tau_j,X}$}Q_\ell^{\tau_j,X} = (T^{\tau_j}_\ell)^{-1}(\cB^{\tau_j,X}). 
\end{displaymath}
(Thus if a function $F$ is measurable with respect to
$Q_\ell^{\tau_j,X}$, then $F$ depends only on what happened at least
$\ell$ time units in the past, where $\ell$ is measured using the time
change $\tau_j$.) 

Let 
\begin{displaymath}
Q_\infty^{\tau_j,X} = \bigcap_{\ell > 0} Q_\ell^{\tau_j,X}. 
\end{displaymath}
The $Q_\ell^{\tau_j,X}$ are a decreasing family of
$\sigma$-algebras, and then, 
by the Martingale Convergence Theorem, for $\beta^{\tau_j,X}$-almost
all $(b,x,s) \in B^X$, 
\begin{equation}
\label{eq:abstract:martingate:convergence}
\lim_{\ell \to \infty} 
\mathbb{E}_j( 1_{K} \mid Q_{\ell}^{\tau_j,X})(b,x,s)  =
\mathbb{E}_j( 1_{K} \mid Q_{\infty}^{\tau_j,X})(b,x,s)
\end{equation}
where $\mathbb{E}_j$ denotes expectation with respect to the measure
$\beta^{\tau_j,X}$. 

\bold{The set $S'$.} In view of
(\ref{eq:abstract:martingate:convergence}) and the 
condition (K2)  we can choose $S' = S'(\delta) \subset
B^{X}$ to be such that for all $\ell > \ell_0$, all $j$, and all 
$(b,x,s) \in S'$, 
\begin{equation}
\label{eq:Kprime:Sprime}
\mathbb{E}_j( 1_{K} \mid Q_{\ell}^{\tau_j,X})(b,x,s) > 1-c_2(\delta).
\end{equation}
By using Lemma~\ref{lemma:measures:time:changes} as in the
proof of Lemma~\ref{lemma:exists:K} we 
may assume that (by possibly making $\ell_0$ larger) we have  
for all $j$, 
\begin{equation}
\label{eq:Sprime:large:measure}
\beta^{\tau_j,X}(S') > 1-c_2(\delta). 
\end{equation}

\bold{The set $E_{good}$.} By Lemma~\ref{lemma:zero:one:law:star} we
may choose a subset $E_{good} \subset \tilde{B}^X$ (which is actually
of the form $\tilde{B} \cross E_{good}'$ for some subset $E_{good}'
\subset X \cross [0,1])$, with $\beta^X(E_{good}) >
1-c_3(\delta)$, and a number $\sigma(\delta) > 0$  
such that for any $(b,x,s) \in E_{good}$, any $j$ and 
any unit vector $w \in \cL_j(b,x)$, \mcc{check delta prime: think ok.}
\begin{equation}
\label{eq:random2:subspace}
\beta\left( \{ a \in B \st d_Y(w,\cV_{>1}^{(j)}(a,x)) > \sigma(\delta) \}
\right) > 1-c_3'(\delta).  
\end{equation}
We may assume that $E_{good} \subset K$. By the Osceledets
multiplicative ergodic theorem and
Lemma~\ref{lemma:Lyapunov:relative:absolute}, 
we may also assume that there exists
$\alpha > 0$ (depending only on the Lyapunov spectrum), and
$\ell_0 = \ell_0(\delta)$ such that for $(b,x,s) \in E_{good}$, 
$\ell > \ell_0$, at least $1-c_3''(\delta)$ 
measure of $a \in B$, and all $\bar{v} \in
\bar{\cV}^{(j)}_{>1}(a,x)$, 
\begin{equation}
\label{eq:upper:bound:hat:V:j}
\| T^{\tau_j}_{-\ell}(a \vee b, x,s)_* \bar{v} \| \le e^{(1-\alpha) \ell}
\|\bar{v}\|. 
\end{equation}
\mcc{explain that the endpoint has to also be in a good set? I think
  not needed}

\bold{The sets $\Omega_\rho$.} In view of
(\ref{eq:Sprime:large:measure}) and
the Birkhoff ergodic theorem, for every $\rho > 0$ there exists a
set $\index{$\Omega_\rho$}\Omega_\rho = \Omega_\rho(\delta) \subset \tilde{B}^X$ such that
\begin{itemize}
\item[{\rm ($\Omega$1)}] $\beta^{X}(\Omega_\rho) > 1 - \rho$. 
\item[{\rm ($\Omega$2)}] There exists $\ell_0' = \ell_0'(\rho)$ 
such that for all $\ell > \ell_0'$,
  and all $(b,x,s) \in \Omega_\rho$, 
$$|\{ t \in [-\ell,\ell] \st T_t(b,x,s) \in S' \cap E_{good} \}| \ge
(1-c_5(\delta)) 2\ell.$$  
\end{itemize}

\begin{lemma}
\label{lemma:can:choose:in:U:perp}
Suppose the measure $\nu$ is not affine. Then there exists $\rho > 0$
so that for every $\delta' > 0$ there exist $(b,x,s) \in \Omega_\rho$,
$(b,y,s) \in \Omega_\rho$ 
with $\|y - x\| < \delta'$ such that $p(y - x) \in
p(U)^{\perp}_{\cx}(x)$,
\begin{equation}
\label{eq:y:minus:x:transverse:U}
d(y-x,U_{\cx}(x)) > \frac{1}{10} \|y-x\|
\end{equation}
and
\begin{equation}
\label{eq:def:general:position}
d(y-x,W^+(b,x)) > \frac{1}{3} \|y-x\|
\end{equation}
(so $y-x$ is in general position with respect to $W^+(b,x)$.) 
\end{lemma}

\bold{Remark.} In view of Theorem~\ref{theorem:no:yeti}, it follows
that for $(b,x,s)$, $(b,y,s)$ satisfying the conditions of
Lemma~\ref{lemma:can:choose:in:U:perp}, 
$p(y-x)$ is orthogonal to the complexification $F_\cx(x)$ of
the Forni subspace $F(x)$.

\bold{Proof.} 
By Fubini's theorem, there exists a subset $\Omega_\rho'
\subset X$ with $\nu(\Omega_\rho') \ge 1-\rho^{1/2}$ such that 
for $x \in \Omega_\rho'$, 
\begin{equation}
\label{eq:def:Omega:rho:prime}
(\beta \cross dt)( \{ (b,s) \st (b,x,s) \in \Omega_\rho \}) \ge
(1-\rho^{1/2}). 
\end{equation}
Let $\cK$ be an arbitrary compact subset of $X_0$ with $\nu(\cK) > 1/2$,  
and let $\tilde{\cK}$
denote its lift to $\tilde{X}_0$. Let $\pi: \tilde{X}_0 \to X_0$ denote the
natural map.  We have
\begin{equation}
\label{eq:measure:Omega:rho:prime:lower}
\nu(\Omega_\rho') \ge (1-2 \rho^{1/2}) \nu(\cK).
\end{equation}

In view of Lemma~\ref{lemma:volume:form:on:U} we can find finitely many 
sets $J_\alpha \subset K_\alpha
\subset \tilde{X}_0$ and
constants $N > 0$ and $\delta_0 > 0$
such
that the following hold:
\begin{itemize}
\item[(i)] For all $\alpha$, $K_\alpha$ is diffeomorphic to an open ball,
  and the restriction of $\pi$ to $K_\alpha$ is injective. 
\item[(ii)] The sets $J_\alpha$ are disjoint, and up to a null set
$\pi(\tilde{\cK}) = \bigsqcup_\alpha \pi(J_\alpha)$.
\item[(iii)] Any point belongs to at most $N$ of the sets $\pi(K_\alpha)$. 
\item[(iv)] Recall that for $x \in \tilde{X}_0$, $U_\cx[x]$ denotes the
  (infinite) affine space whose tangent space is $U_\cx(x)$. We have,
for $\nu$-almost all $x \in J_\alpha$, 
\begin{equation}
\label{eq:lower:bound:lebesgue:Jalpha}
\Vol( U_\cx[x] \cap K_\alpha) \ge \delta_0,
\end{equation}
where $\Vol( \cdot )$ is as in Lemma~\ref{lemma:volume:form:on:U}.
\end{itemize}
Let 
\begin{equation}
\label{eq:def:Omega:rho:twoprime}
\Omega_\rho'' = \{ x \in J_\alpha \st \nu_{U_{\cx}(x)}(\Omega_\rho'
\cap K_\alpha) \ge (1-\rho^{1/4}) \nu_{U_{\cx}(x)}(K_\alpha) \}. 
\end{equation}
In the above equation, $\nu_{U_{\cx}(x)}$ is the conditional measure of
$\nu$ along $U_\cx[x]$ (which is in fact a multiple of the measure
$\Vol$ of Lemma~\ref{lemma:volume:form:on:U}). 
By (\ref{eq:measure:Omega:rho:prime:lower}), properties (ii), (iii) and 
Fubini's theorem, 
$\nu(\Omega_\rho'') \ge (1-2 N \rho^{1/4})\nu(\cK)$. In particular,
$\bigcup_{\rho > 0} \Omega_\rho''$ is conull in $\cK$. 

Note that by the definition of $\Omega_\rho''$, if $x \in
\Omega_\rho'' \cap J_\alpha$ then $U_{\cx}[x] \cap J_\alpha \subset
\Omega_\rho''$. It follows that we may write, for some indexing set
$I_\alpha(\rho)$, 
\begin{displaymath}
\Omega_\rho'' \cap J_\alpha = \bigsqcup_{x \in I_\alpha(\rho)} U_{\cx}[x]
\cap J_\alpha. 
\end{displaymath}
Suppose that for all $\alpha$ and all $\rho > 0$, $I_\alpha(\rho)$ 
is countable. Then, for a positive measure set of $x \in \tilde{X}_0$, 
$x$ has an open
neighborhood in $U_\cx[x]$ whose $\nu$-measure is positive. 
Then by ergodicity of the geodesic flow, this holds for
$\nu$-almost all $x \in \tilde{X}_0$ and without loss of
  generality, for all $x \in I_\alpha(\rho)$.

The restriction of $\nu$ to $U_{\cx}[x]$
is a multiple of the measure $\Vol$ of
Lemma~\ref{lemma:volume:form:on:U}, 
therefore there exists a constant
$\psi(x) \ne 0$ such that for $E \subset U_{\cx}[x]$, $\nu(E) = \psi(x)
\Vol(E)$. Since both $\nu$ and $\Vol$ are
invariant under the $SL(2,\reals)$ action, $\psi(x)$ is invariant, and
thus by ergodicity $\psi$ is constant almost everywhere. 

Let $I_\alpha' = \bigcup_{\rho > 0} I_\alpha(\rho)$. For $x, y \in
I_\alpha'$ write $x \sim y$ if $U_\cx[x] \cap J_\alpha = U_\cx[y]
\cap J_\alpha$, 
and let
$I_\alpha'' \subset I_\alpha'$ be the subset where we keep only one
member of each $\sim$-equivalence class. Note that by properties (i)
and (iv), 
for distinct $x, y \in I_\alpha''$, $U_\cx[x] \cap K_\alpha$ and
$U_\cx[y] \cap K_\alpha$ are disjoint up to a set of measure
$0$. \mcc{justify?} 
Then
(\ref{eq:lower:bound:lebesgue:Jalpha})
implies that for each $\alpha$,
\begin{displaymath}
\nu(K_\alpha) \ge \sum_{x \in I_\alpha''} \nu( U_{\cx}[x] \cap
K_\alpha) = \sum_{x \in I_\alpha''} \psi \Vol( U_{\cx}[x] \cap
K_\alpha) \ge \psi \delta_0 |I_\alpha''|.
\end{displaymath}
where $| \cdot |$
denotes the cardinality of a set. 
Since $\nu$ is a finite measure, we get that each $I_\alpha''$ is
finite. 
Since for a fixed $\cK$, there are only finitely many sets
$K_\alpha$, this implies that the support of 
restriction of $\nu$ to $\cK$  is contained in a finite union of ``affine pieces'' each of the form
$U_\cx[x_j] \cap K_\alpha$ for some $x_j \in \cK$, and the measure $\nu$
restricted to each affine piece coincides with $\psi \Vol$. 
It follows from the
ergodicity of $g_t$ that the affine pieces fit together to form an
(immersed) submanifold. Thus, $\nu$ is affine. \mccc{(more details here?)}

Thus, we may assume that there exist $\alpha$ and $\rho > 0$ such
that $I_\alpha(\rho)$ is not countable. Then we can find $x_1 \in
I_\alpha(\rho)$ and $y_n \in I_\alpha(\rho)$ \mcc{explain?} such that 
\begin{displaymath}
\lim_{n \to \infty} hd(U_{\cx}[x_1] \cap K_\alpha, U_{\cx}[y_n] \cap
K_\alpha)  = 0,  
\end{displaymath}
where $hd$ denotes Hausdorff distance between sets, (using the
distance $d^{X_0}$ defined in \S\ref{sec:semi:markov}).  
Let $f_n:
p(U)_{\cx}[y_n] \to p(U)_{\cx}[x_1]$ denote the function taking $z \in
p(U)_{\cx}[y_n]$ to the unique point in $p(U)_{\cx}[x_1] \cap
p(U)_{\cx}^\perp[z]$. Then, for large $n$, the map $f_n$ is almost
measure preserving, in the sense that for $V \subset p(U)_{\cx}(y_n)$, 
\begin{displaymath}
(0.5) |V| \le |f_n(V)| \le 2 |V|,
\end{displaymath}
where $| \cdot |$ denotes Lebesgue measure. 
Then, in view of
the definition (\ref{eq:def:Omega:rho:twoprime}) of $\Omega_\rho''$, 
for sufficiently large $n$,
there exist $x \in U_{\cx}[x_1] \cap \Omega_\rho'$ and
$y \in U_{\cx}[y_n] \cap \Omega_\rho'$ such that $p(y - x) \in
p(U)^\perp_{\cx}(x)$, and $\|y-x\| < \delta'$. 
Then, by the definition (\ref{eq:def:Omega:rho:prime}) of
$\Omega_\rho'$, we can choose $(b,s)$ so that $(b,x,s) \in
\Omega_\rho$, $(b,y,s) \in \Omega_\rho$, and
(\ref{eq:y:minus:x:transverse:U}) and 
(\ref{eq:def:general:position}) holds. 
\qed\medskip

\bold{Standing Assumption.} We fix $\rho = \rho(\delta)$ so that
Lemma~\ref{lemma:can:choose:in:U:perp} holds. 

The main part of the proof is the following: 
\begin{proposition}
\label{prop:get:more:points:in:Wplus}
There exists $C(\delta) > 1$ such that the following holds:
Suppose for every $\delta' > 0$ there exist 
$(b,x,s), (b,y,s) \in \Omega_\rho$ with  $\|x-y\| \le
\delta'$, $p(x-y) \in p(U)^{\perp}_{\cx}(x)$, and so that
(\ref{eq:y:minus:x:transverse:U}) and (\ref{eq:def:general:position}) hold. 
Then for every $\epsilon > 0$ there exist 
$(b'',x'',s'') \in K_0$, $(b'',y'',s'') \in K_0$, 
such that $y'' - x'' \in U^{\perp}_{\cx}(x'')$, 
\begin{displaymath}
\frac{\epsilon}{C(\delta)}\le \|y'' -x''\| \le C(\delta) \epsilon,
\end{displaymath}
\begin{equation}
\label{eq:y2prime:minus:x2prime:far:U}
d(y''-x'',U_{\cx}(x'')) \ge \frac{1}{C(\delta)} \|y''-x''\|,
\end{equation}
\begin{equation}
\label{eq:y2prime:minus:x2prime:close:Wplus}
d(y''-x'',W^+(b'',x'')) < \delta'',
\end{equation}
where $\delta''$ depends only on $\delta'$, and $\delta'' \to 0$ as
$\delta' \to 0$.
\end{proposition}

\makefig{Proof of
  Proposition~\ref{prop:get:more:points:in:Wplus}. In the figure,
  going ``up'' corresponds to the ``future''. The map $T_m$ for $m >
  0$ takes one   $m$ steps into the
  ``past''.}{fig:martingale}{\includegraphics{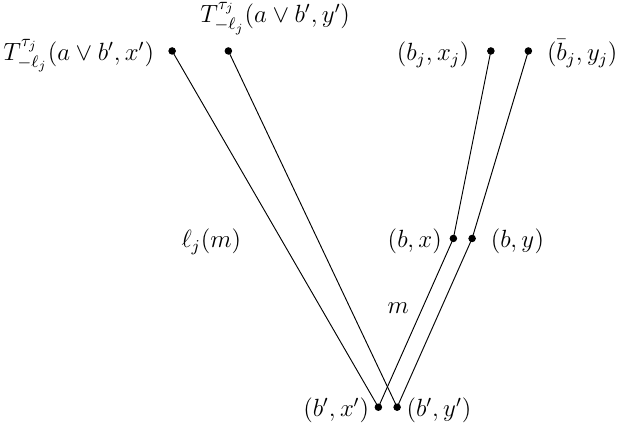}}  

\bold{Proof.} Let $\tilde{\Lambda} \subset \hat{\Lambda}$ denote the
subset $\{ k \st \hat{\lambda}_k \ne 0\}$. We may decompose
\begin{equation}
\label{eq:decomp:Uperp:corrected}
p(U)^{\perp}(x) = \bigoplus_{k \in \tilde{\Lambda}} \cL_k(x) \bigoplus F(x)
\end{equation}
as in \S\ref{sec:random:walks}. For $j \in \tilde{\Lambda}$, 
let $\pi_j$ denote the projection to $\cL_j$, using the decomposition
(\ref{eq:decomp:Uperp:corrected}). Note that by Theorem~\ref{theorem:no:yeti},
the projection of $p(y-x)$ to $F(x)$ is always $0$. 

For $m \in \reals^+$, write
\begin{displaymath}
(b',x',s') = T_m(b,x,s), \quad (b',y',s') =
T_m(b,y,s), 
\end{displaymath}
and let
\begin{displaymath}
w_j(m) = \pi_j(x'-y').
\end{displaymath}
(We will always have $m$ small enough so that the above equation makes
sense). \mcc{say this better}
Let $\ell_j(m)$ be such that
\begin{displaymath}
e^{\ell_j(m)} \| w_j(m) \| = \epsilon. 
\end{displaymath}
We also need to handle the relative homology part (where the action of
the Kontsevich-Zorich cocycle is trivial). Set $\ell_0(m)$ to be the
number such that \mcc{check} 
\begin{displaymath}
e^{\ell_0(m)} \|x'-y'\| = \epsilon.   
\end{displaymath}
Choose $0 < \sigma' \ll \lambda_{min}$ where
$0 < \lambda_{min} = \min_{j \in \tilde{\Lambda} } \hat{\lambda}_j$. 
\mcc{specify exact constants} 
We will be choosing $m$ 
so that 
\begin{equation}
\label{eq:range:m}
\frac{\sigma'}{2} | \log \|y-x\| |  \le  m \le \sigma' |\log \|y-x\||.
\end{equation}
In view of (\ref{eq:def:general:position}) and
Theorem~\ref{theorem:forni}, 
(after some uniformly bounded time), \mcc{fixme:explain} 
$\|w_j(m)\|$ is an increasing
function of $m$ (since the factor of $e^{-t}$ from the geodesic flow beats the
contribution of the Kontsevich-Zorich cocycle). 
Therefore, $\ell_j(m)$ is a decreasing function of $m$.

For a bi-infinite sequence $b \in \tilde{B}$ and $x \in X_0$, let
\begin{displaymath}
  G_j(b,x,s) = \{ m \in \reals_+ \st
  T^{\tau_j}_{-\ell_j(m)}T_m(b,x,s) \in S' \}.
\end{displaymath}
Let $G_{all}(b,x,s) = \bigcap_j G_j(b,x,s) \cap \{ m \st T_m(b,x,s) \in
E_{good} \}$.  
\begin{lemma}
\label{lemma:most:points:good:for:all}
For $(b,x,s) \in \Omega_\rho$ and $N$ sufficiently large, 
\begin{displaymath}
\frac{|G_{all}(b,x,s) \cap [0,N]|}{N} \ge 1-c_6(\delta).
\end{displaymath}
\end{lemma}

\bold{Proof.} We can write $T^{\tau_j}_{-\ell_j(m)}T_m =
T_{-g_j(m)}$. By definition, 
\begin{displaymath}
m \in G_j(b,x,s) \quad \text{ if and only if } \quad
T_{-g_j(m)}(b,x,s) \in S'. 
\end{displaymath}
Since $\ell_j(m)$ is a decreasing function of $m$,
so is $g_j$, and in fact, 
for all $m_2 > m_1$
\begin{displaymath}
g_j(m_1) - g_j(m_2) > m_2 - m_1.   
\end{displaymath}
This implies that 
\begin{equation}
\label{eq:gj:inverse}
g_j^{-1}(m_1) - g_j^{-1}(m_2) < m_1 - m_2. 
\end{equation}
Let $F = \{ t \in [0, g_j(N)] \st T_{-t}(b,x) \not\in S'
\}$. 
By  condition ($\Omega$2), for $N$ large enough, 
$|F| \le (1-c_5(\delta)) g_j(N)$. Note that $G_j^c \cap [0,N] =
g_j^{-1}(F)$.  Then, by (\ref{eq:gj:inverse}), 
\begin{displaymath}
|G_j^c \cap [0,N]| = |g_j^{-1}(F) | \le |F| \le c_5(\delta)
g_j(N) \le c_6(\delta) N,
\end{displaymath}
where as in our convention $c_6(\delta) \to 0$ as $\delta \to 0$. 
\qed\medskip

We now continue the proof of
Proposition~\ref{prop:get:more:points:in:Wplus}. We may assume that
$\delta'$ is small enough so that the right-hand-side of
(\ref{eq:range:m}) is smaller then the $N$ of
Lemma~\ref{lemma:most:points:good:for:all}.  
Suppose $(b,x,s) \in
\Omega_\rho$, $(b,y,s) \in \Omega_\rho$.
By Lemma~\ref{lemma:most:points:good:for:all}, we can  
fix $m \in G_{all}(x)$ such that (\ref{eq:range:m}) holds. 
Write $\ell_j = \ell_j(m)$. Let 
\begin{displaymath}
(b',x',s') = T_m(b,x,s), \quad (b',y',s') =
T_m(b,y,s). 
\end{displaymath}
For $j \in \tilde{\Lambda}$, let 
\begin{displaymath}
(b_j,x_j,s_j) = T^{\tau_j}_{-\ell_j(m)}(b',x',s'), \quad
(\bar{b}_j,y_j,\bar{s}_j) = T^{\tau_j}_{-\ell_j(m)}(b',y',s'). 
\end{displaymath}
Since $m \in G_{all}(b,x,s)$, we have $(b_j,x_j,s_j) \in
S'$, $(\bar{b}_j,y_j,\bar{s}_j) \in S'$. Then, 
by (\ref{eq:Kprime:Sprime}), for all $j$, 
\begin{displaymath}
{\mathbb E}_j( 1_{K} \mid Q_{\ell_j}^{\tau_j,X})(b_j,x_j,s_j) > (1 -
c_2(\delta)), 
\quad {\mathbb E}_j( 1_{K} \mid Q_{\ell_j}^{\tau_j,X})(\bar{b}_j,y_j,
\bar{s}_j) > (1 - c_2(\delta)). 
\end{displaymath}
Since $T^{\tau_j}_{\ell_j}(b_j,x_j,s_j) = (b',x',s')$, 
by \cite[(7.5)]{Benoist:Quint} \mcc{this may need more
  explanation. Not sure what to write here}
we have
\begin{displaymath}
{\mathbb E}_j( 1_{K} \mid Q_{\ell_j}^{\tau_j,X})(b_j,x_j,s_j) = \int_{B}
1_{K}(T^{\tau_j}_{-\ell_j}(a \vee b',x',s')) \, d\beta(a), 
\end{displaymath}
where the notation $a \vee b'$ is as in (\ref{eq:def:vee}). 
Thus, for all $j \in \tilde{\Lambda}$, 
\begin{equation}
\label{eq:T:tau:j:inverse:ell:most:in:K}
\beta\left(\{ a \st T^{\tau_j}_{-\ell_j}(a \vee b', x',s') \in K \}\right) >
1-c_2(\delta). 
\end{equation}
Similarly, for all $j \in \tilde{\Lambda}$, 
\begin{displaymath}
\beta\left(\{ a \st T^{\tau_j}_{-\ell_j}(a \vee b', y',s') \in K \}\right) >
1-c_2(\delta). 
\end{displaymath}

Let $w = x'-y'$, and let $w_j = \pi_j(w)$.  We can write
\begin{equation}
\label{eq:w:bar:wj}
w = \bar{w}_0 + \sum_{j \in \hat{\Lambda}} \bar{w}_j 
\end{equation}
where $\bar{w}_0 \in \ker p$, and for $j > 0$, $\bar{w}_j$ are chosen so that
$\pi_j(\bar{w}_j) = w_j$, and also $\| \bar{w}_j\| \approx \|w_j\|$.

For any $a \in B$, we may write
\begin{displaymath}
w_j = \xi_j(a)+v_j(a),
\end{displaymath}
where $\xi_j(a) \in W_+(b') \tensor \cV_{\le 1}^{(j)}(b',x')$, and 
\begin{displaymath}
v_j(a) \in W_+(b) \tensor \cV_{>1}^{(j)}(a,x') + W_-(a) \tensor
\cL_j(b',x'). 
\end{displaymath}
This decomposition is motivated as follows: if we consider the Lyapunov
decomposition 
\begin{displaymath}
\cx \tensor \cL_j(x) = \bigoplus_{k} \cV_k(a \vee b,x)
\end{displaymath}
then $\xi_j(a)$ belongs to the subspace $\cV_{\le 1}(a \vee b,x)$
corresponding to the top Lyapunov exponent $\sigma_0 +
\hat{\lambda}_j$ for the action of $T_{-t}$, and $v_j \in \oplus_{k
  \ge 2} \cV_k(a \vee b,x)$ will grow with a smaller Lyapunov exponent
under $T_{-t}$. Then $v_j(a)$ will also grow with a smaller Lyapunov
exponent then $\xi_j(a)$ under $T^{\tau_j}_{-\ell}$. 

Since $m \in G_{all}(b,x,s)$, we have $(b',x',s') \in
E_{good}$. Then, by (\ref{eq:random2:subspace}), for 
at least $1-c_3'(\delta)$ fraction of $a \in B$, 
\begin{equation}
\label{eq:random:projection}
\|v_j(a)\| \approx \|\xi_j(a)\| \approx \|w_j\| \approx
\epsilon e^{-\ell_j},   
\end{equation}
where the notation $A \approx B$ means that $A/B$ is bounded between
two constants depending only on $\delta$. 
Since $(b',x',s') \in E_{good} \subset K$, by condition (K3) we have
$|\phi_j(b',x',s')| \le C(\delta)$. Also by
(\ref{eq:T:tau:j:inverse:ell:most:in:K}), for at least
$1-c_2(\delta)$ fraction of $a \in B$, we
have $T^{\tau_j}_{-\ell_j}(a \vee b', x',s') \in K$, so again by
condition (K3) we have
\begin{displaymath}
|\phi_j(T^{\tau_j}_{-\ell_j}(a \vee b', x',s'))| \le C(\delta). 
\end{displaymath}
Thus, by (\ref{eq:random:projection}), 
(\ref{eq:T:tau:k:growth:general}) and
(\ref{eq:upper:bound:hat:V:j}),  we have, for all $j \in
\tilde{\Lambda}$, and at least $1-c_4(\delta)$ fraction of $a \in B$, 
\begin{equation}
\label{eq:estimate:winning:index}
\|T^{\tau_j}_{-\ell_j}(a \vee b',x',s')_* \xi_j(a) \| \approx
\epsilon, \text{ and } 
\| T^{\tau_j}_{-\ell_j}(a \vee b',x',s')_* v_j(a) \| = O( e^{-\alpha \ell_j}), 
\end{equation}
where $\alpha > 0$ depends only on the Lyapunov spectrum. 
(The notation in (\ref{eq:estimate:winning:index}) is defined in
(\ref{eq:def:T:star})). Hence, for at least $1-c_4(\delta)$ fraction
of $a \in B$,  
\begin{displaymath}
\|T^{\tau_j}_{-\ell_j}(a \vee b',x',s')_* w_j \| \approx \epsilon.
\end{displaymath}
Since $\lambda_j \ge 0$ (and by Theorem~\ref{theorem:no:yeti}, 
if $\lambda_j = 0$ then $j=0$, and $\bar{w}_0 \in \ker p$ where the
action of the Kontsevich-Zorich cocycle is trivial), we have for at least
$1-c_4(\delta)$ fraction of  $a \in B$, 
\begin{equation}
\label{eq:estimate:T:minus:tau:j:bar:w:j}
\|T^{\tau_j}_{-\ell_j}(a \vee b',x',s')_* \bar{w}_j \| \approx \epsilon.
\end{equation}

Let 
\begin{displaymath}
t_j(a) = \sup\{t > 0 \st \|T_{-t}(a \vee b',x',s')_* \bar{w}_j\| \le
\epsilon \},  
\end{displaymath}
and let $j(a)$ denote a $j \in \tilde{\Lambda}$
such that $t_j(a)$ is as small as possible as $j$ varies over
$\tilde{\Lambda}$. Then, if $j = j(a)$, then by
(\ref{eq:estimate:T:minus:tau:j:bar:w:j}), 
\begin{equation}
\label{eq:clarify:etimate:winning:index}
\|T_{-t_j(a)}(a \vee b',x',s')_* \bar{w}_j \| \approx 
\|T^{\tau_j}_{-\ell_j}(a \vee b',x',s')_* \bar{w}_j \| \approx
\epsilon. 
\end{equation}
Also, for at least $1-c_4(\delta)$-fraction of $a \in B$, 
if $j = j(a)$ and $k \ne j$, then by
(\ref{eq:estimate:T:minus:tau:j:bar:w:j}),  
\begin{equation}
\label{eq:estimate:losing:index}
\| T^{\tau_j}_{-\ell_j}(a \vee b',x',s')_* \bar{w}_k \| \le
C_1(\delta)\epsilon, 
\end{equation}
where $C_1(\delta)$ depends only on $\delta$. Therefore, by (\ref{eq:w:bar:wj}),
(\ref{eq:clarify:etimate:winning:index}), and 
(\ref{eq:estimate:losing:index}), for at least
$1-c_4(\delta)$-fraction of $a \in B$,  if $j = j(a)$, 
\begin{equation}
\label{eq:xprime:matching:yprime}
\| T^{\tau_j}_{-\ell_j}(a \vee b',x',s')_* (y'-x') \| \approx \epsilon. 
\end{equation}
We now choose $\delta > 0$ so that $c_4(\delta) + 2c_2(\delta) < 1/2$,
and using (\ref{eq:T:tau:j:inverse:ell:most:in:K}) we 
choose $a \in B$ so that (\ref{eq:xprime:matching:yprime}) holds,
and also
\begin{displaymath}
T^{\tau_j}_{-\ell_j}(a \vee b',x',s') \in K, \qquad
T^{\tau_j}_{-\ell_j}(a \vee b',y',s') \in K. 
\end{displaymath}
We may write
\begin{displaymath}
T^{\tau_j}_{-\ell_j}(a \vee b',x',s') = T_{-t}(a \vee b, x', s'), 
\qquad T^{\tau_j}_{-\ell_j}(a \vee b',y',s') = T_{-t'}(a \vee b, y', s')
\end{displaymath}
Then, $|t' - t|  \le C(\delta)$. Therefore by condition (K1), there exists $t''$
with $|t''-t|  \le C(\delta)$ such that
\begin{displaymath}
(b'',x'',s'') = T_{-t''}(a \vee b',x',s') \in K_0, \qquad
(b'',y'',s'') = T_{-t''}(a \vee b',y',s') \in K_0. 
\end{displaymath}

Since $\|w \| \approx \epsilon e^{-\ell_j}$, and $\ell_j \to \infty$
as $\delta' \to 0$, we have $\|w\| = \|x'-y'\| \to 0$ as $\delta' \to 0$. 
Since $T_{-t''}$ does not expand the $W^-$ components, the $W^-$ component
of $x''-y''$ is bounded by the $W^-$ component of $x'-y'$. Thus, the
size of the $W^-$ component of $x''-y''$ tends to $0$ as $\delta' \to
0$. Thus (\ref{eq:y2prime:minus:x2prime:close:Wplus}) holds. 

It remains to prove (\ref{eq:y2prime:minus:x2prime:far:U}). If
\mcc{check constants}
\begin{equation}
\label{eq:p:y2prime:x2prime:big}
\|p(y''-x'')\| \ge \frac{1}{C(\delta)} \|y''-x''\|
\end{equation}
then (\ref{eq:y2prime:minus:x2prime:far:U}) holds since $p(y''-x'')
\in p(U)^\perp(x'')$. This automatically holds for the case where
$|\Sigma| =1$ (and thus, in particular, 
there are no marked points). If not, we may write
\begin{displaymath}
y''-x'' = w_+'' + \bar{w}_0''
\end{displaymath}
where $\|w_+'' \| \le c(\delta) \|\bar{w}_0''\|$ and $\bar{w}_0'' \in
\ker p$. We will need to rule out the case where $\bar{w}_0''$ is very
close to $U^+(x'') \cap \ker p$. 
We will show that this contradicts the assumption
(\ref{eq:y:minus:x:transverse:U}). 

Let $w_+'$, $\bar{w}_0'$ be such that 
\begin{displaymath}
w_+'' = T_{-t''}(a \vee b, x', s')_* w_+', \quad \bar{w}_0'' = T_{-t''}(a \vee b, x', s')_* \bar{w}_0'.
\end{displaymath}
Then $y'-x' = w_+' + \bar{w}_0'$ and in view of
(\ref{eq:tmp:K0:uniform:osceledets}) and
(\ref{eq:random:projection}), \mcc{explain more?}
\begin{displaymath}
\|w_+'\| \le e^{-\lambda_{min}t''/2} \|\bar{w}_0'\| \approx
e^{-\lambda_{min}t''/2} \|y'-x'\|. 
\end{displaymath}
Applying $T_{-m}(b,x',s')$ to both sides we get
\begin{displaymath}
y - x = w_+ + \bar{w}_0,
\end{displaymath}
where $\bar{w}_0 \in \ker p$, and 
\begin{displaymath}
\|w_+\| \le e^{2 m} \|w_+'\| \le e^{2m - \frac{\lambda_{min}t''}{2}}
  \|x-y\|. 
\end{displaymath}
By (\ref{eq:range:m}), $2m - \frac{\lambda_{min}t''}{2} \le -
\frac{\lambda_{min}t''}{4}$. \mcc{make constants explicit}
Thus, $\|w_+\| \le (1/100) \|y-x\|$. Therefore, by
(\ref{eq:y:minus:x:transverse:U}), we have
\begin{displaymath}
d(\bar{w_0}, \ker p \cap U_{\cx}(x)) > \frac{1}{20} \|w_0\|. 
\end{displaymath}
Since the action of the cocycle on $\ker p$ is trivial (and we have
shown that in our situation the component in $\ker p$ dominates
throughout the process), this implies
\begin{displaymath}
d(\bar{w}_0'', \ker p \cap U_{\cx}(x'')) > \frac{1}{20} \|w_0''\| \ge
\frac{1}{40} \|y''-x''\|.
\end{displaymath}
This, together with the assumption that
(\ref{eq:p:y2prime:x2prime:big}) does not hold, 
implies (\ref{eq:y2prime:minus:x2prime:far:U}).
\qed\medskip

\bold{Proof of Theorem~\ref{theorem:P:measures}.} It was already
proved in Theorem~\ref{theorem:main:partone} that $\nu$ is
$SL(2,\reals)$-invariant. 
Now suppose $\nu$ is not affine. We can apply
Lemma~\ref{lemma:can:choose:in:U:perp}, and then
iterate Proposition~\ref{prop:get:more:points:in:Wplus} with
$\delta' \to 0$ and fixed $\epsilon$ and $\delta$. Taking a limit
along a subsequence we get points
$(b_\infty, x_\infty,s_\infty) \in K_0$ and $(b_\infty,
y_\infty,s_\infty)\in K_0$ such that
$\|x_\infty - y_\infty \| \approx \epsilon$, $y_\infty \in
W^+(b_\infty,x_\infty)$ and $y_\infty \in
(U^\perp)^+(b_\infty,x_\infty)$. 
This contradicts Lemma~\ref{lemma:exists:psi} since $K_0 \subset \Psi$. 
Hence $\nu$ is affine. 
\qed\medskip

\appendix
\renewcommand{\thesubsection}{\Alph{section}.\arabic{subsection}}

\section{Forni's results on the $SL(2,\reals)$ action}
\label{sec:app:forni}

In this appendix, we summarize the results we use from the fundamental
work of Forni \cite{Forni:Deviation}. The recent preprint
\cite{Forni:Matheus:Zorich} contains an excellent presentation of these
ideas and also some additional results which we will use as well.

\subsection{The Hodge norm and the geodesic flow}
\label{sec:subsec:app:hodge}
Let $\cM_g$ denote the
moduli space of genus $g$ curves.
Fix a point $S$ in $\cH(\alpha)$; then $S$ is a pair
$(M,\omega)$ where $M \in \cM_g$ and 
$\omega$ is a holomorphic $1$-form on $M$.  Let
$\|\cdot\|_{H,t}$ denote the Hodge norm 
(see e.g.\ \cite{ABEM})
at the surface $M_t = \pi(g_t S)$. 
Here 
$\pi: \cH(\alpha) \to \cM_g$  is the natural map taking
$(M,\omega)$ to $M$. We recall that the Hodge norm is a norm on $H^1(M,\reals)$.

The following fundamental result is due to Forni \cite[\S{2}]{Forni:Deviation}:
\begin{theorem}
\label{theorem:forni}
For any $\lambda \in H^1(M,\reals)$ and any $t \ge 0$, 
\begin{displaymath}
\|\lambda\|_{H,t} \le e^{t} \| \lambda\|_{H,0}.
\end{displaymath}
If in addition $\lambda$ is orthogonal to $\omega$, and 
for some compact subset $\cK$ of $\cM_g$, 
the geodesic segment $[S,g_t S]$ spends at least half the
time in $\pi^{-1}(\cK)$, then we have
\begin{displaymath}
\|\lambda\|_{H,t} \le e^{(1-\alpha)t} \| \lambda\|_{H,0},
\end{displaymath}
where $\alpha > 0$ depends only on $\cK$.
\end{theorem}

\bold{The Hodge norm on relative cohomology.}
Let $\Sigma$ denote the set of zeroes of $\omega$. 
Let $p: H^1(M,\Sigma,\reals) \to H^1(M,\reals)$ denote the
natural map. We define a norm $\| \cdot \|'$ on the relative cohomology group
$H^1(M,\Sigma,\reals)$ as follows:

\begin{equation}
\label{eq:df:relative:hodge}
\|\lambda\|' = \| p(\lambda)\|_H + \sum_{(z,w) \in \Sigma\cross \Sigma}
\left|\int_{\gamma_{z,w}} (\lambda - h)\right|,
\end{equation}
where $\| \cdot \|_H$ denotes the Hodge norm on $H^1(M,\reals)$,
$h$ is the harmonic representative of the cohomology class
$p(\lambda)$ and $\gamma_{z,w}$ is any path connecting the zeroes $z$
and $w$. Since $p(\lambda)$ and $h$ represent the same class in
$H^1(M,\reals)$, the equation (\ref{eq:df:relative:hodge}) does not
depend on the choice of $\gamma_{z,w}$.

Let $\| \cdot\|'_t$ denote the norm (\ref{eq:df:relative:hodge}) on the
surface $M_t$. Then, up to a fixed multiplicative constant, 
the analogue of Theorem~\ref{theorem:forni}
holds, for $\| \cdot \|'_t$, as long as $S\equiv(M,\omega)$ and $g_t S$
belong to a fixed compact set. This assertion is 
essentially Lemma~4.4 from \cite{Athreya:Forni}.
For a self-contained proof in this notation see
\cite[\S{8}]{Eskin:Mirzakhani:Rafi}. 




\mcc{define the AGY norm}

\bold{The Avila-Gou\"ezel-Yoccoz (AGY) norm.} 
The Hodge norm on relative cohomology behaves badly in the thin part
of Teichm\"uller space. Therefore, we will use instead the
Avila-Gou\"ezel-Yoccoz norm $\| \cdot \|_Y$ defined in
\cite{Avila:Gouezel:Yoccoz}, some properties of which were further developed in
\cite{Avila:Gouezel}.  The norms $\| \cdot \|_Y$ and $\| \cdot \|'$ are
equivalent on compact subsets of the strata $\cH_1(\alpha)$, and
therefore the decay estimates on $\| \cdot \|'$  in the style of
Theorem~\ref{theorem:forni} also apply to
the Avila-Gou\"ezel-Yoccoz norm. Furthermore, we have the following:
\begin{theorem}
\label{theorem:decay:relative:norm}
Suppose $S = (M,\omega) \in \cH(\alpha)$.  Let
$\|\cdot\|_{t}$ denote the Avila-Gou\"ezel-Yoccoz (AGY) norm 
on the surface $g_t S$. Then,
\begin{itemize}
\item[{\rm (a)}] For all $\lambda \in H^1(M,\Sigma,
\reals)$ and all $t > 0$, 
\begin{displaymath}
\|\lambda\|_t \le e^{t} \| \lambda\|_0.
\end{displaymath}
\item[{\rm (b)}] 
Suppose for some compact subset $\cK$ of $\cM_g$, the geodesic
segment $[S,g_t S]$ spends at least half the time in $\pi^{-1}(\cK)$. 
Suppose $\lambda \in H^1(M,\Sigma,\reals)$ with $p(\lambda)$
orthogonal to $\omega$. Then we have
\begin{displaymath}
\|\lambda\|_t \le C e^{(1-\alpha)t} \| \lambda\|_0,
\end{displaymath}
where $\alpha > 0$ depends only on $\cK$. 
\end{itemize}
\end{theorem}

\subsection{The Kontsevich-Zorich cocycle}
\label{sec:subsec:app:forni:cocycle}
We recall that $X_0$ denotes a finite cover of a stratum which is a
manifold (see \S\ref{sec:semi:markov}). 
In the sequel, a subbundle $L$ of the Hodge bundle is called {\em
  isometric} if the action of the Kontsevich-Zorich cocycle restricted
to $L$ is by isometries in the Hodge metric. We say that a subbundle
is {\em isotropic} if the symplectic form vanishes identically on the
sections, and {\em symplectic} if the symplectic form is
non-degenerate on the sections. A subbundle is {\em irreducible} if it
cannot be decomposed as a direct sum, and {\em strongly irreducible} if
it cannot be decomposed as a direct sum on any (measurable) 
finite cover of $X_0$. 

\begin{theorem}
\label{theorem:P:invariant:positive:sum}
Let $\nu$ be a $P$-invariant measure on $X_0$, and suppose $L$ is a
$P$-invariant $\nu$-measurable subbundle of the Hodge bundle. Let
$\lambda_1, \dots, \lambda_n$ be the Lyapunov exponents of the
restriction of the Kontsevich-Zorich cocycle to $L$. Then,
\begin{displaymath}
\sum_{i=1}^n \lambda_i \ge 0. 
\end{displaymath}
\end{theorem}

\bold{Proof.} Let the symplectic complement
$L^\dagger$ of $L$ be defined by
\begin{equation}
\label{eq:def:L:dagger}
L^\dagger(x) = \{ v \st v \wedge u = 0 \quad\text{ for all $u \in L(x)$} \}.
\end{equation}
Then, $L^\dagger$ is a $P$-invariant subbundle, and
we have the short exact sequence
\begin{displaymath}
0 \to L \cap L^\dagger \to L \to L/(L \cap L^\dagger) \to 0.
\end{displaymath}
The bundle $L/(L \cap L^\dagger)$ admits an invariant non-degenerate
symplectic form, and therefore, the sum of the Lyapunov exponents on
$L/(L \cap L^\dagger)$ is $\ge 0$. Therefore, it is enough to show that
the sum of the Lyapunov exponents on the isotropic subspace $L \cap
L^\dagger$ is $0$. Thus, without loss of generality, we may assume
that $L$ is isotropic.

Let $\{ c_1, \dots, c_n \}$ be a Hodge-orthonormal basis for the
bundle $L$ at the point $S = (M,\omega)$, where $M$ is a Riemann
surface and $\omega$ is a holomorphic $1$-form on $M$. For $g \in
SL(2,\reals)$, let $V_S(g)$ denote the Hodge norm of the polyvector $c_1
\wedge \dots \wedge c_n$ at the point $gS$, where the vectors $c_i$
are transported following a path from the identity to $g$ 
using the Gauss-Manin connection. (The result does not depend on the
path since the Gauss-Manin connection is flat, and $X_0$ has no
orbifold points). Since $V_S(k g) =
V_S(g)$ for $k \in SO(2)$, we can think of $V_S$ as a function on the
upper half plane $\half$. From the definition of $V_S$ and the
multiplicative ergodic theorem, we see that for
$\nu$-almost all $S \in X_0$, 
\begin{equation}
\label{eq:tmp:lim:log:VS}
 \lim_{t \to \infty} \frac{\log V_S(g_t)}{t} =
\sum_{i=1}^n \lambda_i, 
\end{equation}
where the $\lambda_i$ are as in the statement of
Theorem~\ref{theorem:P:invariant:positive:sum}. 

Let $\Delta_{hyp}$ denote the hyperbolic Laplacian operator (along the
Teichm\"uller disk). By \cite[Lemma~2.8]{Forni:Matheus:Zorich} (see
also \cite[Lemma 5.2 and Lemma 5.2']{Forni:Deviation})
there exists a non-negative
function $\Phi: X_0 \to \reals$ such that for all $S \in X_0$ and all $g
\in SL(2,\reals)$, 
\begin{displaymath}
(\Delta_{hyp} \log V_S)(g) = \Phi(g S). 
\end{displaymath}
We now claim that the Kontsevich-Forni type formula
\begin{equation}
\label{eq:formula:for:sum:of:exponents:subbundle}
\sum_{i=1}^n \lambda_i = \int_{X_0} \Phi(S) \, d\nu(S)
\end{equation}
holds, which clearly implies the
theorem. The formula 
(\ref{eq:formula:for:sum:of:exponents:subbundle}) is proved
in \cite{Forni:Matheus:Zorich} (and for the case of the entire stratum
in \cite{Forni:Deviation}) under the assumption that the measure
$\nu$ is invariant under $SL(2,\reals)$. However, in the proofs, only
averages over ``large circles'' in $\half = SO(2)\backslash
SL(2,\reals)$ are used. Below we show that a slightly modified version
of the proof works under the a-priori weaker assumption that $\nu$
is invariant under $P = AN \subset SL(2,\reals)$. This is not at all
surprising, since large circles in $\half$ are approximately
horocircles (i.e.\ orbits of $N$). 

We now begin the proof of
(\ref{eq:formula:for:sum:of:exponents:subbundle}), following the proof
of \cite[Theorem 1]{Forni:Matheus:Zorich}.

Since (\ref{eq:tmp:lim:log:VS}) holds for $\nu$-almost all $S$ and
$\nu$ is $N$-invariant, (\ref{eq:tmp:lim:log:VS}) also holds for  
almost all $S_0 \in X_0$ and almost all $S \in \Omega_N S_0$, 
where 
\begin{displaymath}
\Omega_N = \left\{ \begin{pmatrix} 1 & s \\ 0 & 1 \end{pmatrix} \st |s|
\le 1 \right\} \subset N. 
\end{displaymath}
We identify $SO(2) \backslash SL(2,\reals) S_0$ with $\half$ so that
$SO(2) g S_0$
corresponds to $g^{-1} \cdot i$. Then $\Omega_N S_0$ corresponds to the 
horizontal line segment connecting $-1+i$ to $1+i$. Let $\epsilon =
e^{-4t}$. \mcc{check normalization}
Then, $g_t \Omega_N S_0$ corresponds to the line segment connecting
$-1+i\epsilon$ to $1+i\epsilon$.

\makefig{Proof of Theorem~\ref{theorem:P:invariant:positive:sum}.}{fig:rectangle}{\includegraphics{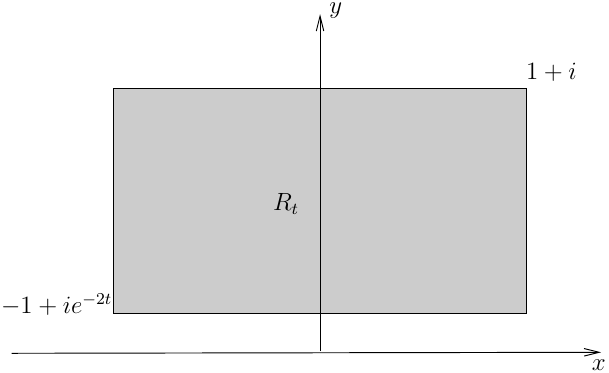}}

Let $f(z) = \log V_{S_0}(SO(2)z)$. Note that $\nabla_{hyp} f$ is
bounded (where $\nabla_{hyp}$ is the gradient with respect to the
hyperbolic metric on $\half$). Then,
(\ref{eq:tmp:lim:log:VS}) implies that for almost all $x \in [-1,1]$, 
\begin{displaymath}
\sum_{i=1}^n \lambda_i = \lim_{T \to \infty} \frac{
f(x+ie^{-2T})-f(x+i)}{T} = \lim_{T \to \infty} \frac{1}{T}\int_0^T
\frac{\partial}{\partial t} \left[ f(x+i e^{-2t}) \right] \, dt
\end{displaymath}
Integrating the above formula from $x=-1$ to $x=1$, we get (using the
bounded convergence theorem), 
\begin{displaymath}
\sum_{i=1}^n \lambda_i =  \lim_{T \to \infty} \frac{1}{T}\int_0^T
\left( \int_{-1}^1  \frac{\partial}{\partial t} \left[ f(x+i e^{-2t}) \right]
\, dx \right) \, dt
\end{displaymath}
Let $R_t$ denote the rectangle with corners at $-1 + i e^{-2t}$,
$1+ i e^{-2t}$, $1+i$ and $-1+i$, see Figure~\ref{fig:rectangle}. 
We now claim that
\begin{equation}
\label{eq:tmp:claim:rectangle}
\int_{-1}^1  \frac{\partial}{\partial t} \left[ f(x+i e^{-2t}) \right]
\, dx  = e^{-4t} \int_{\partial R_t } \frac{\partial f}{\partial n} + O(t e^{-4t}),
\end{equation}
where $\frac{\partial f}{\partial n}$ denotes the (outgoing) 
normal derivative of
$f$ with respect to the hyperbolic metric. 
Indeed, the integral over
the bottom edge of the rectangle $R_t$ on the left hand side
of (\ref{eq:tmp:claim:rectangle}) coincides with the right hand side
of (\ref{eq:tmp:claim:rectangle}) (the factor of $e^{-4t}$ appears
because the hyperbolic length element is $dx/y^2 = e^{-4t} \, dx$.) 
The partial derivative $\frac{\partial f}{\partial n}$ is uniformly
bounded, and the hyperbolic lengths of the other three sides of
$\partial R_t$ are $O(t)$. Therefore (\ref{eq:tmp:claim:rectangle})
follows. 

Now, by Green's
formula (in the hyperbolic metric),
\begin{displaymath}
\int_{\partial {R_t}}
\frac{\partial f}{\partial   n} = \int_{R_t} \Delta_{hyp} f =
\int_{R_t} \Phi, 
\end{displaymath}
We get, for almost all $S_0$,
\begin{displaymath}
\sum_{i=1}^n \lambda_i = \lim_{T \to \infty} \frac{1}{T} \int_0^T
\left(e^{-4t} \int_{R_t} \Phi \right) \, dt \ge 0. 
\end{displaymath}
This completes the proof of the Theorem. It is also easy to conclude
(by integrating over $S_0$) that
(\ref{eq:formula:for:sum:of:exponents:subbundle}) holds.  
\qed\medskip

\begin{theorem}
\label{theorem:zero:lyapunov:isometric}
Let $\nu$ be an ergodic $SL(2,\reals)$-invariant measure, and suppose $L$ is
an $SL(2,\reals)$-invariant $\nu$-measurable subbundle of the Hodge
bundle. Suppose all the Lyapunov exponents of the restriction of the
Kontsevich-Zorich cocycle to $L$ vanish. Then, the action of the
Kontsevich-Zorich cocycle on $L$ is isometric with respect to the
Hodge inner product, and the orthogonal complement $L^\perp$ of $L$
with respect to the Hodge inner product is also an
$SL(2,\reals)$-invariant subbundle. 
\end{theorem}

\bold{Proof.} The first assertion is the content of \cite[Theorem
3]{Forni:Matheus:Zorich}. The second assertion then follows from
\cite[Lemma 4.3]{Forni:Matheus:Zorich}. 
\qed\medskip

\begin{theorem}
\label{theorem:isotropic:zero:lyapunov}
Let $\nu$ be an ergodic $SL(2,\reals)$-invariant measure, and suppose $L$ is
an $SL(2,\reals)$-invariant $\nu$-measurable subbundle of the Hodge
bundle. Suppose $L$ is isotropic. Then all the Lyapunov
exponents of the restriction of the Kontsevich-Zorich cocycle to $L$
vanish (and thus Theorem~\ref{theorem:zero:lyapunov:isometric}
applies to $L$). 
\end{theorem}

\bold{Proof.} For a point $x \in X_0$ and an isotropic $k$-dimensional
subspace $I_k$, let $\Phi_k(x,I_k)$ be as in
\cite[(2.46)]{Forni:Matheus:Zorich} (or \cite[Lemma
5.2']{Forni:Deviation}). We 
have from \cite[Lemma 2.8]{Forni:Matheus:Zorich} 
that
\begin{displaymath}
\Phi_k(x,I_k) \le \Phi_j(x,I_j) \quad \text{ if $i < j$ and $I_k
  \subset I_j$}. 
\end{displaymath}
Let $\lambda_1 \ge \dots \ge \lambda_n$ be the Lyapunov exponents of
the restriction of the Kontsevich-Zorich cocycle to $L$. 
Let $\cV_{\le j}(x)$  denote the direct sum of all the
Lyapunov subspaces corresponding to exponents $\lambda_i \ge
\lambda_j$. By definition, $V_n(x) = L(x)$. Suppose $j=n$ or 
$\lambda_j \ne \lambda_{j+1}$. Then, by
\cite[Corollary~3.1]{Forni:Matheus:Zorich}
the following formula holds:
\begin{displaymath}
\lambda_1 + \dots + \lambda_j = \int_{X_0} \Phi_j(x,\cV_{\le j}(x)) \, d\nu(x)
\end{displaymath}
(This formula is proved in \cite{Forni:Deviation} for the case where
$\nu$ is Lebesgue measure and $L$ is the entire Hodge bundle). 

We will first
show that all the $\lambda_j$ have the same sign. Suppose not, then we
must have $\lambda_n < 0$ but not all $\lambda_j < 0$. Let $k$ be
maximal such that $\lambda_k \ne \lambda_n$. Then
\begin{displaymath}
\lambda_1 + \dots + \lambda_k = \int_{X_0} \Phi_k(x,V_k(x)) \, d\nu(x)
\end{displaymath}
and
\begin{displaymath}
\lambda_1+ \dots + \lambda_n  = \int_{X_0} \Phi_n(x,L(x)) \, d\nu(x)
\end{displaymath}
But $\Phi_k(x,V_k(x)) \le \Phi_n(x,L(x))$ since $V_k(x) \subset
L(x)$. Thus, 
\begin{equation}
\label{eq:lambda:k:plus:1:to:n}
\lambda_{k+1} + \dots + \lambda_n \ge 0. 
\end{equation}
But by the choice of $k$, all the terms in
(\ref{eq:lambda:k:plus:1:to:n}) are equal to each other.  
This implies that $\lambda_n \ge 0$, contradicting our
assumption that $\lambda_n < 0$. Thus all the $\lambda_j$, $1 \le j
\le n$ have the same sign. Since $\nu$ is assumed to be
$SL(2,\reals)$-invariant, and any diagonalizable $g \in SL(2,\reals)$
is conjugate to its inverse, we see that e.g.\ the  $\lambda_j$ cannot
all be positive. Hence, all the Lyapunov exponents $\lambda_j$ are
$0$. 
\qed\medskip

\bold{Algebraic Hulls.} The algebraic hull of a cocycle is defined in
\cite{ZimmerBook}. We quickly recap the definition: Let $G$ be a group
acting on a space $X$, preserving an ergodic measure $\nu$.  
Suppose $H$ is an $\reals$-algebraic group,
and let $A: G \cross X \to H$ be a measurable cocycle. We say that the
$\reals$-algebraic subgroup $H'$ of $H$ is the {\em algebraic hull} of
$A$ if $H'$ is the smallest $\reals$-algebraic subgroup of $H$ such
that there exists a $\nu$-measurable map $C: X \to H$ such that
\begin{displaymath}
C( g x)^{-1} A(g,x) C(x) \in H' \qquad \text{ for almost all $g \in G$
  and $\nu$-almost all $x \in X$. }
\end{displaymath}
It is shown in \cite{ZimmerBook} (see also
\cite[Theorem~3.8]{Morris:Zimmer}) that the algebraic hull exists and
is unique up to conjugation.

\begin{theorem}
\label{theorem:KZ:semisimple}
Let $\nu$ be an ergodic $SL(2,\reals)$-invariant measure. Then, 
\begin{itemize}
\item[{\rm (a)}] The $\nu$-algebraic hull $H'$ of the
  Kontsevich-Zorich cocycle is semisimple.
\item[{\rm (b)}] Every $\nu$-measurable
  \mcc{strongly???}$SL(2,\reals)$-invariant irreducible subbundle of the Hodge
  bundle is either symplectic or isotropic. 
\end{itemize}
\end{theorem}

\bold{Remark.} The fact that the algebraic hull is semisimple for
$SL(2,\reals)$-invariant measures is key to our approach.

\bold{Proof.}
Suppose $L$ is an invariant subbundle. It is enough to show that there
exists an invariant complement to $L$. Let the symplectic complement
$L^\dagger$ of $L$ be defined as in (\ref{eq:def:L:dagger}). 
Then, $L^\dagger$ is also an $SL(2,\reals)$-invariant subbundle, and $K = L
\cap L^\dagger$ is isotropic. By Theorem~\ref{theorem:isotropic:zero:lyapunov},
$K$ is isometric, and $K^\perp$ is also
$SL(2,\reals)$-invariant. Then, 
\begin{displaymath}
L = K \oplus (L \cap K^\perp), \qquad L^\dagger = K \oplus (L^\dagger \cap K^\perp),
\end{displaymath}
and
\begin{displaymath}
H^1(M, \reals) = K \oplus (L \cap K^\perp ) \oplus (L^\dagger \cap K^\perp)
\end{displaymath}
Thus, $L^\dagger \cap K^\perp$ is an $SL(2,\reals)$-invariant complement to
$L$. This proves (a). To prove (b), let $L$ be any irreducible
$SL(2,\reals)$-invariant $\nu$-measurable irreducible subbundle of the Hodge
bundle, and let $K = L \cap L^\dagger$. Since $K \subset L$ and $L$ is
irreducible, we have either $K = 0$ (so $L$ is symplectic), or $K = L$
and so $L$ is isotropic. 
The same could be done on any finite cover. 
\qed\medskip

\bold{The Forni subspace.} 

\begin{definition}[Forni Subspace]
\label{theorem:Forni:analytic:def}
Let
\begin{equation}
\label{eq:Forni:analytic:def}
F(x) = \bigcap_{g \in SL(2,\reals)} g^{-1} (\operatorname{Ann}
B^{\reals}_{gx}), 
\end{equation}
where for $\omega \in X_0$ the quadratic form
$B_\omega^\reals( \cdot, \cdot )$ is as defined in
\cite[(2.33)]{Forni:Matheus:Zorich}. 
\end{definition}

\bold{Remark.} It is clear from the definition, that as long as its
dimension remains constant, $F(x)$ varies real-analytically with $x$.  

\begin{theorem}
\label{theorem:forni:subspace:ergodic:def}
Suppose $\nu$ is an ergodic
$SL(2,\reals)$-invariant measure. Then the subspaces $F(x)$ where $x$
varies over the support of $\nu$ form the maximal 
$\nu$-measurable $SL(2,\reals)$-invariant 
isometric subbundle of the Hodge
bundle. 
\end{theorem}


\bold{Proof.} Let $F(x)$ be as defined in
(\ref{eq:Forni:analytic:def}). Then, $F$ is an
$SL(2,\reals)$-invariant subbundle of the Hodge bundle, 
and the restriction
of $B_x^\reals$ to $F(x)$ is identically $0$. Then, by \cite[Lemma
1.9]{Forni:Matheus:Zorich}, $F$ is isometric. 

Now suppose $M$ is any other $\nu$-measurable 
isometric $SL(2,\reals)$-invariant subbundle of
the Hodge bundle. Then by \cite[Theorem
2]{Forni:Matheus:Zorich}, $M(x) \subset \operatorname{Ann}
B_x^\reals$. Since $M$ is $SL(2,\reals)$-invariant, we have $M \subset
F$. Thus $F$ is maximal. 
\qed\medskip

\begin{theorem}
\label{theorem:properties:forni}
Let $\nu$ be an ergodic $SL(2,\reals)$-invariant measure on 
any finite cover of $X_0$. 
\begin{itemize}
\item[{\rm (a)}] For $\nu$-almost all $x \in X_0$, 
the Forni subspace $F(x)$ is symplectic, and its symplectic
  complement $F^\dagger(x)$ coincides with its Hodge complement $F^\perp(x)$. 

\item[{\rm (b)}] Any $\nu$-measurable $SL(2,\reals)$-invariant
  subbundle of  $F^\perp$ is symplectic,
  and the restriction of the Kontsevich-Zorich cocycle to
  any invariant subbundle of $F^\perp$ has at least one non-zero
  Lyapunov exponent. 
\end{itemize}
\end{theorem}

\bold{Proof.} 
Suppose the subspace $F^\perp$ is not symplectic. Let 
$L = F^\perp \cap (F^\perp)^\dagger$. Then
$L$ is isotropic, and therefore by
Theorem~\ref{theorem:isotropic:zero:lyapunov} and
Theorem~\ref{theorem:zero:lyapunov:isometric}, $L$ is an 
$SL(2,\reals)$-invariant isometric subspace. Hence $L \subset F$ by
Theorem~\ref{theorem:forni:subspace:ergodic:def}.  As $L \subset
F^\perp$ we get $L=0$. Therefore $F^\perp$ is symplectic.

Let $M$ be an irreducible
subbundle of $F^\perp$. Then, in view of
Theorem~\ref{theorem:zero:lyapunov:isometric} and the maximality of
$F$,  $M$ must have at least one non-zero Lyapunov exponent. In
particular, in view of Theorem~\ref{theorem:isotropic:zero:lyapunov},
$M$ cannot be isotropic, so it must be symplectic in view of
Theorem~\ref{theorem:KZ:semisimple} (b). This proves the statement
(b). 

Since $F^\perp$ is symplectic, $(F^\perp)^\dagger$ 
is $SL(2,\reals)$-invariant and complementary to $F^\perp$. Note that
$F$ is also $SL(2,\reals)$-invariant and complementary to $F^\perp$.
In order to conclude that $(F^\perp)^\dagger = F$, it is enough to show
that there is a unique $SL(2,\reals)$-invariant complement to
$F^\perp$. 

Note that another complement to $F^\perp$ would be the graph of an equivariant linear map $A: F \to F^\perp$.
If $A$ is nonzero, then an invariant complement of its kernel in $F$ exists
by Theorem~\ref{theorem:KZ:semisimple}, 
and it even contains an irreducible subbundle $M_2$. Then
$A$ induces an equivariant isomorphism between $M_2$
and its image, an irreducible subbundle $M_1$ of $F^\perp$.
Now, to get a contradiction, it is enough to show that for any 
irreducible subbundles $M_1 \subset F^\perp$ and $M_2 \subset F$, the
algebraic hulls $H'(M_i)$ of the restriction of the Kontsevich-Zorich
cocycle to $M_i$ are not isomorphic to each other. 
But the later statement is clear,
since $H'(M_2)$ is compact and $H'(M_1)$ is not (since it has at
least one non-zero Lyapunov exponent by (b)). Thus, $(F^\perp)^\dagger
= F$. Since we already showed that $F^\perp$ is symplectic, this
implies that so is $F$, which completes the proof of (a). 
\qed\medskip

\section{Entropy and the Teichm\"uller geodesic flow}
\label{sec:appendix:margulis:tomoanov}
The contents of this section are well-known, see
e.g.\ \cite{Ledrappier:Young}, \cite{Margulis:Tomanov:Ratner} and also
\cite{Bufetov:Gurevich}. However, for technical reasons, 
the statements we need do not formally follow
from the results of any of the above papers. Our setting
is intermediate between the homogeneous dynamics setting of \cite{Margulis:Tomanov:Ratner}
and the general $C^2$-diffeomorphism on a compact manifold setup of
\cite{Ledrappier:Young}, but it is closer to the former than the latter. What
follows is a lightly edited but almost verbatim reproduction of
\cite[\S{9}]{Margulis:Tomanov:Ratner}, adapted to the setting of
Teichm\"uller space. It is 
included here primarily for the convenience of the reader. The (minor)
differences between our presentation and that of
\cite{Margulis:Tomanov:Ratner}   
are related to the lack of  uniform hyperbolicity outside of
compact subsets of the space, and some notational changes due to the
fact that our space is not homogeneous. 

\bold{Notation.} We recall some notation from \S\ref{sec:subsec:notation}.
Let $X_0$ denote the finite cover of $\cH_1(\alpha)$ defined in
\S\ref{sec:semi:markov} (which has no orbifold points). 
Let $g_t$ denote the Teichm\"uller geodesic flow.
In this section,  $\nu$ is an ergodic $g_t$-invariant probability measure on
$X_0$. Let $V(x)$ denote a subset of $H^1(M,\Sigma,\reals^2)$. Then we
denote 
\begin{displaymath}
V[x] = \{ y \in X_0 \st y - x \in V(x) \}. 
\end{displaymath}
This makes sense in a neighborhood of $x$. 

Let $d^{X_0}(\cdot,\cdot)$ denote the AGY distance on $X_0$, defined in
\S\ref{sec:semi:markov}. 
Fix a point $p \in X_0$ (so $p$ is not an orbifold point), and 
such that every neighborhood of $p$ in $X_0$ has positive
$\nu$-measure. 
Fix relatively compact neighborhoods  $C'(p)$
and $Q(p)$ of $0$ in $W^+(p)$ and $\reals$ respectively.
Let 
\begin{displaymath}
C = \bigcup_{t \in Q(p)} g_t C'[p], 
\end{displaymath}
For each $c \in \overline{C}$ 
choose a relatively compact neighborhood $B'(c)$ of
$0$ in $W^-(c)$ with diameter in the AGY distance at most $1/200$
so that the $B'(c)$ vary continuously with $c$. 
For $c \in C$, let
\begin{displaymath}
B'[c] = \{ c +v  \st v \in B'(c)\}, \qquad D = \bigsqcup_{c \in C} B'[c].
\end{displaymath}
We assume that $C'(p)$, $Q(p)$ and the $B'(c)$ are
sufficiently small so that $D$ is open and contractible. 

\begin{lemma}(cf.\ \cite[Lemma 9.1]{Margulis:Tomanov:Ratner})
\label{lemma:mt:91}
There exists $s > 0$, 
$C_1 \subset C$ and for each $c \in C_1$ there exists a
subset $E[c] \subset W^-[c]$ such that
\begin{itemize}
\item[{\rm (1)}] $E[c] \subset B'[c]$.
\item[{\rm (2)}] $E[c]$ is open in $W^-[c]$, and the subset $E \equiv \bigcup_{c
    \in C_1} E[c]$ satisfies $\nu(E) > 0$. 
\item[{\rm (3)}] Let $T = g_s$ denote the time $s$ map
    of the geodesic flow. Then whenever 
\begin{displaymath}
T^n E[c] \cap E \ne 0, \qquad c \in C_1, \qquad n > 0, 
\end{displaymath}
we have $T^n E[c] \subset E$.    
\end{itemize}
\end{lemma}

\bold{Proof.} 
Fix a compact subset $K_1 \subset X_0$, with $\nu(K_1^c) < 0.01$. Then by
the Birkhoff ergodic theorem, for
every $\delta > 0$ there exists $R > 0$ and a subset $E_1$ with
$\nu(E_1) > 1-\delta$ such that for all $x \in E_1$ and all $N > R$, 
\begin{displaymath}
|\{ n \in [1,N] \st g_n x \in K_1 \}| \ge (1/2) N. 
\end{displaymath}
By choosing $\delta > 0$ small enough, we may assume that $\nu(D \cap
E_1) > 0$. Let
\begin{displaymath}
C_1 = \{ c \in C \st c + v \in D \cap E_1 \text{ for some $v \in
  B'(c)$ } \}. 
\end{displaymath}
Then there exists a compact $K \supset K_1$ such that for all $c \in
C_1$ and all $x \in B'[c]$, 
\begin{displaymath}
|\{ n \in [1,N] \st g_n x \in K \}| \ge (1/2) N. 
\end{displaymath}
By Lemma~\ref{lemma:forni} there exists $\alpha > 0$ such that for all
$c \in C_1$ and all $x \in B'[c]$, 
\begin{displaymath}
d^{X_0}(g_n x, g_n c) \le \begin{cases} d^{X_0}(x,c) & \text{if $n \le R$ } \\
d^{X_0}(x,c) e^{-\alpha(n-R)} & \text{ if $n > R$ } \end{cases}
\end{displaymath}
Therefore we may choose $s > 0$ such that if we let $T
  = g_s$ denote the time $s$ map of the geodesic flow, 
then for all $c \in C_1$ and all $x \in B'[c]$,
\begin{displaymath}
d^{X_0}(Tx,Tc) \le \frac{1}{10} d^{X_0}(x,c). 
\end{displaymath}
There exists $a > 0$ so that for all $c \in C_1$, 
$B'[c]$ contains the intersection with $W^-[c]$ of a ball in the
AGY metric of
radius $a$ and centered at $c$. 
Let 
\begin{displaymath}
a_0 = \frac{a}{10}
\end{displaymath}
Let $B_0'[c] \subset W^-[c]$ denote the ball in the AGY metric of
radius $a_0$ and centered at $c$.
Let $E^{(0)}[c] = B_0'[c]$, and for $j > 0$ let
\begin{displaymath}
E^{(j)}[c] = E^{(j-1)}[c] \cup \{ T^n B_0'[c'] \st c' \in C_1, n > 0
\text{ and } T^n B_0'[c'] \cap E^{(j-1)}[c] \ne 0 \}. 
\end{displaymath}
Let 
\begin{displaymath}
E[c] = \bigcup_{j\ge 0} E^{(j)}[c], \qquad \text{ and } E = \bigcup_{c
  \in C_1} E[c]. 
\end{displaymath}
 It easily follows from the above definition that $E[c]$ 
has the properties (2) and (3). 
To show (1),  it is enough to show
that for each $j$, 
\begin{equation}
\label{eq:mt:induction}
d^{X_0}(x,c) < a/2, \quad \text{ for all } x \in E^{(j)}[c]. 
\end{equation}
This is done by induction on $j$. 
The case $j=0$ holds since $a_0 = a/10 < a/2$. Suppose
(\ref{eq:mt:induction}) holds for $j-1$, and 
suppose $x \in E^{(j)}[c]\setminus E^{(j-1)}[c]$. Then there exist
$c_0 =c, c_1, \dots, c_j = x$ in $C_1$ 
and non-negative integers $n_0 = 0, \dots, n_j$ such that for
all $1 \le k \le j$, 
\begin{equation}
\label{eq:pairwise:intersect}
T^{n_k}(B_0'[c_k]) \cap T^{n_{k-1}}(B_0'[c_{k-1}]) \ne \emptyset. 
\end{equation}
Let $1\le k\le j$ be such that $n_k$ is minimal. 
Recall that $B'[y] \cap B'[z] = \emptyset$ if $y \ne z$, $y \in C_1$, $z \in
C_1$. Therefore, in view of the inductive assumption, $n_k \ge 1$. 
Applying $T^{-n_k}$ to (\ref{eq:pairwise:intersect}) we get
\begin{displaymath}
\left(\bigcup_{i=1}^{k-1} T^{n_i-n_k} B_0'[c_i] \right) \cap B_0'[c_k] \ne
  \emptyset, \quad \text{
  and } \quad
\left(\bigcup_{i=k+1}^j T^{n_i-n_k} B_0'[c_i] \right) \cap B_0'[c_k] \ne
\emptyset. 
\end{displaymath}
Therefore, in view of (\ref{eq:pairwise:intersect}), and the
definition of the sets $E^{(j)}[c]$, 
\begin{displaymath}
\left(\bigcup_{i=1}^k T^{n_i-n_k} B_0'[c_i] \right) \subset
E^{(k-1)}[c_k], \quad \text{
  and } \quad
\left(\bigcup_{i=k}^j T^{n_i-n_k} B_0'[c_i] \right) \subset E^{(j-k)}[c_k]
\end{displaymath}
By the induction hypothesis, $\diam(E^{(k-1)}[c_k]) < a/2$, and
$\diam(E^{(j-k)}[c_k]) < a/2$. Therefore, 
\begin{displaymath}
\diam\left(\bigcup_{i=1}^j T^{n_i-n_k} B_0'[c_i] \right) \le a.  
\end{displaymath}
Then, applying $T^{n_k}$ we get,
\begin{displaymath}
\diam\left(\bigcup_{i=1}^j T^{n_i} B_0'[c_i] \right) \le \frac{a}{10}
\end{displaymath}
Since $\diam(B_0'[c]) \le a/10$, we get
\begin{displaymath}
\diam\left(\bigcup_{i=0}^j T^{n_i} B_0'[c_i] \right) \le
\diam(B_0'[c_0]) + \diam\left(\bigcup_{i=1}^j T^{n_i} B_0'[c_i] \right) \le
\frac{a}{10}+ \frac{a}{10} < \frac{a}{2}. 
\end{displaymath}
But the set on the left-hand-side of the above equation contains both
$c = c_0$ and $x = c_j$. Therefore $d^{X_0}(c,x) < a/2$, proving
(\ref{eq:mt:induction}). Thus (1) holds. 
\qed\medskip
\begin{lemma}(Ma\~n\'e)
\label{lemma:mane}
Let $E$ be a measurable subset of $X_0$, with $\nu(E) > 0$. If $\nu$ is
a compactly supported
measure on $E$ and $q: E \to (0,1)$ is such that $\log q$ is
$\nu$-integrable, 
then there exists a countable partition $\cP$ of $E$ with entropy
$H(\cP) < \infty$ \mccc{(define!)}
such that, if  $\cP(x)$ denotes the atom of $\cP$ containing $x$, then
$\diam \cP(x)<q(x)$. 
\end{lemma}
\bold{Proof.} See \cite{Mane:article} or \cite[Lemma 13.3]{Mane:book}
\qed\medskip

Let $V(x)$ be a system of real-algebraic subsets of $W^-(x)$. 

\begin{definition}
\label{def:admissible:system}
The system $V(x)$ is admissible if it is $T$-equivariant 
and also for almost all $x \in X_0$, $x$ is a smooth point of
$V[x]$. 
\end{definition}

\begin{definition}
\label{def:subordinate}
We say that a measurable partition $\xi$ of the measure space $(X_0, \nu)$
is {\em subordinate} to an admissible system of real-algebraic subsets $V(x) \subset W^-(x)$ if for almost all (with respect to
$\nu$) $x \in X_0$, we have
\begin{itemize}
\item[{\rm (a)}]   $\xi[x] \subset V[x]$ where $\xi[x]$ denotes, as
  usual, the element of $\xi$ containing $x$.
\item[{\rm (b)}]  $\xi[x]$ is relatively compact in $V[x]$. 
\item[{\rm (c)}]   $\xi[x]$ contains a neighborhood of $x$ in $V[x]$.
\end{itemize}
\end{definition}

Let $\eta$ and $\eta'$ be measurable partitions of $(X_0, \nu)$. 
We write $\eta \le  \eta'$ if $\eta[x] \supset \eta'[x]$ for 
almost all (with respect to $\nu$) $x \in X_0$. 
We define a partition $T \eta$ by $(T\eta)[x] = T (\eta [T^{-1}(x)])$.

\begin{proposition}
\label{prop:mt:93}
Assume that $\nu$ is $T$-ergodic (where $T$ is as in
Lemma~\ref{lemma:mt:91}(3)). Then there exists a measurable partition
$\eta$ of the measure space $(X_0, \nu)$ with the following properties:
\begin{itemize}
\item[{\rm (i)}] $\eta$ is subordinate to $W^-$.
\item[{\rm (ii)}] $\eta$ is $T$-invariant, i.e.\ $\eta \le  T \eta$.
\item[{\rm (iii)}] The mean conditional entropy $H(T\eta \mid \eta)$ 
\mccc{(define!)} is equal to the entropy $h(T, \nu)$ of the
 automorphism $x \to T x$ of the measure space $(X_0, \nu)$.
\end{itemize}
\end{proposition}

\bold{Proof.}
Let $E[c]$ and $E$ be as in Lemma~\ref{lemma:mt:91}. Denote by $\pi: E
\to C_1$ 
the natural projection ($\pi(x)=c$ if $x \in E[c]$). 
We set $\eta[x]=E(\pi(x)$) for every $x \in E$.

We claim that it is enough to find a countable measurable partition
$\xi$ of $(X_0,\nu)$ such that
$H(\xi)< \infty$ and $\eta[x]=\xi^-[x]$ for almost all $x \in E$ where
$\xi^- = \bigvee_{n=0}^\infty T^{-n} \xi$ is the product
of the partitions $T^{-n}\xi$, $0 \le n < \infty$.                 

Indeed, suppose the claim holds. Then it is clear that $\eta$ is
$T$-invariant. The set of $x \in X_0$ for which properties (a) and (b)
(resp.\ (c)) in the
definition of a subordinate partition are satisfied is $T^{-1}$-invariant (resp.\ $T$-invariant)
and contains $E$. But $\nu(E)>0$ and $\nu$ is $T$-ergodic. 
Therefore, $\eta$ is subordinate
to $W^-$. 
To check the property (iii) it is enough to show that the partition
$\xi_s = \bigvee_{k=-\infty}^{\infty} T^k \xi$ is the partition into
points, see \cite[\S{9}]{Rokhlin:book}, or \cite[\S{4.3}]{Katok:Hasselblatt}.  
By \cite{Forni:Deviation} or \cite[Theorem 8.12]{ABEM}
$\xi_s(x) = \{ x \}$ if $T^{-n} x \in E$
for infinitely many $n$. (Recall that by the construction of $E$, any
such geodesic will spend at least half the time in the compact set
$K$).  But $\nu(E) > 0$ and $\nu$ is
$T$-ergodic. Hence $\xi_s[x] = \{x\}$ for almost all $x$, which
completes the proof of the claim. 

Let us construct the desired partition $\xi$. For $x\in E$, let $n(x)$
be the smallest positive integer $n$ such that $T^n x \in E$. 
We have the classical Kac formula \cite{Kac:formula}
\begin{equation}
\label{eq:roklin:tower}
\int_{E} n(x) \, d\nu(x) = 1. 
\end{equation}
Define a probability measure $\nu'$ on $C_1$ by
\begin{equation}
\label{eq:def:nu:prime}
\nu'(F) = \frac{\nu(\pi^{-1}(F))}{\nu(E)}, \qquad F \subset C_1. 
\end{equation}
Property (3) of the family $\{E[c] \st c \in C_1 \}$ 
implies that $n(x)$ is constant on every $E[c]$,
$c \in C_1$. Therefore, in view of (\ref{eq:roklin:tower}) and
(\ref{eq:def:nu:prime}), 
\begin{displaymath}
\int_{C_1} n(c) \, d\nu'(c) < \infty. 
\end{displaymath}
By Lemma~\ref{lemma:forni:upper}, 
there exists $\kappa > 1$ such that for all $x$, $y \in X_0$, 
\begin{displaymath}
d^{X_0}(Tx,Ty) \le \kappa d^{X_0}(x,y). 
\end{displaymath}
Since the
function $n(c)$ is $\nu'$-integrable, 
one can find a positive function $q(c)<\kappa^{-2n(c)}$, $c \in C_1$
such that $\log q$ is $\nu'$-integrable, and the $\nu'$-essential
infimum 
$\operatorname{ess} \inf_{c \in C_1} q(c)$ is $0$.

After replacing, if necessary, $C'(p)$, $Q(p)$ and the $B'(c)$ for $c
\in \overline{C}$ by smaller
subsets we can find $\epsilon > 0$ such that the minimum distance
between lifts of $E$ is at most $\epsilon/10$ and also
\begin{itemize}
\item[{\rm (a)}] $d^{X_0}(x,y) < 2d (\pi(x), \pi(y))$  whenever $x, y \in E$
  and $d^{X_0}(x,y) < \epsilon$, and
\item[{\rm (b)}] if $x, y \in C_1$ then $d^{X_0}(x,y) < \epsilon$. 
\end{itemize}
Since the function $\log q(c)$ is $\nu'$-integrable, 
there exists a countable measurable
partition $\cP$ of $C_1$ such that $H(\cP ) < \infty$ and $\diam \cP(x) <
\tfrac{\epsilon}{2}q(x)$ for almost all $x \in C_1$
(see Lemma~\ref{lemma:mane}). After possibly replacing $\cP$ by a
countable refinement, we may assume that the function $x \to n(x)$ is
constant on the atoms of $\cP$. 
Now we define a countable measurable
partition $\xi$  of $X_0$ by
\begin{displaymath}
\xi(x) = \begin{cases} \pi^{-1}(\cP(\pi(x))) & \text{ if $x \in E$} \\
 X_0\setminus E & \text{ if $x \not\in E$.} \end{cases} 
\end{displaymath}
Since $H (\cP) < \infty$ we get using (\ref{eq:def:nu:prime}) that $H
(\xi) < \infty$. It remains to show that $\xi^-[x]= \eta[x]$ 
for almost all $x \in E$.  It follows from the property (3) of the
family $\{E[c]\}$ that $\eta[z] \subset \xi^-[z]$.
Let $x$ and $y$ be elements in $E$ with $\xi^-[x]=\xi^-[y]$. Since
$\eta[z] \subset \xi[z]$, we can assume that $x$, $y \in C_1$. Then
$d^{X_0}(x,y) < \epsilon$. 
Set $x_1 = x$ , $y_1 =y$ and define by induction 
\begin{displaymath}
x_{k+1} = T^{n(x_k)} x_k, \quad y_{k+1} = T^{n(y_k)} y_k. 
\end{displaymath}
Then, the sequence $\{x_k\}_{k \in \natls}$ (resp.\ $\{y_k\}_{k \in
  \natls}$) is the part of the
$T$-orbit of $x$ (resp.\ $T$-orbit of $y$) which lies in $E$. 

Let
$\tilde{x_1}$, $\tilde{y_1}$ be the lifts of $x_1=x$ and $y_1=y$ 
to Teichm\"uller
space,  \mcc{fix notation?}
and let $\tilde{x}_k$, $\tilde{y}_k$ be defined inductively by
\begin{displaymath}
\tilde{x}_{k+1} = T^{n(x_k)} \tilde{x}_k, \quad \tilde{y}_{k+1} = T^{n(y_k)} \tilde{y}_k. 
\end{displaymath}
Then $\tilde{x}_k$ and $\tilde{y}_k$ are lifts of $x_k$ and $y_k$
respectively. We now claim that for all $k \ge 0$, 
\begin{equation}
\label{eq:mt:equation5}
d^{X_0}(\tilde{x}_k, \tilde{y}_k) < \epsilon q(\pi(x_k)). 
\end{equation}
If $k = 1$, the inequality (\ref{eq:mt:equation5}) is true because
$\diam \cP(x) < \frac{\epsilon}{2} q(\pi(x))$ and $\cP(x) = \cP(y)$. 
Assume that (\ref{eq:mt:equation5}) is proved for $k$. Then
\begin{displaymath}
d^{X_0}(\tilde{x}_{k+1}, \tilde{y}_{k+1}) = d^{X_0}(T^{n(x_k)} \tilde{x}_k,
  T^{n(x_k)} \tilde{y}_k) \le \kappa^{n(x_k)}
  d^{X_0}(\tilde{x}_k,\tilde{y}_k) \le \kappa^{n(x_k)} \epsilon q(\pi(x_k))
  \le \epsilon. 
\end{displaymath}
Then since $x_{k+1}$ and $y_{k+1}$ belong to the same element of the
partition $\xi$ (because $\xi^-[x] = \xi^-[y]$) and $\diam(\cP(x_{k+1}))
\le \frac{\epsilon}{2} q \pi(x_{k+1})$, we get from condition (b) in 
the definition of
$\epsilon > 0$ that (\ref{eq:mt:equation5}) is
true for $k + 1$. 

Since the measure $\nu$ is $T$-ergodic and $\operatorname{ess} \inf q
(c) = 0$  we may assume that \linebreak 
$\liminf_{k \to \infty} 
q(\pi(x_k))=0$ (since this holds for almost all $x \in E$). 
Then (\ref{eq:mt:equation5}) implies that 
\begin{displaymath}
\liminf_{k \to \infty} d^{X_0}(\tilde{x}_k,\tilde{y}_k) = 0.  
\end{displaymath}
By the definition of $\tilde{x}_k$, $\tilde{y}_k$, there exists
a sequence $m_k \to +\infty$ such that $\tilde{x}_k = T^{m_k}
\tilde{x}$, $\tilde{y}_k = T^{m_k} \tilde{y}$. Thus, 
\begin{displaymath}
d^{X_0}(T^{m_k} \tilde{x}, T^{m_k} \tilde{y}) = 0. 
\end{displaymath}
But, by construction $\tilde{x}$ and $\tilde{y}$ are on the same leaf of
$W^{0+}$. This contradicts the non-contraction property of the Hodge
distance \cite[Theorem 8.2]{ABEM}, unless $\tilde{x}=\tilde{y}$. Thus we must
have $x=y$.  
\qed\medskip

\begin{lemma}(see \cite[Proposition 2.2]{Ledrappier:Strelcyn}.) 
\label{lemma:Ledrappier:Strelcyn}
Let $T$ be an automorphism of a measure
space $(X_0, \nu)$, $\nu(X_0)< \infty$, and let $f$ 
be a positive finite measurable function defined
on $X_0$ such that
\begin{displaymath}
\log^- \frac{f \circ T}{f} \in L^1(X,\nu), \quad\text{where }
\log^-(a)= \min(\log a,0).
\end{displaymath}
Then
\begin{displaymath}
\int_{X_0} \log \frac{f \circ T}{f} \, d\nu = 0. 
\end{displaymath}
\end{lemma}

\mcc{missing lemma here (in the comments)}


Suppose $V^-(x) \subset W^-(x)$ is an admissible $T$-equivariant family of
real-algebraic subsets. Let $(T_\reals V^-)(x) \subset W^-(x)$ denote
the tangent space to smooth manifold $V^-[x]$ at $x$. (Recall that
since $V^-$ is admissible, for almost every $x$, $V^-[x]$ is smooth at
$x$). 

\begin{definition}[Margulis Property]
\label{def:margulis:property}
Suppose $V^-(x) \subset W^-(x)$ is an admissible $T$-equivariant family of 
real-algebraic subsets. Let $\tau = \tau(x)$ be a measure on each
$V^-[x]$. We say that $\tau$ has the {\em Margulis Property} if for
almost all $x$, $\tau(x)$ is in the Lebesgue measure class on $V^-[x]$, and also $T_* \tau(x)$ agrees with $\tau(Tx)$ up to
normalization. (In other words the Radon-Nykodym derivative
$\frac{dT_*\tau(x)}{d\tau(Tx)}$ is locally constant along $V^-[x]$). 
\end{definition}

\begin{proposition}
\label{prop:mt:prop96}
Let $T = g_s$ as in Lemma~\ref{lemma:mt:91}(iii). 
Let $V^-(x) \subset W^-(x)$ be a $T$-equivariant family of
real-algebraic subsets. 
Suppose there exists a $T$-invariant measurable partition
$\eta$ of $(X_0,\nu)$ subordinate to $V^-$. Then the following hold:
\begin{itemize}
\item[{\rm (a)}] We have
\begin{displaymath}
H(T\eta \mid \eta) \le s \Delta(V^-),
\end{displaymath}
where $H(T\eta \mid
  \eta)$ is the mean conditional entropy, and 
\begin{displaymath}
\Delta(V^-) = \sum_{i \in I(V)} (1-\lambda_i), 
\end{displaymath}
where $I(V)$ are the Lyapunov subspaces in $T_\reals V$ (counted with multiplicity),
\mcc{improve notation?} and $\lambda_i$ are the corresponding
Lyapunov exponents of the Kontsevich-Zorich cocycle.  
\item[{\rm (b)}] Suppose that for almost all $x$ there exists a
  measure $\tau = \tau(x)$ on 
each $V^-[x]$ with the Margulis property. Then 
\begin{itemize}
\item[{\rm (b1)}] If the conditional measures of $\nu$ along $V^-[x]$
  agree with $\tau(x)$ (up to normalization), then   
\begin{displaymath}
H(T\eta \mid \eta) = s \Delta(V^-) 
\end{displaymath}
\item[{\rm (b2)}] If $H(T\eta \mid \eta) = s \Delta(V^-)$ then the
  conditional measures of $\nu$ along $V^-[x]$ agree with $\tau(x)$ (up
  to normalization). 
\end{itemize}
\end{itemize}
\end{proposition}

\bold{Proof.} Since $\eta \le T \eta$ for almost all $x \in X_0$ we have
a partition $\eta_x$ of $\eta[x]$ such that $\eta_x[y]=(T\eta)[y]$ for
almost all $y \in \eta[x]$. Let $\tau(x)$ be a measure on $V^-(x)$ in 
the Lebesgue measure class. (To simplify notation, we will sometimes
denote $\tau(x)$ simply by $\tau$).  
(Here we pick some  
normalization of the Lebesgue measure on 
the connected components of the intersections of the leaves of $V^-$
with a fixed 
fundamental domain). \mccc{(We must pick a good
  normalization or else $\log J$ is not integrable.)}
\mcc{explain more}  Since $\eta[x] \subset
V^-[x]$, $\tau$ induces a measure on $\eta[x]$ which we will denote
also by $\tau$.  Let $J(x)$ denote the Jacobian of the restriction of
the map $T$ to $V^-[x]$ at $x$
(with respect to the Lebesgue measure class measures $\tau$ 
on $V^-[x]$ and $V^-[Tx]$). Then, by the
Osceledets multiplicative ergodic theorem, for almost all $x \in X_0$, 
\begin{displaymath}
-s \Delta(V^-) = \lim_{N \to \infty} \frac{1}{N} \log \frac{d (T^{-N} \tau)(x)}{d\tau(x)} =
-\lim_{N \to \infty} \frac{1}{N} \sum_{n=0}^{N-1} \log J(T^{-n}x). 
\end{displaymath}
Integrating both sides over $X_0$, we get
\begin{equation}
\label{eq:mt:jacobian}
-\int_{X_0} \log J(x) \, d\nu(x) = s \Delta(V^-). 
\end{equation}
Put $L(x) = \tau(\eta[x])$ and $\tau_x = \tau/L(x)$,
$x \in X_0$.  Note that on $\eta[x]$ we have a conditional probability
measure $\nu_x$ induced by $\nu$. Put $p(x)= \tau_x(\eta_x[x])$ and
$r(x)= \nu_x(\eta_x[x])$.  

Let
\begin{equation}
\label{eq:def:etaprime}
\eta' = \eta \vee T \eta \vee \dots \vee T^k \eta. 
\end{equation}
Then, $\eta'$ is also $T$-invariant, and $H(T \eta' \mid \eta') = H(T
\eta \mid \eta)$. Thus, we can replace $\eta$ by $\eta'$. 

Suppose $\epsilon > 0$ is given. Then, we can choose $k$ large enough
in (\ref{eq:def:etaprime}) so that (after replacing $\eta$ by
$\eta'$), on a set of measure at least
$(1-\epsilon)$, we have 
\begin{equation}
\label{eq:jacobian:estimate}
(1-\epsilon) \le \frac{p(x) L(x)}{J(T^{-1} x) L(T^{-1}x)} \le (1+\epsilon)
\end{equation}

From its definition, $p(x) \le 1$. Also
\begin{equation}
\label{eq:mt:equation7}
-\int_{X_0} \log r(x) \, d\nu(x) = H(T\eta \mid \eta).
\end{equation}

Let $Y_i(x)$, $1 \le i < \infty$ denote the elements of the 
countable partition $\eta_x$ of $\eta[x]$.
Then we have
\begin{equation}
\label{eq:mt:equation8}
\int_{\eta(x)} \log p(y) \, d\nu_x(y) - 
\int_{\eta(x)} \log r(y) \, d\nu_x(y) = \sum_{i=1}^\infty \log
\frac{\tau_x(Y_i(x))}{\nu_x(Y_i(x))} \nu_x(Y_i(x)).  
\end{equation}
We have that
\begin{equation}
\label{eq:mt:equation9}
\sum_{i=1}^\infty \tau_x(Y_i(x)) \le 1, 
\end{equation}
and
\begin{equation}
\label{eq:mt:equation10}
\sum_{i=1}^\infty \nu_x(Y_i(x)) =1. 
\end{equation}
(In (\ref{eq:mt:equation9}), 
we can have strict inequality because apriori it is possible that the measure
$\tau_x$ of $\eta[x]\setminus \bigcup_{i=1}^\infty Y_i(x)$ is
positive). 
From (\ref{eq:mt:equation8}), (\ref{eq:mt:equation9}) and
(\ref{eq:mt:equation10}), using the convexity of log we get that
\begin{displaymath}
\int_{\eta(x)} \log p(y) \, d\nu_x(y) \le \int_{\eta(x)} \log r(y) \,
d\nu_x(y). 
\end{displaymath}
and the equality holds if and only if $p(y)=r(y)$
i.e.\ $\tau_x(\eta_x[y])=\nu_x(\eta_x[y])$ for all $y \in \eta[x]$. 
Now using integration over the quotient space $(X_0, \nu)/\eta$ 
of the measure space $(X_0, \nu)$ by $\eta$, 
we get from (\ref{eq:mt:equation7}) that
\begin{equation}
\label{eq:mt:key:estimate}
H(T \eta \mid \eta) \le - \int_{X_0} \log p(x) \, d\nu(x), 
\end{equation}
and the equality 
holds if and only if $\tau_x((T \eta) [x]) = \nu_x ((T \eta) [x])$
for almost all $x \in X_0$.

In view of (\ref{eq:jacobian:estimate}) and the fact that $p(x) \le
1$, 
\begin{displaymath}
- \int_{X_0} \log p(x) \, d\nu(x) \le 2 \epsilon - \int_{X_0} \log J(x) \,
dv(x) + \int_{X_0} \log_- (L(T^{-1}x)/L(x)) \, d\nu(x).  
\end{displaymath}
The last term vanishes by Lemma~\ref{lemma:Ledrappier:Strelcyn}. 
Since $\epsilon > 0$ is arbitrary, we have, by
(\ref{eq:mt:key:estimate}) and (\ref{eq:mt:jacobian}) that (a) holds. 

Now suppose that $\tau$ is as in (b). 
Then since $\eta_x[x] = T(\eta[T^{-1}x])$
one easily sees that $p(x)= J(T^{-1} x)L(T^{-1}x)/L(x)$. Therefore, by
(\ref{eq:mt:jacobian}) and Lemma~\ref{lemma:Ledrappier:Strelcyn}, 
\begin{displaymath}
-\int_{X_0} \log p(x) \, d\nu(x) = s \Delta(V^-). 
\end{displaymath}
If the conditional measures of $\nu$ along $V^-$ coincide with $\tau$,
then $p(x) = r(x)$ and therefore equality in
(\ref{eq:mt:key:estimate}) holds. This proves (b1). Conversely,
assume that $H(T \eta \mid \eta ) = s \Delta(V^-)$. 
Then $H(T^k \eta \mid \eta) = k s \Delta(V^-)$ for every
$k > 0$. Using the same argument as above and replacing $T$ by $T^k$, 
we get that $\tau_x((T^k \eta) [x]) = \nu_x ((T^k \eta) [x])$
for any $k > 0$ and almost all $x \in X_0$. On the other hand
since $\eta$ is subordinate to $V^-$ and $T$ is contracting on $V^-$,
we have that $\bigvee_{k=1}^\infty T^k \eta$ is the partition into
points. Hence the conditional measures of $\nu$ along $V$ agree with
$\tau$. This proves (b2). 
\qed\medskip

\begin{theorem}
\label{theorem:mt:theorem97}
Let $T=g_s$ denote the time $s$ map of the geodesic flow. 
Assume that $T$ acts ergodically on $(X_0,\nu)$. Let $V^-(x)$ be an
admissible $T$-equivariant system of real-algebraic subsets 
of $W^-(x)$, and let $\Delta(V^-)$
be as in Proposition~\ref{prop:mt:prop96}.  
\begin{itemize}
\item[{\rm (i)}] Suppose $V^-$ has a system of measures $\tau$ with
  the Margulis property, and suppose that for almost all $x$,
the conditional measures of $\nu$ along $V^-[x]$ agree with $\tau(x)$ up to
normalization. Then, $h(T,\nu) \ge s \Delta(V^-)$. 

\item[{\rm (ii)}] Assume that there exists a subset $\Psi \subset X_0$
  with $\nu$-measure $1$ such that $\Psi \cap W^-[x] \subset V^-[x]$
  for every $x \in \Psi$. Then $h(T,\nu) \le s \Delta(V^-)$.

\item[{\rm (iii)}] Assume that there exists a subset $\Psi \subset X_0$
  with $\nu$-measure $1$ such that $\Psi \cap W^-[x] \subset V^-[x]$
  for every $x \in \Psi$. Also assume that $V^-$ has a system
  of measures $\tau$ with the Margulis property,  and that $h(T,\nu) =
  s \Delta(V^-)$. Then, for almost all $x$,
  the conditional measures of $\nu$ along $V^-[x]$ agree with $\tau(x)$ up
  to normalization.  
\end{itemize}
\end{theorem}
\bold{Proof.} According to Proposition~\ref{prop:mt:93}, 
there exists a measurable $T$-invariant partition $\eta$ of $(X_0,\nu)$,
subordinate to $W^-$, such that $H(T\eta \mid \eta)=h(T,\nu)$. By
Lemma~\ref{lemma:embedded:leaves}, we may assume that 
the affine exponential map $W^-(x) \to W^-[x]$ is one-to-one and onto,
and thus $W^-[x]$ has an affine structure. 
Set $\eta'(x)
= V^-[x] \cap \eta[x]$. 

Suppose the assumptions of (i) hold. Then, \mccc{reference?}
\begin{equation}
\label{eq:tmp:h:T:nu:ge:H:T:eta:prime}
h(T,\nu) \ge H(T\eta' \mid \eta'). 
\end{equation}
By Proposition~\ref{prop:mt:prop96} (b1), $H(T \eta' \mid \eta') =
s\Delta(V^-)$. This, together with (\ref{eq:tmp:h:T:nu:ge:H:T:eta:prime})
implies the conclusion of (i). 

Now suppose the assumptions of (ii) or (iii) hold. 
Then $\eta$ and $\eta'$ coincide on $\Psi$,
i.e.\ $\eta[x] \cap \Psi = \eta'[x] \cap \Psi$. Hence $H(T\eta \mid
\eta) = H(T \eta' \mid \eta')$. By Proposition~\ref{prop:mt:93} (iii),
$h(T,\nu) = H(T\eta \mid \eta)$. Using
Proposition~\ref{prop:mt:prop96} (a) we obtain (ii), and using 
Proposition~\ref{prop:mt:prop96} (b2) we obtain (iii). 
\qed\medskip

\section{Semisimplicity of the Lyapunov spectrum}
\label{sec:appendixC}
In this section we work with a bit more generality than we need.
Let $X$ be a space on which $SL(2,\reals)$ acts. 
Let $\mu$ be a compactly supported probability measure on $SL(2,\reals)$
and let $\nu$ be an ergodic $\mu$-stationary probability measure on $X$. 
Let $L$ be a finite dimensional real vector
space, and suppose $A: SL(2,\reals) \cross X \to SL(L)$ is a cocycle,
such that for any $g \in SL(2,\reals)$, 
the map $x \to \log^+\|A(g,x)\|$ is in $L^1(X,\nu)$. 
Let $H'$ be the algebraic hull of the cocycle $A$ (see
\S\ref{sec:subsec:app:forni:cocycle} for the definition).
We may assume that a basis at every point is chosen so that for all $g
\in SL(2,\reals)$ and all $x \in X$, $A(g,x) \in H'$.

\begin{definition}
\label{def:invariant:subspace}
We say that a measurable map $W: X \to L$ is an {\em invariant system of subspaces} for $A(\cdot,
\cdot)$ if for
$\mu$-a.e.\ $g \in SL(2,\reals)$ and $\nu$-a.e.\ $x \in X$,
$A(g,x)W(x) = W(gx)$. 
\end{definition}

\begin{definition}[Strongly Irreducible]
\label{def:strongly:irreducible}
We say that $A$ is {\em strongly irreducible} if on any measurable finite
cover of $X$ there is no nontrivial proper invariant system of
subspaces for $A(\cdot, \cdot)$. 
\end{definition}

\bold{Remark.} If a cocycle is strongly irreducible, then its
algebraic hull is a simple Lie group. 
\medskip




Let $B$ be the space of (one-sided) infinite sequences of
elements of $SL(2,\reals)$. 
We define the measure $\beta$ on $B$ to be $\mu \cross \mu \cdots$. 
Let $\hat{T}: B \cross X \to B \cross X$ be the forward shift, with
$\beta \cross \nu$ as the invariant measure. We denote elements of $B$
by the letter $a$ (following the convention that these refer to
``future'' trajectories). If we write $a = (a_1, a_2, \dots )$ then 
\begin{displaymath}
\hat{T}(a,x) = (Ta, a_1 x)
\end{displaymath}
(and we use the letter $T$ to denote the shift $T(a_1,a_2, \dots) =
(a_2,a_3, \dots)$.) By the Osceledets multiplicative ergodic
theorem, for $\beta \cross \nu$ almost every $(a,x) \in B \cross X$ there
exists a Lyapunov flag 
\begin{equation}
\label{eq:app:forward:Laypunov:flag}
\{0\} = \cV_{\ge k}(a,x) \subset \cV_{\ge k-1}(a,x) \subset \dots 
\subset \cV_{\ge 0}(a,x) = L.  
\end{equation}

\begin{definition}
\label{def:semisimple:lyapunov}
The map $\hat{T}: B \cross X \to B \cross X$ 
has {\em semisimple Lyapunov spectrum} if (after passing to a
measurable finite cover),  
the algebraic hull of the cocycle $\zed \cross (B \cross X) \to SL(L)$
given by 
\begin{displaymath}
(n,a,x) \to A(a_n \dots a_1,x) 
\end{displaymath}
is block-conformal, see
\S\ref{sec:subsec:jordan:form}. 
In other words, $\hat{T}$ has semisimple
Lyapunov spectrum if all the off-diagonal blocks labelled $\ast$ in
(\ref{eq:jordan:block}) are $0$. 
\end{definition}
\medskip

In Appendix~\ref{sec:appendixC} our aim is to prove the following general fact:
\begin{theorem}
\label{theorem:forward:semisimple:lyapunov}
Suppose $A$ is strongly irreducible and $\nu$ is $\mu$-invariant. 
Then $\hat{T}$ has
semisimple Lyapunov spectrum. 
Furthermore, the restriction of $\hat{T}$ to
the top Lyapunov subspace $\cV_{\ge 1}/\cV_{>1}$ consists of a single conformal block, i.e.\
for $\beta \cross \nu$ almost every $(a,x)$ 
there exists an inner product $\langle \cdot,
\cdot \rangle_{a,x}$ on $\cV_{\ge 1}(a,x)/\cV_{> 1}(a,x)$  and a function
$\lambda: B \cross X \to \reals$  such that 
for all $u,v \in \cV_{\ge 1}(a,x)/{\cV_{>1}(a,x)}$, 
\begin{equation}
\label{eq:forward:conformal}
\langle a_1 u, a_1 v\rangle_{(Ta,a x)} = \lambda(a_1 ,x) \langle u,v\rangle_{a,x}.
\end{equation}
If the algebraic hull $H'$ is all of $SL(L)$, then all the Lyapunov
subspaces consist of a single conformal block, i.e.\ for all $1 \le i
\le k-1$ one can define an inner product $\langle \cdot,
\cdot \rangle_{a,x}$ on $\cV_{\ge i}(a,x)/\cV_{>i}(a,x)$ so that
(\ref{eq:forward:conformal}) holds for some function $\lambda = \lambda_i$.
\end{theorem}

\bold{The backwards shift.}
We will actually use the analogue of
Theorem~\ref{theorem:forward:semisimple:lyapunov} for the backwards
shift. 
Let $T: B \cross X \to B \cross X$ be the (backward) shift as in
\S\ref{sec:random:walks}, with
$\beta^X$ as defined in \cite[Lemma~3.1]{Benoist:Quint}
as the invariant measure. By the Osceledets multiplicative ergodic
theorem, for $\beta^X$ almost every $(b,x) \in B \cross X$ there
exists a Lyapunov flag 
\begin{equation}
\label{eq:app:backward:Lyapunov:flag}
\{0\} = \cV_{\le 0}(b,x) \subset \cV_{\le 1}(b,x) \subset \cV_{\le 2}(b,x)
\subset \cV_{\le k}(b,x) = L.  
\end{equation}
We need the following:
\begin{theorem}
\label{theorem:semisimple:lyapunov}
Suppose $A$ is strongly irreducible and $\nu$ is $\mu$-invariant. Then $T$ has
semisimple Lyapunov spectrum. 
Furthermore, the restriction of $T$ to
the top Lyapunov subspace $\cV_{\le 1}$ consists of a single conformal block, i.e.\
for $\beta^X$ almost every $(b,x)$ 
there exists an inner product $\langle \cdot,
\cdot \rangle_{b,x}$ on $\cV_{\le 1}(b,x)$  and a function
$\lambda: B \cross X \to \reals$  such that 
for all $u,v \in \cV_{\le 1}(b,x)$, 
\begin{equation}
\label{eq:conformal}
\langle b_0^{-1} u, b_0^{-1} v\rangle_{(Tb,b_0^{-1} x)} = \lambda(b_0 ,x) \langle u,v\rangle_{b,x}.
\end{equation}
If the algebraic hull $H'$ is all of $SL(L)$, then all the Lyapunov
subspaces consist of a single conformal block, i.e.\ for all $1 \le i
\le k-1$ one can define an inner product $\langle \cdot,
\cdot \rangle_{b,x}$ on $\cV_{\le i}(b,x)/\cV_{< i}(b,x)$ so that
(\ref{eq:conformal}) holds for some function $\lambda = \lambda_i$.
\end{theorem}

\bold{The two-sided shift.}
As in \S\ref{sec:random:walks}, let $\tilde{B}$ be the space of
bi-infinite sequences of elements of $SL(2,\reals)$, and 
we consider the two-sided random walk as a shift map on $\tilde{B}
\cross X$. We abuse notation by using the same letter
  $T$ both for the backwards shift and the bi-infinite shift.
We denote a point in $\tilde{B}$ by $a \vee b$ where $a$
denotes the ``future'' of the trajectory and $b$ denotes the
``past''. Let $\tilde{\beta}^X$ denote the $T$-invariant measure on
$\tilde{B} \cross X$ which projects to the measure $\beta \cross
\nu$ on the  future trajectories, and to the measure $\beta^X$ on the past
trajectories. 
Then, at $\tilde{\beta}^X$ almost all points $(a \vee b, x)$ we have both
the flags (\ref{eq:app:forward:Laypunov:flag}) and
(\ref{eq:app:backward:Lyapunov:flag}). The two flags are generically
in general position (see e.g.\ \cite[Lemma 1.5]{Goldsheid:Margulis})
and thus we can intersect the flags to
define the (shift-invariant) Lyapunov subspaces $\cV_i(a \vee b, x)$ so that
\begin{displaymath}
\cV_{\le i}(b,x) = \bigoplus_{j=1}^i \cV_j(a \vee b,x),\qquad \cV_{\ge
  i}(a,x) =
\bigoplus_{j=i}^m \cV_j(a \vee b, x). 
\end{displaymath}
Then
\begin{equation}
\label{eq:intersection:of:flags}
\cV_{\le i}(b,x)/\cV_{< i}(b,x) \isom \cV_i(a \vee b,x) \isom
\cV_{\ge i}(a,x)/\cV_{> i}(a,x).
\end{equation}

We will prove the following:
\begin{theorem}
\label{theorem:twosided:semisimple:lyapunov}
Suppose $A$ is strongly irreducible and $\nu$ is $\mu$-invariant. 
Then $T$ has semisimple Lyapunov spectrum. 
Furthermore, the restriction of $T$ to
the top Lyapunov subspace $\cV_{\le 1}$ consists of a single conformal block, i.e.\
for $\tilde{\beta}^X$ almost every $(a \vee b,x)$ 
there exists an inner product $\langle \cdot,
\cdot \rangle_{a \vee b,x}$ on $\cV_1(a \vee b,x)$  and a function
$\lambda: \tilde{B} \cross X \to \reals$  such that 
for all $u,v \in \cV_1(a \vee b,x)$, 
\begin{equation}
\label{eq:two:sided:conformal}
\langle a_1 u, a_1 v\rangle_{(T(a \vee b),a_1 x)} = \lambda(a \vee b
,x) \langle u,v\rangle_{a \vee b,x}.
\end{equation}
If the algebraic hull $H'$ is all of $SL(L)$, then all the Lyapunov
subspaces consist of a single conformal block, i.e.\ for all $1 \le i
\le k-1$ one can define an inner product $\langle \cdot,
\cdot \rangle_{a \vee b,x}$ on $\cV_i(b,x)$ so that
(\ref{eq:two:sided:conformal}) holds for some function $\lambda = \lambda_i$.
\end{theorem}

\bold{Remark 1.} The proof of
Theorems~\ref{theorem:forward:semisimple:lyapunov}-\ref{theorem:twosided:semisimple:lyapunov}
we give is essentially taken from
\cite{Goldsheid:Margulis}, and is originally from
\cite{Guivarch:Raugi:Frontiere} and \cite{Guivarch:Raugi:Contraction}. 

For most of the proof, we assume only that $\nu$ is $\mu$-stationary
(and not necessarily $\mu$-invariant). The
exceptions are Lemma~\ref{lemma:non:atomic:lyapunov} and
Claim~\ref{claim:poitwise:martingale}.


We follow \cite{Goldsheid:Margulis} and present the proof of
Theorems~\ref{theorem:forward:semisimple:lyapunov}-\ref{theorem:twosided:semisimple:lyapunov}
for the easier to read 
case where the algebraic hull $H'$ of the cocycle $A$ is
all of $SL(L)$. The general case of semisimple $H'$
  is treated in \cite{Eskin:Matheus:semisimple}. 


\bold{Remark 2.} It is possible to define semisimplicity of the
Lyapunov spectrum in the
context of the action of $g_t = \begin{pmatrix} e^t & 0 \\ 0 & e^{-t}
\end{pmatrix} \subset SL(2,\reals)$ (instead of the random walk). 
Then the analogue of
Theorems~\ref{theorem:forward:semisimple:lyapunov}-\ref{theorem:twosided:semisimple:lyapunov}
remains true; the proof
would use an argument similar to the proof of
Proposition~\ref{prop:sublyapunov:locally:constant}. Since we will
not use this statement we will omit the details.

\subsection{An ergodic lemma}

We recall the following well-known lemma:
\begin{lemma}
\label{lemma:atkinson}
Let $T: \Omega \to \Omega$ be a transformation preserving a 
probability measure $\beta$. Let $F: \Omega \to \reals$ be an 
$L^1$ function. Suppose that for $\beta$-a.e.\ $x \in \Omega$, 
\begin{displaymath}
\liminf \sum_{i=1}^n F(T^i x) = + \infty.
\end{displaymath}
Then $\int_\Omega F \, d \beta > 0$. 
\end{lemma}

\bold{Proof.} This lemma is due to Atkinson \cite{Atkinson}
and Kesten 
\cite{Kesten:sums}. See also
\cite[Lemma 5.3]{Goldsheid:Margulis}, and the
references quoted there. 
\qed\medskip

We will need the following variant:
\begin{lemma}
\label{lemma:K}
Let $T: \Omega \to \Omega$ be a transformation preserving an ergodic  
probability measure $\beta$. Let $F: \Omega \to \reals$ be an 
$L^1$ function. 
Suppose there exists $K' \subset \Omega$ with 
$\beta(K') > 0 $ such that for $\beta$-a.e.\ $x \in \Omega$,
\begin{equation}
\label{eq:assumption}
\liminf \left\{ \sum_{i=1}^n F(T^i x)  \st T^n x \in K' \right\} = + \infty.
\end{equation}
Then $\int_\Omega F \, d \beta > 0$. 
\end{lemma}

\bold{Proof.} After passing to the natural extension, we 
may assume that $T$ is invertible.
We can choose a subset $K \subset K'$ with $\beta(K) >0$, and $C>0$ such that 
for all $x \in K$, we have  
\begin{displaymath}
|F(x)| < C.
\end{displaymath}
Since $K \subset K'$, (\ref{eq:assumption}) holds with $K'$ replaced
by $K$. 

Let 
$A_{-1}=\{x \st x \not \in K\},$ $A_{0}=\{x \st x \in K,\, Tx \in K\},$ and for $n\geq 0,$
\begin{displaymath}
A_{n+1}=\{ x \st  x \in K, \ Tx\not \in K, \dots, \ T^{n}x \not \in
K, \  T^{n+1}x \in K\}.
\end{displaymath}
Also let  $A=\bigsqcup\limits_{n=-1}^{\infty}  A_{n}$.
Note that by the ergodicity of $T$, for almost every $x \in \Omega, $
$$|\{i \st i\geq 0, T^{i}(x) \in K\}| =\infty.             \;\;\;\; (*).$$

Define $G: \Omega \rightarrow \reals$ defined on $A$ (which has full measure) by  
\begin{itemize}
\item $G(x)=0$ if $x \in A_{-1}$.
\item $G(x)=F(x)$ if $x \in A_{0}$.
\item $G(x)=F(x)+F(Tx)+\cdots+F(T^{n} x)$ if $x \in A_{n+1}.$ 
\end{itemize}
We now claim the following hold:
\begin{enumerate}
\item For almost every $x \in \Omega$ we have 
\begin{equation}
\label{eq:one}
\lim_{n\rightarrow \infty} G(x)+G(Tx)+....+G(T^{n}x)=\infty.
\end{equation}
\item $\int_{\Omega} |G| \, d \beta \leq \int_{\Omega} |F| \, d\beta <\infty .$
\item $\int_{\Omega} G(x)\,  d\beta(x) =  \int_{\Omega} F(x) \, d\beta(x).$
\end{enumerate}

\bold{Proof of (1).}
Note that almost every $x \in \Omega$ satisfies
(\ref{eq:assumption}) (with $K'$ replaced by $K$). Also, we have, 
$$G(x)+G(Tx)+....+G(T^{n}x)=\sum_{i=m_{0}}^{m-1} F(T^{i}x), $$
where $m_{0}=\inf \{ k \st T^{k}x \in K\},$ and $m=\inf\{k \st k\geq
n, T^{k}x \in K\}.$
Thus, 
\begin{displaymath}
\sum_{j=0}^n G(T^j x) = \sum_{i=1}^m F(T^i x) - \sum_{i=0}^{m_0-1}
F(T^i x) - F(T^m x).
\end{displaymath}
Since $m_0$ is independent of $n$, $T^m x \in K$ and
for every $x \in K$, we have $|F(x)| <C$, the equation (\ref{eq:assumption})
implies (\ref{eq:one}).
\qed\medskip

\bold{Proof of (3) assuming (2).}
By the definition of $G$ we can use the dominated convergence theorem, and get that 
$$\int_{\Omega}G \, d\beta= \int_K F \, d\beta + \sum_{i=1}^{\infty}
\int_{ A^{i}} F(T^{i} x) \, d\beta(x)$$
where $A^{i}=\bigcup_{j\geq i} A_{i}.$ Then $$T^{i} A^{i}=T^{i}K - (K\cup\cdots T^{i-1}K).$$ 
 Also $K\cup \bigcup_{i=1}^{\infty} T^{i} A^{i}$ has full measure in
 $\Omega$, and for $i\ne j$,
 $T^{i}A^{i} \cap T^{j} A^{j}$ and $K \cap T^{i}A^{i}$ have measure zero. 
 Note that $A^{i}= T^{-i}(T^{i}A^{i}).$
Since $\beta$ is $T$ invariant, we have 
$$  \int_{A^{i}} F(T^{i}x) \, d\beta(x) =\int_{T^{i} A^{i}} F(x) \, d\beta(x), $$
and hence
$$\int_{\Omega} G \, d\beta= \int_K F \, d\beta + \sum_{i=1}^{\infty}
\int_{T^{i} A^{i}} F(x) \, d\beta(x)=
\int_{\Omega} F \, d\beta.$$
\qed\medskip

\bold{Proof of (2).} This follows by applying (3) to $|F|$ instead of
$F$, and then using the triangle inequality. 
\qed\medskip
 
\bold{Proof of Lemma~\ref{lemma:K}.}
Now by (1), and (2), the function $G$ satisfies the assumptions of Lemma~\ref{lemma:atkinson}. Hence we have
$\int_{\Omega} F \, d\beta=\int _{\Omega} G \, d\beta > 0. $
\qed\medskip

\subsection{A zero-one law}
\label{sec:subsec:zero:one}

\mccc{Recall the setting here.}

\begin{lemma}
\label{lemma:no:positive:subharmonic}
Suppose $h$ is a bounded non-negative $\mu$-subharmonic function, i.e.\ for
$\nu$-almost all $x \in X$, 
\begin{equation}
\label{eq:def:subharmonic}
h(x) \le \int_G h(gx) \, d\mu(g). 
\end{equation}
Then $h$ is constant $\nu$-almost everywhere. 
\end{lemma}

\bold{Proof.}
By the random ergodic theorem \cite[Theorem 3.1]{Furman:Survey}, 
for $\nu$-almost all
$x \in X$, 
\begin{displaymath}
\lim_{N \to \infty} \frac{1}{N} \sum_{n=0}^{N-1} \int_G h(g x) \,
d\mu^n(g) = \int_X h \, d\nu
\end{displaymath}
Therefore, by (\ref{eq:def:subharmonic}), for $\nu$-almost all $x \in X$,
\begin{equation}
\label{eq:h:smaller:then:average}
h(x) \le \int_X h \, d\nu.
\end{equation}
Let $s_0 \ge 0$ denote the essential supremum of $h$, i.e.\ 
\begin{displaymath}
s_0 =  \inf \{ s \in \reals \st \nu(\{ h > s \}) = 0 \}.
\end{displaymath}
Suppose $\epsilon > 0$ is arbitrary. We can pick $x \in X$ such that
(\ref{eq:h:smaller:then:average}) holds and $h(x) > s_0 -
\epsilon$. Then, 
\begin{displaymath}
s_0 - \epsilon \le h(x) \le \int_X h \,d\nu \le s_0.
\end{displaymath}
Since $\epsilon > 0$ is arbitrary, $\int_X h \, d\nu = s_0$. Thus
$h(x)=s_0$ for $\nu$-almost all $x$. 
\qed\medskip


Let $\nu$ be an ergodic stationary measure on $X$. Fix $1 \le s <
\dim(L)$, and let $Gr_s$ denote the Grassmannian of $s$-dimensional subspaces
in $L$. Let $\hat{X} = X \cross Gr_s$. We then have an
action of $SL(2,\reals)$ on $\hat{X}$, by 
\begin{displaymath}
g \cdot (x, W) = (gx, A(g,x) W). 
\end{displaymath}
Let $\hat{\nu}$ be a $\mu$-stationary measure on $\hat{X}$
which projects to $\nu$ under the natural map $\hat{X} \to X$. 
We may write
\begin{displaymath}
d\hat{\nu}(x,U) = d\nu(x) \, d\eta_x(U),
\end{displaymath}
where $\eta_x$ is a measure on $Gr_s$. 

Let $m = \dim(L)$. For a subspace $W$ of $L$, let 
\begin{displaymath}
I(W) = \{ U \in Gr_s \st \dim(U \cap W) > \max(0,m - \dim(U) - \dim(W))\}
\end{displaymath}
Then $U \in I(W)$ if and only if $U$ and $W$ intersect more than
general position subspaces of dimension $\dim(U)$ and $\dim(W)$. 

\begin{lemma}(cf.\ \cite[Lemma 4.2]{Goldsheid:Margulis},
\cite[Theorem~2.6]{Guivarch:Raugi:Frontiere})
\label{lemma:non:atomic:lyapunov}
$ $
\begin{itemize}
\item[(i)] Suppose the cocycle is strongly irreducible on $L$. 
  Then for almost all $x \in X$, and any
    $1$-dimensional subspace $W_x \subset L$, $\eta_x(I(W_x)) = 0$. 
\item[(ii)] Suppose the algebraic hull $H'$ of the cocycle is
  $SL(L)$. Then for almost all $x \in X$, for any nontrivial proper
  subspace $W_x \subset L$, $\eta_x(I(W_x)) = 0$.
\end{itemize}
\end{lemma}

\bold{Proof of Lemma~\ref{lemma:non:atomic:lyapunov}.}  We give the
proof under the extra assumption that $\nu$ is $\mu$-invariant (and
not just $\mu$-stationary). The general case is proved in
\cite{Eskin:Matheus:semisimple}.

Suppose there exists  
a subset $E \subset X$ with $\nu(E) > 0$ 
and for all $x \in E$,
a nontrivial subspace $W_x \subset L$ such that
$\eta_x(I(W_x)) > 0$.  
Let $\vec{W} = (W_1, \dots, W_k)$ denote a finite
collection of subspaces of $L$. 
If the assumptions of (i) hold, we are requiring the $W_i$ to be
one-dimensional; if the assumptions of (ii) hold, the $W_i$ are
allowed to be any dimension. Write
\begin{displaymath}
I(\vec{W}) = I(W_1) \cap \dots \cap I(W_k). 
\end{displaymath}
For $x \in E$, let $\cS_x$ denote the set of $I(\vec{W}_x)$ such that 
for any $\vec{W}_x'$ 
so that $I(\vec{W}_x')$ is a proper subset of $I(\vec{W}_x)$, we have 
$\nu_x(I(\vec{W}_x')) = 0$. For $x \in E$, $\cS_x$ is
non-trivial since the subsets $I(\vec{W})$ are algebraic and thus
there cannot  be an infinite descending chain of them.
For $\vec{W} \in \cS_x$, let 
\begin{displaymath}
f_{I(\vec{W})}(x) = \eta_x(I(\vec{W})).
\end{displaymath}
Since $\hat{\nu}$ is $\mu$-stationary and $\nu$ is
  assumed to be $\mu$-invariant, we have 
\begin{equation}
\label{eq:fW:onestep}
f_{I(\vec{W})}(x) = \int_G f_{I(A(g,x) \vec{W})}(g x) \, d\mu(g)
\end{equation}
Let $\cS(x) = \{ I(\vec{W}) \in \cS_x
\st f_{I(\vec{W})}(x) > 0 \}$. 
Then, for $I(\vec{W}_1) \in \cS(x)$, $I(\vec{W}_2) \in \cS(x)$, 
\begin{displaymath}
\eta_x(I(\vec{W}_1) \cap I(\vec{W}_2)) = 0.
\end{displaymath}
Thus
\begin{displaymath}
\sum_{I(\vec{W}) \in \cS(x)} f_{I(\vec W)}(x) \le 1. 
\end{displaymath}
Therefore $\cS(x)$ is at most countable. Let 
\begin{equation}
\label{eq:fmax}
f(x) = \max_{I(\vec{W}) \in \cS(x)} f_{I(\vec W)}(x). 
\end{equation}
Applying (\ref{eq:fW:onestep}) to some $I(\vec{W})$ for which the max is
achieved, we get
\begin{displaymath}
f(x) \le \int_G f(g x) \, d\mu(g)
\end{displaymath}
i.e.\ $f$ is a subharmonic function on $X$. By
Lemma~\ref{lemma:no:positive:subharmonic}, 
$f$ is constant almost everywhere. Now substituting again
into (\ref{eq:fW:onestep}) we get that the cocycle $A$ permutes
the finite set of $I(\vec W)$ where the maximum (\ref{eq:fmax}) is achieved.
Therefore the same is true for the algebraic hull $H'$. 
If the assumptions of (ii) hold, this is a contradiction since $H'$ acts
transitively on subspaces of $L$. If the assumptions of (i) hold then, 
for $\vec{W} = (W_1, \dots, W_k)$, 
since the $W_i$ are $1$-dimensional, we have
\begin{displaymath}
I(\vec{W}) \equiv I(W_1) \cap \dots I(W_k) = \{
\text{subspaces $M \subset L$ such that $W_1 + \dots + W_k \subset M$. } \}
\end{displaymath}
Since $H'$ must permute some finite set of $I(\vec{W})$ it must thus 
permute a finite set of subspaces of $L$ which contradicts the strong
irreducibility assumption. 
\qed\medskip

\subsection{Proof of Theorem~\ref{theorem:twosided:semisimple:lyapunov}}
Recall that we are  assuming that the algebraic hull of the cocycle is
$SL(L)$ for some vector space $L$. Let $m = \dim L$.

\begin{definition}[{\bf $(\epsilon,\delta)$-regular}]
\label{def:epsilon:regular}
Suppose $\epsilon > 0$ and $\delta > 0$ are fixed. 
A measure $\eta$ on $Gr_k(L)$ is {\em 
  $(\epsilon,\delta)$-regular} if for 
any subspace $U$ of $L$, 
\begin{displaymath}
\eta(Nbhd_\epsilon(I(U))) < \delta.
\end{displaymath}
\end{definition}

\begin{lemma}
\label{lemma:strongly:regular}
Suppose $g_n \in GL(L)$ is a sequence of linear transformations, and
$\eta_n$ is a sequence of uniformly $(\epsilon,\delta)$-regular
measures on $Gr_k(L)$ for some $k$. Suppose $\delta
  \ll 1$.  Write
\begin{displaymath}
g_n = K(n) D(n) K'(n), 
\end{displaymath}
where $K(n)$ and $K'(n)$ are orthogonal relative to the standard basis
$\{e_1, \dots e_{m} \}$, and $D(n) = \diag(d_1(n), \dots, d_m(n) \}$
with $d_1(n) \ge \dots \ge d_{m}(n)$.
\begin{itemize}
\item[{\rm (a)}] Suppose 
\begin{equation}
\label{eq:dk:dkplus1}
\frac{d_k(n)}{d_{k+1}(n)} \to \infty
\end{equation}
Then, for any subsequential limit $\lambda$ of $g_n \eta_n$ there
exists a subspace $W \in Gr_k(L)$ such that
\begin{equation}
\label{eq:tmp:K(n):Span:e1es}
K(n) \Span\{ e_{1}, \dots, e_k \} \to W,
\end{equation}
and $\lambda(\{ W\}) \ge 1-\delta$. 
\item[{\rm (b)}] Suppose  $g_n \eta_n \to \lambda$ where $\lambda$ is
  some measure on $Gr_k(L)$. Suppose also that there exists
a subspace $W \in Gr_k(L)$ such that $\lambda(\{W\}) > 5\delta$.
Then, as $n \to \infty$, 
(\ref{eq:dk:dkplus1}) holds. As a consequence, by part (a),
(\ref{eq:tmp:K(n):Span:e1es}) holds and  $\lambda(\{W\}) \ge
1-\delta$. \mcc{check this proof}
\end{itemize}
\end{lemma}

\bold{Proof of (a).} This statement is standard. Suppose $g_n \eta_n
\to \lambda$. Without loss of
generality, $K'(n)$ is the identity (or else we replace $\eta_n$ by
$K'(n) \eta_n$). By our assumptions, for $j_1 < \dots < j_k$,
\begin{displaymath}
\frac{\|g_n(e_{j_1} \wedge \dots \wedge e_{j_k}) \|}{\|g_n(e_1 \wedge \dots
  \wedge e_k)\|} \to 0 \quad \text{unless $j_i = i$ for $1 \le i \le k$.}  
\end{displaymath}
Therefore, if $U \not\in I(\Span\{e_{k+1},\dots, e_m\})$, 
\begin{displaymath}
d(g_n U, K(n)\Span\{e_1, \dots, e_k\}) \to 0,
\end{displaymath}
where $d( \cdot, \cdot)$ denotes some distance in $Gr_k(L)$. 
After passing to a further subsequence, we may assume that for some $W
\in Gr_k(L)$, (\ref{eq:tmp:K(n):Span:e1es}) holds. 
It follows from the $(\epsilon,\delta)$-regularity of
$\eta_n$ that $\lambda(W) \ge 1-\delta$. Since $\delta < 1/2$, $W$ is
uniquely determined by $\lambda$, 
and therefore (\ref{eq:tmp:K(n):Span:e1es}) holds without
passing to a further subsequence (but only assuming $g_n \eta_n \to
\lambda$).  

\bold{Proof of (b).} This is similar to 
\cite[Lemma 3.9]{Goldsheid:Margulis}. 
Suppose $d_k(n)/d_{k+1}(n)$ does not go to $\infty$. 
Then, there is a subsequence of the $g_n$
(which we again denote by $g_n$)  that $K(n) \to K_*$ and 
that for every $j$, either $d_j(n)/d_{j+1}(n)$ converges as $n \to
\infty$ or
$d_{j}(n)/d_{j+1}(n) \to \infty$ as $n \to \infty$.  Also without loss
of generality we may assume that $K'(n)$ is the identity (or else we
replace $\eta_n$ by $K'(n) \eta_n$). 

Let $1 \le s \le k < r \le m$ 
be such that $s$ is as small as possible, $r$ is as
large as possible, and $d_j(n)/d_{j+1}(n)$ is bounded for $s \le j \le r-1$.
Then, for $j_1 < \dots < j_k$,
\begin{multline}
\label{eq:tmp:gn:ej1:ejk}
\frac{\|g_n(e_{j_1} \wedge \dots \wedge e_{j_k}) \|}{\|g_n(e_1 \wedge \dots
  \wedge e_k)\|} \to 0 \quad \text{unless $j_i = i$ for $1 \le i \le
  s-1$} \\ \text{and $s \le j_i \le r$ for $s \le i \le k$.}  
\end{multline}
Let 
\begin{displaymath}
V_- = \Span\{e_1,\dots, e_{s-1}\}, \qquad V_+ = \Span\{e_1,\dots,
e_r\}. 
\end{displaymath}
Let $D_* = \diag(d_*(1), \dots, d_*(m))$ 
be any diagonal matrix such that for $s \le j \le r-1$, 
\begin{displaymath}
d_*(j)/d_*(j+1) = \lim_{n \to \infty} d_j(n)/d_{j+1}(n). 
\end{displaymath}
Then, in view of (\ref{eq:tmp:gn:ej1:ejk}), for $U$ such that $U
\not\in I(V_+^\perp) \cup I(V_-^\perp)$, if along some subsequence
$g_n U \to U'$, we have
\begin{displaymath}
K_* V_- \subset U' \subset K_* V_+.
\end{displaymath}
Therefore, we must have $V_- \subset K_*^{-1}
W \subset V_+$. Furthermore, for $U \not \in I(V_+^\perp) \cup I(V_-^\perp)$,
\begin{displaymath}
\text{if $g_n U \to W$ then $U \in I(D_*^{-1} K_*^{-1}W\cap
V_-^\perp+V_+^\perp)$.} 
\end{displaymath}
\mcc{explain and check} 
But, since $\eta_n$ is $(\epsilon,\delta)$-regular,
\begin{displaymath}
\eta_n(Nbhd_\epsilon(I(V_+^\perp) \cup I(V_-^\perp) \cup
I(D_*^{-1} K_*^{-1}W\cap V_-^\perp+V_+^\perp))) < 3\delta.
\end{displaymath}
Therefore $\lambda(\{W\}) < 3\delta$
which is a contradiction. Thus $d_k(n)/d_{k+1}(n) \to \infty$. 
Now by part (a)
(\ref{eq:tmp:K(n):Span:e1es}) holds, and $\lambda(\{W\}) \ge 1-\delta$.
\qed\medskip

Let $\cF = \cF(L)$ denote the space of full flags on $L$. Let $\hat{X} = X
\cross \cF$. 
The cocycle
$A$ satisfies the cocycle relation
\begin{displaymath}
A(g_1 g_2,x) = A(g_1, g_2 x) A(g_2,x). 
\end{displaymath}
The group $SL(2,\reals)$ acts on the space $\hat{X}$ by
\begin{equation}
\label{eq:straight:action}
g \cdot (x,f) = (g x, A(g,x) f). 
\end{equation}


Let $\hat{\nu}$ be an ergodic $\mu$-stationary measure on $\hat{X}$ 
which projects to $\nu$ under the natural map $\hat{X} \to X$.
(Note there is always at least one such: one chooses $\hat{\nu}$ to be
an extreme point among the measures which project to $\nu$. If
$\hat{\nu} = \hat{\nu}_1 + \hat{\nu}_2$ where the $\hat{\nu}_i$ are
$\mu$-stationary measures then $\nu = \pi_*(\hat{\nu}) =
\pi_*(\hat{\nu}_1) + \pi_*(\hat{\nu}_2)$. Since $\nu$ is
$\mu$-ergodic, this implies that $\pi_*(\hat{\nu}_1) =
\pi_*(\hat{\nu}_2) = \nu$, hence the $\hat{\nu}_i$ also project to
$\nu$. Since $\hat{\nu}$ is an extreme point among such measures, 
we must have $\hat{\nu}_1 = \hat{\nu}_2 = \hat{\nu}$. Thus $\hat{\nu}$
is $\mu$-ergodic.)

\begin{lemma}[Furstenberg]
\label{lemma:furstenberg:formula}
For $1 \le s \le \dim L$, let $\bar{\sigma}_s:
SL(2,\reals) \cross \hat{X}\to \reals$ 
be given by
\begin{displaymath}
\bar{\sigma}_s(g,x,f) = \log \frac{\| A(g,x) \xi_s(f)
  \|}{\|\xi_s(f)\|} 
\end{displaymath}
where  $\xi_s(f)$ is the $s$-dimensional
component of the flag $f$. (The norms in the above equation are on
$\bigwedge^s(V)$, and here and in the following we make sense of such expressions by picking the same basis for the $\xi_s(f)$ in the
numerator and denominator). 
Then,  we have
\begin{displaymath}
\lambda_1 + \dots + \lambda_s = \int_{SL(2,\reals)} 
\int_{\hat{X}} \bar{\sigma}_s( g, x, f) \,
d\hat{\nu}(x,f) \, d\mu(g). 
\end{displaymath}
where $\lambda_i$ is the $i$'th Lyapunov exponent of the cocycle
$A$. 
\end{lemma}

\bold{Proof.} See the proof of \cite[Lemma
5.2]{Goldsheid:Margulis}.
\qed\medskip

We may disintegrate
\begin{displaymath}
d\hat{\nu}(x,f) = d\nu(x) \, d\eta_x(f). 
\end{displaymath}
Note that Lemma~\ref{lemma:non:atomic:lyapunov} applies to the
projections of the measures $\eta_x$ to the various Grassmannians which
are components of $\cF$. 

For $a \in \tilde{B}$, let the measures $\nu_a$,
$\hat{\nu}_a$ be as
defined in \cite[Lemma 3.2]{Benoist:Quint}, i.e.\ 
\begin{displaymath}
\nu_a = \lim_{n \to \infty} (a_n \dots a_1)^{-1}_* \nu
\end{displaymath}
\begin{displaymath}
\hat{\nu}_a = \lim_{n \to \infty} (a_n \dots a_1)^{-1}_*
\hat{\nu}.
\end{displaymath}
The limits exist by the martingale convergence theorem. 
We disintegrate
\begin{displaymath}
d\hat{\nu}_a(x,f) = d\nu_a(x) \, d\eta_{a,x}(f).  
\end{displaymath}
For $1 \le k \le m$, let $\eta_x^k = (\xi_k)_* \eta_x$ and 
$\eta_{a,x}^k = (\xi_k)_* \eta_{a,x}$, 
where $\xi_k: \cF(L) \to Gr_k(L)$ is the natural projection.  
Then, $\eta_x^k$ and
$\eta_{a,x}^k$ are measures on $Gr_k(L)$. 

\begin{claim}
\label{claim:poitwise:martingale}
On a set of $\beta \cross \nu$ full measure,
\begin{displaymath}
\lim_{n \to \infty} (a_n \dots a_1)_*^{-1} \eta_{a_n \dots a_1 x} = \eta_{a,x}.
\end{displaymath}
Equivalently, using (\ref{eq:straight:action}), 
\begin{displaymath}
\lim_{n \to \infty} A((a_n \dots a_1)^{-1}, a_n \dots a_1 x) 
\eta_{a_n \dots a_1 x} = \eta_{a,x}.
\end{displaymath}
\end{claim}
\bold{Proof of claim.} In this claim, we use the invariance of
$\nu$. 
Let $C \subset X$ and $D \subset \cF$ be
measurable,  and let $\chi_C$ denote the characteristic
functions of $C$. 
Recall that $d\hat{\nu}(x,z)=d\nu(x) d\eta_x(z)$ is $\mu$-stationary, so that
\begin{eqnarray*}
\int_C \eta_x(D) \, d\nu(x) &=& \hat{\nu}(C\times D) = (\mu\ast\hat{\nu})(C\times D) \\
&=& \int \chi_C(g y) A(g,y)\eta_y(D) \, d\nu(y) \, d\mu(g) \\
&=& \int \chi_C(x) A(g, g^{-1}x)\eta_{g^{-1}x}(D) d\nu(x) d\mu(g) \\
&=& \int_C \left(\int_G A(g, g^{-1}x)\eta_{g^{-1}x}(D) \,
  d\mu(g)\right) \, d\nu(x)
\end{eqnarray*}
Since $C$ and $D$ are arbitrary, we see that
\begin{displaymath}
\eta_x = \int_G A(g, g^{-1}x)\eta_{g^{-1}x} \, d\mu(g)
\end{displaymath}
Therefore, (replacing $x$ by $a_{n-1} \dots a_1 x$ and $g$ by
$a_n^{-1}$), we have
\begin{displaymath}
\eta_{a_{n-1} \dots a_1 x} = \int_G A(a_n^{-1}, a_n \dots a_1 x)
\eta_{a_n \dots a_1 x} \, d\mu(a_n). 
\end{displaymath}
Multiplying both sides on the left by $A((a_{n-1} \dots a_1)^{-1},
a_{n-1} \dots a_1 x)$ and using the cocycle identity
\begin{displaymath}
A((a_n \dots
a_1)^{-1}, a_n \dots a_1 x) = 
A((a_{n-1} \dots a_1)^{-1},
a_{n-1} \dots a_1 x) A(a_n^{-1}, a_n \dots a_1 x), 
\end{displaymath}
we get
\begin{multline}
\label{eq:etax-stat}
A((a_{n-1} \dots a_1)^{-1},
a_{n-1} \dots a_1 x) 
\eta_{a_{n-1} \dots a_1 x} = \\ = \int_G A((a_n \dots
a_1)^{-1}, a_n \dots a_1 x) \eta_{a_n \dots a_1 x} \, d\mu(a_n). 
\end{multline}
In view of (\ref{eq:etax-stat}), the expression
\begin{displaymath}
A((a_n \dots a_1)^{-1}, a_n \dots a_1 x) 
\eta_{a_n \dots a_1 x}
\end{displaymath}
is a (measure-valued) martingale. Therefore, the claim follows from
the martingale convergence theorem.
\qed\medskip

If the Lyapunov spectrum is simple, we expect the measures $\eta_{a,x}$
to be supported at one point. In the general case, let
\begin{displaymath}
\lambda_1 \ge \lambda_2 \ge \dots \ge \lambda_m
\end{displaymath}
denote the Lyapunov exponents, and let 
\begin{displaymath}
I = \{ 1 \le r \le m-1 \st
\lambda_r = \lambda_{r+1}\}. 
\end{displaymath}
Then, by the multiplicative ergodic
theorem, Lemma~\ref{lemma:non:atomic:lyapunov} and
Lemma~\ref{lemma:strongly:regular} (a), for $r \not\in I$, 
we have $\eta_{a,x}^{m-r}$ is supported
at one point.  (This point is the part of the flag
(\ref{eq:app:forward:Laypunov:flag}) corresponding to the Lyapunov
exponents $\lambda_{r+1}, \dots, \lambda_m$.)

\begin{claim}
\label{claim:non:singular}
For any $r \in I$ and $\beta \cross \nu$-almost all $(a,x)$,  for any
subspace $W(a,x) \in Gr_{m-r}(L)$, we have  
$\eta_{a,x}^{m-r}(\{W(a,x)\}) = 0$. 
\end{claim}

\bold{Proof of claim.} 
Suppose there exists $\delta > 0$ so that 
for some $r \in I$ for a set $(a,x)$ of positive measure, 
there exists $W(a,x)\in Gr_{m-r}(L)$ with
$\eta_{a,x}^r(\{W(a,x)\}) > \delta$. 
Then this happens for a subset of full measure
by ergodicity. 

Note that by the cocycle relation, 
\begin{displaymath}
A(g^{-1} ,g x) = A(g,x)^{-1}. 
\end{displaymath}
Therefore, 
\begin{displaymath}
A((a_n \dots a_1)^{-1}, a_n \dots a_1 x) = A(a_n \dots a_1,x)^{-1}.
\end{displaymath}
Hence, on a set of $\beta \cross \nu$-full measure, 
\begin{displaymath}
\lim_{n \to \infty} A(a_n \dots a_1,x)^{-1} \eta_{a_n \dots a_1 x} = \eta_{a,x}.
\end{displaymath}

In view of Lemma~\ref{lemma:non:atomic:lyapunov} (cf.\ the proof of
Lemma~\ref{lemma:zero:one:law:star}),
there exists $\epsilon > 0$ and a
compact  $\cK_{\delta} \subset X$ with $\nu(\cK_{\delta}) > 1-\delta$
such that the family of
measures $\{\eta_x \}_{x \in \cK_{\delta}}$ is uniformly
$(\epsilon,\delta/5)$-regular. 
Let 
\begin{displaymath}
\cN_\delta(a,x) = \{ n \in \natls \st a_n \dots a_1 x \in
\cK_{\delta} \}. 
\end{displaymath}
Write 
\begin{equation}
\label{eq:decomp:Ainverse}
A( a_n \dots a_1, x)^{-1} = K_n(a,x) D_n(a,x) K_n'(a,x)
\end{equation}
where $K_n$ and $K_n'$ are orthogonal, and $D_n$ is diagonal with
non-increasing entries. We also
write
\begin{equation}
\label{eq:decomp:A:noinverse}
A( a_n \dots a_1, x) = \bar{K}_n(a,x) \bar{D}_n(a,x) \bar{K}_n'(a,x).
\end{equation}
where $\bar{K}_n$ and $\bar{K}_n'$ are orthogonal, and $\bar{D}_n$ 
is diagonal with non-increasing entries.
Let $d_1(n,a,x)
\ge ... \ge d_m(n,a,x)$ be the entries of $D_n(a,x)$, and let
$\bar{d}_1(n,a,x) \ge \bar{d}_2(n,a,x) \ge \bar{d}_m(n,a,x)$ be the
entries of $\bar{D}_n(a,x)$. Then,
\begin{multline}
\label{eq:d:dbar:conversion}
\bar{d}_j(n,a,x) = d_{m+1-j}^{-1}(n,a,x), \\ \bar{K}'_n(a,x) = w_0
K_n(a,x)^{-1} w_0^{-1}, \quad \bar{K}_n(a,x) = w_0 K'_n(a,x)^{-1} w_0, 
\end{multline}
where $w_0 = w_0^{-1}$ is the permutation matrix mapping $e_j$ to $e_{m+1-j}$. 
Then, by Lemma~\ref{lemma:strongly:regular} (b), 
for $\beta \cross \nu$ almost all $(a,x)$, 
  $\eta^{m-r}_{a,x}(\{W(a,x)\}) \ge 1-\delta$  (and thus $W(a,x)$ is
  unique) and as $n \to
\infty$ along $\cN_\delta(a,x)$ we have: 
\begin{displaymath}
d_{m-r}(n,a,x)/d_{m+1-r}(n,a,x) \to \infty,  
\end{displaymath}
and
\begin{equation}
\label{eq:K:converges}
K_n(a,x) \Span\{e_{1}, \dots, e_{m-r}\} \to W(a,x),
\end{equation}
where the $e_i$ are the standard basis for
$L$.  Then, by (\ref{eq:d:dbar:conversion}), 
\begin{equation}
\label{eq:drds}
\bar{d}_{r}(n,a,x)/\bar{d}_{r+1}(n,a,x) \to \infty,
\end{equation}
and
\begin{displaymath}
\bar{K}_n'(a,x)^{-1} \Span\{e_{r+1}, \dots, e_{m}\} \to w_0 W(a,x)
\end{displaymath}
Therefore for any $\epsilon_1 > 0$ there exists a subset
$H_{\epsilon_1} \subset B 
\cross X$ of $\beta \cross \nu$-measure at least $1-\epsilon_1$ such that the convergence in 
(\ref{eq:drds}) and (\ref{eq:K:converges})
is uniform over $(a,x) \in H_{\epsilon_1}$.
Hence there exists $M > 0$ such that for any $(a,x) \in H_{\epsilon_1}$, 
and any $n \in \cN_\delta(a,x)$ with $n > M$, 
\begin{equation}
\label{eq:betaX:small:nbhd:er}
\bar{K}'_n(a,x)^{-1} \Span\{e_{r+1}, \dots, e_m\} \in
Nbhd_{\epsilon_1}(w_0 W(a,x)).
\end{equation}
By Lemma~\ref{lemma:non:atomic:lyapunov} (cf.\ the proof of
Lemma~\ref{lemma:zero:one:law:star}) 
there exists a subset $H''_{\epsilon_1} \subset X$ with $\nu(H''_{\epsilon_1}) >
1-c_2(\epsilon_1)$ with $c_2(\epsilon_1) \to 0$ as $\epsilon_1 \to 0$ such
that for all $x \in H''_{\epsilon_1}$, and any $U \in Gr_{m-r}(L)$,
\begin{displaymath}
\eta^r_x( Nbhd_{2\epsilon_1}(I(U))) < c_3(\epsilon_1),
\end{displaymath}
where $c_3(\epsilon_1) \to 0$ as $\epsilon_1 \to 0$. 
Let
\begin{equation}
H'_{\epsilon_1} = \{ (a,x,f) \st (a,x) \in H_{\epsilon_1}, \quad x \in H''_{\epsilon_1} \quad \text{and}
\quad d(\xi_r(f),
I(w_0 W(a,x))) > 2 \epsilon_1 \}. 
\end{equation}
Then, $(\beta \cross \hat{\nu})(H'_{\epsilon_1}) > 1-\epsilon_1 - c_2(\epsilon_1) -
c_3(\epsilon_1)$, hence $(\beta \cross \hat{\nu})(H'_{\epsilon_1}) \to
1$ as $\epsilon_1 \to 0$. Furthermore, by
(\ref{eq:betaX:small:nbhd:er}) 
and the definition of $H'_{\epsilon_1}$, for $(a,x,f) \in
H'_{\epsilon_1}$ and $n \in N_\delta(a,x)$ with $n > M$, we have
\begin{displaymath}
d(\xi_r(f), I(\bar{K}'_n(a,x)^{-1} \Span\{e_{r+1}, \dots, e_m\})) > \epsilon_1.
\end{displaymath}
Therefore, in view of (\ref{eq:decomp:A:noinverse}) 
there exists $C = C(\epsilon_1)$, such that
for any $(a,x,f) \in H'_{\epsilon_1}$, any $n \in N_\delta(a,x)$ with $n > M$, 
\begin{equation}
\label{eq:r:two:sided:bound}
C  > \frac{\|A(a_n \dots a_1 , x) \xi_r(f)\|}{\|\xi_r(f)\|} \prod_{i=1}^r
\bar{d}_i(n,a,x)^{-1} > \frac{1}{C}, 
\end{equation}
(c.f  \cite[Lemma 5.1]{Goldsheid:Margulis}). 
Note that for all $(a,x,f) \in B \cross \hat{X}$, all $n \in \natls$ and  
$j=r-1$ or $j = r+1$ we have
\begin{equation}
\label{eq:r:plus:minus:one:onesided:bound}
\frac{\|A(a_n \dots a_1, x) \xi_{j}(f)\|}{\|\xi_r(f)\|} \le \|A(a_n \dots
a_1,x)\|_{\bigwedge^j(L)}  \le 
\prod_{i=1}^{j} \bar{d}_i(n,a,x).
\end{equation}
Then, in view of
(\ref{eq:r:two:sided:bound}) and
(\ref{eq:r:plus:minus:one:onesided:bound}), for all $(a,x,f) \in
H'_{\epsilon_1}$, as $n \to \infty$ in $\cN_\delta(a,x)$,  
\begin{multline}
\label{eq:infinite:limit}
\log \frac{\|(A(a_n \dots a_1, x)) \xi_r(f) 
  \|^2}{\|\xi_r(f)\|^2} \frac{\|\xi_{r-1}(f)\|} 
{\|(A(a_n \dots a_1, x)) \xi_{r-1}(f) \|}
\frac{\|\xi_{r+1}(f)\|}
{\|(A(a_n \dots a_1, x)) \xi_{r+1}(f) \| } \ge \\ \ge
\log \frac{\bar{d}_r(n,a,x)}{\bar{d}_{r+1}(n,a,x)} \to \infty 
\end{multline}
Since $(\beta \cross \hat{\nu})(H'_{\epsilon_1}) \to 1$ as $\epsilon_1
\to 0$, (\ref{eq:infinite:limit}) holds as $n
\to \infty$ along $\cN_\delta(a,x)$ for $\beta \cross \hat{\nu}$ almost all
$(a,x,f) \in B \cross \hat{X}$.

For $1 \le s \le m$,
let $\sigma_s: B \cross \hat{X} \to \reals$ be defined by
$\sigma_s(a,x,f) = \bar{\sigma}_s(a_1,x,f)$,
where $\bar{\sigma}$ is as in Lemma~\ref{lemma:furstenberg:formula}. 
Then, the left
hand side of (\ref{eq:infinite:limit}) is exactly
\begin{displaymath}
\sum_{j=0}^{n-1} (2\sigma_{r}-\sigma_{r-1}-\sigma_{r+1})(\hat{T}^j(a,x,f)).
\end{displaymath}
Also, we have $n \in \cN_\delta(a,x)$ if and only if $\hat{T}^n(a,x) \in
\cK_{\delta}$. 
Then, by Lemma~\ref{lemma:K}, 
\begin{displaymath}
\int_{B \cross \hat{X}} (2 \sigma_{r} -\sigma_{r-1} -\sigma_{r+1})(q)
\, d(\beta \cross \hat{\nu})(q) > 0. 
\end{displaymath}
By Furstenberg's formula Lemma~\ref{lemma:furstenberg:formula}, the
left hand side of the above equation is $\lambda_r - \lambda_{r+1}$. 
Thus $\lambda_r > \lambda_{r+1}$, contradicting
our assumption that $r \in I$.  
This completes the proof of the claim. 
\qed\medskip



\bold{Proof of Theorem~\ref{theorem:twosided:semisimple:lyapunov}.} 
Pick an orthonormal basis at each point of $X$, and let $C(a \vee b,
x): L \to L$ be a map which makes the subspaces $\cV_i(a \vee b, x)$
orthonormal. Let $\tilde{A}$ denote the cocycle obtained by
\begin{displaymath}
 \tilde{A}(n, a \vee b, x) = C(T^n(a \vee
b,x))^{-1} A(a_n \dots a_1, x) C(a \vee b, x). 
\end{displaymath}
Then $\tilde{A}$ is cohomologous to $A$. Let
\begin{displaymath}
\hat{\eta}(a \vee b,x)= C( a \vee b,x)_* \eta_x, \qquad
\tilde{\eta}_{a \vee b,x} = C(a \vee b,x)_* \eta_{a,x}.
\end{displaymath}
We have, on a set of $\tilde{\beta}^X$ full measure,
\begin{displaymath}
\tilde{\eta}_{a \vee b ,x} = \lim_{n \to \infty} \tilde{A}(n, a
\vee b, x)^{-1}_* \hat{\eta}(T^n(a \vee b,x)).
\end{displaymath}
In view of Lemma~\ref{lemma:non:atomic:lyapunov} 
there exists $\epsilon > 0$ and a
compact  $\cK_{\delta} \subset \tilde{B} \cross X$ with 
$\tilde{\beta}^X(\cK_{\delta}) > 1-\delta$
such that the family of
measures $\{\hat{\eta}(a \vee b,x) \}_{(a \vee b,x) \in \cK_{\delta}}$
is uniformly $(\epsilon,\delta/5)$-regular. 
Write 
\begin{displaymath}
\tilde{A}(n, a \vee b, x)^{-1} = K_n(a \vee b,x) D_n(a \vee b,x)
K_n'(a \vee b,x)
\end{displaymath}
where $K_n$ and $K_n'$ are orthogonal, and $D_n$ is diagonal with
non-increasing entries. 
Let $d_1(n,a \vee b,x)
\ge ... \ge d_m(n,a,x)$ be the entries of $D_n(a \vee b,x)$.

By Claim~\ref{claim:non:singular}, for $r \in I$ and 
almost all $(a \vee b,x)$
$\tilde{\eta}^{m-r}_{a \vee b,x}$ has no atoms. It follows that for every
$\delta >0$ there exists $\cK_1 = \cK_1(\delta) \subset \tilde{B}
\cross X$ and $\epsilon_1 = \epsilon_1(\delta) > 0$, such that for $(a
\vee b,x) \in \cK_1$, $\eta^{m-r}{a \vee b,x}$ gives measure at most
$\delta$ to the $\epsilon_1$-neighborhood of any point. Then, by
Lemma~\ref{lemma:strongly:regular}(a), there exists $C_1 =
C_1(\delta)$ such that if $(a \vee b,x) \in \cK_1(\delta)$ and $T^n(a
\vee b,x) \in \cK_\delta$ then for $r \in I$ 
\begin{equation}
\label{eq:tmp:d:m:minus:r}
d_{m-r}(n,a \vee b,x)/d_{m+1-r}(n, a \vee b, x) \le C_1.
\end{equation}
Note that the matrix of $\tilde{A}(n, a \vee b,x)$ is block diagonal.
We can write each
block as a scaling factor times a determinant one matrix which we
denote by $\tilde{A}_i(n, a \vee b, x)$. (Thus $\tilde{A}_i(n, a \vee
b, x)$ is, up to a scaling factor, a conjugate of the restriction of
$A(n, a \vee b,x)$ to $\cV_i(a \vee b, x)$.) 
Since the subspaces defining the blocks are by construction
orthogonal, the $KAK$ decomposition of $\tilde{A}(n, a \vee b,x)^{-1}$
is compatible with the $KAK$ decompositions of each $\tilde{A}_i(n, a
\vee b, x)^{-1}$. 
Then, (\ref{eq:tmp:d:m:minus:r}) for all $r \in I$ implies that
for all $(a \vee b,x) \in \cK_1(\delta)$ such that $T^n(a
\vee b,x) \in \cK_\delta$ we have
\begin{displaymath}
\|\tilde{A}_i(n,a \vee b, x)\| \le C_1'(\delta) \qquad\text{ for all $i$.}
\end{displaymath}
It follows that for all $n \in \zed$
\begin{displaymath}
\tilde{\beta}^X(\{ (a \vee b,x) \in B \cross X \st
\|\tilde{A}_i(n,a \vee b,x) \| >
C_1'(\delta) \}) \le 2 \delta.
\end{displaymath}
Since $\delta > 0$ is arbitrary, this means (by definition) that the
cocycle $\tilde{A}_i$ is bounded in the sense of Schmidt, see
\cite{Schmidt:Etds1981}. It is proved in \cite{Schmidt:Etds1981}
that any bounded cocycle is conjugate to a cocycle taking values in an
orthogonal group. Therefore the same holds for the determinant one part of
the cocycle $A|_{\cV_i}$.
\qed\medskip

\bold{Proof of Theorem~\ref{theorem:forward:semisimple:lyapunov} and
  Theorem~\ref{theorem:semisimple:lyapunov}.}

To prove Theorem~\ref{theorem:forward:semisimple:lyapunov}, for the
case where the algebraic hull is all of $SL(L)$, it is
enough to show that for almost all $(a,x)$, 
the inner product $\langle \cdot, \rangle_{a \vee
  b,x}$ does not depend on $b$. The proof is similar to the proof of 
(\ref{eq:locally:constant:quadratic:form}). 

For any $\epsilon > 0$  exists a compact set $K \subset \tilde{B}
\cross X$ of measure $1-\epsilon$ such that the map $(a \vee b, x) \to
\langle \cdot, \cdot \rangle_{a \vee b, x}$ is uniformly continuous on
$K$. Then there exists $\Omega \subset \tilde{B} \cross X$ such that 
$\tilde{\beta}^X(\Omega) = 1$ and $T^n(a \vee b, x) \in K$ for set of
$n$ of asymptotic density at least $1/2$.

For $(a \vee b, x) \in \tilde{B} \cross X$ and $v, w \in \cV_{\ge i}(a,x)/\cV_{>i}(a,x)$, let
\begin{displaymath}
[v,w]_{i,(a \vee b,x)} = \frac{\langle v,w\rangle_{i,(a \vee
    b,x)}}{\langle v,v\rangle_{i,(a \vee b,x)}^{1/2} \langle w,w\rangle_{i,(a \vee b,x)}^{1/2}}
\end{displaymath}

Now suppose $(a \vee b, x) \in \Omega$, and $(a \vee b', x) \in \Omega$. 
Consider the points $T^n( a\vee b, x)$ and
$T^n(a \vee b',x)$, as $n \to \infty$. Then $d(T^n(a \vee b, x), T^n(a
\vee b',x) \to 0$. 
Let
\begin{displaymath}
v_n = A(a_n \dots a_1) v, \quad 
w_n = A(a_n \dots a_1) w. 
\end{displaymath}
Then, by Theorem~\ref{theorem:twosided:semisimple:lyapunov},
we have
\begin{equation}
\label{eq:app:rescaled:un}
[v_n, w_n]_{i,T^n(a \vee b, x)} = [v, w ]_{i,x}, \quad 
[v_n, w_n ]_{i,T^n(a \vee b', x)} = [v, w]_{i,(a \vee b',x)}.
\end{equation}

Now take a sequence $n_k \to \infty$ with $T^n(a \vee b, x) \in K$, 
$T^n(a \vee b',x)  \in K$ (such a sequence exists by the definition of
$\Omega$). Then, 
\begin{displaymath}
[v_{n_k}, w_{n_k}]_{i,T^{n_k}(a \vee b,x)} - 
[v_{n_k}, w_{n_k}]_{i,T^n(a \vee b',x)} \to 0.
\end{displaymath}
Now from (\ref{eq:app:rescaled:un}),  we get
\begin{displaymath}
[v, w]_{i,(a \vee b,x)} = [v, w]_{i,(a \vee b',x)}, 
\end{displaymath}
Therefore, for all $v, w \in \cV_{\ge i}(a,x)/\cV_{>i}(a,x)$
\begin{displaymath}
\langle v, w \rangle_{i, (a \vee b,x)} = c(a,b,b',x) \langle v, w
\rangle_{i, (a \vee b',x)}, 
\end{displaymath}
where $c(a,b,b',x) \in \reals^+$. 
We can (measurably) \mccc{(check!)}
choose, for almost all $(a,x)$ some $b_0 \in B$ so
that $(a \vee b_0, x) \in \Omega$, and then replace $\langle \cdot, \cdot
\rangle_{i, (a \vee b, x)}$ by 
\begin{displaymath}
\langle v, w \rangle'_{i,(a, x)} = \langle v, w \rangle_{i, a \vee b_0, x}.
\end{displaymath}
Then $\langle \cdot, \cdot \rangle'_{i,(a,x)}$ satisfies all the
conditions of Theorem~\ref{theorem:forward:semisimple:lyapunov}. 
This concludes the proof of
Theorem~\ref{theorem:forward:semisimple:lyapunov} for the case where
the algebraic hull is all of $SL(L)$.

The proof of
Theorem~\ref{theorem:semisimple:lyapunov} is identical. 
\qed\medskip

\mccc{ issue of covers in 13-end. actually only need to adjust
  notation.}

\section{Dense subgroups of nilpotent groups}
\label{sec:nilpotent}
The aim of this appendix is to prove
Proposition~\ref{prop:measure:inv:closed:subgroup} which is used
in \S\ref{sec:inductive:step}. 

Let $N$ be a nilpotent Lie group. 
For a subgroup $\Gamma \subset N$, let $\bar{\Gamma}$ denote the
topological closure of $\Gamma$, and let $\bar{\Gamma}^0$ denote the
connected component of $\bar{\Gamma}$ containing the identity $e$ of
$N$. Let $B(x,\epsilon)$ denote the ball of radius $\epsilon$ centered
at $x$ in some left-invariant metric on $N$.

\begin{lemma}
\label{lemma:cartan}
Suppose $N$ is a Lie group, and $S \subset N$ is a 
subset. For $\epsilon > 0$, let $\Gamma_\epsilon$ denote the subgroup
generated by $S \cap B(e,\epsilon)$. Then there exists $\epsilon_1 >
0$ and a connected closed Lie subgroup $N_1$ of $N$ such that for $\epsilon <
\epsilon_1$, $\overline{\Gamma}_\epsilon = N_1$. 
\end{lemma}

\bold{Proof.} By Cartan's theorem (see e.g.\ \cite[\S{0.4}]{Knapp:beyond}), 
any closed subgroup of a Lie group is a closed Lie subgroup.  
Let $\epsilon > 0$ be arbitrary. Since we have $\bar{\Gamma}_{\epsilon'}^0
\subset \bar{\Gamma}_{\epsilon}^0$ for $\epsilon' < \epsilon$,
there exists $\epsilon_0 > 0$ such that for $\epsilon \le \epsilon_0$,
the dimension of the Lie algebra of  
$\bar{\Gamma}_\epsilon^0$ (and thus $\bar{\Gamma}_{\epsilon}^0$ itself) 
is independent of $\epsilon$. Thus there exists a connected closed 
subgroup $N_1 \subset N$ such that for
$\epsilon \le \epsilon_0$, $\bar{\Gamma}_\epsilon^0 = N_1$. In
particular, 
\begin{equation}
\label{eq:tmp:bar:Gamma:epsilon}
\bar{\Gamma}_\epsilon \supset N_1. 
\end{equation}

From the definition it is immediate that
$\bar{\Gamma}_{\epsilon_0}$ is a closed subgroup of $N$. 
Thus, by Cartan's theorem, 
$\bar{\Gamma}_{\epsilon_0}$ and $N_1 = \bar{\Gamma}_{\epsilon_0}^0$ are
closed submanifolds of
$N$. Therefore, there exists $\epsilon_1 < \epsilon_0$ such that
\begin{displaymath}
B(e,\epsilon_1) \cap \bar{\Gamma}_{\epsilon_0} = B(e,\epsilon_1) \cap
\bar{\Gamma}_{\epsilon_0}^0 = B(e,\epsilon_1) \cap N_1. 
\end{displaymath}
Then, for $\epsilon < \epsilon_1 < \epsilon_0$, 
\begin{displaymath}
\Gamma_\epsilon \cap B(e,\epsilon_1) \subset \bar{\Gamma}_{\epsilon_0}
\cap B(e,\epsilon_1) \subset N_1.
\end{displaymath}
Therefore, $\Gamma_\epsilon \subset N_1$, and hence
$\bar{\Gamma}_\epsilon \subset N_1$. In view of
(\ref{eq:tmp:bar:Gamma:epsilon}), the lemma follows.  
\qed\medskip

\begin{lemma}
\label{lemma:replace:by:goodword}
Suppose $N$ is a simply connected 
nilpotent Lie group, and let $S \subset N$ be an
(infinite) subset. For each $\epsilon > 0$ let $\Gamma_\epsilon
\subset N$ denote the subgroup of $N$ generated by the elements
$\gamma \in S \cap B(e,\epsilon)$. Suppose that for all $\epsilon >
0$, $\Gamma_\epsilon$ is dense in $N$. 

Then, for every
$\epsilon > 0$ there exist $0 < \theta < \epsilon$ 
(depending on $\epsilon$ and $S$) such that for
every $\gamma \in \Gamma_\epsilon$ with $d(\gamma,e) < \theta$
there exists $n \in \natls$ and for $1 \le i \le n$ 
elements $\gamma_i \in S$ with 
\begin{equation}
\label{eq:replace:by:goodword:one}
\gamma = \gamma_n \dots \gamma_1
\end{equation}
and for each $1 \le j \le n$, 
\begin{equation}
\label{eq:replace:by:goodword:two}
d(\gamma_j \dots \gamma_1, e) < \epsilon.
\end{equation}
\end{lemma}

\bold{Proof.} We will proceed by induction on $\dim N$. Let $N' =
[N,N]$ 

For $k \in \natls$, let $S^k_\epsilon$ 
be the product of at most $k$ elements in $(S \cup S^{-1}) \cap
B(e, \epsilon)$. Let $T^k_\epsilon =  [S^k_\epsilon, S^k_\epsilon]$. 
This decreases with $\epsilon$, so a variant of Lemma~\ref{lemma:cartan}
shows that, for small enough $\epsilon$, the closure of the group
generated by $T^k_\epsilon$ is 
a closed connected
group $N_k$ (and $N_k$ is independent of $\epsilon$ for $\epsilon$
small enough). 
 Since $N_k$ increases with $k$, 
it is constant for large $k$. Fix $k$ so
that $N_k = N_{k+2}$. 
We will show that $N_k = N'$.  

First, we show that $N_k$ is normal. For $a, b \in  S^k_\epsilon$ 
and $s  \in S_\epsilon$, we have
$s[a, b]s^{-1} = [sas^{-1}, sbs^{-1}] \in T^{k+2}_\epsilon$. 
So, $s T^k_\epsilon s^{-1} \subset T^{k+2}_\epsilon$. Taking the
closure of the generated groups, we get $sN_ks^{-1} \subset  N_{k+2} = N_k$. Hence, $N_k$ is
normalized by $S_\epsilon$. 
Since $S_\epsilon$ generates a dense subset of $N$, $N_k$ is normal.

We have $[ab, c] = a[b, c]a^{-1} [a, c]$. 
This shows that, if $[a, c]$ and $[b, c]$
both belong to $N_k$, then $[ab, c]$ also belongs to $N_k$, by normality. For
$x, y \in  S^k_\epsilon$, we have $[x, y] \in N_k$. 
Taking products, and since $S^k_\epsilon$ generates a dense subgroup of
$N$, we get $[z, y] \in  N_k$ for all $z \in  N$. 
Doing the same argument with the
other variable, we finally have $[z, z'] \in N_k$ for all $z,z' \in
N$, and therefore $N_k = N'$ as desired.

Let $S' = T^k_{\epsilon/4k} \subset N'$. 
For $\delta > 0$ let $\Gamma'_{\delta}$ denote the subgroup of
$N'$ generated by $S' \cap B(e,\delta)$. Since (for sufficiently
small $\delta$) $[B(e,\delta), B(e,\delta)] \subset B(e,\delta)$, we
have, for $\delta < \epsilon/{4k}$, 
\begin{displaymath}
\overline{\Gamma_{\delta}'} \supset
\overline{\{ \text{the subgroup generated by }T^k_{\delta/{4k}}\}} = N'.
\end{displaymath}
\mcc{justify? probably have to say something}
Therefore, $S' \subset N'$ satisfies the conditions of the Lemma.  
Let $\epsilon' > 0$ be such that \mcc{check}
\begin{equation}
\label{eq:tmp:choice:epsilonprime}
B(e,\epsilon')^k \subset B(e,\epsilon/100). 
\end{equation}
Since $\dim N' < \dim N$, by the inductive assumption 
there exist $0 < \theta' < \epsilon'$ 
such that for any $\gamma' \in \Gamma'_{\theta'}$ with
$d(\gamma',e) < \theta'$, there exist $\gamma_i' \in
S'$ such that 
(\ref{eq:replace:by:goodword:one}) holds, and 
(\ref{eq:replace:by:goodword:two}) holds with $\epsilon'$ in place
of $\epsilon$.



Suppose $\epsilon > \eta > 0$. 
By construction, $N/N'$ is abelian. Note that $N$ is connected and
simply connected. Then, since
$\bar{\Gamma}_{\eta} = N$, 
there exists a finite set 
\begin{displaymath}
S_0 \equiv \{ \lambda_1, \dots, \lambda_k \} \subset
\Gamma_{\eta} \cap S 
\end{displaymath}
with $d(\lambda_i,e) < \eta$ for $1 \le i \le k$ so that 
${\lambda_1 N', \dots, \lambda_k N'}$ form a basis over $\reals$ for
the vector space $N/N'$.  Let $\Lambda$ denote the subgroup 
generated by the $\lambda_i$, and let $F' \subset N/N'$ denote the
parallelogram 
centered at the origin whose sides are parallel to the vectors
$\lambda_i N'$. Then $F'$ is a  fundamental domain for the action of
$\Lambda$ on $N/N'$, and  
\begin{displaymath}
\diam F' = O(\eta). 
\end{displaymath}
Let $N_0$ be a local complement to $N'$ in $N$ near the identity $e$.  
We can choose $N_0$ to be a
smooth manifold transversal to $N'$ ($N_0$ need not be a
subgroup). Let $\pi: N \to N/N'$ be the natural map, and let $\pi^{-1}: N/N'
\to N_0$ be the inverse. Let $F = \pi^{-1}(F')$. 
We can now choose $\eta$ sufficiently small so 
that $F \subset B(e,\rho)$, where $\theta' > \rho > \eta > 0$ is such that
\begin{displaymath}
B(e,\rho)^5 \cap N' = [B(e,\rho) B(e,\rho) B(e,\rho) B(e,\rho)
B(e,\rho)] \cap N' \subset B(e, \theta') \cap N'. 
\end{displaymath}
 
We now choose $\theta > 0$ so that $B(e,\theta)  \subset F \cO$
where $\cO \subset N' \cap B(e,\rho)$ 
is some neighborhood of the origin. We now claim that for any $x \in F
\cO$ and any $s \in B(e,\theta)$, there exist $\lambda' \in S_0
\cup S_0^{-1}$ and $\gamma' \in \Gamma'_{\theta'}$ 
such that $\gamma' \lambda' s x
\in F \cO$. Indeed, since $B(e,\theta)N' \subset F N'$, for any $x \in F
N'$, 
\begin{displaymath}
B(x,\theta)N' \subset \bigcup_{\lambda \in S_0 \cup S_0^{-1}}
\lambda B(x,\theta) N'.
\end{displaymath}
Thus, we can find $\lambda' \in S_0 \cup S_0^{-1}$ such that $\lambda'
s x \in F N'$. Since $\Gamma'_{\theta'}$ is dense in $N'$, there
exists $\gamma' 
\in \Gamma'_{\theta'}$ such that $\gamma' \lambda' s x \in F \cO$, completing the
proof of the claim. 

Now suppose $\gamma \in \Gamma_{\theta}$ and $\gamma \in
B(e,\theta) \subset F \cO$. 
Then, we have
\begin{displaymath}
\gamma = s_n \dots s_1, \text{ where $s_i \in S \cap
  B(e,\theta)$. }
\end{displaymath}
Note that $s_1 \in F \cO$. 
We now define elements $\lambda_j' \in S_0 \cup S_0^{-1}$ and $\gamma_j' \in
\Gamma'_{\theta'}$ inductively as follows. At every stage of the
induction, we will have $x_j \equiv \gamma_j' \lambda_j' 
s_j \dots \gamma_1' \lambda_1' 
s_1 \in F \cO$. Suppose $\gamma_1',\dots, \gamma_{j-1}'$ and
$\lambda_1', \dots \lambda_{j-1}'$ have already been chosen. Now
choose $\lambda_j' \in S_0 \cup S_0^{-1}$ and $\gamma_j' \in
\Gamma'_{\theta'}$ so that $x_j = \gamma_j'
\lambda_j' s_j x_{j-1} \in F \cO$. Such $\lambda_j'$ and $\gamma_j'$ exist
by the claim.

Note that 
\begin{displaymath}
\gamma_j' = x_j x_{j-1}^{-1} s_j^{-1} (\lambda_j')^{-1} \in (F \cO ) (F
\cO)^{-1} B(e,\theta)^{-1} (S_0 \cup S_0^{-1}) \subset B(e,\rho)^5 \subset
B(e,\theta').  
\end{displaymath}
Since $x_n = \lambda_n' \gamma_n' s_n \dots \lambda_1' \gamma_1' s_1 \in F
N'$, we have $\lambda_n' s_n \dots \lambda_1' s_1 \in F N'$. Also 
$\gamma = s_n \dots s_1 \in B(x,\theta) \subset F N'$. Since
$F N'$ is a fundamental
domain for the action of $\Lambda$ on $N/N'$, 
$\lambda_n' \dots \lambda_1'  \in N'$. Thus,  
\begin{equation}
\label{eq:preword:gamma}
\gamma = \gamma'  \gamma_n' \lambda_n' s_n \dots  \gamma_1' \lambda_1' s_1,
\end{equation}
where $\gamma' \in N'$. We have
\begin{displaymath}
\gamma' = \gamma x_n^{-1} \in B(e,\theta) (F \cO)^{-1}
\subset B(e,\theta'). 
\end{displaymath}
For notational convenience, denote $\gamma'$ by $\gamma'_{n+1}$. 
By the inductive assumption, for $1 \le i \le n+1$, 
we can express $\gamma_i' = s'_{i1}
\dots s'_{in_i}$ such that $s'_{ij} \in S' \cap B(e,\theta')$ and
so that for all $i$, $j$, 
\begin{displaymath}
d(s_{ij}' \dots s_{i1}', e) \le \epsilon'.
\end{displaymath}
We now substitute this into (\ref{eq:preword:gamma}). Finally, we
express each $s_{ij}'$ as a commutator of a product of at most $k$
elements of $S \cap B(e,\epsilon/4k)$. Then, in view of (\ref{eq:tmp:choice:epsilonprime}),
the resulting word satisfies (\ref{eq:replace:by:goodword:two}). 
\qed\medskip

\begin{proposition}
\label{prop:measure:inv:closed:subgroup}
Suppose $N$ is a simply connected 
nilpotent Lie group, $\cO$ a neighborhood of the
identity in $N$, and $\mu$ a measure on $N$ supported on
$\cO$. Suppose $S \subset N$ is a subset containing elements
arbitrarily close to (and distinct from) $e$, and suppose for each
$\gamma \in S$, 
\begin{equation}
\label{eq:gamma:star:mu:is:mu}
\gamma_* \mu \; \propto \; \mu
\end{equation}
on $\cO \cap \gamma^{-1} \cO$ where both sides make sense. 
Then, there exists a nontrivial connected subgroup $H$ of $N$ and a
neighborhood $\cO'$ of the identity in $H$ 
  such that for all $h \in \cO'$, $h_* \mu \propto \mu$ on $\cO
  \cap h^{-1} \cO$. Furthermore, if
$U$ is a connected subgroup of $N$ and $S$ contains arbitrarily small
elements not contained in $U$, then $H$ is not contained in $U$. 
\end{proposition}

\bold{Proof.} Let $N_1$ and $\epsilon_1$ be as in
Lemma~\ref{lemma:cartan}. By our assumptions on $S$, $N_1$ is
non-trivial (and also $N_1$ is not contained in $U$). Now suppose
$\epsilon > 0$ is such that $B(e,\epsilon) \subset \cO$, and let
$\theta > 0$ be as in Lemma~\ref{lemma:replace:by:goodword}, with $N$
replaced by $N_1$. Without loss of generality, we may assume that
$\theta < \epsilon_1$. 
Let $\Gamma_\theta$ be the subgroup of $N_1$
generated by $S \cap B(e,\theta)$. Since $\theta < \epsilon_1$,
$\Gamma_\theta$ is dense in $N_1$. Now suppose $\bar{\gamma} \in N_1$, and
$d(\bar{\gamma},e) < \theta$. Then, there exists $\gamma_k \in
\Gamma_\theta$ such that $\gamma_k \to \gamma$,  and $d(\gamma_k, e)
< \theta$. We can write each $\gamma_k = \gamma_{k,n} \dots
\gamma_{k,1}$ as in Lemma~\ref{lemma:replace:by:goodword}. Then, by
applying (\ref{eq:gamma:star:mu:is:mu}) repeatedly, we get that
$(\gamma_k)_* \mu \,\propto\, \mu$. Then, taking the
limit as $k \to \infty$ we see 
that $(\bar{\gamma})_* \mu \, \propto \,\mu$. Thus, $\mu$
is invariant (up to normalization) under a 
neighborhood of the origin in $N_1$. 
\qed\medskip


\printindex


\begin{thebibliography}{SCh95b}


\bibitem[ABEM]{ABEM}
J.~Athreya, A.~Bufetov, A.~Eskin and M.~Mirzakhani.
Lattice Point Asymptotics and Volume Growth on Teichm\"uller
  space. 	
\emph{Duke Math. J.}  \textbf{161}  (2012),  no. 6, 1055--1111.


\bibitem[At]{Atkinson}
G.~Atkinson. Recurrence of co-cycles and random walks.
\emph{ J. London Math. Soc.} (2)  \textbf{13}  (1976),  no. 3, 486--488.



\bibitem[Ath]{Jayadev:thesis}
J.~Athreya. Quantitative recurrence and large deviations for  
Teichm\"uller geodesic flow.  \emph{Geom. Dedicata}  \textbf{119}
(2006), 121--140.

\bibitem[AthF]{Athreya:Forni}
J.~Athreya, G.~Forni. 
Deviation of ergodic averages for rational polygonal billiards. \emph{Duke
Math. J.} \textbf{144} (2008), no. 2, 285–-319.  


\bibitem[ACO]{Arnold} 
L.~Arnold, N.~Cong, V.~Osceledets.
Jordan Normal Form for Linear Cocycles.
\emph{Random Operators and Stochastic Equations} \textbf{7} (1999), 301-356.


\bibitem[AEZ]{Athreya:Eskin:Zorich} 
J.~Athreya, A.~Eskin, A.~Zorich.
   Rectangular  billiards and volumes of spaces of quadratic
    differentials on $\mathbb{C}P^1$ (with an appendix by Jon Chaika).
\emph{Ann. Sci. \'Ec. Norm. Sup\'er.} (4)  \textbf{49}  (2016),  no. 6,
  1311--1386. 


\bibitem[ANW]{Aulicino:Nguyen:Wright}
D.~Aulicino, D.~Nguyen, A.~Wright.
Classification of higher rank orbit closures in H\^{}\{odd\}(4).
\emph{J. Eur. Math. Soc. (JEMS)}  \textbf{18}  (2016),  no. 8, 1855--1872.

\bibitem[AEM]{Avila:Eskin:Moeller:yeti}
A.~Avila, A.~Eskin, M.~Moeller.
Symplectic and Isometric SL(2,R) invariant subbundles of the Hodge
bundle.
\emph{J. Reine Angew. Math.}  \textbf{732}  (2017), 1--20.



\bibitem[AF]{Avila:Forni}
A.~Avila, G.~Forni.
Weak mixing for interval exchange transformations and translation
flows. 
\emph{Ann. of Math.} (2) \textbf{165} (2007), no. 2, 637--664.

\bibitem[AG]{Avila:Gouezel}
A.~Avila, S.~Gou\"ezel. 
Small eigenvalues of the Laplacian for algebraic measures in moduli
space, and mixing properties of the Teichmüller flow.  
\emph{Ann. of Math.} (2)  \textbf{178}  (2013),  no. 2, 385--442.

\bibitem[AGY]{Avila:Gouezel:Yoccoz}
A.~Avila, S.~Gou\"ezel, J-C.~Yoccoz. 
Exponential mixing for the Teichm\"uller flow.
\emph{Publ. Math. Inst. Hautes \'Etudes Sci.}  \textbf{104}  (2006), 143--211.


\bibitem[AV2]{Avila:Viana:extremal} 
A.~Avila, M.~Viana. Extremal Lyapunov exponents: an invariance
principle and applications. 
\emph{Invent. Math.}  \textbf{181}  (2010),  no. 1, 115--189.

\bibitem[ASV]{Avila:Santamaria:Viana} 
A.~Avila, J.~Santamaria, M.~Viana. 
Cocycles over partially hyperbolic maps.
\emph{Ast\'erisque}  No. \textbf{358}  (2013), 1--12.

\bibitem[AV1]{Avila:Viana} 
A.~Avila, M.~Viana. 
Simplicity of Lyapunov spectra: proof of the Zorich-Kontsevich conjecture,
    \emph{Acta Math.} \textbf{198} , No. 1, (2007), 1--56.


\bibitem[Ba]{Bainbridge:L:shaped} 
M.~Bainbridge.
Billiards in L-shaped tables with barriers.
\emph{Geom. Funct. Anal.}  {\bf 20 } (2010),  no. 2, 299--356.

\bibitem[BaM]{Bainbridge:Moeller}
M. Bainbridge and M. M\"oller.
Deligne-Mumford compactification of the real multiplication locus and
Teichm\"uller curves in genus 3.
\emph{Acta Math.}  \textbf{208}  (2012),  no. 1, 1--92.
		
\bibitem[BoM]{Bouw:Moeller}
I. Bouw and M. M\"oller. 
Teichm\"uller curves, triangle groups, and Lyapunov exponents.
 \emph{Ann. of Math.} (2)  \textbf{172}  (2010),  no. 1, 139--185.



\bibitem[BQ]{Benoist:Quint}
Y.~Benoist and J-F~Quint.
Mesures Stationnaires et Ferm\'es Invariants des espaces
  homog\`enes. 
(French)  [Stationary measures and invariant subsets of homogeneous spaces]  \emph{Ann. of Math.} (2)  \textbf{174}  (2011),  no. 2, 1111--1162.


\bibitem[BG]{Bufetov:Gurevich}
A.~Bufetov, B.~Gurevich.
Existence and Uniqueness of the Measure of Maximal Entropy for the
Teichm\"uller Flow on the Moduli Space of Abelian Differentials. 
(Russian)  \emph{Mat. Sb.}  \textbf{202}  (2011),  no. 7, 3--42; translation in  \emph{Sb. Math.}  \textbf{202}  (2011),  no. 7-8, 935--970




\bibitem[Ca]{Calta:thesis}
K.~Calta.
Veech surfaces and complete periodicity in genus two.  
\emph{J. Amer. Math. Soc.}  \textbf{17}  (2004),  no. 4, 871--908 

\bibitem[CK]{Climenhaga:Katok}
V.~Climenhaga, A.~Katok.
Measure theory through dynamical eyes.
\texttt{arXiv:1208.4550 [math.DS]}.




\bibitem[CW]{Calta:Wortman}
K.~Calta, K.~Wortman.
On unipotent flows in H(1,1).
\emph{Ergodic Theory Dynam. Systems}  \textbf{30}  (2010),  no. 2, 379--398.



\bibitem[Dan1]{Dani:invariant:1979}
S.G.\ Dani,
\newblock {On invariant measures, minimal sets and a lemma of {M}argulis},
\newblock {\em Invent. Math.} {\bf 51} (1979), 239--260.

\bibitem[Dan2]{Dani:horospherical:1981}
S.G.\ Dani,
\newblock {Invariant measures and minimal sets of horoshperical flows},
\newblock {\em Invent. Math.} {\bf 64} (1981), 357--385.

\bibitem[Dan3]{Dani:orbits:I}
S.G. Dani,
\newblock {On orbits of unipotent flows on homogeneous spaces},
\newblock {\em Ergod.\ Theor.\ Dynam.\ Syst.} {\bf 4} (1984), 25--34.

\bibitem[Dan4]{Dani:orbits:II}
S.G.\ Dani,
\newblock {On orbits of unipotent flows on homogenous spaces {II}},
\newblock {\em Ergod.\ Theor.\ Dynam.\ Syst.} {\bf 6} (1986), 167--182.


\bibitem[De]{Encyclopaedia:distances}
M.~Deza and E.~Deza,
\newblock{Encyclopaedia of Distances, 3rd edition}, Springer-Verlag, 2014. 


\bibitem[DM1]{Dani:Margulis:values}
S.G.\ Dani and G.A.\ Margulis,
\newblock {Values of quadratic forms at primitive integral points},
\newblock {\em Invent.\ Math.} {\bf 98} (1989), 405--424.

\bibitem[DM2]{Dani:Margulis:unipotent}
S.G.\ Dani and G.A.\ Margulis,
\newblock {Orbit closures of generic unipotent flows on homogeneous spaces of
  $SL(3,\reals)$},
\newblock {\em Math.\ Ann.} {\bf 286} (1990), 101--128.


\bibitem[DM3]{Dani:Margulis:asymptotic}
S.G.\ Dani and G.A.\ Margulis.
\newblock {Asymptotic behaviour of trajectories of unipotent flows on
  homogeneous spaces},
\newblock {\em Indian.\ Acad.\ Sci.\ J.} {\bf 101} (1991), 1--17.


\bibitem[DM4]{Dani:Margulis:distribution}S.G.\ Dani and G.A.~Margulis,
  \newblock {Limit distributions of orbits of unipotent
    flows and values of quadratic forms}, \newblock in: \emph{ I.\ M.\
    Gelfand Seminar}, Amer.\ Math.\ Soc., Providence, RI, 1993, pp.\
  91--137.


\bibitem[Ef]{Effros} Effros, Edward G.  
Transformation groups and $C^*$-algebras. 
\emph{Ann. of Math.} (2) \textbf{81} (1965) 38–55. 

\bibitem[EKL]{EKL} M. Einsiedler, A. Katok, E. Lindenstrauss. 
Invariant measures and the set of exceptions to Littlewood's
conjecture.  
\emph{Ann. of Math.} (2)  \textbf{164}  (2006),  no. 2, 513--560.




\bibitem[EL]{Einsiedler:Lindenstrauss:Pisa}
M.~Einsiedler, E.~Lindenstrauss. 
``Diagonal actions on locally homogeneous spaces.'' 
In Homogeneous flows, moduli spaces and arithmetic, 155–241, Clay
Math. Proc., 10, Amer. Math. Soc., Providence, RI, 2010. 



\bibitem[EMa]{Eskin:Masur} A.~Eskin, H.~Masur. 
Asymptotic formulas on flat surfaces.
\emph{Ergodic Theory  and Dynamical
    Systems,} {\bf 21} (2) (2001), 443--478.



\bibitem[EMM]{Eskin:Marklof:Morris}
A.~Eskin, J.~Marklof, D.~Morris.
Unipotent flows on the space of branched covers of Veech
  surfaces.
\emph{Ergodic Theory Dynam. Systems}  \textbf{26}  (2006),  no. 1, 129--162.


\bibitem[EMM1]{Eskin:Margulis:Mozes:31}
A.~Eskin, G.~Margulis, and S.~Mozes.
\newblock {Upper bounds and asymptotics in a quantitative version of
the   {O}ppenheim conjecture.}
\newblock {\em  Ann. of Math. (2)} {\bf 147} (1998), no. 1, 93--141.


\bibitem[EMM2]{Eskin:Margulis:Mozes:22}
A.~Eskin, G.~Margulis, and S.~Mozes.
\newblock {Quadratic forms of signature $(2,2)$ and eigenvalue spacings
on flat 2-tori.} 
\newblock{{\em Ann. of Math. (2)}  {\bf 161} (2005),  no. 2, 679--725.}


\bibitem[EMiMo]{Eskin:Mirzakhani:Mohammadi}
A. Eskin, M.~Mirzakhani and A. Mohammadi.
Isolation, equidistribution, and orbit closures for the $SL(2,\reals)$
action on Moduli space.
\emph{Ann. of Math.} (2)  \textbf{182}  (2015),  no. 2, 673--721.


\bibitem[EMR]{Eskin:Mirzakhani:Rafi}
A.~Eskin, M.~Mirzakhani and K.~Rafi.
Counting closed geodesics in strata. 
\texttt{arXiv:1206.5574 [math.GT]} (2012).




\bibitem[EMS]{Eskin:Masur:Schmoll} A.~Eskin, H.~Masur, M.~Schmoll, 
Billiards in rectangles with barriers.
\emph{Duke Math. J.} \textbf{118} No. 3 (2003), 427--463.

\bibitem[EMZ]{Eskin:Masur:Zorich} A.~Eskin, H.~Masur,  A.~Zorich,
Moduli  spaces  of  Abelian  differentials:  the  principal
boundary, counting problems and  the  Siegel--Veech constants.
\emph{Publ. Math. Inst. Hautes \'Etudes Sci.}, {\bf 97} (1) (2003),
61--179.

\bibitem[EMat]{Eskin:Matheus:semisimple}
A.~Eskin, C.~Matheus.
Semisimplicity of the Lyapunov spectrum for irreducible cocycles.
\emph{Preprint.}






\bibitem[Fo]{Forni:Deviation} G.~Forni. 
Deviation  of  ergodic averages  for area-preserving  flows  on
surfaces of higher     genus.  
\emph{Ann. of Math.}, \textbf{155}, no. 1, 1--103  (2002)

\bibitem[Fo2]{Forni:handbook}
G.~Forni. {On the {L}yapunov exponents of the {K}ontsevich-{Z}orich
              cocycle}. 
\emph{Handbook of dynamical systems. {V}ol. 1{B}}, pp 549--580. 
Elsevier B. V., Amsterdam, 2006.




\bibitem[FoM]{Forni:Matheus}  G.~Forni, C.~Matheus,
An example of a Teichm\"uller disk in genus 4 with degenerate
Kontsevich-Zorich spectrum, \texttt{arXiv:0810.0023} (2008), 1--8.


\bibitem[FoMZ]{Forni:Matheus:Zorich}
G.~Forni, C.~Matheus, A.~Zorich.
Lyapunov Spectrum of Invariant Subbundles of the Hodge Bundle.
\emph{Ergodic Theory Dynam. Systems}  \textbf{34}  (2014),  no. 2, 353--408.






\bibitem[Fu]{Furman:Survey}
A.~Furman.
Random walks on groups and random transformations.
\emph{Handbook of dynamical systems}, Vol. 1A, 931 - 1014, 
North-Holland, Amsterdam, 2002. 



\bibitem[F1]{F1}
H.~Furstenberg.
A Poisson formula for semi-simple Lie groups. 
\emph{Ann. of Math.} \textbf{77}(2),
335~-386 (1963).

\bibitem[F2]{F2}
H.~Furstenberg. 
Non commuting random products. 
\emph{Trans. Amer. Math. Soc.} \textbf{108},
377~-428 (1963). 

\bibitem[Fi1]{Filip:semisimplicity}
S.~Filip.
``Semisimplicity and rigidity of the Kontsevich-Zorich cocycle.''
 \textit{Invent. Math.}  \textbf{205}  (2016),  no. 3, 617--670.

\bibitem[Fi2]{Filip:splitting}
S.~Filip.
``Splitting mixed Hodge structures over affine invariant manifolds.''
 \textit{Ann. of Math.} (2)  \textbf{183}  (2016),  no. 2, 681--713.

\bibitem[GM]{Goldsheid:Margulis}
I.Ya. Gol'dsheid and G.A. Margulis. 
\newblock{Lyapunov indices of a product of random matrices.}
\emph{Russian Math. Surveys} \textbf{44:5} (1989), 11-71.

\bibitem[GR1]{Guivarch:Raugi:Frontiere}
Y. Guivarc'h and A. Raugi, 
\newblock{Frontiere de Furstenberg, propriotes de contraction et
     theoremes de convergence,}
\emph{ Z. Wahrsch. Verw. Gebiete} \textbf{69} (1985), 187-242.


\bibitem[GR2]{Guivarch:Raugi:Contraction}
Y.~Guivarc'h and A.~Raugi. 
\newblock{Propri\'et\'es de contraction d'un semi-groupe de matrices
  inversibles. Coefficients de Liapunoff d'un produit de matrices
  al\'eatoires ind\'ependantes. (French) [Contraction properties of an
  invertible matrix semigroup. Lyapunov coefficients of a product of
  independent random matrices]} 
\emph{Israel J. Math.} \textbf{65} (1989), no. 2, 165–196.


\bibitem[HLM]{Hubert:Lanneau:Moeller:Ratner}
P.~Hubert, E.~Lanneau, M. M\"oller. 
${\rm GL}^+_2(\mathbb{R})$-orbit closures via topological
splittings. 
Surveys in differential geometry. Vol. XIV. Geometry of Riemann surfaces and their moduli spaces, 145~169,
Surv. Differ. Geom., 14, Int. Press, Somerville, MA, 2009. 




\bibitem[HST]{Hubert:Schmoll:Trubetzkoy}
P. Hubert, M. Schmoll and S. Troubetzkoy,
Modular fibers and illumination problems.  
\emph{Int. Math. Res. Not. IMRN}  \textbf{2008},  no. 8,

\bibitem[Ka]{Kac:formula}
M.~Kac.  On the notion of recurrence in discrete stochastic processes.
 \emph{Bull. Amer. Math. Soc.}  \textbf{53},  (1947). 1002--1010.

\bibitem[Ke]{Kesten:sums}
H.~Kesten. Sums of stationary sequences cannot grow slower than linearly.
 \emph{Proc. Amer. Math. Soc.}  \textbf{49},  (1975), 205--211.

\bibitem[Kn]{Knapp:beyond}
A.~Knapp. 
Lie Groups, Beyond an Introduction. Second Edition. Progress in
Mathematics Vol. \textbf{140}, Birkh\"auser, Boston, 2002. 


\bibitem[KS]{Kalinin:Sadovskaya}
B.~Kalinin, V.~Sadovskaya.
Cocycles with one exponent over partially hyperbolic systems.
\emph{Geom. Dedicata}  \textbf{167}  (2013), 167--188.

\bibitem[KH]{Katok:Hasselblatt}
A.~Katok, B.~Hasselblat.
Introduction to the Modern Theory of Dynamical Systems. 
Cambridge University Press, 1995.



\bibitem[KuSp]{Lojasiewicz}
K.~Kurdyka, S.~Spodzieja. 
Separation of real algebraic sets and the \L{}ojasiewicz exponent.
 \emph{Proc. Amer. Math. Soc.}  \textbf{142}  (2014),  no. 9, 3089--3102.


\bibitem[LN1]{Lanneau:Nguyen:Teichmuller:curves}
E.~Lanneau, D.~Nguyen.
Teichmueller curves generated by Weierstrass Prym eigenforms in genus
three and genus four.
\emph{J. Topol.}  \textbf{7}  (2014),  no. 2, 475--522.

\bibitem[LN2]{Lanneau:Nguyen:periodicity}
E.~Lanneau, D.~Nguyen.
Complete periodicity of Prym eigenforms.
Ann. Sci. \'Ec. Norm. Sup\'er. (4)  \textbf{49}  (2016),  no. 1, 87--130.

\bibitem[LN3]{Lanneau:Nguyen:GL2R}
E.~Lanneau, D.~Nguyen.
$GL^+(2,R)$-orbits in Prym eigenform loci.
\emph{Geom. Topol.}  \textbf{20}  (2016),  no. 3, 1359--1426.


\bibitem[L]{Ledrappier:Positivity}
F.~Ledrappier. Positivity of the exponent for stationary sequences of matrices.
Lyapunov exponents (Bremen, 1984), 
56--73, Lecture Notes in Math., 1186, Springer, Berlin,  1986. 

\bibitem[LS]{Ledrappier:Strelcyn}
F.~Ledrappier and J.~M.~Strelcyn. 
A proof of the estimation from below in Pesin's entropy formula. 
\emph{Ergodic Theory Dyn. Syst.} \textbf{2} (1982), 203--219.

\bibitem[LY]{Ledrappier:Young}
F.~ Ledrappier, and L.~S.~Young: 
The metric entropy of diffeomorphisms. I. 
\emph{Ann. Math.} \textbf{122} (1985), 503--539.



\bibitem[M1]{Mane:article}
R. Ma\~n\'e. A proof of Pesin's formula. \emph{Ergodic Theory
  Dyn. Syst.} \textbf{1} (1981), 95--102.

\bibitem[M2]{Mane:book}
R. Ma\~n\'e. Ergodic Theory and Differentiable Dynamics. Springer-Verlag
1987. 




\bibitem[Mar1]{Margulis:return}
G.A.\ Margulis,
\newblock {\em On the action of unipotent groups in the space of lattices},
\newblock In \emph{Lie Groups and their representations}, Proc.\ of Summer School
  in Group Representations, Bolyai Janos Math.\ Soc., Akademai Kiado, Budapest,
  1971, p.\ 365--370, Halsted, New York, 1975.

\bibitem[Mar2]{Margulis:formes}
G.A.\ Margulis,
\newblock {Formes quadratiques ind\`efinies et flots unipotents sur les spaces
  homog\`enes},
\newblock {\em C.\ R.\ Acad.\ Sci.\ Paris Ser.\ I} {\bf 304} (1987), 247--253.

\bibitem[Mar3]{Margulis:number}
G.A.\ Margulis,
\newblock {em Discrete Subgroups and Ergodic Theory},
\newblock in: \emph{Number theory, trace formulas and discrete subgroups}, a
  symposium in honor of A {S}elberg, pp.\ 377--398. Academic Press, Boston,
  MA, 1989.

\bibitem[Mar4]{Margulis:indefinite}
G.A.\. Margulis,
\newblock {\em Indefinite quadratic forms and unipotent flows on homogeneous
  spaces},
\newblock In \emph{Dynamical systems and ergodic theory}, Vol.\ 23, pp.\
  399--409, Banach Center Publ., PWN -- Polish Scientific Publ., Warsaw, 1989.


\bibitem[MaT]{Margulis:Tomanov:Ratner}
G.~A.~Margulis and G.~M.~Tomanov.
Invariant measures for actions of unipotent groups over local fields
on homogeneous spaces.
\emph{Invent. Math.} \textbf{116} (1994), 347--392.




\bibitem[Mas1]{masur:interval} H.~Masur. 
Interval exchange transformations and measured foliations. 
\emph{Ann. of Math.} (2)  \textbf{115}
  (1982), no. 1, 169--200.

\bibitem[Mas2]{Masur:upper}
  H.~Masur.
  The growth rate of trajectories of a quadratic differential.
  \emph{Ergodic Theory Dynam. Systems} \textbf{10} (1990), 151--176.

\bibitem[Mas3]{Masur:lower}
  H.~Masur.
  Lower bounds for the number of saddle connections and closed
trajectories of a quadratic differential.
  In \emph{Holomorphic Functions and Moduli}, Vol.~1, D.~Drasin, ed.,
  Springer-Verlag: New York, 1988, pp.~215--228.

\bibitem[MW]{Matheus:Wright}
C.~Matheus, A.~Wright. 
Hodge-Teichmueller planes and finiteness results for Teichmueller curves.
\emph{Duke Math. J.}  \textbf{164}  (2015),  no. 6, 1041--1077.

\bibitem[Mc1]{McMullen:Billiards}
C. McMullen. 
Billiards and Teichmüller curves on Hilbert modular surfaces. 
\emph{J. Am. Math. Soc.} \textbf{16}, 857--885 (2003) 

\bibitem[Mc2]{McMullen:geodesics}
C. McMullen. 
Teichmüller geodesics of infinite complexity. 
\emph{Acta Math.} \textbf{191}, 191--223 (2003)

\bibitem[Mc3]{McMullen:discriminant}
 C. McMullen. 
Teichmüller curves in genus two: Discriminant and spin. 
\emph{Math. Ann.} \textbf{333}, 87--130 (2005)

\bibitem[Mc4]{McMullen:decagon} 
C. McMullen. 
Teichmüller curves in genus two: The decagon and beyond. 
\emph{J. Reine Angew. Math.} \textbf{582}, 173--200. 

\bibitem[Mc5]{McMullen:torsion}
C. McMullen. 
Teichmüller curves in genus two: torsion divisors and ratios of sines.
\emph{Invent. Math.}  \textbf{165}  (2006),  no. 3, 651--672.


\bibitem[Mc6]{McMullen:SL2R}
  C.~McMullen.
  Dynamics of $SL_2(\reals)$  over moduli space in genus two. 
 \emph{Ann. of Math.} (2)  \textbf{165}  (2007),  no. 2, 397--456.

\bibitem[M\"o1]{Moeller:Variations}
M. M\"oller. 
Variations of Hodge structures of a Teichmüller curve. 
\emph{J. Amer. Math. Soc.} \textbf{19} (2006), no. 2, 327--344.

\bibitem[M\"o2]{Moeller:Periodic}
M. M\"oller. 
Periodic points on Veech surfaces and the Mordell-Weil group over a
Teichmüller curve. 
\emph{Invent. Math.} \textbf{165} (2006), no. 3, 633--649.

\bibitem[M\"o3]{Moeller:Finiteness}
M. M\"oller. 
Finiteness results for Teichmüller curves. 
\emph{Ann. Inst. Fourier (Grenoble)} \textbf{58} (2008), no. 1, 63--83. 


\bibitem[M\"o4]{Moeller:Linear} M.~M\"oller,
    Linear manifolds in the moduli space of one-forms,
    \emph{Duke Math. J.} \textbf{144} No. 3 (2008), 447--488.



\bibitem[Mor]{Morris-Ratner}
      D.~Witte~Morris.
      \emph{Ratner's Theorems on Unipotent Flows.}
  Univ.\ of Chicago Press: Chicago, 2005. 
      \texttt{http://arxiv.org/math.DS/0310402}.


\bibitem[Moz]{Mozes:epimorphic}
S. Mozes. 
Epimorphic subgroups and invariant measures.  
\emph{Ergodic Theory Dynam. Systems}  {\bf 15}  (1995),  no. 6, 1207--1210. 


\bibitem[MoSh]{Mozes:Shah}
S. Mozes, N. Shah.
On the space of ergodic invariant measures of unipotent flows. 
\emph{Ergodic Theory Dynam. Systems} {\bf 15} (1995), no. 1, 149~-159. 

\bibitem[MZ]{Morris:Zimmer} R.~Zimmer, D.~Witte Morris. Ergodic
  theory, groups, and geometry.  CBMS Regional Conference Series in
  Mathematics, 109. Published for the Conference Board of the
  Mathematical Sciences, Washington, DC; by the American Mathematical
  Society, Providence, RI, 2008. x+87 pp. ISBN: 978-0-8218-0980-8


\bibitem[NW]{Nguyen:Wright}
D.~Nguyen, A.~Wright. 
Non-Veech surfaces in H\^{}hyp(4) are generic.
\emph{Geom. Funct. Anal.}  \textbf{24}  (2014),  no. 4, 1316--1335.

\bibitem[NZ]{Nevo:Zimmer}
A. Nevo and R. Zimmer. 
Homogeneous Projective Factors for actions of semisimple Lie groups. 
\emph{Invent. Math.}  \textbf{138}  (1999),  no. 2, 229--252.

\bibitem[Ra1]{RatnerRig}
       M.~Ratner.
      Rigidity of horocycle flows.
      \emph{Ann. of Math.} \textbf{115} (1982), 597--614.

\bibitem[Ra2]{RatnerQuot}
      M.~Ratner.
      Factors of horocycle flows.
      \emph{Ergodic Theory Dynam. Systems} \textbf{2} (1982), 465--489.

\bibitem[Ra3]{RatnerJoin}
      M.~Ratner.
      Horocycle flows, joinings and rigidity of products.
      \emph{Ann. of Math.} \textbf{118} (1983), 277--313.

\bibitem[Ra4]{RatnerSolv}
      M.~Ratner.
      Strict measure rigidity for unipotent subgroups of solvable groups.
      \emph{Invent. Math.} \textbf{101} (1990), 449--482.

\bibitem[Ra5]{RatnerSS}
      M.~Ratner.
      On measure rigidity of unipotent subgroups of semisimple groups.
      \emph{Acta Math.} \textbf{165} (1990), 229--309.

\bibitem[Ra6]{RatnerMeas}
      M.~Ratner.
      On Raghunathan's measure conjecture.
      \emph{Ann. of Math.} \textbf{134} (1991), 545--607.

\bibitem[Ra7]{RatnerEqui}
  M.~Ratner.
  Raghunathan's topological conjecture and distributions of unipotent 
flows.
  \emph{Duke Math. J.} \textbf{63} (1991), no.~1, 235--280.


\bibitem[R]{Rokhlin:book}
V.A.~Rokhlin. 
Lectures on the theory of entropy of transformations with invariant
         measures. \emph{Russ. Math. Surv.} \textbf{22:5} (1967),
         1--54. 

\bibitem[Sch]{Schmidt:Etds1981}
K.~Schmidt. 
Amenability, Kazhdan's property T, strong ergodicity and invariant
 means for ergodic group-actions.
\emph{Ergodic Theory Dynamical Systems}  \textbf{1}  (1981),  no. 2, 223--236.








\bibitem[Ve1]{veech:gauss}
W.~Veech, Gauss measures 
for transformations on the space of interval exchange maps, 
\emph{Ann. of Math.}, \textbf{15} (1982),
201--242.

\bibitem[Ve2]{Veech:Siegel}
  W. Veech.
  Siegel measures.
  \emph{Ann. of Math.} \textbf{148} (1998), 895--944.

\bibitem[Wr1]{Wright:field}
A.~Wright.
The field of definition of affine invariant submanifolds of the moduli
space of abelian differentials. 
\emph{Geom. Topol.}  \textbf{18}  (2014),  no. 3, 1323--1341.

\bibitem[Wr2]{Wright:cylinder}
A.~Wright.
Cylinder deformations in orbit closures of translation surfaces.
\emph{Geom. Topol.}  \textbf{19}  (2015),  no. 1, 413--438.

\bibitem[WWF]{WWF}
L.~Wang, X.~Wang, J.~Feng.
Subspace Distance Analysis with Application to
Adaptive Bayesian Algorithm for Face Recognition. \emph{Pat. Rec.}
{\bf 39}(3) (2006) 456-464.


\bibitem[ZK]{Zemlyakov:Katok}
A.N.~Zemlyakov, A.B.~Katok.
\newblock{\em Topological transitivity of billiards in  
polygons.}
\newblock{\em Matem. Zametki} {\bf 18}(2) (1975) pp.~291--300.
English translation in {\em Math. Notes} {\bf 18}(2) (1976) pp.~760--764. 



\bibitem[Zi1]{Zimmer:amenable:reduction}
R. J. Zimmer, Induced and amenable ergodic actions of Lie groups, 
\emph{Ann. Sci. Ecole Norm. Sup.}
\textbf{11} (1978), no. 3, 407~428.


\bibitem[Zi2]{ZimmerBook}
  R.~J.~Zimmer.
  \emph{Ergodic Theory and Semisimple Groups.}
  Birkh\"auser: Boston, 1984.





\bibitem[Zo]{Zorich:survey}
A. Zorich, A. 
Flat Surfaces. 
\emph{Frontiers in Number Theory, Physics, and Geometry. I. Berlin:
  Springer, 2006.} 437--583. 



\end{thebibliography}
\end{document}